\documentclass[12pt,oneside]{book}
\usepackage{geometry}                		
\usepackage{graphicx}									
\usepackage{amssymb}\vspace{0.25cm}
\usepackage{mathtools}
\usepackage{setspace}
\usepackage{tikz-cd}
\usepackage[mathscr]{euscript}
\newcommand{\abs}[1]{\left\vert#1\right\vert}
\usepackage[normalem]{ulem}
\usepackage{manfnt}
\usepackage{pgfplots}
\usepackage{pifont}
\usepackage[safe]{tipa}
\usepackage{marvosym}
\usepackage{scalerel,stackengine}
\usepackage{titlesec}
\usepackage{makeidx}
\usepackage{centernot}
\usepackage{xspace}
\usepackage[noauto]{chappg}
\usepackage[OT2,T1]{fontenc}
\usepackage{hyperref}
\usepackage{amsmath}
\usepackage[amsmath,thmmarks]{ntheorem}
\usepackage{tabularx}
\usepackage{yfonts}
\usepackage{epstopdf}
\usepackage{pdfpages}
\usepackage{accents}
\newcolumntype{Y}{>{\raggedright\arraybackslash}X}

\newcommand\fg{%
\mathfrak{g}
}

\newcommand\gA{%
\textgoth{A}
}

\newcommand\gC{%
\textgoth{C}
}
\newcommand\gf{%
\textgoth{f}
}

\newcommand\ggS{%
\textgoth{S}
}
\newcommand\gT{%
\textgoth{T}
}

\newcommand\Cx{%
\mathbb{C}
}

\newcommand\Q{%
\mathbb{Q}
}
\newcommand\R{%
\mathbb{R}
}

\newcommand\Z{%
\mathbb{Z}
}
\newcommand\invh{%
\Large\textnormal{\textturnh} 
}

\newcommand\sB{%
\mathcal{B}
}

\newcommand\sE{%
\mathcal{E}
}
\newcommand\sF{%
\mathcal{F}
}
\newcommand\sH{%
\mathcal{H}
}
\newcommand\sI{%
\mathcal{I}
}
\newcommand\sL{%
\mathcal{L}
}

\newcommand\sM{%
\mathcal{M}
}

\newcommand\sP{%
\mathcal{P}
}

\newcommand\sS{%
\mathcal{S}
}

\newcommand\gLP{%
\sL - \sP
}

\newcommand\Bo{%
\text{$\sB$o}
}

\newcommand\Reg{%
\text{$\R$e}\hspace{0.05cm}
}
\newcommand\Img{%
\text{Im}\hspace{0.05cm}
}

\newcommand\id{%
\text{id}
}
\newcommand\itr{%
\text{int}\hspace{0.05cm}
}
\newcommand\nth{%
\text{th}
}

\newcommand\td{%
\text{d}
}

\newcommand\tI{%
\text{I}
}
\newcommand\tII{%
\text{II}
}

\newcommand\tJ{%
\text{J}
}

\newcommand\tO{%
\text{O}
}
\newcommand\txo{%
\text{o}
}
\newcommand\Lp{
\text{L}
}

\newcommand\subsetx{%
\ \subset \ 
}

\newcommand\gen{%
\text{gen}\hsx
}
\newcommand\len{%
\text{len}\hsx
}
\newcommand\lext{%
\text{lext}\hsx
}
\newcommand\rext{%
\text{rext}\hsx
}
\newcommand\ent{%
\text{ent}\hsx
}
\newcommand\PW{%
\text{PW}
}

\newcommand\tRH{%
\text{RH}
}

\newcommand\tcf{%
\text{cf}
}

\newcommand\ML{%
\text{$M\hspace{-.07cm}L$}
}

\def\thereforex{\boldsymbol{\text{ }
\leavevmode
\lower0.4ex\hbox{\textbullet}
\kern-.9em\raise1.1ex\hbox{\textbullet}
\kern-0.9em\lower0.4ex\hbox{\textbullet}
\hspace{0.1cm}\thinspace\text{ }}}
\def\thereforez{\boldsymbol{\text{ }
\leavevmode
\lower0.4ex\hbox{$\circ$}
\kern-.9em\raise1.1ex\hbox{$\circ$}
\kern-0.9em\lower0.4ex\hbox{$\circ$}
\hspace{0.1cm}\thinspace\text{ }}}
\newcommand\ra{%
\rightarrow
}
\newcommand\lra{%
\longrightarrow
}

\newcommand\ds{%
\displaystyle
}

\newcommand\erf{%
\text{erf}\hspace{0.05cm}
}

\newcommand\sptx{%
\text{spt}\hspace{0.05cm}
}

\newcommand\sgn{%
\text{sgn}\hspace{0.05cm}
}

\newcommand\un[1]{%
\underline{#1}\xspace
}

\newcommand\restr[2]{%
{#1}|{#2}
}

\newcommand\restrBig[2]{%
{#1}\Big|{#2}
}

\newcommand\hsx{%
\hspace{0.05cm}
}
\newcommand\hsy{%
\hspace{0.03cm}
}
\newcommand\hsz{%
\hspace{0.1cm}
}

\newcommand\thetafcn{
\theta
}

\newcommand\bcdot{%
{\boldsymbol{\ \cdot \ }}
}

\newcommand\bxi{%
{\boldsymbol{\xi}}
}
\newcommand\bzeta{%
{\boldsymbol{\zeta}}
}
\newcommand\bXi{%
{\boldsymbol{\Xi}}
}

\newcommand{\norm}[1]{\left\lVert #1 \right\rVert}

\newcommand{\normx}[1]{\big|\hspace{-.05cm}\big| #1 \big|\hspace{-.05cm}\big|}

\newcommand\limsupx{%
\ov{\text{lim}}
\hspace{0.05cm}
}

\newcommand\liminfx{%
\un{\text{lim}}
\hspace{0.05cm}
}

\newcommand\ov[1]{%
\overline{#1}
}

\newcommand\thickbar[1]{\accentset{\rule{.4em}{0.8pt}}{#1}}

\newcommand\ovs[1]{
\mkern 1.5mu
\overline{\mkern-2.75mu\mbox{$#1$}\raisebox{3.1mm}{}\mkern-1.5mu}
\mkern 1.5mu
}

\newcommand\ovhc[1]{
\mkern 3.5mu
\overline{\mkern-2.5mu\mbox{$#1$}\raisebox{3.2mm}{}\mkern1.5mu}
\mkern 3.5mu
}

\newcommand\chisubab{\chi\raisebox{-.1cm}{$\scaleto{_{a,b}}{6.5pt}$}}
\newcommand\whchisubab{\widehat{\chi}\raisebox{-.1cm}{$\scaleto{_{a,b}}{6.5pt}$}}

\newcommand\chisubminusOneOne{\chi\raisebox{-.1cm}{$\scaleto{_{-1, 1}}{8.5pt}$}}

\newcommand\chisubBminusAAB{\chi\raisebox{-.1cm}{$\scaleto{_{[-A, A]}}{8.5pt}$}}
\newcommand\chisubBminusnnB{\chi\raisebox{-.1cm}{$\scaleto{_{[-n, n]}}{8.5pt}$}}

\definecolor{ultramarine}{RGB}{0, 32, 96}

\definecolor{darkcerulean}{rgb}{0.3, 0.27, 0.49}

\definecolor{forestgreen}{rgb}{0.0, 0.27, 0.13}
\definecolor{forestgreenweb}{rgb}{0.13, 0.55, 0.13}
\definecolor{deepjunglegreen}{rgb}{0.0, 0.29, 0.29}

\definecolor{midnightblue}{rgb}{0.1, 0.1, 0.44}
\definecolor{midnightgreen}{rgb}{0.0, 0.29, 0.33}
\definecolor{myrtle}{rgb}{0.13, 0.26, 0.12}
\definecolor{darkviolet}{rgb}{0.58, 0.0, 0.83}
\definecolor{darkgreen}{rgb}{0.0, 0.2, 0.13}
\definecolor{officegreen}{rgb}{0.0, 0.5, 0.0}

\definecolor{harvardcrimson}{rgb}{0.79, 0.0, 0.09}
\definecolor{hollywoodcerise}{rgb}{0.96, 0.0, 0.63}
\definecolor{debianred}{rgb}{0.84, 0.04, 0.33}
\definecolor{darkturquoise}{rgb}{0.0, 0.81, 0.82}

\definecolor{darktangernine}{rgb}{1.0, 0.35, 0.35}
\definecolor{aureolin}{rgb}{0.99, 0.93, 0.0}
\definecolor{canaryyellow}{rgb}{1.0, 0.94, 0.0}
\definecolor{amber}{rgb}{1.0, 0.75, 0.0}

\definecolor{urobilin}{rgb}{0.88, 0.68, 0.13}
\definecolor{uscgold}{rgb}{1.0, 0.8, 0.0}

\makeatletter

\def\env@sqcases{%
  \let\@ifnextchar\new@ifnextchar
  \left\lbrack
  \def\arraystretch{1.2}%
  \array{@{}l@{\quad}l@{}}%
}
\makeatother

\usepackage{scalerel,stackengine}
\stackMath
\newcommand\reallywidehat[1]{%
\savestack{\tmpbox}{\stretchto{%
  \scaleto{%
    \scalerel*[\widthof{\ensuremath{#1}}]{\kern-.6pt\bigwedge\kern-.6pt}%
    {\rule[-\textheight/2]{1ex}{\textheight}}
  }{\textheight}%
}{0.5ex}}%
\stackon[1pt]{#1}{\tmpbox}%
}
\makeatletter
\DeclareRobustCommand\widecheck[1]{{\mathpalette\@widecheck{#1}}}
\def\@widecheck#1#2{%
    \setbox\z@\hbox{\m@th$#1#2$}%
    \setbox\tw@\hbox{\m@th$#1%
       \widehat{%
          \vrule\@width\z@\@height\ht\z@
          \vrule\@height\z@\@width\wd\z@}$}%
    \dp\tw@-\ht\z@
    \@tempdima\ht\z@ \advance\@tempdima2\ht\tw@ \divide\@tempdima\thr@@
    \setbox\tw@\hbox{%
       \raise\@tempdima\hbox{\scalebox{1}[-1]{\lower\@tempdima\box
\tw@}}}%
    {\ooalign{\box\tw@ \cr \box\z@}}}
\makeatother

\newtheoremstyle{xx}
  {4pt}
  {0pt}
  {\upshape}
  {\bfseries}
  {}
  { }
  {}
  
\makeatletter 
 \newtheoremstyle{myu}%
  {\upshape\item[ \indent\indent\bf\underline{\theorem@headerfont ##2:}]}%
\makeatother
\makeatletter 
 \newtheoremstyle{myn}%
  {\item[\hskip\labelsep \ \bf ##1 \theorem@headerfont ##2.]}%
\makeatother
\theoremstyle{myn}
\newtheorem{theoremn}{Theorem} 
\theoremstyle{myu}
{\upshape}

\newtheorem{x}[theoremn]{}
\newtheorem{rf}[theoremn]{}

\DeclareFontFamily{U}{wncy}{}
\DeclareFontShape{U}{wncy}{m}{n}{<->wncyr10}{}
\DeclareSymbolFont{mcy}{U}{wncy}{m}{n}
\DeclareMathSymbol{\Sh}{\mathord}{mcy}{"58} 

\title{ZEROS}
\author{Garth Warner\\
Department of Mathematics\\
University of Washington}
\date{}									

\titleformat{\chapter}[display]
{\normalfont\filcenter\huge\bfseries}{}{0pt}{\large}

\titleformat{\chapter}[display]
{\normalfont\filcenter\huge\bfseries}{}{0pt}{\large}

\usepackage[OT2,OT1]{fontenc} 
\newcommand\cyr
{
\renewcommand\rmdefault{wncyr} 
\renewcommand\sfdefault{wncyss} 
\renewcommand\encodingdefault{OT2} 
\normalfont
\selectfont
}

\DeclareTextFontCommand{\textcyr}{\cyr}

\makeindex 
\begin{document}

\maketitle                              

\titlespacing*{\chapter}{0pt}{-50pt}{40pt}
\setlength{\parskip}{0.1em}
\renewcommand{\thepage}{\roman{page}}
\begingroup
\fontsize{11pt}{11pt}\selectfont

\endgroup 

\[
\textbf{ABSTRACT}
\]
\\

The purpose of this book is two fold.
\\[-.5cm]

\qquad
(1) \  
To give a systematic account of classical ``zero theory'' as developed by 
Jensen, 
P\'olya, 
Titchmarsh, 
Cartwright, 
Levinson 
and others. 
\\[-.5cm]

\qquad
(2) \  
To set forth developments of a more recent nature with a view toward their possible application to the Riemann Hypothesis.
\\[2cm]

\[
\textbf{ACKNOWLEDGEMENTS}
\]

\textbullet \quad
My thanks to Judith Clare for her superb job of difficult technical typing and her meticulous proofreading.
\\[-.5cm]

\textbullet \quad 
My thanks to Saundra Martin for her expert help with library matters.
\\[-.5cm]

\textbullet \quad 
My thanks to M. Scott Osborne for his mathematical input and willingness to discuss technicalities.
\\[-.5cm]

\textbullet \quad 
My thanks to David Clark for his rendering the original transcript into AMS-LaTeX.  
\\[-.5cm]

\newpage

\centerline{\textbf{\LARGE CONTENTS}}
\vspace{0.5cm}

\allowdisplaybreaks
\begin{align*}
\S1.  \qquad &\text{INFINITE PRODUCTS} 
\\[10pt]
\S2.  \qquad &\text{ORDER} 
\\[10pt]
\S3.  \qquad &\text{TYPE}  
\\[10pt]
\S4.  \qquad &\text{CONVERGENT EXPONENT}
\\[10pt]
\S5.  \qquad &\text{CANONICAL PRODUCTS}
\\[10pt]
\S6.  \qquad &\text{EXPONENTIAL FACTORS }
\\[10pt]
\S7.  \qquad &\text{REPRESENTATION THEORY}
\\[10pt]
\S8.  \qquad &\text{ZEROS}
\\[10pt]
\S9.  \qquad &\text{JENSEN CIRCLES}
\\[10pt]
\S10.  \qquad &\text{CLASSES OF ENTIRE FUNCTIONS}
\\[10pt]
\S11.  \qquad &\text{DERIVATIVES}
\\[10pt]
\S12.  \qquad &\text{JENSEN POLYNOMIALS}
\\[10pt]
\S13.  \qquad &\text{CHARACTERIZATIONS}
\\[10pt]
\S14.  \qquad &\text{SHIFTED SUMS}
\\[10pt]
\S15.  \qquad &\text{JENSEN CIRCLES (BIS)}
\\[10pt]
\S16.  \qquad &\text{STURM CHAINS}
\\[10pt]
\S17.  \qquad &\text{EXPONENTIAL TYPE}
\\[10pt]
\S18.  \qquad &\text{THE BOREL TRANSFORM}
\\[10pt]
\S19.  \qquad &\text{THE INDICATOR FUNCTION}
\\[10pt]
\S20.  \qquad &\text{DUALITY}
\\[10pt]
\S21.  \qquad &\text{FOURIER TRANSFORMS}
\\[10pt]
\S22.  \qquad &\text{PALEY-WIENER}   
\\[10pt]
\S23.  \qquad &\text{DISTRIBUTION FUNCTIONS}   
\\[10pt]
\S24.  \qquad &\text{CHARACTERISTIC FUNCTIONS}
\\[10pt]
\S25.  \qquad &\text{HOLOMORPHIC CHARACTERISTIC FUNCTIONS}
\\[10pt]
\S26.  \qquad &\text{ENTIRE CHARACTERISTIC FUNCTIONS}
\\[10pt]
\S27.  \qquad &\text{ZERO THEORY: BERNSTEIN FUNCTIONS}
\\[10pt]
\S28.  \qquad &\text{ZERO THEORY: PALEY-WIENER FUNCTIONS}
\\[10pt]
\S29.  \qquad &\text{INTERMEZZO}
\\[10pt]
\S30.  \qquad &\text{TRANSFORM THEORY: \ JUNIOR GRADE}
\\[10pt]
\S31.  \qquad &\text{TRANSFORM THEORY: \ SENIOR GRADE}
\\[10pt]
\S32.  \qquad &\text{APPLICATION OF INTERPOLATION}
\\[10pt]
\S33.  \qquad &\text{ZEROS OF $W_{A,\alpha}$}
\\[10pt]
\S34.  \qquad &\text{ZEROS OF $f_A$}
\\[10pt]
\S35.  \qquad &\text{MISCELLANEA}
\\[10pt]
\S36.  \qquad &\text{LOCATION, LOCATION, LOCATION}
\\[10pt]
\S37.  \qquad &\text{THE $\sF_0$-CLASS}
\\[10pt]
\S38.  \qquad &\text{$\bzeta$, $\bxi$, AND $\bXi$}
\\[10pt]
\S39.  \qquad &\text{THE de BRUIJN-NEWMAN CONSTANT}
\\[10pt]
\S40.  \qquad &\text{TOTAL POSITIVITY}
\\[10pt]
\S41.  \qquad &\text{CHANGE OF VARIABLE} 
\\[10pt]
\S42.  \qquad &\text{$D(n,2)$}  
\\[10pt]
\S43.  \qquad &\text{POSITIVE QUADRATIC FORMS}
\\[10pt]
\S44.  \qquad &\text{ONE EQUIVALENCE}
\\[10pt]
\S45.  \qquad &\text{SUGGESTED READING}
\\[10pt]
\end{align*}

\pagenumbering{bychapter}
\setcounter{chapter}{0}
\chapter{
$\boldsymbol{\S}$\textbf{1}.\quad  INFINITE PRODUCTS}
\setlength\parindent{2em}
\setcounter{theoremn}{0}
\renewcommand{\thepage}{\S1-\arabic{page}}

\vspace{-.25cm}
\qquad
Let $\{z_n : n = 1, 2, \ldots \}$ be a sequence of complex numbers.
\\[-.25cm]

\begin{x}{\small\bf DEFINITION} \ 
The infinite product 

\[
\prod\limits_{n=1}^\infty \
(1 + z_n)
\]
is \un{convergent} if the following conditions are satisfied.
\\[-.25cm]

\qquad \textbullet \quad 
The partial products

\[
\prod\limits_{n=1}^N \
(1 + z_n)
\]
approach a finite limit as $N \ra \infty$.
\\[-.5cm]

\qquad \textbullet \quad 
From some point on, say $n > N_0$, $z_n \neq -1$, and then

\[
\lim\limits_{N \ra \infty} \ 
\prod\limits_{N_0 + 1}^N \
(1 + z_n)
\ \neq \ 
0.
\]
\\[-1cm]

[Note: \ 
The infinite product
\[
\prod\limits_{n=1}^\infty \
(1 + z_n)
\]
is \un{divergent} if it is not convergent.]
\\[-.25cm]
\end{x}

{\small\bf \un{N.B}} \ 
The convergence  of 
\[
\prod\limits_{n=1}^\infty \
(1 + z_n)
\]
implies that $1 + z_n \ra 1$, hence that $z_n \ra 0$.
\\[-.25cm]

\begin{x}{\small\bf REMARK} \ 
It can happen that
\[
\prod\limits_{n=1}^\infty \
(1 + z_n)
\ = \ 
0
\]
but only when at least one factor is zero.
\\[-.25cm]
\end{x}

\begin{x}{\small\bf EXAMPLE} \ 
On the one hand, 
\[
\prod\limits_{n=2}^\infty \
\bigg(
1 - \frac{1}{n^2}
\bigg)
\ = \ 
\frac{1}{2}.
\]
\allowdisplaybreaks\begin{align*}
\Big[ 
\hspace{2cm} 
\prod\limits_{n=2}^N \
\bigg(
1 - \frac{1}{n^2}
\bigg)\ 
&=\ 
\prod\limits_{n=2}^N \ 
\frac{n^2 - 1}{n^2}
\\[15pt]
&=\ 
\prod\limits_{n=2}^N \ 
\frac{(n-1) \cdot (n+1)}{n^2} 
\\[15pt]
&=\ 
\frac{1 \cdot 3}{2^2} 
\hsx \cdot \hsx 
\frac{2 \cdot 4}{3^2} 
\hsx \cdot \hsx 
\frac{3 \cdot 5}{4^2} 
\cdots
\frac{(N-3) \cdot N}{(N - 1)^2} 
\hsx \cdot \hsx 
\frac{(N-1) \cdot (N+1)}{N^2} 
\\[15pt]
&=\ 
\frac{1}{2} 
\hsx \cdot \hsx 
1 
\hsx \cdot \hsx 
1 
\cdots
1
\hsx \cdot \hsx 
\frac{N + 1}{N} 
\\[15pt]
&\ra \ 
\frac{1}{2} 
\qquad (N \ra \infty).\Big]
\end{align*}
On the other hand, 
\[
\prod\limits_{n=1}^\infty \
\bigg(
1 - \frac{1}{n^2}
\bigg)
\ = \ 
0.
\]
\\[-.75cm]
\end{x}

\begin{x}{\small\bf EXAMPLE} \ 
For all $N_0 > 1$, 

\[
\lim\limits_{N \ra \infty} \ 
\prod\limits_{N_0 + 1}^N \
\bigg(
1 - \frac{1}{n}
\bigg)
\ = \ 
0.
\]
Therefore the infinite product

\[
\prod\limits_{n = 2}^\infty \
\bigg(
1 - \frac{1}{n}
\bigg)
\]
is divergent.
\\[-.25cm]
\end{x}

Turning to the theory, we shall first consider the case of the real numbers.  
\\[-.25cm]

\begin{spacing}{1.75}
\begin{x}{\small\bf LEMMA} \ 
If $\{a_n : n = 1, 2, \ldots\}$ is a sequence of nonnegative real numbers, then \ 
$\ds 
\prod\limits_{n = 1}^\infty \
(1 + a_n)
$
is convergent iff \ 
$\ds
\sum\limits_{n = 1}^\infty \ a_n$
is convergent.
\\[-.5cm]

PROOF \ 
In fact, $\forall \ N$, 

\[
a_1 + a_2 + \cdots + a_N 
\ \leq \ 
\prod\limits_{n = 1}^N \ 
(1 + a_n)
\ \leq \ 
\exp (a_1 + a_2 + \cdots + a_N ).
\]
\\[-1.5cm]
\end{x}
\end{spacing}


\begin{x}{\small\bf EXAMPLE} \ 
The infinite product

\[
\prod\limits_{n = 1}^\infty \
\bigg(1 + \frac{1}{n^p}\bigg)
\]
is convergent for $p > 1$ and divergent for $p \leq 1$.
\\[-.25cm]
\end{x}

\begin{spacing}{1.75}
\begin{x}{\small\bf LEMMA} \ 
If $\{a_n : n = 1, 2, \ldots\}$ is a sequence of nonnegative real numbers, then \ 
$\ds 
\prod\limits_{n = 1}^\infty \
(1 - a_n)
$
is convergent iff \ 
$\ds
\sum\limits_{n = 1}^\infty \ a_n$
is convergent.
\end{x}
\end{spacing} 

PROOF \ 
If $a_n$ does not tend to 0, then both the product and the series are divergent, 
so there is no loss of generality in assuming from the beginning that 
$\ds
a_n < \frac{1}{2}$ 
$\ds \bigg(\implies 1 - a_n > \frac{1}{2}\bigg)$.
\\

\qquad \textbullet \quad
Suppose that 
$\ds
\prod\limits_{n = 1}^\infty \ 
(1 - a_n)
$
is convergent $-$then the partial products
\[
\prod\limits_{n = 1}^N \ 
(1 - a_n)
\]
constitute a monotone decreasing sequence with a positive limit $L$: $\forall \ N$, 
\[
\prod\limits_{n = 1}^N \ 
(1 - a_n)
\ \geq \ 
L 
\ > \ 
0.
\]
But

\[
1 + a_n
\ \leq \ 
\frac{1}{1 - a_n},
\]
thus 

\[
\prod\limits_{n = 1}^N \ 
(1 + a_n)
\ \leq \ 
\prod\limits_{n = 1}^N \ 
\frac{1}{1 - a_n}
\ \leq \ 
\frac{1}{L}.
\]
Since the partial products

\[
\prod\limits_{n = 1}^N \ 
(1 + a_n)
\]
constitute a monotone increasing sequence, it follows that 
$\ds
\prod\limits_{n = 1}^\infty \ 
(1 + a_n)
$
is convergent, 
hence the same is true of 
$\ds
\sum\limits_{n = 1}^\infty \ a_n
$ 
(cf. 1.5).
\\

\qquad \textbullet \quad
Suppose that 
$\ds
\sum\limits_{n = 1}^\infty \ 
a_n
$
is convergent $-$then 
$\ds
\sum\limits_{n = 1}^\infty \ 
2 a_n
$
is convergent, thus
$\ds
\prod\limits_{n = 1}^\infty \ 
(1 + 2 a_n)
$
is convergent (cf. 1.5), so there exists $K > 0$ such that $\forall \ N$, 

\[
\prod\limits_{n = 1}^N \ 
(1 + 2 a_n)
\ \leq \ 
K.
\]
But
\[
0
\ \leq \ 
a_n 
\ < \ 
\frac{1}{2}
\ \implies \ 
1 - a_n 
\ \geq \ 
\frac{1}{1 + 2 \hsy a_n}
\]

\qquad 
$\implies$

\[
\prod\limits_{n = 1}^N \ 
(1 -  a_n)
\ \geq \ 
\prod\limits_{n = 1}^N \ \frac{1}{1 + 2 a_n} 
\ \geq \ 
\frac{1}{K}
\ > \ 
0.
\]
And

\[
\prod\limits_{n = 1}^\infty \ 
(1 -  a_n)
\]
is monotone increasing.
\\

\begin{x}{\small\bf EXAMPLE} \ 
The infinite product

\[
\prod\limits_{n = 1}^\infty \
\bigg(1 - \frac{1}{n^p}\bigg)
\]
is convergent for $p > 1$ and divergent for $p \leq 1$.
\\[-.25cm]
\end{x}

\begin{spacing}{1.75}
\begin{x}{\small\bf LEMMA} \ 
Let $\{a_n : n = 1, 2, \ldots\}$ be a sequence of real numbers.  
\\[-.25cm]

\noindent
Assume: \ 
$
\begin{cases}
\quad
\ds 
\sum\limits_{n = 1}^\infty \
a_n
\\[18pt]
\quad
\ds 
\sum\limits_{n = 1}^\infty \ a_n^2
\end{cases}
$
are convergent $-$then \ 
$\ds 
\prod\limits_{n = 1}^\infty \
(1 + a_n)
$
\ 
is convergent.

\end{x}
\end{spacing}

PROOF \ 
Supposing as we may that $\forall \ n$, 
$\ds \abs{a_n} < \frac{1}{2}$, 
note that 

\[
\log (1 + a_n) 
\ = \ 
a_n + \tO(a_n^2).
\]
Therefore the series

\[
\sum\limits_{n = 1}^\infty \
\log (1 + a_n)
\]
is convergent to $L$, say, hence
\allowdisplaybreaks\begin{align*}
\prod\limits_{n = 1}^N \
(1 + a_n) \ 
&=\ 
\exp 
\Big( \log \hsx 
\prod\limits_{n = 1}^N \
(1 + a_n) 
\Big)
\\[15pt]
&=\ 
\exp \hsy 
\bigg(
\sum\limits_{n = 1}^N \
\log (1 + a_n)
\bigg)
\\[15pt]
&
\ra 
e^L 
\qquad (N \ra \infty)
\\[15pt]
&\neq\ 
0.
\end{align*}
\\[-1.25cm]

\begin{x}{\small\bf EXAMPLE} \ 
The infinite product
\[
\prod\limits_{n = 1}^\infty \
\bigg(
1 + \frac{(-1)^{n-1}}{n}
\bigg)
\]
is convergent.
\\[-.25cm]
\end{x}

\begin{x}{\small\bf LEMMA} \ 
Let $\{a_n : n = 1, 2, \ldots\}$ be a sequence of real numbers.  
\\[-.25cm]

\noindent
Assume: \ 
$
\begin{cases}
\quad
\sum\limits_{n = 1}^\infty \
a_n
\quad \text{is convergent}
\\[18pt]
\quad
\sum\limits_{n = 1}^\infty \
a_n^2
\quad \text{is divergent}
\end{cases}
$
$-$then 
\ 
$\ds 
\prod\limits_{n = 1}^\infty \
(1 + a_n)
$
\ 
is divergent.
\\


[Use the inequality

\[
x \hsx - \hsx  \log (1 + x) \ > \quad
\begin{cases}
\ \ds
\frac{\ds\frac{x^2}{2}}{(1 + x)}
\hspace{01cm} 
(x > 0)
\\[15pt]
\hspace{0.5cm}  
\ds
\frac{x^2}{2} 
\hspace{1.5cm} 
(0 > x > -1)
\end{cases}
.]
\]
\\[-.75cm]
\end{x}

\begin{x}{\small\bf EXAMPLE} \ 
The infinite product
\[
\prod\limits_{n = 1}^\infty \
\bigg(
1 + \frac{(-1)^{n-1}}{\sqrt{n}}
\bigg)
\]
is divergent.
\\[-.25cm]
\end{x}

\begin{x}{\small\bf REMARK} \ 
It can happen that both 
$\ds 
\ 
\sum\limits_{n = 1}^\infty \
a_n
\ 
$
and
$\ds 
\ 
\sum\limits_{n = 1}^\infty \
a_n^2
\ 
$
are divergent, yet 
\ 
$\ds 
\prod\limits_{n = 1}^\infty \ (1 + a_n)
$
\ 
is convergent.
\\

[Consider

\[
\bigg(
1 - \frac{1}{\sqrt{2}}
\bigg)
\hsy
\bigg(
1 + \frac{1}{\sqrt{2}} + \frac{1}{2}
\bigg)
\hsy
\bigg(
1 - \frac{1}{\sqrt{3}}
\bigg)
\hsy
\bigg(
1 + \frac{1}{\sqrt{3}} + \frac{1}{3}
\bigg)
\cdots
.]
\]
\\[-1cm]
\end{x}

Let $\{z_n : n = 1, 2, \ldots\}$ be a sequence of complex numbers.
\\[-.25cm]

\begin{x}{\small\bf CRITERION} \ 
The infinite product

\[
\prod\limits_{n = 1}^\infty \ (1 + z_n)
\]
is convergent iff $\forall \ \varepsilon > 0$, 
$\exists \ N(\varepsilon)$ such that 
$\forall \ N > N(\varepsilon)$ and every $k \geq 1$, 

\[
\abs{(1 + z_{N + 1}) \cdots (1 + z_{N + k})  - 1}
\ < \ 
\varepsilon.
\]
\\[-1cm]

PROOF  \ 
\\[-.5cm]

\qquad \textbullet \quad
\un{Necessity} \ 
Choose $N_0$ per 1.1, put

\[
P_N 
\ = \ 
\prod\limits_{N_0 + 1}^N \ (1 + z_n)
\]
and fix $C > 0$: 

\[
\forall \ N > N_0, \ 
\abs{P_N} > C.
\]
Since $\{P_N\}$ is a Cauchy sequence, by taking $N_0$ large enough, one can arrange that
$\forall \ N > N_0$ and every $k \geq 1$, 

\[
\abs{P_{N + k} - P_N} \ < \ C \hsy \varepsilon.
\]
Therefore 
\[
\abs{\frac{P_{N + k}}{P_N} - 1 }
\ < \ 
\frac{C}{P_N} \hsy \varepsilon
\ < \ 
\varepsilon
\]
or still, 
\[
\abs
{
\big(
1 + z_{N+1}
\big)
\cdots
\big(
1 + z_{N+k}
\big)
-1
}
\ < \ 
\varepsilon.
\]
\\[-.75cm]

\qquad \textbullet \quad
\un{Sufficiency} \ 
First take 
$
\ds
\varepsilon = \frac{1}{2}
$, 
hence 
$
\ds
\forall \ N > N \bigg(\frac{1}{2}\bigg)
$
and every $k \geq 1$, 

\[
\abs
{
\big(
1 + z_{N+1}
\big)
\cdots
\big(
1 + z_{N+k}
\big)
-1
}
\ < \ 
\frac{1}{2}.
\]
So, for all 
$
\ds
n > N_0 \equiv N \bigg(\frac{1}{2}\bigg) + 1
$, 
$z_n \neq -1$, and if 

\[
\lim\limits_{N \ra \infty} \ 
\prod\limits_{N_0 + 1}^N \ (1 + z_n)
\]
exists, it cannot be zero since 

\[
\frac{1}{2} 
\ < \ 
\bigg| 
\prod\limits_{N_0 + 1}^N \ (1 + z_n)
\hsx
\bigg| 
\ < \ 
\frac{3}{2}.
\]
Take now $\varepsilon > 0$ and choose 
$
\ds
N \bigg(\frac{\varepsilon}{2}\bigg)
 > 
 N \bigg(\frac{1}{2}\bigg)
$
$-$then 
$
\ds
\forall \ N > N \bigg(\frac{\varepsilon}{2}\bigg)
$ 
and every $k \geq 1$, 

\[
\abs
{
\big(
1 + z_{N+1}
\big)
\cdots
\big(
1 + z_{N+k}
\big)
-1
}
\ < \ 
\frac{\varepsilon}{2}, 
\]
from which 

\[
\abs{\frac{P_{N + k}}{P_N} - 1 }
\ < \ 
\frac{\varepsilon}{2}
\]
or still, 
\allowdisplaybreaks\begin{align*}
\abs{P_{N + k} - P_N} \ 
&<\ 
\abs{P_N} \hsy \frac{\varepsilon}{2}
\\[11pt]
&<\ 
\bigg(\frac{3}{2} \bigg) \hsy \frac{\varepsilon}{2}
\\[11pt]
&=\ 
\frac{3}{4} \hsy \varepsilon
\\[11pt]
&<\ 
\varepsilon.
\end{align*}
Therefore

\[
\bigg\{\ 
\prod\limits_{N_0 + 1}^N \ 
 (1 + z_n)
\bigg\}
\]
is a Cauchy sequence, thus is convergent.
\\[-.25cm]
\end{x}

\begin{x}{\small\bf DEFINITION} \ 
The infinite product

\[
\prod\limits_{n = 1}^\infty \ (1 + z_n)
\]
is \un{absolutely convergent} if the infinite product

\[
\prod\limits_{n = 1}^\infty \ (1 + \abs{z_n})
\]
is convergent.
\\[-.25cm]
\end{x}

\begin{x}{\small\bf LEMMA} \ 
An absolutely convergent infinite product
\[
\prod\limits_{n = 1}^\infty \ (1 + z_n)
\]
is convergent.
\\[-.5cm]

PROOF \ 
One has only to note that

\[
\abs
{
\big(
1 + z_{N+1}
\big)
\cdots
\big(
1 + z_{N+k}
\big)
-1
}
\ \leq \ 
\big(
1 + \abs{z_{N+1}}
\big)
\cdots
\big(
1 + \abs{z_{N+k}}
\big)
-1
\]
and the apply 1.14.
\\[-.25cm]
\end{x}

\begin{x}{\small\bf REMARK} \ 
In view of 1.5, 
\ 
$
\ds
\prod\limits_{n = 1}^\infty \ (1 + \abs{z_n})
$
\ 
is convergent iff 
\ 
$
\ds
\sum\limits_{n = 1}^\infty \ \abs{z_n})
$
\ 
is convergent.
\\[-.25cm]
\end{x}

\begin{x}{\small\bf EXAMPLE} \ 
The infinite product

\[
\prod\limits_{n = 1}^\infty \ \frac{\sin(z / n)}{(z / n)}
\]
is absolutely convergent for all finite $z$ 
(with the usual convention at $z = 0$).
\\[-.5cm]

[Observe that

\[
\frac{\sin(z / n)}{(z / n)} 
- 1
\ = \ 
\tO_z
\bigg(
\frac{1}{n^2}
\bigg)
\qquad (n \ra \infty).]
\]
\\[-1.75cm]
\end{x}

\begin{spacing}{2}
It is initially tempting to think that absolute convergence should be the demand that 
$
\ds
\prod\limits_{n = 1}^\infty \ \abs{1 + z_n}
$
is convergent but this will not do since then it is no longer true that ``absolute convergence'' implies convergence.
\\[-1.5cm]
\end{spacing}

\begin{x}{\small\bf EXAMPLE} \ 
The infinite product

\[
\prod\limits_{n = 1}^\infty \ 
\bigg(
1 + \frac{\sqrt{-1}}{n}
\bigg)
\]
is divergent but the infinite product 
\[
\prod\limits_{n = 1}^\infty \ 
\abs
{
1 + \frac{\sqrt{-1}}{n}
}
\]
is convergent.
\\[-.25cm]
\end{x}

\begin{x}{\small\bf LEMMA} \ 
If the infinite product

\[
\prod\limits_{n = 1}^\infty \ (1 + z_n)
\]
is absolutely convergent, then it can be rearranged at will without changing its value, 
which is thus independent of the order of the factors.
\\[-.25cm]
\end{x}

\begin{x}{\small\bf EXAMPLE} \ 
The infinite product

\[
P 
\ = \ 
\bigg(
1 - \frac{1}{2}
\bigg)
\hsy 
\bigg(
1 + \frac{1}{3}
\bigg)
\hsy 
\bigg(
1 - \frac{1}{4}
\bigg)
\hsy
\bigg(
1 + \frac{1}{5}
\bigg)
\hsy
\bigg(
1 - \frac{1}{6}
\bigg)
\cdots
\]
is convergent (cf. 1.10) but not absolutely convergent and has value $1/2$, 
while the rearrangement

\[
Q 
\ = \ 
\bigg(
1 - \frac{1}{2}
\bigg)
\hsy 
\bigg(
1 - \frac{1}{4}
\bigg)
\hsy 
\bigg(
1 + \frac{1}{3}
\bigg)
\hsy
\bigg(
1 - \frac{1}{6}
\bigg)
\hsy
\bigg(
1 - \frac{1}{8}
\bigg)
\hsy
\bigg(
1 + \frac{1}{5}
\bigg)
\cdots
\]
has value $1/ 2 \sqrt{2}$.
\\[-.25cm]
\end{x}

\begin{x}{\small\bf EXAMPLE} \ 
Fix a complex number 
$q: \abs{q} < 1$.  
Introduce the absolutely convergent infinite products

\[
\begin{cases}
\ds \ 
q_0 
\ = \ 
\prod\limits_{n = 1}^\infty \ 
\big(
1 - q^{2 n}
\big), 
\hspace{.9cm}
q_1 
\ = \ 
\prod\limits_{n = 1}^\infty \ 
\big(
1 + q^{2 n}
\big), 
\\[26pt]
\ds \ 
q_2 
\ = \ 
\prod\limits_{n = 1}^\infty \ 
\big(
1 + q^{2 n - 1}
\big), 
\hspace{0.5cm}
q_3 
\ = \ 
\prod\limits_{n = 1}^\infty \ 
\big(
1 - q^{2 n - 1}
\big)
\end{cases}
.
\]
\\[-.75cm]

\noindent
Then

\[
\begin{cases}
\ds \ 
q_0 \hsy q_3 
\ = \ 
\prod\limits_{n = 1}^\infty \ 
\big(
1 - q^n
\big)
\\[26pt]
\ds \ 
q_1 \hsy q_2 
\ = \ 
\prod\limits_{n = 1}^\infty \ 
\big(
1 + q^n
\big)
\end{cases}.
\]
In addition, 
\allowdisplaybreaks\begin{align*}
q_0 \ 
&=\ 
\prod\limits_{n = 1}^\infty \ \big(1 - q^{2 n}\big)
\\[15pt]
&=\ 
\prod\limits_{m = 1}^\infty \
\big(1 - q^{4 m}\big)
\prod\limits_{m = 1}^\infty \
\big(1 - q^{4 m - 2}\big)
\\[15pt]
&=\ 
\prod\limits_{m = 1}^\infty \
\bigg(1 - q^{2 m}\bigg)
\prod\limits_{m = 1}^\infty \
\bigg(1 + q^{2m}\bigg)
\prod\limits_{m = 1}^\infty \
\bigg(1 + q^{2 m - 1}\bigg)
\prod\limits_{m = 1}^\infty \
\bigg(1 - q^{2 m - 1}\bigg)
\\[15pt]
&=\ 
q_0 \hsy q_1 \hsy q_2 \hsy q_3,
\end{align*}
so
\[
q_1 \hsy q_2 \hsy q_3
\ = \ 
1.
\]
\\[-.75cm]
\end{x}

\begin{x}{\small\bf EXAMPLE} \ 
The infinite product

\[
\prod\limits_{n = 1}^\infty \ 
\bigg(
1 - \frac{z^2}{n^2}
\bigg)
\]
is absolutely convergent and has value

\[
\frac{\sin \hsy \pi \hsy z}{\pi \hsy z}.
\]
Consider now the infinite product

\[
\bigg(
1 - z
\bigg)
\bigg(
1 + z
\bigg)
\bigg(
1 - \frac{z}{2}
\bigg)
\bigg(
1 + \frac{z}{2}
\bigg)
\cdots \hsx .
\]
Officially, therefore

\[
z_1 = -z, 
\quad
z_2 = z, 
\quad 
z_3 = -\frac{z}{2}, 
\quad 
z_4 = \frac{z}{2}, 
\ldots, 
\]
and the associated series of absolute values is

\[
\abs{z} + \abs{z} + \frac{\abs{z}}{2} + \frac{\abs{z}}{2} + \cdots, 
\]
which is not convergent if $z \neq 0$.  
Nevertheless, our infinite product is convergent and has value

\[
\frac{\sin \hsy \pi \hsy z}{\pi \hsy z},
\]
as can be seen by looking at the sequence of partial products.  
To correct for the failure of absolute convergence, form instead the infinite product

\[
\big\{
(1 - z) e^z
\big\}
\big\{
(1 + z) e^{-z}
\big\}
\big\{
\Big(1 - \frac{z}{2}\Big) e^{z/2}
\big\}
\big\{
\Big(1 + \frac{z}{2}\Big) e^{-z/2}
\big\}
\cdots \hsx .
\]
To place it into the 
$
\prod\limits_{n = 1}^\infty \ 
(1 + z_n)
$
format, note that the 
$(2n - 1)^\nth$ term is  

\[
\bigg(
1 - \frac{z}{n}
\bigg)
\hsx
e^{z / n} - 1
\]
and the $(2 n)^\nth$ term is 

\[
\bigg(
1 + \frac{z}{n}
\bigg)
\hsx
e^{- z / n} - 1.
\]
But

\[
\bigg(
1 \mp \frac{z}{n}
\bigg)
\hsx
e^{\pm z / n}
\ = \ 
1 + \tO_z
\bigg(
\frac{1}{n^2}
\bigg)
\qquad 
(n \ra \infty).
\]
Since

\[
1 + 1 + \frac{1}{2^2} + \frac{1}{2^2} + \frac{1}{3^2} + \frac{1}{3^2} + \cdots
\]
is convergent, it follows that the foregoing infinite product is absolutely convergent and it too has value

\[
\frac{\sin \hsy \pi \hsy z}{\pi \hsy z}.
\]
\\[-1cm]
\end{x}

\begin{x}{\small\bf EXAMPLE} \ 
The infinite product

\[
\bigg(
1 - z
\bigg)
\bigg(
1 - \frac{z}{2}
\bigg)
\bigg(
1 + z
\bigg)
\bigg(
1 - \frac{z}{3}
\bigg)
\bigg(
1 - \frac{z}{4}
\bigg)
\bigg(
1 + \frac{z}{2}
\bigg)
\cdots
\]
is convergent and has value

\[
\exp (-z\hsy \log 2) \hsx
\frac{\sin \hsy \pi \hsy z}{\pi \hsy z}.
\]
\\[-1cm]

[Judiciously insert the appropriate exponential correction factors.]
\\[-.25cm]
\end{x}

Let $\{f_n (z) : n = 1, 2, \ldots\}$ 
be a sequence of complex valued functions defined on some nonempty subset $S$ of the complex plane.
\\[-.25cm]


\begin{x}{\small\bf DEFINITION} \ 
The infinite product
\[
\prod\limits_{n = 1}^\infty \ 
\big(1 + f_n (z)\big)
\]
is \un{uniformly convergent in $S$} 
if $\forall \ \varepsilon > 0$, $\exists \ N(\varepsilon)$ 
such that $\forall \ N > N(\varepsilon)$ 
and every $k \geq 1$ and every $z \in S$, 

\[
\abs
{
\big(
1 + f_{N + 1} (z)
\big)
\cdots
\big(
1 + f_{N + k} (z)
\big)
-1
}
\ < \ 
\varepsilon.
\]
\\[-1cm]
\end{x}

\begin{spacing}{1.75}
\begin{x}{\small\bf LEMMA} \ 
Suppose that $\forall \ n > 0$, $\exists \ M_n > 0$ such that $\forall \ z \in S$, 
$\abs{f_n (z)} \leq M_n$.  
\\[-1.cm]
\end{x}
\end{spacing}

\noindent
Assume: \ 
$
\ds
\sum\limits_{n = 1}^\infty \ M_n
$
is convergent $-$then the infinite product
\[
\prod\limits_{n = 1}^\infty \ 
\big(1 + f_n (z)\big)
\]
is absolutely and uniformly convergent in $S$.
\\[-.5cm]

PROOF \ 
Absolute convergence is immediate (cf. 1.17):
\[
\sum\limits_{n = 1}^\infty \ 
\abs{f_n (z)} 
\ \leq \ 
\sum\limits_{n = 1}^\infty \ M_n 
\ < \ 
\infty.
\]
As for uniform convergence, the assumption on the $M_n$ implies that 
\ 
$
\ds
\prod\limits_{n = 1}^\infty \ 
(1 + M_n)
$
\ 
is convergent (cf. 1.5).  
On the other hand, 
\allowdisplaybreaks\begin{align*}
\abs{(1 + f_{N + 1} (z)) \cdots (1 + f_{N + k}(z))  - 1} \ 
&\leq\ 
(1 + \abs{f_{N + 1} (z)}) \cdots (1 + \abs{f_{N + k}(z)})  - 1
\\[11pt]
&\leq\ 
(1 + M_{N + 1}) \cdots (1 + M_{N + k}) - 1, 
\end{align*}
thus it remains only to quote 1.14.
\\[-.25cm]

\begin{x}{\small\bf REMARK} \ 
It suffices to assume that 
\ 
$
\ds
\sum\limits_{n = 1}^\infty \ 
\abs{f_n(z)}
$
\ 
is uniformly convergent in $S$ with a bounded sum.
\\[-.25cm]
\end{x}

\begin{x}{\small\bf EXAMPLE} \ 
Take for $S$ a compact subset of 
$\{z : \abs{z} < 1\}$ 
$-$then $S$ is contained in 
$\{z : \abs{z} \leq \delta\}$ 
for some $\delta < 1$, so $\forall \ z \in S$, 

\[
\sum\limits_{n = 1}^\infty \ 
\abs{z^n}
\ \leq \ 
\sum\limits_{n = 1}^\infty \ 
\delta^n 
\ = \ 
\frac{\delta}{1 - \delta} \hsy .
\]
Therefore the infinite product

\[
\prod\limits_{n = 1}^\infty \ 
\big(1 + z^n\big)
\]
is absolutely and uniformly convergent in $S$.
\\[-.25cm]
\end{x}

\begin{x}{\small\bf THEOREM} \ 
Let $f_n (z)$ $(n = 1, 2, \ldots)$ be continuous (holomorphic) in a 
region\footnote[2]{a.k.a.: nonempty open connected subset of $\Cx$} 
$D$ and suppose that the infinite product

\[
\prod\limits_{n = 1}^\infty \ 
\big(1 + f_n (z)\big)
\]
is uniformly convergent on compact subsets of $D$ $-$then the function defined by

\[
\prod\limits_{n = 1}^\infty \ 
\big(1 + f_n (z)\big)
\]
is continuous (holomorphic) in $D$. 
\\[-.25cm]
\end{x}

\begin{x}{\small\bf EXAMPLE} \ 
The infinite product

\[
\prod\limits_{n = 1}^\infty \ 
\Big(
1 + \frac{z}{n}
\Big)
\hsy 
\exp
\Big(-\frac{z}{n}\Big)
\]
is uniformly convergent on compact subsets of $\Cx$ and if as usual, $\Gamma (z)$ stands for  
the gamma function, then

\[
\frac{1}{\Gamma (z)} 
\ = \ 
z \hsy e^{\gamma z} \ 
\prod\limits_{n = 1}^\infty \ 
\Big(
1 + \frac{z}{n}
\Big)
\hsy 
\exp
\Big(-\frac{z}{n}\Big), 
\]
where
\[
\gamma  
\ = \ 
\lim\limits_{n \ra \infty} \ 
\big(H_n - \log n\big)
\]
is Euler's constant.
\\[-.5cm]

[Note: \ 
\[
\Gamma (z)
\ = \ 
\frac{1}{z} \ 
\prod\limits_{n = 1}^\infty \ 
\bigg(
1 + \frac{1}{n}
\bigg)^z
\hsy 
\exp
\bigg(
1 + \frac{z}{n}
\bigg)^{-1}
\]
is meromorphic with simple poles at  0 (residue 1) and the negative integers 
\\[-.25cm]

\noindent
$-n = -1, -2, \ldots$ 
$
\bigg( \text{residue} \ 
\ds
\frac{(-1)^n}{n !} 
\bigg)
$
.]
\\
\end{x}

\[
\text{APPENDIX}
\]

Given a complex number $\tau$ whose imaginary part is positive, 
let $q = \exp \big(\pi \hsy \sqrt{-1} \hsx \tau)$, 
thus $\abs{q} < 1$.
\\[-.25cm]

{\small\bf LEMMA} \ 
The \un{theta functions} 
\\[-.25cm]

\[
\begin{cases}
\ 
\thetafcn_1 (\restr{z}{\hsx \tau}) 
\\[8pt]
\ 
\thetafcn_2 (\restr{z}{\hsx \tau}) 
\\[8pt]
\ 
\thetafcn_3 (\restr{z}{\hsx \tau}) 
\\[8pt]
\ 
\thetafcn_4 (\restr{z}{\hsx \tau}) 
\end{cases}
\]
defined by the series
\\[-.25cm]

\[
\begin{cases}
\ds \ 
\thetafcn_1 (\restr{z}{\hsx \tau}) 
\ = \ 
2 \ 
\sum\limits_{n = 0}^\infty \ 
(-1)^n \hsx 
q^{\big(n + \frac{1}{2}\big)^2} \hsx \sin (2n + 1) z
\\[18pt]
\ds \ 
\thetafcn_2 (\restr{z}{\hsx \tau}) 
\ = \ 
2 \ 
\sum\limits_{n = 0}^\infty \ 
q^{\big(n + \frac{1}{2}\big)^2} \hsx \cos (2n + 1) z
\\[18pt]
\ 
\thetafcn_3 (\restr{z}{\hsx \tau}) 
\ = \ 
1 + 2 \ 
\sum\limits_{n = 1}^\infty \ 
q^{n^2} \hsy \cos 2 n z
\\[18pt]
\ds \ 
\thetafcn_4 (\restr{z}{\hsx \tau}) 
\ = \ 
1 + 2 \ 
\sum\limits_{n = 1}^\infty \ 
(-1)^n \hsx 
q^{n^2} \hsy \cos 2 n z
\end{cases}
\]

\noindent
are entire functions of $z$.
\\[-.25cm]

[The defining series are uniformly convergent on compact subsets of $\Cx$.]
\\

RELATIONS
\\

\qquad \textbullet \quad
$
\ds
\thetafcn_1 \big(\restr{z}{\tau}\big) 
\ = \ 
-\sqrt{-1} \hsx 
\exp \big(\sqrt{-1} \hsx z + \frac{1}{4} \hsx \pi \hsy \sqrt{-1} \hsx \tau \big) 
\hsx 
\thetafcn_4 \big( \restrBig{z + \frac{\pi \tau}{2}}{\hsx \tau} \big)
$
\\[.25cm]

\qquad \textbullet \quad
$
\ds
\thetafcn_2 (\restr{z}{\tau}) 
\ = \ 
\thetafcn_1 \Big( \restrBig{z + \frac{\pi}{2}}{\hsx \tau} \Big)
$
\\[.25cm]

\qquad \textbullet \quad
$
\ds
\thetafcn_3 (\restr{z}{\tau}) 
\ = \ 
\thetafcn_4 \Big( \restrBig{z + \frac{\pi}{2}}{\hsx \tau} \Big).
$
\\

ZEROS \ 
Let $m$, $n$ be integers.
\\

\qquad \textbullet \quad
$
\ds
\thetafcn_1 (\restr{m \hsy \pi + n \hsy \pi \tau}{\hsx \tau})
\ = \ 0
$
\\[.25cm]

\qquad \textbullet \quad
$
\ds
\thetafcn_2 \Big(\restr{\frac{\pi}{2} + m \hsy \pi + n \hsy \pi \tau\hsx}{\tau}\Big)
\ = \ 0
$
\\[.25cm]

\qquad \textbullet \quad
$
\ds
\thetafcn_3 \Big(\restr{\frac{\pi}{2} + \frac{\pi \tau}{2} + m \hsy \pi + n \hsy \pi \tau\hsx}{\tau}\Big)
\ = \ 0
$
\\[.25cm]

\qquad \textbullet \quad
$
\ds
\thetafcn_4 \Big(\restr{\frac{\pi \tau}{2} + m \hsy \pi + n \hsy \pi \tau\hsx}{\tau}\Big)
\ = \ 0.
$
\\

\noindent
These formulas give all the zeros of the respective theta functions and each zero is simple.
\\

PRODUCTS \ 
Let

\[
q_0 
\ = \ 
\prod\limits_{n = 1}^\infty \ 
\big(
1 - q^{2 n}
\big)
\qquad \text{(cf. 1.22).}
\]
\\[-1cm]

\qquad \textbullet \quad
$
\ds
\thetafcn_1 (\restr{z}{\tau}) 
\ = \ 
2 \hsy q_0 \hsy q^{1/4} \ \sin z \ 
\prod\limits_{n = 1}^\infty \ 
\big(
1 - 2\hsy q^{2 n} \cos 2z + q^{4n}
\big)
$
\\[.25cm]

\qquad \textbullet \quad
$
\ds
\thetafcn_2 (\restr{z}{\tau}) 
\ = \ 
2 \hsy q_0 \hsy q^{1/4} \ \cos z \ 
\prod\limits_{n = 1}^\infty \ 
\big(
1 + 2\hsy q^{2 n} \cos 2z + q^{4n}
\big)
$
\\[.25cm]

\qquad \textbullet \quad
$
\ds
\thetafcn_3 (\restr{z}{\tau}) 
\ = \ 
q_0 \ 
\prod\limits_{n = 1}^\infty \ 
\big(
1 + 2\hsy q^{2 n - 1} \cos 2z + q^{4n-2}
\big)
$
\\[.25cm]

\qquad \textbullet \quad
$
\ds
\thetafcn_4 (\restr{z}{\tau}) 
\ = \ 
q_0 \ 
\prod\limits_{n = 1}^\infty \ 
\big(
1 - 2\hsy q^{2 n - 1} \cos 2z + q^{4n-2}
\big).
$
\\[.25cm]

TRANSFORMATIONS \ 
\\[-.25cm]

\qquad \textbullet \quad
$
\ds
\thetafcn_1 (\restr{z}{\tau}) 
\ = \ 
\sqrt{-1} \hsx
\big(- \sqrt{-1} \hsx \tau \big)^{-\frac{1}{2}} \ 
\exp
\Big(
\frac{z^2}{\pi \hsy \sqrt{-1} \hsx \tau}
\Big)
\hsx 
\thetafcn_1 \Big(\restrBig{\frac{z}{\tau}}{\hsx -\tau^{-1}}\Big) 
$
\\[.25cm]

\qquad \textbullet \quad
$
\ds
\thetafcn_2 (\restr{z}{\tau}) 
\ = \ 
\big(- \sqrt{-1} \hsx \tau \big)^{-\frac{1}{2}} \ 
\exp
\Big(
\frac{z^2}{\pi \hsy \sqrt{-1} \hsx \tau}
\Big)
\hsx 
\thetafcn_4 \Big(\restrBig{\frac{z}{\tau}}{\hsx -\tau^{-1}}\Big) 
$
\\[.25cm]

\qquad \textbullet \quad
$
\ds
\thetafcn_3 (\restr{z}{\tau}) 
\ = \ 
\big(- \sqrt{-1} \hsx \tau \big)^{-\frac{1}{2}} \ 
\exp
\Big(
\frac{z^2}{\pi \hsy \sqrt{-1} \hsx \tau}
\Big)
\hsx 
\thetafcn_3 \Big(\restrBig{\frac{z}{\tau}}{\hsx -\tau^{-1}}\Big) 
$
\\[.25cm]

\qquad \textbullet \quad
$
\ds
\thetafcn_4 (\restr{z}{\tau}) 
\ = \ 
\big(- \sqrt{-1} \hsx \tau \big)^{-\frac{1}{2}} \ 
\exp
\Big(
\frac{z^2}{\pi \hsy \sqrt{-1} \hsx \tau}
\Big)
\hsx 
\thetafcn_2 \Big(\restrBig{\frac{z}{\tau}}{\hsx -\tau^{-1}}\Big).
$
\\

[Note: \ 
The square root is real and positive when $\tau$ is purely imaginary.]
\\

{\small\bf EXAMPLE} 
Take $z = x$ real and $\tau = \sqrt{-1} \hsx t$ $(t > 0)$ $-$then

\[
\thetafcn_3 (\restr{x}{\sqrt{-1} \hsx t}) 
\ = \ 
\frac{1}{\sqrt{t}} \ 
\exp
\Big(
-\frac{x^2}{\pi \hsy t}
\bigg)
\hsx 
\thetafcn_3 
\Big(
\restrBig
{
\frac{x}{\sqrt{-1} \hsx t}
}
{
\frac{\sqrt{-1}}{t}
}
\Big).
\]
Specializing still further, let $x = 0$, and put

\[
\thetafcn (t) 
\ = \ 
\sum\limits_{n = 1}^\infty \ 
e^{- n^2 \hsy \pi \hsy t}, 
\]
thus
\allowdisplaybreaks\begin{align*}
1 + 2 \hsy \thetafcn (t) \ 
&=\ 
\thetafcn_3 \big(\restr{0}{\sqrt{-1} \hsx t}\big)
\\[15pt]
&=\ 
\frac{1}{\sqrt{t}} \ 
\thetafcn_3 \Big(\restrBig{0}{\frac{\sqrt{-1}}{t}}\Big)
\\[15pt]
&=\ 
\frac{1}{\sqrt{t}} \hsx 
\Big(
1 + 2 \hsx \thetafcn \Big(\frac{1}{2}\Big)
\Big).
\end{align*}



\chapter{
$\boldsymbol{\S}$\textbf{2}.\quad  ORDER}
\setlength\parindent{2em}
\setcounter{theoremn}{0}
\renewcommand{\thepage}{\S2-\arabic{page}}

\qquad 
Given an entire function 
\[
f(z) 
\ = \ 
\sum\limits_{n = 0}^\infty \ 
c_n z^n 
\qquad 
(\implies 
\lim\limits_{n \ra \infty} \ 
\abs{c_n}^{1/n} = 0), 
\]
put 
\[
M(r ; f) 
\ = \ 
\max\limits_{\abs{z} = r} \abs{f(z)}.
\]
\\[-1cm]

\begin{x}{\small\bf LEMMA} \ 
$M(r ; f)$ is a continuous increasing function of $r$.
\\[-.25cm]
\end{x}

\begin{x}{\small\bf LEMMA} \ 
If $f$ is not a constant, then 
\[
M(r ; f)  \ra \infty 
\qquad (r \ra \infty).
\]
\\[-1.25cm]
\end{x}

\begin{x}{\small\bf LEMMA} \ 
If for some $\lambda > 0$, 
\[
\underset{r \ra \infty}{\liminfx} \ 
\frac{M(r ; f)}{r^\lambda}
\ = \ 
0,
\]
then $f$ is a polynomial of degree $\leq \lambda$.
\\[-.5cm]

PROOF  \ 
In general, 
\[
\abs{c_n}
\ \leq \ 
\frac{M(r ; f)}{r^n}, 
\]
so for $n > \lambda$, 
\[
\abs{c_n}
\ \leq \ 
\underset{r \ra \infty}{\liminfx} \ 
\frac{M(r ; f)}{r^\lambda}
\ = \ 
0.
\]
\\[-.25cm]
\end{x}

\begin{x}{\small\bf EXAMPLE} \ 
We have
\[
\begin{cases}
\ 
M(r ; \exp z^n) 
\ = \ 
\exp r^n 
\qquad (n = 1, 2, \ldots)
\\[4pt]
\ 
M(r ; \exp e^z) 
\ = \ 
\exp e^r 
\end{cases}
.
\]
\\[-.25cm]
\end{x}


\begin{x}{\small\bf EXAMPLE} \ 
We have
\[
\begin{cases}
\ \ds
M(r ; \sin z) 
\ = \ 
\frac{e^r - e^{-r}}{2}
\\[11pt]
\ \ds
M(r ; \cos z) 
\ = \ 
\frac{e^r + e^{-r}}{2}
\end{cases}
.
\]
\\[-1.cm]
\end{x}

\begin{x}{\small\bf LEMMA} \ 
Let
\[
p(z)
\ = \ 
a_0 + a_1 z + \cdots + a_n z^n 
\qquad (a_n \neq 0, \ n \geq 1)
\]
be a polynomial of degree $n$ $-$then
\[
M(r; p(z)) \sim \abs{a_n} r^n 
\qquad (r \ra \infty).
\]
\\[-1.25cm]
\end{x}

\begin{x}{\small\bf DEFINITION} \ 
An entire function is said to be \un{transcendental} if it is not a polynomial.
\\[-.25cm]
\end{x}

\begin{x}{\small\bf LEMMA} \ 
If $f$ is transcendental, then for any polynomial $p$, 
\[
\lim\limits_{r \ra \infty} \ 
\frac{M(r; p)}{M(r; f)}
\ = \ 0.
\]
\\[-1.25cm]
\end{x}

\begin{x}{\small\bf DEFINITION} \ 
If $f \not\equiv C$ is an entire function, then its \un{order} $\rho$ $(= \rho (f))$ is given by
\[
\underset{r \ra \infty}{\limsupx} \ 
\frac{\log \log M(r ; f)}{\log r} \hsx .
\]

[Note: \ 
Conventionally, the order of $f \equiv C$ is 0.]
\\[-.25cm]
\end{x}

\begin{x}{\small\bf REMARK} \ 
The reason that one works with $\log \log M(r ; f)$ rather than $\log M(r ; f)$ 
is that if $f$ is transcendental, then 
\[
\lim\limits_{r \ra \infty} \ 
\frac{\log M(r; f)}{\log r}
\ = \ \infty.
\]
\\[-.25cm]
\end{x}
\begin{x}{\small\bf EXAMPLE} \ 
Every polynomial is an entire function of order 0 (cf. 2.6) but 
there are transcendental entire functions of order 0, e.g., 
$
\ds
\sum\limits_{n = 0}^\infty \ 
e^{-n^2} \hsx z^n
$
(cf. 2.27).
\\[-.25cm]
\end{x}

\begin{x}{\small\bf EXAMPLE} \ 
The entire function $\exp z^n$ $(n = 1, 2, \ldots)$ is of order $n$.  
On the other hand, the entire function $\exp e^z$ is of order $\infty$.
\\[-.25cm]
\end{x}

\begin{x}{\small\bf DEFINITION} \ 
$f$ is of \un{finite order} if $\rho$ is finite; otherwise, $f$ is of \un{infinite order}.
\\[-.25cm]
\end{x}

\begin{x}{\small\bf LEMMA} \ 
An entire function $f$ is of finite order iff there exists a positive constant $K$ such that 

\[
M(r; f) 
\ < \ 
\exp r^K 
\qquad (r \gg 0), 
\]
the greatest lower bound of the set of all such $K$ then being the order of $f$. 
\\[-.25cm]
\end{x}

\begin{x}{\small\bf LEMMA} \ 
An entire function $f$ is of finite order iff there exist positive constants $B$, $C$, and $K$ such that
\[
M(r; f) 
\ < \
B \hsx \exp C r^K
\qquad (r \gg 0),
\]
the greatest lower bound of the set of all such $K$ then being the order of $f$.
\\[-.5cm]

[Note: \ 
In general, the constants $B$ and $C$ depend on $K$.]
\\[-.25cm]
\end{x}

\begin{x}{\small\bf APPLICATION} \ 
Suppose that $f$ is an entire function of finite order.  
Given a complex constant $A$, 
let 
$f_A(z) = f(z + A)$ $-$then
$\rho(f) = \rho(f_A)$.  
\\[-.5cm]

[For $\exists \ K > 0$: 
\[
M(r; f) 
\ < \ 
\exp r^K 
\qquad (r \gg 0).
\]
But
\[
\abs{z}
\ < \ 
\abs{A}
\hsx \implies \hsx 
\abs{z + A} < 2 \hsy \abs{z}
\]
\qquad 
$\implies$

\[
M(r; f_A) 
\ < \ 
\exp 2^K \hsy r^K 
\qquad (r \gg 0).]
\]
\\[-1.25cm]
\end{x}

\begin{x}{\small\bf APPLICATION} \ 
Suppose that $f$ is an entire function of finite order.  
Given a nonzero complex constant $A$, 
let 
$f_A(z) = f (A \hsy z)$ $-$then
$\rho(f) = \rho(f_A)$. 
\\[-.5cm]

[For $\exists \ K > 0$: 
\[
M(r; f) 
\ < \ 
\exp r^K 
\qquad (r \gg 0) \hsx .
\]
But
\[
\abs{A z} 
\ \leq \ 
\abs{A} \abs{z}
\]
\qquad 
$\implies$

\[
M(r; f_A) 
\ < \ 
\exp \abs{A}^K r^K
\qquad (r \gg 0) \hsx .]
\]
\\[-1.25cm]
\end{x}

\begin{x}{\small\bf LEMMA} \ 
If $M(r; f) \sim h(r)$ $(r \ra \infty)$, then

\[
\underset{r \ra \infty}{\limsupx} \ 
\frac{\log \log M(r ; f)}{\log r}
\ = \ 
\underset{r \ra \infty}{\limsupx} \ 
\frac{\log \log h(r)}{\log r} \hsx .
\]

PROOF \ 
Assuming that $r \gg 0$, write
\allowdisplaybreaks
\allowdisplaybreaks\begin{align*}
\log M(r ; f) \ 
&=\ 
\log 
\bigg(
\frac{M(r ; f)}{h(r)}
h(r)
\bigg)
\\[15pt]
&=\ 
\log h(r) 
\hsx + \hsx 
\log 
\frac{M(r ; f)}{h(r)}
\\[15pt]
&=\ 
\log h(r) 
\bigg[
1 + 
\frac{1}{\log h(r)}
\hsy
\log 
\frac{M(r ; f)}{h(r)}
\bigg]
\end{align*}

$\implies$

\[
\frac{\log \log M(r ; f)}{\log r}
\ = \ 
\frac{\log \log h(r)}{\log r}
\hsx + \hsx
\frac{
\log 
\Big[
1 + 
\ds\frac{1}{\log h(r)}
\hsx 
\log \frac{M(r ; f)}{h(r)}
\Big]
}
{\log r}, 
\]
from which the assertion.
\\[-.25cm]
\end{x}

\begin{x}{\small\bf EXAMPLE} \ 
If $C$ is a positive constant, then 

\[
\underset{r \ra \infty}{\limsupx} \ 
\frac{\log \log C \hsy e^r}{\log r} 
\ = \ 1.
\]
This said, take now in 2.18
\[
h(r) 
\ = \ 
\frac{e^r}{2}
\]
to conclude that the entire functions $\sin z$ and $\cos z$ are both of order 1 (cf. 2.5).
\\[-.5cm]

[Note: \ 
Define entire functions
\[
\frac{\sin \sqrt{z}}{\sqrt{z}}, 
\quad 
\cos \sqrt{z}
\]
by the appropriate power series $-$then each is of order 
$
\ds
\frac{1}{2}
$.]
\\[-.25cm]
\end{x}

\begin{x}{\small\bf EXAMPLE} \ 
Put
\[
\Gamma_1 (z) 
\ = \ 
\int\limits_1^\infty \ 
t^z \hsy e^{-t} 
\ \td t.
\] 
Then $\Gamma_1$ is entire and  
\[
M(r; \Gamma_1) 
\ = \ 
\sqrt{2 \hsy \pi \hsy r} \ 
\Big(
\frac{r}{e}
\Big)^r
\hsx
\Big(
1 + \tO\Big(\frac{1}{r}\Big)
\Big) \hsx .
\]
Therefore
\[
\log M(r; \Gamma_1) 
\ \sim  \ 
r \log r 
\qquad (r \ra \infty),
\]
so $\rho(\Gamma_1) = 1$.
\\[-.25cm]
\end{x}

Sometimes it is simpler to work directly with 
$\log M(r; f)$.
\\[-.25cm]


\begin{x}{\small\bf EXAMPLE} \ 
Fix $\alpha > 0$ and let
\[
f_\alpha (z) 
\ = \ 
\prod\limits_{n = 1}^\infty \
\Big(
1 + \frac{z^n}{n^{\alpha \hsy n}}
\Big).
\]
Then 
\allowdisplaybreaks
\allowdisplaybreaks\begin{align*}
\log M(r; f_\alpha) 
&= \ 
\sum\limits_{n = 1}^\infty \
\log \Big(1 + \frac{r^n}{n^{\alpha \hsy n}}\Big)
\\[15pt]
&=\ 
\int\limits_0^\infty \ 
\log \Big(1 + \frac{r^u}{u^{\alpha \hsy u}}\Big) 
\ \td u 
\hsx + \hsx 
\tO(r^{\frac{1}{\alpha}})
\\[15pt]
&\sim\ 
r^\frac{2}{\alpha} \ 
\frac{1}{\alpha} \ 
\int\limits_1^\infty \ 
t^{-\frac{2}{\alpha} - 1} \hsy 
\log t 
\ \td t 
\qquad (r \ra \infty), 
\end{align*}
where we made the change of variable 
$
\ds 
t = \frac{r}{u^\alpha}.  
$
In the integral
\[
\int\limits_1^\infty \ 
t^{-\frac{2}{\alpha} - 1} \hsy 
\log t 
\ \td t, 
\]
let 
$
\ds
x = t^{\frac{2}{\alpha}}
$,
hence
\allowdisplaybreaks
\allowdisplaybreaks\begin{align*}
\frac{\alpha}{2} \ 
\int\limits_1^\infty \
\ds\frac{\log x^\frac{\alpha}{2}}{x^2} 
\ \td x \ 
&=\ 
\frac{\alpha^2}{4} \ 
\int\limits_1^\infty \
\frac{\log x}{x^2} 
\ \td x \ 
\\[15pt]
&=\ 
\frac{\alpha^2}{4} 
\hsx 
\Gamma(2)
\\[15pt]
&=\ 
\frac{\alpha^2}{4}.
\end{align*}
Therefore
\[
\log M(r; f_\alpha) 
\ \sim \ 
\frac{\alpha}{4}
\hsx 
r^\frac{2}{\alpha}
\qquad (r \ra \infty), 
\]
so
\[
\rho (f_\alpha) 
\ = \ 
\frac{2}{\alpha}.
\]
\\[-.25cm]
\end{x}


As will now be seen, the order $\rho$ of an entire function $f$ can be computed 
from the coefficients of its power series expansion at the origin. 
\\[-.25cm]

\begin{x}{\small\bf SUBLEMMA} \ 
If there exist positive constants $A$ and $K$ such that 

\[
M(r; f) 
\ < \ 
\exp A \hsy r^k 
\qquad (r \gg 0),
\]
then 
\[
\abs{c_n}
\ < \ 
\Big(
\frac{e \hsy A \hsy K}{n}
\Big)^{n/K}
\qquad (n \gg 0).
\]

PROOF \ 
For $r \gg 0$, say $r \geq r_0$, 
\[
\abs{c_n}
\ \leq \ 
\frac{M(r; f) }{r^n}
\ < \ 
\exp (A \hsy r^K - n \log r).
\]
As a function of $r$, 
\[
A \hsy r^K - n \log r
\]
achieves its minimum at $r_n$, 
where 
$r_n^K = n / (A \hsy K)$.  
But for $n \gg 0$, $r_n \geq r_0$.  
And 
\allowdisplaybreaks
\allowdisplaybreaks\begin{align*}
\exp \Big(A\hsy r_n^K - n \log r_n\Big)
&=\ 
\exp \Big(A \frac{n}{A \hsy K}\Big)
\hsx
\exp \Big(-n \log \Big(\frac{n}{A \hsy K}\Big)^{1/K}\Big)
\\[15pt]
&=\ 
\exp \Big(\frac{n}{ K}\Big)
\hsx
\exp \Big(\log \Big(\frac{n}{A \hsy K}\Big)^{-n/K}\Big)
\\[15pt]
&=\ 
\Big(\frac{e \hsy A \hsy K}{n}\Big)^{n/K}.
\end{align*}
\\[-1.25cm]
\end{x}

\begin{x}{\small\bf LEMMA} \ 
If there exist positive constants $A$ and $K$ such that 

\[
\abs{c_n} 
\ < \ 
\Big(
\frac{e \hsy A \hsy K}{n}
\Big)^{n / K} 
\qquad (r \gg 0),
\]
then $\forall \ \varepsilon > 0$, 
\[
M(r; f) 
\ < \ 
\exp (A + \varepsilon) 
\hsy r^K
\qquad (r \gg 0),
\]
hence
\[
M(r; f) 
\ < \ 
\exp r^{K + \varepsilon} 
\qquad (r \gg 0).
\]

PROOF \ 
We can and will assume that $c_0 = 0$ and

\[
\abs{c_n} 
\ < \ 
\Big(
\frac{e \hsy A \hsy K}{n}
\Big)^{n / K} 
\qquad \forall \ n \geq 1.
\]

Accordingly, 
\allowdisplaybreaks\begin{align*}
M(r; f)  \ 
&\leq \ 
\sum\limits_{n = 1}^\infty \ 
\abs{c_n}  \hsy r^n
\\[15pt]
&\leq \ 
\sum\limits_{n = 1}^\infty \ 
\Big(
\frac{e \hsy A \hsy K}{n}
\Big)^{n / K} 
\hsy r^n
\\[15pt]
&= \ 
\sum\limits_{n = 1}^\infty \ 
\Big(
\frac{e \hsy A \hsy r^K}{n / K}
\Big)^{n / K}.
\end{align*}

Put $m = [n / K]$:
\[
\begin{cases}
\ \ds
m!
\ \sim \ 
\Big(
\frac{m}{e}
\Big)^m
\hsx \sqrt{2 \hsy \pi \hsy m}
\\[15pt]
\ \ds
\sqrt{2 \hsy \pi \hsy m}
\ < \ 
C_1 \hsx 
\Big(
\frac{A + \varepsilon / 2}{A}
\Big)^{m + 1}
\end{cases}
.
\]

Therefore
\allowdisplaybreaks
\allowdisplaybreaks\begin{align*}
\Big(
\frac{e \hsy A \hsy r^K}{m}
\Big)^{m + 1}\ 
&=\ 
\Big(
\frac{e}{m}
\Big)
\hsy 
\Big(
\frac{e}{m}
\Big)^{m}
\hsy
\Big(
A \hsy r^K
\Big)^{m + 1}
\\[26pt]
&=\ 
\Big(\frac{e}{m}\Big)
\hsx
\frac
{\Big(\ds\frac{e}{m}\Big)^{m}}
{\ds\frac{\sqrt{2 \hsy \pi \hsy m}}{m!}}
\hsx
\frac{\ds\sqrt{2 \hsy \pi \hsy m}}{m!}
\hsx
\Big(
A \hsy r^K
\Big)^{m + 1}
\\[26pt]
&<\ 
C_2 \hsx
\frac{\sqrt{2 \hsy \pi \hsy m}}{m!}
\Big(
A \hsy r^K
\Big)^{m + 1}
\\[26pt]
&<\ 
C_3 \hsx \frac{1}{m!} \hsx
\Big(
\frac{A + \varepsilon / 2}{A}
\Big)^{m + 1}
\Big(
A \hsy r^K
\Big)^{m + 1}
\\[26pt]
&=\ 
C_3 \
\frac{(A + \varepsilon / 2)^{m + 1} \hsy r^{K (m + 1)}}{m!}
\end{align*}
$\implies$
\allowdisplaybreaks
\allowdisplaybreaks\begin{align*}
\sum\limits_{m = 1}^\infty \ 
\Big(
\frac{(A + \varepsilon / 2)^{m + 1} \hsy r^{K (m + 1)}}{m!}
\Big)
\ 
&=\ 
(A + \varepsilon / 2) \hsy 
(r^K) \hsy
(\exp(A + \varepsilon / 2) \hsy r^K - 1)
\\[15pt]
&<\ 
(A + \varepsilon / 2) \hsy 
(r^K) \hsy
\exp(A + \varepsilon / 2) \hsy r^K
\\[15pt]
&<\ 
\exp(A + \varepsilon) \hsy r^K
\qquad (r \gg 0).
\end{align*}
\\[-1.25cm]
\end{x}

\begin{x}{\small\bf THEOREM} \ 
The order of the entire function

\[
f(z) 
\ = \ 
\sum\limits_{n = 0}^\infty \ 
c_n \hsy z^n
\]
is given by 
\[
\rho 
\ = \ 
\underset{r \ra \infty}{\limsupx} \ 
\frac
{n \hsy \log n}
{\ds \log \hsy \frac{1}{\abs{c_n}}}
\]
or, equivalently, is given by
\[
\rho 
\ = \ 
\underset{r \ra \infty}{\limsupx} \ 
\frac
{\log n}
{\ds \log \hsy \frac{1}{\abs{c_n}^{1/n}}} \hsx.
\]

[Note: \ 
The terms for which $c_n = 0$ are taken to be 0.]
\\[-.5cm]

PROOF \ 
Suppose first that $\rho$ is finite $-$then for any $K > \rho$, 

\[
M(r; f) 
\ < \ 
\exp r^K
\qquad (r \gg 0), 
\]
thus by 2.22, 
\[
\abs{c_n} 
\ < \ 
\Big(
\frac{e \hsy K}{n}
\Big)^{n / K}
\qquad (n \gg 0).
\]
Therefore
\[
K 
\ > \ 
\frac
{\log n}
{\ds \log \hsy \frac{1}{\abs{c_n}^{1/n}}}
\hsx + \hsx 
\frac
{\raisebox{.2cm}{$\ds \log \frac{1}{e \hsy K}$}}
{\ds \log \hsy \frac{1}{\abs{c_n}^{1/n}}}
\qquad (n \gg 0).
\]
But 
\[
\lim\limits_{n \ra \infty} \ 
\log \frac{1}{\abs{c_n}^{1/n}}
\ = \ 
\infty,
\]
so 

\[
K 
\ \geq \ 
\underset{n \ra \infty}{\limsupx} \ 
\frac
{\log n}
{\ds \log \hsy \frac{1}{\abs{c_n}^{1/n}}}
\]

\qquad 
$\implies$

\[
\rho 
\ \geq \ 
\underset{n \ra \infty}{\limsupx} \ 
\frac
{\log n}
{\ds \log \hsy \frac{1}{\abs{c_n}^{1/n}}} \hsx .
\]
To reverse this, let
\[
K^\prime
\ > \ 
\underset{n \ra \infty}{\limsupx} \ 
\frac
{\log n}
{\ds \log \hsy \frac{1}{\abs{c_n}^{1/n}}}.
\]
Choose a positive integer $N(K^\prime)$:
\[
\frac
{\log n}
{\ds \log \hsy \frac{1}{\abs{c_n}^{1/n}}}
\ < \ 
K^\prime 
\qquad (n > N(K^\prime))
\]
or still, 
\[
\abs{c_n} 
\ < \ 
\Big(
\frac{1}{n}
\Big)^{n / K^\prime}
\qquad (n > N(K^\prime)).
\]
Then, thanks to 2.23 (with 
$
\ds
A = \frac{1}{e \hsy K^\prime}
$
), 
given $\varepsilon > 0$, there is an $R(\varepsilon)$: 
\[
M(r; f) 
\ < \ 
\exp 
\Big(
\frac{1}{e \hsy K^\prime} + \varepsilon
\Big)
\hsx r^{K^\prime}
\ < \ 
\exp r^{K^\prime + \varepsilon}
\qquad (r > R(\varepsilon)),
\]
hence
\[
\rho \leq K^\prime + \varepsilon 
\implies 
\rho \leq K^\prime
\implies
\rho \leq 
\underset{n \ra \infty}{\limsupx} \ 
\frac
{\log n}
{\ds \log \hsy \frac{1}{\abs{c_n}^{1/n}}} \hsx .
\]
In summary: 
For $\rho$ finite, 
\[
\rho  
\ = \
\underset{n \ra \infty}{\limsupx} \ 
\frac
{\log n}
{\ds \log \hsy \frac{1}{\abs{c_n}^{1/n}}}.
\]
Turning to the case of an infinite $\rho$, on the basis of what has been said above, it is clear that if 

\[
\underset{n \ra \infty}{\limsupx} \ 
\frac
{\log n}
{\ds \log \hsy \frac{1}{\abs{c_n}^{1/n}}} \hsx 
\]
is finite, then $\rho$ is finite, i.e., if $\rho$ is infinite, then

\[
\underset{n \ra \infty}{\limsupx} \ 
\frac
{\log n}
{\ds \log \hsy \frac{1}{\abs{c_n}^{1/n}}} 
\]
is infinite.
\\[-.25cm]
\end{x}

\begin{x}{\small\bf APPLICATION} \ 
The order of an entire function is unchanged by differentiation: 
$\rho (f) = \rho (f^\prime)$.
\\[-.25cm]
\end{x}

\begin{x}{\small\bf EXAMPLE} \ 
Let $0 < \rho < \infty$ $-$then the entire function

\[
f (z) 
\ = \ 
\sum\limits_{n = 1}^\infty \ 
\Big(
\frac{\rho \hsy e}{n}
\Big)^{n / \rho} \hsx z^n
\]
is of order $\rho$.
\\[-.25cm]
\end{x}

\begin{x}{\small\bf EXAMPLE} \ 
The entire function
\[
f (z) 
\ = \ 
\sum\limits_{n = 2}^\infty \ 
\Big(
\frac{1}{\log n}
\Big)^n \hsx z^n
\]
is of infinite order and the entire function
\[
f (z) 
\ = \ 
\sum\limits_{n = 0}^\infty \
e^{-n^2} 
\hsx z^n
\]
is of zero order.
\\[-.25cm]
\end{x}

\begin{x}{\small\bf EXAMPLE} \ 
Fix $\alpha > 0$ $-$then the entire function

\[
\ML_\alpha (z) 
\ = \ 
\sum\limits_{n = 0}^\infty \
\frac{z^n}{\Gamma (\alpha n + 1)}
\]
is of order 
$
\ds 
\frac{1}{\alpha}
$.
\\[-.25cm]

[Note: \ 
Obviously, 
\[
\begin{cases}
\ 
\ds
\ML_1 (z) 
\ = \ 
\sum\limits_{n = 0}^\infty \
\frac{z^n}{\Gamma (n + 1)}
\ = \ 
\sum\limits_{n = 0}^\infty \
\frac{z^n}{n !}
\ = \ 
e^z
\\[26pt]
\ds
\ML_2 (z)  
\ = \ 
\sum\limits_{n = 0}^\infty \
\frac{z^n}{\Gamma (2n + 1)}
\ = \ 
\sum\limits_{n = 0}^\infty \
\frac{z^n}{(2n) !}
\ = \ 
\cosh \sqrt{z}
\end{cases}
.]
\]
\\[-.5cm]
\end{x}

\begin{x}{\small\bf EXAMPLE} \ 
The \un{Bessel function} $\tJ_\nu (z)$ of the first kind of real index $\nu > -1$
is defined by the series

\[
\Big(
\frac{z}{2}
\Big)^\nu
\ 
\sum\limits_{n = 0}^\infty \ 
\frac{\ds(-1)^n \hsx  \Big(\frac{z}{2}\Big)^{2 \hsy n}}{n ! \hsy \Gamma (\nu + n + 1)},
\]
where 
$
\ds
\Big(
\frac{z}{2}
\Big)^\nu
\ = \ 
\exp \Big(\nu \log \frac{z}{2}\Big)
$, 
the logarithm having its principal value.  
Multiplying up, 

\[
\Big(
\frac{z}{2}
\Big)^{-\nu} \ \tJ_\nu (z)
\]
is therefore entire and, moreover, it is of order 1.
\\[-.25cm]
\end{x}

\begin{x}{\small\bf EXAMPLE} \ 
Fix $\alpha > 1$ $-$then the entire function

\[
\Phi_\alpha (z)
\ = \ 
\int\limits_0^\infty \ 
\exp (- t^\alpha) 
\hsy 
\cos z \hsy t 
\ \td t
\]
is of order 
$
\ds
\frac{\alpha}{\alpha - 1}
$.
\\[-.25cm]

[One first has to check that $\Phi_\alpha (z)$ really is entire, 
which can be seen by noting that it is uniformly convergent on compact subsets of $\Cx$:

\[
\abs{\cos z t}
\ \leq \ 
e^{t \hsy \abs{z}}
\]
\qquad
$\implies$
\[
\abs{\exp (- t^\alpha) \hsx \cos z t}
\ \leq \ 
\exp (t \hsy \abs{z} - t^\alpha) 
\ \leq \ 
\exp (- t)
\]
for all $t$ such that 
$
\ds
t^{\alpha - 1} > 1 + \abs{z}
$.
This settled, to compute the order, write
\allowdisplaybreaks\begin{align*}
\Phi_\alpha (z) \ 
&=\ 
\int\limits_0^\infty \ 
\exp (- t^\alpha) \ 
\bigg[
\sum\limits_{n = 0}^\infty \
\frac{(-1)^n \hsy z^{2 \hsy n} \hsy t^{2 \hsy n} }{(2 \hsy n)!}
\bigg]
\ \td t
\\[15pt]
&=\  
\sum\limits_{n = 0}^\infty \
\bigg[
\int\limits_0^\infty \ 
\exp (- t^\alpha) \hsx 
t^{2 \hsy n}
\ \td t
\bigg]
\hsx
\frac{(-1)^n \hsy z^{2 \hsy n} }{(2 \hsy n)!}
\\[15pt]
&=\ 
\frac{1}{\alpha} \ 
\sum\limits_{n = 0}^\infty \
\frac{(-1)^n}{(2 \hsy n)!}
\ 
\Gamma\Big(\frac{2 n + 1}{\alpha}\Big) 
\hsx 
z^{2 \hsy n},
\end{align*}
and then proceed \ldots \hsy .]
\\[-.5cm]

[Note: \ 
As a special case, 
\[
\Phi_2 (z) \ 
\ = \ 
\frac{1}{2} \hsx
\sqrt{\pi} \ 
\exp \Big(-\frac{z^2}{4}\Big), 
\]
an entire function of order 2 (by direct inspection).]
\\[-.25cm]
\end{x}

\begin{x}{\small\bf LEMMA} \ 
If $f_1$, $f_2$ are entire functions of respective orders $\rho_1$, $\rho_2$ 
and if $\rho_1 \leq \rho_2$ $(\rho_1 < \rho_2)$, then the order of 
$f_1 + f_2$ is $\leq \rho_2$ $(= \rho_2)$.
\\[-.25cm]
\end{x}

\begin{x}{\small\bf EXAMPLE} \ 
If 
$f_1 = e^z$, 
$f_2 = -e^z$ 
$-$then 
$\rho_1 = \rho_2 = 1$ but the order of 
$f_1 + f_2$ is 0.  
\\[-.25cm]
\end{x}

\begin{x}{\small\bf EXAMPLE} \ 
If $f$ is an entire function of order $\rho$, then for any polynomial $p$, 
the order of $f + p$ is equal to $\rho$.
\\[-.25cm]
\end{x}

\begin{x}{\small\bf LEMMA} \ 
If $f_1$, $f_2$ are entire functions of respective orders $\rho_1$, $\rho_2$ 
and if $\rho_1 \leq \rho_2$ $(\rho_1 < \rho_2)$, then the order of 
$f_1 \hsy f_2$ is $\leq \rho_2$ $(= \rho_2)$.
\\[-.25cm]
\end{x}

\begin{x}{\small\bf EXAMPLE} \ 
If 
$f_1 = e^z$, 
$f_2 = e^{-z}$ 
$-$then 
$\rho_1 = \rho_2 = 1$ but the order of 
$f_1 \hsy f_2$ is 0.  
\\[-.25cm]
\end{x}

\begin{x}{\small\bf EXAMPLE} \ 
If $f$ is an entire function of order $\rho$, 
then for any nonzero polynomial $p$, 
the order of $p \hsy f$ is equal to 
$\rho$.
\\[-.5cm]

[Note: \ 
If the quotient 
$\ds
\frac{f}{p}
$ 
is an entire function, 
then it too is of order 
$\rho$.
\\

Proof: \ 
$
\ds
\rho
\hsy
\Big(\frac{f}{p}\Big) 
\ = \ 
\rho 
\hsy
\Big(p \cdot \frac{f}{p}\Big) 
\ = \ 
\rho (f).]
$
\\[-.25cm]
\end{x}

\begin{x}{\small\bf LEMMA} \ 
If $f$, $g$ are entire functions and if 
$
\ds
\frac{f}{g}
$
is an entire function, then
\[
\rho
\hsy
\Big(\frac{f}{g}\Big) 
\ \leq \ 
\max(\rho (f),\rho (g)).
\]

PROOF \ 
Since 
$
\ds
g \cdot \frac{f}{g} 
= f
$, 
in the event that 
$
\ds
\rho
\hsy
\Big(\frac{f}{g}\Big) 
> 
\rho (g)
$, 
we have
\[
\rho
\hsy
\Big(\frac{f}{g}\Big)
\ = \ 
\rho 
\hsy
\Big(g \cdot \frac{f}{g}\Big) 
\ = \ 
\rho (f) 
\qquad (\tcf. \ 2.34),
\]
leaving the case 
$
\ds
\rho
\hsy
\Big(\frac{f}{g}\Big) 
\leq 
\rho (g).
$
\\[-.25cm]
\end{x}

\begin{x}{\small\bf EXAMPLE} \ 
Consider the theta functions

\[
\begin{cases}
\
\theta_1(\restr{z}{\tau})
\\[4pt]
\
\theta_2(\restr{z}{\tau})
\\[4pt]
\
\theta_3(\restr{z}{\tau})
\\[4pt]
\
\theta_4(\restr{z}{\tau})
\end{cases}
\]
of the Appendix to \S1 $-$then each is of order 2.  
First
\[
\begin{cases}
\ \ds
\theta_2(\restr{z}{\tau}) 
\ = \ 
\theta_1 \Big(\restrBig{z + \frac{\pi}{2}\hsy}{\tau}\Big)
\\[18pt]
\ \ds
\theta_3(\restr{z}{\tau}) 
\ = \ 
\theta_4 \Big(\restrBig{z + \frac{\pi}{2}\hsy}{\tau}\Big)
\end{cases}
.
\]
Therefore
\[
\begin{cases}
\
\rho (\theta_2) 
\ = \ 
\rho (\theta_1) 
\\[4pt]
\
\rho (\theta_3) 
\ = \ 
\rho (\theta_4) 
\end{cases}
,
\]
provided that 
$\theta_1$
and 
$\theta_4$ 
are of finite order (cf. 2.16).  
Next, recall the relation

\[
\theta_1
(\restr{z}{\tau})
\ = \ 
- \sqrt{-1} \hsx 
\exp
\Big(
\sqrt{-1} \hsx  z + \frac{1}{4} \pi \sqrt{-1} \hsx \tau
\Big)
\hsx 
\theta_4 
\Big(
\restrBig{z + \frac{\pi \hsy \tau}{2}}{\tau}
\Big) \hsx .
\]
Granting for the moment that 
$\rho (\theta_1) = 2$, 
the fact that 
$\exp (\sqrt{-1} \hsx z)$ 
is of order 1 in conjunction with 2.34 forces
\[
\rho 
\Big(
\theta_4 
\Big(
\restrBig{z + \frac{\pi \hsy \tau}{2}}{\tau}
\Big)
\Big)
\ = \ 
2
\]
from which 
$\rho (\theta_4) = 2$
(cf. 2.16).  
To deal with 
$\theta_1$, 
given $z$, let
\[
\lambda
\ = \ 
(2 \hsy \abs{z}  + \log 2)
/ 
\log \abs{1 / q}- \frac{1}{2}.
\]
Then
\allowdisplaybreaks
\allowdisplaybreaks\begin{align*}
\abs{\theta_1(\restr{z}{\tau})} \ 
&\leq \ 
2 \ 
\sum\limits_{n = 0}^\infty \ 
\abs{q}^{\big(n + \frac{1}{2}\big)^2}
\hsx
e^{(2 n + 1)\hsy \abs{z}}
\\[15pt]
&\leq \ 
2 \ 
\sum\limits_{n \leq \lambda} \ 
\abs{q}^{\big(n + \frac{1}{2}\big)^2}
e^{(2 n + 1)\hsy \abs{z}}
\hsy + \hsx 
2 \ 
\sum\limits_{n > \lambda} \ 
\Big(\frac{1}{2}\Big)^{n + \frac{1}{2}}
\\[12pt]
&= \  
\tO
\big(
e^{(2 \hsy \lambda + 1) \hsy \abs{z}}
\big)
\\[12pt]
&= \  
\tO
\big(
e^{C\hsy \abs{z}^2}
\big).
\\[-.5cm]
\end{align*}
Therefore 
$\rho (\theta_1) \leq 2$.  
That 
$\rho (\theta_1) = 2$
is established in 4.27.
\\[-.25cm]
\end{x}

\begin{x}{\small\bf EXAMPLE} \ 
Then entire function 
\[
1 + 
\sum\limits_{n = 1}^\infty \ 
\Big(\frac{1}{2}\Big)^{n^2}
\hsx
e^{n z}
\]
is of order 2.
\\[-.25cm]
\end{x}

\begin{x}{\small\bf NOTATION} \ 
Given an entire function  $f$, let
\[
A(r; f) 
\ = \ 
\max\limits_{\abs{z} = r} \ 
\Reg f(z).
\]
\\[-1.25cm]
\end{x}

\begin{x}{\small\bf RAPPEL} \ 
If for some $C > 0$, $d > 0$, 
\[
A(r; f) 
\ < \ 
C \hsy r^d
\qquad (r \gg 0), 
\]
then $f$ is a polynomial of degree $\leq [d]$.
\\[-.25cm]
\end{x}

\begin{x}{\small\bf LEMMA} \ 
If $f$ is an entire function and if the order of $F = e^f$ is finite, 
then $f$ is a polynomial (and the order of $F$ is equal to the degree of $f$).
\\[-.5cm]

PROOF \ 
From the definitions, 
\[
\log \abs{F (z)} 
\ = \ 
\Reg f(z),
\]
hence
\[
\log M(r ; f) 
\ = \ 
A(r; f).
\]
But $\forall \ \varepsilon > 0$, 
\[
\frac{\log \log M(r ; F)}{\log r} 
\ < \ 
\rho (F) + \varepsilon 
\qquad (r \gg 0), 
\]
thus
\[
\log M(r ; F)
\ < \ 
r^{\rho (F) + \varepsilon} 
\qquad (r \gg 0) 
\]
and so
\[
A (r; f)
\ < \ 
r^{\rho (F) + \varepsilon} 
\qquad (r \gg 0) \hsx .
\]
Therefore $f$ is a polynomial of degree $\ \leq [\rho (F) + \varepsilon]\ $ or still, 
$f$ is a polynomial of degree $\leq  [\rho (F)]$.
\\[-.25cm]
\end{x}


\chapter{
$\boldsymbol{\S}$\textbf{3}.\quad  TYPE}
\setlength\parindent{2em}
\setcounter{theoremn}{0}
\renewcommand{\thepage}{\S3-\arabic{page}}

\qquad 
Let $f$ be an entire function of order $\rho$, where $0 < \rho < \infty$.
\\[-.25cm]

\begin{x}{\small\bf DEFINITION} \ 
The \un{type} $\tau$ $(= \tau (f))$ of $f$ is given by 

\[
\underset{r \ra \infty}{\limsupx} \ 
\frac{\log M(r; f)}{r^\rho}.
\]
\\[-1.5cm]
\end{x}

\begin{x}{\small\bf EXAMPLE} \ 
The entire function
\[
\exp (a_0 + a_1 \hsy z + \cdots + a_n \hsy z^n) 
\qquad (a_n \neq 0, \ n \geq 1)
\]
is of order $n$ and type $\abs{a_n}$.
\\[-.25cm]
\end{x}

\begin{x}{\small\bf EXAMPLE} \ 
The entire functions
\[
\begin{cases}
\ 
\sin A\hsy z
\\[4pt]
\ 
\cos A\hsy z
\end{cases}
(A \neq 0)
\]
are of order 1 and type $\abs{A}$.
\\[-.25cm]
\end{x}

\begin{x}{\small\bf DEFINITION} \ 
$f$ is of 
\un{maximal type} 
if 
$\tau = \infty$, 
of 
\un{minimal type} 
if 
$\tau = 0$, 
and of 
\un{intermediate type} 
if 
$0 < \tau < \infty$.
\\[-.25cm]
\end{x}

\begin{x}{\small\bf REMARK} \ 
$f$ is of \un{finite type} if $0 \leq \tau < \infty$, 
which will be the case iff there exists a positive constant $C$ such that

\[
M(r; f) 
\ < \ 
\exp C \hsy r^\rho 
\qquad (r \gg 0),
\]
the greatest lower bound of the set of all such $C$ then being the type of $f$. 
\\[-.25cm]
\end{x}

Here is a formula for the type parallel to that of 2.24 for the order. 
\\[-.25cm]

\begin{x}{\small\bf THEOREM} \ 
The type of the entire function 

\[
f(z) 
\ = \ 
\sum\limits_{n = 0}^\infty \ 
c_n \hsy z^n
\]
is given by
\[
\tau 
\ = \ 
\frac{1}{\rho \hsy e} \ 
\underset{n \ra \infty}{\limsupx} \ 
\big(
n \hsy \abs{c_n}^{\rho / n}
\big).
\]

PROOF \ 
Suppose first that $\tau$ is finite $-$then for any $A > \tau$, 

\[
M(r; f) 
\ < \ 
\exp A \hsy r^\rho 
\qquad (r \gg 0),
\]
thus by 2.22, 
\[
\abs{c_n} 
\ < \ 
\Big(\frac{\rho \hsy e \hsy A}{n}
\Big)^{n / \rho}
\qquad (n\gg 0),
\]
so
\[
A 
\ > \ 
\frac{1}{\rho \hsy e} \ 
n \hsy \abs{c_n}^{\rho / n}
\qquad (n\gg 0).
\]
Therefore
\[
A 
\ \geq \ 
\frac{1}{\rho \hsy e} \ 
\underset{n \ra \infty}{\limsupx} \ 
\big(
n \hsy \abs{c_n}^{\rho / n}
\big)
\]

\qquad 
$\implies$

\[
\tau
\ \geq \ 
\frac{1}{\rho \hsy e} \ 
\underset{n \ra \infty}{\limsupx} \ 
\big(
n \hsy \abs{c_n}^{\rho / n}
\big).
\]
To go the other way, let
\[
K^\prime
\ > \ 
\frac{1}{\rho \hsy e} \ 
\underset{n \ra \infty}{\limsupx} \ 
\big(
n \hsy \abs{c_n}^{\rho / n}
\big).
\]
Choose a positive integer $N(K^\prime)$: \ 
\[
\frac{1}{\rho \hsy e} \ 
n \hsy \abs{c_n}^{\rho / n}
\ < \ 
K^\prime
\qquad (n > N(K^\prime))
\]
or still, 

\[
\abs{c_n}  
\ < \ 
\Big(
\frac{\rho \hsy e \hsy K^\prime}{n}
\Big)^{n / \rho}
\qquad (n > N(K^\prime)).
\]
Then, thanks to 2.23 (with $A = K^\prime$, $K = \rho$), given any $\varepsilon > 0$, there is an $R(\varepsilon)$:

\[
M(r; f) 
\ < \ 
\exp(K^\prime + \varepsilon) r^\rho 
\qquad (r > R(\varepsilon)),
\]
hence

\[
\tau 
\ \leq \ 
K^\prime + \varepsilon 
\hsx \implies \hsx
\tau 
\ \leq \ 
K^\prime
\hsx \implies \hsx
\tau 
\ \leq \
\frac{1}{\rho \hsy e} \ 
\underset{n \ra \infty}{\limsupx} \ 
\big(
n \hsx \abs{c_n}^{\rho / n}
\big).
\]
In summary: \ 
For $\tau$ finite, 

\[
\tau 
\ = \ 
\frac{1}{\rho \hsy e} \ 
\underset{n \ra \infty}{\limsupx} \ 
\big(
n \hsx \abs{c_n}^{\rho / n}
\big) .
\]
Turning to the case of an infinite $\tau$, 
on the basis of what has been said above, 
it is clear that if 

\[
\frac{1}{\rho \hsy e} \ 
\underset{n \ra \infty}{\limsupx} \ 
\big(
n \hsx \abs{c_n}^{\rho / n}
\big)
\]
is finite, then $\tau$ is finite, i.e., if $\tau$ is infinite, then

\[
\frac{1}{\rho \hsy e} \ 
\underset{n \ra \infty}{\limsupx} \ 
\big(
n \hsx \abs{c_n}^{\rho / n}
\big)
\]
is infinite.
\\[-.25cm]
\end{x}
\begin{x}{\small\bf APPLICATION} \ 
They type of an entire function is unchanged by differentiation: \ 
$\tau (f) = \tau(f^\prime)$.
\\[-.25cm]
\end{x}

\begin{x}{\small\bf EXAMPLE} \ 
Let $0 < \rho < \infty$ $-$then the entire function

\[
f (z) 
\ = \ 
\sum\limits_{n = 2}^\infty \ 
\Big(
 \frac{\rho \hsy e}{n \log n}
\Big)^{n / \rho} \hsy z^n
\]
is of order $\rho$ and of minimal type. 
\\[-.25cm]
\end{x}

\begin{x}{\small\bf EXAMPLE} \ 
Let $0 < \rho < \infty$ $-$then the entire function
\[
f (z) 
\ = \ 
\sum\limits_{n = 2}^\infty \ 
\Big(
\rho \hsy e \hsx \frac{\log n}{n}
\Big)^{n / \rho} \hsy z^n
\]
is of order $\rho$ and of maximal type.
\\[-.25cm]
\end{x}

\begin{x}{\small\bf EXAMPLE} \ 
The entire function

\[
z 
\ra 
\int\limits_0^1 \ 
e^{z \hsy t^2}
\ \td t
\]
is of order 1 and of type 1.
\\[-.25cm]
\end{x}

\begin{x}{\small\bf EXAMPLE} \ 
Let 
$0 < \rho < \infty$, 
$0 < \tau < \infty$
$-$then the entire function

\[
f (z) 
\ = \ 
\sum\limits_{n = 1}^\infty \ 
\Big(
\frac{\rho \hsy e \hsx \tau}{n}
\Big)^{n / \rho} \hsy z^n
\]
is of order $\rho$ and type $\tau$ (cf. 2.26).
\\[-.25cm]
\end{x}

\begin{x}{\small\bf EXAMPLE} \ 
Fix 
$\alpha > 0$, 
$A > 0$ 
$-$then the entire function

\[
\ML_{\alpha, A} (z) 
\ = \ 
\sum\limits_{n = 0}^\infty \ 
\frac{(A \hsy z)^n}{\Gamma (\alpha \hsy n + 1)}
\]
is of order 
$\ds
\frac{1}{\alpha}
$
and type $A$ (cf. 2.28).
\\[-.25cm]
\end{x}

\begin{x}{\small\bf EXAMPLE} \ 
Fix 
$t > 0$ 
and let 

\[
\theta_t (z) 
\ = \ 
1 + 
\sum\limits_{n = 1}^\infty \ 
\big(
e^{-\pi \hsy t}
\big)^{n^2} \hsx e^{n \hsy z}.
\]
Then $\theta_t$ is of order 2 and of type
$\ds
\frac{1}{4 \hsy \pi \hsy t}
$.
\\[-.25cm]

[Note: \ 
As a special case, 

\[
\theta_{\frac{\log 2}{\pi}} (z) 
\ = \ 
1 + 
\sum\limits_{n = 1}^\infty \ 
\Big(
\frac{1}{2}
\Big)^{n^2}
e^{n \hsy z},
\]
an entire function of order 2 and of type 
$\ds
\frac{1}{4 \hsy \log 2}
$
(cf. 2.39).]
\\[-.25cm]
\end{x}

\begin{x}{\small\bf LEMMA} \ 
Let $f_1$, $f_2$ be entire functions of respective orders $\rho_1$, $\rho_2$, where 
$0 < \rho_1 < \infty$, 
$0 < \rho_2 < \infty$, 
and respective types
$\tau_1$, $\tau_2$.
\\[-.25cm]

\qquad \textbullet \quad
If 
$\rho_1 < \rho_2$, 
then 
$\rho (f_1 f_2) = \rho (f_2)$
and 
$\tau (f_1 \hsy f_2) = \tau_2$.
\\[-.25cm]

\qquad \textbullet \quad
If 
$\rho_1 = \rho_2$, 
if
$0 < \tau_1 < \infty$, 
if
$\tau_2 = 0$,
then 
$\rho (f_1 f_2) = \rho_1 = \rho_2$ 
and 
$\tau (f_1 \hsy f_2) = \tau_1$.
\\[-.25cm]

\qquad \textbullet \quad
If
$\rho_1 = \rho_2$, 
if
$\tau_1 = \infty$,
if 
$0 \leq \tau_2 < \infty$,
then 
$\rho (f_1 f_2) = \rho_1 = \rho_2$ 
and 
$\tau (f_1 \hsy f_2) = \infty$.
\\[-.25cm]
\end{x}


\chapter{
$\boldsymbol{\S}$\textbf{4}.\quad  CONVERGENCE EXPONENT}
\setlength\parindent{2em}
\setcounter{theoremn}{0}
\renewcommand{\thepage}{\S4-\arabic{page}}

\qquad
Let 
$\{r_n : n = 1, 2, \ldots\}$ 
be a sequence of positive real numbers with 
\[
0 < r_1 \leq r_2 \leq \cdots 
\qquad (r_n \ra \infty),
\]
finite repetitions being permitted.
\\[-.25cm]

\begin{x}{\small\bf DEFINITION} \ 
The greatest lower bound $\kappa$ of the positive $p$ for which the series
\[
\sum\limits_{n = 1}^\infty \ 
\frac{1}{r_n^p}
\]
is convergent is called the \un{convergence exponent} of the sequence 
$\{r_n : n = 1, 2, \ldots\}$.
\\[-.25cm]
\end{x}

\qquad
{\small\bf \un{N.B.}} \ 
If $\forall \ p$, 
\[
\sum\limits_{n = 1}^\infty \ 
\frac{1}{r_n^p}
\ = \ 
\infty,
\]
then take $\kappa = \infty$.
\\[-.5cm]

\begin{x}{\small\bf EXAMPLE} \ 
The sequence $\{e^n\}$ has convergence exponent 0.
\\[-.25cm]
\end{x}

\begin{x}{\small\bf EXAMPLE} \ 
The sequence $\{\log n\}$ has convergence exponent $\infty$.
\\[-.25cm]
\end{x}

\begin{x}{\small\bf REMARK} \ 
Take $\kappa < \infty$ $-$then the series
\[
\sum\limits_{n = 1}^\infty \ 
\frac{1}{r_n^\kappa}
\]
may or may not converge.
\\[-.5cm]

[The sequence $\{n\}$ has convergence exponent 1 and 
$
\ 
\ds
\sum\limits_{n = 1}^\infty \ 
\frac{1}{n}
\ 
$
is divergent while the sequence 
$\{n \hsy (\log n)^2\}$ 
also has convergence exponent 1 but 
$
\ 
\ds
\sum\limits_{n = 2}^\infty \ 
\frac{1}{n \hsy (\log n)^2}
\ 
$
is convergent.]
\\[-.25cm]
\end{x}


\begin{x}{\small\bf LEMMA} \ 
We have
\[
\kappa 
\ = \ 
\underset{n \ra \infty}{\limsupx} \ 
\frac{\log n}{\log r_n}.
\]
\\[-1.5cm]
\end{x}

\begin{x}{\small\bf DEFINITION} \ 
The 
\un{counting function} 
$n(r)$ $(r \geq 0)$ of the sequence
$\{r_n : n = 1, 2, \ldots\}$ 
is the number of $r_n$ such that $r_n \leq r$, i.e., 
\[
n(r) 
\ = \ 
\sum\limits_{r_n \leq r} \ 
1.
\]

[Note: \ 
$n(r) = 0$ for $0 \leq r < r_1$.  
In addition, $n(r)$ is right continuous, increasing, integer valued, and piecewise constant.]
\\[-.25cm]
\end{x}

\begin{x}{\small\bf EXAMPLE} \ 
Take $r_n = n$ $\forall \ n$ $-$then $n (r) = [r]$.
\\[-.25cm]
\end{x}

\begin{x}{\small\bf EXAMPLE} \ 
Let 
$\{r_n : n = 1, 2, \ldots\}$ 
be the sequence derived from the lattice points in the plane (excluding $(0,0)$) 
$-$then 
\[
\sum\limits_{n = 1}^\infty \ 
\frac{1}{r_n^p} 
\ = \ 
\sum\limits_{(m,n) \neq (0,0)} \ 
\frac{1}{(m^2 + n^2)^{p/2}},
\]
the series on the right being convergent if $p > 2$ and divergent if $p \leq 2$, hence $\kappa = 2$.  
And here
\[
n(r) 
\sim 
\pi \hsy r^2 
\qquad (r \ra \infty).
\]
\\[-1.5cm]
\end{x}

\begin{x}{\small\bf LEMMA} \ 
We have
\[
\underset{r \ra \infty}{\limsupx} \ 
\frac{\log n(r)}{\log r}
\ = \ 
\underset{n \ra \infty}{\limsupx} \ 
\frac{\log n}{\log r_n}.
\]
\\[-1.5cm]
\end{x}

\begin{x}{\small\bf APPLICATION} \ 
The convergence exponent $\kappa$ is given by 

\[
\underset{r \ra \infty}{\limsupx} \ 
\frac{\log n(r)}{\log r}
\qquad (\tcf. \ 4.5).
\]
\\[-1.5cm]
\end{x}


\begin{x}{\small\bf DEFINITION} \ 
Take $\kappa < \infty$ $-$then the 
\un{density} 
of the sequence 
$\{r_n : n = 1, 2, \ldots\}$ 
is
\[
\Delta \ 
\ =\ 
\underset{n \ra \infty}{\limsupx} \ 
\frac{n}{r_n^\kappa}.
\]
\\[-1.25cm]
\end{x}

\begin{x}{\small\bf EXAMPLE} \ 
Fix $p > 1$ and let $r_n = n^p$ $-$then 
$\kappa = 1/p$ and $\Delta = 1$.
\\[-.25cm]
\end{x}

\begin{x}{\small\bf LEMMA} \ 
We have
\[
\Delta \ 
\ =\ 
\underset{r \ra \infty}{\limsupx} \ 
\frac{n (r)}{r^\kappa}.
\]
\\[-1.25cm]
\end{x}

\begin{x}{\small\bf DEFINITION} \ 
Take $\kappa < \infty$ $-$then the 
\un{genus} 
of the sequence 
$\{r_n : n = 1, 2, \ldots\}$ 
is the smallest nonnegative integer $\fg$ such that 
\[
\sum\limits_{n = 1}^\infty \ 
\frac{1}{r_n^{\fg + 1}}
\]
is convergent.
\\[-.25cm]
\end{x}

\begin{x}{\small\bf LEMMA} \ 
Assume that $\kappa$ is finite.
\\[-.25cm]

\qquad \textbullet \quad
If $\kappa$ is not an integer, then $\fg = [\kappa]$.
\\[-.25cm]

\qquad \textbullet \quad 
If $\kappa$ is an integer, then $\fg = \kappa - 1$ if 
$
\ 
\ds
\sum\limits_{n = 1}^\infty \ 
\frac{1}{r_n^\kappa}
\ 
$
is convergent while 
$\fg = \kappa$ if 
$
\ 
\ds
\sum\limits_{n = 1}^\infty \ 
\frac{1}{r_n^\kappa}
\ 
$
is divergent.
\\[-.25cm]
\end{x}

Having dispensed with the formalities, 
we shall now come back to complex variable theory.  
So suppose that $f$ is a transcendental entire function of finite order $\rho$.  
Arrange the nonzero zeros of $f$ in a sequence $z_1, z_2, \ldots$ such that 
\[
0 
\ < \ 
\abs{z_1}
\ \leq \ 
\abs{z_2}
\ \leq \ 
\cdots
\]
with multiple zeros counted according to their multiplicities and let 
$r_n = \abs{z_n}$.
\\

\begin{x}{\small\bf THEOREM} \ 
Given $\varepsilon > 0$, 

\[
\underset{r \ra \infty}{\limsupx} \ 
\frac{n (r)}{r^{\rho + \varepsilon}}
\ \leq \ 
e (\rho + \varepsilon).
\]
\\[-1cm]

Before detailing the proof, it will be best to make some initial reductions.
\\[-.25cm]

\qquad \textbullet \quad 
If the number of zeros of $f$ is finite, then $n (r)$ is eventually constant and the result is trivial.  
It will therefore be assumed that 
$r_n = \abs{z_n} \ra \infty$.
\\[-.25cm]

\qquad \textbullet \quad
If $f (0) = 0$, write
$f (z) = z^m \hsy g(z)$ $(g (0) \neq 0)$ $-$then 
the order of $f$  equals the order of $g$ (cf. 2.36) 
so we can just as well assume from the beginning that $f (0) \neq 0$. 
\\[-.25cm]

\qquad \textbullet \quad 
Since multiplication by a nonzero constant does not affect the order of the zeros, 
there is no loss of generality in assuming that $\abs{f (0)} = 1$.
\\[-.25cm]
\end{x}

\begin{x}{\small\bf JENSEN INEQUALITY} \ 
If $\abs{f (0)} = 1$, then $\forall \ r > 0$, 
\[
\int\limits_0^r \ 
\frac{n (t)}{t}
\ \td t 
\ \leq \ 
\log M(r; f).
\]

Proceeding to the proof of 4.16, fix a parameter $\lambda \in ]0, 1[$ $-$then 
\allowdisplaybreaks\begin{align*}
\int\limits_0^r \ 
\frac{n (t)}{t}
\ \td t \ 
&\geq \ 
\int\limits_{\lambda \hsy r}^r \ 
\frac{n (t)}{t}
\ \td t
\\[15pt]
&\geq \ 
n (\lambda \hsy r) \ 
\int\limits_{\lambda \hsy r}^r \ 
\frac{\td t}{t}
\\[15pt]
&=\ 
n (\lambda \hsy r)
\hsx 
\log \frac{1}{\lambda}
\end{align*}
or still,
\[
n (\lambda \hsy r)
\ \leq \ 
\frac{1}{\ds \log \frac{1}{\lambda}}
\hsx 
\log M (r; f)
\]
or still,
\[
\frac{n (\lambda \hsy r)}{\log M (r; f)}
\ \leq \ 
\frac{1}{\ds \log \frac{1}{\lambda}}\hsx .
\]
Therefore
\[
\underset{r \ra \infty}{\limsupx} \ 
\frac{n (\lambda \hsy r)}{\log M (r; f)}
\ \leq \ 
\frac{1}{\ds \log \frac{1}{\lambda}} \hsx .
\]
But
\[
\log M(r; f) 
\ < \ 
r^{\rho + \varepsilon} 
\qquad (r \gg 0),
\]
thus
\[
\underset{r \ra \infty}{\limsupx} \ 
\frac{n (\lambda \hsy r)}{r^{\rho + \varepsilon}}
\ \leq \ 
\frac{1}{\ds \log \frac{1}{\lambda}}
\]
or still, 
\[
\underset{r \ra \infty}{\limsupx} \ 
\frac{n (r)}{r^{\rho + \varepsilon}}
\ \leq \ 
\frac{1}{\lambda^{\rho + \varepsilon}}
\hsx 
\frac{1}{\ds \log \frac{1}{\lambda}} \hsx .
\]
To finish up, simply take
\[
\lambda 
\ = \ 
e^{-1/(\rho + \varepsilon)}.
\]
\\[-1.25cm]
\end{x}

\begin{x}{\small\bf APPLICATION} \ 
If $f$ is a transcendental entire function of finite order $\rho$, then 
$\forall \ \varepsilon > 0$, 

\[
n (r) 
\ = \ 
\tO \big(r^{\rho + \varepsilon}\big).
\]
\\[-1.25cm]
\end{x}

\begin{x}{\small\bf LEMMA} \ 
If $\abs{f(0)} = 1$, then 
\[
n (r) 
\ \leq \ 
\log M(e \hsy r; f).
\]

PROOF \ 
In fact, 
\allowdisplaybreaks\begin{align*}
n(r) \ 
&=\ 
n (r) \ 
\int\limits_r^{e \hsy r} \ 
\frac{\td t}{t}
\\[15pt]
&\leq \ 
\int\limits_r^{e \hsy r} \ 
\frac{n (t)}{t}
\ \td t
\\[15pt]
&\leq \ 
\int\limits_0^{e \hsy r} \ 
\frac{n (t)}{t}
\ \td t
\\[15pt]
&\leq \ 
\log M(e \hsy r; f).
\end{align*}
\\[-1cm]
\end{x}

\begin{x}{\small\bf THEOREM} \ 
If $f$ is a transcendental entire function of finite order $\rho$, then 
the convergence exponent $\kappa$ of the sequence 
$\{r_n = \abs{z_n}\}$ 
is $\leq \rho$.
\\[-.5cm]

PROOF \ 
This, of course, is trivial if $f$ has a finite number of zeros 
(for then $\kappa = 0$), 
so as above it will be assumed that $f$ has an infinite number of zeros 
(hence that $r_n = \abs{z_n} \ra \infty$), 
matters reducing to the case when $\abs{f(0)} = 1$:
\allowdisplaybreaks\begin{align*}
\kappa \ 
&=\ 
\underset{r \ra \infty}{\limsupx} \ 
\frac{\log n(r)}{\log r}
\qquad (\tcf. \ 4.10)
\\[15pt]
&\leq \ 
\underset{r \ra \infty}{\limsupx} \ 
\frac{\log \log M(e \hsy r; f)}{\log r}
\qquad (\tcf. \ 4.19)
\\[15pt]
&\leq \ 
\underset{r \ra \infty}{\limsupx} \ 
\frac{\log \log M(e \hsy r; f)}{\log e \hsy r}
\cdot
\frac{\log e \hsy r}{\log r}
\\[15pt]
&= \ 
\underset{r \ra \infty}{\limsupx} \ 
\frac{\log \log M(r; f)}{\log r}
\\[15pt]
&= \ 
\rho.
\end{align*}
\\[-1.25cm]
\end{x}

\begin{x}{\small\bf COROLLARY} \ 
If $p > \rho$, then
\[
\sum\limits_{n = 1}^\infty \ 
\frac{1}{\abs{z_n}^p}
\ < \ 
\infty.
\]
\\[-1.25cm]
\end{x}

\begin{x}{\small\bf EXAMPLE} \ 
It can happen that $\kappa < \rho$. 
E.g.: \ 
If $f (z) = e^z$, then $\rho = 1$ but
there are no zeros, thus $\kappa = 0$.  
Another ``for instance'' is given by 
$\ds 
e^{z^2} \hsy \sin z
$, 
where $\kappa = 1 < 2 = \rho$. 
\\[-.5cm]

[Note: \ 
The so-called canonical products constitute a class of entire functions of finite order for which $\kappa = \rho$ 
(cf. 5.10).]
\\[-.25cm]
\end{x}

\begin{x}{\small\bf REMARK} \ 
If $\kappa$ is positive, then $f$ has an infinite number of zeros. 
\\[-.25cm]
\end{x}

\begin{x}{\small\bf DEFINITION} \ 
Let $f$ be a transcendental entire function of finite order $\rho$ $-$then 
$f$ is said to be of 
\un{convergence class} 
or 
\un{divergence class} 
according to whether 
\[
\sum\limits_{n = 1}^\infty \ 
\frac{1}{\abs{z_n}^\kappa}
\]
is convergent or divergent.
\\[-.25cm]
\end{x}

\begin{x}{\small\bf EXAMPLE} \ 
The transcendental entire function

\[
f(z) 
\ = \ 
\prod\limits_{n = 2}^\infty \ 
\Big(
1 - \frac{z}{n \hsy (\log n)^2}
\Big)
\]
is of order 1.  
Here $\kappa = 1$ and $f(z)$ is of convergence class (cf. 4.4).
\\[-.25cm]
\end{x}

\begin{x}{\small\bf EXAMPLE} \ 
The transcendental entire functions

\[
\begin{cases}
\ \sin z
\\[4pt]
\ \cos z
\end{cases}
\]
are of order 1 and of divergence class.
\\[-.25cm]
\end{x}

\begin{x}{\small\bf EXAMPLE} \ 
Consider the theta functions

\[
\begin{cases}
\ \theta_1 (\restr{z}{\tau})
\\[4pt]
\ \theta_2 (\restr{z}{\tau})
\\[4pt]
\ \theta_3 (\restr{z}{\tau})
\\[4pt]
\ \theta_4 (\restr{z}{\tau})
\end{cases}
\]
of the Appendix to \S1 $-$then the zeros of each of them are enumerated there and in all four cases, 
\[
\sum\limits_{n = 1}^\infty \ 
\frac{1}{\abs{z_n}^p}
\]
is convergent if $p > 2$ and divergent if $p \leq 2$ (cf. 4.8), 
hence $\kappa = 2$.  
On the other hand, 
it was shown in 2.38 that 
$\rho(\theta_1) \leq 2$, 
so 
$\rho(\theta_1) = 2$ 
$(\implies \rho(\theta_2) = \rho(\theta_3) = \rho(\theta_4) = 2)$. 
Therefore the theta functions are of divergence class.
\\[-.25cm]
\end{x}

\begin{x}{\small\bf LEMMA} \ 
If $\abs{f(0)} = 1$ and if $0 < \rho = \kappa < \infty$, then
\[
\Delta 
\ \leq \ 
e^\rho \tau.
\]

PROOF \ 
In fact, 
\allowdisplaybreaks\begin{align*}
\Delta \ 
&=\ 
\underset{r \ra \infty}{\limsupx} \ 
\frac{n (r)}{r^\kappa}
\qquad (\tcf. \ 4.13)
\\[15pt]
&\leq \ 
\underset{r \ra \infty}{\limsupx} \ 
e^\kappa \ 
\frac{\log M(e\hsy r; f)}{(e \hsy r)^\kappa}
\qquad (\tcf. \ 4.19)
\\[15pt]
&= \ 
\underset{r \ra \infty}{\limsupx} \ 
e^\rho \ 
\frac{\log M(e\hsy r; f)}{(e \hsy r)^\rho}
\\[15pt]
&= \ 
\underset{r \ra \infty}{\limsupx} \ 
e^\rho \ 
\frac{\log M(r; f)}{r^\rho}
\\[15pt]
&= \ 
e^\rho \tau
\qquad (\tcf. \ 3.1).
\end{align*}
\\[-1cm]
\end{x}

Maintaining the assumption that $f$ is a transcendental entire function of finite order $\rho$, 
suppose further that $f$ is of finite type $\tau$ (cf. 3.5), so $\rho > 0$.
\\[-.25cm]

\begin{x}{\small\bf THEOREM} \ 
We have
\[
\underset{r \ra \infty}{\limsupx} \ 
\frac{n(r)}{r^\rho}
\ \leq \ 
\rho \hsy e \hsy \tau.
\]
\\[-1.25cm]
\end{x}

The technical key to proving this is to employ a generalization of 4.17.
\\[-.25cm]

\begin{x}{\small\bf JENSEN INEQUALITY} \ 
If $f$ has a zero of order $m$ at the origin, then
\[
\int\limits_0^r \ 
\frac{n(t)}{t}
\ \td t
\ \leq \ 
\log M(r; f) 
- 
\log \abs{\frac{f^{(m)} (0)}{m !}} \hsy r^m.
\]

[Note: \ 
When $m = 0$, the correction term becomes
\[- \log \abs{f (0)}
\
\]
which disappears if in addition $\abs{f(0)} = 1$.]
\\[-.25cm]
\end{x}

To establish 4.29, 
start by fixing a parameter $\lambda \in \ ]0, 1[$ and then proceed as in the proof of 4.16:
\[
\int\limits_0^r \ 
\frac{n(t)}{t}
\ \td t
\ \geq \ 
n (\lambda \hsy r) \hsx \log \frac{1}{\lambda}
\]
or still, 
\[
n (\lambda \hsy r) 
\ \leq \ 
\frac{1}{\ds\log \frac{1}{\lambda}}
\hsx 
\Big(
\log M (r; f) 
- 
\log \abs{\frac{f^{(m)} (0)}{m !}} \hsy r^m
\Big)
\]
or still, 
\[
\frac
{n (\lambda \hsy r)}
{\log M (r; f)}
\ \leq \ 
\frac{1}{\ds\log \frac{1}{\lambda}}
\ 
\bigg(
1 - 
\frac
{\log \abs{\frac{f^{(m)} (0)}{m !}} \hsy r^m}
{\log M (r; f)}
\bigg)
\hsx .
\]
But
\[
\lim\limits_{r \ra \infty} \ 
\frac
{\log r}
{\log M (r; f)}
\ = \ 
0
\qquad (\tcf. \ 2.10).
\]
Therefore
\[
\underset{r \ra \infty}{\limsupx} \ 
\frac
{n (\lambda \hsy r)}
{\log M (r; f)}
\ \leq \ 
\frac{1}{\ds \log \frac{1}{\lambda}} \hsx .
\]
Since $f$ is of finite type, $\forall \ \varepsilon > 0$, 
\[
\log M(r; f) 
\ < \ 
(\tau + \varepsilon) \hsy r^\rho 
\qquad (r \gg 0).
\]
And this implies that
\[
\underset{r \ra \infty}{\limsupx} \ 
\frac{n (\lambda \hsy r)}{(\tau + \varepsilon) r^\rho}
\ \leq \ 
\frac{1}{\ds \log \frac{1}{\lambda}}
\]
or still, 
\[
\underset{r \ra \infty}{\limsupx} \ 
\frac{n (r)}{r^\rho}
\ \leq \ 
\frac{\tau + \varepsilon}{\ds \lambda^\rho \hsx \log \frac{1}{\lambda}}
\hsx .
\]
Setting $\lambda = e^{-1/\rho}$ then gives
\[
\underset{r \ra \infty}{\limsupx} \ 
\frac{n (r)}{r^\rho}
\ \leq \ 
\rho \hsy e \hsx (\tau + \varepsilon),
\]
so in the limit $(\varepsilon \ra 0)$
\[
\underset{r \ra \infty}{\limsupx} \ 
\frac{n (r)}{r^\rho}
\ \leq \ 
\rho \hsy e \hsy \tau.
\]
\\[-1cm]

\begin{x}{\small\bf REMARK} \ 
It follows that if $f$ has finite order and finite type, then 4.18 can be sharpened to

\[
n(r) 
\ = \ 
\tO (r^\rho).
\]
\\[-.25cm]
\end{x}


\chapter{
$\boldsymbol{\S}$\textbf{5}.\quad  CANONICAL PRODUCTS}
\setlength\parindent{2em}
\setcounter{theoremn}{0}
\renewcommand{\thepage}{\S5-\arabic{page}}

\qquad
Given a nonnegative integer $p$, let
\[
E(z, 0) 
\ = \ 
1 - z
\qquad (p = 0)
\]
and 
\[
E(z, p) 
\ = \ 
(1 - z) \hsx
\exp 
\bigg(
z + \frac{z^2}{2} + \cdots + \frac{z^p}{p}
\bigg)
\qquad (p > 0).
\]
\\[-1.25cm]

[Note:  \ 
The polynomial 
\[
z + \frac{z^2}{2} + \cdots + \frac{z^p}{p}
\]
is the $p^\nth$ partial sum of the expansion
\[
\log \frac{1}{1 - z}
\ = \ 
\sum\limits_{k = 1}^\infty \ 
\frac{z^k}{k}\hsx.]
\]
\\[-1cm]

\begin{x}{\small\bf DEFINITION} \ 
The functions $E(z, p)$ are called \un{primary factors}.
\\[-.25cm]
\end{x}

\begin{x}{\small\bf LEMMA} \ 
If $\abs{z} \leq 1$, then

\[
\abs{E(z, p) - 1} 
\ \leq \ 
\abs{z}^{p+1}.
\]
\\[-1.25cm]

PROOF \ 
Assuming that $p$ is positive, write

\[
E (z, p) 
\ = \ 
1 + 
\sum\limits_{n = 1}^\infty \ 
A_n z^n.
\]
Then
\[
E^\prime (z, p) 
\ = \ 
\sum\limits_{n = 1}^\infty \ 
n \hsy A_n \hsy z^{n-1}.
\]
Meanwhile, 

\[
E^\prime (z, p) 
\ = \ 
- z^p \hsx 
\exp \Big(z + \frac{z^2}{2} + \cdots + \frac{z^p}{p}\Big).
\]
Therefore

\[
A_1 = A_2 = \cdots = A_p = 0
\quad
\text{and} \quad 
A_n < 0 
\quad 
(n > p).
\]
On the other hand, $E(1, p) =0$, so

\[
\sum\limits_{n = p + 1}^\infty \ 
\abs{A_n} 
\ = \ 
1.
\]
Accordingly,
\allowdisplaybreaks\begin{align*}
\abs{z} \leq 1 
\implies 
\abs{E(z, p) - 1} \ 
&\leq \ 
\sum\limits_{n = p + 1}^\infty \ 
\abs{A_n} \hsx \abs{z}^n
\\[15pt]
&=\ 
\abs{z}^{p+1} \ 
\sum\limits_{n = p + 1}^\infty \ 
\abs{A_n} \hsx \abs{z}^{n - p - 1}
\\[15pt]
&\leq \ 
\abs{z}^{p+1}
\sum\limits_{n = p + 1}^\infty \ 
\abs{A_n}
\\[15pt]
&=\ 
\abs{z}^{p+1}.
\end{align*}
\\[-1cm]
\end{x}

Let 
$\{z_n: n = 1, 2, \ldots \}$ 
be a sequence of nonzero complex numbers with 

\[
0 < \abs{z_1} \leq \abs{z_2} \leq \cdots 
\qquad 
(\abs{z_n} \ra \infty),
\]
finite repetitions being permitted.  
Put $r_n = \abs{z_n}$ 
and assume that the convergence exponent $\kappa$ of the sequence 
$\{r_n: n = 1, 2, \ldots \}$ 
is finite.
\\[-.25cm]

Fix a nonnegative integer $p$ such that the series

\[
\sum\limits_{n = 1}^\infty \ 
\frac{1}{r_n^{p + 1}}
\]
is convergent.


\begin{x}{\small\bf NOTATION} \ 
Let
\[
P(z, p) 
\ = \ 
\prod\limits_{n = 1}^\infty \
E (\frac{z}{z_n}, p).
\]
\\[-1.25cm]
\end{x}

{\small\bf \un{N.B.}} \ 
At the origin, 

\[
P(0, p) 
\ = \ 
1.
\]

\begin{x}{\small\bf THEOREM} \ 
$P(z, p)$ is an entire function whose zeros are the $z_n$.
\\[-.5cm]

PROOF \ 
Taking into account 5.2, it is a question of applying 1.26 and 1.29.  
So consider the series

\[
\sum\limits_{n = 1}^\infty \
\Big(E \Big(\frac{z}{z_n}, p\Big) - 1\Big).
\]
Given $R > 0$, choose $N \gg 0$ : $n > N$ $\implies$ $\abs{z_n} > R$ $-$then for $\abs{z} \leq R$, 

\[
\Big|
E \Big(\frac{z}{z_n}, p\Big) - 1
\Big|
\ \leq \ 
\bigg|
\frac{z}{z_n}
\bigg|^{p + 1}
\ \leq \ 
\frac{R^{p + 1}}{\abs{z_n}^{p + 1}}
\]
and by assumption 

\[
\sum\limits_{n > N} \ 
\frac{1}{\abs{z}^{p + 1}}
\ < \ 
\infty.
\]
\\[-1cm]
\end{x}

\begin{x}{\small\bf LEMMA} \ 
For all complex $z$, if $p = 0$, 

\[
\log  \abs{E(z, 0)}
\ \leq \ 
\log ( 1 + \abs{z}),
\]
and if $p > 0$, 

\[
\log  \abs{E(z, p)}
\ \leq \ 
C_p \hsx 
\frac{\ \ \abs{z}^{p + 1}}{1 + \abs{z}},
\]
where $C_p = 3 \hsy e \hsy (2 + \log p)$.
\\[-.5cm]

PROOF \ 
The first inequality is trivial.  
To establish the second inequality, 
consider two cases.
\\[-.25cm]

\qquad \textbullet \quad 
$
\ds
\abs{z} \leq \frac{p}{p + 1}
$ $-$then
\allowdisplaybreaks
\allowdisplaybreaks\begin{align*}
\log  \abs{E(z, p)} \ 
&=\ 
\log  \abs{(E(z, p) - 1)+ 1}
\\[11pt]
&\leq\ 
\log \big(\abs{E(z, p) - 1} + 1\big)
\\[11pt]
&\leq\ 
\abs{E(z, p) - 1}
\\[11pt]
&\leq\ 
\abs{z}^{p + 1}
\qquad (\tcf. \ 5.2), 
\end{align*}
since $\log (x + 1) \leq x$ for $x \geq 0$.
\\[-.25cm]

\qquad \textbullet \quad
$
\ds
\abs{z} > \frac{p}{p + 1}
$ $-$then
\allowdisplaybreaks
\allowdisplaybreaks\begin{align*}
\log \abs{E(z, p)} \ 
&\leq\ 
2 \hsy \abs{z} 
+ \hsy \frac{\abs{z}^2}{2}
+ \cdots + 
\frac{\abs{z}^p}{p} 
\\[11pt]
&=\ 
\abs{z}^p \hsy
\Big(
\frac{1}{p} + \frac{1}{p-1} \hsx \frac{1}{\abs{z}} 
+ \cdots + 
\frac{1}{2} \frac{1}{\abs{z}^{p - 2}} 
+ 2 \hsy \frac{1}{\abs{z}^{p-1}}
\Big)
\\[11pt]
&\leq\ 
\abs{z}^p \hsy
\Big(
\frac{p + 1}{p}
\Big)^{p-1}
\Big(
2 + \frac{1}{2} + \cdots + \frac{1}{p}
\Big)
\\[11pt]
&\leq\ 
\abs{z}^p \hsy
\Big(
1 + \frac{1}{p}
\Big)^p
\hsy 
\bigg(
2 + \hsx 
\int\limits_1^p \ 
\frac{ \td t}{t}
\bigg)
\\[11pt]
&\leq\ 
\abs{z}^p
\hsx
e \hsy (2 + \log p) 
\\[11pt]
&=\ 
e \hsy (2 + \log p) 
\hsx
\abs{z}^p
\hsx
\frac{1 + \abs{z}}{1 + \abs{z}}
\\[11pt]
&=\ 
e \hsy (2 + \log p) 
\hsx
\Big(1 + \frac{1}{\abs{z}}
\Big)
\hsx 
\frac{\abs{z}^{p + 1}}{1 + \abs{z}}
\\[11pt]
&\leq\ 
3 \hsx e \hsy (2 + \log p) 
\hsx 
\frac{\abs{z}^{p + 1}}{1 + \abs{z}}
\\[11pt]
&=\ 
C_p \hsx 
\frac{\abs{z}^{p + 1}}{1 + \abs{z}},
\end{align*}
\\[-.25cm]
since

\[
1 + \frac{1}{\abs{z}} 
\ < \ 
1 + \frac{p + 1}{p}
\ = \ 
1 + 1 + \frac{1}{p} 
\ \leq \ 
3.
\]
\\[-.25cm]
\end{x}

\begin{x}{\small\bf SUBLEMMA} \ 
We have 

\[
\lim\limits_{r \ra \infty} \ 
\frac{n(r)}{r^{p + 1}}
\ = \ 
0.
\]
\\[-.5cm]

PROOF \
In fact, 
\allowdisplaybreaks\begin{align*}
\sum\limits_{n = 1}^\infty \ 
\frac{1}{r_n^{p + 1}} \ 
&=\ 
\int\limits_0^\infty \ 
\frac{\td n(t)}{t^{p + 1}}
\\[15pt]
&=\ 
\lim\limits_{r \ra \infty} \ 
\frac{n(r)}{r^{p + 1}}
\hsx + \hsx 
(p + 1) \ 
\int\limits_0^\infty \ 
\frac{n(t)}{t^{p + 2}} \ \td t .
\end{align*}
And
\allowdisplaybreaks\begin{align*}
\frac{n(r)}{r^{p + 1}} \ 
&=\ 
(p + 1) \hsy n(r) \ 
\int\limits_r^\infty \ 
\frac{ \td t}{t^{p + 2}} \
\\[15pt]
&\leq 
(p + 1) \ 
\int\limits_r^\infty \ 
\frac{n(t)}{t^{p + 2}} \ \td t
\\[15pt]
&\ra 
0 
\qquad (r \ra \infty).
\end{align*}
\\[-1cm]
\end{x}

\begin{x}{\small\bf LEMMA} \ 
Put $r = \abs{z}$ $-$then for $p = 0$, 

\[
\log \abs{P(z, 0)} 
\ \leq \ 
\int\limits_0^r \ 
\frac{n(t)}{t} \ \td t
\hsx + \hsx 
r \ 
\int\limits_r^\infty \ 
\frac{n(t)}{t^{2}} \ \td t, 
\]
and for $p > 0$, 

\[
\log \abs{P(z, p)} 
\ \leq \ 
(p + 1) \hsx C_p \hsy r^p \ 
\bigg(
\int\limits_0^r \ 
\frac{n(t)}{t^{p + 1}} \ \td t
\hsx + \hsx 
r \ 
\int\limits_r^\infty \ 
\frac{n(t)}{t^{p + 2}} \ \td t
\bigg).
\]
\\[-1cm]

PROOF \
If $p = 0$, 
\allowdisplaybreaks
\allowdisplaybreaks\begin{align*}
\log \abs{P(z, 0)} \ 
&\leq \ 
\sum\limits_{n = 1}^\infty \ 
\log \bigg(1 + \frac{r}{r_n}\bigg) 
\qquad (\tcf. \ 5.5)
\\[15pt]
&=\ 
\int\limits_0^\infty \ 
\log \bigg(1 + \frac{r}{t}\bigg) 
\ \td n(t)
\\[15pt]
&=\ 
\log \Big(1 + \frac{r}{t}\Big) \hsy n(t) \bigg|_0^\infty 
\hsx + \hsx 
r \hsx 
\int\limits_0^\infty \ 
\frac{n (t)}{t (t + r)}
\ \td t
\\[15pt]
&=\ 
\log \Big(1 + \frac{r}{t}\Big)\hsy t \hsy \frac{n(t)}{t} \bigg|_0^\infty 
\hsx + \hsx 
r \hsx 
\int\limits_0^\infty \ 
\frac{n (t)}{t (t + r)}
\ \td t
\\[15pt]
&=\ 
r \hsx 
\int\limits_0^\infty \ 
\frac{n (t)}{t (t + r)}
\ \td t
\\[15pt]
&\leq \ 
\int\limits_0^r \ 
\frac{n(t)}{t} \ \td t
\hsx + \hsx 
r \ 
\int\limits_r^\infty \ 
\frac{n(t)}{t^{2}} \ \td t
\end{align*}
and if $p > 0$, 
\allowdisplaybreaks
\allowdisplaybreaks\begin{align*}
\log \abs{P(z, p)} \ 
&\leq \ 
C_p \ 
\sum\limits_{n = 1}^\infty \ 
\frac{r^{p + 1}}{r_n^p (r + r_n)} 
\qquad (\tcf. \ 5.5)
\\[15pt]
&=\ 
C_p \hsx r^{p + 1} \ 
 \int\limits_0^\infty \ 
 \frac{\td n (t)}{t^p (t + r)}
 \\[15pt]
&=\ 
C_p \hsx r^{p + 1} \ 
\frac{n (t)}{t^p (t + r)} \bigg|_0^\infty
\\[15pt]
&
\hspace{1cm}
\hsx + \hsx 
C_p \hsx r^{p + 1} \ 
\int\limits_0^\infty \ 
\Big(
\frac{p}{t^{p+1} (t + r)}
\hsx + \hsx 
\frac{1}{t^p (t + r)^2}
\Big)
n(t) 
\ \td t
\\[15pt]
&=\ 
C_p \hsy r^{p + 1} \hsy 
\frac{n (t)}{t^{p+1} (1 + r/t)} \bigg|_0^\infty
\\[15pt]
&
\hspace{1cm}
\hsy + \hsy 
C_p \hsx r^{p + 1} \hsy 
\int\limits_0^\infty \hsy
\Big(
\frac{p}{t^{p+1} (t + r)}
\hsy + \hsy 
\frac{1}{t^p (t + r)^2}
\Big)
n(t) 
\hsy \td t
\\[15pt]
&=\ 
C_p \hsx r^{p + 1} \ 
\int\limits_0^\infty \ 
\Big(
\frac{p}{t^{p+1} (t + r)}
\hsy + \hsy 
\frac{1}{t^p (t + r)^2}
\Big)
n(t) 
\hsy \td t
\\[15pt]
&=\ 
C_p \hsx r^{p + 1} \ 
\bigg(
\int\limits_0^r \ 
\hsx + \hsx
\int\limits_r^\infty \ 
\bigg)
\bigg(
\frac{p}{t^{p+1} (t + r)}
\hsy + \hsy 
\frac{1}{t^p (t + r)^2} 
\bigg)
n(t) \ \td t
\\[15pt]
&\leq\ 
(p + 1) \hsy C_p \hsx r^p \ 
\bigg(
\int\limits_0^r \ 
\frac{n(t)}{t^{p+1}} \ \td t
\hsx + \hsx 
r \hsx
\int\limits_r^\infty \ 
\frac{n(t)}{t^{p+2}} \ \td t
\bigg).
\end{align*}
\\[-.25cm]
\end{x}

\begin{x}{\small\bf REMARK} \ 
For use below, note that these inequalities involve $z$ only through its modulus $r$, hence provide estimates for
\[
\log M (r; P(z, p)).
\]
\\[-1.25cm]

It has been assumed from the outset that the convergence exponent $\kappa$ of the sequence 
$\{r_n : n = 1, 2, \ldots\}$ 
is finite, 
thus it makes sense to take $p = \fg$, 
the genus of the sequence
$\{r_n : n = 1, 2, \ldots\}$ 
(cf. 4.14).
\\[-.25cm]
\end{x}

\begin{x}{\small\bf DEFINITION} \ 
\[
P (z, \fg) 
\ = \ 
\prod\limits_{n = 1}^\infty \ 
E 
\Big(
\frac{z}{z_n}, \fg
\Big)
\]
is called the \un{canonical product} formed from the $z_n$.
\\[-.5cm]

[Note: \ 
$P(z, \fg)$ is a transcendental entire function and the infinite product defining 
$P(z, \fg)$ is absolutely convergent (cf. 5.4).]
\\[-.25cm]
\end{x}

\begin{x}{\small\bf THEOREM} \ 
The order $\rho$ of $P(z, \fg)$ is equal to $\kappa$.
\\[-.5cm]

PROOF \ 
It suffices to show that $\rho \leq \kappa$, hence is finite 
(for then, on general grounds, $\kappa \leq \rho$ (cf. 4.20)).  
In any event, 

\[
\fg 
\ \leq \ 
\kappa 
\ \leq \ 
\fg + 1
\qquad (\tcf. \ 4.15)
\]
and it will be assumed that $\fg$ is positive.
\\[-.25cm]

\qquad 
\un{Case 1:} \ 
$\kappa 
\ < \ 
\fg + 1$.  
Choose $\varepsilon > 0$ : $\kappa + \varepsilon < \fg + 1$ $-$then

\[
n(t) 
\ < \ 
t^{\kappa +\varepsilon} 
\qquad (t \gg 0) 
\quad (\tcf. \ 4.10), 
\]
so
\allowdisplaybreaks\begin{align*}
\log M (r; P(z, \fg))\ 
&\leq \ 
(\fg + 1) \hsx C_\fg \hsy r^\fg \hsy 
\bigg(
\tO(1) 
\hsx + \hsx 
\int\limits_0^r \hsy 
t^{\kappa + \varepsilon - \fg - 1} 
\ \td t
\hsx + \hsx 
r \ 
\int\limits_r^\infty \ 
t^{\kappa + \varepsilon - \fg - 2} 
\ \td t
\bigg)
\\[15pt]
&\leq \ 
(\fg + 1) \hsx C_\fg \hsy r^\fg \hsy 
\bigg(
\tO(1) 
\hsx + \hsx 
\frac{r^{\kappa + \varepsilon - \fg}}{\kappa + \varepsilon - \fg}
\hsx + \hsx 
\frac{r^{\kappa + \varepsilon - \fg}}{\fg + 1 - \kappa - \varepsilon}
\bigg)
\\[15pt]
&< \ 
r^{\kappa + 2\varepsilon}
\qquad (r \gg 0).
\end{align*}
Therefore $\rho \leq \kappa$. 
\\[-.25cm]

\qquad 
\un{Case 2:} \ 
$\kappa = \fg + 1$.  
Owing to 5.6, 

\[
\lim\limits_{r \ra \infty } \  
\frac{n(r)}{r^{\fg + 1}}
\ = \ 
0.
\]
Fix $\varepsilon > 0$ and choose $r_0$ : 

\[
r > r_0 
\hsx \implies \  
\begin{cases}
\ds \ 
\frac{n(r)}{r^{\fg + 1}} 
\hspace{1.5cm} < \varepsilon
\\[15pt]
\ds \ 
\int\limits_r^\infty \ 
\frac{n(t)}{t^{\fg + 2}} 
\ \td t
\hspace{.3cm} < \varepsilon
\end{cases}
.
\]
Then
\allowdisplaybreaks\begin{align*}
\log M (r; P(z, \fg))\ 
&\leq \ 
(\fg + 1) \hsx C_\fg \hsy r^\fg \hsy 
\bigg(
r \hsx \frac{n(r)}{r^{\fg + 1}} + r \hsy \varepsilon
\bigg)
\\[11pt]
&\leq \ 
(\fg + 1) \hsx C_\fg \hsy r^\fg \hsy (r \varepsilon + r \varepsilon)
\\[11pt]
&= \ 
2 \hsx (\fg + 1) \hsx C_\fg \hsy \varepsilon \hsy r^{\fg + 1}
\\[11pt]
&= \ 
2 \hsx (\fg + 1) \hsx C_\fg \hsy \varepsilon \hsy r^\kappa.
\end{align*}
Restated: \ 
$\forall \ C > 0$, 

\[
\log M (r; P(z, \fg))
\ \leq \ 
C \hsy r^\kappa
\qquad (r \gg 0).
\]
Therefore $\rho \leq \kappa$ (and more (cf. 5.16)).
\\[-.5cm]

[Note: \ 
The discussion when $\fg = 0$ is similar but simpler.]
\\[-.25cm]
\end{x}

\begin{x}{\small\bf LEMMA} \ 
Let $Q$ be a polynomial of degree $q$ and put

\[
f(z) 
\ = \ 
e^{Q(z)} 
\hsx 
P(z, \fg).
\]
Then 

\[
\rho(f) 
\ = \ 
\max (q, \kappa).
\]
\\[-1.25cm]

PROOF \ 
Since $q$ equals the order of $e^Q$ and since $\kappa$ equals the order of $P(z, \fg)$, 
it follows from 2.34 that

\[
\rho(f) 
\ \leq \ 
\max (q, \kappa).
\]
On the other hand, $\kappa \leq \rho (f)$, (cf. 4.20).  
And

\allowdisplaybreaks\begin{align*}
\frac{f}{P} = e^Q 
\implies q = \rho (e^Q) 
&\leq \ 
\max (\rho (f) , \kappa) 
\qquad (\tcf. \ 2.37)
\\[8pt]
&= \ 
\rho(f).
\end{align*}
Therefore

\[
\max (q, \kappa) 
\ \leq \ 
\rho (f).
\]
\\[-1.25cm]

[Note: \ 
It is a corollary that if $\rho (f)$ is not an integer, then $\rho(f) = \kappa$.]
\\[-.25cm]
\end{x}

\begin{x}{\small\bf EXAMPLE} \ 
The canonical product

\[
\{
(1 - z) \hsy e^z
\}
\hsy 
\{
(1 + z) \hsy e^{-z}
\}
\hsy 
\Big\{
\Big(1 - \frac{z}{2}\Big) \hsy e^{z/2}
\Big\}
\hsy 
\Big\{
\Big(1 + \frac{z}{2}\Big) \hsy e^{-z/2}
\cdots
\Big\}
\]
represents

\[
\frac{\sin \pi \hsy z}{\pi \hsy z}
\qquad (\tcf. \ 1.23).
\]
\\[-1cm]
\end{x}

\begin{x}{\small\bf EXAMPLE} \ 
The reciprocal

\[
\frac{1}{z \hsy \Gamma (z)}
\ = \ 
e^{\gamma \hsy z} \ 
\prod\limits_{n = 1}^\infty \ 
\bigg(1 + \frac{z}{n}\bigg) \hsx \exp \bigg(- \frac{z}{n}\bigg)
\]
is a transcendental entire function of order 1.  
To see this, take $z_n = -n$ $(n = 1, 2, \ldots)$ 
$-$then $\kappa = 1$ and $\fg = 1$ (cf. 4.15).  
In view of 5.10, the order of the canonical product

\[
\prod\limits_{n = 1}^\infty \ 
\bigg(1 + \frac{z}{n}\bigg) \hsx \exp \bigg(- \frac{z}{n}\bigg)
\]
is 1, as is the order of $e^{\gamma \hsy z}$.  
Therefore the order of 
$
\ds
\frac{1}{z \hsy \Gamma (z)}
$ 
equals
\[
\max (1, 1) 
\ = \ 
1
\qquad (\tcf. \ 5.11).
\]
\\[-1.25cm]
\end{x}

\begin{x}{\small\bf EXAMPLE} \ 
Let $\omega_1$, $\omega_2$ be two nonzero complex constants whose ratio is not purely real.  
Put

\[
\Omega_{m, n}
\ = \ 
m \hsy \omega_1 + n  \hsy \omega_2
\qquad ((m, n) \neq (0,0))
\]
and consider

\[
\prod\limits_{m, n} \ 
\bigg(
1 - \frac{z}{\Omega_{m, n}}
\bigg)
\hsy 
\exp
\bigg(
\frac{z}{\Omega_{m, n}} + \frac{1}{2} \hsx \big(\frac{z}{\Omega_{m, n}}\big)^2
\bigg).
\]
Then here, $\kappa = 2$ and $\fg = 2$ (cf. 4.15).  
Setting 

\[
\sigma (\restr{z}{\omega_1, \omega_2})
\ = \ 
\prod\limits_{m, n} \ 
\cdots \hsy,
\]
it follows that $\sigma (\restr{z}{\omega_1, \omega_2})$ is a transcendental entire function of order 2.
\\[-.25cm]
\end{x}

The proof of 5.10 fell into two cases: 
\[
\kappa \ < \  \fg + 1
\quad \text{or} \quad 
\kappa \ = \  \fg + 1.
\]
\\[-1.25cm]

\begin{x}{\small\bf RAPPEL} \ 
(cf. 4.15)
\\[-.25cm]

\qquad\qquad \textbullet \quad
If $\kappa$ is not an integer, then $\fg = [\kappa]$.
\\[-.5cm]

\qquad\qquad \textbullet \quad
If $\kappa$ is an integer, then 
$\fg = \kappa - 1$ if \ 
$
\ds 
\sum\limits_{n = 1}^\infty \ \frac{1}{\abs{z_n}^\kappa}
$
is convergent, while 
$\fg = \kappa$ 
if \ 
$
\ds 
\sum\limits_{n = 1}^\infty \ \frac{1}{\abs{z_n}^\kappa}
$
is divergent.
\\[-.25cm]

[Note: \ 
Employing the terminology of 4.24, in this situation

\[
\begin{cases}
\ 
\text{$P(z, \fg)$ of convergence class $\implies \fg = \kappa - 1$}
\\[4pt]
\ 
\text{$P(z, \fg)$ of divergence class $\implies \fg = \kappa$}
\end{cases}
.]
\]
\\[-1cm]

So, if $\kappa$ is not an integer, then $\kappa < \fg + 1$ 
and if $\kappa$ is an integer, then $\kappa < \fg + 1$ if 
$
\ 
\ds 
\sum\limits_{n = 1}^\infty \ \frac{1}{\abs{z_n}^\kappa}
$
is divergent but $\kappa = \fg + 1$ if 
$
\ 
\ds 
\sum\limits_{n = 1}^\infty \ \frac{1}{\abs{z_n}^\kappa}
$
is convergent.
\\

With these points in mind, we shall now proceed to the determination of the type $\tau$ of $P(z, \fg)$.
\\[-.5cm]

[Note: \ 
The very definition of type requires that $0 < \rho < \infty$.  
It is automatic that $\rho$ is finite and it is also automatic that $\rho$ is positive if $\kappa$ is not an integer 
or if $\kappa$ is an integer and $\fg = \kappa - 1$ but if $\kappa$ is an integer and $\fg = \kappa$, 
then it will be assumed that $\kappa (= \rho)$ is positive.]
\\[-.25cm]
\end{x}

\begin{x}{\small\bf THEOREM} \ 
If $\kappa$ is not an integer and if 
$
\ 
\ds 
\sum\limits_{n = 1}^\infty \ \frac{1}{\abs{z_n}^\kappa}
$
is convergent, then $P(z, \fg)$ is of minimal type. 
\\[-.5cm]

[Here $\kappa = \fg + 1$, thus the assertion is implied by the ``Case 2'' analysis in 5.10.]
\\[-.25cm]
\end{x}

\begin{x}{\small\bf LEMMA} \ 
Take $\rho > 0$ $-$then 

\[
\Delta 
\ \leq \ 
e^\rho \hsy \tau.
\]
\\[-1cm]

PROOF  \ 
Since $P(0, \fg) = 1$, in view of 4.19, 

\[
n(r) 
\ \leq \ 
\log M (e \hsy r; P(z, \fg)),
\]
thus
\[
\frac{n(r)}{r^\kappa}
\ \leq \ 
\frac{\log M (e \hsy r; P(z, \fg))}{r^\kappa}
\]

\qquad 
$\implies$
\allowdisplaybreaks
\allowdisplaybreaks\begin{align*}
\Delta \ 
&=\ 
\underset{r \ra \infty}{\limsupx} \ 
\frac{n(r)}{r^\kappa}
\qquad (\tcf. \ 4.13)
\\[15pt]
&\leq \ 
\underset{r \ra \infty}{\limsupx} \ 
e^\kappa \ 
\frac{\log M (e \hsy r; P(z, \fg))}{(e \hsy r)^\kappa}
\\[15pt]
&=\  
\underset{r \ra \infty}{\limsupx} \ 
e^\rho \ 
\frac{\log M (e \hsy r; P(z, \fg))}{(e \hsy r)^\rho}
\\[15pt]
&=\  
e^\rho \ 
\underset{r \ra \infty}{\limsupx} \ 
\frac{\log M (e \hsy r; P(z, \fg))}{(e \hsy r)^\rho}
\\[15pt]
&=\  
e^\rho \hsy \tau.
\end{align*}
\\[-1.25cm]
\end{x}

Suppose that $\kappa$ is not an integer (hence $\rho > 0$ and 
$\fg < \kappa < \fg + 1$).
\\[-.25cm]

\begin{x}{\small\bf LEMMA} \ 
Put 

\[
K_{0, \kappa} 
\ = \ 
\frac{1}{\kappa} + \frac{1}{1 - \kappa}
\]
and 
\[
K_{\fg, \kappa} 
\ = \ 
(\fg + 1) \hsy C_\fg \ 
\bigg[
\frac{1}{\kappa - \fg} + \frac{1}{\fg + 1 - \kappa}
\bigg]
\qquad (\fg > 0).
\]
Then
\[
\tau 
\ \leq \ 
2 \hsy K_{\fg, \kappa}  \hsy \Delta.
\]
\\[-1cm]

PROOF \ 
Given $\varepsilon > 0$, we have 

\[
n(t) 
\ < \ 
(\Delta + \varepsilon) t^\kappa 
\qquad (t \gg 0).
\]
Therefore, taking $\fg > 0$, 

\allowdisplaybreaks\begin{align*}
&
\log M (r; P(z, \fg)) \ 
\\[15pt]
&
\hspace{1cm}
\leq \ 
(\fg + 1) \hsy C_\fg \hsy r^\fg \ 
\bigg(
\int\limits_0^r \ 
\frac{n(t)}{t^{\fg +1}}\ 
\ \td t 
\hsx + \hsx 
r\hsx 
\int\limits_r^\infty\ 
\frac{n(t)}{t^{\fg + 2}} 
\ \td t
\bigg)
\qquad (\tcf. \ 5.7)
\\[15pt]
&
\hspace{1cm}
\leq \ 
(\fg + 1) \hsy C_\fg \hsy r^\fg \ 
\bigg(
\tO (1) + (\Delta + \varepsilon) \ 
\int\limits_0^r \ 
t^{\kappa - \fg - 1} 
\ \td t
\hsx + \hsx 
(\Delta + \varepsilon) r \ 
\int\limits_r^\infty\ 
t^{\kappa - \fg - 2} 
\ \td t
\bigg)
\\[15pt]
&
\hspace{1cm}
\leq \  
(\fg + 1) \hsy C_\fg \hsy r^\fg \ 
\bigg(
\tO (1) + (\Delta + \varepsilon) \ 
\frac{r^{\kappa - \fg}}{\kappa - \fg} 
\hsx + \hsx 
(\Delta + \varepsilon) 
\frac{r^{\kappa - \fg}}{\fg + 1 - \kappa} 
\\[15pt]
&
\hspace{1cm}
< \ 
2 \hsy K_{\fg, \kappa} \hsy (\Delta + \varepsilon) \hsy r^\kappa 
\qquad (r \gg 0).
\end{align*}
Since $\rho = \kappa$ it follows that 

\[
\underset{r \ra \infty}{\limsupx} \ 
\frac{\log M (r; P(z, \fg))}{r^\rho}
\ \leq \ 
K_{\fg, \kappa} \hsy (\Delta + \varepsilon), 
\]
i.e., 
\[
\tau 
\ \leq \ 
2 \hsy K_{\fg, \kappa} \hsy \Delta.
\]
\\[-1.25cm]

[Note: \ 
The discussion when $\fg = 0$ is similar but simpler.]
\\[-.25cm]
\end{x}

\begin{x}{\small\bf THEOREM} \ 
If $\kappa$ is not an integer, then $P(z, \fg)$ is of maximal, minimal, or intermediate type 
according to whether 
$\Delta = \infty$, 
$\Delta = 0$, 
or 
$0 < \Delta < \infty$ 
and conversely.
\\[-.5cm]

[This is implied by 5.17 and 5.18.]
\\[-.25cm]
\end{x}

There remains the case when $\kappa$ is an integer $> 0$ and 
$
\ds 
\sum\limits_{n = 1}^\infty \ 
\frac{1}{\abs{z_n}^\kappa}
$
is divergent (hence $\fg = \kappa$).  
To this end, let 

\[
\delta (r) 
\ = \ 
\bigg|
\frac{1}{\kappa} \ 
\sum\limits_{\abs{z_n} < r} \ 
z_n^{-\kappa}
\bigg|,
\]
put 
\[
\delta
\ = \ 
\underset{r \ra \infty}{\limsupx} \ 
\delta (r) ,
\]
and set
\[
\Gamma 
\ = \ 
\max (\delta, \Delta).
\]

\begin{x}{\small\bf THEOREM} \ 
Under the preceding conditions, $P(z, \fg)$ is of maximal, minimal, or intermediate type 
according to whether 
$\Gamma = \infty$, 
$\Gamma = 0$, 
or 
$0 < \Gamma < \infty$ 
and conversely.
\\[-.25cm]

The proof can be divided into two parts.
\\[-.25cm]

\qquad \textbullet \ \quad
$\exists \ C > 1$: 
\[
\Gamma 
\ \leq \ 
C \hsy e^\rho \hsy \tau.
\]
\\[-1.25cm]

[First, it can be shown that for some $C > 1$, 

\[
\delta (r)
\ < \ 
C \ 
\frac{\log M (e \hsy r; P(z, \fg))}{r^\kappa}
\qquad (r \gg 0).
\]
Thus 

\[
\delta (r)
\ < \ 
(C \hsy e^\rho)\ 
\frac{\log M (e \hsy r; P(z, \fg))}{(e \hsy r)^\rho}
\qquad (r \gg 0)
\]
and so
\[
\delta 
\ \leq \ 
C \hsy e^\rho \hsy \tau.
\]
Meanwhile, 
\[
\Delta 
\ \leq \ 
e^\rho \hsy \tau
\qquad (\tcf. \ 5.17).
\]
Therefore
\[
\Gamma 
\ \leq \ 
C \hsy e^\rho \hsy \tau.
\]
\\[-1cm]

\qquad \textbullet \quad
$\exists \ K > 0$: 
\[
\tau 
\ \leq \ 
K \hsy \Gamma.
\]
\\[-1.25cm]

[Write
\[
P(z, \fg)
\ = \ 
\exp
\bigg(
\big(
\frac{1}{\kappa} \ 
\sum\limits_{\abs{z_n} < r} \ 
z_n^{-\kappa} 
\big) 
z^\kappa
\bigg)
\hsx \times \hsx 
\prod\limits_{\abs{z_n} < r} \ 
E
\bigg(
\frac{z}{z_n}, \fg - 1
\bigg)
\prod\limits_{\abs{z_n} < r} \ 
E
\bigg(
\frac{z}{z_n}, \fg 
\bigg),
\]
where $r = \abs{z}$ and take $\kappa > 1$ $-$then
\allowdisplaybreaks\begin{align*}
\log M (r; P(z, \fg))\ 
&\leq \ 
\delta (r) r^\kappa 
\hsx + \hsx 
C_\fg\hsx
\bigg(
r^\fg \hsx 
\int\limits_0^r \ 
\frac{\td n(t)}{t^{\fg -1} (t + r)} 
\hsx + \hsx 
r^{\fg + 1} \hsx 
\int\limits_r^\infty\ 
\frac{\td n(t)}{t^{\fg} (t + r)} 
\bigg)
\\[15pt]
&\leq \ 
\delta (r) r^\kappa 
\hsx + \hsx 
(\fg + 1)
C_\fg\hsx
\bigg(
r^{\fg - 1} \hsx 
\int\limits_0^r \ 
\frac{n(t)}{t^\fg} 
\hsx + \hsx 
r^{\fg + 1} \hsx 
\int\limits_r^\infty\ 
\frac{n(t)}{t^{\fg + 2}}  
\ \td t
\bigg).
\end{align*}
But $\forall \ \varepsilon > 0$, 
\[
n(t)
\ < \ 
(\Delta + \varepsilon) \hsy t^\kappa 
\qquad (t \gg 0).
\]
Therefore 

\[
\log M (r; P(z, \fg))
\ \leq \ 
\delta (r) r^\kappa 
\hsx + \hsx 
2 \hsy (\fg + 1) C_\fg \hsy (\Delta + \varepsilon) r^\kappa 
\qquad (r \gg 0).
\]
And finally
\allowdisplaybreaks\begin{align*}
\tau \ 
&=\ 
\underset{r \ra \infty}{\limsupx} \ 
\frac{\log M (r; P(z, \fg))}{r^\kappa}
\\[15pt]
&\leq \
\delta + 2 \hsy (\fg + 1) C_\fg \Delta
\\[15pt]
&\leq \
\Gamma + 2 \hsy (\fg + 1) C_\fg \hsy \Gamma
\\[15pt]
&= \
(1 + 2 \hsy (\fg + 1) C_\fg) \Gamma
\\[15pt]
&\equiv \
K \hsy \Gamma.
\end{align*}
\\[-1.25cm]

[Note: \ 
Minor modifications in the argument are needed if $\kappa = 1$.]
\\[-.25cm]
\end{x}

\begin{spacing}{1.65} 
\begin{x}{\small\bf EXAMPLE} \ 
In the setup of 5.12, the zeros are $\pm n$ $(n = 1, 2, \ldots)$, say 
$z_1 = 1$,  
$z_2 = -1$,  
$z_3 = 2$,  
$z_4 = -2$,  
$\ldots, $
hence 
$r_1 = 1$, 
$r_2 = 1$, 
$r_3 = 2$, 
$r_4 = 2, \ldots$
. 
Here $\kappa = 1$ and 
$
\ds
\frac{\sin \pi \hsy z}{\pi \hsy z}
$ 
is of divergence class.  
Moreover, 
\\[-1.5cm]
\end{x}
\end{spacing}
\[
\delta (r) \ = 0 \
 \ (r > 0) \ 
\ \implies \ 
\delta = 0.
\]
On the other hand, 

\[
\Delta
\ = \ 
\underset{n \ra \infty}{\limsupx} \ 
\frac{n}{r_n}
\qquad (\tcf. \ 4.11).
\]
But

\[
\frac{1}{r_1}
\ = \ 
\frac{1}{1}, 
\quad 
\frac{2}{r_2}
\ = \ 
\frac{2}{1}, 
\quad 
\frac{3}{r_3}
\ = \ 
\frac{3}{2}, 
\quad 
\frac{4}{r_4}
\ = \ 
\frac{4}{2}, 
\ldots \hsx . 
\]
Therefore $\Delta = 2$, and 

\[
\Gamma 
\ = \ 
\max (\delta, \Delta) 
\ = \ 
\max (0, 2) 
\ = \ 
2.
\]
I.e.: \ 
$
\ds
\frac{\sin \pi \hsy z}{\pi \hsy z}
$ 
is of intermediate type.
\\[-.25cm]

\begin{x}{\small\bf EXAMPLE} \ 
In the setup of 5.13, the zeros are $-n$ $(n = 1, 2, \ldots)$, say 
$z_n = -n$.  
Here $\kappa = 1$ and 
$
\ds \frac{1}{z \hsy\Gamma (z)}
$
is of divergence class.  
However, in contrast with 5.21, 

\[
\delta
\ = \ 
\underset{n \ra \infty}{\limsupx} \ 
\bigg(
1 + \frac{1}{2} + \cdots +  \frac{1}{n}
\bigg)
\ = \ 
\infty.
\]
Since it is clear that $\Delta = 1$, we thus have

\[
\Gamma 
\ = \ 
\max (\delta, \Delta) 
\ = \ 
\max (\infty, 1) 
\ = \ 
\infty.
\]
Consequently, 

\[
\prod\limits_{n = 1}^\infty \ 
\bigg(
1 + \frac{z}{n} \hsx 
\bigg)
\exp 
\bigg(
-\frac{z}{n}
\bigg)
\]
is of maximal type.  
But the order of $e^{\gamma \hsy z}$ is 1 and the type of $e^{\gamma \hsy z}$ is $\gamma$.  
An appeal to 3.14 then implies that

\[
\frac{1}{z \hsy\Gamma (z)}
\ = \ 
\prod\limits_{n = 1}^\infty \ 
\bigg(
1 + \frac{z}{n} \hsx 
\bigg)
\exp 
\bigg(
-\frac{z}{n}
\bigg)
\]
is of maximal type.
\\[-.25cm]
\end{x}


\chapter{
$\boldsymbol{\S}$\textbf{6}.\quad  EXPONENTIAL FACTORS}
\setlength\parindent{2em}
\setcounter{theoremn}{0}
\renewcommand{\thepage}{\S6-\arabic{page}}

\qquad 
Take a canonical product $P(z, \fg)$ per \S5, 
let $Q$ be a polynomial of degree $q \geq 1$ and put
\[
f(z)
\ = \ 
e^{Q(z)} \hsx P(z, \fg).
\]
Then
\[
\rho \quad (= \rho(f))
\ = \ 
\max (q, \kappa) 
\qquad (\tcf. \ 5.11).
\]
\\[-1cm]

[Note: \ 
Recall that it is always true that $\kappa \leq \rho$ (cf. 4.20).]
\\[-.25cm]

\begin{x}{\small\bf DEFINITION} \ 
The \un{genus} of $f$ is the nonnegative integer
\[
\gen f 
\ = \ 
\max (q, \fg) .
\]
\\[-1.5cm]
\end{x}

\begin{x}{\small\bf LEMMA} \ 
We have
\[
\gen f 
\ \leq \ 
\rho.
\]
\\[-1.25cm]

[This is because $\fg \leq \kappa$ (cf. 5.15).]
\\[-.25cm]
\end{x}

\begin{x}{\small\bf LEMMA} \ 
If $\rho$ is not an integer, then the genus of $f$ is $[\rho]$.
\\[-.5cm]

PROOF \ 
For here $\rho = \kappa$ (and $\rho > q)$.  
But in general, 
\[
\fg 
\ \leq \ 
\kappa 
\ \leq \ 
\fg + 1, 
\]
so in this case
\[
\fg 
\ < \ 
\rho 
\ < \ 
\fg + 1, 
\]
thus 
\[
\gen f 
\ = \ 
\max (q, \fg) 
\ = \ 
\max (q, [\rho]) 
\ = \ 
[\rho].
\]
\\[-1.25cm]
\end{x}

\begin{x}{\small\bf LEMMA} \ 
If $\rho$ is an integer, then the genus of $f$ is either equal to $\rho$ or to $\rho - 1$.
\\[-.5cm]

PROOF \ 
The genus of $f$ is necessarily less than or equal to $\rho$ (cf. 6.2).  
If
it is less than $\rho$, then $q < \rho$ $(\implies q \leq \rho - 1)$ and $\rho = \kappa$, hence

\[
\fg 
\ \leq \ 
\rho 
\ \leq \ 
\fg + 1.
\]
But by assumption, $\fg < \rho$.  Therefore $\fg = \rho - 1$ and 

\[
\gen f 
\ = \ 
\max (q, \fg) 
\ = \ 
\max (q, \rho - 1) 
\ = \ 
\rho - 1.
\]
\\[-1.25cm]
\end{x}

\begin{x}{\small\bf REMARK} \ 
When $\rho$ is an integer, there are five possibilities.
\\[-.25cm]

\qquad (i) \quad 
$\kappa < \rho$, 
$\fg \leq \kappa$, 
$q = \rho$, 
$\gen f = \rho$
\\[-.5cm]

\qquad (ii) \quad 
$\kappa = \rho$, 
$\fg = \rho$, 
$q = \rho$, 
$\gen f = \rho$
\\[-.5cm]

\qquad (iii) \quad 
$\kappa = \rho$, 
$\fg = \rho$, 
$q < \rho$, 
$\gen f = \rho$
\\[-.5cm]

\qquad (iv) \quad 
$\kappa = \rho$, 
$\fg = \rho - 1$, 
$q = \rho$, 
$\gen f = \rho$
\\[-.5cm]

\qquad (v) \quad 
$\kappa = \rho$, 
$\fg = \rho - 1$, 
$q < \rho$, 
$\gen f = \rho - 1$.
\\[-.25cm]

\noindent
And examples illustrating the various possibilities can be constructed.
\\[-.25cm]
\end{x}

\begin{x}{\small\bf THEOREM} \ 
Suppose that $\rho$ is nonintegral $-$then $f$ is of maximal, minimal, or intermediate type 
according to whether 
$\Delta = \infty$, 
$\Delta = 0$, 
or 
$0 < \Delta < \infty$ 
and conversely.  
\\[-.5cm]

PROOF \ 
In this situation, $\rho = \kappa$ 
(the order of $P$ (cf. 5.10)), while $\rho > q$ ($q$ the order of $e^Q$).  
Therefore the type of $f$ equals the type of $P$ (cf. 3.14), so we can quote 5.19.
\\[-.5cm]
\end{x}

\begin{x}{\small\bf THEOREM} \ 
Suppose that $\rho$ is integral.  
Assume: $\fg < \rho$ $-$then $f$ is either of minimal type or of intermediate type.  
\\[-.5cm]

PROOF \ 
The assumption that $\fg$ is less than $\rho$ puts us in cases (i), (iv), or (v) above.  
Since the series 
$
\textcolor{white}{\bigg|_{\Big|}^{\big|}}
\ds
\sum\limits_{n = 1}^\infty \ 
\frac{1}{\abs{z_n}^\rho}
\ 
$
is convergent, 
one can replace $\kappa$ by $\rho$ in 5.16 and conclude that $P(z, \fg)$ is of minimal type.  
\\[-.75cm]


\qquad \textbullet \quad
In case (i), the order of $e^Q$ is strictly greater than the order of $P$ : $q > \kappa$.  
Therefore
\[
\tau(f) 
\ = \ 
\tau(e^Q) 
\ = \ 
\abs{a_q}
\ \neq \ 
0
\qquad (\tcf. \ 3.14),
\]
so $f$ is of intermediate type.  
\\[-.25cm]

\qquad \textbullet \quad
In case (iv), the order of $e^Q$ and the order of $P$ are one and the same: $q = \kappa$.  
Since 
$0 < \tau(e^Q) = \abs{a_q} < \infty$, \ 
$0 = \tau(P)$, the conclusion is that $\tau (f) = \abs{a_q}$ 
(cf. 3.14), thus $f$ is of intermediate type.
\\[-.25cm]

\qquad \textbullet \quad
In case (v), the order of $e^Q$ is strictly smaller than the order of $P$ : $q > \kappa$.  
Therefore
\[
\tau(f) 
\ = \ 
\tau(P) 
\ = \ 
0
\qquad (\tcf. \ 3.14),
\]
i.e., $f$ is of minimal type.  
\\[-.25cm]

Assuming still that $\rho$ is integral, it remains to deal with cases (ii) and (iii) 
$(\implies \fg = \rho)$.  
Agreeing to write

\[
\begin{cases}
\ 
a_\rho = a_q 
\hspace{0.5cm} \text{if} \ q = \rho
\\[4pt]
\ 
a_\rho = 0
\hspace{0.65cm} \text{if} \ q < \rho
\end{cases}
,
\]
let
\[
\delta (r)
\ = \ 
\Big|
a_\rho + \frac{1}{\rho} \ 
\sum\limits_{\abs{z_n} < r}  \ 
z_n^{-\rho}
\Big|,
\]
put
\[
\delta
\ = \ 
\underset{r \ra \infty}{\limsupx} \ 
\delta (r),
\]
and set
\[
\Gamma 
\ = \ 
\max (\delta, \Delta).
\]
\\[-1.cm]
\end{x}

\begin{x}{\small\bf THEOREM} \ 
Suppose that $\rho$ is integral.  
Assume: \ 
$\fg = \rho$ $-$then $f$ is of maximal, minimal, or intermediate type 
according to whether 
$\Gamma = \infty$, 
$\Gamma = 0$, 
or 
$0 < \Gamma < \infty$ 
and conversely.  
\\[-.5cm]

PROOF \ 
The case (iii) scenario is straightforward: 
$q < \kappa = \rho$, 
hence 
$\tau(f) = \tau (P)$, 
the latter being controlled by 5.20 ($a_\rho = 0$, so the $\Gamma$ there is the $\Gamma$ here).  
As for what happens in case (ii), simply repeat the proof of 5.20 subject to the complication resulting from the 
presence of $a_q \neq 0$ in the definition of $\delta$, the trick being to write

\allowdisplaybreaks\begin{align*}
f(z)
&= \ 
\exp
\Big(
\Big(
a_\rho + \frac{1}{\rho} \ 
\sum\limits_{\abs{z_n} < r} \ 
z_n^{-\rho}
\Big)
z^\rho
\Big)
\hsx
\exp
\big(
Q(z) - a_\rho z^\rho
\big)
\\[15pt]
&
\hspace{3.5cm}
\times \ 
\prod\limits_{\abs{z_n} < r} \ 
E
\Big(
\frac{z}{z_n}, \hsx \fg - 1
\Big)
\ 
\prod\limits_{\abs{z_n} \geq r} \ 
E
\Big(
\frac{z}{z_n}, \hsx \fg
\Big).
\end{align*}
\\[-1.25cm]
\end{x}

\begin{x}{\small\bf REMARK} \ 
Under the preceding assumptions, if $f$ is of minimal type, then

\[
\frac{1}{\rho}\ 
\sum\limits_{n = 1}^\infty \ 
\frac{1}{z_n^\rho} 
\ = \ 
- a_\rho.
\]
\\[-.25cm]
\end{x}


\chapter{
$\boldsymbol{\S}$\textbf{7}.\quad  REPRESENTATION THEORY}
\setlength\parindent{2em}
\setcounter{theoremn}{0}
\renewcommand{\thepage}{\S7-\arabic{page}}

\qquad
Let $f$ be an entire function $-$then as regards its zeros, 
there are three possibilities.
\\[-.25cm]

\qquad 1. \quad 
$f$ has no zeros.
\\[-.5cm]

\qquad 2. \quad 
$f$ has a finite number of zeros.
\\[-.5cm]

\qquad 3. \quad 
$f$ has an infinite number of zeros.
\\[-.25cm]

\begin{x}{\small\bf THEOREM} \ 
If $f$ has no zeros, then there is an entire function $g$ such that $f = e^g$.
\\[-.5cm]

PROOF \ 
Since $f$ has no zeros, 
$\ds
\frac{1}{f}
$
is entire, as is
$\ds
\frac{f^\prime}{f}
$.  
Define $g$ by the prescription
\[
g(z) 
\ = \ 
\int\limits_0^z \ 
\frac{f^\prime (t)}{f(t)}
\ \td t,
\]
the path of integration being immaterial $-$then 
$
\ds
g^\prime = \frac{f^\prime}{f}$.  
And 
\allowdisplaybreaks\begin{align*}
\big(
f \hsy e^{-g}
\big)^\prime \ 
&=\ 
f^\prime \hsy e^{-g} - f \hsy g^\prime e^{-g}
\\[11pt]
&=\ 
e^{-g}
\Big(
f^\prime - f \hsy \frac{f^\prime}{f}
\Big)
\\[11pt]
&=\ 
0.
\end{align*}
Therefore
\[
f(z) \hsx e^{-g(z)}
\ = \ 
f(0) \hsx e^{-g(0)}
\ = \ 
f(0)
\]
\qquad\qquad 
$\implies$
\[
f(z) 
\ = \ 
f(0) \hsx e^{g(z)}.
\]
Conclude by absorbing $f(0)$ into the exponential.
\\[-.25cm]
\end{x}

\begin{x}{\small\bf REMARK} \ 
If $f$ has no zeros, if $f = e^g$, and if $f$ is of finite order, then $g$ is a polynomial (cf. 2.42).
\\[-.25cm]
\end{x}

Suppose now that $f$ is an entire function with finitely many zeros \hsx 
$z_1 \neq 0, \hsx \ldots,$\hsx$z_n \neq 0$  
(each counted with multiplicity), 
as well as a zero of order $m \geq 0$ at the origin $-$then the entire function 
\[
\frac{f(z)}{\ds z^m \hsy  \prod\limits_{k = 1}^n \ \bigg(1 - \frac{z}{z_k}\bigg)}
\]
has no zeros, hence equals
\[
e^{g(z)}, 
\]
where $g(z)$ is entire, so

\[
f(z) 
\ = \ 
z^m \hsy 
e^{g(z)} \hsy
\prod\limits_{k = 1}^n \ 
\bigg(
1 - \frac{z}{z_k}
\bigg) \hsx .
\]

\qquad
{\small\bf \un{N.B.}} \ 
If $f$ is of finite order, then $g$ is a polynomial (cf. 7.2).  
\\[-.25cm]

Assume henceforth that $f$ is a transcendental entire function of finite order $\rho$ with an infinite 
number of nonzero zeros 
$\{z_n : n \geq 1\}$ 
and a zero of order $m \geq 0$ at the origin.  
Set $\Pi (z) = P(z, \fg)$. 
\\[-.25cm]

\begin{x}{\small\bf HADAMARD FACTORIZATION} \ 
We have
\[
f(z)
\ = \ 
z^m \hsy e^{Q(z)} \hsy
\Pi (z),
\]
where $Q(z)$ is a polynomial of degree $q \leq \rho$.
\\[-.5cm]

PROOF \ 
The quotient
\[
\frac{f(z)}{z^m \hsy \Pi (z)}
\]
is entire and has no zeros, 
thus can be written as $e^{Q(z)}$, where $Q(z)$ is entire.  
Owing to 2.37, the order of 

\[
\frac{f(z)}{z^m \hsy \Pi (z)}
\]
is $\leq$ the maximum of $\rho$ and the order of $z^m \hsy \Pi (z)$, 
the order of the latter being that of $\Pi(z)$ (cf. 2.36), 
which in turn is equal to $\kappa$ (cf. 5.10).
But $\kappa$ is $\leq \rho$ (cf. 4.20).  
Therefore the order of $e^{Q(z)}$ is $\leq \rho$, so $Q(z)$ is a polynomial of degree $q \leq \rho$ (cf. 2.42).
\\[-.25cm]
\end{x}

\begin{x}{\small\bf REMARK} \ 
If $f$ is a transcendental entire function of finite nonintegral order $\rho$, 
then it is automatic that $f$ has an infinity of zeros.
\\[-.5cm]

[In fact, 
\[
\rho 
\ = \ 
\max (q, \kappa) 
\quad (\tcf. \ 5.11) 
\quad
\implies
\quad
\rho = \kappa.
\]
But if $f$ had finitely many zeros, then of necessity, $\kappa = 0, \ldots$ .]
\\[-.25cm]
\end{x}

By definition (cf. 6.1), 
\[
\gen f 
\ = \ 
\max (q, \fg) 
\]
and the simplest cases
\[
\gen f \ =\ 
\begin{cases}
\ 
0
\\[4pt]
\ 
1
\end{cases}
\]
are of special interest.
\\[-.25cm]

\begin{x}{\small\bf LEMMA} \ 
If $\gen f = 0$ or 1, then $\rho \leq 2$. 
\\[-.5cm]

PROOF \ 
if $\rho$ is not an integer, 
then $\gen f = [\rho]$ (cf. 6.3), 
hence $\rho < 2$.  
On the other hand, 
if $\rho$ is an integer, 
then $\gen f = \rho$ or $\rho - 1$ (cf. 6.4), 
hence $\rho \leq 2$. 
\\[-.25cm]

\qquad \textbullet \quad 
$\gen f = 0$.  
Here $q = 0$, so $Q(z) = C$, and 
\[
f(z) \ = \ 
z^m \hsy e^C \ 
\prod\limits_{n = 1}^\infty \ 
\bigg(
1 - \frac{z}{z_n}
\bigg),
\]
where
\[
\sum\limits_{n = 1}^\infty \ 
\frac{1}{\abs{z_n}} 
\ < \ 
\infty.
\]
\\[-.25cm]

\qquad \textbullet \quad 
$\gen f = 1$.  
\[
\begin{cases}
\ 
q = 1
\\[4pt]
\ 
\fg = 1
\end{cases}
\implies 
f(z) \ = \ 
z^m \hsy e^{a z + b} \ 
\prod\limits_{n = 1}^\infty \ 
\Big(
1 - \frac{z}{z_n}
\Big)
\hsy e^{z / z_n}, 
\]
where $a \neq 0$ and
\[
\sum\limits_{n = 1}^\infty \ 
\frac{1}{\abs{z_n}^2} 
\ < \ 
\infty
\]
but 
\[
\sum\limits_{n = 1}^\infty \ 
\frac{1}{\abs{z_n}} 
\ = \ 
\infty.
\]

\[
\begin{cases}
\ 
q = 0
\\[4pt]
\ 
\fg = 1
\end{cases}
\implies 
f(z) \ = \ 
z^m \hsy e^C \ 
\prod\limits_{n = 1}^\infty \ 
\Big(
1 - \frac{z}{z_n}
\Big)
\hsy e^{z / z_n}, 
\]
\\[-.75cm]

\noindent
where 
\[
\sum\limits_{n = 1}^\infty \ 
\frac{1}{\abs{z_n}^2} 
\ < \ 
\infty
\]
but 
\[
\sum\limits_{n = 1}^\infty \ 
\frac{1}{\abs{z_n}} 
\ = \ 
\infty.
\]

\[
\begin{cases}
\ 
q = 1
\\[4pt]
\ 
\fg = 0
\end{cases}
\implies 
f(z) \ = \ 
z^m \hsy e^{a z + b} \ 
\prod\limits_{n = 1}^\infty \ 
\Big(
1 - \frac{z}{z_n}
\Big), 
\]
\\[-.75cm]

\noindent
where $a \neq 0$ and 
\[
\sum\limits_{n = 1}^\infty \ 
\frac{1}{\abs{z_n}} 
\ < \ 
\infty.
\]
\end{x}


\chapter{
$\boldsymbol{\S}$\textbf{8}.\quad  ZEROS}
\setlength\parindent{2em}
\setcounter{theoremn}{0}
\renewcommand{\thepage}{\S8-\arabic{page}}

\qquad 
Let $f$ be an entire function.
\\[-.25cm]

\begin{x}{\small\bf DEFINITION} \ 
A 
\un{critical point} 
of $f$ is a zero of $f^\prime$.
\\[-.25cm]
\end{x}

Suppose that
\[
f (z) 
\ = \ 
\prod\limits_{i = 1}^k \ 
(z - z_i)^{m_i}
\]
is a polynomial of degree $n$, 
thus 
$
\ 
\ds
\sum\limits_{i = 1}^k \ 
m_i = n
\ 
$
and the $z_i$ are distinct.  
There are then two kinds of critical points. 
\\[-.25cm]

\qquad \textbullet \quad
A zero $z_i$ of multiplicity $m_i > 1$ is said to be of the \un{first kind}.
Counting it $m_i - 1$ times (its multiplicity as a zero of $f^\prime$), 
it follows that there are $n - k$ critical points of the first kind.
\\[-.25cm]

\qquad \textbullet \quad
Since the degree of 
$f^\prime$ 
is $n - 1$, there are $k - 1$ additional critical points, 
these being termed of the \un{second kind}.  
They are not zeros of $f$ but are zeros of 
$
\ds
\frac{f^\prime}{f}
$ 
(defined on $\Cx - \{z_1, \ldots, z_k\}$), i.e., are zeros of 

\[
\sum\limits_{i = 1}^k \ 
\frac{m_i}{z - z_i}.
\]
\\[-1.25cm]

\begin{x}{\small\bf REMARK} \ 
There is no simple relation between the number of distinct zeros of a polynomial and its derivative.
\\[-.25cm]

\qquad (1) \
The polynomial 
$
\ds
\prod\limits_{i = 1}^k \ 
(z - i)^2
$
has $k$ distinct zeros while its derivative has $2 k - 1$ distinct zeros.
\\[-.25cm]

\qquad (2) \
The polynomial $z^n - 1$ has $n$ distinct zeros but its derivative has just one. 
\\[-.25cm]

\qquad (3) \
The polynomial 
$z^{n - 1} (z - 1)$ has two distinct zeros as does its derivative.
\\[-.25cm]
\end{x}

\begin{x}{\small\bf THEOREM} \ 
The zeros of 
$f^\prime$ 
belong to the convex hull of the zeros of $f$.  
\\[-.5cm]

PROOF \ 
It suffices to consider a zero $z_0$ of the second kind:
\[
\sum\limits_{i = 1}^k \ 
\frac{m_i}{z_0 - z_i}
\ = \ 0
\ \implies \ 
\sum\limits_{i = 1}^k \ 
\frac{m_i}{\bar{z}_0 - \bar{z}_i}
\ = \ 0
\]
\qquad 
$\implies$
\[
\sum\limits_{i = 1}^k \ 
m_i \ 
\frac{z_0 - z_i}{\abs{z_0 - z_i}^2}
\ = \ 
0
\]
\qquad 
$\implies$
\[
z_0 \ 
\sum\limits_{i = 1}^k \ 
\frac{m_i}{\abs{z_0 - z_i}^2}
\ = \ 
\sum\limits_{i = 1}^k \ 
m_i \hsx
\frac{z_i}{\abs{z_0 - z_i}^2}
\]
\qquad 
$\implies$
\[
z_0
\ = \ 
\sum\limits_{i = 1}^k \ 
\lambda_i \hsy z_i,
\]
where

\[
\lambda_i 
\ = \ 
\frac
{
\ds
\frac{m_i}{\abs{z_0 - z_i}^2}}
{
\ds
\sum\limits_{j = 1}^k \ 
\frac{m_j}{\abs{z_0 - z_j}^2}
}
\ > \ 
0
\]
and
\[
\sum\limits_{i = 1}^k \ 
\lambda_i 
\ = \ 1.
\]
\\[-1.cm]
\end{x}


\begin{x}{\small\bf EXAMPLE} \ 
There are transcendental entire functions for which this result is false.
\\[-.5cm]

[Take
\[
f (z) 
\ = \ 
z \hsy \exp \frac{z^2}{2}.
\]
It has one zero, viz. $z = 0$, but its derivative
\[
f^\prime (z) 
\ = \ 
(1 + z^2) \hsy \exp \frac{z^2}{2}
\]
has two zeros, viz $\pm \sqrt{-1}$.]
\\[-.25cm]
\end{x}

\begin{x}{\small\bf NOTATION} \ 
Given a nonempty closed subset $T$ of $\Cx$, 
let
$\langle T \rangle$ 
stand for its closed convex hull.
\\[-.25cm]
\end{x}

\begin{x}{\small\bf LEMMA} \ 
Let $f$ be a transcendental entire function of finite order $\rho$ with $\gen f = 0$.  
Assume: \ 
The zeros of $f$ lie in $T$ $-$then the zeros of 
$f^\prime$ 
lie in 
$\langle T \rangle$.
 \\[-.5cm]

PROOF \ 
Decompose $f$ per 7.3:
\[
f (z) 
\ = \ 
C \hsy z^m \ 
\prod\limits_{n = 1}^\infty \ 
\Big(
1 - \frac{z}{z_n}
\Big),
\]
and put
\[
f_N (z) 
\ = \ 
C \hsy z^m \ 
\prod\limits_{n = 1}^N \ 
\Big(
1 - \frac{z}{z_n}
\Big).
\]
Then
\[
f_N \ra f
\qquad (N \ra \infty)
\]
uniformly on compact subsets of $\Cx$, so
\[
f_N^\prime \ra f^\prime
\qquad (N \ra \infty)
\]
uniformly on compact subsets of $\Cx$.  
But  the zeros of 
$f^\prime$ 
are limits of zeros of the 
$f_N^\prime$, 
these in turn being elements of 
$\langle T \rangle$
(cf. 8.3).
\\[-.5cm]
 
[Note: \ 
In terms of $\rho$, 
\[
0 
\ \leq \ 
\rho 
\ < \ 
1 
\ \implies \ 
\gen f \ = \ [\rho] \ = \ 0
\qquad (\tcf. \ 6.3)
\]
or
\[
\rho \ = \ 1
\ \text{and} \ 
\gen f = \rho - 1 \ = \ 1 - 1 \ = \ 0
\qquad (\tcf. \ 6.4).]
\]
\\[-1.25cm]
\end{x}

\begin{x}{\small\bf EXAMPLE} \ 
The transcendental entire function
\[
f (z)
\ = \ 
\prod\limits_{k = 0}^K \ 
\cos (z - k \sqrt{-1})^{1/2}
\]
is of order $1/2$ and its zeros lie in the set
\[
T : \Reg z \geq 0 
\quad \& \quad 
0 \leq \Img z \leq K.
\]
Since here 
$T = \langle T \rangle$, 
the zeros of its derivative also lie in $T$.
\\[-.25cm]
\end{x}

\begin{x}{\small\bf REMARK} \ 
Take $\rho = 1$ and suppose that the conditions of 6.8 are in force with $f$ of minimal type, 
hence $\Gamma = 0$ and 
\allowdisplaybreaks
\allowdisplaybreaks\begin{align*}
\sum\limits_{n = 1}^\infty \ 
\frac{1}{t_n} \ 
&=\ 
- a_1
\qquad (\tcf. \ 6.9)
\\[8pt]
&\equiv \ 
-a.
\end{align*}
Then 8.6 still goes through.  
Thus write
\[
f (z) 
\ = \ 
C \hsy z^m \hsy e^{a \hsy z} \ 
\prod\limits_{n = 1}^\infty \ 
\Big(
1 - \frac{z}{z_n}
\Big) \hsy e^{-z/z_n}
\qquad (\tcf. \ 7.3)
\]
and let
\[
f_N (z) 
\ = \ 
C \hsy z^m \hsy e^{a \hsy z}\ 
\prod\limits_{n = 1}^N \ 
\Big(
1 - \frac{z}{z_n}
\Big) \hsy e^{-z/z_n}.
\]
Since
\[
\sum\limits_{n = 1}^N \ 
\frac{1}{z_n} - a
\ra 0
\qquad (N \ra \infty),
\]
it follows that
\[
f_N \ra f
\qquad (N \ra \infty)
\]
uniformly on compact subsets of $\Cx$.  
\\[-.25cm]
\end{x}

\begin{x}{\small\bf EXAMPLE} \ 
Fix $\tau > 0$ $-$then 
\[
f(z) 
\ = \ 
(z^2 - 1)^m \hsy e^{\tau \hsy z}
\]
is a transcendental entire function of order 1 and type $\tau$ and its zeros lie in the convex set 
$[-1, 1]$.  
On the other hand, $f$ has a critical point at
\[
-\frac{1}{\tau} \hsx 
\big(
m + \sqrt{m^2 + \tau^2} 
\hsx 
\big)
\ \notin \ [-1, 1].
\]
Therefore the assumption of minimal type cannot be dropped in 8.8.
\\[-.25cm]
\end{x}

Before proceeding further, it will be best to recall some standard generalities.
\\[-.25cm]

\begin{x}{\small\bf LEMMA} \ 
Suppose that $f$ is a real analytic function $-$then in any finite interval $I$, 
$f$ has at most a finite number of distinct zeros.
\\[-1.25cm]
\end{x}

\begin{spacing}{1.75}
[Note: \ 
This is false if $f$ is merely $C^\infty$: \ 
Take $I = [0,1]$ and consider 
$
\ds
f (x) = 
x \sin \Big(\frac{1}{x}\Big)
$.]
\\[-1cm]
\end{spacing}

\begin{x}{\small\bf ROLLE'S THEOREM} \ 
Suppose that $f$ is a real analytic function $-$then between any two consecutive zeros of $f$, 
say $f(a) = 0$, $f(b) = 0$, $(a < b)$, 
$f^\prime$ 
has an odd number of zeros in $]a, b[$ counted according to multiplicity.
\\[-.25cm]
\end{x}


\begin{x}{\small\bf LEMMA} \ 
Suppose that $f$ is a real analytic function and let $I$ be a finite interval.  
Assume: \ 
$f^\prime$ 
has $Z^\prime$ 
zeros in $I$ counted according to multiplicity $-$then 
$f$ has at most $Z^\prime + 1$ zeros in $I$ counted according to multiplicity.
\\[-.5cm]

PROOF  \ 
Let $d$ be the number of distinct zeros of $f$ in $I$ and let $D$ denote the number of zeros of $f$ in $I$ 
counted according to multiplicity.  
At a zero of $f$ of multiplicity $m_k$, 
$f^\prime$ 
has a zero of multiplicity $m_k - 1$.  
In addition, by Rolle's theorem, 
$f^\prime$ 
has at least one zero between two consecutive zeros of $f$.  
Therefore 
\allowdisplaybreaks
\allowdisplaybreaks\begin{align*}
Z^\prime \ 
&\geq \
\sum\limits_{k = 1}^d \ 
(m_k - 1) + d  - 1
\\[15pt]
&=\ 
D - d + d - 1 
\\[15pt]
&=\ 
D - 1
\end{align*}
\qquad\qquad
$\implies$
\[
D 
\ \leq \ 
Z^\prime + 1.
\]

[Note: \ 
It is thus a corollary that if $f$ has $Z$ zeros in $I$ counted according to multiplicity, 
then 
$f^\prime$ 
has at least
$Z - 1$ zeros in $I$ counted according to multiplicity.]
\end{x}

\begin{x}{\small\bf DEFINITION} \ 
An entire function is said to be \un{real} if it assumes real values on the real axis.
\\[-.5cm]

[Note: \ 
The restriction of a real entire function to the real axis is a real analytic function.]
\\[-.5cm]
\end{x}

\qquad
{\small\bf \un{N.B.}} \ 
If 
\[
f (z) 
\ = \ 
\sum\limits_{n = 0}^\infty \ 
c_n \hsy z^n, 
\]
then $f$ is real iff $\forall \ n$, $c_n$ is real.  
\\[-.25cm]


\begin{x}{\small\bf EXAMPLE} \ 
If $f$ is a polynomial and if the zeros of $f$ are real, then $f$ is real 
(to within a multiplicative constant) but not conversely.
\\[-.5cm]
\end{x}

\begin{x}{\small\bf REMARK} \ 
If $f$ is a transcendental entire function of finite order and if $\gen f = 0$, 
then the reality of its zeros forces the reality of $f$ (up to a constant factor) 
but this need not be true if $\gen f > 0$ (although it will be if $f$ is a canonical product with real zeros).
\\[-.5cm]
\end{x}

\begin{x}{\small\bf THEOREM} \ 
If $f$ is a polynomial and if the zeros of $f$ are real, then the zeros of $f^\prime$ are real.
\\[-.5cm]

[In view of 8.3, this is immediate.]
\\[-.5cm]

[Note: \ 
Suppose that 
$z_1 < \cdots < z_k$ are the distinct zeros of $f$ $-$then by Rolle's theorem, 
$f$ has at least one critical point in each of the intervals 
\hsy
$]z_i, z_{i + 1}[$ 
\hsy
$(i = 1, \ldots, k - 1)$ 
and these critical points are of the second kind.  
Since there are $k - 1$ critical points of the second kind, there is but one critical point in 
\hsy
$]z_i, z_{i + 1}[$ 
\hsy
and it is simple.  
Finally, all critical points of $f$ are to be found in 
$[z_1, z_k]$.]
\\[-.25cm]
\end{x}

\begin{x}{\small\bf EXAMPLE} \ 
The zeros of the following polynomials are real and simple.
\\[-.25cm]

\qquad \textbullet \quad
The 
\un{Legendre polynomials}:
\[
P_n (x) 
\ = \ 
\frac{1}{2^n \hsy n !} \hsx
\frac{\td^n}{\td x^n} \hsx
(x^2 - 1)^n.
\]

\qquad \textbullet \quad
The 
\un{Laguerre polynomials}:
\[
L_n (x) 
\ = \ 
\frac{e^x}{n !} \hsx 
\frac{\td^n}{\td x^n} \hsx
e^{-x} \hsy x^n.
\]


\qquad \textbullet \quad
The \un{Hermite polynomials}:
\[
H_n (x) 
\ = \ 
(-1)^n \hsx
e^{x^2} \hsx 
\frac{\td^n}{\td x^n} \hsx
e^{-x^2}.
\]

A polynomial 
\[
f (z) 
\ = \ 
\prod\limits_{n = 1}^N \ 
E
\Big(
\frac{z}{z_n}, 0
\Big)
\ = \ 
\prod\limits_{n = 1}^N \ 
\Big(
1 - \frac{z}{z_n}
\Big)
\]
of degree $N$ is, in particular, a canonical product, so 8.16 is a special case of the next result (compare to 8.6).
\\[-.25cm]
\end{x}

\begin{x}{\small\bf THEOREM} \ 
Let
\[
f (z) 
\ = \ 
\prod\limits_{n = 1}^\infty \ 
E
\Big(
\frac{z}{z_n}, \fg
\Big)
\]
be a canonical product whose zeros are real $-$then the zeros of 
$f^\prime$ 
are real.
\\[-.5cm]

PROOF \ 
Working with the zeros of 
$f^\prime$ 
that are not zeros of $f$, pass to
\[
\frac{f^\prime (z)}{f (z)} 
\ = \ 
z^\fg \ 
\sum\limits_{n = 1}^\infty \ 
\frac{1}{z_n^\fg \hsy (z - z_n)},
\]
which shows that the origin is a zero of multiplicity $\fg$ of 
$f^\prime (z)$. 
Let
\[
F(z) 
\ = \ 
z^{-\fg} \hsx \frac{f^\prime (z)}{f (z)} 
\]
and write 
$z_n = x_n + \sqrt{-1} \ 0$, hence
\[
F(z) 
\ = \ 
\sum\limits_{n = 1}^\infty \ 
\frac{1}{x_n^\fg \hsy (z - x_n)}.
\]
Suppose now that 
\[
f^\prime (c) 
\ = \ 
f^\prime (a + \sqrt{-1} \hsx b) 
\ = \ 0,
\]
the claim being that $b = 0$.  
To see this, separate the real and imaginary parts in $F (c) = 0$ to get
\[
a \ 
\sum\limits_{n = 1}^\infty \  
\frac{1}{x_n^\fg \hsy \abs{c - x_n}^2}
\ - \ 
\sum\limits_{n = 1}^\infty \ 
\frac{1}{x_n^{\fg - 1} \hsy \abs{c - x_n}^2}
\ = \ 
0
\]
and 
\[
b \ 
\sum\limits_{n = 1}^\infty \ 
\frac{1}{x_n^\fg \hsy \abs{c - x_n}^2} 
\ = \ 
0.
\]
\\[-.75cm]

\qquad \textbullet \quad
If $\fg$ is even or if $\forall \ n$, $x_n > 0$ $(x_n < 0)$, then $b = 0$. 
\\[-.25cm]

\qquad \textbullet \quad
If $\fg$ is odd and there are positive as well as negative $x_n$, then 
\allowdisplaybreaks
\allowdisplaybreaks\begin{align*}
b \neq 0 
&\implies \ 
\sum\limits_{n = 1}^\infty \ 
\frac{1}{x_n^\fg \hsy \abs{c - x_n}^2} 
\ = \ 0
\\[15pt]
&\implies \ 
\sum\limits_{n = 1}^\infty \ 
\frac{1}{x_n^{\fg - 1} \hsy \abs{c - x_n}^2}
\ = \ 0.
\end{align*}
But this is impossible since $\fg - 1$ is even.
\\[-.25cm]
\end{x}

\begin{x}{\small\bf ADDENDUM} \ 
Let
$\zeta^\prime < \zeta^{\prime\prime}$ 
be consecutive zeros of $f$ of the same sign $-$then 
there is exactly one distinct zero of 
$f^\prime$ in 
$]\zeta^\prime, \zeta^{\prime\prime}[$.
\\[-.5cm]

[By Rolle's theorem, there is at least one $\zeta$ in 
$]\zeta^\prime, \zeta^{\prime\prime}[$ 
such that 
$f^\prime (\zeta) = 0$ 
(bear in mind that $f$ is real).  
As for its uniqueness, if $\fg$ is even or if $\forall \ n$, 
$x_n > 0$ $(x _n < 0)$, then the sign of 
\[
F^\prime (x) 
\ = \ 
- 
\sum\limits_{n = 1}^\infty \ 
\frac{1}{x_n^\fg \hsy (x - x_n)^2}
\]
is constant, thus $F(x)$ is monotonic between 
$\zeta^\prime$
and 
$\zeta^{\prime\prime}$, 
thus cannot vanish more than 
once in 
$]\zeta^\prime, \zeta^{\prime\prime}[\hsy$.  
So, if $\alpha \neq \beta$ were distinct zeros of 
$f^\prime$ 
in 
$]\zeta^\prime, \zeta^{\prime\prime}[\hsy$, 
then $\fg$ would 
have to be odd and there would have to be both positive and negative $x_n$.  
But
\[
\begin{cases}
\ 
0 \ = \  F(\alpha) + F(\beta) \ = \  (\alpha + \beta) X - 2Y
\\[8pt]
\ 
0 \ = \  F(\alpha) - F(\beta) \ = \  (\beta - \alpha) X
\end{cases}
\implies 
X = 0 
\qquad (\alpha \neq \beta)
\]
\\[-1cm]

$\implies$
\[
- 2 \ 
\sum\limits_{n = 1}^\infty \ 
\frac{1}{x_n^{\fg - 1} \hsy (\alpha - x_n) \hsy (\beta - x_n)} 
\ = \ 0.
\]

\noindent
This, however, is impossible: \ 
$\fg - 1$ is even and $\forall \ n$, $(\alpha - x_n) \hsy (\beta - x_n) > 0$.]
\\[-.25cm]
\end{x}

\begin{x}{\small\bf REMARK} \ 
It can be shown that the genus of 
$f^\prime$ 
is equal to the genus of $f$.  
\\[-.5cm]

[This is obvious if the order $\rho$ of $f$ is not an integer (for $\rho = \rho^\prime$ (the order of 
$f^\prime$) 
(cf. 2.25) and 
$\gen f = [\rho] = [\rho^\prime] = \gen f^\prime$  
(cf. 6.3)) but not so obvious otherwise.]
\\[-.25cm]
\end{x}

\begin{x}{\small\bf EXAMPLE} \ 
Let 
\[
f_\alpha (z) 
\ = \ 
\prod\limits_{n = 2}^\infty \ 
\Big(
1 + \frac{1}{n (\log n)^\alpha}
\Big)
\qquad (1 < \alpha < 2).
\]
Then $\rho (f_\alpha) = 1$, 
$\gen f_\alpha = 0$, and 
$\gen f_\alpha^\prime = 0$. 
On the other hand, 
\[
A \neq 0 \implies \gen (f_\alpha - A) = 1
\]
\qquad \qquad
$\implies$
\[
\gen (f_\alpha - A)^\prime 
\ = \ 
\gen f_\alpha^\prime = 0.
\]
\\[-1cm]
\end{x}

If $f$ is a nonconstant real entire function, 
then the zeros of $f$ are either real or, if nonreal, occur in conjugate pairs 
$(z_0, \bar{z}_0)$.
\\[-.25cm]

\qquad
{\small\bf \un{N.B.}} \ 
The multiplicity of $z_0$ is the same as the multiplicity of $\bar{z}_0$.
\\[-.25cm]

\begin{x}{\small\bf LEMMA} \ 
If $f$ is a nonconstant real polynomial, 
then the number of nonreal zeros of 
$f^\prime$ 
counted according to multiplicity is $\leq$ the number of nonreal zeros of $f$ 
counted according to multiplicity.
\\[-.5cm]

PROOF\ 
Suppose that the degree of $f$ is $n$, 
the number of real zeros of $f$ counted according to multiplicity is $r$, 
and the number of nonreal zeros of $f$ counted according to multiplicity is $n - r$, 
then for 
$f^\prime$ 
they are $= n - 1, \ \geq r - 1$, (cf. 8.12), and 
$\leq n - 1 - (r - 1) = n - r$.
\\[-.25cm]
\end{x}

Let $f$ be a nonconstant real entire function of finite order $\rho$ and suppose that $f$ has 
$0 \leq C = 2 D < \infty$ nonreal zeros counted according to multiplicity $-$then 
$f^\prime$ 
has 
$0 \leq C^\prime = 2 D^\prime \leq  C = 2 D < \infty$ 
nonreal zeros counted according to multiplicity (see 8.24 below).
\\[-.25cm]

\un{Extra Zeros} \quad
This refers to 
$f^\prime$ 
and there are two kinds.
\\[-.25cm]

\qquad \textbullet \quad
If $\zeta^\prime < \zeta^{\prime\prime}$ are consecutive real zeros of $f$, 
then by Rolle's theorem, 
$f^\prime$ 
has an odd number of zeros in 
$]\zeta^\prime, \zeta^{\prime\prime}[$ 
counted according to multiplicity, say $2 k + 1$.  
One then says that 
$f^\prime$ 
has $2 k$ extra zeros between $\zeta^\prime $and $\zeta^{\prime\prime}$.
\\[-.25cm]

\qquad \textbullet \quad
If $f$ has a largest real zero $x_L$ or a smallest real zero $x_S$, 
then any zero of 
$f^\prime$  
in 
$]x_L, +\infty[$ 
or 
$]-\infty, x_S[$ 
is called extra and will be counted according to multiplicity.
\\[-.25cm]

Let $E^\prime$ denote the total number of extra zeros of 
$f^\prime$ .
\\[-.25cm]


\begin{x}{\small\bf EXAMPLE} \ 
Take for $f$ a canonical product whose zeros are real (cf. 8.18) $-$then
it might be that 0 is extra as in 
\[
\begin{tikzpicture}[scale=0.5]
      \draw[-] (-3,0) -- (3,0) node[below] {$$};
      \draw[] (0,-.4) -- (0,-.4) node[below] {$x_L$};
      \draw[] (0,0) -- (0,0) node[] {$|$};
      \draw[] (3,-.4)-- (3,-.4) node[below] {$0$};
      \draw[] (3,0) -- (3,0) node[] {$|$}; 
\end{tikzpicture}
\qquad
\begin{tikzpicture}[scale=0.5]
\draw[] (0,0) -- (0,0) node[above] {$\text{or}$};
\draw[] (0,-.4)-- (0,-.4) node[below] {$$};
\end{tikzpicture}
\qquad
\begin{tikzpicture}[scale=0.5]
      \draw[-] (-3,0) -- (3,0) node[below] {$$};
      \draw[] (0,-.4) -- (0,-.4) node[below] {$x_S$};
      \draw[] (0,0) -- (0,0) node[] {$|$};
      \draw[] (-3,-.4)-- (-3,-.4) node[below] {$0$};
      \draw[] (-3,0) -- (-3,0) node[] {$|$};
      
      \draw[] (3.75,0) -- (3.75,0) node[below] {$.$};
\end{tikzpicture}
\]
\\[-1.25cm]
\end{x}

\begin{x}{\small\bf THEOREM}\footnote[2]{\vspace{.11 cm}
E. Borel,  
\textit{Lecons sur les Fonctions Enti\`eres}, Gauthier-Villars, 1900, pp. 37-47.} \ 
Under the preceding assumptions on $f$, 
\[
E^\prime + C^\prime 
\ \leq \ 
C + \gen f,
\]
and 
\[
\gen f
\ = \ 
\gen f^\prime.
\]
\\[-1.25cm]
\end{x}

\begin{x}{\small\bf SCHOLIUM} \ 
If $f$ has a canonical product whose zeros are real, then $E^\prime \leq \fg$ (cf. 8.18).
\\[-.75cm]

[Note: \ 
As a special case, if $f$ is a polynomial and if the zeros of $f$ are real, 
then $E^\prime = 0$ 
(the critical points guaranteed by Rolle's theorem are simple (cf. 8.16)).]
\\[-.25cm]
\end{x}

\begin{x}{\small\bf EXAMPLE} \ 
Take
\[
f (z) 
\ = \ 
(z + 1) \hsy \exp \frac{z^2}{2}.
\]
It has one real zero, viz. $z = -1$, and its derivative
\[
f^\prime (z) 
\ = \ 
(1 + z + z^2) \hsy \exp \frac{z^2}{2} 
\]
has two nonreal zeros, viz.

\[
z 
\ = \ 
\frac{-1 \pm \sqrt{-3}}{2} \hsx .
\]
Here
\[
\begin{cases}
\ 
E^\prime = 0
\\[4pt]
\ 
C^\prime = 2
\end{cases}
,\quad
\begin{cases}
\ 
C = 0
\\[4pt]
\ 
\gen f = 2
\end{cases}
.
\]
\\[-.5cm]
\end{x}

\begin{x}{\small\bf EXAMPLE} \ 
Take 
\[
f (z) 
\ = \ 
(z^2 - 4) \hsy \exp \frac{z^2}{3}.
\]
It has two real zeros, viz. $z = \pm2$, and its derivative
\[
f^\prime (z) 
\ = \ 
\frac{2}{3} \hsy z \hsy (z^2 - 1) \hsy \exp \frac{z^2}{3}
\]
has three real  zeros, viz. $z = -1, 0, 1$.  
Here
\[
\begin{cases}
\ 
E^\prime = 2
\\[4pt]
\ 
C^\prime = 0
\end{cases}
,\quad
\begin{cases}
\ 
C = 0
\\[4pt]
\ 
\gen f = 2
\end{cases}
.
\]

[Note: \ 
The three zeros between -2 and 2 are per Rolle and $3 = 2 + 1$, so $E^\prime = 2$.]
\\[-.25cm]
\end{x}

\begin{x}{\small\bf EXAMPLE} \ 
Take
\[
f (z) 
\ = \ 
(z^2 - 1) \hsy e^z.
\]
It has two real zeros, viz. $z = \pm 1$, and its derivative
\[
f^\prime (z) 
\ = \ 
(z^2 + 2 z - 1) \hsy e^z
\]
has two nonreal zeros, viz. $z = -1 \pm \sqrt{2}$.  
Here
\[
\begin{cases}
\ 
E^\prime = 1
\\[4pt]
\ 
C^\prime = 0
\end{cases}
,\quad
\begin{cases}
\ 
C = 0
\\[4pt]
\ 
\gen f = 1
\end{cases}
.
\]

[Note: \ 
The zero 
$-1 + \sqrt{2} \hsx $ lies between -1 and 1 and is per Rolle but the zero 
$-1 - \sqrt{2} \hsx $ lies to the left of -1, hence is extra.]
\\[-.25cm]
\end{x}


\begin{x}{\small\bf REMARK} \ 
If $f$ is a nonconstant real polynomial, then
\[
E^\prime + C^\prime \ = \quad 
\begin{cases}
\ 
C \hspace{1.25cm} \text{if deg $f > C$}
\\[4pt]
\ 
C  - 1\hspace{0.5cm} \text{if deg $f = C$}
\end{cases}
.
\]

[Note: \ 
In particular, $C^\prime \leq C$ (cf. 8.22).]
\\[-.25cm]
\end{x}

\begin{x}{\small\bf THEOREM} \ 
Let $f$ be a nonconstant real entire function of finite order $\rho$.  
Assume: \ 
The zeros of $f$ are real and 
$\gen f = 0$ or 1 $-$then the zeros of 
$f^\prime$ 
are real and 
\[
\gen f 
\ = \ 
\gen f^\prime.
\]

PROOF \ 
In this situation, 
\[
E^\prime + C^\prime 
\ \leq \ 
\gen f 
\qquad (\tcf. \ 8.24), 
\]
so
\[
\gen f = 0 
\hsx \implies \hsx 
C^\prime = 0.
\]
And
\allowdisplaybreaks
\allowdisplaybreaks\begin{align*}
\gen f = 1 
&\implies
E^\prime + C^\prime \leq 1
\\[11pt]
&\implies
C^\prime \leq 1.
\end{align*}
But $C^\prime$ is even.   
Therefore $C^\prime = 0$ (although $E^\prime$ might be 1 (cf. 8.28)).
\\[-.5cm]

[Note: \ 
It follows that 
$f^\prime$ 
satisfies the same general conditions as $f$.]
\\[-.25cm]
\end{x}


\chapter{
$\boldsymbol{\S}$\textbf{9}.\quad  JENSEN CIRCLES}
\setlength\parindent{2em}
\setcounter{theoremn}{0}
\renewcommand{\thepage}{\S9-\arabic{page}}

\qquad 
We begin with a computation.
\\[-.25cm]

\begin{x}{\small\bf LEMMA} \ 
Let $c = a + \sqrt{-1} \hsx b$ $-$then $\forall \ z = x + \sqrt{-1} \hsx y$, 
\\[-1.cm]

\allowdisplaybreaks\begin{align*}
\Img \hsx
\left[
\frac{1}{z - c}
\hsx + \hsx 
\frac{1}{z - \bar{c}}
\right] \ 
&=\ 
-\Img 
\left[
\frac{z - c}{\abs{z - c}^2}
\hsx + \hsx 
\frac{z - \bar{c}}{\abs{z - \bar{c}}^2}
\right] 
\\[15pt]
&=\ 
-\Img \hsx
\left[
\frac
{(z - c) \hsy (z - \bar{c}) \hsy (\bar{z} - c) 
\hsx + \hsx 
 (z - \bar{c}) \hsy (z - c) \hsy (\bar{z} - \bar{c}) 
}
{\abs{z - c}^2 \hsy \abs{z - \bar{c}}^2}
\right] 
\\[15pt]
&=\ 
- 2 \hsy \Img \hsx
\left[
\frac
{(z - c) \hsy (z - \bar{c}) \hsy (\bar{z} - a)}
{\abs{z - c}^2 \hsy \abs{z - \bar{c}}^2}
\right] 
\\[15pt]
&=\ 
- 2 \hsy \Img 
\left[
\frac
{(z - a - \sqrt{-1} \hsx b) \hsx (z - a + \sqrt{-1} \hsx b) \hsx (\bar{z} - a)}
{\abs{z - c}^2 \hsy \abs{z - \bar{c}}^2}
\right] 
\\[15pt]
&=\ 
- 2 \hsy y \  
\frac
{\abs{z - a}^2 - b^2}
{\abs{z - c}^2 \hsy \abs{z - \bar{c}}^2}
\\[15pt]
&=\ 
- 2 \hsy y \  
\frac
{(x - a)^2 \hsx + \hsx y^2 \hsx - \hsx b^2}
{\abs{z - c}^2 \hsy \abs{z - \bar{c}}^2}
.
\end{align*}
\\[-1.25cm]
\end{x}

Given a real polynomial $f$, denote by 
$z_1, \ldots, z_\ell$ 
those zeros of $f$ which lie in the open upper half-plane.
\\[-.25cm]

\begin{x}{\small\bf DEFINITION} \ 
Put
\[
\gC_j 
\ = \ 
\{z \in \Cx : \abs{z - \Reg z_j} \leq \Img z_j \ (j = 1, \ldots \ell)\}.
\]
Then the $\gC_j$ are called the \un{Jensen circles} of $f$. 
\\[-.5cm]

[Note: \ 
The line segment joining the pair $z_j$, $\bar{z}_j$ is the vertical diameter of $\gC_j$.]

E.g.: \ 
If $f$ is a real polynomial with roots $z$ and $w$ in the upper half-plane, say
\\[-.25cm]

$
\begin{cases}
\ \ds
z =  -3 \hsx \frac{\sqrt{2}}{2}\ + \sqrt{2}\sqrt{-1}
\\[11pt]
\ \ds 
w =  4 \hsx \frac{\sqrt{3}}{2}\ + 2\sqrt{-1}
\end{cases}
$
then the Jensen circles of $f$ per $z$ and $w$ are: 
\[
\begin{tikzpicture}[scale=0.75]
      \draw[<->] (-6,0) -- (6,0) node[below] {$$};
      \draw[<->] (0,-3.75) -- (0,3.75) node[below] {$$};

      \draw[] (-2.1213, 1.4142) -- (-2.1213, 1.4142) node[above] {$z$};
      \draw[] (-2.1213, 1.4142)-- (-2.1213, 1.4142) node[] {\textbullet};
      \draw[] (-2.1213, -1.4142) -- (-2.1213, -1.4142) node[below] {$\bar{z}$};
      \draw[] (-2.1213, -1.4142) -- (-2.1213, -1.4142) node[] {\textbullet};
      \draw[blue!100] (-2.1213,0) circle(1.4142cm);
      
      \draw[] (-.1,-.3)-- (-.1,-.3) node[right] {\scriptsize $(0,0)$};
      
      \draw[] (3.4641, 2) -- (3.4641, 2) node[above] {$w$}; 
      \draw[] (3.4641, 2) -- (3.4641, 2) node[] {\textbullet};
      \draw[] (3.4641,-2) -- (3.4641,-2) node[below] {$\bar{w}$}; 
      \draw[] (3.4641,-2) -- (3.4641,-2) node[] {\textbullet};
      \draw[red!100] (3.4641,0) circle(2.cm);

\end{tikzpicture}
\]
\\[-1.25cm]
\end{x}

\begin{x}{\small\bf THEOREM} \ 
Let $f$ be a real polynomial $-$then the nonreal critical points of $f$ lie in the union
\[
\bigcup\limits_{j = 1}^\ell \ 
\gC_j
\]
of the Jensen circles of $f$. 
\\[-.5cm]

PROOF  \ 
Take $f$ monic of degree $n$, so 
\allowdisplaybreaks\begin{align*}
f (z) \ 
&=\ 
\prod\limits_{i = 1}^k \ 
(z - z_i)^{m_i} 
\\[15pt]
&=\ 
\prod\limits_{\Img z_i = 0} \ 
(z - z_i)^{m_i} 
\hsx \cdot \hsx
\prod\limits_{\Img z_i > 0} \
(z - z_i)^{m_i} 
\hsx 
(z - \bar{z}_i)^{m_i} 
\\[15pt]
&=\ 
\prod\limits_{\Img z_i = 0} \ 
(z - z_i)^{m_i} 
\hsx \cdot \hsx 
\prod\limits_{j = 1}^\ell \ 
(z - z_j)^{m_j} 
\hsx 
(z - \bar{z}_j)^{m_j}.
\end{align*}
Since the only issue is the position of the critical points of the second kind, pass to
\[
\frac{f^\prime (z)}{f (z)}
\ = \ 
\sum\limits_{\Img z_i = 0} \ 
\frac{m_i}{z - z_i}
\hsx + \hsx 
\sum\limits_{j = 1}^\ell \ 
m_j 
\hsx
\bigg[
\frac{1}{z - z_j} + \frac{1}{z - \bar{z}_j}
\bigg].
\]
Write
\[
z = x + \sqrt{-1} \hsx y
\quad \text{and} \quad 
z_j = x_j + \sqrt{-1} \hsx y_j
\qquad (j = 1, \ldots, \ell).
\]
Then
\[
\Img \frac{f^\prime (z)}{f (z)} 
\ = \ 
-y \ 
\left[
\
\sum\limits_{\Img z_i = 0} \ 
\frac{m_i}{\abs{z - z_i}^2} 
\ + \ 
2 \hsx 
\sum\limits_{j = 1}^\ell \ 
m_j \ 
\frac{(x - x_j)^2 + y^2 - y_j^2}{\abs{z - z_j}^2 \abs{z - \bar{z}_j}^2}
\hsx
\right] 
\quad (\tcf. \ 9.1)
.
\]
To say that $z \in \gC_j$ means that
\[
\abs{ x + \sqrt{-1} \hsx y - x_j}
\ \leq \ 
y_j
\]
or still, that
\[
(x - x_j)^2 + y^2 
\ \leq \ 
y_j^2.
\]
Therefore
\[
z \not\in \gC_j
\hsx \implies \hsx
(x - x_j)^2 + y^2 - y_j^2 
\ > \ 
0.
\]
Accordingly, outside the union of the $\gC_j$, at a $z$ with $y \neq 0$ we have
\[ 
\sgn \hsx \Img \frac{f^\prime (z)}{f (z)} 
\ = \ 
- \sgn \hsx y 
\ \neq \ 
0
\]
\qquad 
$\implies$
\[
f^\prime (z) 
\ \neq \ 
0.
\]
\\[-1.25cm]
\end{x}

Inspection of the preceding proof then leads to the following conclusion.
\\[-.25cm]

\begin{x}{\small\bf SCHOLIUM} \ 
A nonreal critical point of the second kind lies in the interior of at least one of the Jensen circles of $f$ 
unless it is a boundary point of each of them (in which case $f$ has no real zeros).
\\[-.25cm]
\end{x}

\begin{x}{\small\bf EXAMPLE} \ 
Let 
$
f 
=
z^5 - 5 z^4 + 5 z^3 + 35 z^2 + 130 z,
$
then the zeros of $f$ are 
\[ 
0, \ 
\begin{cases}
\ 
\alpha \ \approx \ -1.5509 \hsx + \hsx 1.6771 \hsy \sqrt{-1}
\\[4pt]
\ 
\bar{\alpha}  \ \approx \ -1.5509 \hsx - \hsx 1.6771 \hsy \sqrt{-1}
\end{cases}
\quad
\begin{cases}
\ 
\beta   \ \approx \ 4.0509 \hsx + \hsx 2.9160  \hsy \sqrt{-1}
\\[4pt]
\ 
\bar{\beta}   \ \approx \ 4.0509 \hsx - \hsx 2.9160 \hsy \sqrt{-1}
\end{cases}
.
\] 
And 
\allowdisplaybreaks\begin{align*}
f^\prime \ 
&= \ 
5 z^4 - 20z^3 + 15 z^2 + 70z + 130
\\[8pt]
&=\ 
5 \hsy(z^2 + 2 z + 2) \hsy (z^2 - 6 z + 13).
\end{align*}
\\[-.25cm]
Thus the critical points of $f$ are 
$
\ 
\begin{cases}
\ 
\mu \ = \ - 1 + \sqrt{-1}
\\[4pt]
\
\bar{\mu} \ = \ -1 -  \hsy\sqrt{-1}
\end{cases}
,\ 
\begin{cases}
\ 
\nu \ = 3 + 2 \hsy  \sqrt{-1}
\\[4pt]
\
\bar{\nu} \ = \ 3 - 2 \hsy\sqrt{-1}
\end{cases}
.
$
\\

\noindent
The Jensen circles, zeros, and critical points of $f$ are shown in the diagram:
\[
\begin{tikzpicture}[scale=0.65]
      \draw[<->] (-6,0) -- (9,0) node[below] {$$};
      \draw[<->] (0,-4.) -- (0,4.) node[below] {$$};

      \draw[] (-1.5509, 1.6771) -- (-1.5509, 1.6771) node[above] {$\alpha$};
      \draw[] (-1.5509, 1.6771)-- (-1.5509, 1.6771) node[] {\textbullet};
      \draw[] (-1.5509, -1.6771) -- (-1.5509, -1.6771) node[below] {$\bar{\alpha}$};
      \draw[] (-1.5509, -1.6771) -- (-1.5509, -1.6771) node[] {\textbullet};
      \draw[harvardcrimson, thick] (-1.5509,0) circle(1.6771cm);
      
      \draw[] (-.1,-.3)-- (-.1,-.3) node[right] {\scriptsize $(0,0)$};
      
      \draw[] (4.0509, 2.9160) -- (4.0509, 2.9160) node[above] {$\beta$}; 
      \draw[] (4.0509, 2.9160) -- (4.0509, 2.9160) node[] {\textbullet};
      \draw[] (4.0509,-2.9160) -- (4.0509,-2.9160) node[below] {$\bar{\beta}$}; 
      \draw[] (4.0509,-2.9160) -- (4.0509,-2.9160) node[] {\textbullet};
      \draw[harvardcrimson, thick] (4.0509,0) circle(2.9160cm);
      
      \draw[] (-1, 1) -- (-1, 1) node[] {{\color{red}\small \textbullet}};
      \draw[] (-1, 1) -- (-1, 1) node[left] {\small $\mu$};
      
      \draw[] (-1, -1) -- (-1, -1) node[] {{\color{red}\small \textbullet}};
      \draw[] (-1, -1) -- (-1, -1) node[left] {\small $\bar{\mu}$};

      \draw[] (3, 2) -- (3, 2) node[] {{\color{red}\small \textbullet}};
      \draw[] (3, 2) -- (3, 2) node[left] {\small $\nu$};
      
      \draw[] (3, -2) -- (3, -2) node[] {{\color{red}\small \textbullet}};
      \draw[] (3, -2) -- (3, -2) node[left] {\small $\bar{\nu}$};

      \draw[forestgreenweb, thick, -] (-1.5509, 1.6771) -- (4.0509, 2.9160);      
      \draw[forestgreenweb, thick, -] (-1.5509, -1.6771) -- (4.0509, -2.9160);
      \draw[forestgreenweb, thick, -] (-1.5509, 1.6771) -- (-1.5509, -1.6771);
      \draw[forestgreenweb, thick, -] (4.0509, +2.9160) -- (4.0509, -2.9160);
      
\end{tikzpicture}
\]
\\[-1.5cm]
\end{x}

\begin{x}{\small\bf LEMMA} \ 
Let $x_0$ be a point on the real line lying outside all the Jensen

circles of $f$.  
Assume: \ 
$f(x_0) = 0$ $-$then in each of the half-planes
\[
\begin{cases}
\ 
\{z \in \Cx : \Reg z < x_0\}
\\[4pt]
\ 
\{z \in \Cx : \Reg z > x_0\}
\end{cases}
,
\]
the number of zeros is the same as the number of critical points. 
\\[-.25cm]
\end{x}

\begin{x}{\small\bf LEMMA} \ 
Let $x_0$ be a point on the real line lying outside all the Jensen circles of $f$.  
Assume: \ 
$f (x_0) \neq 0$ $-$then in each of the half-planes
\[
\begin{cases}
\ 
\{z \in \Cx : \Reg z < x_0\}
\\[4pt]
\ 
\{z \in \Cx : \Reg z > x_0\}
\end{cases}
,
\]
the number of zeros is at least as large as the number of critical points 
(but can exceed it by at most one).
\\[-.25cm]
\end{x}

\begin{x}{\small\bf THEOREM} \ 
Let $a < b$ be two real numbers lying outside all the Jensen circles of $f$.  
Denote by $M$ the number of zeros and by $M^\prime$ the number of critical points in the strip
\[
\{z \in \Cx : a < \Reg z < b\}.
\]
Then
\\[-.25cm]

\qquad \textbullet \quad
$f (a) = 0$ and $f (b) = 0 \implies M^\prime = M + 1$.
\\[-.25cm]

\qquad \textbullet \quad
$f (a) = 0$ or $f (b) = 0 \implies M \leq M^\prime \leq M + 1$.
\\[-.25cm]

\qquad \textbullet \quad
$f (a) \neq 0$ and $f (b) \neq 0 \implies M - 1 \leq M^\prime \leq M + 1$.
\\[-.25cm]
\end{x}

\begin{x}{\small\bf EXAMPLE} \ 
The assumption that $a$ and $b$ lie outside all the Jensen circles of $f$ cannot be dropped.
\\[-.5cm]

[Take
\[
f (z) 
\ = \ 
z^4 + 4
\]
and let
\[
\begin{cases}
\ 
a = -1
\\[4pt]
\ 
b = 1
\end{cases}
, 
\quad \text{so} \quad 
\begin{cases}
\ 
f(a) \neq 0
\\[4pt]
\ 
f(b) \neq 0
\end{cases}
.
\]
Then $M = 0$ but $M^\prime = 3.$]
\\[-.25cm]
\end{x}


\chapter{
$\boldsymbol{\S}$\textbf{10}.\quad  CLASSES OF ENTIRE FUNCTIONS}
\setlength\parindent{2em}
\setcounter{theoremn}{0}
\renewcommand{\thepage}{\S10-\arabic{page}}

\qquad
Let $T$ be a nonempty closed subset of $\Cx$.
\\[-.25cm]

\begin{x}{\small\bf DEFINITION} \ 
A \un{$T$-polynomial} is a polynomial whose zeros are in $T$.
\\[-.25cm]
\end{x}

\begin{x}{\small\bf DEFINITION} \ 
A \un{$T$-function} is an entire function $\not\equiv 0$ which is the uniform limit 
on compact subsets of $\Cx$ of a sequence of $T$-polynomials.
\\[-.25cm]
\end{x}

\begin{x}{\small\bf NOTATION} \ 
Let 
\[
\ent (T)
\]
stand for the class of $T$-functions.
\\[-.25cm]
\end{x}

\qquad
{\small\bf \un{N.B.}} \ 
The product of two $T$-functions is a $T$-function.
\\[-.25cm]

\begin{x}{\small\bf LEMMA} \ 
If $f \in \ent (T)$, then all its zeros lie in $T$.  
\\[-.5cm]

[Note: \ 
As will be seen below (cf. 10.14), the converse of this assertion is false: 
An entire function whose zeros are in $T$ need not belong to $\ent(T)$.]
\\[-.25cm]
\end{x}

\begin{x}{\small\bf LEMMA} \ 
If $T$ is bounded, then $\ent(T)$ is the set of $T$-polynomials.  
\\[-.5cm]

PROOF \ 
Let $f \in \ent(T)$, and suppose that $f_n \ra f$ uniformly on compact subsets of $\Cx$, 
where $\{f_n\}$ is a sequence of $T$-polynomials.  
Since all the zeros of $f$ lie in $T$ and since $T$ is bounded, 
their number is finite, call it $N$.  
By Rouch\'e's theorem, 
the number of zeros of $f_n$ is also $N$ provided $n \gg 0$, 
thus the $f_n$ are of degree $N$ provided $n \gg 0$.  
But the Taylor coefficients of $f$ are the limits of the Taylor coefficients of the $f_n$, 
hence $f$ is a polynomial of degree $N$.
\\[-.25cm]
\end{x}

Abstractly, the problem then is to characterize $\ent(T)$ in terms of the properties
of $T$.  
This can be done (more or less) but instead of delving into the general theory, 
we shall consider only those special cases that will be needed later on, 
namely: \ 
\[
\begin{cases}
\ 
T = \ ]-\infty, 0] 
\quad \text{or} \quad [0, +\infty[
\\[8pt]
\ 
T = \ ]-\infty, +\infty[
\end{cases}
\]
{subject to the restriction that here 
``$T$-polynomials'' and ``$T$-functions'' are real (so, e.g., $\sqrt{-1} \hsx (z^2 - 1)$ is not a 
$T$-polynomial even though its zeros are real).
\\[-.25cm]

\begin{x}{\small\bf LEMMA} \ 
We have 
\[
\begin{cases}
\ 
\ent(\hsx]-\infty, 0]\hsx)
\\[8pt]
\ 
\ent(\hsx[0, +\infty[\hsx)
\end{cases}
\subsetx \ent(\hsx]-\infty, +\infty[\hsx).
\]

[This is obvious.]
\\[-.25cm]
\end{x}

\begin{x}{\small\bf EXAMPLE} \ 
If $f = C$ $(C \neq 0)$, then $f \in \ent(\hsx[0, +\infty[\hsx)$.
\\[-.5cm]

[Consider
\[
C \hsx 
\Big(
1 - \frac{z}{k}
\Big)
\qquad (k = 1, 2, \ldots).]
\]
\\[-1.cm]
\end{x}

\begin{x}{\small\bf EXAMPLE} \ 
Since
\[
e^{-z} 
\ = \ 
\lim\limits_{n \ra \infty} \ 
\Big(
1 - \frac{z}{n}
\Big)^n, 
\]
it follows that
\[
e^{-z} \in \ent(\hsx[0, +\infty[\hsx).
\]
\\[-1.25cm]
\end{x}

\begin{x}{\small\bf EXAMPLE} \ 
The zeros of 
\[
\Big(
1 - \frac{z^2}{n^2}
\Big)
\]
are $z = \pm n$, so
\[
\prod\limits_{n = 1}^N \ 
\Big(
1 - \frac{z^2}{n^2}
\Big)
\ \in  \ 
\ent(\hsx]-\infty, +\infty[\hsx), 
\]
which implies that
\[
\frac{\sin \pi z}{\pi z}
\in  \ent(\hsx]-\infty, +\infty[\hsx) 
\qquad (\tcf. \ 1.23).
\]
\\[-1.cm]
\end{x}

\begin{x}{\small\bf EXAMPLE} \ 
The zeros of the Laguerre polynomials (cf. 8.17) are real and positive, hence $\forall \ n$, 
\[
L_n \in \ent(\hsx[0, +\infty[\hsx).
\]
Consider now the Bessel function of index 0:
\[
\tJ_0 (z) 
\ = \ 
1 
- 
\frac{1}{1 ! \hsy 1 !} 
\hsx
\Big(
\frac{z}{2}
\Big)^2
+ 
\frac{1}{2 ! \hsy 2 !} 
\hsx
\Big(
\frac{z}{2}
\Big)^4
-
\frac{1}{3 ! \hsy 3 !} 
\hsx
\Big(
\frac{z}{2}
\Big)^6
+ \cdots \hsx .
\]
Then
\[
\tJ_0 (z) 
\ = \ 
\lim\limits_{n \ra \infty} \ 
L_n
\Big(
\frac{z^2}{4 \hsy n}
\Big)
\]
uniformly on compact subsets of $\Cx$, thus 
\[
\tJ_0 (z) \in \ent(\hsx[0, +\infty[\hsx).
\]

[In fact, 
\[
L_n 
\Big(
\frac{z^2}{4 \hsy n}
\Big)
\ = \ 
1 
- 
\frac{z^2}{2 \cdot 2}
+ 
\frac{z^4}{2 \cdot4 \cdot 2 \cdot 4}
\hsx
\Big(
1 - \frac{1}{n}
\Big)
- 
\frac{z^6}{2 \cdot4 \cdot 6 \cdot 2 \cdot4 \cdot 6}
\hsx
\Big(
1 - \frac{1}{n}
\Big)
\Big(
1 - \frac{2}{n}
\Big)
+ \cdots \hsx .]
\]
\\[-1.cm]
\end{x}

\begin{x}{\small\bf THEOREM} \ 
Let $f \not\equiv 0$ be a real entire function $-$then $f \in \ent(\hsx[0, +\infty[\hsx)$
iff $f$ has a representation of the form
\[
f(z) 
\ = \ 
C \hsy z^m \hsy e^{a z} \ 
\prod\limits_{n = 1}^\infty \ 
\Big(
1 - \frac{z}{\lambda_n}
\Big),
\]
where $C \neq 0$ is real, 
$m$ is a nonnegative integer, 
a is real and $\leq 0$, 
the $\lambda_n$ are
real and $> 0$ with 
$
\ 
\ds 
\sum\limits_{n = 1}^\infty \ 
\frac{1}{\lambda_n} 
< 
\infty.
$
\\[-.5cm]
\end{x}

\begin{spacing}{1.75}
[Note: \ 
Functions having finitely many zeros are accommodated by the convention that 
$\lambda_n = \infty$ and 
$
\ds
0 = \frac{1}{\lambda_n}
$ 
$(n \geq n_0)$ and an empty product is taken to be 1.]
\end{spacing}

\begin{x}{\small\bf REMARK} \ 
$\ent(\hsx[0, +\infty[\hsx)$ 
is closed under differentiation (cf. 8.16).
\\[-.25cm]
\end{x}

\begin{x}{\small\bf REMARK} \ 
Let 
$f \in \ent(\hsx[0, +\infty[\hsx)$ $-$then $\fg = 0$, so
\[
\gen f \ = \quad
\begin{cases}
\ 
0 \ \text{if} \ a = 0
\\[4pt]
\ 
1 \ \text{if} \ a \neq 0
\end{cases}
\]
and $\rho \leq 1$. 
\\[-.25cm]
\end{x}

\begin{x}{\small\bf EXAMPLE} \ 
The real entire function 
\[
e^{-z^2} \ 
\prod\limits_{n = 1}^\infty \ 
\Big(
1 - \frac{z}{n^2}
\Big)
\]
has zeros in $\ [0, +\infty[\ $ but does not belong to $\ent(\hsx[0, +\infty[\hsx)$.
\\[-.25cm]
\end{x}

That the conditions of 10.11 are necessary is straightforward: 
Consider
\[
p_k (z) 
\ = \ 
C \hsx 
\Big(
1 - \frac{z}{k}
\Big)
\hsx
\Big(
z - \frac{1}{k}
\Big)^m
\hsx
\Big(
1 + \frac{a z}{k}
\Big)^k
\hsx
\prod\limits_{n = 1}^k \ 
\Big(
1 - \frac{z}{\lambda_n}
\Big)
\hsx .
\]

This said, suppose now that 
$f \in \ent(\hsx[0, +\infty[\hsx)$ 
and write
\[
f (z) 
\ = \ 
a_0 - a_1 z + a_2 z^2 - \cdots \hsy .
\]
Let
\[
p_k (z) 
\ = \ 
a_{k \hsy 0} -a_{k \hsy 1} z + a_{k \hsy 2} z^2 - \cdots + (-1)^k a_{k \hsy k}z^k
\]
be a sequence of polynomials whose zeros are real and positive such that $p_k \ra f$
uniformly on compact subsets of $\Cx$ $-$then
\[
\lim\limits_{k \ra \infty} \ 
a_{k \hsy \ell} 
\ = \ 
a_\ell.
\]
\\[-1.25cm]

\begin{x}{\small\bf REDUCTION} \ 
There is no loss of generality in assuming that $a_0 \neq 0$.
\\[-.5cm]

[Fix a positive real number $\alpha$ which is smaller than the smallest positive zero of $f$ 
(cf. 10.4), 
pass to $f (z + \alpha)$, and note that $f (\alpha) \neq 0$.]
\\[-.25cm]

Therefore one can work instead with 

\[
\frac{f(z)}{a_0}, \ 
\frac{p_k (z)}{a_{k \hsy 0}}
\qquad 
\text{(since \ $\lim\limits_{k \ra \infty} \ a_{k \hsy 0} \ = \ a_0 \ \neq \ 0$)}.
\]
So, recast, 
\[
f(z) 
\ = \ 
1 - a_1 z + a_2 z^2 - \cdots
\]
and 
\allowdisplaybreaks
\allowdisplaybreaks\begin{align*}
p_k (z) \
&= \ 
1 - a_{k \hsy 1}  z + a_{k \hsy 2} z^2 - \cdots + (-1)^k a_{k \hsy k}z^k
\\[15pt]
&\equiv 
\Big(1 - \frac{z}{\lambda_{k \hsy 1}}\Big)
\hsx
\Big(1 - \frac{z}{\lambda_{k \hsy 2}}\Big)
\cdots 
\Big(1 - \frac{z}{\lambda_{k \hsy k}}\Big),
\end{align*}
where the zeros 
$\lambda_{k \hsy \ell} \neq 0$
are positive and 

\[
0 
\ < \ 
\lambda_{k \hsy 1}
\ \leq \ 
\lambda_{k \hsy 2}
\ \leq \ 
\cdots
\ \leq \ 
\lambda_{k \hsy k}.
\]
\\[-1.25cm]
\end{x}

\qquad
{\small\bf \un{N.B.}} \ 
The $a_k$ and the $a_{k \hsy \ell}$ are nonnegative.
\\

\begin{x}{\small\bf LEMMA}\footnote[2]{\vspace{.11 cm}
O. Schl\"omilch,  
\textit{Zeitschr. f. Math. und Physik} \textbf{3} (1858), pp. 301-308, 
(see page 308, formula 15).} \ 
Let 

\[
\Phi (z) 
\ = \ 
1 - c_1 z + c_2 z^2 - \cdots + (-1)^n c_n z^n
\]
be a real polynomial whose zeros are real and positive $-$then
\[
\frac{c_1}{n} 
\ \geq \ 
\left[
\raisebox{.1cm}
{
$
\ds
\frac{c_2}{\raisebox{-.1cm}{$\binom{n}{2}$}}
$
}
\right]^{1/2}
\ \geq \ 
\cdots
\ \geq  
\left[
\raisebox{.1cm}
{
$
\ds
\frac{c_p}{\raisebox{-.1cm}{$\binom{n}{p}$}}
$
} 
\right]^{1/2}
\ \geq \ 
\cdots
\ \geq \ 
(c_n)^{1/n}.
\]
\\[-1cm]

Take $\Phi = p_k$, thus

\[
\frac{a_{k \hsy \ell}}{k} 
\ \geq \ 
\left[
\raisebox{.1cm}
{
$
\ds
\frac{a_{k \hsy \ell}}{\raisebox{-.1cm}{$\binom{k}{\ell}$}}
$
} 
\right]^{1/\ell}
\]
\qquad 
$\implies$
\[
(a_{k \hsy \ell})^\ell \ 
\frac{k (k -1) \cdots (k - \ell + 1)}{k^\ell}
\ \geq \ 
a_{k \hsy \ell}, 
\]

\noindent
so in the limit as $k \ra \infty$, 

\[
\frac{(a_1)^\ell}{\ell !}
\ \geq \ 
a_\ell.
\]
\\[-1.cm]
\end{x}

\begin{x}{\small\bf LEMMA} \ 
$f$ is of finite order $\rho \leq 1$. 
\\[-.25cm]

PROOF \ 
In fact,
\allowdisplaybreaks\begin{align*}
\abs{f(z)} \ 
&\leq \ 
\sum\limits_{\ell = 0}^\infty \ 
a_\ell 
\hsx \abs{z}^\ell
\\[15pt]
&\leq 
\sum\limits_{\ell = 0}^\infty \ 
\frac{(a_1)^\ell}{\ell !} 
\hsx \abs{z}^\ell
\\[15pt]
&= 
\exp (a_1 \abs{z})
\end{align*}

\qquad 
$\implies$

\[
M(r; f) 
\ \leq \ 
\exp a_1 r,
\]
from which the assertion (cf. 2.15).
\\[-.5cm]
\end{x}

Enumerate the zeros of $f$ in the usual way:
\[
0 
\ < \ 
\lambda_1 
\ \leq \ 
\lambda_2 
\ \leq \ 
\cdots \hsx .
\]
Then
\[
\lim\limits_{k \ra \infty} \ 
\lambda_{k \hsy \ell} 
\ = \ 
\lambda_\ell. 
\]
But
\allowdisplaybreaks
\allowdisplaybreaks\begin{align*}
a_{k \hsy 1} \ 
&=\ 
\frac{1}{\lambda_{k \hsy 1}}
\hsx + \hsx
\frac{1}{\lambda_{k \hsy 2}}
\hsx + \hsx
\cdots
\hsx + \hsx
\frac{1}{\lambda_{k \hsy k}}
\\[11pt]
&\geq\
\frac{1}{\lambda_{k \hsy 1}}
\hsx + \hsx
\frac{1}{\lambda_{k \hsy 2}}
\hsx + \hsx
\cdots
\hsx + \hsx
\frac{1}{\lambda_{k \hsy \ell}}
\end{align*}
\qquad
$\implies$
\allowdisplaybreaks\begin{align*}
a_1 \ 
&=\ 
\lim\limits_{k \ra \infty} \ 
a_{k \hsy 1} 
\\[15pt]
&\geq\
\lim\limits_{k \ra \infty} \ 
\Big(
\frac{1}{\lambda_{k \hsy 1}}
\hsx + \hsx
\frac{1}{\lambda_{k \hsy 2}}
\hsx + \hsx
\cdots
\hsx + \hsx
\frac{1}{\lambda_{k \hsy \ell}}
\Big)
\\[15pt]
&=\
\frac{1}{\lambda_1}
\hsx + \hsx
\frac{1}{\lambda_2}
\hsx + \hsx
\cdots
\hsx + \hsx
\frac{1}{\lambda_\ell}.
\end{align*}
Therefore the series 
$
\ 
\ds
\frac{1}{\lambda_1}
\hsx + \hsx
\frac{1}{\lambda_2}
+ \cdots
\ 
$ 
converges and 
\[
\sum\limits_{n = 1}^\infty \ 
\frac{1}{\lambda_n}
\ \leq \ 
a_1.
\]

Proceeding, write
\[
f(z) 
\ = \ 
e^{Q(z)} \ 
\prod\limits_{n = 1}^\infty \ 
\Big(
1 - \frac{z}{\lambda_n}
\Big)
\qquad (\tcf. \ 7.3), 
\]
where $q \leq \rho \leq 1$ and $\fg = 0$, hence
\[
\gen f 
\ = \ 
\max (q, \fg) 
\ = \ q.
\]
And
\[
Q(z) 
\ = \ 
a \hsy z + b, 
\]
the final claim being that $a$ is real and $\leq 0$.
\\[-.25cm]

[Note: \ 
$
\ds
1 = f(0) = 
e^b \ 
\prod\limits_{n = 1}^\infty \ 
1
= 
e^b.]
$
\\[-.25cm]

However
\[
1 - a_1 z + \cdots 
\ = \ 
(1 + a z + \cdots ) 
\hsy 
\Big(
1 - 
\Big(
\sum\limits_{n = 1}^\infty \ 
\frac{1}{\lambda_n} 
\Big)
\hsy z 
+ \cdots
\Big)
\]
\qquad 
$\implies$
\[
- a_1 
\ = \ 
a - \ 
\sum\limits_{n = 1}^\infty \ 
\frac{1}{\lambda_n} 
\]
\qquad 
$\implies$
\allowdisplaybreaks\begin{align*}
a \
&= \ 
-a_1 +\ 
\sum\limits_{n = 1}^\infty \ 
\frac{1}{\lambda_n} 
\\[11pt]
&\leq \ 
0,
\end{align*}
thereby completing the proof of 10.11.
\\[-.25cm]

\begin{x}{\small\bf REMARK} \ 
The fact that $f$ is of finite order $\rho \leq 1$ was established by appealing to 10.16.
This can be avoided.  
Indeed, $
\{a_{k \hsy 1} : k = 1, 2, \ldots\}$ 
converges to $a_1$, hence is bounded, say $0 \leq a_{k \hsy 1} \leq M$, hence
\allowdisplaybreaks\begin{align*}
\abs{p_k (z)} \ 
&\leq \ 
\sum\limits_{\ell = 1}^k \ 
\abs{1 - \frac{z}{\lambda_{k \hsy \ell}}}
\\[11pt]
&\leq \ 
\sum\limits_{\ell = 1}^k \ 
\Big(
1 - \frac{\abs{z}}{\lambda_{k \hsy \ell}}
\Big)
\\[11pt]
&\leq \ 
\exp
\Big(
\abs{z} \ 
\sum\limits_{\ell = 1}^k \ 
\frac{1}{\lambda_{k \hsy \ell}}
\Big)
\\[11pt]
&\leq \ 
\exp (\abs{z} \hsy a_{k \hsy 1})
\\[11pt]
&\leq \ 
\exp (M \abs{z}).
\end{align*}
And then 
\[
\abs{f(z)}
\ = \ 
\lim\limits_{k \ra \infty} \ 
\abs{p_k (z)}
\ \leq \ 
\exp (M \abs{z}).
\]
\\[-1.cm]
\end{x}

\begin{x}{\small\bf THEOREM} \ 
Let $f \not\equiv 0$ be a real entire function 
$-$then 
$f \in \ent(\hsx]-\infty$, $+\infty[\hsx)$
iff $f$ has a representation of the form 
\[
f(z) 
\ = \ 
C \hsy z^m \hsy e^{a z^2 + b z} \ 
\prod\limits_{n = 1}^\infty \ 
\Big(
1 - \frac{z}{\lambda_n}
\Big)^{z/\lambda_n},
\]
where $C \neq 0$ is real, $m$ is a nonnegative integer, $a$ is real and $\leq 0$, $b$ is real, 
the $\lambda_n$ are real with 
$
\ 
\ds
\sum\limits_{n = 1}^\infty \ 
\frac{1}{\lambda_n^2}
< \infty.
$
\\[-.25cm]

[Note: \ 
Functions having finitely many zeros are accommodated by the convention that 
$\lambda_n = \infty$ 
and 
$
\ds
0 = \frac{1}{\lambda_n}
$ 
$(n \geq n_0)$ 
and an empty product is taken to be 1.]
\\[-.25cm]
\end{x}

\begin{x}{\small\bf REMARK} \ 
$\ent(\hsx ]-\infty, +\infty[\hsx)$  
is closed under differentiation (cf. 8.16).
\\[-.25cm]
\end{x}

\begin{x}{\small\bf REMARK} \ 
Let 
$f \in \ent(\hsx]-\infty, +\infty[\hsx)$.
\\[-.25cm]

\qquad \textbullet \quad
$\fg = 0 \implies \gen f = 0, \hsx 1, \hsx 2$.
\\[-.25cm]

\qquad \textbullet \quad
$\fg = 1 \implies \gen f = 1, \hsx 2$.
\\[-.25cm]

To see that the conditions of 10.19 are necessary, introduce
\[
\Lambda_k
\ = \ 
b 
\hsx + \hsx 
\sum\limits_{n = 1}^k \ 
\frac{1}{\lambda_n}
\]
and let
\[
p_k (z) 
\ = \ 
C \hsy
\Big(1 - \frac{z}{k}\Big)
\hsx 
\Big(z - \frac{1}{k}\Big)^m
\hsx 
\Big(1 + \frac{a z^2}{k}\Big)^k
\hsx 
\Big(1 + \frac{\Lambda_k \hsy z}{n_k}\Big)^{n_k}
\ 
\prod\limits_{n = 1}^k \ 
\Big(1 - \frac{z}{\lambda_n}\Big),
\]
where the $n_k \ra \infty$ $(k \ra \infty)$ are chosen subject to
\[
\abs{z} \leq k 
\hsx \implies \hsx 
\
\abs{
\Big(1 + \frac{\Lambda_k \hsy z}{n_k}\Big)^{n_k}
- e^{\Lambda_k z}
}
\ < \  
\frac{1}{k} \hsx
\exp
\Big(
-k \ 
\sum\limits_{n = 1}^k \ 
\frac{1}{\abs{\lambda_n}}
\Big)\hsx .
\]
\\[-1cm]

Turning to the sufficiency, let 
$f \in \ent(\hsx]-\infty, +\infty[\hsx)$ 
and normalize the situation so that as before
\[
f(z) 
\ = \ 
1 - a_1 z + a_2 z^2 - \cdots
\]
and
\allowdisplaybreaks\begin{align*}
p_k (z) \ 
&=\ 
1 - a_{k \hsy 1} \hsy z + a_{k \hsy 2} \hsy z^2 - \cdots + (-1)^k a_{k \hsy k} \hsy z^k
\\[11pt]
&\equiv \ 
\Big(1 - \frac{z}{\lambda_{k \hsy 1}}\Big)
\hsx
\Big(1 - \frac{z}{\lambda_{k \hsy 2}}\Big)
\cdots 
\Big(1 - \frac{z}{\lambda_{k \hsy k}}\Big),
\end{align*}
where the zeros $\lambda_{k \hsy \ell} \neq 0$ are real and 
\[
0 
\ < \ 
\abs{\lambda_{k \hsy 1}}
\ \leq \ 
\abs{\lambda_{k \hsy 2}}
\ \leq \ 
\cdots
\ \leq \ 
\abs{\lambda_{k \hsy k}}.
\]
\\[-1.25cm]
\end{x}

\begin{x}{\small\bf SUBLEMMA} \ 
$\forall$ complex $z$,
\[
\abs{(1 + z) e^{-z}}
\ \leq \ 
e^{4 \hsy \abs{z}^2}.
\]

PROOF \ 
If 
$
\ds
\abs{z} \leq \frac{1}{2}
$, 
then
\[
\abs{(1 + z) e^{-z}}
\ \leq \ 
e^{\abs{z}^2}
\ \leq \ 
e^{4 \hsy \abs{z}^2}.
\]
On the other hand, if 
$
\ds 
\abs{z} \geq \frac{1}{2}
$, 
then
\allowdisplaybreaks\begin{align*}
\abs{(1 + z) e^{-z}} \ 
&\leq \ 
(1 + \abs{z}) e^{\abs{z}}
\\[11pt]
&\leq \ 
e^{2 \hsy \abs{z}} 
\\[11pt]
&\leq \ 
e^{4 \hsy \abs{z}^2}.
\end{align*}

From the definitions, 
\[
a_1 
\ = \ 
\lim\limits_{k \ra \infty} \ 
a_{k 1} 
\ = \ 
\lim\limits_{k \ra \infty} \ 
\sum\limits_{\ell= 1}^k \ 
\frac{1}{\lambda_{k \hsy \ell}}.
\]
Next
\allowdisplaybreaks\begin{align*}
a_{k \hsy 2} \ 
&=\ 
\sum\limits_{i < j} \ 
\frac{1}{\lambda_{k \hsy i}}
\hsx 
\frac{1}{\lambda_{k \hsy j}}
\\[15pt]
&=\ 
\frac{1}{2} \ 
\sum\limits_{i \neq j} \ 
\frac{1}{\lambda_{k \hsy i}}
\hsx 
\frac{1}{\lambda_{k \hsy j}}.
\end{align*}
But
\allowdisplaybreaks\begin{align*}
\Big(
\sum\limits_{i= 1}^k \ 
\frac{1}{\lambda_{k \hsy i}}
\Big)
\hsx 
\Big(
\sum\limits_{j= 1}^k \ 
\frac{1}{\lambda_{k \hsy j}}
\Big)
&=\ 
\sum\limits_{\ell= 1}^k \ 
\frac{1}{\lambda_{k \hsy \ell}^2}
\hsx + \hsx 
\sum\limits_{i \neq j}  \
\frac{1}{\lambda_{k \hsy i}}
\hsx 
\frac{1}{\lambda_{k \hsy j}}
\\[15pt]
&=\ 
\sum\limits_{\ell= 1}^k \ 
\frac{1}{\lambda_{k \hsy \ell}^2}
\hsx + \hsx 
2 \ 
\sum\limits_{i < j} \
\frac{1}{\lambda_{k \hsy i}}
\hsx 
\frac{1}{\lambda_{k \hsy j}}.
\end{align*}
So, upon letting $k \ra \infty$, we get
\[
a_1^2 
\ = \ 
\lim\limits_{k \ra \infty} \ 
\sum\limits_{\ell= 1}^k \ 
\frac{1}{\lambda_{k \hsy \ell}^2}
\ + \ 
2 a_2
\]
or still, 
\[
a_1^2 - 2 \hsy a_2
\ = \ 
\lim\limits_{k \ra \infty} \ 
\sum\limits_{\ell= 1}^k \ 
\frac{1}{\lambda_{k \hsy \ell}^2}.
\]
\\[-.25cm]

Fix constants
$
\begin{cases}
\ 
U > 0
\\[4pt]
\ 
V > 0
\end{cases}
$
such that $\forall \ k$, 
\[
\begin{cases}
\ \ 
\ds
\bigg| \hsx
\sum\limits_{\ell= 1}^k \ \frac{1}{\lambda_{k \hsy \ell}}
\hsx
\bigg| 
\ \leq \ U
\\[18pt]
\ \ 
\ds
\sum\limits_{\ell= 1}^k \ 
\frac{1}{\lambda_{k \hsy \ell}^2}
\ \leq \ V
\end{cases}
.
\]
\\[-.25cm]
\end{x}

\begin{x}{\small\bf LEMMA} \ 
We have
\[
\abs{p_k (z)}
\ \leq \ 
\exp (U \abs{z}\hsx + \hsx 4 \hsy V \abs{z}^2).
\]

PROOF \ 
Write
\allowdisplaybreaks
\allowdisplaybreaks\begin{align*}
\abs{p_k (z) \hsy e^{a_{k \hsy 1} \hsy z}} \ 
&=\ 
\bigg| \hsx
p_k (z)
\hsy 
\exp
\Big(
\sum\limits_{\ell= 1}^k \ 
\frac{z}{\lambda_{k \hsy \ell}}
\Big)
\hsx
\bigg| 
\\[15pt]
&=\ 
\bigg| \hsx
\prod\limits_{\ell= 1}^k \ 
\Big(
1 - \frac{z}{\lambda_{k \hsy \ell}}
\Big)
\exp
\Big(
\frac{z}{\lambda_{k \hsy \ell}}
\Big)
\hsx
\bigg| 
\\[15pt]
&\leq\ 
\prod\limits_{\ell= 1}^k \ 
\abs
{
\Big(
1 - \frac{z}{\lambda_{k \hsy \ell}}
\Big)
\exp
\Big(
\frac{z}{\lambda_{k \hsy \ell}}
\Big)
}
\\[15pt]
&\leq\ 
\prod\limits_{\ell= 1}^k \ 
\exp 
\Big( 
4 
\hsx 
\Big| \hsx
\frac{z}{\lambda_{k \hsy \ell}}
\hsx \Big|^2
\hsx
\Big)
\qquad (\tcf. \ 10.22)
\\[15pt]
&\leq\ 
\exp 
\Big( 
4 
\hsx 
\Big(
\sum\limits_{\ell= 1}^k \ 
\frac{1}{\lambda_{k \hsy \ell}^2}
\Big)
\hsx
\abs{z}^2
\Big)
\\[15pt]
&\leq\ 
\exp (4 \hsy V \abs{z}^2).
\end{align*}
Therefore
\allowdisplaybreaks
\allowdisplaybreaks\begin{align*}
\abs{p_k (z)} \ 
&=\ 
\abs{p_k (z) \hsx e^{a_{k \hsy 1} z} \hsx e^{-a_{k \hsy 1} z}}
\\[15pt]
&\leq\ 
\abs{p_k (z) \hsx e^{a_{k \hsy 1} z}}
\hsx 
\abs{e^{-a_{k \hsy 1} z}}
\\[15pt]
&\leq\ 
\exp (4 \hsy V \abs{z}^2)
\hsx 
\exp (\abs{a_{k \hsy 1}} \hsy  \abs{z})
\\[15pt]
&\leq\ 
\exp (U  \abs{z} \hsx + \hsx 4 \hsy V \abs{z}^2).
\end{align*}

Consequently, $f$ is of finite order $\rho \leq 2$ (cf. 10.18).
\\[-.25cm]
\end{x}

\begin{x}{\small\bf LEMMA} \ 
If 
$\lambda_1, \lambda_2, \ldots$ 
are the zeros of $f$ and if 
\[
0 
\ \leq \ 
\abs{\lambda_1}
\ \leq \ 
\abs{\lambda_2}
\ \leq \ 
\cdots, 
\]
then
\[
\lim\limits_{k \ra \infty} \ 
\lambda_{k \hsy \ell}
\ = \ 
\lambda_{\ell}
\]
and
\[
a_1^2 - 2 \hsy a_2 \ 
\ \geq \ 
\sum\limits_{n =1}^\infty \
\frac{1}{\lambda_n^2}.
\]

PROOF \ 
Start by writing 
\[
\frac{1}{\lambda_{k \hsy 1}^2} + \frac{1}{\lambda_{k \hsy 2}^2} + \cdots + \frac{1}{\lambda_{k \hsy k}^2}
\ \geq \ 
\frac{1}{\lambda_{k \hsy 1}^2} + \frac{1}{\lambda_{k \hsy 2}^2} + \cdots + \frac{1}{\lambda_{k \hsy \ell}^2}
\]
and then let $k \ra \infty$, hence
\allowdisplaybreaks\begin{align*}
a_1^2 - 2 \hsy a_2 \ 
&=\ 
\lim\limits_{k \ra \infty} \ 
\Big(
\sum\limits_{\ell =1}^k \
\frac{1}{\lambda_{k \hsy \ell}^2}
\Big)
\\[15pt]
&\geq \ 
\lim\limits_{k \ra \infty} \ 
\Big(
\frac{1}{\lambda_{k \hsy 1}^2} + \frac{1}{\lambda_{k \hsy 2}^2} + \cdots + \frac{1}{\lambda_{k \hsy \ell}^2}
\Big)
\\[15pt]
&=\ 
\frac{1}{\lambda_1^2} + \frac{1}{\lambda_2^2} + \cdots + \frac{1}{\lambda_\ell^2},
\end{align*}
which implies that 
\[
a_1^2 - 2 \hsy a_2 
\ \geq \ 
\sum\limits_{n =1}^\infty \
\frac{1}{\lambda_n^2} \hsx .
\]

Accordingly, 
\[
\sum\limits_{n =1}^\infty \
\frac{1}{\lambda_n^2}
\ < \ 
\infty 
\qquad (\implies \fg = 0 \ \text{or} \ 1)
\]
and the product
\[
\prod\limits_{n =1}^\infty \
\Big(
1 - \frac{z}{\lambda_n}
\Big)
e^{z /  \lambda_n}
\]
is an entire function whose zeros are the $\lambda_n$ (cf. 5.4).
To see that its order is also $\leq 2$, write
\allowdisplaybreaks\begin{align*}
\bigg|
\prod\limits_{n =1}^\infty \
\Big(
1 - \frac{z}{\lambda_n}
\Big)
e^{z /  \lambda_n}
\bigg|
&\leq \ 
\prod\limits_{n =1}^\infty \
\abs{
\Big(
1 - \frac{z}{\lambda_n}
\Big)
e^{z /  \lambda_n}}
\\[15pt]
&\leq \ 
\prod\limits_{n =1}^\infty \
\exp
\Big(
4 \hsx 
\frac{\abs{z}^2}{\lambda_n^2}
\Big)
\qquad (\tcf. \ 10.22)
\\[15pt]
&\leq \ 
\exp 
\Big(
4 \ 
\Big(
\sum\limits_{n =1}^\infty \
\frac{1}{\lambda_n^2}
\Big)
\abs{z}^2
\Big).
\end{align*}

Thanks to 2.37, the order of 
\[
\frac{f(z)}
{\ds\prod\limits_{n =1}^\infty \
\Big(
1 - \frac{z}{\lambda_n}
\Big)
e^{z /  \lambda_n}}
\]
is $\leq$ the maximum of $\rho$ and the order of 
\[
\prod\limits_{n =1}^\infty \
\Big(
1 - \frac{z}{\lambda_n}
\Big)
e^{z /  \lambda_n},
\]
thus is $\leq 2$, so 
\[
\frac{f(z)}
{
\ds\prod\limits_{n = 1}^\infty \
\Big(
1 - \frac{z}{\lambda_n}
\Big)
e^{z /  \lambda_n}
}
\ = \ 
e^{Q(z)},
\]
where 
\[
Q(z) 
\ = \ 
a  z^2 + bz + c
\]
is a polynomial of degree $\leq 2$ (cf. 2.42).
\\[-.25cm]

[Note: \ 
$
\ds
1 = f(0) = 
e^c \ 
\prod\limits_{n = 1}^\infty \
1 
= e^c.]
$
\\

There remain the claims that 
(1) $b$ is real and 
(2) $a$ is real and $\leq 0$.  
To this end, compare coefficients: 
\\[-.5cm]

\qquad (1) \quad 
$b = -a_1 = \lim\limits_{k \ra \infty} \ a_{k \hsy 1}$, 
\
which is real.
\\[-.25cm]

\qquad (2) \quad 
$\ds
a = -\frac{1}{2} \hsx 
\Big(
a_1^2 - 2 \hsy a_2 
\hsx - \hsx 
\sum\limits_{n = 1}^\infty \ 
\frac{1}{\lambda_n^2}
\Big)
$
\\[-.25cm]

\noindent
and
\[
a_1^2 - 2 \hsy a_2 
\hsx - \hsx 
\sum\limits_{n = 1}^\infty \ 
\frac{1}{\lambda_n^2}
\ \geq \  
0
\qquad (\tcf. \ 10.24).
\]

The proof of 10.19 is therefore complete.
\\[-.25cm]
\end{x}

\qquad
{\small\bf \un{N.B.}} \ 
Take an $f \in \ent(\hsx[0, +\infty[\hsx)$ and write 
\[
f(z) 
\ = \ 
C \hsy z^m \hsy e^{a z} \ 
\prod\limits_{n = 1}^\infty \ 
\Big(
1 - \frac{z}{\lambda_n}
\Big)
\qquad (\tcf. \ 10.11).
\]
Then since the $\lambda_n$ are real and $> 0$ with 
$
\ds 
\sum\limits_{n = 1}^\infty \ 
\frac{1}{\lambda_n}
< \infty
$, 
we have
\[
f(z) 
\ = \ 
C \hsy z^m
\exp
\Big(
\Big(
a 
\hsx - \hsx 
\sum\limits_{n = 1}^\infty \ 
\frac{1}{\lambda_n}
\Big)
z
\Big)
\ 
\prod\limits_{n = 1}^\infty \ 
\Big(
1 - \frac{z}{\lambda_n}
\Big)
\hsy e^{z/\lambda_n}
\]
and 
$
\ 
\ds 
\sum\limits_{n = 1}^\infty \ 
\frac{1}{\lambda_n^2}
< \infty
$.
\\

\begin{x}{\small\bf DEFINITION} \ 
The \un{Laguerre-Polya} class of entire functions is comprised of the elements of 
$\ent(\hsx ]-\infty, +\infty[\hsx)$.
\\[-.25cm]
\end{x}

\begin{x}{\small\bf DEFINITION} \ 
The \un{type-\tI} Laguerre-Polya class of entire functions is comprised of the elements of 
\[
\ent(\hsx ]-\infty, 0]) \hsx \cup \hsx \ent(\hsx[0, +\infty[\hsx).
\]
\\[-1.5cm]
\end{x}

\begin{x}{\small\bf DEFINITION} \ 
The \un{type-\tII} Laguerre-Polya class of entire functions is comprised of the elements of 
$\ent(\hsx ]-\infty, +\infty[\hsx)$ 
which are not of type I.
\\[-.25cm]
\end{x}

\begin{x}{\small\bf NOTATION} \ 
$\gLP$, $\sI - \sL - \sP$, $\sI\sI - \sL - \sP$.
\\[-.25cm]
\end{x}


\begin{x}{\small\bf EXAMPLE} \ 
Let $p$ be a real polynomial with real zeros only.
\\[-.25cm]

\qquad \textbullet \quad
If all the nonzero zeros of $p$ are either positive or negative, then 
$p \in \sI - \sL - \sP$.
\\[-.5cm]

\qquad \textbullet \quad
If $p$ has both positive and negative zeros, the $p \in \sI\sI - \sL - \sP$.
\\[-.5cm]
\end{x}


\begin{x}{\small\bf EXAMPLE} \ 
The function
\[
\frac{1}{\Gamma (z)}
\ = \ 
z \hsy e^{\gamma \hsy z} \ 
\prod\limits_{n = 1}^\infty \ 
\Big(
1 + \frac{z}{n}
\Big)
\hsx 
\exp
\Big(
-\frac{z}{n}
\Big)
\]
is in 
$\tII - \sL - \sP$ 
(cf. 1.30).
\\[-.25cm]
\end{x}

Given $A \geq 0$ $(A < \infty)$, put
\[
S(A)
\ = \ 
\{z : \abs{\Img z} \leq A\}.
\]
\\[-1.5cm]

\begin{x}{\small\bf NOTATION} \ 
$A - \sL - \sP$ stands for the class of real entire functions $f \not\equiv 0$ 
that have a representation of the form

\[
f(z) 
\ = \ 
C \hsy z^m \hsy e^{a z^2 + b\hsy z} \ 
\prod\limits_{n = 1}^\infty \ 
\Big(
1 - \frac{z}{z_n}
\Big)
\hsx 
e^{z / z_n},
\]
where $C \neq 0$ is real, $m$ is a nonnegative integer, $a$ is real and $\leq 0$, $b$ is real, 
the $z_n \in S(A) - \{0\}$ with 
$
\ 
\ds
\sum\limits_{n = 1}^\infty \ 
\frac{1}{\abs{z_n}^2} 
\hsx < \hsx 
\infty
$.
\\[-.25cm]

[Note: \ 
Therefore
\[
0 - \sL - \sP
\ = \ 
\sL - \sP.]
\]
\\[-1.25cm]
\end{x}

\begin{x}{\small\bf THEOREM} \ 
$f \in A - \sL - \sP$ iff $f$ is the uniform limit on compact subsets of $\Cx$ of a sequence of real polynomials 
whose only zeros are in $S(A)$.
\\[-.25cm]
\end{x}

\begin{x}{\small\bf REMARK} \ 
Take $T \in S(A)$ $-$then 
\[
A - \sL - \sP \subset \ent (S(A)),
\]
the containment being proper if $A > 0$.
\\[-.5cm]

[Note: \ 
It is possible to characterize $\ent (S(A))$ but we shall omit the details as they will not be needed.]
\\[-.25cm]
\end{x}


\begin{x}{\small\bf EXAMPLE} \ 
The real polynomial $z (z^2 + 1)$ belongs to 
$1 - \sL - \sP$.
\\[-.25cm]
\end{x}

\begin{x}{\small\bf LEMMA} \ 
$A - \sL - \sP$ is closed under differentiation.
\\[-.5cm]

[This is because $S(A)$ is convex, so 8.3 is applicable.]
\\[-.25cm]
\end{x}

\begin{x}{\small\bf NOTATION} \ 
Denote by 
\\[-.25cm]
\[
* - \sL - \sP
\]
the class of real entire functions of the form 
\[
\phi(z) 
\ = \ 
p(z) \hsy f(z),
\]
where $p$ is a real polynomial and $f \in \sL - \sP$.
\\[-.25cm]
\end{x}

\begin{x}{\small\bf LEMMA} \ 
$\phi \in * - \sL - \sP$ iff $\phi \in A - \sL - \sP$ for some $A$ and $\phi$ has at most a finite number of 
nonreal zeros. 
\\[-.25cm]
\end{x}

\begin{x}{\small\bf LEMMA} \ 
$* - \sL - \sP$ is closed under differentiation.
\\[-.5cm]

PROOF \ 
Take a $\phi \in * - \sL - \sP$ and fix an $A: \phi \in A - \sL - \sP$  
$-$then 
$\phi^\prime \in A - \sL - \sP$ 
(cf. 10.35) and has at most a finite number of nonreal zeros (cf. 8.24).
\\[-.5cm]
\end{x}

Let $\phi \in * - \sL - \sP$ and suppose that 
$a \pm \sqrt{-1} \hsx b$ 
is a pair of conjugate nonreal zeros of $\phi$.  
\\[-.25cm]

\begin{x}{\small\bf DEFINITION} \ 
Given $k \geq 1$, the ellipse whose minor axis has 
$a + \sqrt{-1} \hsx b$ 
and
$a - \sqrt{-1} \hsx b$ 
as endpoints and whose major axis has length 
$2 \hsy b \hsy \sqrt{k}$
is called the \un{Jensen ellipse} of order $k$ of $\phi$.
\\[-.5cm]
\end{x}

The notion of ``Jensen ellipse'' generalizes that of ``Jensen circle'' 
(in the context of a real polynomial) 
and the proof of the following result is a computation similar to that used in 9.3.
\\[-.5cm]


\begin{x}{\small\bf THEOREM} \ 
Let 
$\phi \in * - \sL - \sP$
$-$then every nonreal zero of 
$\phi^{(k)}$ 
lies in the union of the Jensen ellipses of order $k$ of $\phi$.
\\[-.5cm]

[Note: \ 
Restated, if 
$a_j \pm \sqrt{-1} \hsx b_j$ $(j = 1, \ldots, d)$ 
are the nonreal zeros of $\phi$ and if 
$z = x + \sqrt{-1} \hsx y$ 
is a nonreal zero of 
$\phi^{(k)}$, 
then for some $j$,
\[
\frac{(x - a_j)^2}{k} + y^2
\ \leq \ 
b_j^2.]
\]
\\[-1.5cm]

The symbols $C$, $C^\prime$, $E^\prime$ 
employed in 8.24 make sense in the present setting 
(replace the ``$f$'' there by the ``$\phi$'' here).  
Therefore
\[
E^\prime  +  C^\prime
\ \leq \ 
C + \gen \phi 
\]
and
\[
\gen \phi 
\ = \ 
\gen \phi^\prime.
\]
\\[-.75cm]
\end{x}

\begin{x}{\small\bf LEMMA} \ 
Let 
$\phi \in * - \sL - \sP$ 
$-$then $C^\prime \leq C$ (cf. 8.22).
\end{x}


\chapter{
$\boldsymbol{\S}$\textbf{11}.\quad  DERIVATIVES}
\setlength\parindent{2em}
\setcounter{theoremn}{0}
\renewcommand{\thepage}{\S11-\arabic{page}}

\begin{x}{\small\bf DEFINITION} \ 
An entire function $\phi$ is said to be of \un{growth $(2, A)$} $(0 \leq A < \infty)$ 
if its order is $< 2$ or is of order 2 with type not exceeding $A$.
\\[-.25cm]
\end{x}

Denote by 
\[
\ent (2, A)
\]
the class of entire functions of growth $(2, A)$ $-$then 
\[A < A^\prime 
\ \implies \ 
\ent(2, A) \subset \ent(A, A^\prime).
\]
In particular:
\[
\ent(2, 0) \subset \ent(2, A).
\]
\\[-1.25cm]

\begin{x}{\small\bf LEMMA} \ 
The class $\ent(2, A)$ is closed under differentiation (cf. 2.25 and 3.7).
\\[-.25cm]
\end{x}

{\small\bf \un{N.B.}} \ 
If $\phi \in \ent(2, A)$, then for every $a > A$, 
\[
M(r; \phi) 
\ < \ 
e^{ar^2} 
\qquad (r \gg 0).
\]
\\[-1.25cm]

We shall now establish some technicalities that will be needed for the proof of the main result (viz. 11.9 infra).
\\[-.25cm]

\begin{x}{\small\bf NOTATION} \ 
Given positive real numbers $A > 0$, $B > 0$, let 
\[
C 
\ = \ 
\big(B + \sqrt{B^2 + 2 \hsy A^{-1}}\hsx\big) / 2,
\]
thus
\[
2\hsy A \hsy C (C - B) 
\ = \ 1.
\]
\\[-1.25cm]
\end{x}

\begin{x}{\small\bf LEMMA} \ 
If $\phi \in \ent(2, A)$, then
\[
\underset{n \ra \infty}{\limsupx} \ 
\sqrt{n} \ 
\bigg[
\frac{M\big(B \hsy \sqrt{n}; \phi^{(n)}\big)}{n !}
\bigg]^{1/n}
\ \leq \ 
2\hsy A \hsy C e^{A C^2}.
\]

PROOF \ 
Take $a > A$ and let 
\[
c 
\ = \ 
\big(
B + \sqrt{B^2 + 2 a^{-1}}
\hsx
\big) 
 / 2, 
\]
so that
\[
2 \hsy a \hsy c (c - B) 
\ = \ 1.
\]
Determine $r_0$: 
\[
r \geq r_0 
\implies
M(r ; \phi) < e^{a r^2}.
\]
Then for $n = 1, 2, \ldots$, 
\[
\log 
\bigg[
\frac{M\big(B \hsy \sqrt{n}; \phi^{(n)}\big)}{n !}
\bigg]^{1/n}
\ \leq \ 
\frac{a r^2}{n} - \log (r - B \hsy \sqrt{n}\hsy)
\]
if $r > \max (r_0, B \hsy \sqrt{n}\hsx )$.  
Since the RHS attains its minimum
\[
\log \frac{2 \hsy a \hsy c \hsy e^{a c^2}}{\sqrt{n}}
\]
at $r = c \hsy\sqrt{n}$, it follows that 
\[
\underset{n \ra \infty}{\limsupx} \ 
\sqrt{n} \ 
\bigg[
\frac{M\big(B \hsy \sqrt{n}; \phi^{(n)}\big)}{n !}
\bigg]^{1/n}
\ \leq \ 
2 \hsy a \hsy c \hsy e^{a c^2}.
\]
To finish, let $a \downarrow A$.
\\[-.25cm]
\end{x}

Let $f$ be an entire function and suppose that 
$z_0, z_1, \ldots$ is a sequence of complex numbers such that 
$\forall \ n \geq 0$, $f^{(n)} (z_n) = 0$ $-$then $\forall \ n > 0$, 
\[
f(z)
\ = \ 
\int\limits_{z_0}^z \ 
\int\limits_{z_1}^{\zeta_1} \ 
\cdots
\int\limits_{z_{n-1}}^{\zeta_{n-1}} \ 
f^{(n)} (\zeta_n) 
\ \td \zeta_n 
\cdots 
\td \zeta_2 \hsx
\td \zeta_1 \hsx.
\]
\\[-1.25cm]

\begin{x}{\small\bf SUBLEMMA} \ 
We have
\[
\abs{f(z)} 
\ \leq \ 
\frac{1}{n !} \ 
\sup\limits_{w \in H_n} \ 
\abs{f^{(n)} (w)} \ 
\big(
\abs{z - z_0} + \abs{z_0 - z_1}  + \cdots + \abs{z_{n-2} - z_{n-1}}  
\big)^n,
\]
where $H_n$ is the convex hull of the set 
$\{z, z_0, z_1, \ldots, z_{n-1}\}$.
\\[-.25cm]
\end{x}

\begin{x}{\small\bf SUBLEMMA} \ 
If $w \in H_n$, then 
\[
\abs{w} 
\ \leq \ 
\abs{z} + 
\abs{z - z_0} + \abs{z_0 - z_1}  + \cdots + \abs{z_{n-2} - z_{n-1}} . 
\]

PROOF \ 
Let $D_n$ be the closed disk of radius the RHS centered at the origin: 
$z \in D_n$.  
Next, 
\allowdisplaybreaks\begin{align*}
&
\abs{z_0} 
\leq 
\abs{z} + \abs{z - z_0} 
\implies z_0 \in D_n
\\[11pt]
&\abs{z_1} 
\leq 
\abs{z} + \abs{z - z_0} + \abs{z_0 - z_1} 
\implies z_1 \in D_n
&
\\[11pt]
&\hspace{1.5cm}
\vdots
\end{align*}
Therefore $D_n$ contains 
$z, z_0, z_1, \ldots, z_{n-1}$, 
hence being convex, $D_n$ contains $w$.
\\[-.25cm]

Accordingly, $H_n \subset D_n$, and 
\[
\abs{f(z)} 
\ \leq \ 
\frac{1}{n !} \ 
\sup\limits_{w \in D_n} \ 
\abs{f^{(n)} (w)} \ 
\big(
\abs{z - z_0} + \abs{z_0 - z_1}  + \cdots + \abs{z_{n-2} - z_{n-1}}  
\big)^n.
\]
\\[-1.25cm]
\end{x}

\begin{x}{\small\bf LEMMA} \ 
Maintaining the notation and assumptions of 11.4, 
suppose further that
\[
2 \hsy A \hsy B \hsy C \hsy e^{A C^2} 
\ < \ 
1.
\]
Impose the following conditions: \ 
$\exists$ a sequence 
$z_0, z_1, \ldots$ 
of complex numbers such that 
$\forall \ n \geq 0$, $\phi^{(n)} (z_n) = 0$ and 
\[
\underset{n \ra \infty}{\limsupx} \ 
\big(
\abs{z_0 - z_1} 
+ 
\abs{z_1 - z_2} 
+ \cdots + 
\abs{z_{n-1} - z_n} 
\big) 
/ \sqrt{n}
\ < \ 
B.
\]
Then
\[
\phi 
\ \equiv \ 
0.
\]
\\[-1.25cm]

PROOF \ 
In fact, 
\allowdisplaybreaks\begin{align*}
\underset{n \ra \infty}{\limsupx} \ 
B \hsy \sqrt{n} \ 
\bigg[
\frac{M\big(B \hsy \sqrt{n}; \phi^{(n)}\big)}{n !}
\bigg]^{1/n} \ 
&\leq \ 
2 \hsy A \hsy B \hsy C \hsy e^{A C^2} 
\qquad (\tcf. \ 11.4)
\\[15pt]
&<\ 
1
\end{align*}

\qquad $\implies$

\[
\lim\limits_{n \ra \infty} \ 
\frac{M\big(B \hsy \sqrt{n}; \phi^{(n)}\big)}{n !}
\hsx 
(B \hsy \sqrt{n}\hsx)^n
\ = \ 
0.
\]
Fix $z$ and determine $n_0$:
\[
n \geq n_0 
\implies
\abs{z} + 
\abs{z - z_0} + \abs{z_0 - z_1}  + \cdots + \abs{z_{n-2} - z_{n-1}}  
\ \leq \ 
B \hsy \sqrt{n},
\]
so $n \geq n_0$, 
\allowdisplaybreaks\begin{align*}
\implies \ 
\abs{\phi(z)} \ 
&\leq 
\frac{M\big(B \hsy \sqrt{n}; \phi^{(n)}\big)}{n !}
\hsx 
(B \hsy \sqrt{n}\hsx)^n
\qquad (\tcf. \ 11.5 \ \text{and} \ 11.6)
\\[11pt]
&\implies 
\abs{\phi(z)} = 0 
\\[11pt]
&\implies 
\phi(z) = 0.
\end{align*}
\\[-1.25cm]
\end{x}

\begin{x}{\small\bf SUBLEMMA} \ 
Let $\gamma_k = \alpha_k + \sqrt{-1} \hsx \beta_k$ $(\beta_k > 0)$ $(k = 0, 1, \ldots, n)$ be complex numbers such that

\[
\abs{\gamma_{k+1} - \alpha_k} 
\ \leq \ 
\beta_k
\qquad (k = 0, 1, \ldots, n - 1).
\]
Then 
\[
0 \leq \beta_n \leq \beta_{n - 1} \leq \cdots \leq \beta_0
\]
and 
\allowdisplaybreaks\begin{align*}
\abs{\gamma_0 - \gamma_1} 
&+ 
\abs{\gamma_1 - \gamma_2} 
+ \cdots + 
\abs{\gamma_{n-1} - \gamma_n} 
\leq \ 
\beta_0 - \beta_n
+
\sqrt{n}\hsx 
\big(
\beta_0^2 - \beta_n^2
\big)^{1/2}.
\end{align*}

PROOF \ 
The decrease of the $\beta_k$ is immediate and induction on $n$ leads to the inequality
\[
\abs{\alpha_0 - \alpha_1} 
+ 
\abs{\alpha_1 - \alpha_2} 
+ \cdots + 
\abs{\alpha_{n-1} - \alpha_n} 
\ \leq \ 
\sqrt{n}\hsx 
\big(
\beta_0^2 - \beta_n^2
\big)^{1/2}, 
\]
from which
\allowdisplaybreaks\begin{align*}
\abs{\gamma_0 - \gamma_1} 
&+ 
\abs{\gamma_1 - \gamma_2} 
+ \cdots + 
\abs{\gamma_{n-1} - \gamma_n} \ 
\\[11pt]
&\leq \ 
\abs{\alpha_0 - \alpha_1} 
+ 
\abs{\alpha_1 - \alpha_2} 
+ \cdots + 
\abs{\alpha_{n-1} - \alpha_n} 
\\[11pt]
&\hspace{1.5cm}
\hsx + \hsx 
(\beta_0 -  \beta_1) 
+ 
(\beta_1 -  \beta_2) 
+ \cdots + 
(\beta_{n-1} -  \beta_n) 
\\[11pt]
&\leq \ 
\sqrt{n}\hsx 
\big(
\beta_0^2 - \beta_n^2
\big)^{1/2}
 + 
\beta_0 - \beta_n.
\end{align*}

[Note: \ 
Extending the setup to infinity, let 
$
\ds
\beta 
= 
\lim\limits_{n \ra \infty} \ 
\beta_n$, 
hence
\[
\underset{n \ra \infty}{\limsupx} \ 
\big(
\abs{\gamma_0 - \gamma_1} 
+ 
\abs{\gamma_1 - \gamma_2} 
+ \cdots + 
\abs{\gamma_{n-1} - \gamma_n} 
\big)
/ 
\sqrt{n}
\ \leq \ 
\big(
\beta_0^2 - \beta^2
\big)^{1/2}.]
\]
\\[-1.25cm]
\end{x}

To see how data of this type is going to arise, 
take a 
$\phi \in * \hsy - \hsy \sL \hsy - \hsy \sP$ 
$-$then $\forall \ n \geq 0$,  
$\phi^{(n)} \in * \hsy - \hsy \sL \hsy - \hsy \sP$ 
(cf. 10.38) and given a nonreal zero $z_{n+1}$ of 
$\phi^{(n+1)}$ in the open upper half-plane, there is a nonreal zero of 
$z_{n}$ of 
$\phi^{(n)}$ in the open upper half-plane such that 
\[\abs{z_{n+1} - \Reg z_n} 
\ \leq \ 
\Img z_n.
\]
\\[-1.25cm]

[Note: \ 
This is a consequence of 10.40 (use Jensen circles, replacing the $\phi$ there by $\phi^{(n)}$ 
and then applying the theory to the pair 
$(\phi^{(n)},\phi^{(n+1)}))$.]
\\[-.25cm]


\begin{x}{\small\bf THEOREM} \ 
Let $\phi \in * - \sL - \sP$ $-$then there is a positive integer $N_0$ such that 
$\forall \ N \geq N_0$, $\phi^{(N)}$ has only real zeros, thus is in $\sL - \sP$. 
\\[-.25cm]
\end{x}

In order to utilize the machinery developed above, there is one crucial preliminary to be dealt with.
\\[-.25cm]

Let 
$\phi \in * - \sL - \sP$ and let $c_1, \bar{c}_1, \ldots, c_J, \bar{c}_J$ 
denote the nonreal zeros of $\phi$ $-$then $\phi$ has a representation of the form

\[
C \hsx 
\prod\limits_{j = 1}^J \ 
(z - c_j) \hsy (z - \bar{c}_j)
\hsy 
z^m 
\hsy
e^{az^2 + bz}
\ 
\prod\limits_{n = 1}^\infty \ 
\Big(1 - \frac{z}{\lambda_n}\Big) 
\hsy
e^{z/\lambda_n},
\]
where the various parameters are subject to the conditions enumerated in 10.19.
\\[-.25cm]

\begin{x}{\small\bf LEMMA} \ 
A given $\phi \in * - \sL - \sP$ is of growth $(2, \abs{a})$.
\\[-.25cm]

PROOF \ 
It is simply a matter of examining the various possibilities.
\\[-.25cm]

[Note: \ 
The polynomial
\[
C \hsx 
\prod\limits_{j = 1}^J \ 
(z - c_j) \hsy (z - \bar{c}_j)
\hsy 
z^m 
\]
can be safely ignored.]
\\[-.75cm]
\end{x}

\begin{spacing}{2}
\qquad 1. \
If 
$a = 0$, 
$b = 0$, 
and if the product 
$
\ds
\prod\limits_{n = 1}^\infty \ 
\Big(
1 - \frac{z}{\lambda_n} 
\Big)\hsx 
e^{z / \lambda_n}
$
is finite (recall the conventions set forth in 10.19), 
then the order of $\phi$ is 0.
\\[-.75cm]
\end{spacing}

\begin{spacing}{2}
\qquad 2. \
If 
$a = 0$, 
$b \neq 0$, 
and if the product 
$
\ds
\prod\limits_{n = 1}^\infty \ 
\Big(
1 - \frac{z}{\lambda_n} 
\Big)\hsx 
e^{z / \lambda_n}
$
is finite, 
then the order of $\phi$ is 1 (cf. 2.36).
\\[-1cm]
\end{spacing}

\begin{spacing}{2}
\qquad 3. \
If 
$a \neq 0$, 
$b = 0$ or $\neq 0$, 
and if the product 
$
\ds
\prod\limits_{n = 1}^\infty \ 
\Big(
1 - \frac{z}{\lambda_n} 
\Big)\hsx 
e^{z / \lambda_n}
$
is finite, 
then the order of $\phi$ is 2 and its type is $\abs{a}$ (cf. 3.2).
\\[-1cm]
\end{spacing}

\begin{spacing}{2}
\qquad 4. \
If 
$a = 0$, 
$b = 0$, 
and if the product 
$
\ds
\prod\limits_{n = 1}^\infty \ 
\Big(
1 - \frac{z}{\lambda_n} 
\Big)\hsx 
e^{z / \lambda_n}
$
is infinite, then there are two possibilities.
\\[-1cm]

\qquad \textbullet \quad
$
\ds
\sum\limits_{n = 1}^\infty \ 
\frac{1}{\lambda_n^2} 
\ < \ 
\infty
$
\ and \ 
$
\ds
\sum\limits_{n = 1}^\infty \ 
\frac{1}{\abs{\lambda_n}} 
\ = \ 
\infty
$
$-$then $\fg = 1$ is the genus of the sequence 
$\{\abs{\lambda_n} : n = 1, 2, \ldots\}$ 
(cf. 4.14), 
hence
$
\ds
\prod\limits_{n = 1}^\infty \ 
\Big(
1 - \frac{z}{\lambda_n} 
\Big)\hsx 
e^{z / \lambda_n}
$
is the associated canonical product (cf. 5.9).  
As such, its order is $\kappa$ (the convergence exponent of the sequence
$\{\abs{\lambda_n} : n = 1, 2, \ldots\}$) (cf. 5.10).  
But 
$1 \leq \kappa \leq 1 + 1$.  
(cf. 4.15), so the order of the product
$
\ds
\prod\limits_{n = 1}^\infty \ 
\Big(
1 - \frac{z}{\lambda_n} 
\Big)\hsx 
e^{z / \lambda_n}
$
is $\leq 2$.  
It remains to analyze the situation when $\kappa = 2$.  
This, however, is immediate: 
$
\ 
\ds
\prod\limits_{n = 1}^\infty \ 
\Big(
1 - \frac{z}{\lambda_n} 
\Big)\hsx 
e^{z / \lambda_n}
\ 
$
is of minimal type (cf. 5.16), thus is of growth $(2, 0)$ or still, 
is of growth $(2, \abs{a})$ (since here $a = 0$). 
\\[-.5cm]

\qquad \textbullet \quad
$
\ds
\sum\limits_{n = 1}^\infty \ 
\frac{1}{\lambda_n^2} 
\ < \ 
\infty
$
\ and \ 
$
\ds
\sum\limits_{n = 1}^\infty \ 
\frac{1}{\abs{\lambda_n}} 
\ < \ 
\infty
$
$-$then $\fg = 0$ is the genus of the sequence 
$\{\abs{\lambda_n} : n = 1, 2, \ldots\}$ 
(cf. 4.14) and we can write
\[
\prod\limits_{n = 1}^\infty \ 
\Big(
1 - \frac{z}{\lambda_n} 
\Big)\hsx 
e^{z / \lambda_n}
\ = \ 
\exp
\bigg( 
\Big( 
\sum\limits_{n = 1}^\infty \ 
\frac{1}{\lambda_n} 
\Big)
z
\bigg)
\prod\limits_{n = 1}^\infty \ 
\Big(
1 - \frac{z}{\lambda_n} 
\Big).
\]
Thanks to 5.11, the order of the RHS is 
$\max (1, \kappa) \leq \max (1, 1) = 1$
if 
$
\ 
\ds
\sum\limits_{n = 1}^\infty \ 
\frac{1}{\lambda_n} 
\ \neq \ 
0
\ 
$
or $\kappa \leq 1$ if 
$
\ 
\ds
\sum\limits_{n = 1}^\infty \ 
\frac{1}{\lambda_n} 
\ = \ 
0.
$
\\[-1cm]
\end{spacing}

\begin{spacing}{2}
\qquad 5. \quad
If 
$a = 0$, 
$b \neq 0$, 
and if the product 
$
\ds
\prod\limits_{n = 1}^\infty \ 
\Big(
1 - \frac{z}{\lambda_n} 
\Big)\hsx 
e^{z / \lambda_n}
$
is infinite, then there are two possibilities.
\\[-.75cm]


\qquad \textbullet \quad
$
\ds
\sum\limits_{n = 1}^\infty \ 
\frac{1}{\lambda_n^2} 
\ < \ 
\infty
$
\ and \ 
$
\ 
\ds
\sum\limits_{n = 1}^\infty \ 
\frac{1}{\abs{\lambda_n}} 
\ = \ 
\infty .
\ 
$
Suppose first that $\kappa$ is $< 2$
$-$then the order of
\[
e^{b\hsy z} \ 
\prod\limits_{n = 1}^\infty \ 
\Big(
1 - \frac{z}{\lambda_n} 
\Big)
\hsx e^{z / \lambda_n}
\]
is $\max (1, \kappa) < 2$ (cf. 5.11).  
On the other hand, if $\kappa = 2$, then the order of 
\[
e^{b\hsy z} \ 
\prod\limits_{n = 1}^\infty \ 
\Big(
1 - \frac{z}{\lambda_n} 
\Big)
\hsx e^{z / \lambda_n}
\]
is $\max (1, 2) = 2$ (cf. 5.11).  
As for its type, use 3.14 in the ``$\rho_1 < \rho_2$'' scenario to see that it is minimal, thus
\[
e^{b\hsy z} \ 
\prod\limits_{n = 1}^\infty \ 
\Big(
1 - \frac{z}{\lambda_n} 
\Big)
\hsx e^{z / \lambda_n}
\]
is of growth $(2, 0)$ or still, is of growth $(2, \abs{a})$ (since here $a = 0$).
\\[-.5cm]

\qquad \textbullet \quad
$
\ds
\sum\limits_{n = 1}^\infty \ 
\frac{1}{\lambda_n^2} 
\ < \ 
\infty
$
\ and \ 
$
\ds
\sum\limits_{n = 1}^\infty \ 
\frac{1}{\abs{\lambda_n}} 
\ < \ 
\infty
$
$-$then the order of the product
$
\ds
\prod\limits_{n = 1}^\infty \ 
\Big(
1 - \frac{z}{\lambda_n} 
\Big)
$
is $\leq 1$, 
hence the order of 
\[
e^{b\hsy z} \ 
\prod\limits_{n = 1}^\infty \ 
\Big(
1 - \frac{z}{\lambda_n} 
\Big)\hsx 
e^{z / \lambda_n}
\ = \ 
\exp
\Big(
\Big(
b + 
\sum\limits_{n = 1}^\infty \ 
\frac{1}{\lambda_n} 
\Big)
z
\Big)
\prod\limits_{n = 1}^\infty \ 
\Big(
1 - \frac{z}{\lambda_n} 
\Big)
\]
is $\leq 1$ (cf. 5.11).
\\[-1cm]
\end{spacing}

\begin{spacing}{2}
\qquad 6. \quad
If 
$a \neq 0$, 
$b = 0$, or $\neq 0$,
and if the product 
$
\ds
\prod\limits_{n = 1}^\infty \ 
\Big(
1 - 
\frac{z}{\lambda_n} 
\Big)
e^{z / \lambda_n}
$
is infinite, 
then there are two possibilities.  
\\[-.5cm]

\qquad \textbullet \quad
$
\ds
\sum\limits_{n = 1}^\infty \ 
\frac{1}{\lambda_n^2} 
\ < \ 
\infty
$
\ and \ 
$
\ 
\ds
\sum\limits_{n = 1}^\infty \ 
\frac{1}{\abs{\lambda_n}} 
\ = \ 
\infty.
\ 
$
Suppose first that $\kappa$ is $< 2$
$-$then the order of
\[
e^{a z^2 + bz} \ 
\prod\limits_{n = 1}^\infty \ 
\Big(
1 - \frac{z}{\lambda_n} 
\Big)
e^{z / \lambda_n}
\]
is $\max (2, \kappa) = 2$ (cf. 5.11) and its type is $\abs{a}$ (apply 3.14 (first bullet point)).  
As for what happens when $\kappa = 2$, the product
$
\ds
\prod\limits_{n = 1}^\infty \ 
\Big(
1 - \frac{z}{\lambda_n} 
\Big)\hsx 
e^{z / \lambda_n}
$ 
is of minimal type (see above), 
so another appeal to 3.14 (second bullet point) allows one to conclude that the type of 
\[
e^{a z^2 + bz} \ 
\prod\limits_{n = 1}^\infty \ 
\Big(
1 - \frac{z}{\lambda_n} 
\Big)
e^{z / \lambda_n}
\]
is again $\abs{a}$.
\\[-.5cm]

\qquad \textbullet \quad
$
\ds
\sum\limits_{n = 1}^\infty \ 
\frac{1}{\lambda_n^2} 
\ < \ 
\infty
$
\ and \ 
$
\ds
\sum\limits_{n = 1}^\infty \ 
\frac{1}{\abs{\lambda_n}} 
\ < \ 
\infty
$
$-$then the order of the product
$
\ds
\prod\limits_{n = 1}^\infty \ 
\Big(
1 - \frac{z}{\lambda_n} 
\Big)
$
is $\leq 1$, 
hence the order of 

\[
e^{a z^2 + bz} \ 
\prod\limits_{n = 1}^\infty \ 
\Big(
1 - \frac{z}{\lambda_n} 
\Big)\hsx 
e^{z / \lambda_n}
\ = \ 
\exp
\Big(
a z^2 + 
\Big(
b + 
\sum\limits_{n = 1}^\infty \ 
\frac{1}{\lambda_n} 
\Big)
z
\Big)
\prod\limits_{n = 1}^\infty \ 
\Big(
1 - \frac{z}{\lambda_n} 
\Big)
\]
is 2 (cf. 5.11) and its type is $\abs{a}$ 
(use 3.14 in the ``$\rho_1 < \rho_2$'' scenario).
\end{spacing}

Passing now to the proof of 11.9, it suffices to show that there is a positive $N_0$ such that 
$\phi^{(N_0)}$ has only real zeros (cf. 10.38 and 10.41).  
Proceeding by contradiction,
suppose that $\forall \ n \geq 0$, 
$\phi^{(n)}$ 
has a nonreal zero and let $X_n$ denote the set of nonreal zeros of 
$\phi^{(n)}$  
in the open upper half-plane $\Img z > 0$ $-$then each $X_n$ is finite and the product
$
\ds
X = \prod\limits_{n = 0}^\infty \ X_n
$
is a nonempty compact set.  
Given $n = 1, 2, \ldots,$ put
\[
E_n 
\ = \ 
\Big\{
\big(\zeta_0, \zeta_1, \ldots\big) 
\in X: 
\abs{\zeta_{j + 1} - \Reg \zeta_j} 
\hsx \leq \hsx 
\Img \zeta_j, 
\hsx 
j = 0, 1, \ldots, n
\Big\}.
\]
Then $E_n$ is a closed subset of $X$ and $E_1 \supset E_2 \supset \cdots$.  
Furthermore, $E_n$ is nonempty, so 
$
\
\ds
\bigcap\limits_{n = 1}^\infty \ 
E_n 
\ \neq \ 
\emptyset
$, 
thus one can find a sequence $z_0, z_1, \ldots$ of complex numbers such that 
\[
\Img z_n > 0, 
\quad
\phi^{(n)}(z_n) = 0,
\quad
\abs{z_{n + 1} - \Reg z_n} \leq \Img z_n.
\]
Write
$z_n = a_n + \sqrt{-1} \hsx b_n$ $(b_n > 0)$ $-$then $\{b_n\}$ is a decreasing sequence and 
\[
\abs{z_m - z_{m + 1}} 
\hsx + \hsx
\abs{z_{m + 1} - z_{m + 2}} 
\hsx + \hsx
\cdots
\hsx + \hsx
\abs{z_{m + n -  1} - z_{m + n}} 
\ \leq \ 
b_m - b_{m + n} 
+ 
\sqrt{n} \hsx
\big(
b_m^2 - b_{m + n}^2
\big)^{1/2}.
\]
Here $m = 0, 1, \ldots$ and $n = 1, 2, \ldots \hsx .$ 
Therefore
\[
\underset{n \ra \infty}{\limsupx} \ 
\frac
{
\abs{z_m - z_{m + 1}} 
\hsx + \hsx
\abs{z_{m + 1} - z_{m + 2}} 
\hsx + \hsx
\cdots
\hsx + \hsx
\abs{z_{m + n -  1} - z_{m + n}} 
}
{
\sqrt{n}
}
\ \leq \ 
\big(
b_m^2 - b^2
\big)^{1/2},
\]
where we have set 
$
\ds 
b = \lim\limits_{n \ra \infty} \ b_n.
$  
Fix $A > \abs{a}$, hence

{\ }
\\[-1.5cm]
\[
\phi \in \ent(2, A) 
\qquad (\tcf. \ 11.10).
\]
Choose $B > 0$: 
\[
2 \hsy A \hsy B \hsy C \hsy e^{A c^2} 
\ < \ 
1
\]
and choose $m$: 
\[
\big(
b_m^2 - b^2
\big)^{1/2} 
\ < \ 
B. 
\]
Then
\[
\underset{n \ra \infty}{\limsupx} \ 
\frac{
\abs{z_m - z_{m + 1}} 
\hsx + \hsx
\abs{z_{m + 1} - z_{m + 2}} 
\hsx + \hsx
\cdots
\hsx + \hsx
\abs{z_{m + n -  1} - z_{m + n}} 
}
{
\sqrt{n} 
}
\ < \ 
B.
\]
But
\[
\phi \in \ent(2, A) 
\implies 
\phi^{(m)} \in \ent(2, A) 
\qquad (\tcf. \ 11.2).
\]
And this means that 11.7 is applicable to $\phi^{(m)}$: 

\[
\implies \ 
\phi^{(m)} 
\ \equiv \ 
0.
\]
Contradiction $\ldots$ \hsy.
\\[-.25cm]
\begin{x}{\small\bf EXAMPLE} \ 
The real entire function 
$
e^{z^2}
$
belongs to $\ent (2, 1)$.  
However, it is not in $* - \sL - \sP$ and 11.9 does not obtain.
\end{x}


\chapter{
$\boldsymbol{\S}$\textbf{12}.\quad  JENSEN POLYNOMIALS}
\setlength\parindent{2em}
\setcounter{theoremn}{0}
\renewcommand{\thepage}{\S12-\arabic{page}}

\qquad 
Given a real entire function 
\[
f (z) 
\ = \ 
\sum\limits_{n = 0}^\infty \ 
c_n \hsy z^n, 
\]
put 
$\gamma_n = f^{(n)} (0)$, 
thus
\[
f (z) 
\ = \ 
\sum\limits_{n = 0}^\infty \ 
\frac{\gamma_n}{n !}  \hsy z^n.
\]
\\[-.75cm]

\begin{x}{\small\bf DEFINITION} \ 
The 
\un{$n^\nth$ Jensen polynomial}
$\tJ_n$ associated with $f$ is defined by
\[
\tJ_n (f; z) 
\ = \ 
\sum\limits_{k = 0}^n \ 
\binom{n}{k} \hsy \gamma_k \hsy z^k.
\]
\\[-1.cm]
\end{x}

\begin{x}{\small\bf LEMMA} \ 
The sequence $\{\tJ_n (f; t)\}$ is generated by $e^x \hsy f (x \hsy t)$, i.e., 
\[
e^x \hsy f (x \hsy t)
\ = \ 
\sum\limits_{n = 0}^\infty \ 
\tJ_n (f; t) \hsy \frac{x^n}{n !}
\qquad (x, t \in \R).
\]
\\[-1.25cm]
\end{x}

\begin{x}{\small\bf LEMMA} \ 
We have 
\[
z \hsy \tJ_n^\prime  (f; z)  
\ = \ 
n \hsy \tJ_n (f; z) - n \hsy \tJ_{n - 1} (f; z) 
\qquad (n \geq 1).
\]
\\[-1.5cm]
\end{x}

\begin{x}{\small\bf DEFINITION} \ 
The 
\un{$n^\nth$ Appell polynomial}
$\tJ_n^*$ associated with $f$ is defined by
\[
\tJ_n^* (f; z) 
\ = \ 
\sum\limits_{k = 0}^n \ 
\binom{n}{k} \hsy \gamma_k \hsy z^{n - k}.
\]
\\[-1.25cm]
\end{x}

\begin{x}{\small\bf LEMMA} \ 
The sequence $\{\tJ_n^* (f; t)\}$ is generated by $e^{x \hsy t} \hsy f (x)$, i.e., 
\[
e^{x \hsy t} \hsy f (x)
\ = \ 
\sum\limits_{n = 0}^\infty \ 
\tJ_n^* (f; t) \hsy \frac{x^n}{n !}
\qquad (x, t \in \R).
\]
\\[-1.25cm]
\end{x}


\begin{x}{\small\bf LEMMA} \ 
We have 
\[
\frac{\td}{\td z} \tJ_n^* (f; z)
\ = \ 
n \hsy \tJ_{n - 1}^* (f; z)
\qquad (n \geq 1).
\]
\\[-1.25cm]
\end{x}

\qquad
{\small\bf \un{N.B.}} \ 
Obviously, 
\[
\begin{cases}
\ \ds
\tJ_n (f; z) 
\ = \ 
z^n \hsy 
\tJ_n^* \Big(f; \frac{1}{z}\Big)
\\[18pt]
\ \ds
\tJ_n^* (f; z)
\ = \ 
z^n \hsy 
\tJ_n \Big(f; \frac{1}{z}\Big)
\end{cases}
.
\]
Therefore the zeros of $\tJ_n$ are real iff the zeros of $\tJ_n^*$ are real.
\\[-.25cm]

\begin{x}{\small\bf DEFINITION} \ 
The 
\un{$(n, m)^\nth$ Jensen polynomial} 
associated with $f$ is defined by
\[
\tJ_{n, m} (f; z) 
\ = \ 
\sum\limits_{k = 0}^n \ 
\binom{n}{k} \hsy \gamma_{k+m} \hsy z^k.
\]
\\[-1.25cm]
\end{x}

\qquad
{\small\bf \un{N.B.}} \ 
Therefore
\[
\tJ_{n, m} (f; z) 
\ = \ 
\tJ_n (f^{(m)}; z) .
\]
\\[-1.25cm]

\begin{x}{\small\bf LEMMA} \ 
We have 
\allowdisplaybreaks\begin{align*}
\tJ_n^{(m)} (f; z) \ 
&=\ 
\frac{n !}{(n - m)!} \  
\tJ_{n - m, m} (f; z)
\\[11pt]
&=\ 
\frac{n !}{(n - m)!} \  
\tJ_{n - m} (f^{(m)}; z).
\end{align*}
\\[-1cm]
\end{x}

\begin{x}{\small\bf THEOREM} \ 
On compact subsets of $\Cx$, 
\[
\tJ_n \Big(f; \frac{z}{n}\Big) \ra f (z)
\]
uniformly.
\\[-.5cm]

PROOF \ 
Fix a compact set $K \subset \Cx$.  
Given $\varepsilon > 0$, choose $N > 2$: 
\[
\sum\limits_{n = N + 1}^\infty \ 
\abs
{
\frac{\gamma_n}{n !} \hsx z^n 
}
\ < \ 
\frac{\varepsilon}{4}
\quad (z \in K).
\]
Next, choose $N^\prime > N$: 
\[
n \geq N^\prime 
\implies
\bigg| \hsy
\sum\limits_{k = 2}^N \ 
\Big(
\frac{\gamma_k}{k !} 
- 
\Big(
1 - \frac{1}{n}
\Big)
\cdots
\Big(
1 - \frac{k-1}{n}
\Big)
\frac{\gamma_k}{k !}
\Big)
z^k
\hsx
\bigg|
\ < \ 
\frac{\varepsilon}{2}
\quad (z \in K).
\]
Then $\forall \ z \in K$ and $\forall \ n \geq N^\prime$ : 
\allowdisplaybreaks
\allowdisplaybreaks\begin{align*}
&\abs{f (z) - \tJ_n \Big(f; \frac{z}{n}\Big)}
\\[15pt]
&=\ 
\bigg|
\sum\limits_{n = N + 1}^\infty \ 
\frac{\gamma_n}{n !} \hsy z^n 
\hsx + \hsx 
\sum\limits_{k = 0}^N \ 
\frac{\gamma_k}{k !} \hsy z^k 
\hsy - \hsy 
\Big(
\gamma_0 + \gamma_1 z 
+ 
\sum\limits_{k = 2}^N \ 
\Big(
1 - \frac{1}{n}
\Big)
\cdots
\Big(
1 - \frac{k-1}{n}
\Big)
\frac{\gamma_k}{k !} \hsy z^k 
\Big)
\bigg|
\\[15pt]
&=\ 
\bigg|
\sum\limits_{n = N + 1}^\infty \ 
\frac{\gamma_n}{n !} \hsy z^n 
\hsy + \hsy
\sum\limits_{k = 2}^N \ 
\frac{\gamma_k}{k !} \hsy z^k 
\hsy - \hsy 
\sum\limits_{k = 2}^n \ 
\Big(
1 - \frac{1}{n}
\Big)
\cdots
\Big(
1 - \frac{k-1}{n}
\Big)
\frac{\gamma_k}{k !} \hsy z^k 
\bigg|
\\[15pt]
&=\ 
\bigg|
\sum\limits_{n = N + 1}^\infty \ 
\frac{\gamma_n}{n !} \hsy z^n 
\hsy + \hsy
\sum\limits_{k = 2}^N \ 
\Big(
\frac{\gamma_k}{k !} 
- 
\Big(
1 - \frac{1}{n}
\Big)
\cdots
\Big(
1 - \frac{k-1}{n}
\Big)
\frac{\gamma_k}{k !}
\Big)
z^k
\\[15pt]
&
\hspace{3.5cm}  
- \ 
\sum\limits_{k = N + 1}^n \ 
\Big(
1 - \frac{1}{n}
\Big)
\cdots
\Big(
1 - \frac{k-1}{n}
\Big)
\frac{\gamma_k}{k !} \hsy z^k 
\bigg|
\\[15pt]
&\leq\ 
\sum\limits_{n = N + 1}^\infty \ 
\abs{\frac{\gamma_n}{n !} \hsy z^n} 
\hsy + \hsy
\bigg|
\sum\limits_{k = 2}^N \ 
\Big(
\frac{\gamma_k}{k !} 
- 
\Big(
1 - \frac{1}{n}
\Big)
\cdots
\Big(
1 - \frac{k-1}{n}
\Big)
\frac{\gamma_k}{k !}
\hsx
z^k
\bigg|
\\[15pt]
&
\hspace{3.5cm} 
+ \ 
\sum\limits_{k = N + 1}^n \ 
\bigg|
\Big(
1 - \frac{1}{n}
\Big)
\cdots
\Big(
1 - \frac{k-1}{n}
\Big)
\frac{\gamma_k}{k !} \hsy z^k 
\bigg|
\\[11pt]
&<\ 
\frac{\varepsilon}{4} + \frac{\varepsilon}{2} + \frac{\varepsilon}{4}
\\[11pt]
&=\ 
\varepsilon.
\end{align*}
\\[-1.25cm]
\end{x}

In what follows, certain classical facts from the theory of equations will be admitted without proof.  
To begin with: 
\\[-.25cm]

\begin{x}{\small\bf HERMITE-POULAIN CRITERION} \ 
Suppose that the real polynomial 
\[
a_0 + a_1 \hsy z + \cdots + a_n \hsy z^n
\]
has real zeros only.  
Let $p(z)$ be a real polynomial $-$then the polynomial
\[
P(z) 
\ = \ 
a_0 \hsy p(z) + a_1 \hsy p^\prime(z) + \cdots + a_n \hsy p^{(n)}(z)
\]
has at least as many real zeros as $p(z)$ does.
\\[-.5cm]

[Note: \ 
By taking limits, one can extend 12.10, viz. replace the real polynomial
\[
a_0 + a_1 \hsy z + \cdots + a_n \hsy z^n
\]
by an element $f \in \sL - \sP$ $-$then for any real polynomial $p(z)$, the polynomial 
\[
\sum\limits_{k = 0}^d \ 
\frac{f^{(k)} (0)}{k !} \hsy p^{(k)} (z)
\qquad (d = \deg p)
\]
has at least as many real zeros as $p(z)$ does.]
\\[-.25cm]
\end{x}

\begin{x}{\small\bf APPLICATION} \ 
A real polynomial has real zeros only iff its Jensen polynomials have real zeros only.
\\[-.5cm]

[Suppose that 
\[
f (z) 
\ = \ 
\gamma_0 + \frac{\gamma_1}{1 !} \hsy z + \cdots + \frac{\gamma_d}{d !} \hsy z^d
\]
is a real polynomial of degree $d$. 
\\[-.25cm]

\qquad \textbullet \quad
If $f(z)$ has real zeros only, take $p (z) = z^n$ in 12.10 to see that $\forall \ n = 1, 2, \ldots$,
\[
\tJ_n^* (f; z) 
\ = \ 
\gamma_0 \hsy z^n + \binom{n}{1} \hsy \gamma_1 \hsy z^{n - 1} + \cdots
\]
has real zeros only, so the same is true of $\tJ_n (f; z)$.
\\[-.25cm]

\qquad \textbullet \quad
If $\forall \ n = 1, 2, \ldots$, $\tJ_n (f; z)$ has real zeros only, then 
\[
f (z) 
\ = \ 
\lim\limits_{n \ra \infty} \ 
\tJ_n \Big(f; \frac{z}{n}\Big)
\]
has  real zeros only (cf. 12.9).]
\\[-.25cm]
\end{x}

\begin{x}{\small\bf MALO-SCHUR CRITERION} \ 
Suppose that the zeros of
\[
a_0 + a_1 \hsy z + \cdots + a_n \hsy z^n
\]
are real and the zeros of
\[
b_0 + b_1 \hsy z + \cdots + b_m \hsy z^m
\]
are real and of the same sign.  
Put $k = \min (n, m)$ $-$then the zeros of 
\[
a_0 \hsy b_0 + 1! \hsy a_1 \hsy b_1 \hsy z + \cdots + k ! \hsy a_k \hsy b_k \hsy z^k 
\]
are real.
\\[-.25cm]
\end{x}

\begin{x}{\small\bf EXAMPLE} \ 
Suppose that the zeros of
\[
a_0 + a_1 \hsy z + \cdots + a_n \hsy z^n
\]
are real $-$then the zeros of 
\[
a_n + a_{n - 1} \hsy z + \cdots + a_0 \hsy z^n
\]
are real.  
Working now with 
\[
(1 + z)^n 
\ = \ 
1 + \binom{n}{1} \hsy z + \cdots + z^n, 
\]
it follows that the zeros of 
\[
a_n + n \hsy a_{n -1} \hsy z + \cdots + n! \hsy a_0 \hsy z^n
\]
are real, or still, that the zeros of 
\\
\[
\frac{a_n}{n !} + \frac{a_{n - 1}}{(n - 1)!} \hsy z + \cdots + a_0 \hsy z^n
\]
\\[-.25cm]
are real, or still, that the zeros of 
\[
a_0 + \frac{a_1}{1 !} z+ \cdots + \frac{a_n}{n !} \hsy z^n
\]
are real.  
Consequently, if the zeros of 
\[
b_0 + b_1 \hsy z + \cdots + b_m \hsy z^m
\]
are real and of the same sign, then the zeros of
\[
a_0 \hsy b_0 + a_1 \hsy b_1 \hsy z + \cdots + a_k \hsy b_k \hsy z^k 
\qquad (k = \min(n, m))
\]
are real.
\\[-.25cm]
\end{x}

\begin{x}{\small\bf THEOREM} \ 
Let $f \not\equiv 0$ be a real entire function.  
$-$then $f \in \sL - \sP$ iff its Jensen polynomials have real zeros only.
\\[-.5cm]

PROOF\ 
In view of 12.9, it is clear that the condition is sufficient.  
Turning to the necessity, given that $f \in \sL - \sP$, choose a sequence 
$\{p_k : k = 1, 2, \ldots\}$ of real polynomials having real zeros only such that 
$p_k \ra f$ uniformly on compact subsets of $\Cx$, say
\[
p_k (z) 
\ = \ 
\gamma_{k \hsy 0} + \frac{\gamma_{k \hsy 1}}{1 !} + \cdots \hsy .
\]
Then the Jensen polynomials $\tJ_n (p_k ; z)$ have real zeros only (cf. 12.11).  
But for fixed $n$, 
\[
\lim\limits_{k \ra \infty} \ 
\tJ_n (p_k ; z) 
\ = \ 
\tJ_n (f ; z)
\]
uniformly on compact subsets of $\Cx$.
\\[-.25cm]
\end{x}


\begin{x}{\small\bf REMARK} \ 
If $f \in \sL - \sP$, then
\[
\tJ_n \Big(f; \frac{z}{n}\Big) \ra f (z) 
\]
uniformly on compact subsets of $\Cx$ and the zeros of 
$
\ds
\tJ_n \Big(f; \frac{z}{n}\Big)
$
are real.  
By comparison, the partial sums
\[
\sum\limits_{k = 0}^n \ 
\frac{\gamma_k}{k !} \hsy z^k, 
\]
while uniformly convergent on compact subsets of $\Cx$, may very well have nonreal zeros.  

\noindent
E.g.: Take $f (z) = e^z$ $-$then
\[
\sum\limits_{k = 0}^n \ 
\frac{z^k}{k !}
\]
has no real zeros if $n$ is even and has one real zero if $n$ is odd.
\\[-.25cm]
\end{x}

\begin{x}{\small\bf DEFINITION} \ 
A sequence 
$\gamma_0, \gamma_1, \ldots$ 
of real numbers is said to be a 
\un{multiplier sequence} 
if $\forall \ n = 1, 2, \ldots$, 
the real polynomial
\[
\sum\limits_{k = 0}^n \ 
\binom{n}{k} \hsy \gamma_k \hsy z^k
\]
has real zeros only or, equivalently, if $\forall \ n = 1, 2, \ldots$, the real polynomial 
\[
\sum\limits_{k = 0}^n \ 
\binom{n}{k} \hsy \gamma_k \hsy z^{n-k}
\]
has real zeros only. 
\\[-.25cm]
\end{x}

If $f \in \sL - \sP$, then the associated sequence 
$\gamma_0, \gamma_1, \ldots$ 
is a multiplier sequence (cf. 12.14).
\\[-.25cm]

\begin{x}{\small\bf EXAMPLE} \ 
Take 
\[
f (z) \ = \ 
\begin{cases}
\ 
e^z
\\[4pt]
\ 
e^{-z}
\end{cases}
\]
to see that
\[
\begin{cases}
\ 
1, \hsx 1, \hsx 1, \ldots
\\[4pt]
\ 
1, -1, \hsx1, \ldots
\end{cases}
\]
are multiplier sequences.
\\[-.25cm]
\end{x}

\begin{x}{\small\bf EXAMPLE} \ 
Let $p$ be a positive integer and take $f (z) = z^p \hsy e^z$ $-$then
\[
z^p \hsy e^z 
\ = \ 
p ! \ \frac{z^p}{p !} 
\ + \  
\frac{(p + 1)!}{1 !} \hsy \frac{z^{p + 1}}{(p + 1)!} 
\ + \  
\cdots \hsy .
\]
Therefore the sequence 
\[
0 , \hsx 0 , \ldots, 0,  \hsx p !, \hsx \frac{(p + 1)!}{1 !}, \dots
\]
is a multiplier sequence.
\\[-.5cm]

[Note: \ 
Specialize and let $p = 1$, thus $0, 1, 2, \ldots$ is a multiplier sequence.]
\\[-.25cm]
\end{x}

\begin{x}{\small\bf EXAMPLE} \ 
Take $f (z) = e^{-z^2/2}$ $-$then
\[
e^{-z^2/2}
\ = \ 
1 
\ - \  
\frac{z^2}{2 !} 
\ + \  
1 \cdot 3 \hsy \frac{z^4}{4 !} 
\ - \  
1 \cdot 3 \cdot 5 \hsy \frac{z^6}{6 !} 
\ + \  
\cdots .
\]
Therefore the sequence 
\[
1, \ 0,\  -1, \ 0, \ 1 \cdot 3, \ 0, \ -1 \cdot 3 \cdot 5, \ 0, \ldots
\]
is a multiplier sequence.
\\[-.25cm]
\end{x}

\begin{x}{\small\bf EXAMPLE} \ 
Take
\[
f (z) \ = \ 
\begin{cases}
\ 
\cos z
\\[4pt]
\ 
\sin z
\end{cases}
\]
$-$then
\[
\begin{cases}
\ 
1, \hsx 0, -1,  \hsx 0, 1, \hsx 0, -1, \ldots
\\[4pt]
\ 
0, \hsx 1, \hsx 0, -1, \hsx 0, \hsx 1, \hsx 0, \ldots
\end{cases}
\]
are multiplier sequences.
\\[-.25cm]
\end{x}

\begin{x}{\small\bf THEOREM} \ 
Let 
$\gamma_0, \gamma_1, \ldots$ 
be a multiplier sequence and put 
$
\ds
c_n  = \frac{\gamma_n}{n !}
$
$-$then 
\[
f (z) 
\ = \ 
\sum\limits_{n = 0}^\infty \ 
c_n \hsy z^n
\]
is a real entire function and, as such, is in $\sL - \sP$.
\\[-.5cm]

PROOF \ 
The objective is to find an estimate for $\abs{c_n}$ that suffices to ensure the covergence of the series at every $z$.  
This said, let $\gamma_r$ be the first nonzero entry in the sequence 
$\gamma_0, \gamma_1, \ldots \hsx .$ 
Take $n > r$: 
\allowdisplaybreaks\begin{align*}
\sum\limits_{k = 0}^n \ 
\binom{n}{k} \hsy \gamma_k \hsy z^{n-k} \ 
&=\ 
\sum\limits_{k = 0}^n \ 
\frac{n !}{(n - k)!} \hsy \frac{\gamma_k}{k !} \hsy z^{n-k}
\\[15pt]
&=\ 
\sum\limits_{k = 0}^n \ 
\frac{n !}{(n - k)!} \hsy c_k \hsy z^{n-k}
\\[15pt]
&=\ 
c_0 \hsy z^n + n \hsy c_1 \hsy z^{n-1} + \cdots + n! \hsy c_n
\\[15pt]
&=\ 
n (n - 1) \cdots (n - r + 1) c_r z^{n-1} + \cdots + n! \hsy c_n
\end{align*}
and denote by 
$
\lambda_1, \lambda_2, \ldots, \lambda_{n - r} 
$
its (necessarily real) zeros $-$then 
\[
\lambda_1^2 + \lambda_2^2 + \cdots +\lambda_{n-r}^2
\ = \ 
(n - r)^2 \hsy 
\Big(
\frac{c_{r + 1}}{c_r}
\Big)^2 
- 2 (n - r) (n - r - 1) \hsy 
\frac{c_{r + 2}}{c_r}
\]
and 
\[
\lambda_1 \hsy \lambda_2 \cdots \lambda_{n-r}
\ = \ 
(-1)^{n-r} \hsx (n - r)! \hsx \frac{c_n}{c_r}.
\]
But
\[
\frac{\lambda_1^2 + \lambda_2^2 + \cdots +\lambda_{n-r}^2}{n-r}
\ \geq \ 
((\lambda_1 \hsy \lambda_2 \cdots \lambda_{n-r})^2)^{\frac{1}{n-r}}.
\]

\noindent
Therefore
\[
\abs{c_n}
\ < \ 
C \hsx 
\frac{(M n)^{(n-r)/2}}{(n-r)!},
\]
where $C$ and $M$ are positive constants independent of $n$.  
And this estimate will do the trick.
\\[-.25cm]
\end{x}

\begin{x}{\small\bf LEMMA} \ 
Let $\gamma_0, \gamma_1, \ldots$ be a multiplier sequence. Suppose that 
\[
c_0 + c_1 z
+ \cdots + 
c_d \hsy z^d
\]
is a real polynomial whose zeros are real and of the same sign $-$then the zeros of the real polynomial 
\[
\gamma_0 \hsy c_0 \hsy  +  \gamma_1\hsy  c_1  \hsy z 
+ \cdots + 
\gamma_d \hsy  c_d  \hsy z^d
\]
are real.
\\[-.5cm]

PROOF \ 
Thanks to 12.12, the zeros of the real polynomial 
\[
\gamma_0 \hsy c_0 + 1! \binom{n}{1}  \gamma_1 \hsy  c_1  \hsy z 
+ \cdots + 
d ! \hsy \binom{n}{d} \hsy\gamma_d \hsy  c_d  z^d
\qquad (n > d)
\]
are real.  
Replacing $z$ by 
$\ds\frac{z}{n}
$, 
it follows that the zeros of the real polynomial 
\[
\gamma_0 \hsy c_0 \hsx + \hsx   \gamma_1 \hsy  c_1  \hsy z 
\hsx + \hsx
\cdots 
\hsx + \hsx 
\Big(1  -  \frac{1}{n}\Big) 
\hsy 
\Big(1 - \frac{2}{n}\Big) 
\cdots 
\Big(1 - \frac{d - 1}{n}\Big) 
\hsy\gamma_d \hsy  c_d  \hsy z^d
\] 
are real so, upon letting $n \ra \infty$, we conclude that the zeros of the real polynomial
\[
\gamma_0 \hsy c_0 \hsy  \hsx + \hsx  \gamma_1 \hsy  c_1  \hsy z 
\hsx + \hsx
\cdots 
\hsx + \hsx  
\gamma_d  \hsy  c_d  \hsy z^d
\]
are real.
\\[-.5cm]

[Note: \ 
The stated property is characteristic.  
Proof: \ 
The zeros of the real polynomial
\[
(1 + z)^n 
\ = \ 
\sum\limits_{k = 0}^n \ 
\binom{n}{k} \hsy z^k
\]
are real and of the same sign.]
\\[-.25cm]
\end{x}

\begin{x}{\small\bf APPLICATION} \ 
Let $\gamma_0, \gamma_1, \ldots$ be a multiplier sequence $-$then the 
\un{Turan inequalities} obtain: 
\[
\gamma_n^2 - \gamma_{n-1} \hsy \gamma_{n + 1} 
\ \geq \ 
0 
\qquad (n = 1, 2, \ldots).
\]

[The zeros of the real polynomial 
\[
z^{n-1} + 2 \hsy z^n + z^{n + 1}
\]
are real and $\leq 0$.  
Therefore the zeros of the real polynomial 
\[
\gamma_{n-1} \hsy z^{n-1} + 2 \hsy \gamma_n \hsy z^n + \gamma_{n+1} \hsy z^{n + 1}
\]
are real, from which the assertion.]
\\[-.25cm]
\end{x}

\begin{x}{\small\bf LAGUERRE CRITERION} \ 
\ 
Let $Q(x)$ be a real polynomial whose zeros are real and lie outside the interval $[0,d]$ 
$-$then for any real sequence 
$c_0, c_1, \ldots, c_d$, 
the number of nonreal zeros of the real polynomial 
\[
 Q(0) \hsy c_0 \hsy  +  Q(1) \hsy  c_1  \hsy z 
+ \cdots + 
 Q(d) \hsy  c_d  \hsy z^d
\]
is $\leq$ the number of nonreal zeros of the real polynomial 
\[
c_0 + c_1 z
+ \cdots + 
c_d \hsy z^d.
\]
[Note: \ 
Accordingly, if the zeros of 
\[
c_0 + c_1 z
+ \cdots + 
c_d \hsy z^d
\]
are real, then the zeros of 
\[
 Q(0) \hsy c_0 \hsy  \ +  \ Q(1) \hsy  c_1  \hsy z 
\ + \cdots + \ 
 Q(d) \hsy  c_d  \hsy z^d
\]
are also real.]
\\[-.25cm]
\end{x}

\begin{x}{\small\bf THEOREM} \ 
Let $f \in \sL - \sP$ and assume that the zeros of $f$ are negative.  
Suppose that 
\[
c_0 + c_1 z
+ \cdots + 
c_d \hsy z^d
\]
is a real polynomial whose zeros are real $-$then the zeros of the real polynomial
\[
f(0) \hsy c_0 \hsy  +  f(1) \hsy  c_1  \hsy z 
+ \cdots + 
f(d) \hsy  c_d   \hsy z^d
\]
are real.
\\[-.5cm]

PROOF \ 
Take $f (0) = 1$ and write
\[
f (z) 
\ = \ 
e^{a \hsy z^2 + b \hsy z} \ 
\prod\limits_{n = 1}^\infty \ 
\Big(
1 - \frac{z}{\lambda_n}
\Big) 
\hsx
e^{z / \lambda_n}
\qquad (\tcf. \ 10.19).
\]
Choose $k \hsx > \hsx  0$: $\sqrt{k} \hsx > \hsx d \hsy \sqrt{-a}$ \ $(a \leq 0)$ and put
\[
Q_k(z) 
\ = \ 
\Big(
1 + \frac{a z^2}{k}
\Big)^k
\hsy
\Big(
1 - \frac{z}{\lambda_1}
\Big)
\hsy
\cdots
\hsy
\Big(
1 - \frac{z}{\lambda_k}
\Big),
\]
the interval of exclusion thus being $[0,d]$.  
Let
\[
B_k 
\ = \ b + \frac{1}{\lambda_1} + \cdots + \frac{1}{\lambda_k}.
\]
Then the zeros of the real polynomial 
\[
c_0  + c_1 \hsy  e^{B_k} \hsy z 
+ \cdots + 
c_d \hsy  e^{d \hsy B_k} \hsy z^d
\]
are real, hence the zeros of the real polynomial 
\[
c_0 \hsy Q_k(0)  + c_1 \hsy Q_k(1) \hsy e^{B_k} \hsy z 
+ \cdots + 
c_d \hsy Q_k(d) \hsy e^{d \hsy B_k} \hsy z^d
\]
are also real.  
Now let $k \ra \infty$.
\\[-.25cm]
\end{x}

\qquad
{\small\bf \un{N.B.}} \ 
An additional assumption to the effect that the zeros of 
\[
c_0 + c_1 z
+ \cdots + 
c_d \hsy z^d
\]
are of the same sign is inutile.
\\[-.25cm]

\begin{x}{\small\bf SCHOLIUM} \ 
If $f \in \sL - \sP$ and if the zeros of $f$ are negative, then the sequence 
$f(0), f(1), \ldots$ 
is a multiplier sequence.
\\[-.25cm]
\end{x}

\begin{x}{\small\bf EXAMPLE} \ 
Take 
$
\ds
f (z)  = e^{z^2 \log q}
$
$(0 < q \leq 1)$ 
$-$then 
$
\ds
f (n)  = q^{n^2}
$, 
so 
$
\ds
\Big\{q^{n^2}\Big\}
$
is a multiplier sequence.
\\[-.25cm]
\end{x}

\begin{x}{\small\bf EXAMPLE} \ 
Take
$
\ds
f (z)  = \frac{1}{\Gamma (z + 1)}
$
(cf. 10.30) $-$then 
$
\ds
f (n) = \frac{1}{n !}
$, 
so 
$
\ds
\Big\{
\frac{1}{n !} : n = 0, 1, \ldots
\Big\}
$
is a multiplier sequence.
\\[-.25cm]

[Note: \ 
Given $\alpha > 0$, put $(\alpha)_0 = 1$ and 
\[
(\alpha)_n
\ = \ 
\alpha \hsy (\alpha + 1) \hsy \cdots (\alpha + n - 1) 
\qquad (n \geq 1).
\]
Take now
\[
f (z)
\ = \ 
\frac{\Gamma (\alpha)}{\Gamma ( z + \alpha)}.
\]
Then
\[
f (n)
\ = \ 
\frac{\Gamma (\alpha)}{\Gamma ( n + \alpha)} 
\ = \ 
\frac{1}{(\alpha)_n},
\]
so 
$
\ds
\Big\{
\frac{1}{(\alpha)_n} : n = 0, 1, \ldots
\Big\}
$ is a multiplier sequence.]
\\
\end{x}

\begin{x}{\small\bf THEOREM} \ 
Let $f \in \sL - \sP$ and assume that the zeros of $f$ are negative.  
Suppose that 
\[
F (z) 
\ = \ C_0 + C_1 \hsy z + \cdots
\]
is in $\sL - \sP$ $-$then the series 
\[
f (0) \hsy C_0 + f (1) \hsy C_1 \hsy z + \cdots
\]
is a real entire function and, as such, is in $\sL - \sP$ .
\\[-.5cm]

PROOF \ 
The initial claim is that the series
\[
f (0) \hsy C_0 + f (1) \hsy C_1  z + \cdots
\]
is convergent for every $z$.  
Thus decompose $f$ per 10.19: 
\[
f (z) 
\ = \ 
C \hsy e^{a \hsy z^2 + bz} \ 
\prod\limits_{n = 1}^\infty \ 
\Big(
1 - \frac{z}{\lambda_n}
\Big) 
\hsx
e^{z / \lambda_n}.
\]
Then
\[
(1 + t) \hsy e^{-t} 
\ \leq \ 1
\qquad (t \geq 0)
\]
\qquad 
$\implies$
\[
\Big(
1 - \frac{t}{\lambda_n}
\Big) 
\hsx
e^{t / \lambda_n}
\ = \ 
\Big(
1 + \Big(\frac{t}{-\lambda_n}\Big)
\Big) 
\hsx
e^{-(t / -\lambda_n)}
\ \leq \ 
1 
\qquad (\lambda_n < 0).
\]
So, for $k$ a nonnegative integer, 
\[
\abs{f (k)} 
\ \leq \ 
\abs{C} 
\hsx 
e^{a \hsy k^2}
\hsx
e^{b \hsy k}
\ \leq \ 
\abs{C} 
\hsx 
e^{b \hsy k}
\qquad (a \leq 0).
\]
Therefore
\[
\underset{k \ra \infty}{\limsupx} \ 
\abs{f (k)}^{1/k}
\hsx 
\abs{C_k}^{1/k}
\ = \ 
0,
\]
which settles the convergence issue.  
To verify the $\sL - \sP$ contention, note first that the zeros of 
\[
\tJ_n (F; z) 
\ = \ 
C_0 + n \hsy C_1 \hsy z + n (n -1) \hsy C_2 \hsy z^2 + \cdots
\]
are real (cf. 12.14).  
Therefore the zeros of the real polynomial 
\[
f(0) \hsy C_0  + n \hsy f(1) \hsy C_1 \hsy z + n (n -1) \hsy f(2) \hsy C_2 \hsy z^2 + \cdots
\]
are real (cf. 12.25).  
But this polynomial is the $n^\nth$ Jensen polynomial of the series 
\[
f(0) \hsy C_0  + f(1) \hsy C_1 \hsy z  + \cdots,
\]
so another application of 12.14 finishes the argument.
\\[-.25cm]
\end{x}

\begin{x}{\small\bf EXAMPLE} \ 
Take $F (z) = e^z$ $-$then 
\[
\sum\limits_{n = 0}^\infty \ 
\frac{f (n)}{n !} \hsy z^n
\]
is in $\sL - \sP$.
\\[-.25cm]
\end{x}

\begin{x}{\small\bf EXAMPLE} \ 
Take $F (z) = e^{-z^2}$ $-$then 
\[
\sum\limits_{n = 0}^\infty \ 
(-1)^n \hsx
\frac{f (2 \hsy n)}{n !} \hsy z^{2 \hsy n}
\]
is in $\sL - \sP$.
\\[-.25cm]
\end{x}

\begin{x}{\small\bf EXAMPLE} \ 
Fix a positive integer $m$ and take
\[
f (z)
\ = \ 
\frac{\Gamma (z + 1)}{\Gamma (m \hsy z + 1)}.
\]
Then
\[
f (n)
\ = \ 
\frac{n !}{(m \hsy n)!} \hsx , 
\]
hence
\[
\sum\limits_{n = 0}^\infty \ 
\frac{z^n}{(m \hsy n)!}
\ \equiv \ 
\ML_m (z) 
\qquad (\tcf. \ 2.28)
\]
is in $\sL - \sP$.
\\[-.5cm]

[Note: \ 
The poles of the numerator, viz. 
$
-1, -2, \ldots
$, 
are absorbed by the poles of the denominator, viz.
$
\ds
-\frac{1}{m}, \hsy -\frac{2}{m}, \hsy \ldots, -\frac{m}{n},  \ldots \ 
$.]
\\[-.25cm]
\end{x}

\begin{x}{\small\bf EXAMPLE} \ 
Recall that the Bessel function $\tJ_\nu (z)$ of the first kind of real index $\nu > -1$ is defined by the series
\[
\Big(
\frac{z}{2}
\Big)^\nu \ 
\sum\limits_{n = 0}^\infty \ 
\frac{\ds (-1)^n \hsy \Big(\frac{z}{2}\Big)^{2 \hsy n}}{n ! \hsy \Gamma (\nu + n  + 1)}
\qquad (\tcf. \ 2.29).
\]
To apply the foregoing machinery, rewrite this as 
\[
\tJ_\nu (z)
\ = \ 
\Big(
\frac{z}{2}
\Big)^\nu 
\hsx
\invh_\nu 
\big(
\frac{z}{2}
\big),
\]
where
\[
\invh_\nu (z)
\ = \ 
\sum\limits_{n = 0}^\infty \ 
 (-1)^n \ \frac{f_\nu (2 \hsy n)}{n !} \hsx z^{2 \hsy n}.
\]
Here
\[
f_\nu (z)
\ = \ 
\ds\frac{1}
{
\raisebox{-.25cm}
{$
\ds
\Gamma 
\Big(
\nu + \ds\frac{z}{2}+ 1
\Big)
$}
}
\]
is in $\sL - \sP$ and its zeros are negative (since $\nu > -1$).  
Therefore the zeros of $\tJ_\nu (z)$ are real\footnote[2]{\vspace{.11 cm}
E. Lommel,  
\textit{Studien \"uber die Bessel'schen Functionen}, Teubner, Leipzig, 1868, \S19.}.
\\[-.25cm]
\end{x}


\begin{x}{\small\bf EXAMPLE} \ 
Given $p = 1, 2, \ldots$, 
\[
\Phi_{2 \hsy p} (z) 
\ = \ 
\int\limits_0^\infty \ 
\exp (-t^{2 \hsy p})
\hsy 
\cos z \hsy t 
\ \td t
\qquad (\tcf. \ 2.30)
\]
is in $\sL - \sP$.
\\[-.5cm]

[In fact, 
\[
2 \hsx p \ 
\Phi_{2 \hsy p} (z) 
\ = \ 
\sum\limits_{n = 0}^\infty \ 
(-1)^n \hsx 
\frac{f_p (2 \hsy n)}{n !} \hsx 
z^{2 \hsy n},
\]
where
\[
f_p (z)
\ = \ 
\frac{\ds\Gamma\Big(\frac{z}{2} + 1\Big) \hsx \Gamma\Big(\frac{z + 1}{2 \hsy p}\Big) }{\Gamma (z + 1)},
\]
the poles of the numerator, viz. 
\[
-2,\hsy -4,\hsy -6,\hsy \ldots, -1, -(1 + 2\hsy p), \hsy -(1 + 4 \hsy p), \hsy \ldots,
\]
being absorbed by the poles of the denominator, viz. $-1, -2, -3, \ldots \ $.
\\[-.5cm]

[Note: \ 
$\Phi_2 (z)$ has no zeros but 
$\Phi_4 (z), \Phi_6 (z), \ldots$, 
have an infinity of zeros.  
\\[-.5cm]

\noindent
Proof: \ 
The order of 
$\Phi_{2 \hsy p} (z)$ 
is 
$
\ds
\frac{2p}{2p - 1}
$, 
which lies strictly between 1 and 2 if $p > 1$, so one can cite 7.4.]
\\[-.25cm]
\end{x}

If $f \in \sL - \sP$ then $f^\prime \in \sL - \sP$ (cf. 10.20 and 10.25).
\\[-.5cm]

[Note: \ 
Letting 
$\gamma_0, \gamma_1, \ldots$ 
be the multiplier sequence associated with $f$, 
it follows that 
$\gamma_0^\prime = \gamma_1, 
\gamma_1^\prime = \gamma_2, 
\ldots$ 
is a multiplier sequence (namely the one associated with $f^\prime$).]
\\

\begin{x}{\small\bf EXAMPLE} \ 
The $n^\nth$ Hermite polynomial is, by definition, 
\[
H_n (z) 
\ = \ 
(-1)^n \ e^{z^2} \ \frac{\td^n}{\td z^n} \hsx e^{-z^2} 
\qquad (\tcf. \ 8.17),
\]
so
\[
\frac{\td^n}{\td z^n} \hsy e^{-z^2}
\ = \ 
(-1)^n \hsx 
H_n (z) \hsy e^{-z^2}.
\]
The fact that $e^{-z^2}$ is in $\sL - \sP$ then implies that 
$
\ds
\frac{\td^n}{\td z^n} \hsy e^{-z^2}$ is in $\sL - \sP$, 
thus the zeros of $H_n (z)$ must be real.
\\[-.25cm]
\end{x}

While $\sL - \sP$ is not a vector space, there are circumstances in which it is closed under addition.
\\[-.25cm]

\begin{x}{\small\bf LEMMA} \ 
If $f \in \sL - \sP$, then $\forall \ a \in \R$, 
\[
a \hsy f + f^\prime \in \sL - \sP
\qquad (\tcf. \ 12.10).
\]

PROOF \
The product 
$
f(z) \hsy e^{a \hsy z}
$
is in 
$\sL - \sP$, 
as is the derivative 
$
\ds
\frac{\td}{\td z} (f(z) \hsy e^{a \hsy z})
$, 
as is the product
$
\ds
e^{-a \hsy z} \hsy 
\frac{\td}{\td z} (f(z) \hsy e^{a \hsy z})
$, 
thus
\[
a \hsy f (z) + f^\prime  (z)
\] 
is in $\sL - \sP$.
\\[-.25cm]
\end{x}

\begin{x}{\small\bf EXAMPLE} \ 
Let $p$ be a real polynomial with real zeros only.  
Take $\alpha > 0$, $\beta \in \R$, and define $F$ by
\[
F (z) 
\ = \ 
\int\limits_{-\infty}^\infty \ 
p(\sqrt{-1} \hsx t) \hsy 
\exp (-\alpha \hsy t^2 + \sqrt{-1} \hsx \beta \hsy t + \sqrt{-1} \hsx z \hsy t ) 
\ \td t.
\]
Then 
$F \in \sL - \sP$.
\\[-.5cm]

[Supposing that $p$ is monic, write
\[
p (z) 
\ = \ 
(z + a_1) \cdots (z + a_n) 
\qquad  (a_1, \ldots, a_n \in \R).
\]
Put
\[
F_0 (z) 
\ = \ 
\int\limits_{-\infty}^\infty \ 
\exp (-\alpha \hsy t^2 + \sqrt{-1} \hsx \beta \hsy t + \sqrt{-1} \hsx z \hsy t ) 
\ \td t.
\]
Then
\[
F_0 (z) 
\ = \ 
\Big(
\frac{\pi}{\alpha}
\Big)^{1/2}
\hsx 
\exp 
\Big(
\frac{- (z + \beta)^2}{4 \hsy \alpha}
\Big), 
\]
so $F_0 \in \sL - \sP$.  
Now define $F_k$ $(k = 1, \ldots, n)$ by
\[
F_k (z) 
\ = \ 
\int\limits_{-\infty}^\infty \ 
p_k (\sqrt{-1} \hsx t) 
\hsx
\exp (-\alpha \hsy t^2 + \sqrt{-1} \hsx \beta \hsy t + \sqrt{-1} \hsx z \hsy t ) 
\ \td t,
\]
where
\[
p_k (z) 
\ = \ 
(z + a_1) \cdots (z + a_k).
\]
Then
\allowdisplaybreaks\begin{align*}
F_1 \ 
&= \ a_1 \hsy F_0 + F_0^\prime
\\[11pt]
&
\ \vdots
\\[11pt]
F \ 
= \ F_n \ 
&= \ a_n \hsy F_{n -1 } + F_{n -1 }^\prime,
\end{align*}
so $F \in \sL - \sP$.]
\\[-.25cm]
\end{x}

\[
\text{APPENDIX}
\]

A multiplier sequence 
$\gamma_0, \gamma_1, \ldots$
is said to be \un{strict} if it has the following property: \ 
Given any real polynomial
\[
c_0 + c_1 z
+ \cdots + 
c_d \hsy z^d
\]
whose zeros are real, the zeros of the real polynomial
\[
\gamma_0  \hsy c_0 
+ 
\gamma_1 \hsy  c_1 \hsy z
+ \cdots + 
\gamma_d \hsy  c_d \hsy z^d
\]
are also real (cf. 12.22).
\\[-.25cm]

\begin{spacing}{1.75}
\qquad
{\small\bf EXAMPLE} \ 
Let $f \in \sL - \sP$ and assume that the zeros of $f$ are negative 
$-$then the sequence 
$f(0), f(1), \ldots$ 
is a strict multiplier sequence (cf. 12.25).  
In particular: \ 
$
\ds
\Big\{
\frac{1}{n !} : n = 0, 1, \ldots
\Big\}
$
is a strict multiplier sequence (cf. 12.28 (or 12.13)).
\\[-.25cm]
\end{spacing}

\qquad
{\small\bf LEMMA} \ 
A strict multiplier sequence acting on a polynomial whose zeros are real and of the same sign preserves 
the reality and the sign of the zeros.
\\[-.25cm]\\

\qquad
{\small\bf EXAMPLE} \ 
Take 
$f (z) = (z^2 + 2 z - 1)\hsy e^z$ 
and consider the corresponding multiplier sequence 
$\{-1 + n + n^2 : n = 0, 1, \ldots \}$ 
$-$then its action on 
$(z + 1)^2$ is 
\[
-1 \hsy (1) + 1 \hsy (2)z + 5 \hsy (2) z^2.
\]
\begin{spacing}{1.75}
\noindent
The zeros of this polynomial are 
$\ds
\frac{-1 \pm \sqrt{11}}{10}
$, 
hence are real but of opposite sign.  
Therefore the multiplier sequence 
$\{- 1 + n + n^2 : n = 0, 1, \ldots\}$
is not strict.
\\[-.25cm]
\end{spacing}

\qquad
{\small\bf DEFINITION} \ 
Given two sequences
\[
\begin{cases}
a_0, a_1, \ldots
\\[4pt]
\ 
b_0, b_1, \ldots
\end{cases}
\]
of real numbers, 
their \un{component wise product} is the sequence 
$a_0 \hsy b_0, a_1 \hsy b_1, \ldots \ $.
\\

\qquad
{\small\bf LEMMA} \ 
If
\[
\begin{cases}
\ 
\alpha_0, \alpha_1, \ldots
\\[4pt]
\ 
\beta_0, \beta_1, \ldots
\end{cases}
\]
are strict multiplier sequences, then so is their component wise product.
\\

\qquad
{\small\bf LEMMA} \ 
If

\[
\begin{cases}
\ 
\alpha_0, \alpha_1, \ldots
\\[4pt]
\ 
\beta_0, \beta_1, \ldots
\end{cases}
\]
are multiplier sequences and if $\alpha_0, \alpha_1, \ldots$ is strict, 
then their component wise product is a multiplier sequence.
\\[-.5cm]

PROOF \ 
Let
\[
c_0 + c_1 z
+ \cdots + 
c_d \hsy z^d
\]
be a real  polynomial whose zeros are real and of the same sign $-$then
\[
\alpha_0  \hsy c_0 
+ 
\alpha_1 \hsy  c_1 \hsy z
+ \cdots + 
\alpha_d \hsy  c_d \hsy z^d
\]
is a real polynomial whose zeros are real and of the same sign, 
thus the zeros of the real polynomial
\[
\alpha_0 \hsy \beta_0 \hsy c_0 
+ 
\alpha_1 \hsy \beta_1 \hsy c_1 \hsy z
+ \cdots + 
\alpha_d \hsy \beta_d \hsy c_d \hsy z^d
\]
are real (cf. 12.22), which implies that 
$\alpha_0 \hsy \beta_0, \alpha_1 \hsy \beta_1, \ldots$ is a multiplier sequence 
(see the comment appended to 12.22).
\\

\qquad
{\small\bf APPLICATION} \ 
For $f \in \sL - \sP$, say
\[
f (z) 
\ = \ 
\sum\limits_{n = 0}^\infty \ 
c_n \hsy z^n.
\]
Then $c_0, c_1, \ldots$ is a multiplier sequence.
\\[-.5cm]

[For
\[
c_n 
\ = \ 
\frac{\gamma_n}{n !}
\]
\begin{spacing}{1.75}
\noindent
and
$
\ds
\Big\{\frac{1}{n !} : n = 0, 1, \ldots \Big\}
$
is a strict multiplier sequence while 
$\gamma_0, \gamma_1, \ldots$ 
is a multiplier sequence (cf. 12.14).]
\\[-.75cm]
\end{spacing}

[Note: \ 
A priori, 
\[
c_n^2 - c_{n - 1} \hsy c_{n + 1} 
\ \geq \ 
0
\qquad (n = 1, 2, \ldots) \quad (\tcf. \ 12.23)
\]
but this can be sharpened:
\[
\gamma_n^2 - \gamma_{n - 1} \hsy \gamma_{n + 1} 
\ \geq \ 
0
\]
\qquad 
$\implies$
\[
(n !)^2 \hsy c_n^2 - (n - 1)! \hsy (n + 1)! \hsy c_{n - 1} \hsy c_{n + 1} 
\ \geq \ 
0
\]
\qquad 
$\implies$
\[
n \hsy c_n^2 - (n + 1) \hsy c_{n - 1} \hsy c_{n + 1} 
\ \geq \ 
0
\]
\qquad 
$\implies$
\[
c_n^2 - \Big(1  + \frac{1}{n} \Big)c_{n - 1} \hsy c_{n + 1} 
\ \geq \ 
0
\]
\qquad 
$\implies$
\[
c_n^2 - c_{n - 1} \hsy c_{n + 1} 
\ \geq \ 
0.]
\]


\chapter{
$\boldsymbol{\S}$\textbf{13}.\quad  CHARACTERIZATIONS}
\setlength\parindent{2em}
\setcounter{theoremn}{0}
\renewcommand{\thepage}{\S13-\arabic{page}}

\qquad
Let

\[
f (z) 
\ = \ 
\sum\limits_{n = 0}^\infty \ 
c_n \hsy z^n 
\]
be in $\sL - \sP$ $-$then
\[
c_n
\ = \ 
\frac{\gamma_n}{n !} 
\qquad (\gamma_n = f^{(n)} (0))
\]
and 
$\gamma_0, \gamma_1, \ldots$
is a multiplier sequence (cf. 12.14).  
Therefore (cf. 12.23)

\[
\gamma_n^2 - \gamma_{n-1} \hsy \gamma_{n+1} 
\ \geq \ 
0
\qquad (n = 1, 2, \ldots).
\]
\\[-1.cm]

\begin{x}{\small\bf EXAMPLE} \ 
Consider the Hermite polynomials 
$\{H_n : n = 0, 1, \ldots\}$ (cf. 12.35) $-$then for real $t$ and complex $z$, 

\[
\exp (2 t z - z^2) 
\ = \ 
\sum\limits_{n = 0}^\infty \ 
\frac{H_n (t)}{n !}
\hsy z^n.
\]
Since $\forall \ t$, the function 

\[
z \ra \exp (2 t z - z^2) 
\]
is in $\sL - \sP$, it follows that 

\[
H_n^2 (t) - H_{n - 1} (t) \hsy H_{n + 1} (t) 
\ \geq \ 
0
\qquad (n = 1, 2, \ldots).
\]
\\[-1.25cm]
\end{x}

\begin{x}{\small\bf EXAMPLE} \ 
Consider the Laguerre polynomials 
$\{L_n^{(\alpha)} : n = 0, 1, \ldots \}$ of index $\alpha > -1$ and degree $n$, thus

\[
L_n^{(\alpha)} (t) 
\ = \ 
\frac{t^{-\alpha} \hsy e^t}{n !} \hsx
\frac{\td^n}{\td t^n} \hsy
e^{-t} \hsy 
t^{n + \alpha} 
\qquad (\tcf. \ 8.17 \ (L_n^{(0)} \equiv L_n)),
\]
where
\[
L_n^{(\alpha)} (0)
\ = \ 
\frac{(1 + \alpha)_n}{n !}\hsx .
\]
In terms of the Bessel function $\tJ_\alpha$, for real $t > 0$ and complex $z$, 
\allowdisplaybreaks
\allowdisplaybreaks\begin{align*}
\Gamma (1 + \alpha) (t z)^{-\alpha / 2} \hsy 
\tJ_\alpha (2 \sqrt{t \hsy z}) \ 
&=\ 
\sum\limits_{n = 0}^\infty \ 
\frac{L_n^{(\alpha)}  (t)}{(1 + \alpha)_n} \hsx z^n
\\[15pt]
&=\ 
\sum\limits_{n = 0}^\infty \ 
\frac{n !}{(1 + \alpha)_n} 
\
L_n^{(\alpha)} (t) 
\
\frac{z^n}{n !}
\\[15pt]
&=\ 
\sum\limits_{n = 0}^\infty \ 
\frac{L_n^{(\alpha)}  (t)}{L_n^{(\alpha)} (0)}
\ 
\frac{z^n}{n !}.
\end{align*}
Since $\forall \ t > 0$, the function 

\[
z \ra (t z)^{-\alpha / 2}
\tJ_\alpha (2 \hsx \sqrt{t \hsy z})
\]
is in $\sL - \sP$ (cf. 12.33), it follows that 

\[
\bigg[
\frac{L_n^{(\alpha)}  (t)}{L_n^{(\alpha)} (0)}
\bigg]^2
\ - \
\frac{L_{n-1}^{(\alpha)}  (t)}{L_{n-1}^{(\alpha)} (0)} 
\
\frac{L_{n+1}^{(\alpha)}  (t)}{L_{n+1}^{(\alpha)} (0)} 
\ \geq \ 
0
\qquad (n = 1, 2, \ldots).
\]

[Note: \ 
As we know, 

\[
\Big(\frac{z}{2}\Big)^{-\alpha} \hsy
\tJ_\alpha (z) 
\in \sL - \sP, 
\]
so by evenness, 
\[
\Big(\frac{\sqrt{z}}{2}\Big)^{-\alpha} \hsy
\tJ_\alpha (\sqrt{z}) 
\in \sL - \sP
\]

\qquad 
$\implies$
\[
2^\alpha \hsy 
z^{-\alpha/2}  \hsy
 \tJ_\alpha (\sqrt{z}) 
\in \sL - \sP
\]

\qquad 
$\implies$
\[
2^\alpha \hsy
(4 z)^{-\alpha / 2}
\tJ_\alpha (2 \sqrt{z}) 
\in \sL - \sP
\]

\qquad 
$\implies$
\[
z^{-\alpha/2}  \hsy
\tJ_\alpha (2 \sqrt{z}) 
\in \sL - \sP.]
\]
\\[-1.75cm]
\end{x}

\begin{x}{\small\bf LEMMA} \ 
If $f \in \sL - \sP$, then for all real $t$, 

\[
(f^{(n)} (t))^2 
\hsy - \hsy 
f^{(n-1)} (t) f^{(n+1)} (t) 
\ \geq \ 
0 
\qquad (n \geq 1),
\]
with equality iff $f^{(n-1)} (z)$ is of the form $C \hsy e^{b \hsy z}$ or $t$ is a multiple zero of $f^{(n-1)} (z)$.
\\[-.5cm]

PROOF \ 
Decompose $f$ per 10.19:

\[
f (z) 
\ = \ 
C \hsy z^m \hsy e^{a z^2 + b z} \ 
\prod\limits_{n = 1}^\infty \ 
\big(
1 - \frac{z}{\lambda_n} 
\big) \hsx e^{z / \lambda_n}.
\]
Then

\[
\frac{f^\prime (t)}{f (t)}
\ = \ 
\frac{m}{t} \hsx + \hsx 2 a t \hsx + \hsx b \hsx + \hsx
\sum\limits_{n = 1}^\infty \ 
\big(
\frac{1}{t - \lambda_n} 
\hsx + \hsx
\frac{1}{\lambda_n} 
\big)
\]

\qquad 
$\implies$

\allowdisplaybreaks\begin{align*}
\frac{\td}{\td t} \hsy \Big(\frac{f^\prime (t)}{f (t)}\Big)
&= \ 
\frac{f (t) \hsy f^{\prime\prime} (t) - (f^\prime (t))^2 }{(f (t))^2}
\\[15pt]
&= \ 
- \frac{m}{t^2} \hsx + \hsx 2 a \hsx - \hsx
\sum\limits_{n = 1}^\infty \ 
\frac{1}{(t - \lambda_n)^2}.
\end{align*}
If $f (z) = C \hsy e^{b \hsy z}$ or if $t$ is a multiple zero of $f (z)$, then 

\[
f (t) \hsy f^{\prime\prime} (t) - (f^\prime (t))^2 
\ = \ 
0.
\]
On the other hand, if 
$f (z) \neq C \hsy e^{b \hsy z}$ 
and if $c$ is not a zero of $f (z)$, then 

\[
- \frac{m}{c^2} + 2 a - 
\sum\limits_{n = 1}^\infty \ 
\frac{1}{(c - \lambda_n)^2}
\ < \ 
0
\]

\qquad 
$\implies$

\[
f (c) \hsy f^{\prime\prime} (c) - (f^\prime (c))^2
\ < \ 
0,
\]
so by continuity, 
\[
f (t) \hsy f^{\prime\prime} (t) - (f^\prime (t))^2 
\ \leq \ 
0
\]
for all real $t$.  
If equality obtains and if $f (z) \neq C \hsy e^{b \hsy z}$, then $t$ must be a zero of $f(z)$ (cf. supra), 
hence $t$ must be a multiple zero of $f(z)$: 

\[
(f^\prime (t))^2
\ = \ 
0
\implies 
f^\prime (t)
\ = \ 
0.
\]
Proceed from here by iteration (bear in mind that $\sL - \sP$ is closed under differentiation (cf. 10.20 and 10.25)).
\\[-.25cm]

[Note: \ 
In particular, 
\[
(f^{(n)} (0))^2 
\hsy - \hsy 
f^{(n-1)} (0) f^{(n+1)} (0) 
\ \geq \ 
0, 
\]
i.e., 
\[
\gamma_n^2 
\hsy - \hsy 
\gamma_{n-1} \hsy \gamma_{n+1}
\ \geq \ 
0 
\qquad (n = 1, 2, \ldots) .]
\]
\\[-1.cm]
\end{x}

\begin{x}{\small\bf EXAMPLE} \ 
Take 
\[
f (z) 
\ = \ 
z (z^2 + 1).
\]
Then
\[
f^\prime (t)^2  - f (t) \hsy  f^{\prime\prime} (t) 
\ = \ 
3 t^4 + 1 
\ > \ 
0.
\]
Still, $f \notin \sL - \sP$ (because it has the nonreal zeros $\pm \sqrt{-1}$).
\\[-.25cm]
\end{x}

\begin{x}{\small\bf EXAMPLE} \ 
Take 
\[
f (z) 
\ = \ 
e^z - e^{2 z}.
\]
Then
\[
(f^{(n)} (t))^2 
\hsy - \hsy 
f^{(n-1)} (t) f^{(n+1)} (t) 
\ = \ 
2^{n-1} \hsy e^{3 \hsy t} 
\ > \ 
0 
\qquad (n \geq 1).
\]
Still, $f \notin \sL - \sP$ (because it has the nonreal zeros 
$2 \hsy \pi \hsy  \sqrt{-1} \hsx k$ $(k = \pm 1, \pm 2, \ldots)$).
\\[-.25cm]
\end{x}

Therefore the inequalities 

\[
(f^{(n)} (t))^2 
\hsy - \hsy 
f^{(n-1)} (t) f^{(n+1)} (t) 
\ \geq \ 
0 
\qquad (n \geq 1)
\] 
do not serve to characterize the elements of $\sL - \sP$ (even if they are strict).
\\

\begin{x}{\small\bf NOTATION} \ 
Given a real entire function $f$, 
let $L_0 (f) (t) = f(t)^2$ and for $n = 1, 2, \ldots$,  let
\[
L_n (f) (t) 
\ = \ 
\sum\limits_{k = 0}^{2 n} \ 
\frac{(-1)^{k+n}}{(2 n) !} \hsy 
\binom{2 n}{k} \hsy 
f^{(k)} (t) \hsy
f^{(2 n-k)} (t)
\qquad (t \in \R).
\]
\\[-1cm]
\end{x}

\qquad
{\small\bf \un{N.B.}} \ 
For the record, 
\allowdisplaybreaks
\allowdisplaybreaks\begin{align*}
L_1 (f) (t) \
&= \ 
\sum\limits_{k = 0}^{2} \ 
\frac{(-1)^{k+1}}{2} \hsy 
\binom{2}{k} \hsy 
f^{(k)} (t) \hsy
f^{(2-k)} (t)
\\[15pt]
&= \ 
-
\frac{f (t) \hsy f^{\prime\prime} (t)}{2} 
\ + \ 
(f^\prime (t))^2 
\ - \ 
\frac{f^{\prime\prime} (t) \hsy f (t)}{2} 
\\[15pt]
&= \ 
(f^\prime (t))^2 
\hsx - \hsx 
f (t) \hsy f^{\prime\prime} (t).
\end{align*}
\\[-1.cm]

\begin{x}{\small\bf THEOREM} \ 
Let $f \in A - \sL - \sP$ (cf. 10.31) $-$then $f \in 0 - \sL - \sP$ $(= \sL - \sP)$ iff 
$\forall \ n \geq 0$ and $\forall \ t \in \R$, 

\[
L_n (f) (t) 
\ \geq \ 
0.
\]
\\[-1.5cm]
\end{x}


Some preparation will help ease the way.
\\[-.25cm]

\begin{x}{\small\bf NOTATION} \ 
Given a real entire function $f$, for fixed $x \in \R$, let
\allowdisplaybreaks\begin{align*}
f_x (y) \ 
&=\ 
\abs{f (x + \sqrt{-1} \hsx y)}^2 
\\[11pt]
&\equiv \ 
f (x + \sqrt{-1} \hsx y) \hsy f (x - \sqrt{-1} \hsx y).
\end{align*}
Then $f_x$ is an even function of $y$ and 

\[
f_x (y)
\ = \ 
\sum\limits_{n = 0}^\infty \ 
\Lambda_n (f) (x) y^{2 n}, 
\]
where
\[
\Lambda_n (f) (x) 
\ = \ 
\frac{f_x^{(2 n)} (0)}{(2 n)!}.
\]
\\[-1.cm]
\end{x}

\begin{x}{\small\bf LEMMA} \ 
We have

\[
\Lambda_n (f) (x) 
\ = \  
L_n (f) (x).
\]

PROOF \ 
In fact, 
\allowdisplaybreaks
\allowdisplaybreaks\begin{align*}
(2 n)! \hsx \Lambda_n (f) (x) \ 
&=\ 
f_x^{(2 n)} (0)
\\[11pt]
&=\ 
\frac{\td}{\td y} \hsy
\abs{f (x + \sqrt{-1} \hsx y)}^2 \hsy
\bigg|_{y = 0}
\\[11pt]
&=\ 
\frac{\td}{\td y} \hsy
(f (x + \sqrt{-1} \hsx y) \hsy f (x - \sqrt{-1} \hsx y))
\bigg|_{y = 0}
\\[11pt]
&=\ 
\sum\limits_{k = 0}^n \ 
\binom{2n}{k}\hsy
\frac{\td^k}{\td y^k} \hsy
f (x + \sqrt{-1} \hsx y) \hsy
\bigg|_{y = 0}
\hsx \cdot \hsx
\frac{\td^{2n-k}}{\td y^{2n-k}} \hsy
f (x - \sqrt{-1} \hsx y) \hsy
\bigg|_{y = 0}
\\[11pt]
&=\ 
\sum\limits_{k = 0}^n \ 
(-1)^{k+n} \hsy
\binom{2n}{k}\hsy
f^{(k)} (x) \hsy
f^{(2n - k)} (x) 
\\[11pt]
&=\ 
(2 n)! \hsx L_n (f) (x).
\end{align*}
\\[-.75cm]

When convenient to do so, write
\[
\begin{cases}
\ 
L_n (f) (t)
\ = \ 
L_n (f (t))
\\[4pt]
\ 
\Lambda_n (f) (t)
\ = \ 
\Lambda_n (f (t))
\end{cases}
.
\]
\\[-.75cm]
\end{x}

\begin{x}{\small\bf LEMMA} \ 
For every real $a$, 
\[
L_n ((x + a) f (x)) 
\ = \ 
(x + a)^2 \hsy 
L_n (f(x)) + L_{n-1} (f(x))
\qquad (n = 1, 2, \ldots).
\]

PROOF \ 
From the definitions, 
\allowdisplaybreaks
\allowdisplaybreaks\begin{align*}
\sum\limits_{n = 0}^\infty \ 
L_n 
&((x + a) f (x)) \hsy 
y^{2 \hsy n}
\\[15pt]
&=\ 
\sum\limits_{n = 0}^\infty \ 
\Lambda_n ((x + a) (f (x)) \hsy y^{2 \hsy n}
\\[15pt]
&=\
\abs{(x + a + \sqrt{-1} \hsx y) \hsy 
f(x + \sqrt{-1} \hsx y)}^2
\\[15pt]
&=\ 
((x + a)^2 + y^2)  \
\sum\limits_{n = 0}^\infty \ 
\Lambda_n (f (x))
\hsy
y^{2 \hsy n}
\\[15pt]
&=\ 
(x + a)^2 \
\sum\limits_{n = 0}^\infty \ 
\Lambda_n (f (x))
\hsy
y^{2 \hsy n}
\ + \
\sum\limits_{n = 0}^\infty \ 
\Lambda_n (f (x))
\hsy
y^{2 \hsy n + 2}
\\[15pt]
&=\ 
(x + a)^2 \
\sum\limits_{n = 0}^\infty \ 
\Lambda_n (f (x))
\hsy
y^{2 \hsy n}
\ + \
\sum\limits_{n = 1}^\infty \ 
\Lambda_{n-1} (f (x))
\hsy
y^{2 \hsy n}
\\[15pt]
&=\ 
(x + a)^2 \hsy
\Lambda_0 (f (x)) 
\ + \
\sum\limits_{n = 1}^\infty \ 
\big[
(x + a)^2 \hsy
\Lambda_n (f (x)) \hsy
\hsx + \hsx
\Lambda_{n-1} (f (x)) 
\big]
\hsy
y^{2 \hsy n}
\\[15pt]
&=\ 
(x + a)^2 \hsy
L_0 (f (x)) 
\hsx + \hsx
\sum\limits_{n = 1}^\infty \ 
\big[
(x + a)^2 \hsy
L_n (f (x)) \hsy
\hsx + \hsx
L_{n-1} (f (x)) 
\big]
\hsy
y^{2 \hsy n}.
\end{align*}
\\[-.75cm]

To establish the necessity in 13.7, it can be assumed that $f$ is a real polynomial with real zeros only.  
For this purpose, proceed by induction on the degree
of $f$, the assertion being clear when $\deg f = 0$.  
If $\deg f > 0$, write $f (x) = (x + a) g(x)$, where $a \in R$ and $g(x)$ is a real polynomial with real zeros only.  
By the induction hypothesis, $L_n (g(x)) \geq 0$ for all $n \geq 0$.  
Now apply 13.10 to see that the same is true of $f$.

Turning to the sufficiency in 13.7, if $f \not\equiv 0$ is not in $\sL - \sP$, then $f$ has a nonreal zero 
$z_0 = x_0 + \sqrt{-1} \hsx y_0$, so

\[
0 
\ = \ 
\abs{f(z_0)}^2 
\ = \ 
\sum\limits_{n = 0}^\infty \ 
L_n (f) (x_0) \hsy y_0^{2 n}
\qquad (y_0 \neq 0).
\]
Since each term in the sum on the right is nonnegative, it follows that $L_n (f) (x_0) = 0$ $\forall \ n \geq 0$, 
hence $\forall \ y \in \R$,

\[
0 
\ = \ 
\abs{f (x_0 + \sqrt{-1} \hsx y)}^2 
\ = \ 
\sum\limits_{n = 0}^\infty \ 
L_n (f) (x_0) \hsy y^{2 \hsy n},
\]
implying thereby that $f \equiv 0$.
\\[-.5cm]

[Note: \ 
The assumption that $f \in A - \sL - \sP$ serves to ensure that if $f \notin 0 - \sL - \sP$ ($= \sL - \sP$), 
then $f$ has a nonreal zero.]
\\[-.25cm]
\end{x}

\begin{x}{\small\bf EXAMPLE} \ 
Take $f (z) = (z^2 + 1) e^z$ $-$then 

\[
\begin{cases}
\ 
L_1 (f) (t) 
\ = \ 2 (t^2 - 1) \hsy e^{2 \hsy t}
\\[4pt]
\ 
L_2 (f) (t) 
\ = \ 
e^{2 \hsy t}
\end{cases}
\]
and $L_n (f) (t) = 0$ $(n > 2)$.  
Here

\[
t^2 < 1 
\implies L_1 (f) (t) < 0
\]
and, of course, $f \notin \sL - \sP$ (but $f \in * - \sL - \sP$).
\\[-.25cm]
\end{x}


\begin{x}{\small\bf THEOREM} \ 
Let $f \in A - \sL - \sP$  (cf 10.31) $-$then $f \in 0 - \sL - \sP$ $(= \sL - \sP)$ iff $\forall \ z$,
\[
\abs{f^\prime (z)}^2 
\ \geq \ 
\Reg (f (z) \hsy \ov{f^{\prime\prime} (z)}).
\]

PROOF \ 
Suppose first that $f \in \sL - \sP$:
\[
\abs{f(x + \sqrt{-1} \hsx y)}^2
\ = \ 
\sum\limits_{n = 0}^\infty \ 
L_n (f) (x) \hsy y^{2 \hsy n}
\]
\qquad 
$\implies$
\allowdisplaybreaks\begin{align*}
\frac{\partial^2}{\partial y^2} \hsx \abs{f(x + \sqrt{-1} \hsx y)}^2 \ 
&=\ 
\sum\limits_{n = 0}^\infty \ 
(2n + 2) \hsy
(2n + 1) \hsy 
L_{n + 1} (f) (x) \hsy y^{2 \hsy n}
\\[15pt]
&\geq \ 
0 
\qquad (\tcf. \ 13.7).
\end{align*}
On the other hand, 

\[
\frac{\partial^2}{\partial y^2} \hsx \abs{f(x + \sqrt{-1} \hsx y)}^2
\ = \ 
2 \abs{f^\prime (z)}^2  
\hsx - \hsx 
2 \Reg (f (z) \hsy \ov{f^{\prime\prime} (z)}).
\]
As for the converse, 
let $z_0 = x_0 + \sqrt{-1} \hsx y_0$ be a zero of $f$ and consider 

\[
f_0 (y) 
\ \equiv \ 
f_{x_0} (y) 
\ = \ 
\abs{f(x_0 + \sqrt{-1} \hsx y)}^2.
\]
Then 
\[
\frac{\td^2}{\td y^2}\hsx  f_0 (y) 
\ \geq \ 
0,
\]
so $f_0 (y)$ is a convex even function of $y$, thus has a unique minimum, 
which must be taken on at $y = 0$.  
But
\[
0 
\ = \ 
f (z_0) 
\ = \ 
f(x_0 + \sqrt{-1} \hsx y_0)
\ \implies \ 
y_0 
\ = \ 
0.
\]
Therefore the zeros of $f$ are real, hence 
$f \in 0 - \sL - \sP$ $(= \sL - \sP)$.
\\[-.25cm]
\end{x}

\begin{x}{\small\bf THEOREM} \ 
Let $f \in A - \sL - \sP$  (cf. 10.31) $-$then $f \in 0 - \sL - \sP$ $(= \sL - \sP)$ iff 
$\forall \ z = x + \sqrt{-1} \hsx y$ $(y \neq 0)$, 

\[
\frac{1}{y} \ 
\Img \big(- f^\prime (z) \ov{f (z)}\big)
\ \geq \ 
0.
\]

[This is a simple consequence of the canonical computation $\ldots$ .]
\\[-.5cm]
\end{x}

\[
\text{APPENDIX}
\]

Let 
$f \in \sL - \sP$ 
\ 
be transcendental. 
If 
$\hsx f (t_0) \hsx \neq \hsx 0 \hsx$ 
and 
$\hsx f^\prime (t_0) = 0$, 
then 
$f (t_0) \hsy f^{\prime\prime} (t_0) < 0$ 
(cf. 13.3), 
so 
$t_0$ is a simple zero of 
$f^\prime \in \sL - \sP$.  
\\

\qquad
{\small\bf LEMMA} \ 
Let $f \in \sL - \sP$ be transcendental.  
Suppose that $f^{(n)}$ has a multiple zero at $t_0$ $-$then

\[
f (t_0) 
\ = \ 
f^\prime (t_0) 
\ = \ 
\cdots
\ = \ 
f^{(n)} (t_0) 
\ = \ 
0.
\]
\\[-1.25cm]

\qquad
{\small\bf SCHOLIUM} \ 
If the zeros of $f$ are simple, then the zeros of all its derivatives are simple.
\\[-.25cm]

\qquad
{\small\bf THEOREM} \ 
Let $f \in \sL - \sP$ be transcendental.  
Assume: $f$ satisfies the differential equation 
\[
f^{(n)} (z) 
\ = \ 
A (z) \hsy f (z),
\]
where $\restr{A}{\R}$ is real analytic $-$then the zeros of $f$ are simple.
\\[-.25cm]

PROOF \ 
Proceeding by contradiction, suppose that at some $t_0$, 
$f (t_0) = f^\prime (t_0)  = 0$, 
thus 
$f^{(n)} (t_0) = 0$.
Since
\[
f^{(n+1)} (z) 
\ = \ 
A^\prime (z) \hsy f (z) + A (z) f^\prime (z),
\]
it follows that 
$f^{(n+1)} (t_0) = 0$.  
Owing now to the lemma, 
\[
f (t_0) 
\ = \ 
f^\prime (t_0) 
\ = \ 
\cdots
\ = \ 
f^{(n)} (t_0)
\ = \ 
f^{(n+1)} (t_0) 
\ = \ 
0.
\]
But 
\[
f^{(n+k)} (z) 
\ = \ 
\sum\limits_{\ell = 0}^k \ 
\binom{k}{\ell} \hsy 
A^{(k - \ell)} (z) \hsy
f^{(\ell)} (z).
\]
Therefore $f$ and all its derivatives vanish at $t_0$, a non sequitur.


\chapter{
$\boldsymbol{\S}$\textbf{14}.\quad  SHIFTED SUMS}
\setlength\parindent{2em}
\setcounter{theoremn}{0}
\renewcommand{\thepage}{\S14-\arabic{page}}

\qquad 
Let $f \not\equiv 0$ be a real entire function.
\\[-.5cm]

\begin{x}{\small\bf NOTATION} \ 
Given a real number $\lambda$, put
\[
f_\lambda (z) 
\ = \ 
f (z + \sqrt{-1} \hsx \lambda) 
\hsx + \hsx 
f (z - \sqrt{-1} \hsx \lambda).
\]

[Note: \ 
$f_\lambda$ is again a real entire function.]
\\[-.5cm]
\end{x}

Obviously,
\[
f_\lambda 
\ = \ f_{-\lambda}.
\]
\\[-1.25cm]

\begin{x}{\small\bf EXAMPLE} \ 
Take $f (z) = z^n$ $-$then 

\[
f_\lambda (z) 
\ = \ 
2 \ 
\prod\limits_{k = 0}^{n-1} \ 
\Big(
z - \lambda \cot \Big[\frac{(2 k + 1)\hsy \pi}{2 n}\Big]
\Big).
\]
\\[-.75cm]
\end{x}

\begin{x}{\small\bf EXAMPLE} \ 
Take 
$
f (z) \ = \ 
\begin{cases}
\ 
\sin z
\\[4pt]
\ 
\cos z
\end{cases}
$
$-$then
\[
f_\lambda (z) \ = \ 
2 \cosh \lambda \ 
\begin{cases}
\ 
\sin z
\\[4pt]
\ 
\cos z
\end{cases}
.
\]
\\[-1.25cm]
\end{x}

Let $EX_f$ denote the set of $\lambda$ such that $f_\lambda \equiv 0$ or for which $f_\lambda$ 
has the form 
$C_\lambda \hsy \exp (b_\lambda z)$, 
where 
$C_\lambda \neq 0$ and $b_\lambda$ are real constants.
\\

\begin{x}{\small\bf LEMMA} \ 
Suppose that $f$ is not of the form $C \hsy e^{b \hsy z}$, where $C \neq 0$ and $b$ are real constants 
$-$then $EX_f$ is a discrete subset of $\R$ (if not empty).  
\\[-.5cm]

[In fact, 
\[
EX_f
\ = \ 
\{\lambda : L_1 (f_\lambda) \equiv 0\}.]
\]
\\[-1.25cm]
\end{x}

\begin{x}{\small\bf EXAMPLE} \ 
Take $f (z) = e^z$ $-$then 

\[
f_\lambda (z) 
\ = \ 
2 (\cos \lambda) \hsy e^z, 
\]
so $EX_f = \R$.
\\[-.5cm]

[Note: \ 
$f$ is in $\sL - \sP$ but technically the zero function 
(e.g., $f_\frac{\pi}{2}$) is not in $\sL - \sP$.]
\\[-.25cm]
\end{x}

\begin{x}{\small\bf EXAMPLE} \ 
Take $f (z) = e^z (a_0 + a_1 z)$, where $a_0$ and $a_1 \neq 0$ are real $-$then 
\[
f_\lambda (z) 
\ = \ 
e^z (A_1 z + A_0), 
\]
where
\[
A_1 
\ = \ 
2 \hsy a_1 \hsy \cos \lambda
\]
and
\[
A_0 
\ = \ 
2 \hsy a_0 \hsy \cos \lambda - 2 \hsy a_1\hsy  \lambda \hsy \sin \lambda.
\]
Therefore
\[
EX_f 
\ = \ 
\{(2 k + 1) \hsx \frac{\pi}{2} : k = 0, \pm 1, \ldots\}.
\]
And
\allowdisplaybreaks
\allowdisplaybreaks\begin{align*}
\lambda \in EX_f \ (\lambda \neq 0) 
&\implies 
A_0 = - 2 \hsy a_1 \hsy \lambda \hsy \sin \lambda \neq 0
\\[11pt]
&\implies 
f_\lambda \not\equiv 0.
\end{align*}
\\[-1cm]
\end{x}

\begin{x}{\small\bf EXAMPLE} \ 
Take 
\[
f (z) 
\ = \  
e^{b z}\hsy p(z)
\qquad (\text{$b$ real}), 
\]
where
\[
p (z) 
\ = \ 
a_0 + a_1 z + \cdots + a_n z^n 
\qquad (a_n \neq 0)
\]
is a real polynomial of degree $n \geq 2$ with real zeros only $-$then 
\[
f_\lambda (z)
\ = \ 
e^{b \hsy z}
\big(
A_n z^n + A_{n-1} z^{n-1} + \cdots + A_0
\big).
\]
Here
\[
A_n 
\ = \ 
2 \hsy a_n \hsy \cos \lambda b 
\]
and
\[
A_{n-1} 
\ = \ 
2 \hsy a_{n-1} \hsy \cos \lambda b 
\hsx - \hsx
2 \hsy \lambda \hsy n \hsy a_n \hsy \sin \lambda b.
\]

\qquad \textbullet \quad
If $\cos \lambda \hsy b \neq 0$, then $A_n \neq 0$ and $f_\lambda$ has $n$ zeros.
\\[-.25cm]

\qquad \textbullet \quad
If $\cos \lambda \hsy b = 0$, then $A_n = 0$ but if in addition $\lambda \neq 0$, 
then $A_{n-1} \neq 0$, thus $f_\lambda$ has $n - 1$ zeros.
\\[-.5cm]

Since $n \geq 2$, the conclusion is that $EX_f = \emptyset$.
\\[-.25cm]
\end{x}

\begin{x}{\small\bf REMARK} \ 
It is clear that if $\forall \ \lambda$, $f_\lambda\not\equiv 0$ has a zero, 
then $EX_f = \emptyset$.
\\[-.5cm]

[For instance, if $f \in \sL - \sP$ and if 

\[
f (z) 
\ = \ 
C \hsy z^m \hsy e^{b \hsy z} \ 
\prod\limits_{n = 1}^\infty\ 
\Big(1 - \frac{z}{\lambda_n}\Big) \hsx 
e^{z / \lambda_n}
\qquad (\tcf. \ 10.19)
\]
has an infinite number of zeros, then $\forall \ \lambda$, $f_\lambda \not\equiv 0$ 
has an infinite number of zeros, hence $EX_f = \emptyset$.]
\\[-.25cm]
\end{x}

\begin{x}{\small\bf LEMMA} \ 
If $f \in \sL - \sP$, then $\forall \ \lambda \in \R$, 
either $f_\lambda \in \sL - \sP$ or $f_\lambda \equiv 0$.
\\[-.5cm]

PROOF  \ 
By the usual approximation argument, it will be enough to consider the case when $f$ is a 
real polynomial with real zeros only, say
\[
f (z) 
\ = \ 
C \hsy z^m \ 
\prod\limits_{n = 1}^N \ 
\Big(1 - \frac{z}{\lambda_n}\Big)
\qquad (C \neq 0).
\]
So take $\lambda > 0$ and suppose that $f_\lambda (z) = 0$ $(z = x + \sqrt{-1} \hsx y)$ $-$then

\[
\abs{f (z+ \sqrt{-1} \hsx \lambda)}
\ = \ 
\abs{f (z- \sqrt{-1} \hsx \lambda)}
\] 
$\implies$
\allowdisplaybreaks
\allowdisplaybreaks\begin{align*}
1 \ 
&=\ 
\frac
{\abs{f (z + \sqrt{-1} \hsx \lambda)}^2}
{\abs{f (z - \sqrt{-1} \hsx \lambda)}^2}
\\[15pt]
&=\ 
\frac
{\abs{(z + \sqrt{-1} \hsx \lambda)^2}^m}
{\abs{(z - \sqrt{-1} \hsx \lambda)^2}^m}
\ \cdot \ 
\frac
{
\ds
\prod\limits_{n = 1}^N \ 
\abs{\lambda_n - (z + \sqrt{-1} \hsx \lambda)}^2
}
{
\ds
\prod\limits_{n = 1}^N \ 
\abs{\lambda_n - (z - \sqrt{-1} \hsx \lambda)}^2
}
\\[15pt]
&=\ 
\bigg[
\frac
{x^2 + (y + \lambda)^2}
{x^2 + (y -  \lambda)^2}
\bigg]^m
\ \cdot \ 
\prod\limits_{n = 1}^N \ 
\frac
{
(x - \lambda_n)^2 + (y + \lambda)^2
}
{
(x - \lambda_n)^2 + (y - \lambda)^2
}
\hsx .
\end{align*}
\\[-1cm]

\noindent
If $y > 0$, then all factors on the RHS are $> 1$, while if $y <0$, then all factors on the RHS are $< 1$.  
As this is impossible, it follows that $y = 0$.
\\[-.5cm]

[Note: \ 
More generally, the same argument can be used to show that the polynomial
\[
f (z + \sqrt{-1} \hsx \lambda) 
-
\gamma \hsy 
f (z - \sqrt{-1} \hsx \lambda)
\qquad (\gamma \in \Cx,  \ \abs{\gamma} = 1)
\]
has real zeros only.]
\\[-.25cm]
\end{x}

\qquad
{\small\bf \un{N.B.}} \ 
Consequently, $\forall \ \lambda \in \R$, 
\[
f \in \sL - \sP
\implies 
L_1 (f_\lambda) (t) 
\ \geq \ 
0
\qquad (t \in \R) 
\quad \text{(cf. 13.3).}
\]
\\[-1cm]

\begin{x}{\small\bf EXAMPLE} \ 
Take $f (z) = z (1 + z^2)$ $-$then 
\[
L_1 (f_\lambda) (t) 
\ = \ 
12 t^4 + (6 \lambda^2 - 2)^2 
\ \geq \ 
0,
\]
yet $f \notin \sL - \sP$.
\\[-.5cm]

[Note: \ 
\[
L_1 (f_\lambda) (0) 
\ = \ 
(6 \lambda^2 - 2)^2
\]
and the expression on the right vanishes at 
$
\ds
\lambda
\ = \ 
\pm\frac{1}{\sqrt{3}}
\hsx$.]
\\
\end{x}

\begin{x}{\small\bf THEOREM} \ 
If $f \in \sL - \sP$ and if $EX_f = \emptyset$, then $\forall \ \lambda \neq 0$, 
the zeros of $f_\lambda$ are simple.
\\[-.5cm]

PROOF \ 
Take $\lambda > 0$ and suppose that $t_0$ is a multiple zero of $f_\lambda$: 

\[
\begin{cases}
\ 
f_\lambda (t_0) 
= 0 
\implies 
f (t_0 + \sqrt{-1} \hsx \lambda)
= 
-
f (t_0 - \sqrt{-1} \hsx \lambda)
\\[8pt]
\ 
f_\lambda^\prime (t_0)
= 0 
\implies 
f^\prime (t_0 - \sqrt{-1} \hsx \lambda)
=
-
f^\prime (t_0 + \sqrt{-1} \hsx \lambda)
\end{cases}
.
\]
Now 
\[
f (t_0 - \sqrt{-1} \hsx \lambda)
\hsy
f^\prime (t_0 + \sqrt{-1} \hsx \lambda)
\]
is real iff 
\[
f (t_0 - \sqrt{-1} \hsx \lambda)
\hsy
f^\prime (t_0 + \sqrt{-1} \hsx \lambda)
\ = \ 
\ov
{
f (t_0 - \sqrt{-1} \hsx \lambda)
\hsy
f^\prime (t_0 + \sqrt{-1} \hsx \lambda)
}\hsx .
\]
But
\allowdisplaybreaks
\allowdisplaybreaks\begin{align*}
\ov
{
f (t_0 - \sqrt{-1} \hsx \lambda)
\hsy
f^\prime (t_0 + \sqrt{-1} \hsx \lambda)
}
\ 
&=\ 
f (t_0 + \sqrt{-1} \hsx \lambda)
\hsy
f^\prime (t_0 - \sqrt{-1} \hsx \lambda)
\\[15pt]
&=\ 
(-
f (t_0 - \sqrt{-1} \hsx \lambda))
\hsx
(-f^\prime (t_0 + \sqrt{-1} \hsx \lambda)
)
\\[15pt]
&= 
f (t_0 - \sqrt{-1} \hsx \lambda)
\hsy
f^\prime (t_0 + \sqrt{-1} \hsx \lambda).
\end{align*}
On the other hand, for $\Img z > 0$, 
\allowdisplaybreaks
\allowdisplaybreaks\begin{align*}
\Img \frac{f^\prime (z)}{f (z)} \ 
&=\ 
\Img 
\Big(
\frac{m}{z} + 2 \hsy a \hsy z + b + \ 
\sum\limits_{n = 1}^\infty \ 
\Big(
\frac{1}{z - \lambda_n} + \frac{1}{\lambda_n}
\Big)
\Big)
\\[15pt]
&< \ 
0.
\end{align*}
Setting $z = t_0 + \sqrt{-1} \hsx \lambda$ then leads to a contradiction: \
\allowdisplaybreaks
\allowdisplaybreaks\begin{align*}
\Img \hsx \frac{f^\prime (t_0 + \sqrt{-1} \hsx \lambda)}{f (t_0 + \sqrt{-1} \hsx \lambda)} \ 
&=\ 
\Img \hsx
\frac
{f^\prime (t_0 + \sqrt{-1} \hsx \lambda) 
\hsy 
\ov{f (t_0 + \sqrt{-1} \hsx \lambda)}}
{\abs{f (t_0 + \sqrt{-1} \hsx \lambda)}^2}
\\[15pt]
&=\ 
\frac{1}{\abs{f (t_0 + \sqrt{-1} \hsx \lambda)}^2}
\
\Img \hsx
\big(
f^\prime (t_0 + \sqrt{-1} \hsx \lambda)
\hsy
f (t_0 - \sqrt{-1} \hsx \lambda)
\big)
\\[15pt]
&=\ 
0.
\end{align*}

[Note: \ 
This point is illustrated by 14.2 and 14.3.]
\\[-.25cm]
\end{x}

\begin{x}{\small\bf THEOREM} \ 
If $f \in \sL - \sP$ and if $EX_f = \emptyset$, then $\forall \ \lambda \neq 0$, 
\[
L_1 (f_\lambda) (t)\  > \ 0
\qquad (t \in \R) \quad \text{(cf. 13.3)}. 
\]
\\[-1.5cm]
\end{x}

\begin{x}{\small\bf REMARK} \ 
Suppose that $f \in A - \sL - \sP$ has the property that 
$\forall \ \lambda \neq 0$, 
\[
L_1 (f_\lambda) (t)\  > \ 0
\qquad (t \in \R) \quad \text{(cf. 13.3)}. 
\]
Then $EX_f = \emptyset$ and it is an open question as to whether $f \in \sL - \sP$.
\\[-.5cm]

[Note: \ 
If specialized to the case when $f \in * - \sL - \sP$, 
the stated condition does indeed imply that $f \in \sL - \sP$.  
In passing, observe that the strict inequality 
$L_1 (f_\lambda) (t) \hsx >  \hsx 0$ is necessary (cf. 14.10).]
\end{x}


\chapter{
$\boldsymbol{\S}$\textbf{15}.\quad  JENSEN CIRCLES (BIS)}
\setlength\parindent{2em}
\setcounter{theoremn}{0}
\renewcommand{\thepage}{\S15-\arabic{page}}

\qquad 
Given a real polynomial $f$, denote by 
$z_1, \ldots, z_\ell$ 
those zeros of $f$ which lie in the open upper half-plane.
\\[-.25cm]

\begin{x}{\small\bf NOTATION} \ 
Given a real polynomial $f$ and a real number $\lambda$, 
for $j = 1, \ldots, \ell$, put
\[
\gC_j (\lambda) 
\ = \ 
\{z \in \Cx : \abs{z - \Reg z_j}^2 \leq (\Img z_j)^2 - \lambda^2\}.
\]

[Note: \ 
Take 
$\gC_j (\lambda) = \emptyset$ 
if 
$\abs{\lambda} > \abs{\Img z_j}\hsx$.]
\\[-.25cm]
\end{x} 

\qquad
{\small\bf \un{N.B.}} \ 
In particular: 
\[
\gC_j (0) 
\ = \ 
\gC_j
\qquad (\tcf. \ 9.2).
\]
\\[-1.25cm]

\begin{x}{\small\bf THEOREM} \ 
For any $\lambda \neq 0$, the nonreal zeros of the polynomial
\[
f (z + \sqrt{-1} \hsx \lambda) 
- 
\gamma
f (z - \sqrt{-1} \hsx \lambda) 
\qquad (\gamma \in \Cx, \ \abs{\gamma} = 1)
\]
lie in the union of the $\gC_j (\lambda)$.
\\[-.25cm]

PROOF \ 
Take $f$ monic of degree $n$, so

\[
f (z) 
\ = \ 
\prod\limits_{\Img  z_i = 0} \ 
\big(z - z_i\big)^{m_i} 
\hsx \cdot \hsx
\prod\limits_{j = 1}^\ell \ 
\big(z - z_j\big)^{m_j} 
\big(z - \bar{z}_j\big)^{m_j} 
\qquad (\tcf. \ 9.3).
\]
Write
\[
z = x + \sqrt{-1} \hsx y
\quad \text{and} \quad 
z_j = x_j + \sqrt{-1} \hsx y_j
\qquad (j = 1, \ldots, \ell).
\]
Then

\qquad \textbullet \quad
$
\abs{z + \sqrt{-1} \hsx \lambda - z_i}^2
-
\abs{z - \sqrt{-1} \hsx \lambda - z_i}^2
\ = \ 
4 \hsy \lambda \hsy y
\qquad (\Img z_i = 0).
$
\\[-.25cm]


\qquad \textbullet \quad
$
\abs{z + \sqrt{-1} \hsx \lambda - z_j}^2
\hsx
\abs{z + \sqrt{-1} \hsx \lambda - \bar{z}_j}^2
$
\\[-.25cm]

\hspace{6cm}
$
\hsx - \hsx
\abs{z - \sqrt{-1} \hsx \lambda - z_j}^2
\hsx
\abs{z - \sqrt{-1} \hsx \lambda - \bar{z}_j}^2
$
\\[-.25cm]

\hspace{3cm}
$
\ = \ 
8 \hsy \lambda \hsy y
\hsx
\big[
(x - x_j)^2 + y^2 + \lambda^2 - y_j^2
\big].
$
\\[-.25cm]

\noindent
If now $z$ is nonreal and lies outside all the $\gC_j (\lambda)$, then

\[
(x - x_j)^2 + y^2 + \lambda^2 - y_j^2 \ > \ 0.
\]
Therefore every factor in the product representation of 
$\abs{f (z + \sqrt{-1} \hsx \lambda)}^2$ 
is larger than the corresponding factor in the product representation of 
$\abs{f (z - \sqrt{-1} \hsx \lambda)}^2$ 
if $\lambda \hsy  y > 0$ and vice-versa if $\lambda \hsy y <  0$.  
To recapitulate: \ 
\[
\begin{cases}
\ 
\lambda \hsy y \hsx > \hsx 0 
\ \implies \ 
\abs{ f (z + \sqrt{-1} \hsx \lambda)} 
\hsx > \hsx
\abs{ f (z - \sqrt{-1} \hsx \lambda)} 
\\[4pt]
\ 
\lambda \hsy y \hsx < \hsx 0 
\ \implies \ 
\abs{ f (z + \sqrt{-1} \hsx \lambda)} 
\hsx < \hsx
\abs{ f (z - \sqrt{-1} \hsx \lambda)} 
\end{cases}
.
\]
Accordingly, at such a $z$, the polynomial
\[
f (z + \sqrt{-1} \hsx \lambda)
\hsx - \hsx 
\gamma \hsy
f (z - \sqrt{-1} \hsx \lambda)
\]
cannot vanish.
\\[-.25cm]
\end{x}

\qquad
{\small\bf \un{N.B.}} \ 
If \ 
$\abs{\lambda} = \abs{\Img z_j} = \abs{y_j}$, then 
\[
\gC_j (\lambda) 
\ = \ 
\{z \in \Cx : (x - x_j)^2 + y^2 
\ \leq \ 
y_j^2 - \lambda^2  
\ = \ 
0\},
\]
so in this situation, $x = x_j$ and $y = 0$, thus
\[
\gC_j (\lambda) 
\ = \ 
\{(x_j, 0)\}.
\]
\\[-1.25cm]

\begin{x}{\small\bf COROLLARY} \ 
For any $\lambda \neq 0$, the nonreal zeros of the polynomial

\[
f_\lambda (z) 
\ = \ 
f (z + \sqrt{-1} \hsx \lambda)
\hsx + \hsx 
f (z - \sqrt{-1} \hsx \lambda)
\]
lie in the union of the 
$\gC_j (\lambda)$.
\\[-.5cm]

[Simply take $\gamma = -1$.]
\\[-.25cm]
\end{x}

\begin{x}{\small\bf COROLLARY} \ 
For any $\lambda \neq 0$ and for any $\xi \in \Cx$ $(\xi \neq 0)$, the nonreal zeros of the polynomial

\[
\xi \hsy
f (z + \sqrt{-1} \hsx \lambda)
\hsx + \hsx 
\bar{\xi} \hsy
f (z - \sqrt{-1} \hsx \lambda)
\]
lie in the union of the 
$\gC_j (\lambda)$.
\\[-.25cm]

[Simply take $\ds\gamma = -\frac{\bar{\xi}}{\xi}$.]
\\[-.25cm]
\end{x}

\begin{x}{\small\bf REMARK} \ 
One can recover 9.3 from 15.2.  
Thus let 
$
\ds 
\lambda_n = \frac{1}{n}
$ 
and consider
\[
f_n (z) 
\ = \ 
\frac
{
f (z + \sqrt{-1} \hsx \lambda_n)
\hsx - \hsx 
f (z - \sqrt{-1} \hsx \lambda_n)
}
{2 \hsy \lambda_n }.
\]
Then

\[
\lim\limits_{n \ra \infty} \ 
f_n (z)  
\ = \ 
f^\prime (z)
\]
uniformly on compact subsets of $\Cx$.  
Moreover, the zeros of $f_n (z)$ are contained in the union of the 
$\gC_j (\lambda_n)$ 
and the real line which is a subset of the union of the Jensen circles of $f$ and the real line.
\\[-.25cm]
\end{x}

\begin{x}{\small\bf LEMMA} \ 
Let $f$ be a real polynomial whose zeros lie in the strip 

\[
S (A) 
\ = \ 
\{z : \abs{\Img z} \leq A\}
\qquad (A > 0).
\]
Then for all 
$\lambda \neq 0$, 
the zeros of the polynomial 

\[
f (z + \sqrt{-1} \hsx \lambda)
\hsx - \hsx 
\gamma \hsy
f (z - \sqrt{-1} \hsx \lambda)
\qquad (\gamma \in \Cx, \ \abs{\gamma} = 1)
\]
lie in 
$S \big( \sqrt{A^2 - \lambda^2}\hsx \big)$
if 
$\abs{\lambda} < A$ 
and lie in $S(0) = \R$ 
if 
$A \leq \abs{\lambda}$.
\\[-.5cm]

PROOF \ 
If $z = x + \sqrt{-1} \hsx y \in \gC_j(\lambda)$ is a nonreal zero and if $\abs{\lambda} < A$ , then 
\[
y^2 
\ \leq \ 
(x - x_j)^2 + y^2 
\ \leq \ 
y_j^2 - \lambda^2
\ \leq \ 
A^2 - \lambda^2,
\]
hence 
$z \in S \big( \sqrt{A^2 - \lambda^2} \hsx \big)$.  
Meanwhile, at the transition point 
$A = \abs{\lambda}$, 
there is no nonreal zero in any of the 
$\gC_j(\lambda)$ and on the other side
$A < \abs{\lambda}$, 
all the $\gC_j(\lambda)$ are empty.
\\[-.25cm]
\end{x}

\begin{x}{\small\bf REMARK} \ 
If $A = 0$, hence if $f \in \sL - \sP$, then 
$\forall \ \lambda \neq 0$, 
the zeros of the polynomial
\[
f (z + \sqrt{-1} \hsx \lambda)
\hsx - \hsx 
\gamma \hsy
f (z - \sqrt{-1} \hsx \lambda)
\qquad (\gamma \in \Cx, \ \abs{\gamma} = 1)
\]
are real (cf. 14.9) and this persists to $\lambda = 0$:
\[
f (z) - \gamma \hsy f (z) 
\ = \ 
(1 - \gamma) \hsy f (z).
\]
\\[-1.25cm]
\end{x}

\begin{x}{\small\bf THEOREM} \ 
Let $f \in A - \sL - \sP$ (cf. 10.31) $-$then the zeros of $f_\lambda$ lie in 
$S \big( \sqrt{A^2 - \lambda^2} \hsx \big)$ 
if 
$\abs{\lambda} < A$ 
and lie in 
$S (0) = \R$ if $A \leq \abs{\lambda}$.
\\[-.5cm]

[Taking into account 15.6 and 15.7, apply 10.32.]
\\[-.5cm]

[Note: \ 
It is a corollary that 
\[
f_\lambda \in A_\lambda - \sL - \sP,
\]
where
\[
A_\lambda
\ = \ 
\big(\max (A^2 - \lambda^2, 0)\big)^{1/2}.]
\]
\\[-.25cm]
\end{x}


\chapter{
$\boldsymbol{\S}$\textbf{16}.\quad  STURM CHAINS}
\setlength\parindent{2em}
\setcounter{theoremn}{0}
\renewcommand{\thepage}{\S16-\arabic{page}}

\qquad
Given nonconstant real polynomials $P$ and $Q$, put
\[
F (z) 
\ = \ 
P (z) + \sqrt{-1} \hsx Q(z).
\]

\begin{x}{\small\bf LEMMA} \ 
Suppose that $F (z)$ has all its zeros in either the open upper half-plane or the open lower half-plane 
$-$then $P$ and $Q$ have real zeros only.  
\\[-.5cm]

PROOF \ 
Working under the open lower half-plane supposition, write
\[
F (z) 
\ = \ 
C_n (z - z_1) \cdots (z - z_n) 
\qquad (C_n \neq 0).
\]
Then for $\Img z > 0$, 
\[
\abs{z - z_k} 
\ > \ 
\abs{\bar{z} - z_k} 
\qquad (\Img z_k < 0, \ k = 1, \ldots, n)
\]
\qquad
$\implies$
\[
\abs{F (z)} 
\ > \ 
\abs{F (\bar{z})} 
\]
\qquad 
$\implies$
\allowdisplaybreaks\begin{align*}
2 \sqrt{-1} \hsx (P(\bar{z}) \hsy Q (z) \hsx - \hsx P (z) \hsy Q (\bar{z})) \ 
&=\ 
F (z) \hsy \ovs{F(z)} \hsx - \hsx F (\bar{z}) \hsy\ovs{F(\bar{z})}
\\[11pt]
&> \
0.
\end{align*}
Therefore $P$ and $Q$ have real zeros only (nonreal zeros of either $P$ or $Q$ would occur in conjugate pairs).
\\[-.5cm]

[Note: \ 
$P$ and $Q$ have no common zero (otherwise $F$ would have a real zero: \ 
$\abs{F (x)}^2 = P (x)^2 + Q (x)^2$).]
\\[-.25cm]
\end{x}

Here is an application.  
Let $f$ be a nonconstant real polynomial with real
zeros only, so $f \in \sL - \sP$, thus taking $\lambda > 0$, the zeros of 
$f (z + \sqrt{-1} \hsx\lambda)$ lie in the  open lower half-plane.  
Define nonconstant real polynomials $P$ and $Q$ by writing 
\[
f (z + \sqrt{-1} \hsx \lambda) 
\ = \ 
P (z) + \sqrt{-1} \hsx Q (z).
\]
Then $P, \hsy Q \in \sL - \sP$ and $\forall \ x \in \R$, 
\allowdisplaybreaks\begin{align*}
f_\lambda (x) \ 
&=\ 
f (x + \sqrt{-1} \hsx \lambda) 
-
\ovs{f (x + \sqrt{-1} \hsx \lambda) }
\\[11pt]
&=\ 
2 \hsy P (x)
\end{align*}
\qquad \qquad
$\implies$
\[
f_\lambda \in \sL - \sP 
\qquad (\tcf. \ 14.9).
\]
\\[-.75cm]

\begin{x}{\small\bf REMARK} \ 
If $\mu$ and $\nu$ are real and if $\mu^2 + \nu^2 > 0$, then the zeros of $F$ and
\[
(\mu - \sqrt{-1} \hsx \nu)F 
\ = \ 
(\mu \hsy P + \nu \hsy Q) + \sqrt{-1} \hsx  (\mu \hsy Q - \nu \hsy P)
\]
are the same.  
Therefore
\[
\begin{cases}
\ 
\mu \hsy P + \nu \hsy Q
\\[4pt]
\ 
\mu \hsy Q - \nu \hsy P
\end{cases}
\]
have real zeros only.
\\[-.25cm]
\end{x}

\begin{x}{\small\bf SUBLEMMA} \ 
The zeros of 
\[
\Big(
1 + \frac{\sqrt{-1} \hsx \lambda \hsy z}{n}
\Big)^n
\qquad (\lambda > 0)
\]
lie in the open upper half-plane, hence the zeros of 
\[
1 
- 
\binom{n}{2} \hsx \frac{\lambda^2 \hsy z^2}{n^2} 
+ 
\binom{n}{4} \hsx \frac{\lambda^4 \hsy z^4}{n^4} 
- \cdots
\]
are real (cf. 16.1).
\\[-.25cm]
\end{x}

\begin{x}{\small\bf LEMMA} \ 
Let $f$ be a real polynomial $-$then $f_\lambda$ has at least as many real zeros as $f$ does.
\\[-.5cm]

PROOF \ 
Take $\lambda > 0$ $-$then the polynomial 
\[
f (z) 
- 
\binom{n}{2} \hsx \frac{\lambda^2}{n^2} \hsx f^{\prime\prime} (z) 
+
\binom{n}{4} \hsx \frac{\lambda^4}{n^4} \hsx f^{\prime\prime\prime\prime} (z) 
-
\cdots
\]
has at least as many zeros as $f(z)$ does (cf. 12.10).  
But there is an expansion
\[
\frac{f_\lambda (z)}{2}
\ = \ 
f (z)
- 
\frac{\lambda^2}{2 !} \hsx f^{\prime\prime} (z) 
+
\frac{\lambda^4}{4 !} \hsx f^{\prime\prime\prime\prime} (z) 
-
\cdots, 
\]
so it remains only to let $n \ra \infty$.
\\[-.25cm]
\end{x}

\begin{x}{\small\bf LEMMA} \ 
Assume: \ 
\\[-.25cm]

\qquad \textbullet \quad
$F (z)$ has $n$ zeros in the closed lower half-plane
\\[-.5cm]

\noindent
or

\qquad \textbullet \quad
$F (z)$ has $n$ zeros in the closed upper half-plane.
\\[-.5cm]

\noindent
Then $P$ and $Q$ have $n$ pairs of nonreal zeros at most.
\\[-.5cm]

[Note: \ 
The case $n = 0$ is 16.1.]
\\[-.25cm]
\end{x}

There is more to be said about $(P, Q)$ and $F$ but for this it will be best to first introduce some machinery.
\\[-.25cm]

Let 
\[
P_n (x), P_{n-1} (x), \ldots, P_1 (x), P_0 (x)
\]
be a sequence of real polynomials such that 
$\deg P_k = k$ and $P_k^{(k)} (0) > 0$ $(k = 0, \ldots, n)$. 
\\[-.5cm]

[Note: \ 
Therefore $P_0 (x)$ is a positive constant.]
\\

\begin{x}{\small\bf DEFINITION} \ 
The $P_k$ are a \un{Sturm chain} if the following conditions are satisfied.
\\[-.5cm]


\qquad \textbullet \quad
Two consecutive terms $P_k$, $P_{k+1}$ cannot vanish simultaneously.  
\\[-.25cm]

\qquad \textbullet \quad
Whenever one of the $P_{n-1}, \ldots, P_1$ vanishes, the neighboring terms have opposite signs.
\\[-.25cm]
\end{x}

\begin{x}{\small\bf EXAMPLE} \ 
Consider the Legendre polynomials

\[
P_n (x) 
\ = \ 
\frac{1}{2^n \hsy n !} 
\hsx 
\frac{\td^n}{\td x^n} 
\hsx 
(x^2 - 1)^n 
\qquad (\tcf. \ 8.17).
\]
Then
\[
P_0 (x) = 1, 
\quad 
P_1 (x) = x, 
\quad 
P_2 (x) = \frac{3}{2} \hsy x^2 - \frac{1}{2},
\]
and for $k > 2$, 
\[
P_k (x) 
\ = \ 
\frac{2^k \hsy \Big(\ds\frac{1}{2}\Big) \hsy k}{k !} 
\hsx
x^k 
\hsx + \hsx 
\pi_{k-2} (x),
\]
where $\pi_{k-2}$ is a polynomial of degree $(k - 2)$ in $x$.  
Furthermore, there is a recurrence relation 
\[
(k+1) P_{k+1} (x) 
\ = \ 
(2k + 1) \hsy x \hsy P_k (x) - k \hsy P_{k-1} (x).
\]
Thus, in consequence, the sequence
\[
P_n (x), P_{n-1} (x), \ldots, P_1 (x), P_0 (x)
\]
is a Sturm chain.
\\[-.5cm]

[Note: \ 
This setup is the tip of the iceberg: \ 
Consider a weight function $w(x) > 0$ $(a < x < b)$ ($a$ or $b$ potentially infinite) and an associated sequence 
$\{P_n (x)\}$ of orthogonal real polynomials.]
\\[-.25cm]
\end{x}

\begin{x}{\small\bf EXAMPLE} \ 
Fix $\lambda > -1$ and let

\[
P_{\lambda, n} (x) 
\ = \ 
\int\limits_{-1}^1 \ 
(1 - t^2)^\lambda 
\hsx 
(x + \sqrt{-1} \hsx t)^n 
\ \td t
\qquad (n = 0, 1, \ldots).
\]
Then the sequence 
\[
P_{\lambda, n} (x), P_{\lambda, n-1} (x), \ldots, P_{\lambda, 1} (x), P_{\lambda, 0} (x)
\]
is a Sturm chain.
\\[-.25cm]
\end{x}

\begin{x}{\small\bf STURM CRITERION} \ 
Suppose that 
\[
P_n (x), P_{n-1} (x), \ldots, P_1 (x), P_0 (x)
\]
is a Sturm chain $-$then the zeros of the $P_k$ $(k = 1, \ldots, n)$ are real and simple.  
\\[-.25cm]
\end{x}

Return now to 
\[
F (z) 
\ = \ 
P (z) +  \sqrt{-1} \hsx Q (z).
\]
\\[-1.25cm]

\begin{x}{\small\bf LEMMA} \ 
Under the assumptions of 16.1, $P$ and $Q$ have real zeros only and, in addition, these zeros are simple.
\\[-.5cm]

[Note: \ 
The new information is the assertion of simplicity.]
\\[-.25cm]
\end{x}

It suffices to work with $P$ (since
$- \sqrt{-1} \hsx F = Q - \sqrt{-1} \hsx P$), 
the idea being to exhibit a Sturm chain
\[
P (x) 
\ = \ 
P_n (x), P_{n-1} (x), \ldots, P_1 (x), P_0 (x), 
\]
thereby enabling one to quote 16.9.
\\[-.5cm]

As before, write
\[
F (z) 
\ = \ 
C_n (z - z_1) \cdots (z - z_n) 
\qquad (C_n \neq 0),
\]
take $C_n = 1$, and let 
\[
z_1 = a_1 + \sqrt{-1} \hsx b_1 \ (b_1 < 0), 
\ldots, 
z_n = a_n + \sqrt{-1} \hsx b_n \ (b_n < 0).
\]
Put
\allowdisplaybreaks\begin{align*}
F_k (x) \ 
&=\ 
(x - a_1 - \sqrt{-1} \hsx b_1) 
\cdots
(x - a_k - \sqrt{-1} \hsx b_k) 
\\[11pt]
&\equiv \ 
P_k (x) + \sqrt{-1} \hsx Q_k (x).
\end{align*}
Then
\[
\begin{cases}
\ 
P_k (x) 
\ = \ 
(x - a_k)\hsy P_{k-1} (x) + b_k \hsy Q_{k-1} (x)
\\[8pt]
\ 
Q_k (x)
\ = \ 
-b_k \hsy P_{k-1} (x) + (x - a_k) \hsy Q_{k-1} (x)
\end{cases}
.
\]
Replacing $k$ by $k + 1$ gives 
\[
P_{k+1} (x)
\ = \ 
(x - a_{k+1}) P_k (x) + b_{k+1} Q_k (x)
\]
from which (by elimination of $Q_k (x)$) 
\[
b_k \hsy P_{k+1} (x)
\ = \ 
(b_k \hsy (x - a_{k+1}) + b_{k+1} \hsy (x - a_k) ) P_k (x) 
- b_{k+1} (b_k^2 + (x - a_k)^2) P_{k-1} (x).
\]
Setting $P_0 (x) = 1$ and noting that by construction, the $P_k$ are monic, it thus follows that 
\[
P (x) 
\ = \ 
P_n (x), P_{n-1} (x), \ldots, P_1 (x), P_0 (x)
\] 
is a Sturm chain, as desired.
\\[-.5cm]

At this juncture, return to the inequality
\[
2 \sqrt{-1} \hsx \big(P (\bar{z}) \hsy Q(z) - P (z) \hsy Q (\bar{z})) 
\ > \ 
0
\qquad (\Img z > 0)
\]
and divide it by $- 2 \sqrt{-1} \hsx (z - \bar{z})$ to get 
\[
- 
\hsx
\frac
{P (\bar{z}) (Q (z) - Q (\bar{z})) 
\hsx - \hsx 
Q (\bar{z}) (P (z) - P (\bar{z}))}
{z - \bar{z}}
\ > \ 
0
\qquad (\Img z > 0).
\]
Letting $z$ approach the real axis, we conclude that 
\[
Q(x) \hsy P^\prime (x) - P (x) \hsy Q^\prime (x)
\ \geq \ 
0.
\]
\\[-1.25cm]

\begin{x}{\small\bf REMARK} \ 
Recall that $P$ and $Q$ have no common zeros, so if $P (x_0) = 0$,  
then $Q (x_0) \neq 0$.  
On the other hand, $x_0$ is simple (cf. 16.10), hence $P^\prime (x_0) \neq 0$.  
Therefore
\[
Q(x_0) \hsy P^\prime (x_0) - P (x_0) \hsy Q^\prime (x_0)
\ = \ 
Q(x_0) \hsy P^\prime (x_0)
\ > \ 
0.
\]
Accordingly, 
\[
Q(x) \hsy P^\prime (x) - P (x) \hsy Q^\prime (x)
\ \geq \ 
0
\]
whenever $P (x) = 0$ (and, analogously, whenever $Q (x) = 0$).
\\[-.25cm]
\end{x}

\begin{x}{\small\bf LEMMA} \ 
Between any two consecutive zeros of $Q$ there is one and only one zero of $P$ and 
between any two consecutive zeros of $P$ there is one and only one zero of $Q$, 
i.e., $P$ and $Q$ have \un{interlacing zeros}. 
\\[-.5cm]

PROOF \ 
The rational function 
\[
R (x) 
\ = \ 
\frac{P (x)}{Q (x)}
\]
has a nonnegative derivative at all $x$ except at the zeros of $Q (x)$.  
Moreover, between any two consecutive zeros of $Q(x)$, $R (x)$ climbs from $-\infty$ to $+\infty$ and, 
in so doing, determines a unique zero of $P (x)$.
\\[-.25cm]
\end{x}

\begin{x}{\small\bf REMARK} \ 
This property of the data forces an after the fact restriction on the degrees of $P$ and $Q$, viz.

\[
\deg P 
\ = \ 
\deg Q 
\quad \text{or} \quad
\begin{cases}
\ 
\deg P 
\ = \ 
\deg Q + 1
\\[4pt]
\ 
\deg Q 
\ = \ 
\deg P + 1
\end{cases}
.
\]
\\[-1.25cm]
\end{x}

The preceding considerations can be turned around.  
Spelled out, make the following assumptions.
\\[-.25cm]

\qquad \textbullet \quad
The zeros of $P$ and $Q$ are real and simple.
\\[-.25cm]

\qquad \textbullet \quad
The zeros of $P$ and $Q$ are interlacing.
\\[-.25cm]

\qquad \textbullet \quad
There exists an $x_0$ such that 
\[
Q(x_0) \hsy P^\prime (x_0) - P (x_0) \hsy Q^\prime (x_0)
\ > \ 
0.
\]
Then
\[
F (z) 
\ = \ 
P(z) + \sqrt{-1} \hsx Q (z)
\]
has all its  zeros in the open lower half-plane.
\\[-.5cm]

To begin with, it is clear that $P$ and $Q$ do not have a common zero 
(their zeros being interlacing), thus $F$ cannot have a real zero.  
Suppose, therefore, that $F (z_0) = 0$, where 
$z_0 = x_0 + \sqrt{-1} \hsx y_0$ $(y_0 \neq 0)$ $-$then 
\[
\frac{P (z_0)}{Q (z_0)} + \sqrt{-1} 
\ = \ 
0.
\]
Denoting by 
$a_1 < a_2 < \cdots < a_n$ the zeros of $Q$, pass to the decomposition
\[
\frac{P (z)}{Q (z)}
\ = \ 
A + \frac{A_1}{z - a_1} + \frac{A_2}{z - a_2} + \cdots + \frac{A_n}{z - a_n}, 
\]
where $A$ is a real constant and 
\[
A_k 
\ = \ 
\frac{P (a_k)}{Q^\prime (a_k)}
\qquad (k = 1, 2, \ldots, n).
\]
Here
\[
\begin{cases}
\ 
P (a_k) \hsy P (a_{k+1}) \ < \ 0
\\[8pt]
\ 
Q^\prime (a_k) \hsy Q^\prime (a_{k+1}) \ < \ 0
\end{cases}
, 
\]
so
\[
A_1, A_2, \ldots, A_n
\]
have one and the same sign.  
But
\[
-\sqrt{-1} 
\ = \ 
A 
\hsx + \hsx
\frac{A_1}{z_0 - a_1}
\hsx + \hsx
\frac{A_2}{z_0 - a_2}
\hsx + \hsx \cdots \hsx + \hsx
\frac{A_n}{z_0 - a_n}
\]
\qquad 
$\implies$
\[
-1
\ = \ 
-y_0 \ 
\sum\limits_{k = 1}^n \ 
\frac{A_k}{(x_0 - a_k)^2 + y_0^2}
\]
\qquad 
$\implies$
\[
1
\ = \ 
y_0 \ 
\sum\limits_{k = 1}^n \ 
\frac{A_k}{(x_0 - a_k)^2 + y_0^2} \hsx.
\]
\\[-.75cm]

\noindent
There are then two possibilities: \ 
All the $A_k$ are $> 0$, in which case $y_0$ is positive, 
or all the $A_k$ are negative, in which case $y_0$ is negative.  
And this means that $F (z)$ has all its zeros either in the open upper half-plane or the open lower half-plane.  
\\[-.5cm]

It remains to eliminate the first contingency.  
However, if it held, then, arguing as before, we would have
\[
Q(x) \hsy P^\prime (x) - P (x) \hsy Q^\prime (x)
\ \leq \ 
0, 
\]
contradicting the assumption that there exists an $x_0$ such that 
\[
Q(x_0) \hsy P^\prime (x_0) - P (x_0) \hsy Q^\prime (x_0)
\ > \ 
0.
\]

[Note: \ 
\allowdisplaybreaks\begin{align*}
\forall \ k, \hsx A_k < 0 
&\implies
\bigg(\frac{P (x)}{Q (x)}\bigg)^\prime  > 0 
\qquad (x \neq a_k)
\\[15pt]
&\implies
Q(x) \hsy P^\prime (x) - P (x) \hsy Q^\prime (x) > 0.]
\end{align*}

In summary: 
\[
F (z) 
\ = \ 
P (z) + \sqrt{-1} \hsx Q (z)
\]
has all its zeros in the open lower-half plane.
\\
\begin{x}{\small\bf REMARK} \ 
The developments in this $\S$  are known collectively as \un{Hermite-Bieler theory}
\\[-.25cm]
\end{x}

\chapter{
$\boldsymbol{\S}$\textbf{17}.\quad  EXPONENTIAL TYPE}
\setlength\parindent{2em}
\setcounter{theoremn}{0}
\renewcommand{\thepage}{\S17-\arabic{page}}

\qquad
Given an entire function
\[
f (z) 
\ = \ 
\sum\limits_{n = 0}^\infty \ 
c_n \hsy z^n, 
\]
put
\[
T(f) 
\ = \ 
\underset{r \ra \infty}{\limsupx} \ 
\frac{\log M(r; f)}{r}.
\]
\\[-1cm]

\begin{x}{\small\bf DEFINITION} \ 
$f$ is of \un{exponential type} if $T (f) < \infty$, in which case $T (f)$ is called the \un{exponential type} of $f$. 
\\[-.25cm]
\end{x}

\qquad
{\small\bf \un{N.B.}} \ 
$f$ is of exponential type iff there exists a positive constant $K$: 
\[
f (z) 
\ = \ 
\tO \big(e^{K \hsy \abs{z}}\big),
\]
the greatest lower bound of the set of $K$ for which the relation holds then being the exponential type of $f$. 
\\[-.25cm]

\begin{x}{\small\bf LEMMA} \ 
If $f$ is of exponential type, then its order $\rho (f)$ is $\leq 1$. 
\\[-.25cm]
\end{x}

\begin{x}{\small\bf LEMMA} \ 
If $f$ is of exponential type and if $T (f) > 0$, then its order $\rho (f)$ is $= 1$ and $T (f) = \tau (f)$.
\\[-.25cm]
\end{x}

\begin{x}{\small\bf LEMMA} \ 
If $f$ is of exponential type and if $T (f) = 0$, then there are two possibilities: \ 
$\rho (f) < 1$ 
or 
$\rho (f) = 1$ and $\tau (f) = 0$.
\\[-.25cm]
\end{x}

\begin{x}{\small\bf SCHOLIUM} \ 
The set of entire functions of exponential type is comprised of the entire functions of order $< 1$ and the entire functions of order 1 
and of finite type.
\\[-.25cm]
\end{x}


\begin{x}{\small\bf EXAMPLE} \ 
The entire function 
\[
\frac{\sin \sqrt{z}}{\sqrt{z}}
\]
is of order 
$
\ds
\frac{1}{2}
$.  
It is of type 1 but of exponential type 0.  
\\[-.25cm]
\end{x}

\begin{x}{\small\bf EXAMPLE} \ 
The entire function 
\[
\frac{1}{z \hsy \Gamma (z)}
\]
is of order 1 (cf. 5.13).  
However, it is of maximal type (cf. 5.22), hence is not of exponential type.
\\[-.25cm]
\end{x}

\begin{x}{\small\bf LEMMA} \ 
If $f$ is of exponential type, then $f^\prime$ is of exponential type and $T (f) = T(f^\prime)$ 
(cf. 2.25 and 3.7).
\\[-.25cm]
\end{x}

\begin{x}{\small\bf LEMMA} \ 
If $f$, $g$ are of exponential type and if 
$
\ 
\ds
\frac{f}{g}
\ 
$
is entire, then 
$
\ 
\ds
\frac{f}{g}
\ 
$
is of exponential type.
\\[-.5cm]

PROOF \ 
On general gounds, 
\allowdisplaybreaks\begin{align*}
\rho \Big(\frac{f}{g}\Big) \ 
&\leq \ 
\max(\rho (f),\rho (g))
\qquad (\tcf. \ 2.37)
\\[8pt]
&\leq \ 
\max (1, 1) 
\\[8pt]
&= \ 
1.
\end{align*}
There is nothing to prove if 
$
\ds
\rho \Big(\frac{f}{g}\Big)  < 1
$, 
so assume that 
$
\ds
\rho \Big(\frac{f}{g}\Big)  = 1
$ 
and distinguish two cases.
\\[-.25cm]

\un{Case 1:} \quad
$\rho (g) < 1$ $-$then $\rho (f) = 1$ 
\\[-.5cm]

\noindent
and
\[
\tau (f) 
\ = \ 
\tau \Big(g \cdot\frac{f}{g}\Big)
\ = \ 
\tau \Big(\frac{f}{g}\Big)
\qquad (\tcf. \ 3.14),
\]
thus
$
\ds
\frac{f}{g}
$
is of finite type. 
\\[-.25cm]

\un{Case 2:} \quad
$\rho (g) = 1$ $-$then $0 \leq \tau (g) < \infty$ 
and if 
$
\ds
\tau \Big(\frac{f}{g}\Big)  = \infty
$, 
it would follow that 
\[
\tau (f) 
\ = \ 
\tau \Big(g \cdot\frac{f}{g}\Big)
\ = \ 
\infty
\qquad (\tcf. \ 3.14),
\]
contradicting $0 \leq \tau (f) < \infty$. 
\\[-.25cm]
\end{x}

\begin{x}{\small\bf THEOREM} \ 
Suppose that $f$ is an entire function $-$then
\[
T (f) 
\ = \ 
\frac{1}{e} \ 
\underset{n \ra \infty}{\limsupx} \ 
n \hsy \abs{a_n}^{1/n} 
\qquad (\tcf. \ 3.6).
\]

[Note: \ 
In terms of the $\gamma_n$, 
\[
T (f) 
\ = \ 
\underset{n \ra \infty}{\limsupx} \ 
\abs{\gamma_n}^{1/n}.
\]

Proof: 
\allowdisplaybreaks\begin{align*}
\frac{1}{e} \ 
\underset{n \ra \infty}{\limsupx} \ 
n \hsy \abs{a_n}^{1/n} \ 
&=\ 
\frac{1}{e} \ 
\underset{n \ra \infty}{\limsupx} \ 
n \hsy \abs{\frac{\gamma_n}{n !}}^{1/n} 
\\[15pt]
&=\ 
\underset{n \ra \infty}{\limsupx} \ 
\bigg(
\frac{n^n \hsy e^{-n} \hsx \sqrt{2 \hsy \pi \hsy n}}{n!}
\bigg)^{1/n}
\ 
\frac{n}{e \big(n^n \hsy e^{-n} \hsx \sqrt{2 \hsy \pi \hsy n} \hsx \big)^{1/n}}
\ 
\abs{\gamma_n}^{1/n}
\\[15pt]
&=\ 
\underset{n \ra \infty}{\limsupx} \ 
\abs{\gamma_n}^{1/n}.]
\end{align*}
\\[-1.25cm]
\end{x}

\begin{x}{\small\bf APPLICATION} \ 
An entire function $f$ is of exponential type iff
\[
\underset{n \ra \infty}{\limsupx} \ 
n \hsy \abs{a_n}^{1/n}
\ < \ 
\infty.
\]
\\[-1.5cm]
\end{x}


\begin{x}{\small\bf NOTATION} \ 
$\sE_0$ is the set of entire functions of exponential type.
\\[-.25cm]
\end{x}

\begin{x}{\small\bf LEMMA} \ 
$\sE_0$ is a vector space.
\\[-.5cm]

PROOF \ 
Let 
\\[-.75cm]
\[
\begin{cases}
\ \ds
f (z) \ = \ 
\sum\limits_{n = 0}^\infty \ 
a_n \hsy z^n
\\[15pt]
\ \ds
g (z) \ = \ 
\sum\limits_{n = 0}^\infty \ 
b_n \hsy z^n
\end{cases}
\]
be elements of $\sE_0$ $-$then 
\allowdisplaybreaks\begin{align*}
\abs{a_n + b_n}^{1/n} \ 
&\leq \ 
(2 \max \big(\abs{a_n}, \abs{b_n}\big)^{1/n}
\\[11pt]
&\leq \ 
2^{1/n}  \big(\abs{a_n}^{1/n}\hsx + \hsx  \abs{b_n}^{1/n}\big)
\end{align*}
\qquad 
$\implies$
\allowdisplaybreaks\begin{align*}
\underset{n \ra \infty}{\limsupx} \ 
n \hsx \abs{a_n + b_n}^{1/n} \ 
&\leq \ 
\underset{n \ra \infty}{\limsupx} \ 
2^{1/n} 
\hsx
n 
\hsx 
\big(
\abs{a_n}^{1/n} 
\hsx + \hsx 
\abs{b_n}^{1/n}
\big)
\\[11pt]
&\leq \ 
\lim\limits_{n \ra \infty} \ 
2^{1/n} 
\hsx \cdot \hsx
\underset{n \ra \infty}{\limsupx} \ 
n \hsy 
\big(
\abs{a_n}^{1/n} 
\hsx + \hsx 
n \hsy \abs{b_n}
\big)^{1/n}
\\[11pt]
&\leq \ 
\underset{n \ra \infty}{\limsupx} \ 
n \hsx \abs{a_n}^{1/n} 
\hsx + \hsx 
\underset{n \ra \infty}{\limsupx} \ 
n \hsx \abs{b_n}^{1/n} 
\\[11pt]
&< \ 
\infty.
\end{align*}
\\[-1.5cm]
\end{x}

\begin{x}{\small\bf EXAMPLE} \ 
A trigonometric polynomial
\[
\sum\limits_{k \hsy = \hsy  -n}^n \ 
c_k 
\hsx 
e^{\sqrt{-1} \hsx k \hsy z}
\]
is an entire function of exponential type $n$. 
\\[-.25cm]
\end{x}

\begin{x}{\small\bf LEMMA} \ 
$\sE_0$ is an algebra.
\\[-.5cm]

PROOF \ 
Given 

\[
\begin{cases}
\ 
f \in \sE_0
\\[4pt]
\ 
g \in \sE_0
\end{cases}
,
\]
choose positive constants
\\[-.5cm]
\[
\begin{cases}
\ 
(K, M)
\\[4pt]
\ 
(L, N)
\end{cases}
:  \ 
\begin{cases}
\ 
\abs{f (z)} \ \leq \  M \hsy e^{K \hsy \abs{z}}
\\[6pt]
\ 
\abs{g (z)} \ \leq \ N \hsy e^{L \hsy \abs{z}}
\end{cases}
.
\]
Then

\[
\abs{f (z) \hsy g(z)} 
\ \leq \ 
M  N \hsy e^{(K + L) \hsy \abs{z}}.
\]
\\[-1.25cm]
\end{x}

\begin{x}{\small\bf LEMMA} \ 
$\sE_0$ is closed under translation: 
If $f (z)$ is of exponential type $T (f)$ and if $A$, $B$ are complex constants, 
then 
$f (A z + B)$ is of exponential type $\abs{A} T (f)$.  
\\[-.25cm]
\end{x}

Embedded in the theory are a variety of estimates, a sampling of the simplest of these being given below..
\\[-.25cm]

\begin{x}{\small\bf LEMMA} \ 
Let $f \in \sE_0$, say
\[
\abs{f (z)} 
\ \leq \ 
C_K \hsy e^{K \hsy\abs{z}}.
\]
Assume: \ 
$\forall$ real $x$, 
\[
\abs{f (x)} 
\ \leq \
M.
\]
Then $\forall$ real $y$, 

\[
\abs{f (x + \sqrt{-1} \hsx y)}
\ \leq \ 
M \hsy e^{K \hsy\abs{y}}.
\]

[This is a standard application of Phragm\'en-Lindel\"of \ldots $\hsx.$]
\\[-.25cm]
\end{x}

\begin{x}{\small\bf THEOREM} \ 
Let $f \in \sE_0$.  
Assume: \ 
$\forall$ real $x$, 
\[
\abs{f (x)} 
\ \leq \ 
M.
\]
Then $\forall$ real $y$, 
\[
\abs{f (x + \sqrt{-1} \hsx y)}
\ \leq \ 
M \hsy e^{T (f) \hsy \abs{y}}.
\]

PROOF \ 
Given $\varepsilon > 0$, $\exists \ C_\varepsilon > 0 \hsy$: 
\[
\abs{f (z)}
\ \leq \ 
C_\varepsilon \hsy \exp((T (f) + \varepsilon) \abs{z}).
\]
So,  $\forall$ real $y$,
\[
\abs{f (x + \sqrt{-1} \hsx y)}
\ \leq \ 
M \hsy \exp((T (f) + \varepsilon) \abs{y}).
\]
Now let $\varepsilon \ra 0$: 
\\[-.5cm]

\qquad \qquad
$\implies$
\[
\abs{f (x + \sqrt{-1} \hsx y)}
\ \leq \ 
M \hsy e^{T (f) \hsy \abs{y}}.
\]

[Note: \ 
Accordingly, if $T (f) = 0$, then $f$ is a constant.  
In particular: 
Every entire function of order less than one which is bounded on the real axis must be a constant.]
\\[-.25cm]
\end{x}

\begin{x}{\small\bf EXAMPLE} \ 
Given $\phi \in \Lp^1 [-A, A]$ $(0 < A < \infty)$, put
\[
f (z) 
\ = \ 
\frac{1}{\sqrt{2 \hsy \pi}} \ 
\int\limits_{-A}^A \ 
\phi (t) \hsy e^{\sqrt{-1} \hsx z \hsy t} 
\ \td t.
\]
Then $f (z)$ is entire and 
\allowdisplaybreaks\begin{align*}
\abs{f (z) } \ 
&\leq \ 
\frac{1}{\sqrt{2 \hsy \pi}} \ 
\int\limits_{-A}^A \ 
\abs{\phi (t)} 
\hsx 
e^{-y \hsy t}
\ \td t
\qquad (z = x + \sqrt{-1} \hsx y)
\\[15pt]
&\leq \ 
\frac{1}{\sqrt{2 \hsy \pi}} \ 
e^{A \hsy \abs{y}} \ 
\int\limits_{-A}^A \ 
\abs{\phi (t)} 
\ \td t
\end{align*}
\qquad 
$\implies$
\[
T (f) 
\ \leq \ 
A,
\]
thus $f (z)$ is of exponential type.  
And: 
\allowdisplaybreaks
\allowdisplaybreaks\begin{align*}
\abs{f (x) } \ 
&\leq \ 
\frac{1}{\sqrt{2 \hsy \pi}} \ 
\int\limits_{-A}^A \ 
\abs{\phi (t)} 
\ \td t
\\[15pt]
&\equiv \ 
M,
\end{align*}
thereby realizing the assumption of 17.18.
\\[-.25cm]
\end{x}

\begin{x}{\small\bf LEMMA} \ 
Let $f \in \sE_0$.  Suppose that 
\[
f (x) \ra 0 
\quad \text{as} \quad 
\abs{x} \ra \infty.
\]
Then
\[
f (x + \sqrt{-1} \hsx y) \ra 0 
\quad \text{as} \quad 
\abs{x} \ra \infty
\]
uniformly in every horizontal strip.
\\[-.5cm]

[On the basis of the foregoing, this follows from Montel's theorem.]
\\[-.25cm]
\end{x}

\begin{x}{\small\bf EXAMPLE} \ 
Take the data as in 17.19 $-$then by the Riemann-Lebesgue lemma (cf. 21.6), 
\[
f (x) \ra 0 
\quad \text{as} \quad 
\abs{x} \ra \infty.
\]
\\[-1.25cm]
\end{x}

\begin{x}{\small\bf LEMMA} \ 
Let $f \in \sE_0$ with $T (f) > 0$.  
Assume: 
$\forall$ real $x$, 
\[
\abs{f (x)} 
\ \leq \ 
M.
\]
Then
\[
f^\prime (x) 
\ = \ 
\frac{4 \hsy T (f)}{\pi^2} \ 
\sum\limits_{k \hsy = \hsy -\infty}^\infty \ 
(-1)^k 
\hsx 
\frac{1}{(2 k + 1)^2}
\hsx
f \Big(x + \frac{2k + 1}{2 T (f)} \hsx \pi\Big), 
\]
the convergence being uniform on compact subsets of $\R$.
\\[-.5cm]

PROOF \ 
Suppose initially that $T (f) = 1$ and consider the meromorphic function
\[
F (z) 
\ = \ 
\frac{f (z)}{z^2 \hsy \cos z}.
\]
Let $\Gamma_n$ be the square contour with corners at 
$(1 + \sqrt{-1}) \hsx \pi \hsy n$, 
$(-1 + \sqrt{-1}) \hsx \pi \hsy n$, 
$(-1 - \sqrt{-1}) \hsx \pi \hsy n$, 
$(1 - \sqrt{-1}) \hsx \pi \hsy n$, 
$-$then $F$ has no singularities on $\Gamma_n$ but inside $\Gamma_n$ it might have a pole at the origin or at the points
$
\ds
\frac{2k + 1}{2} \hsx \pi$ 
$(-n \leq  n - 1)$.
So, from residue theory,
\[
\frac{1}{2 \hsy \pi \sqrt{-1}} \ 
\int\limits_{\Gamma_n} \ 
F (z) 
\td z
\ = \ 
f^\prime (0) 
\hsx - \hsx 
\sum\limits_{k \hsy = \hsy -n}^{n-1} \ 
(-1)^k \ 
\frac{4}{\pi^2 \hsy (2k + 1)^2} 
\hsx 
f \Big(\frac{2k + 1}{2} \hsx \pi\Big).
\]
Next
\[
z \in \Gamma_n 
\implies
\abs{\cos z} 
\ >  \ 
\frac{e^{\abs{y}}}{4} 
\qquad (y = \Img z).
\]
Meanwhile (cf. 17.18), 
\[
\abs{f(x + \sqrt{-1} \hsx y)}
\ \leq \ 
M \hsy e^{\abs{y}}
\qquad (T (f) = 1).
\]
Therefore
\allowdisplaybreaks\begin{align*}
z \in \Gamma_n \implies 
\abs{F (z)} \ 
&= 
\frac{\abs{f (z)}}{\abs{z^2 \cos z}}
\\[11pt]
&< \ 
4 \hsy M \hsy \abs{z}^{-2}
\end{align*}
\qquad 
$\implies$
\[
\int\limits_{\Gamma_n} \ 
F (z) 
\ \td z 
\ra 0 
\qquad (n \ra \infty)
\]
\qquad 
$\implies$

\[
f^\prime (0) 
\ = \
\frac{4}{\pi^2} \ 
\sum\limits_{k \hsy = \hsy -\infty}^\infty \ 
(-1)^k \ 
\frac{1}{(2k + 1)^2} 
\ 
f \Big(\frac{2k + 1}{2} \hsx \pi\Big).
\]
Working now with $f (z + x_0)$ at a fixed $x_0 \in \R$ 
(the exponential type of this function is still 1 (f. 17.16)), we conclude that 
\[
f^\prime (x_0) 
\ = \
\frac{4}{\pi^2} \ 
\sum\limits_{k \hsy = \hsy -\infty}^\infty \ 
(-1)^k \ 
\frac{1}{(2k + 1)^2} 
\ 
f \Big(x_0 + \frac{2k + 1}{2} \hsx \pi\Big).
\]
Finally, to eliminate the restriction that $T (f) = 1$, consider the function 
$
\ds
f\Big(\frac{z}{T (f)}\Big)$
of exponential type 1 (cf. 17.16) $-$then
\[
f^\prime\Big(\frac{x}{T (f)}\Big) \ 
\frac{1}{T (f)} \ 
\ = \ 
\frac{4}{\pi^2} \ 
\sum\limits_{k \hsy = \hsy -\infty}^\infty \ 
(-1)^k 
\ 
\frac{1}{(2 \hsy k + 1)^2}
\
f \Big(\frac{x}{T (f)} + \frac{2k + 1}{2\hsy T(f)} \hsx \pi\Big),
\]
i.e., $\forall$ real $x$, 
\[
f^\prime(x)
\ = \ 
\frac{4 \hsy T (f)}{\pi^2} \ 
\sum\limits_{k \hsy = \hsy -\infty}^\infty \ 
(-1)^k
\frac{1}{(2k + 1)^2}
\ 
f \Big(x + \frac{2k + 1}{2\hsy T(f)} \hsx \pi\Big).
\]
\\[-1,25cm]
\end{x}

\begin{x}{\small\bf APPLICATION} \ 
Take $f (z) = \sin z$ and evaluate at $x = 0 \hsy$: 
\\[-.25cm]

\qquad
$\implies$
\[
1
\ = \ 
\frac{4}{\pi^2} \ 
\sum\limits_{k \hsy = \hsy -\infty}^\infty \ 
\frac{1}{(2k + 1)^2} .
\]
\\[-1.cm]
\end{x}

\begin{x}{\small\bf THEOREM} \ 
Let $f \in \sE_0$ with $T (f) > 0$.  
Assume: 
$\forall$ real $x$, 
\[
\abs{ f(x)}
\ \leq \ 
M.
\]
Then
\[
\abs{f^\prime(x)}
\ \leq \ 
M \hsy T (f).
\]

PROOF \ 
In fact, 
\allowdisplaybreaks\begin{align*}
\abs{f^\prime (x) } \ 
&\leq \
T (f) 
\ 
\frac{4}{\pi^2} \ 
\sum\limits_{k \hsy = \hsy -\infty}^\infty \ 
\frac{1}{(2k + 1)^2} 
\hsx 
\Big|
\hsx
f \Big(x + \frac{2k + 1}{2 \hsy T (f)} \hsx \pi\Big)
\hsx
\Big|
\\[15pt]
&\leq \
M \hsy T (f) \ 
\frac{4}{\pi^2} \ 
\sum\limits_{k \hsy = \hsy -\infty}^\infty \ 
\frac{1}{(2k + 1)^2}
\\[15pt]
&= \
M \hsy T (f).
\end{align*}
\\[-1.25cm]
\end{x}


\begin{x}{\small\bf COROLLARY} \ 
Let $f \in \sE_0$ with $T (f) > 0$.  
Assume: 
$\forall$ real $x$, 
\[
\abs{f (x)} 
\ \leq \ 
M.
\]
Then (cf. 17.8)
\[
\big|
f^{(n)} (x)
\big|
\ \leq \ 
M \hsy T (f)^n
\qquad (n = 1, 2, \ldots).
\]
\\[-1.25cm]
\end{x}

\begin{x}{\small\bf EXAMPLE} \ 
Take
\[
f (z) 
\ = \ 
\sum\limits_{k \hsy = \hsy -n}^n \
c_k 
\hsy
e^{\sqrt{-1} \hsx k \hsy z} 
\qquad (\tcf. \ 17.14)
\]
and let $M$ be the maximum of $\abs{f (x)}$ $-$then
\[
\abs{f^\prime (x)} 
\ \leq \ 
M  n.
\]
\\[-1.25cm]
\end{x}

\begin{x}{\small\bf REMARK} \ 
Here is a suggestive way to write the assumption and the conclusion of 17.24: 
\[
\abs{f (x)} 
\ \leq \ 
\Big|
M \hsy e^{\sqrt{-1} \hsx T (f) \hsy x}
\Big|
\ \implies \ 
\abs{f^\prime (x)} 
\ \leq \ 
\Big|\big(M \hsy e^{\sqrt{-1} \hsx T (f) \hsy x}\big)^\prime\Big|.
\]
\\[-1.25cm]
\end{x}

Working on the real axis, let $\norm{\hsx \cdot \hsx}_p$ be the $\Lp^p$-norm: 
\[
\norm{f}_p 
\ = \ 
\bigg[
\int\limits_{-\infty}^\infty \ 
\abs{f (x)}^p 
\td x
\bigg]^{1/p}
\qquad (p \geq 1).
\]

[Note: \ 
$\norm{\hsx \cdot \hsx}_p$ is translation invariant: \ 
$\forall \ f$, $\forall  \ t$, $\norm{f_t}_p = \norm{f}_p$, 
where
$f_t (x) = f (x + t)$.]
\\

\begin{x}{\small\bf THEOREM} \ 
Let $f \in \sE_0$. 
Assume: 
\[
\norm{f}_p 
\ < \ 
\infty.
\]
Then $\forall$ real $y$, 
\[
\int\limits_{-\infty}^\infty \ 
\abs{f (x + \sqrt{-1} \hsx y}^p 
\td x
\ \leq \ 
\norm{f}_p^p 
\hsx 
e^{p \hsy T(f) \hsy \abs{y}}.
\]

PROOF \ 
It suffices to consider the case when $y > 0$.  
To this end, let
\[
F_A (z) 
\ = \ 
\int\limits_{-A}^A \ 
\abs{f (z + t)}^p 
\ \td t.
\]
Then
\allowdisplaybreaks\begin{align*}
\abs{F_A (x) } \ 
&\leq \ 
\int\limits_{-\infty}^\infty \ 
\abs{f (x + t)}^p 
\ \td t
\\[15pt]
&= \ 
\norm{f}_p^p 
\\[15pt]
&< \ 
\infty.
\end{align*}
In addition, $\abs{f (z)}^p$ is subharmonic, thus $\abs{F_A (z)}$ is subharmonic.  
Using Phragm\'en-Lindel\"of in its subharmonic formulation, it follows that
\[
\abs{F_A ( x + \sqrt{-1} \hsx y} 
\ \leq \ 
\norm{f}_p^p 
\ 
e^{p \hsy T(f) \hsy \abs{y}}.
\]
Finish by sending $A$ to infinity.
\\[-.25cm]
\end{x}

\begin{x}{\small\bf LEMMA} \ 
Let $f \in \sE_0$.  
Assume: 
\[
\norm{f}_p
\ < \ 
\infty.
\]
Then $f$ is bounded on the real axis: 
$\forall$ real $x$, 
\[
\abs{f (x)} 
\ \leq \ 
M.
\]

PROOF \ 
Because $\abs{f (z)}^p$ is subharmonic, we have 
\[
\abs{f (x)}^p
\ \leq \ 
\frac{1}{2 \hsy \pi} \ 
\int\limits_0^{2 \hsy \pi} \ 
\big| f(x + r \hsy e^{\sqrt{-1} \hsx \theta}\hsx )\big| ^p 
\ \td \theta
\]
$\implies$
\allowdisplaybreaks
\allowdisplaybreaks\begin{align*}
\abs{f (x)}^p \ 
\int\limits_0^1 \ 
r
\ \td r
&\leq \ 
\frac{1}{2 \hsy \pi} \ 
\int\limits_0^{2 \hsy \pi} \ 
\int\limits_0^1 \
\big| f(x + r \hsy e^{\sqrt{-1} \hsx \theta}\hsx )\big| ^p 
\hsx
r 
\hsx \td r \hsx \td \theta
\\[15pt]
&\leq \ 
\frac{1}{2 \hsy \pi} \ 
\iint\limits_{s^2 + t^2 \leq 1} \ 
\big| f (x + s + \sqrt{-1} \hsx t)\big| ^p 
\ \td s \hsx \td t
\\[15pt]
&\leq \ 
\frac{1}{2 \hsy \pi} \ 
\int\limits_{-1}^1
\ \td t
\ 
\int\limits_{-1}^1 
\abs{f (x + s + \sqrt{-1} \hsx t)}^p 
\ \td s
\end{align*}
$\implies$
\allowdisplaybreaks\begin{align*}
\abs{f (x)}^p \ 
&\leq \ 
\frac{1}{\pi} \ 
\int\limits_{-1}^1
\ \td t
\ 
\int\limits_{-\infty}^\infty \ 
\abs{f (x + s + \sqrt{-1} \hsx t)}^p 
\ \td s
\\[15pt]
&= \ 
\frac{1}{\pi} \ 
\int\limits_{-1}^1
\ \td t
\ 
\int\limits_{-\infty}^\infty \ 
\abs{f (s + \sqrt{-1} \hsx t)}^p 
\ \td s
\\[15pt]
&\leq \ 
\frac{1}{\pi} \ 
\int\limits_{-1}^1
\norm{f}_p^p 
\hsx 
e^{p \hsy T (f) \hsy \abs{t}}
\ \td t
\\[15pt]
&= \ 
\frac{2}{\pi} \ 
\norm{f}_p^p 
\hsx 
\int\limits_0^1 \
e^{p \hsy T (f) \hsy t}
\ \td t
\\[15pt]
&\equiv \ 
M^p.
\end{align*}
\\[-1cm]
\end{x}

\begin{x}{\small\bf REMARK} \ 
If 
$
\norm{f}_p
\ < \ 
\infty
$
and if $T (f) = 0$, then arguing as above, 
\allowdisplaybreaks\begin{align*}
\abs{f (x + \sqrt{-1} \hsx y)}^p
&\leq \ 
\frac{1}{\pi} \ 
\int\limits_{y-1}^{y+1} \ 
\td t
\ 
\int\limits_{-\infty}^\infty \ 
\abs{f (s + \sqrt{-1} \hsx t)}^p
\td s
\\[15pt]
&\leq \ 
\frac{1}{\pi} \ 
\int\limits_{y-1}^{y+1} \ 
\norm{f}_p^p 
\td t
\qquad (\tcf. \ 17.28)
\\[15pt]
&= \ 
\frac{2}{\pi} \ 
\norm{f}_p^p 
\\[15pt]
&< \ 
\infty.
\end{align*}
Therefore $f$ is a constant, hence $f$ is identically zero (cf. 17.34).
\\[-.25cm]
\end{x}

\begin{x}{\small\bf THEOREM} \ 
Let $f \in \sE_0$ with $T (f) > 0$.  
Assume: 
\[
f \in \Lp^p (-\infty, \infty).
\]
Then 
$f^\prime \in \Lp^p (-\infty, \infty)$ 
and 
\[
\norm{f^\prime}_p 
\ \leq \ 
\norm{f}_p 
\hsx
T (f).
\]

PROOF \
Apply 17.22 in the obvious way (legal in view of 17.29).
\\[-.25cm]
\end{x}

\begin{x}{\small\bf SUBLEMMA} \ 
If $f \in \Lp^1(-\infty, \infty)$ and if $f$ is uniformly continuous, then the
limit of $f (x)$ as $x$ approaches plus or minus infinity is zero.  
\\[-.5cm]

PROOF \
Given 
$\varepsilon > 0$, 
choose
$\delta > 0$: 

\[
\abs{x - y} < \delta 
\implies 
\abs{f (x) - f( y)} 
< 
\frac{\varepsilon}{2}.
\]
Choose $R > 0$: 

\[
\int\limits_{R}^\infty \ 
\abs{f} 
\hsx + \hsx 
\int\limits_{-\infty}^{-R} \ 
\abs{f} 
\ < \ 
\varepsilon \hsy \delta.
\]
Claim: 
\[
\begin{cases}
\ 
x \hsx > \hsx R + \delta
\implies
\abs{f (x)} < \varepsilon
\\[4pt]
\ 
x \hsx < \hsx -R - \delta
\implies
\abs{f (x)} < \varepsilon
\end{cases}
.
\]
Consider the first of these assertions and to get a contradiction, assume instead that 
$\abs{f (x)} \geq \varepsilon$ $-$then
\[
x - \delta 
\ < \  
y
\ < \  
x + \delta 
\]
\qquad \qquad
$\implies$
\allowdisplaybreaks
\allowdisplaybreaks\begin{align*}
\abs{f (y)} \ 
&=\ 
\abs{f (x) + f (y) - f (x)}
\\[11pt]
&\geq\ 
\abs{f (x)} - \abs{f (y) - f (x)}
\\[11pt]
&=\
\abs{f (x)} - \abs{f (x) - f (y)}
\\[11pt]
&>\ 
\varepsilon - \frac{\varepsilon}{2} 
\\[11pt]
&=\
\frac{\varepsilon}{2} 
\end{align*}
\qquad \qquad
$\implies$

\[
\int\limits_{x - \delta }^{x + \delta } \ 
\abs{f} 
\ > \ 
\frac{\varepsilon}{2} \hsy (2 \hsy \delta) 
\ = \ 
\varepsilon \hsy \delta.
\]
But

\[
\int\limits_{x - \delta }^{x + \delta } \ 
\abs{f} 
\ < \ 
\int\limits_R^\infty \ 
\abs{f} 
\ < \ 
\varepsilon \hsy \delta.
\]
\\[-.75cm]
\end{x}

\begin{x}{\small\bf LEMMA} \ 
Let 
\[
\Phi 
\ = \ 
\phi * \chisubminusOneOne, 
\]
where 
$\phi \in \Lp^1 (-\infty, \infty)$
and 
$\chisubminusOneOne$ 
is the characteristic function of $[-1,1]$ $-$then 
$\phi \in \Lp^1 (-\infty, \infty)$
is uniformly continuous and 
\[
\begin{cases}
\ 
\lim\limits_{x \ra +\infty} \ 
\Phi (x) \ = \ 0
\\[15pt]
\ 
\lim\limits_{x \ra -\infty} \ 
\Phi (x) \ = \ 0
\end{cases}
.
\]

[Note: \ 
The $*$ stands, of course, for convolution.]
\\[-.25cm]
\end{x}

\begin{x}{\small\bf THEOREM} \ 
Let $f \in \sE_0$.  
Assume: 
\[
\norm{f}_p
\ < \ 
\infty.
\]
Then
\[
f (x) \ra 0 
\quad \text{as} \quad 
\abs{x} \ra \infty.
\]

PROOF \
Proceeding as in 17.29, 
\[
\pi \abs{f (x)}^p 
\ \leq \ 
\int\limits_{-1}^1 \ 
\ \td t
\ 
\int\limits_{-1}^1 \ 
\abs{f (x + s + \sqrt{-1} \hsx t)}^p 
\ \td s.
\]
Let
\[
\phi (s) 
\ = \ 
\int\limits_{-1}^1 \ 
\abs{f (s + \sqrt{-1} \hsx t)}^p 
\ \td t.
\]
Then
\\[-.75cm]
\allowdisplaybreaks
\allowdisplaybreaks\begin{align*}
\int\limits_{-\infty}^\infty \  
\abs{\phi (s) } 
\ \td s \ 
&=\ 
\int\limits_{-\infty}^\infty \  
\bigg(
\int\limits_{-1}^1 \ 
\abs{f (s + \sqrt{-1} \hsx t)}^p 
\ \td t
\bigg)
\ \td s 
\\[15pt]
&=\ 
\int\limits_{-1}^1 \ 
\ \td t \
\int\limits_{-\infty}^\infty \  
\abs{f (s + \sqrt{-1} \hsx t)}^p 
\ \td s
\\[15pt]
&<\ 
\infty.
\end{align*}
I.e.: 
$\phi \in \Lp^1 (-\infty, \infty)$.  
And
\allowdisplaybreaks
\allowdisplaybreaks\begin{align*}
\phi * \chisubminusOneOne (x) \ 
&=\ 
\int\limits_{-\infty}^\infty \  
\phi (x - s) 
\hsy 
\chisubminusOneOne (s) 
\ \td s 
\\[18pt]
&=\ 
\int\limits_{-1}^1 \ 
\phi (x - s)
\ \td s 
\\[18pt]
&=\ 
\int\limits_{-1}^1 \ 
\phi (x + s)
\ \td s 
\\[18pt]
&=\ 
\int\limits_{-1}^1 \ 
\bigg(
\int\limits_{-1}^1 \ 
\abs{f (x + s + \sqrt{-1} \hsx t) }^p 
\ \td t
\bigg)
\hsx \td s 
\\[18pt]
&=\ 
\int\limits_{-1}^1 \ 
\ \td t
\int\limits_{-1}^1 \
\abs{f (x + s + \sqrt{-1} \hsx t) }^p 
\ \td s .
\end{align*}
\\[-.75cm]

\noindent
Now quote 17.33.
\\[-.25cm]
\end{x}

Let $\{\lambda_n \}$ be a real increasing sequence such that 
$\lambda_{n+1} - \lambda_n \hsx \geq \hsx 2 \hsy \delta \hsx > \hsx 0$.
\\[-.25cm]

[Note: \ 
The intervals
$]\lambda_n - \delta, \lambda_n+ \delta[$
are then pairwise disjoint:
\\[-.5cm]

\[
\begin{cases}
\ 
x < \lambda_n  + \delta
\\[4pt]
\ 
x > \lambda_{n+1} - \delta
\end{cases}
\implies \ 
\lambda_n + \delta 
> 
\lambda_{n+1} - \delta
\ \implies \ 
2 \hsy \delta 
>
\lambda_{n+1}  - \lambda_n .]
\]
\\[-1cm]

\begin{x}{\small\bf THEOREM} \ 
Let $f \in \sE_0$.  
Assume: 
\[
\norm{f}_p 
\ < \ 
\infty.
\]
Then
\[
\sum\limits_n \ 
\abs {f (\lambda_n)}^p 
\ \leq \ 
2 \ 
\frac{e^{\delta \hsy p \hsy T (f)}}{\delta \hsy \pi}
\hsx 
\norm{f}_p ^p.
\]

PROOF \ 
We have
\allowdisplaybreaks
\allowdisplaybreaks\begin{align*}
\sum\limits_n \ 
\abs {f (\lambda_n)}^p \ 
&\leq \ 
\frac{1}{\delta^2 \hsy \pi} \ 
\sum\limits_n \ 
\iint\limits_{\abs{z} \leq \delta} \ 
\abs{f (\lambda_n + z)}^p
\ \td x \hsy \td y
\\[15pt]
&\leq \ 
\frac{1}{\delta^2 \hsy \pi} \ 
\sum\limits_n \ 
\int\limits_{-\delta}^\delta \
\int\limits_{-\delta}^\delta \
\abs{f (\lambda_n + x + \sqrt{-1} \hsx y)}^p
\ \td x \hsy \td y
\\[15pt]
&= \
\frac{1}{\delta^2 \hsy \pi} \ 
\sum\limits_n \ 
\int\limits_{-\delta}^\delta \
\int\limits_{\lambda_n-\delta}^{\lambda_n+\delta} \
\abs{f (x + \sqrt{-1} \hsx y)}^p
\ \td x  \hsy \td y
\\[15pt]
&\leq \
\frac{1}{\delta^2 \hsy \pi} \ 
\int\limits_{-\delta}^\delta \
\int\limits_{-\infty}^\infty \ 
\abs{f (x + \sqrt{-1} \hsx y)}^p
\ \td x  \hsy \td y
\\[15pt]
&\leq \
\frac{1}{\delta^2 \hsy \pi} \ 
\int\limits_{-\delta}^\delta \
\norm{f}_p ^p
\hsx 
e^{p \hsy T (f) \hsy \abs{y}}
\ \td y
\qquad (\tcf. \ 17.28)
\\[15pt]
&\leq \
\frac{2}{\delta^2 \hsy \pi} \ 
\bigg(
\int\limits_0^\delta \ 
e^{p \hsy T (f) \hsy y}
\ \td y
\bigg)
\norm{f}_p ^p
\\[15pt]
&\leq \
2 \ 
\frac{e^{\delta \hsy p \hsy T (f)}}{\delta \hsy \pi}
\hsx 
\norm{f}_p ^p.
\end{align*}
\\[-.25cm]
\end{x}

\chapter{
$\boldsymbol{\S}$\textbf{18}.\quad  THE BOREL TRANSFORM}
\setlength\parindent{2em}
\setcounter{theoremn}{0}
\renewcommand{\thepage}{\S18-\arabic{page}}

\qquad 
Let $K$ be a nonempty convex compact subset of $\Cx$.
\\[-.25cm]

\begin{x}{\small\bf DEFINITION} \ 
Put
\[
H_K (z) 
\ = \ 
\sup\limits_{w \in K} \ 
\Reg (w \hsy z).
\]
Then
\[
H_K : \Cx \ra \Cx
\]
is called the \un{support function} of $K$.
\\[-.25cm]
\end{x}

\qquad
{\small\bf \un{N.B.}} \ 
$H_K$ is homogeneous of degree 1: 
\[
H_K (t \hsy z) 
\ = \ 
t \hsy H_K (z) 
\qquad (t > 0).
\]
Therefore 
\[
H_K (z) 
\ = \ 
H_K \big(\abs{z} \hsy e^{\sqrt{-1} \hsx \theta}\big)
\ = \ 
\abs{z} \hsx 
H_K \big(e^{\sqrt{-1} \hsx \theta}\big).
\]

[Note: \ 
Of course, $H_K(0) = 0$.]
\\[-.25cm]

\qquad
{\small\bf \un{N.B.}} \ 
$H_K$ is convex: 
\[
H_K  (\lambda \hsy z_1 + (1 - \lambda) z_2) 
\ \leq \ 
\lambda \hsy H_K (z_1) + (1 - \lambda) \hsy H_K(z_2)
\qquad 
(0 < \lambda < 1).
\]

[Note: \ 
It thus follows that $H_K$ is continuous.]
\\[-.25cm]

\begin{x}{\small\bf EXAMPLE} \ 
Take $K = \{x_0 + \sqrt{-1} \hsx y_0\}$ (a singleton) $-$then 
\[
H_K (z) 
\ = \ 
\abs{z} \hsx 
(x_0 \hsy \cos \theta - y_0 \hsy \sin \theta).
\]
\\[-1.25cm]
\end{x}

\begin{x}{\small\bf EXAMPLE} \ 
Take $K = \{z : \abs{z} \leq R\}$ $-$then 
\[
H_K (z) 
\ = \ 
R \hsy\abs{z}.
\]
\\[-1.25cm]
\end{x}


\begin{x}{\small\bf EXAMPLE} \ 
Take $K = [-a, a]$ $(a > 0)$ $-$then 
\[
H_K (z) 
\ = \ 
a \hsy \abs{z} \hsy \abs{\cos \theta}.
\]
\\[-1.5cm]
\end{x}

\begin{x}{\small\bf EXAMPLE} \ 
Take $K = [- \sqrt{-1} \hsx a, \sqrt{-1} \hsx a]$ $(a > 0)$ $-$then 
\[
H_K (z) 
\ = \ 
a \hsy \abs{z} \hsy \abs{\sin \theta}.
\]
\\[-1.5cm]
\end{x}

\begin{x}{\small\bf LEMMA} \ 
$\forall \ w \in K$, 
\[
(\Reg w) \cos \theta - (\Img w) \sin \theta 
\ = \ 
\Reg \big(w \hsy e^{\sqrt{-1} \hsx \theta}\big)
\ \leq \ 
H_K \big(e^{\sqrt{-1} \hsx \theta}\big).
\]
\\[-1.25cm]
\end{x}

\begin{x}{\small\bf APPLICATION} \ 
\\[-.25cm]

\qquad \textbullet \quad
Take $\theta = 0$ to get 
\[
\Reg w
\ \leq \ 
H_K (1).
\]

\qquad \textbullet \quad
Take $\theta = \pi$ to get 
\[
-\Reg w
\ \leq \ 
H_K (-1).
\]
Therefore 
\[
-H_K (-1) 
\ \leq \ 
\Reg w
\ \leq \ 
H_K (1) .
\]
\\[-1.25cm]
\end{x}

\begin{x}{\small\bf APPLICATION} \ 
\\[-.25cm]

\qquad \textbullet \quad
Take $\ds\theta = \frac{\pi}{2}$ to get 
\[
-\Img w
\ \leq \ 
H_K (\sqrt{-1}) .
\]

\qquad \textbullet \quad
Take $\ds\theta = \frac{3 \hsy \pi}{2}$ to get 
\[
-\Img w \hsy (-1)
\ \leq \ 
H_K (-\sqrt{-1}) .
\]
Therefore 
\[
- H_K (\sqrt{-1}) 
\ \leq \ 
\Img w
\ \leq \ 
H_K (-\sqrt{-1}).
\]
\\[-1.25cm]
\end{x}


\begin{x}{\small\bf EXAMPLE} \ 
Suppose that
\[
\begin{cases}
\ 
H_K (1) \leq 0
\\[4pt]
\ 
H_K (-1) \leq 0
\end{cases}
.
\]

Then
\[
0
\ \leq \ 
-H_K (-1)
\leq
\Reg w
\ \leq \ 
H_K (1)
\ = \ 
0
\]
\qquad \qquad
$\implies$
\[
\Reg w 
 \ = 0.
\]
Therefore $K$ is contained in the imaginary axis.
\\[-.25cm]
\end{x}

\begin{x}{\small\bf DEFINITION} \ 
Suppose that

\[
f (z) 
\ = \ 
\sum\limits_{n = 0}^\infty \ 
\frac{\gamma_n}{n !} \hsx z^n
\]
is of exponential type $-$then its \un{Borel transform} $\sB_f$ is defined by the prescription 
\[
\sB_f (w)
\ = \ 
\sum\limits_{n = 0}^\infty \ 
\frac{\gamma_n}{w^{n + 1}}.
\]

[Note: \ 
The series converges if $\abs{w} > T (f)$ and diverges if $\abs{w} < T(f)$.]
\\[-.25cm]
\end{x}

\begin{x}{\small\bf EXAMPLE} \ 
Take $f (z) = e^z$ $-$then
\[
\sB_f (w)
\ = \ 
\frac{1}{w - 1}.
\]
\\[-1.25cm]
\end{x}

\begin{x}{\small\bf EXAMPLE} \ 
Take $f (z) = e^{\sqrt{-1} \hsx z}$ $-$then 
\[
\sB_f (w)
\ = \ 
\frac{1}{w - \sqrt{-1}}.
\]
\\[-1.25cm]
\end{x}

\begin{x}{\small\bf LEMMA} \ 
Fix $T^\prime > T (f)$ and suppose that $\Reg w > 2 \hsy T^\prime$ $-$then
\[
\sB_f (w)
\ = \ 
\int\limits_0^\infty \ 
f (t) \hsy e^{-w t}
\ \td t.
\]
PROOF \ 
First of all, 
\allowdisplaybreaks
\allowdisplaybreaks\begin{align*}
\Big|
\hsx
f (z) - \sum\limits_{k = 0}^n \ c_k \hsy z^k
\hsx
\Big|
\ 
&\leq\ 
\sum\limits_{k = n+1}^\infty \ \abs{c_k }\hsy \abs{r}^k\ 
\\[15pt]
&=\ 
\sum\limits_{k = n + 1}^\infty \ 
\abs{c_k }\hsx 
R^k \hsx 
\Big(\frac{r}{R}\Big)^k
\qquad (R > r)
\\[15pt]
&\leq\ 
M (R; f) \ 
\sum\limits_{k = n + 1}^\infty \ 
\Big(\frac{r}{R}\Big)^k
\\[15pt]
&=\ 
\Big(\frac{r}{R}\Big)^{n + 1}
M (R; f) \hsx
\frac{1}
{
\raisebox{-.2cm}
{
$
\ds
1 - \frac{r}{R}
$
}
}
\\[15pt]
&\leq\ 
\Big(\frac{r}{R}\Big)^{n + 1} \hsx
e^{R \hsy T^\prime} \hsx
\frac{R}{R - r}.
\end{align*}
Now take $R = 2 \hsy r$ to get
\[
\Big|
\hsx
f (z) - \sum\limits_{k = 0}^n \ c_k \hsy z^k
\hsx
\Big|
\ 
\ \leq \ 
\Big(\frac{1}{2}\Big)^n \hsy e^{2 \hsy r \hsy T^\prime}.
\]
Since 
\[
\abs{e^{- w \hsy t}} 
\ = \ 
\exp (- (\Reg w) \hsy t),
\]
it then follows that
\allowdisplaybreaks
\allowdisplaybreaks\begin{align*}
\bigg|
\int\limits_0^\infty \ 
f (t) \hsy e^{-w \hsy t} 
\ \td t
&\hsx - \hsx
\int\limits_0^\infty \ 
\bigg(
\sum\limits_{k = 0}^n \
c_k \hsy t^k
\bigg)
\hsy e^{-w \hsy t} 
\ \td t
\bigg| \ 
\\[15pt]
&\leq
\int\limits_0^\infty \ 
\Big|
\hsx
f (t) 
\
-
\
\sum\limits_{k = 0}^n \
c_k \hsy t^k
\hsx
\Big|
\ 
\exp (- (\Reg w) \hsy t)
\ \td t
\\[15pt]
&\leq
\Big(\frac{1}{2}\Big)^n
\int\limits_0^\infty \ 
\exp((2 \hsy T^\prime  - \Reg w) \hsy t)
\ \td t.
\end{align*}
But
\[
\Reg w 
\ > \ 
2 \hsy T^\prime 
\implies 
(2 \hsy T^\prime  - \Reg w) 
\ < \ 
0
\]
\qquad
$\implies$
\[
\int\limits_0^\infty \ 
\exp((2 \hsy T^\prime  - \Reg w) \hsy t)
\ \td t
\ < \ 
\infty.
\]
Therefore the infinite series

\[
\sum\limits_{n = 0}^\infty \
c_n \ 
\int\limits_0^\infty \ 
\hsy t^n \hsy e^{-w \hsy t}
\ \td t
\]
is convergent and has sum
$
\ds
\int\limits_0^\infty \ 
f (t) \hsy e^{-w \hsy t}
\ \td t.  
$  
And finally
\allowdisplaybreaks\begin{align*}
\sum\limits_{n = 0}^\infty \ 
c_n \ 
\int\limits_0^\infty \ 
t^n \hsy e^{-w \hsy t}
\ \td t \ 
&=\ 
\sum\limits_{n = 0}^\infty \ 
\gamma_n \ 
\int\limits_0^\infty \ 
\frac{t^n}{n !} \ e^{-w \hsy t}
\ \td t \ 
\\[15pt]
&=\ 
\sum\limits_{n = 0}^\infty \ 
\frac{\gamma_n}{w^{n + 1}} \ 
\\[15pt]
&=\ 
\sB_f (w).
\end{align*}

[Note: \ 
The constant implicit in the asymptotics has been set equal to 1.  
To proceed in general, break 
$
\ 
\ds
\int\limits_0^\infty \ 
\ldots
\ \td t
\ 
$
into 
$
\ 
\ds
\int\limits_0^{t_0} \ 
\ldots
\ \td t
\hsx + \hsx
\int\limits_{t_0}^\infty \ 
\ldots
\ \td t
$\hsx .]
\\[-.25cm]
\end{x}

Keeping still to the assumption that $f$ is of exponential type, 
let $K_f$ denote the intersection of all the convex compact subsets of $\Cx$ outside of which $\sB_f$ is holomorphic. 
\\

\qquad
{\small\bf \un{N.B.}} \ 
Therefore $K_f$ is the smallest convex compact subset of $\Cx$ outside of which $\sB_f$ is holomorphic.
\\[-.25cm]

\begin{x}{\small\bf DEFINITION} \ 
$K_f$ is the \un{indicator diagram} of $f$.
\\[-.25cm]
\end{x}

\begin{x}{\small\bf LEMMA} \ 
The extreme points of $K_f$ are singular points of $\sB_f$.
\\[-.5cm]

PROOF \ 
If $p \in K_f$ were an extreme point of $K_f$ which was not a singular point of $\sB_f$, 
then upon removing a certain neighborhood of $p$ from $K_f$ one would be led to a smaller 
convex compact subset of $\Cx$ outside of which $\sB_f$ is holomorphic.
\\[-.25cm]
\end{x}

\begin{x}{\small\bf EXAMPLE} \ 
Let 
\[
f (z) 
\ = \ 
\sum\limits_{k = 1}^n \ 
P_k (z) \hsy e^{c_k \hsy z}
\]
be an 
\un{exponential polynomial} 
(meaning that the $P_k$ are polynomials and the $c_k$ are complex numbers).  
Since the Borel transform of a monomial 
$z^p \hsy e^{c_k \hsy z}$ 
equals 
$p ! \hsx (w - c_k)^{-p-1}$, 
the poles at the $c_k$ are the only singularities of the Borel transform of $f$, 
so the indicator diagram of $f$ is the convex hull of the set 
$\{c_1, \ldots, c_n\}$.
\\[-.25cm]
\end{x}

\begin{x}{\small\bf NOTATION} \ 
Write $H_f$ in place of $H_{K_f}$.
\\[-.25cm]
\end{x}

\begin{x}{\small\bf EXAMPLE} \ 
Take $f (z) = \sin \pi \hsy z$ $-$then 
\[
\sB_f (w) 
\ = \ 
\frac{1}{2 \hsy \sqrt{-1}} \ 
\bigg[
\frac{1}{w -\sqrt{-1} \hsx \pi }
\hsx - \hsx 
\frac{1}{w +\sqrt{-1} \hsx \pi }
\bigg]
\]
and 
\[
K_f 
\ = \ 
[-\sqrt{-1} \hsx \pi, \sqrt{-1} \hsx \pi].
\]
Here
\[
H_f (z) 
\ = \ 
\pi \hsy \abs{z} \hsy \abs{\sin \theta}
\qquad (\tcf. \ 18.5), 
\]
so
\[
H_f (\pm \sqrt{-1}) 
\ = \ 
\pi 
\ = \ 
\tau (f).
\]
\\[-1.25cm]
\end{x}

Let $\Gamma$ be a rectifiable Jordan curve containing $K_f$ in its interior.
\\

\begin{x}{\small\bf THEOREM} \ 
We have

\[
f (z) 
\ = \ 
\frac{1}{2 \hsy \pi \hsy \sqrt{-1}} \
\int\limits_\Gamma \ 
\sB_f (w) \hsy e^{z \hsy w} 
\ \td w.
\]

PROOF
Take for $\Gamma$ the circle $\abs{w} = T (f) + \varepsilon$ $(\varepsilon > 0)$ $-$then 
\allowdisplaybreaks
\allowdisplaybreaks\begin{align*}
\frac{1}{2 \hsy \pi \hsy \sqrt{-1}} \
\int\limits_\Gamma \ 
\sB_f (w) \hsy e^{z \hsy w} 
\ \td w
&=\ 
\frac{1}{2 \hsy \pi \hsy \sqrt{-1}} \
\int\limits_\Gamma \ 
\Big(
\sum\limits_{n = 0}^\infty \ 
\frac{n ! \hsy c_n}{w^{n+1}}
\Big)
\hsy e^{z \hsy w} 
\ \td w
\\[15pt]
&=\ 
\sum\limits_{n = 0}^\infty \ 
n! \hsy c_n
\ 
\frac{1}{2 \hsy \pi \hsy \sqrt{-1}} \
\int\limits_\Gamma \ 
\frac{e^{z \hsy w}}{w^{n+1}}
\ \td w
\\[15pt]
&=\ 
\sum\limits_{n = 0}^\infty \ 
c_n \hsy z^n 
\\[15pt]
&=\ 
f (z).
\end{align*}
\\[-1.25cm]
\end{x}

\begin{x}{\small\bf LEMMA} \ 
$K_f = \emptyset$ iff $f \equiv 0$.
\\[-.5cm]

PROOF 
If $K_f = \emptyset$, then $\sB_f$ is everywhere holomorphic (including $\infty$), 
thus $\sB_f$ is a constant.  
But $\sB_f (\infty) = 0$, so $\sB_f \equiv 0$ $\implies$ $f \equiv 0$ (cf. 18.19).  
Conversely, if $f \equiv 0$, then $\forall \ n$, $\gamma_n = 0$, hence $\sB_f \equiv 0$.
\\[-.25cm]
\end{x}

\begin{x}{\small\bf EXAMPLE} \ 
Suppose that 
\[
\begin{cases}
\ 
H_f (\sqrt{-1}) \ < \  0
\\[4pt]
\ 
H_f (-\sqrt{-1}) \ < \  0
\end{cases}
.
\]
Then $K_f = \emptyset$, implying that $f \equiv 0$.
\\[-.5cm]

[From 18.8,

\[
\begin{cases}
\ 
-H_f (\sqrt{-1}) > 0
\implies \Img w > 0
\\[4pt]
\ 
H_f (-\sqrt{-1}) < 0
\implies \Img w < 0
\end{cases}
.]
\]
\\[-1.cm]
\end{x}

\begin{x}{\small\bf NOTATION} \ 
$\sH_0 (\infty)$ is the set of functions that are holomorphic near $\infty$ and vanish at $\infty$. 
\\[-.5cm]

[Note: \ 
If $\Phi \in \sH_0 (\infty)$, then there is an expansion 
\[
\Phi (z) 
\ = \ 
\sum\limits_{n = 0}^\infty \ 
\frac{A_n}{z^{n + 1}}, 
\]
where
\[
A_n 
\ = \ 
\frac{1}{2 \hsy \pi \hsy \sqrt{-1}} \ 
\int\limits_\Gamma \ 
\Phi (w) \hsy w^n 
\ \td w
\qquad (n = 0, 1, \ldots),
\]
$\Gamma$ a suitable contour.]
\\[-.75cm]
\end{x}

E.g.: 
\[
f \in \sE_0 \implies \sB_f \in \sH_0 (\infty).
\]
\\[-1cm]

\begin{x}{\small\bf LEMMA} \ 
The arrow
\[
\sB : \sE_0 \ra \sH_0 (\infty)
\]
that sends $f$ to $\sB_f$ is a linear injection.
\\[-.5cm]

PROOF \ 
Using the inversion formula for the Laplace transform, if 
$\sB_f = \sB_g$, then for $u = \Reg w \gg 0$ (cf. 18.13), 
\allowdisplaybreaks\begin{align*}
f (t) \ 
&= \
\frac{1}{2 \hsy \pi \hsy \sqrt{-1}} \ 
\int\limits_{u - \sqrt{-1} \hsx \infty}^{u + \sqrt{-1} \hsx \infty} \ 
e^{t \hsy w} \hsy \sB_f (w)
\ \td w
\\[15pt]
&=\ 
\frac{1}{2 \hsy \pi \hsy \sqrt{-1}} \ 
\int\limits_{u - \sqrt{-1} \hsx \infty}^{u + \sqrt{-1} \hsx \infty} \ 
e^{t \hsy w} \hsy \sB_g (w)
\ \td w
\\[15pt]
&=\ 
g (t).
\end{align*}
\\[-1.25cm]
\end{x}

\qquad
{\small\bf \un{N.B.}} \ 
The inverse

\[
\sB^{-1} : \sB \sE_0 \ra \sE_0
\]
is constructed via 18.19: 

\[
\sB^{-1} (\sB_f) (z) 
\ = \ 
\frac{1}{2 \hsy \pi \hsy \sqrt{-1}} \ 
\int\limits_\Gamma \ 
\sB_f (w) \hsy e^{z \hsy w}
\ \td w.
\]
\\[-1.25cm]

\begin{x}{\small\bf LEMMA} \ 
The arrow

\[
\sB : \sE_0 \ra \sH_0 (\infty)
\]
that sends $f$ to $\sB_f$ is a linear surjection.
\\[-.5cm]

PROOF \ 
Fix $\Phi \in \sH_0 (\infty)$ and let $S (\Phi)$ be the smallest convex compact subset of $\Cx$ 
in whose complement $\Phi$ is holomorphic.  
Put
\[
N (S (\Phi), r) 
\ = \ 
\{w \in \Cx : d (w, S(\Phi)) < r\}
\]
and let $\Gamma$ be a rectifiable Jordan curve containing $S (\Phi)$ in its interior: 
\[
S (\Phi) \subsetx \itr \Gamma \subsetx N (S (\Phi), r).
\]
Consider now the holomorphic function 
\[
f (z) 
\ = \ 
\frac{1}{2 \hsy \pi \hsy \sqrt{-1}} \ 
\int\limits_\Gamma \ 
\Phi (w) \hsy e^{z \hsy w}
\ \td w.
\]
Then
\allowdisplaybreaks\begin{align*}
\sup\limits_{w \in \Gamma} \ 
\Reg (z \hsy w) \ 
&\leq \ 
\sup\limits_{w \in S(\Phi)} \ 
(\Reg (z \hsy w) + r \hsy \abs{z})
\\[11pt]
&=\ 
H_{S (\Phi)} (z) + r \hsy \abs{z}
\end{align*} 
\qquad 
$\implies$
\[
\abs{f (z)} 
\ \leq \ 
C \hsy \exp(H_{S (\Phi)} (z) + r \hsy \abs{z}), 
\]
where 
\[
C 
\ = \ 
\frac{\len \Gamma}{2 \hsy \pi} 
\hsx
\sup\limits_{w \in \Gamma} \ 
\abs{\Phi (w)}.
\]
Choose $R \gg 0$: 
\[
S (\Phi) \subsetx \{z : \abs{z} \leq R\}
\]
\qquad 
$\implies$
\[
\abs{f (z)} 
\ \leq \ 
C \hsy \exp(R \abs{z} + r \abs{z})
\qquad (\tcf. \ 18.3).
\]
Therefore $f \in \sE_0$.  
And $\sB_f = \Phi$ (details below).
\\[-.5cm]

[Let $T$ be the analytic functional  defined by the rule
\[
\langle F, T \rangle 
\ = \ 
\frac{1}{2 \hsy \pi \hsy \sqrt{-1}} \ 
\int\limits_\Gamma \ 
\Phi (w) \hsy F (w)
\ \td w.
\]
Then by definition its FL-transform $\widehat{T}$ is the function
\[
\langle e^{z \hsy w}, \widehat{T} \rangle 
\ = \ 
\frac{1}{2 \hsy \pi \hsy \sqrt{-1}} \ 
\int\limits_\Gamma \ 
\Phi (w) \hsy e^{z \hsy w}
\ \td w,
\]
thus here
\[
\langle e^{z \hsy w}, \widehat{T} \rangle 
\ = \ f (z). 
\]
On the other hand, the prescription 

\[
F 
\ra 
\frac{1}{2 \hsy \pi \hsy \sqrt{-1}} \ 
\int\limits_\Gamma \ 
\sB_f (w) \hsy F (w) 
\ \td w
\]
defines an analytic functional $S$ whose FL-transform is also $f(z)$ (cf. 18.19).
But

\[
f (z) \ = \ 
\begin{cases}
\ 
\langle e^{z \hsy w}, \widehat{T} \rangle 
\ = \ 
\ds
\sum\limits_{n = 0}^\infty \ 
\frac{\langle w^n, T \rangle}{n !} \hsy z^n
\\[26pt]
\ 
\langle e^{z \hsy w}, \widehat{S} \rangle 
\ = \ 
\ds
\sum\limits_{n = 0}^\infty \ 
\frac{\langle w^n, S \rangle}{n !} \hsy z^n
\end{cases}
\]
\\[-1cm]

\qquad 
$\implies$
\[
\langle w^n, T \rangle 
\ = \ 
\langle w^n, S \rangle 
\qquad (n = 0, 1, \ldots)
\]

\qquad 
$\implies$
\[
\Phi 
\ = \ 
\sB_f.]
\]

[Note: \ 
See 20.2 for the definition of ``analytic functional''.]
\\[-.25cm]
\end{x}


\chapter{
$\boldsymbol{\S}$\textbf{19}.\quad  THE INDICATOR FUNCTION}
\setlength\parindent{2em}
\setcounter{theoremn}{0}
\renewcommand{\thepage}{\S19-\arabic{page}}

\qquad 
Let $f$ be an entire function of exponential type. 
\\[-.25cm]

\begin{x}{\small\bf DEFINITION} \ 
The \un{indicator function}
\[
h_f : \Cx^\times \ra \Cx
\]
of $f$ is defined by
\[
h_f (z)
\ = \ 
\underset{r \ra \infty}{\limsupx} \ 
\frac{\log \abs{f (r \hsy z)}}{r}.
\]

[Note: \ 
Sometimes
\[
h_f (e^{\sqrt{-1} \hsx \theta})
\ = \ 
\underset{r \ra \infty}{\limsupx} \ 
\frac{\log 
\big|
f (r \hsy e^{\sqrt{-1} \hsx \theta})
\big|
}{r}
\]
is referred to as the exponential type of $f$ in the direction $\theta$.  
Obviously, 

\[
h_f \big(e^{\sqrt{-1} \hsx \theta}\big)
\ \leq \ 
T (f).]
\]
\\[-1.5cm]
\end{x}

\begin{x}{\small\bf EXAMPLE} \ 
Take $f (z) = \exp (a + \sqrt{-1} \hsx b) \hsy z$ 
\ ($a, b \in \R$) 
$-$then
\[
h_f (z) 
\ = \ 
\abs{z} \hsy (a \cos \theta - b \sin \theta)
\qquad (z = \abs{z} e^{\sqrt{-1} \hsx \theta}).
\]
\\[-1.25cm]
\end{x}

\begin{x}{\small\bf LEMMA} \ 
If $f \equiv 0$, then 
$h_f \equiv -\infty$ 
and if 
$h_f \equiv -\infty$, 
then 
$f \equiv 0$.
\\[-.25cm]
\end{x}

\begin{x}{\small\bf LEMMA} \ 
If $f \not\equiv 0$, then 
$h_f \big(e^{\sqrt{-1} \hsx \theta}\big) > -\infty$ everywhere.
\\[-.25cm]
\end{x}

\begin{x}{\small\bf LEMMA} \ 
If $f \not\equiv 0$, then $h_f (z)$ is a continuous function of $z \in \Cx$ if 
$h_f (0)$ is defined to be 0.
\\[-.25cm]
\end{x}

\qquad
{\small\bf \un{N.B.}} \ 
$h_f$ ($f \not\equiv 0$) is homogeneous of degree 1: 
\[
h_f (t \hsy z) 
\ = \ 
t \hsx h_f (z) 
\qquad (t > 0).
\]
Therefore
\[
h_f (z) 
\ = \ 
h_f (\abs{z} \hsy e^{\sqrt{-1} \hsx \theta})
\ = \ 
\abs{z} 
\hsy 
h_f (e^{\sqrt{-1} \hsx \theta}).
\]
\\[-1.25cm]

\begin{x}{\small\bf REMARK} \ 
It can be shown that $h_f$ ($f \not\equiv 0)$ is subharmonic. 
\\[-.25cm]
\end{x}

\begin{x}{\small\bf THEOREM} \ 
If $f \not\equiv 0$, then $H_f = h_f$.
\\[-.5cm]

PROOF \ 
It will be enough to prove that $\forall \ \theta$, 
\[
H_f (e^{\sqrt{-1} \hsx \theta})
\ = \ 
h_f (e^{\sqrt{-1} \hsx \theta}).
\]
To this end, we shall first show that
\[
h_f (e^{\sqrt{-1} \hsx \theta})
\ \leq \ 
H_f (e^{\sqrt{-1} \hsx \theta}).
\]
Thus write
\[
f (z) 
\ = \ 
\frac{1}{2 \hsy \pi \sqrt{-1}} \ 
\int\limits_{\Gamma_\varepsilon} \ 
\sB_\gf (w) 
\hsx
e^{z \hsy w} 
\ \td w
\qquad (\tcf. \ 18.19),
\]
choosing $\Gamma_\varepsilon$ so as to remain within the $\varepsilon$-neighborhood of $K_f$ subject to 
$
K_f \subset \itr \Gamma_\varepsilon
$ 
$-$then
\[
\big|
f \big(r \hsy e^{\sqrt{-1} \hsx \theta} \big)
\big|
\ \leq \ 
\frac{\len \Gamma_\varepsilon}{2 \hsy \pi} 
\ \cdot \ 
\sup\limits_{w \in \Gamma_\varepsilon}  \ 
\abs{\sB_\gf (w)}
\ \cdot \ 
\sup\limits_{w \in \Gamma_\varepsilon}  \ 
\exp \big( r \hsy \Reg \big(w \hsy e^{\sqrt{-1} \hsx \theta}\big)\big)
\]
\\[-1cm]

\qquad 
$\implies$
\allowdisplaybreaks\begin{align*}
h_f (e^{\sqrt{-1} \hsx \theta})\ 
&\leq \ 
\sup\limits_{w \in \Gamma_\varepsilon}  \ 
\Reg \big(w \hsy e^{\sqrt{-1} \hsx \theta}\big)
\\[15pt]
&\leq \ 
H_f \big(e^{\sqrt{-1} \hsx \theta}\big) + \varepsilon
\end{align*}
\qquad
\qquad 
$\implies$
\[
h_f (e^{\sqrt{-1} \hsx \theta})
\ \leq \ 
H_f (e^{\sqrt{-1} \hsx \theta}).
\]
As for the opposite direction, it suffices to work at $\theta = 0$, the claim being that
\[
H_f (1) 
\ \leq \ 
h_f (1).
\]
But $\forall \ \varepsilon > 0$, 
\[
\abs{f (t)} 
\ < \ 
\exp( (h_f (1) + \varepsilon) t) 
\qquad (t \gg 0).
\]
Therefore the integral 
\[
\int\limits_0^\infty \ 
f (t) 
\hsx 
e^{-w \hsy t} 
\ \td t
\]
is a holomorphic function of $w$ in the half-plane $\Reg w > h_f (1)$.  
Since $h_f (1) \leq T (f)$, it follows from 18.13 that $\sB_\gf$ has no singularities to the right 
of the line $x = h_f (1)$, so $H_f (1) \leq h_f (1)$.
\\[-.25cm]
\end{x}

\begin{x}{\small\bf APPLICATION} \ 
\\[-.25cm]

\qquad \textbullet \quad
$H_f$ convex $\implies$ $h_f$ convex.
\\[-.25cm]

\qquad \textbullet \quad
$h_f$ subharmonic $\implies$ $H_f$ subharmonic.
\\[-.25cm]
\end{x}

\begin{x}{\small\bf REMARK} \ 
Any complex valued function with domain $\Cx$ which is subharmonic and homogeneous of degree 1 is necessarily convex. 
\\[-.25cm]
\end{x}

\begin{x}{\small\bf LEMMA} \ 
If $T (f) > 0$, then $T (f) = \tau (f)$ (cf. 17.3) and 
\[
\tau (f) 
\ = \ 
\sup\limits_{0 \leq \theta \leq 2 \hsy \pi} \ 
h_f \big(e^{\sqrt{-1} \hsx \theta} \big).
\]
\\[-1.25cm]
\end{x}

\begin{x}{\small\bf LEMMA} \ 
Assume that $f \not\equiv 0$ $-$then $T (f) = 0$ iff $h_f = 0$. 
\\[-.5cm]

PROOF \ 
If $T (f) = 0$, then $\sB_f$ is holomorphic in the region $\abs{w} > 0$, 
so $K_f = \{0\}$ (cf. 18.20), hence $H_f = 0$, hence $h_f = 0$.  
Conversely, if $h_f = 0$ then $T (f) = 0$
($T (f) > 0$ being ruled out by 19.10).
\\[-.25cm]
\end{x}

\begin{x}{\small\bf LEMMA} \ 
If $f$, $g \in \sE_0$ and if $g$ is an exponential polynomial, then 
\[
h_{f \hsy g} 
\ = \ 
h_f + h_g.
\]

[Note: \ 
Recall that $\sE_0$ is an algebra (cf. 17.15), thus $f  g \in \sE_0$.]
\\[-.25cm]
\end{x}

\begin{x}{\small\bf COROLLARY} \ 
If $f$, $g \in \sE_0$ and if $g$ is an exponential polynomial, and 
\\[.25cm]
if 
$
\ 
\ds
\frac{f}{g}
\ 
$
is entire, then 
$
\ 
\ds
\frac{f}{g}
\
$
is of exponential type (cf. 17.9) and 
\[
h_{\frac{f}{g}}
\ = \ 
h_f - h_g.
\]
\\[-1.5cm]
\end{x}

\begin{x}{\small\bf THEOREM} \ 
Suppose that $f \in \sE_0$ has the property that $h_f (\pm \sqrt{-1}) < \pi$.  
Assume further that $f (n) = 0$ for 
$n = 0, \pm 1, \pm 2, \ldots$ $-$then $f \equiv 0$.
\\[-.5cm]

PROOF \ 
Let 
\[
\phi (z)
\ = \ 
\frac{f (z)}{g (z)},
\]
where $g (z) = \sin \pi z$ $-$then $\phi \in \sE_0$.  
But $g$ is an exponential polynomial, so
\[
h_\phi
\ = \ 
h_f - h_g
\]
\qquad 
$\implies$
\allowdisplaybreaks\begin{align*}
h_\phi (\pm \sqrt{-1}) \
&=\ 
h_f (\pm \sqrt{-1}) - h_g (\pm \sqrt{-1})
\\[11pt]
&=\ 
h_f (\pm \sqrt{-1}) - \pi
\qquad (\tcf. \ 18.5)
\\[11pt]
&<\ 
\pi - \pi
\\[11pt]
&=\ 
0
\\[4pt]
\implies
\hspace{2cm}
&
\\[4pt]
\phi \ 
&\equiv \ 0
\qquad (\tcf. \ 18.21 \quad (h_\phi = H_\phi))
\\[4pt]
\implies
\hspace{2cm}
&
\\[4pt]
f \ 
&\equiv 0.
\end{align*}
\\[-1cm]
\end{x}

\begin{x}{\small\bf REMARK} \ 
One cannot replace 
$h_f (\pm \sqrt{-1}) < \pi$ 
by 
$h_f (\pm \sqrt{-1}) = \pi$ 
(consider $\sin \pi z$).
\\[-.25cm]
\end{x}

\begin{x}{\small\bf LEMMA} \ 
If $f \in \sE_0$, then $\forall$ complex constant $c$, $f_c \in \sE_0$ (cf. 17.16) and 
\[
K_f 
\ = \ 
K_{f_c }.
\]

[Note: \ 
Here
\[
f_c (z) 
\ = \ 
f (z + c).]
\]
\\[-1.25cm]
\end{x}

\qquad
{\small\bf \un{N.B.}} \ 
Therefore
\[
H_f 
\ = \ 
H_{f_c }
\]
or still, 
\[
h_f 
\ = \ 
h_{f_c }.
\]
\\[-1.25cm]

\begin{x}{\small\bf THEOREM} \ 
Suppose that $f \in \sE_0$ has the property that 
$h_f (\pm \sqrt{-1} \hsx) < \pi$.  
Assume further that $f (n) = 0$ for $n = 0, 1, 2, \ldots$ $-$then $f \equiv 0$.
\\[-.5cm]

PROOF \ 
\allowdisplaybreaks\begin{align*}
0 
\ = \ 
f (n) 
&= \ 
\frac{1}{2 \hsy \pi \sqrt{-1}} \ 
\int\limits_\Gamma \ 
\sB_\gf (w) 
\hsx
e^{n \hsy w} 
\ \td w
\qquad (\tcf. \ 18.19)
\\[4pt]
\implies
\hspace{2cm}
&
\\[4pt]
0 
&= \ 
\frac{1}{2 \hsy \pi \sqrt{-1}} \ 
\int\limits_\Gamma \ 
\sB_\gf (w) 
\hsx
\frac{1}{1 - z \hsy e^w}
\ \td w
\\[4pt]
\implies
\hspace{2cm}
&
\\[4pt]
0 
&= \ 
\frac{1}{2 \hsy \pi \sqrt{-1}} \ 
\int\limits_\Gamma \ 
\sB_\gf (w) 
\hsx
\frac{z}{1 - z \hsy e^w}
\ \td w
\\[4pt]
\implies
\hspace{2cm}
&
\\[4pt]
0 
&= \ 
-
\frac{1}{2 \hsy \pi \sqrt{-1}} \ 
\int\limits_\Gamma \ 
\sB_\gf (w) 
\hsx
e^{-w}
\ \td w
\qquad (z \ra \infty)
\\[4pt]
\implies
\hspace{2cm}
&
\\[4pt]
f (-1) \ 
&= \ 
0.
\end{align*}
Now apply the same argument to $f_{-1}$ to see that 
\[
f_{-1} (-1) 
\ = \ 
f (-2) 
\ = \ 
0.
\]
ETC.  One may then quote 19.14.
\\[-.5cm]

[Note: \ 
In view of 19.16, $\forall \ n$, $h_{f_n} (\pm \sqrt{-1} \hsx) < \pi$, and so $\forall \ w \in K_{f_n}$, 
\[
- \pi 
\ < \ 
- H_{f_n} (\sqrt{-1}) 
\ \leq \ 
\Img w 
\ \leq \ 
H_{f_n} (- \sqrt{-1}) 
\ < \
\pi,
\]
as follows from 18.8.]
\\[-.25cm]
\end{x}

\begin{x}{\small\bf FACT} \ 
$\forall \ f \in \sE_0$, 
\[
h_{f^\prime} 
\ \leq \ 
h_f.
\]

[In fact, 
\[
K_{f^\prime} \subset K_f
\implies 
H_{f^\prime} 
\ \leq \ 
H_f.]
\]
\\[-.25cm]
\end{x}

\chapter{
$\boldsymbol{\S}$\textbf{20}.\quad  DUALITY}
\setlength\parindent{2em}
\setcounter{theoremn}{0}
\renewcommand{\thepage}{\S20-\arabic{page}}

\qquad 
We shall provide here a description of the three standard realizations of the dual of the entire functions.  
\\[-.25cm]

\begin{x}{\small\bf NOTATION} \ 
$\sE$ is the set of entire functions.
\\[-.25cm]
\end{x}

By definition, the $C^0$-topology on $\sE$ is the topology of uniform convergence on compact subsets of $\Cx$.  
Denote its dual by $\sE^*$.  
Since $\sE$ is a closed subspace of $C^0 (\R^2)$, every continuous linear functional 
$\Lambda \in \sE^*$ extends to a continuous linear functional on $C^0 (\R^2)$, 
hence determines a compactly supported Radon measure.
\\[-.25cm]

\begin{x}{\small\bf DEFINITION} \ 
The elements of $\sE^*$ are called \un{analytic functionals}.
\\[-.25cm]
\end{x}

\begin{x}{\small\bf EXAMPLE} \ 
The compactly supported Radon measures 
\[
F \ra F (0)
\]
and
\[
F \ra 
\frac{1}{2 \hsy \pi \sqrt{-1}} \ 
\int\limits_{\abs{z} = 1} \ 
\frac{F (z)}{z} 
\ \td z
\]
restrict to the same analytic functional.
\\[-.25cm]
\end{x}

\begin{x}{\small\bf REMARK} \ 
The $C^0$-topology on $\sE$ coincides with the $C^\infty$-topology on $\sE$.  
Since $\sE$ is a closed subspace of $C^\infty (\R^2)$, every continuous linear functional 
$\Lambda \in \sE^*$ extends to a continuous linear functional on $C^\infty (\R^2)$, 
hence determines a compactly supported distribution. 
\\[-.5cm]

[Note: \ 
Recall that if $F_1, F_2, \ldots$ is a sequence in $\sE$ and if $F_n \ra F$ 
uniformly on compact subsets of $\Cx$, then $F_n^\prime \ra F^\prime$ 
uniformly on compact subsets of $\Cx$.]
\\[-.25cm]
\end{x}


\begin{x}{\small\bf NOTATION} \ 
$\sM_0$ is the set of compactly supported Radon measures on $\R^2$.
\\[-.25cm]
\end{x}

\begin{x}{\small\bf DEFINITION} \ 
Given $\mu \in \sM_0$, its \un{FL-transform} 
$\widehat{\mu}$ 
is defined by 
\[
\widehat{\mu} (z) 
\ = \ 
\int\limits \ 
e^{z \hsy w} 
\td \mu (w).
\]
\\[-1.25cm]
\end{x}

\begin{x}{\small\bf LEMMA} \ 
$\widehat{\mu}  (z)$ is an entire function of exponential type.
\\[-.5cm]

PROOF \ 
To see that $\widehat{\mu}$ is entire, simply observe that
\[
\frac{\td}{\td z} 
\hsx 
\widehat{\mu} (z) 
\ = \ 
\int\limits \ 
(w) 
\hsx 
e^{z \hsy w} 
\td \mu (w).
\]
Next choose $R \gg 0$: $\sptx \mu$ is contained in the circle of radius $R$ centered at the origin 
$-$then
\allowdisplaybreaks
\allowdisplaybreaks\begin{align*}
\abs{\widehat{\mu} (z)} \ 
&=\ 
\int \ 
\abs{e^{z \hsy w}} 
\hsx
\abs{\td \mu (w)} 
\\[15pt]
&=\ 
e^{R \hsy \abs{z}}
\ 
\int \ 
\
\abs{\td \mu (w)} .
\end{align*}
\\[-1cm]
\end{x}

\begin{x}{\small\bf NOTATION} \ 
Given 
$\mu, \nu \in \sM_0$, write $\mu \sim \nu$ if $\widehat{\mu} = \widehat{\nu}$.
\\[-.25cm]
\end{x}

\begin{x}{\small\bf LEMMA} \ 
$\mu \sim \nu$ iff $\forall \ F \in \sE$, 
\[
\langle F, \mu \rangle 
\ = \ 
\langle F, \nu \rangle.
\]
Therefore $\sim$ is an equivalence relation on $\sM_0$.
\\[-.25cm]
\end{x}

\begin{x}{\small\bf EXAMPLE} \ 
Take 
$\td \mu = \restr{\td z}{\Gamma}$, where $\Gamma$ is a circle $-$then
\[
\widehat{\mu} (z) 
\ = \ 
\int\limits_\Gamma \ 
e^{z \hsy w} 
\td w
\ = \ 
0.
\]
So $\mu \sim 0$ but $\mu \neq 0$.
\\[-.25cm]
\end{x}

\begin{x}{\small\bf NOTATION} \ 
Given $\mu \in \sM_0$, let $[\mu]$ be its associated equivalence class.
\\[-.25cm]
\end{x}


\begin{x}{\small\bf LEMMA} \ 
The arrow
\[
\sM_0 / \sim \ \  \ra  \ \sE_0
\]
that sends 
$[\mu]$ to $\widehat{\mu}$ is a  linear bijection. 
\\[-.5cm]

PROOF \ 
Injectivity is manifest while surjectivity is an application of 18.19.
\\[-.25cm]
\end{x}

\begin{x}{\small\bf RAPPEL} \ 
The arrow
\[
\sB : \sE_0 \ra \sH_0 (\infty)
\]
that sends $f$ to $\sB_f$ is a linear bijection (cf. 18.23 and 18.24).
\\[-.25cm]
\end{x}

\begin{x}{\small\bf NOTATION} \ 
Let $F \in \sE$.
\\[-.25cm]

\qquad \textbullet \quad
Given $f \in \sE_0$, put
\[
\langle F, f \rangle \ 
\ = \ 
\sum\limits_{n = 0}^\infty \ 
\frac{\gamma_n}{n!} 
\hsx
F^{(n)} (0)
\qquad (\gamma_n = f^{(n)} (0)).
\]

\qquad \textbullet \quad
Given $\Phi \in \sH_0 (\infty)$, put
\[
\langle F, \Phi \rangle \ 
\ = \ 
\frac{1}{2 \hsy \pi \sqrt{-1}} \ 
\int\limits_{\Gamma} \ 
\Phi (w) 
\hsx
\hsx F (w) 
\ \td w.
\]

\qquad \textbullet \quad
Given $[\mu] \in \sM_0 / \sim$, put
\[
\langle F, [\mu] \rangle \ 
\ = \ 
\int \ 
F (w) 
\ \td \mu (w) 
\qquad (= \langle F, \mu\rangle).
\]
\\[-1.25cm]
\end{x}

\begin{x}{\small\bf LEMMA} \ 
Each of these prescriptions defines an analytic functional.
\\[-.25cm]
\end{x}

\begin{x}{\small\bf LEMMA} \ 
Suppose given a triple 
$(f, \Phi, [\mu])$.  
Assume: \ 
$\Phi = \sB_f$ 
and 
$\widehat{\mu} = f$ 
$-$then these three data points give rise to the same analytic functional.
\\[-.5cm]

PROOF \ 
By definition (cf. 20.6), 
\allowdisplaybreaks
\allowdisplaybreaks\begin{align*}
\widehat{\mu} (z) 
&=\ 
\int \ 
e^{z \hsy w} 
\ \td \mu (w) 
\\[15pt]
&=\ 
\int \ 
\sum\limits_{n = 0}^\infty \ 
\frac{(z w)^n}{n !} 
\ \td \mu (w)
\\[15pt]
&= \ 
\sum\limits_{n = 0}^\infty \ 
\frac{\langle w^n, \mu \rangle}{n !}
\hsx 
z^n
\end{align*}
$\implies$
\allowdisplaybreaks
\allowdisplaybreaks\begin{align*}
\langle F, f \rangle \ 
&=\ 
\langle F, \widehat{\mu} \rangle
\\[15pt]
&= \ 
\sum\limits_{n = 0}^\infty \ 
\frac{\langle w^n, \mu \rangle}{n !}
\ 
F^{(n)} (0)
\\[15pt]
&= \ 
\Big\langle
\sum\limits_{n = 0}^\infty \ 
\frac{F^{(n)} (0)}{n !} \hsx w^n, \mu
\Big\rangle
\\[15pt]
&= \ 
\langle F, \mu \rangle
\\[15pt]
&= \ 
\langle F, [\mu] \rangle.
\end{align*}
On the other hand, 
\allowdisplaybreaks
\allowdisplaybreaks\begin{align*}
\langle F, \sB_{\widehat{\mu}}  \rangle \ 
&=\ 
\frac{1}{2 \hsy \pi \sqrt{-1}} \ 
\int\limits_{\Gamma} \ 
\sB_{\widehat{\mu}} (w) 
\hsx F (w) 
\ \td w \ 
\\[15pt]
&= \ 
\frac{1}{2 \hsy \pi \sqrt{-1}} \ 
\int\limits_{\Gamma} \ 
\sB_{\widehat{\mu}} (w) 
\ 
\sum\limits_{n = 0}^\infty \ 
\frac{F^{(n)} (0)}{n !}
\hsx 
w^n 
\ \td w \ 
\\[15pt]
&= \ 
\sum\limits_{n = 0}^\infty \ 
\frac{F^{(n)} (0)}{n !}\ 
\frac{1}{2 \hsy \pi \sqrt{-1}} \ 
\int\limits_{\Gamma} \ 
\sB_{\widehat{\mu}} (w) 
w^n 
\ \td w 
\\[15pt]
&= \ 
\sum\limits_{n = 0}^\infty \ 
\frac{F^{(n)} (0)}{n !} 
\hsy 
(\widehat{\mu})^{(n)} (0)
\qquad (\tcf. \ 18.19)
\\[15pt]
&= \ 
\sum\limits_{n = 0}^\infty \ 
\frac{(\widehat{\mu})^{(n)} (0)}{n !} 
\hsy 
F^{(n)} (0)
\\[15pt]
&= \ 
\langle F, \widehat{\mu} \rangle
\\[15pt]
&= \ 
\langle F, f \rangle.
\end{align*}
\\[-1.25cm]
\end{x}

\begin{x}{\small\bf SCHOLIUM} \ 
Each of the spaces
$\sE_0$, 
$\sH_0 (\infty)$, 
$\sM_0 / \sim$ 
can be viewed as 
$\sE^*$.
\\[-.5cm]

[Note: \ 
If $\Lambda \in \sE^*$, then there is a $\mu \in \sM_0$$\hsy:$ $\forall \ F \in \sE$, 
\[
\langle F, \Lambda \rangle 
\ = \ 
\langle F, \mu \rangle.
\]
And if $\nu \in \sM_0$ has the same property, then $\mu \sim \nu$ (cf. 20.9).]
\\[-.25cm]
\end{x}

\begin{x}{\small\bf EXAMPLE} \ 
Take $\mu = \delta_1$ $-$then 
$\widehat{\mu} (z) = e^z$ 
and 
$\ds\sB_{\widehat{\mu}} (w) = \frac{1}{w - 1}$.  
Here
\[
\langle F, \delta_1 \rangle 
\ = \ 
F (1)
\]
while
\allowdisplaybreaks
\allowdisplaybreaks\begin{align*}
\langle F, \widehat{\mu} \rangle  \ 
&=\ 
\sum\limits_{n = 0}^\infty \ 
\frac{(\widehat{\mu})^{(n)} (0)}{n !}
\hsx 
F^{(n)} (0)
\\[15pt]
&=\ 
\sum\limits_{n = 0}^\infty \ 
\frac{F^{(n)} (0)}{n !}
\\[15pt]
&=\
F (1)
\end{align*}
and
\allowdisplaybreaks\begin{align*}
\frac{1}{2 \hsy \pi \sqrt{-1}} \ 
\int\limits_{\Gamma} \ 
\sB_{\widehat{\mu}} (w) \hsx F (w) 
\ \td w \ 
&= 
\frac{1}{2 \hsy \pi \sqrt{-1}} \ 
\int\limits_{\Gamma} \ 
\frac{F (w)}{w - 1} 
\ \td w \ 
\\[15pt]
&=\
F (1).
\end{align*}
\end{x}


\chapter{
$\boldsymbol{\S}$\textbf{21}.\quad  FOURIER TRANSFORMS}
\setlength\parindent{2em}
\setcounter{theoremn}{0}
\renewcommand{\thepage}{\S21-\arabic{page}}

\qquad 
Working on the real axis, 
the sign convention of the Fourier transform of an $f \in \Lp^1(-\infty, \infty)$ is ``plus'':
\[
\hat{f} (x) 
\ = \ 
\frac{1}{\sqrt{2 \hsy \pi}} \ 
\int\limits_{-\infty}^\infty \ 
f (t) \hsy e^{\sqrt{-1} \hsx x \hsy t} 
\ \td t.
\]

\begin{spacing}{1.75}
[Note: \ 
From the point of view of harmonic analysis, 
the ambient Haar measure is 
$
\ds
\frac{1}{\sqrt{2 \hsy \pi}}
$
times Lebesgue measure.]
\\[-1cm]
\end{spacing}

\begin{x}{\small\bf LEMMA} \ 
Let $f \in \Lp^1 (-\infty, \infty)$ $-$then $\hat{f} (x)$ is a uniformly continuous function of $x$. 
\\[-.5cm]

PROOF \ 
Write
\allowdisplaybreaks
\allowdisplaybreaks\begin{align*}
\abs{\hat{f} (x + y)  - \hat{f} (x) } \ 
&= \ 
\frac{1}{\sqrt{2 \hsy \pi}} \ 
\bigg|\hsx
\int\limits_{-\infty}^\infty \ 
f (t) 
\hsx
e^{\sqrt{-1} \hsx x \hsy t}
\hsx
(e^{\sqrt{-1} \hsx y \hsy t} - 1)
\ \td t
\hsx\bigg|
\\[15pt]
&\leq \ 
\frac{1}{\sqrt{2 \hsy \pi}} \ 
\int\limits_{-\infty}^\infty \ 
\abs{f (t)} 
\hsx
\big|
\hsx
e^{\sqrt{-1} \hsx y \hsy t} - 1
\hsx 
\big|
\ \td t
\\[15pt]
&\leq \ 
\frac{1}{\sqrt{2 \hsy \pi}} \ 
\int\limits_{-\infty}^\infty \ 
\abs{f (t)} 
\hsx
(2 (1 - \cos y \hsy t))^{1/2}
\ \td t
\\[15pt]
&\leq \ 
\frac{1}{\sqrt{2 \hsy \pi}} \ 
\int\limits_{-\infty}^\infty \ 
\abs{f (t)} 
\hsx
2
\hsx
\big| 
\sin \Big(\frac{y \hsy t}{2}\Big)
\big|
\ \td t
\\[15pt]
&= \ 
\frac{2}{\sqrt{2 \hsy \pi}} \ 
\bigg[
\int\limits_{-\infty}^{-R} \ 
\ + \ 
\int\limits_{R}^\infty \ 
\ + \ 
\int\limits_{-R}^R \ 
\bigg]
\cdots
\\[15pt]
&\leq \ 
\frac{2}{\sqrt{2 \hsy \pi}} \ 
\bigg[
\int\limits_{-\infty}^{-R} \ 
\ + \ 
\int\limits_{R}^\infty \ 
\bigg]
\abs{f (t)} 
\ \td t
\ + \ 
\frac{1}{\sqrt{2 \hsy \pi}} \ 
\int\limits_{-R}^R \ 
\abs{f (t)} 
\abs{y \hsy t} 
\ \td t
\\[15pt]
&\leq \ 
\frac{2}{\sqrt{2 \hsy \pi}} \ 
\bigg[
\int\limits_{-\infty}^{-R} \ 
\ + \ 
\int\limits_{R}^\infty \ 
\bigg]
\abs{f (t)} 
\ \td t
\ + \ 
\frac{\abs{y}}{\sqrt{2 \hsy \pi}} 
\
R
\hsx
\int\limits_{-R}^R \ 
\abs{f (t)} 
\ \td t.
\end{align*}
Given $\varepsilon > 0$, choose $R$ large enough to render
\[
\frac{2}{\sqrt{2 \hsy \pi}} \ 
\bigg[
\int\limits_{-\infty}^{-R} \ 
\ + \ 
\int\limits_{R}^\infty \ 
\bigg]
\abs{f (t)} 
\ \td t
\ < \ 
\frac{\varepsilon}{2}.
\]
This done, choose $y$ small enough to render
\[
\frac{\abs{y}}{\sqrt{2 \hsy \pi}} \ 
R
\hsx
\int\limits_{-R}^R \ 
\abs{f (t)} 
\ \td t
\ < \ 
\frac{\varepsilon}{2}.
\]
So, with these choices, 

\[
\abs{\hat{f} (x + y)  - \hat{f} (x) }
\ < \ 
\varepsilon.
\]
\\[-.75cm]
\end{x}

\begin{x}{\small\bf EXAMPLE} \ 
Take 
$
\scalebox{1.11}{$f (t) = e^{-\abs{t}}$}
$ 
$-$then 
\[
\hat{f} (x) 
\ = \ 
\Big(
\frac{2}{\pi}
\Big)^{1/2}
\hsx 
\frac{1}{1 + x^2}.
\]
\\[-1.25cm]
\end{x}

\begin{x}{\small\bf EXAMPLE} \ 
Take 
$
\ds
\scalebox{1.11}{$f (t) = e^{-\frac{1}{2} \hsy t^2}$}
$ 
$-$then 
\[
\scalebox{1.11}{$\hat{f} (x) \ = \ e^{-\frac{1}{2} \hsy x^2}$}
.
\]
\\[-1.25cm]
\end{x}


\begin{x}{\small\bf EXAMPLE} \ 
Take 
$
\scalebox{1.11}{$f (t) = e^{-e^t} \hsy e^t$}
$ 
$-$then 

\[
\hat{f} (x) 
\ = \ 
\frac{1}{\sqrt{2 \hsy \pi}} \ 
\Gamma (1 + \sqrt{-1} \hsx x).
\]
\\[-.25cm]
\end{x}

\begin{x}{\small\bf NOTATION} \ 
Let 
\[
C_0 (-\infty, \infty)
\]
stand for the set of continuous functions $F$ on $\R$ such that 
\[
F(x) \ra 0 
\quad \text{as} \quad 
\abs{x} \ra \infty.
\]

[Note: \ 
When equipped with the supremum norm, 
$C_0 (-\infty, \infty)$
is a Banach algebra and 
$C_c (-\infty, \infty)$
is a dense subalgebra.]
\\[-.25cm]
\end{x}

\begin{x}{\small\bf RIEMANN-LEBESGUE LEMMA} \ 
Let $f \in \Lp^1 (-\infty, \infty)$ $-$then $\hat{f} \in C_0 (-\infty, \infty)$.
\\[-.25cm]
\end{x}

\qquad
{\small\bf \un{N.B.}} \ 
The arrow
\[
\Lp^1 (-\infty, \infty) \ra C_0(-\infty, \infty)
\]
that sends $f$ to $\hat{f}$ is a bounded linear transformation: 
\[
\normx{\hat{f}}_\infty 
\ = \ 
\sup\limits_{-\infty \hsy < \hsy x \hsy < \hsy \infty} \big|\hat{f} (x)\big| 
\ \leq \ 
\frac{1}{\sqrt{2 \hsy \pi}} \ \norm{f}_1.
\]
\\[-1.25cm]

\begin{x}{\small\bf REMARK} \ 
Not every $F \in C_0 (-\infty, \infty)$ is the Fourier transform of a function in  $\Lp^1 (-\infty, \infty)$.
\\[-.5cm]

[Consider the function defined for $x \geq 0$ by the rule

\[
F (x) \ = \ 
\begin{cases}
\ \ \
\ds
\frac{x}{e}
\hspace{1.cm} (0 \leq x \leq e)
\\[18pt]
\ 
\ds
\frac{1}{\log x}
\hspace{0.65cm} (x > e)
\end{cases}
\]
and put

\[
F(x) 
\ = \ 
-F(-x) 
\qquad (x \leq 0).]
\]
\\[-1.cm]
\end{x}

\begin{x}{\small\bf RAPPEL} \ 
Let \gA be a subalgebra of $C_0 (-\infty, \infty)$.  
Assume: 
\\[-.25cm]

\qquad \textbullet \quad
\gA is selfadjoint: \ $F \in \gA \implies \ovs{F} \in \gA$.
\\[-.25cm]

\qquad \textbullet \quad
\gA separates points: $\forall \ x, y \in \R$ with $x \neq y$, $\exists \ F \in \gA$: $F(x) \neq F(y)$.
\\[-.25cm]

\qquad \textbullet \quad
\gA vanishes at no point: \ $\forall \ x \in \R$, $\exists \ F \in \gA$: $F(x) \neq 0$.
\\[-.25cm]

Then $\gA$ is dense in $C_0 (-\infty, \infty)$.
\\[-.25cm]
\end{x}

\begin{x}{\small\bf NOTATION} \ 
Let 
\[
\gA (-\infty, \infty)
\]
stand for the set of all $\widehat{f}$ \  $(f \in \Lp^1 (-\infty, \infty))$.
\\[-.25cm]
\end{x}

\begin{x}{\small\bf LEMMA} \ 
$\gA (-\infty, \infty)$ is an algebra.
\\[-.5cm]

PROOF \ 
It is clear that 
$\gA (-\infty, \infty)$ 
is a vector space.  
If now $\hat{f}, \hsx \hat{g} \in \gA (-\infty, \infty)$, then
\[
\hat{f} \cdot \hat{g} 
\ = \ 
\frac{1}{\sqrt{2 \hsy \pi}}
\hsx 
(f * g)^{\widehat{\hsx}},
\]
the $*$ being convolution.
\\[-.25cm]
\end{x}

\begin{x}{\small\bf THEOREM} \ 
$\gA (-\infty, \infty)$ is dense in $C_0 (-\infty, \infty)$.
\\[-.5cm]

PROOF \ 
\\[-.5cm]

\qquad \textbullet \quad
$\gA (-\infty, \infty)$  is selfadjoint.
\\[-.5cm]

[Given $f \in \Lp^1 (-\infty, \infty)$, 
\allowdisplaybreaks
\allowdisplaybreaks\begin{align*}
(\ovhc{\widehat{f}\hsx}) (x) \ 
&= \ 
\frac{1}{\sqrt{2 \hsy \pi}}
\int\limits_{-\infty}^\infty \ 
\ovs{f (t)} 
\hsx e^{-\sqrt{-1} \hsx x \hsy t} 
\ \td t
\\[15pt]
&= \
\frac{1}{\sqrt{2 \hsy \pi}}
\int\limits_{-\infty}^\infty \ 
\ovs{f (-t)}
\hsx
\hsx e^{\sqrt{-1} \hsx x \hsy t} 
\ \td t
\\[15pt]
&= \
\widehat{g} (x)
\qquad (g (t) = \ovs{f (-t)}).]
\end{align*}

\qquad \textbullet \quad
$\gA (-\infty, \infty)$  separates points.
\\[-.5cm]

[In fact, 
\[
C_c^\infty (-\infty, \infty)
\subsetx 
\sS (-\infty, \infty)
\subsetx
\gA (-\infty, \infty).]
\]

\qquad \textbullet \quad
$\gA (-\infty, \infty)$  vanishes at no point (obvious).
\\[-.25cm]
\end{x}

\begin{x}{\small\bf THEOREM} \ 
If 
$f_1, f_2 \in \Lp^1 (-\infty, \infty)$ 
and if 
$\hat{f}_1 = \hat{f}_2$ 
everywhere, then 
$f_1 = f_2$ 
almost everywhere.
\\[-.25cm]
\end{x}

In general, the Fourier transform $\widehat{f}$ of $f$ need not belong to $\Lp^1 (-\infty, \infty)$.
\\[-.25cm]

\begin{x}{\small\bf EXAMPLE} \ 
Take 
\[
f (t) 
\ = \ 
\begin{cases}
\ 
1 
\hspace{0.5cm} 
(\abs{t} \leq 1)
\\[4pt]
\ 
0
\hspace{0.5cm} 
(\abs{t} > 1)
\end{cases}
.
\]
Then
\[
\widehat{f} (x) 
\ = \ 
\Big(\frac{2}{\pi}\Big)^{1/2} 
\hsx
\frac{\sin x}{x}
\]
is not in $\Lp^1 (-\infty, \infty)$.
\\[-.25cm]
\end{x}

Accordingly, it cannot be expected that Fourier inversion will hold on the nose.  
Still, there are summability results.
\\[-.25cm]

\begin{x}{\small\bf THEOREM} \ 
If 
$f \in \Lp^1 (-\infty, \infty)$, 
then for almost all $t$, 
\[
f (t) 
\ = \ 
\lim\limits_{R \ra \infty} \ 
\frac{1}{\sqrt{2 \hsy \pi}} \ 
\int\limits_{-R}^R \ 
\widehat{f} (x)
\hsx 
\Big(
1 - \frac{\abs{x}}{R}
\Big)
\hsx e^{-\sqrt{-1} \hsx t \hsy x}
\ \td x.
\]

[Note: \ 
This relation is also valid at every continuity point of $f$.]
\\[-.25cm]
\end{x}

\begin{x}{\small\bf REMARK} \ 
If 
$f \in \Lp^1 (-\infty, \infty)$, 
then as $R \ra \infty$, 
\[
\frac{1}{\sqrt{2 \hsy \pi}} \ 
\int\limits_{-R}^R \ 
\widehat{f} (x)
\hsx 
\Big(
1 - \frac{\abs{x}}{R}
\Big)
\hsx e^{-\sqrt{-1} \hsx t \hsy x}
\ \ra  \ 
f(t)
\]
in the $\Lp^1$ norm.
\\[-.25cm]
\end{x}

\begin{x}{\small\bf THEOREM} \ 
If 
$f \in \Lp^1 (-\infty, \infty)$ 
and if 
$\hat{f} \in \Lp^1 (-\infty, \infty)$,
then 
\[
f (t) 
\ = \ 
\frac{1}{\sqrt{2 \hsy \pi}} \ 
\int\limits_{-\infty}^\infty \ 
\widehat{f} (x)
\hsx e^{-\sqrt{-1} \hsx t \hsy x}
\ \td x
\]
almost everywhere.
\\[-.25cm]
\end{x}

\begin{x}{\small\bf THEOREM} \ 
If 
$f \in \Lp^1 (-\infty, \infty)$ 
and if 
$\widehat{f} \in \Lp^1 (-\infty, \infty)$, 
then 
\[
f (t) 
\ = \ 
\frac{1}{\sqrt{2 \hsy \pi}} \ 
\int\limits_{-\infty}^\infty \ 
\widehat{f} (x)
\hsx e^{-\sqrt{-1} \hsx t \hsy x}
\ \td x
\]
everywhere provided $f$ is continuous everywhere.
\\[-.25cm]
\end{x}

\begin{x}{\small\bf EXAMPLE} \ 
Take
\[
f (t) 
\ = \ 
\begin{cases}
\ 
1 - \abs{t} 
\hspace{0.5cm} 
(\abs{t} \leq 1)
\\[4pt]
\ 
\ 0
\hspace{1.35cm}
(\abs{t} > 1)
\end{cases}
.
\]
Then
\[
\widehat{f} (x)
\ = \ 
\frac{1}{\sqrt{2 \hsy \pi}} \ 
\frac{\sin^2 (x/2)}{(x/2)^2},
\]
so here the assumptions of 21.17 are met, thus $\forall \ t$, 

\allowdisplaybreaks\begin{align*}
\frac{1}{\sqrt{2 \hsy \pi}} \ 
\int\limits_{-\infty}^\infty \ 
\frac{1}{\sqrt{2 \hsy \pi}} \ 
\frac{\sin^2 (x/2)}{(x/2)^2}
\hsx e^{-\sqrt{-1} \hsx t \hsy x}
\ \td x \ 
&=\ 
\frac{1}{2 \hsy \pi} \ 
\int\limits_{-\infty}^\infty \ 
\frac{\sin^2 (x/2)}{(x/2)^2}
\hsx e^{\sqrt{-1} \hsx t \hsy x}
\ \td x \ 
\\[15pt]
&= \
\begin{cases}
\ 
1 - \abs{t} 
\hspace{0.5cm} 
(\abs{t} \leq 1)
\\[4pt]
\ 
\ 0
\hspace{1.35cm}
(\abs{t} > 1)
\end{cases}
.
\end{align*}
In particular: \ 
At $t = 0$, 
\[
\frac{1}{2 \hsy \pi} \ 
\int\limits_{-\infty}^\infty \ 
\frac{\sin^2 (x/2)}{(x/2)^2}
\ \td x \ 
\ = \ 
1
\]
\qquad 
$\implies$
\[
\int\limits_{-\infty}^\infty \ 
\frac{\sin^2 x}{x^2}
\ \td x \ 
\ = \ 
\pi.
\]
\\[-.75cm]
\end{x}

\begin{x}{\small\bf EXAMPLE} \ 
Take
\[
f (t)  \ = \ 
\begin{cases}
\ 
t \hsy e^{-t}
\hspace{0.5cm}
(t \geq 0)
\\[4pt]
\ 
\ 0
\hspace{.9cm}
(t < 0)
\end{cases}
.
\]
Then 
$f \in \Lp^1 (-\infty, \infty)$ .  
Moreover,
\[
\widehat{f} (x)
\ = \ 
\frac{1}{\sqrt{2 \hsy \pi}} \ 
\frac{1}{(1 - \sqrt{-1} \hsx x)^2}
\]
is also in 
$\Lp^1 (-\infty, \infty)$.  
Therefore at every $t$ (cf. 21.17),
\allowdisplaybreaks
\allowdisplaybreaks\begin{align*}
f (t) \ 
&=\ 
\frac{1}{\sqrt{2 \hsy \pi}} \ 
\int\limits_{-\infty}^\infty \ 
\widehat{f} (x)
\hsx
e^{- \sqrt{-1} \hsx t \hsy x}
\ \td x
\\[15pt]
&=\ 
\frac{1}{\sqrt{2 \hsy \pi}} \ 
\int\limits_{-\infty}^\infty \ 
\frac{1}{\sqrt{2 \hsy \pi}} \ 
\frac{1}{(1 + \sqrt{-1} \hsx x)^2}
e^{\sqrt{-1}  \hsx t \hsy x}
\ \td x
\\[15pt]
&=\ 
\widehat{\phi} (t),
\end{align*}
where
\[
\phi (x) 
\ = \ 
\frac{1}{\sqrt{2 \hsy \pi}} \ 
\frac{1}{(1 + \sqrt{-1} \hsx x)^2}.
\]
\\[-1cm]
\end{x}

\begin{x}{\small\bf THEOREM} \ 
If 
$f \in \Lp^1 (-\infty, \infty)$ 
is continuously differentiable and if 
$f^\prime \in \Lp^1 (-\infty, \infty)$, 
then $\forall \ x$, 
\[
\big({f}^\prime\big)^{\widehat{ }} (x)
\ = \ 
- \sqrt{-1} \hsx x \hsx \widehat{f} (x).
\]

PROOF \ 
Write
\[
f(x) - f(0) 
\ = \ 
\int\limits_{0}^x \ 
f^\prime (t)
\ \td t.
\]
Then 
\[
\begin{cases}
\ 
\ds
\lim\limits_{x \hsy \ra \hsy \infty} \ 
f(x) 
\ = \ 
f(0) + 
\int\limits_0^\infty \ 
f^\prime (t)
\ \td t
\ = \ 0
\\[26pt]
\ 
\ds
\lim\limits_{x \hsy \ra \hsy -\infty} \ 
f(x) 
\ = \ 
f(0) + 
\int\limits_0^{-\infty} \ 
f^\prime (t)
\ \td t
\ = \ 0
\end{cases}
,
\]
$f$ being $\Lp^1$.  
But for $x \neq 0$, 
\[
\int\limits_{-R}^R \ 
f (t) \hsx e^{\sqrt{-1} \hsx x \hsy t} 
\ \td t
\ = \ 
\frac{e^{\sqrt{-1} \hsx x \hsy t} }{\sqrt{-1} \hsx x} 
\ 
f (t) 
\ 
\hsx \bigg|_{t \hsy = \hsy -R}^{t \hsy = \hsy R} 
\ - \ 
\int\limits_{-R}^R \ 
\frac{e^{\sqrt{-1} \hsx x \hsy t} }{\sqrt{-1} \ x} 
\hsx 
f^\prime (t)
\ \td t.
\]
Therefore, upon letting $R \ra \infty$, we have
\[
\int\limits_{-\infty}^\infty \ 
f (t) \hsx e^{\sqrt{-1} \hsx x \hsy t} 
\ \td t
\ = \ 
-
\int\limits_{-\infty}^\infty \ 
\frac{e^{\sqrt{-1} \hsx x \hsy t} }{\sqrt{-1} \hsx x} 
\hsx 
f^\prime (t)
\ \td t
\]
\qquad 
$\implies$
\[
- \sqrt{-1} \hsx x \hsx \widehat{f} (x)
\ = \ 
\big({f}^\prime\big)^{\widehat{ }} (x)
\qquad (x \neq 0).
\]
This relation is also valid at $x = 0$.  
In fact, both sides are continuous and the LHS is zero at $x = 0$ whereas the RHS at $x = 0$ equals
\allowdisplaybreaks\begin{align*}
\int\limits_{-\infty}^\infty \ 
f^\prime (t) 
\ \td t\ 
&=\ 
f(\infty) - f(-\infty)
\\
&=\ 
0 - 0 
\\[4pt]
&=\ 
0.
\end{align*}

[Note: \ 
By iteration, if $f$ is continuously differentiable $n$ times and if 
$f^{(k)} \in \Lp^1 (-\infty, \infty)$ 
$(0 \leq k \leq n)$,  
then $\forall \ x$, 
\[
\big({f}^{(n)}\big)^{\widehat{ }} (x)
\ = \ 
(-\sqrt{-1} \hsx x)^n 
\hsx \hat{f} (x).]
\]
\\[-1.5cm]
\end{x}

\begin{x}{\small\bf RAPPEL} \ 
If $0 < A < \infty$, then
\[
\Lp^2 [-A, A] 
\subsetx
\Lp^1 [-A, A] 
\]
\\[-1cm]
but this is false if $A = \infty$: \ 
The function 
\[
f(x) 
\ = \ 
\frac{1}{1 + \abs{x}}
\]
is in 
$\Lp^2 (-\infty, \infty)$ but is not in 
$\Lp^1 (-\infty, \infty)$.
\\[-.5cm]
\end{x}

We shall now turn to the $\Lp^2$-theory of the Fourier transform.
\\[-.5cm]

\begin{x}{\small\bf PLANCHEREL THEOREM} \ 
If 
$f \in \Lp^1 (-\infty, \infty) \cap \Lp^2 (-\infty, \infty)$, 
then 
$\widehat{f} \in \Lp^2 (-\infty, \infty)$ 
and 
$\restr{\wedge }{\Lp^1 (-\infty, \infty) \cap  \Lp^2 (-\infty, \infty)}$ 
extends uniquely to an isometric isomorphism
\[
\wedge : 
\Lp^2 (-\infty, \infty) \ra  \Lp^2 (-\infty, \infty).
\]
It is of period 4 (i.e., 
$
\wedge^4 = \id
$) and has pure point spectrum 
$1, \sqrt{-1}, -1, -\sqrt{-1}$.
\\[-.5cm]

[Note: \ 
For the record, given 
$f_1, f_2 \in \Lp^2 (-\infty, \infty)$ , 
\[
\int\limits_{-\infty}^\infty \ 
f_1 (t)  \hsx \ov{f_2 (t)} 
\ \td t 
\ = \ 
\int\limits_{-\infty}^\infty \ 
\widehat{f_1} (x)  \hsx \ov{\widehat{f_2} (x)} 
\ \td x.
\]
In particular: 
$\forall \ f \in \Lp^2 (-\infty, \infty)$, 
\[
\norm{f}_2 
\ = \ 
\normx{\widehat{f} \hsx}_2.]
\]
\\[-1.5cm]
\end{x}

\qquad
{\small\bf \un{N.B.}} \ 
Computationally, if 
$f \in \Lp^2 (-\infty, \infty)$, then as $R \ra \infty$, 
\[
\frac{1}{\sqrt{2 \hsy \pi}} \ 
\int\limits_{-R}^R \ 
f (t) \hsx e^{\sqrt{-1} \hsx x \hsy t} 
\ \td t 
\ \lra \ 
\widehat{f} (x) 
\]
in the $\Lp^2$-norm and 
\\[-1.5cm]

\[
\frac{1}{\sqrt{2 \hsy \pi}} \ 
\int\limits_{-R}^R \ 
\widehat{f} (x) 
\hsx
e^{- \sqrt{-1} \hsx t \hsy x}
\ \td x 
\ \lra \ 
f(t)
\]
in the $\Lp^2$-norm.
\\

\begin{x}{\small\bf REMARK} \ 
Let 
\[
h_n (x) 
\ = \ (2^n \hsy n !)^{-1/2} 
\hsx
\pi^{-1/4} 
\hsx 
e^{-x^2 / 2} 
\hsx 
H_n (x), 
\]
where
\[
H_n (x)
\ = \ 
(-1)^n \hsx 
e^{x^2} 
\hsx 
\frac{\td^n}{\td x^n} 
\hsx
e^{-x^2 / 2} 
\]
is the $n^\nth$ Hermite polynomial (cf. 8.17) $(n \geq 0)$ $-$then $\{h_n\}$ is an orthonormal basis for 
$\Lp^2 (-\infty, \infty)$ 
and
\[
\wedge (h_n) 
\ = \ 
\widehat{h}_n 
\ = \ 
(\sqrt{-1})^n 
\hsx 
h_n.]
\]
\\[-1.25cm]
\end{x}

\begin{x}{\small\bf RAPPEL} \ 
If 
$f, g \in \Lp^2 (-\infty, \infty)$, 
then their convolution $f * g$ belongs to 
$C_0 (-\infty, \infty)$
and 
\[
\norm{f * g}_\infty
\ \leq \ 
\norm{f}_2 
\hsx
\norm{g}_2 .
\]

[Note: \ 
The same cannot be said if 
$f, g \in \Lp^1 (-\infty, \infty)$.  
For example, take 

\[
f (t) \ = \ 
\begin{cases}
\ 
\ds
\frac{1}{\sqrt{t}} 
\hspace{0.5cm} (0 < t < 1)
\\[8pt]
\ \
0 
\hspace{0.75cm} (t \leq 0 \ \text{or} \ t \geq 1)
\end{cases}
, 
\quad
g (t) \ = \ 
\begin{cases}
\ 
\ds
\frac{1}{\sqrt{1-t}}
\hspace{0.5cm} (0 < t < 1)
\\[8pt]
\hspace{0.75cm}
0 
\hspace{1cm} (t \leq 0 \ \text{or} \ t \geq 1)
\end{cases}
.
\]
Then
\[
(f * g) (1) 
\ = \ 
\int\limits_{-\infty}^\infty \ 
f (t) \hsy g(1 - t) 
\ \td t 
\ = \ 
\int\limits_0^1 \ 
\frac{\td t}{t}
\]
is undefined.]
\\[-.5cm]
\end{x}

Let 
$f, g \in \Lp^2 (-\infty, \infty)$ 
$-$then 
$f  \cdot  g \in \Lp^1 (-\infty, \infty)$ 
and
\[
\int\limits_{-\infty}^\infty \ 
f (t) \hsy g(t) 
\ \td t \ 
\ = \ 
\int\limits_{-\infty}^\infty \ 
\widehat{f} (x) \hsy \widehat{g}(-x) 
\ \td x.
\]
So, $\forall \ x_0$, 
\allowdisplaybreaks
\allowdisplaybreaks\begin{align*}
\int\limits_{-\infty}^\infty \ 
f (t) \hsy g(t) \hsx e^{\sqrt{-1} \hsx x_0 \hsy t}
\ \td t \ 
&= \ 
\int\limits_{-\infty}^\infty \ 
\widehat{f} (x) \hsy \widehat{g}(x_0-x) 
\ \td x
\\[15pt]
&= \ 
(\widehat{f} * \widehat{g}) (x_0)
\end{align*}
\qquad 
$\implies$
\[
\hsx (f \cdot g)^{\widehat{\ }}
\ = \ 
\frac{1}{\sqrt{2 \hsy \pi}} \ 
(\widehat{f} * \widehat{g}).
\]
\\[-.5cm]

\begin{x}{\small\bf THEOREM} \ 
$\gA (-\infty, \infty)$ consists precisely of the convolutions 
$F * G$, where 
$F, G \in \Lp^2  (-\infty, \infty)$.
\\[-.5cm]

PROOF \ 
Given 
$F, G \in \Lp^2  (-\infty, \infty)$, write 
\[
\begin{cases}
\ 
F = \widehat{f}
\\[15pt]
\ 
G = \widehat{g}
\end{cases}
\qquad 
(f, g \in \Lp^2  (-\infty, \infty)).
\]
Then
\[
F * G 
\ = \ 
\widehat{f} * \widehat{g}
\ = \ 
\sqrt{2 \hsy \pi} 
\hsx (f \cdot g)^{\widehat{\ }} \in \gA (-\infty, \infty).
\]
Conversely, every 
$\phi \in \Lp^1  (-\infty, \infty)$ 
is a product $f \cdot g$ with 
$f, g \in \Lp^2  (-\infty, \infty)$, 
thus matters can be turned around.
\\[-.5cm]

[Note: \ 
Let 
$f = \sqrt{\abs{\phi}}$ 
and take  
$g = \phi / \sqrt{\abs{\phi}}$ 
when $f$ is not zero but take $g = 0$ when $f = 0$.]
\\[-.25cm]
\end{x}

\begin{x}{\small\bf THEOREM} \ 
If 
$f \in \Lp^2 (-\infty, \infty)$, 
then for almost all $t$, 
\[
f (t) 
\ = \ 
\lim\limits_{R \ra \infty} \ 
\frac{1}{\sqrt{2 \hsy \pi}} \ 
\int\limits_{-R}^R \ 
\widehat{f} (x) \ 
\Big(
1 - \frac{\abs{x}}{R}
\Big)
\hsx
e^{- \sqrt{-1} \hsx t \hsy x}
\ \td x.
\]
\\[-1.25cm]
\end{x}

\begin{x}{\small\bf APPLICATION} \ 
If 
$f_1 \in \Lp^1 (-\infty, \infty)$
and
$f_2 \in \Lp^2 (-\infty, \infty)$ 
and if 
$\widehat{f}_1 = \widehat{f}_2$
almost everywhere, then 
$f_1 = f_2$ 
almost everywhere.
\\[-.5cm]

[Use the preceding result in conjunction with 21.14.]
\\[-.25cm]
\end{x}

\begin{x}{\small\bf LEMMA} \ 
Let 
$f \in \Lp^2 (-\infty, \infty)$
$-$then the restriction of $f$ to $[a,b]$ is $\Lp^2$, hence is $\Lp^1$, and 
\[
\int\limits_a^b \ 
f (t) 
\ \td t
\ = \ 
\frac{1}{\sqrt{2 \hsy \pi}} \ 
\int\limits_{-\infty}^\infty \ 
\widehat{f} (x) \ 
\frac{e^{-\sqrt{-1} \hsx b \hsx x} - e^{-\sqrt{-1} \hsx a \hsy x}}{-\sqrt{-1} \hsx x}
\td x.
\]

[If $\hsx \chisubab \hsx $ is the characteristic function of $[a,b]$, then 
\[
\whchisubab (x) 
\ = \ 
\frac{1}{\sqrt{2 \hsy \pi}} \ 
\frac{e^{\sqrt{-1} \hsx b \hsx x} - e^{\sqrt{-1} \hsx a \hsy x}}{\sqrt{-1} \hsx x}.]
\]
\\[-1.cm]
\end{x}

\begin{x}{\small\bf THEOREM} \ 
If 
$f \in \Lp^2 (-\infty, \infty)$ 
is continuously differentiable and if 
$f^\prime \in \Lp^2 (-\infty, \infty)$, 
then
\[
(f^\prime)^{\widehat{\ }} (x) 
\ = \ 
-\sqrt{-1} \hsx x \hsy \widehat{f} (x)
\]
almost everywhere (cf. 21.20).
\\[-.5cm]

PROOF \ 
Start by writing 
\[
f (t + h) - f(t) 
\ = \ 
\int\limits_t^{t + h} \ 
f^\prime (s) 
\ \td s.
\]
Next apply 21.28 to the integral on the right (replacing $f$ by $f^\prime$): 
\[
\int\limits_t^{t + h} \ 
f^\prime (s) 
\ \td s
\ = \ 
\frac{1}{\sqrt{2 \hsy \pi}} \ 
\int\limits_{-\infty}^\infty \ 
(f^\prime)^{\widehat{\ }} (x) 
\ 
\Big(
\frac{e^{-\sqrt{-1} \hsx h \hsy x} - 1}{-\sqrt{-1} \hsx x}
\Big)
\hsx
e^{-\sqrt{-1} \hsx t \hsy x} 
\ \td x.
\]
On the other hand, 
\[
f (t + h) - f(t) 
\ = \ 
\frac{1}{\sqrt{2 \hsy \pi}} \ 
\int\limits_{-\infty}^\infty \ 
\widehat{f} (x)
\ 
\big(
e^{-\sqrt{-1} \hsx h \hsy x} - 1
\big)
\hsx
e^{-\sqrt{-1} \hsx t \hsy x} 
\ \td x
\]
in the $\Lp^2$-sense.  
But 
\[
(f^\prime)^{\widehat{\ }} (x) \in \Lp^2 (-\infty, \infty), 
\quad
\frac{e^{-\sqrt{-1} \hsx h \hsy x} - 1}{-\sqrt{-1} \hsx x} \in \Lp^2 (-\infty, \infty)
\]
\qquad 
$\implies$
\[
(f^\prime)^{\widehat{\ }} (x) 
\ 
\Big(
\frac{e^{-\sqrt{-1} \hsx h \hsy x} - 1}{-\sqrt{-1} \hsx x}
\Big)
\in \Lp^1 (-\infty, \infty).
\]
Meanwhile
\[
\widehat{f} (x)
\hsx
\big(
e^{-\sqrt{-1} \hsx h \hsy x} - 1
\big) 
\in \Lp^2 (-\infty, \infty).
\]
Therefore (cf. 21.27)
\[
(f^\prime)^{\widehat{\ }} (x) 
\ 
\Big(
\frac{e^{-\sqrt{-1} \hsx h \hsy x} - 1}{-\sqrt{-1} \hsx x}
\Big)
\ = \ 
\widehat{f} (x)
\hsx
\big(
e^{-\sqrt{-1} \hsx h \hsy x} - 1
\big) 
\]
almost everywhere.  
Take $h = 1$ and $x \neq 2 \hsy \pi \hsy n$: 
\\[-.25cm]

\qquad 
$\implies$
\[
(f^\prime)^{\widehat{\ }} (x) 
\ = \ 
-\sqrt{-1} \hsx x \hsy \widehat{f} (x)
\]
almost everywhere.  
\\[-.5cm]

[Note: \ 
It follows that $x \hsy \widehat{f} (x)$ belongs to $\Lp^2 (-\infty, \infty)$.]
\\[-.25cm]
\end{x}

\[
\text{APPENDIX}
\]

Assuming that 
$
\ds
\nu > -\frac{1}{2}
$, 
take
\[
f_\nu (t) 
\ = \ 
0
\quad \text{if} \ \abs{t} \geq 1
\]
and take
\[ 
f_\nu (t) 
\ = \ 
(1 - t^2)^{\nu - \frac{1}{2}}
\quad \text{if} \ \abs{t} < 1.
\]
Then 
$f_\nu \in \Lp^1 (-\infty, \infty)$
and
\allowdisplaybreaks
\allowdisplaybreaks\begin{align*}
\widehat{f}_\nu (x) \ 
&=\ 
\Big(
\frac{2}{\pi}
\Big)^{1/2}
\int\limits_0^1 \
(1 - t^2)^{\nu - \frac{1}{2}}
\hsx 
\cos x t
\ \td t
\\[15pt] 
&=\ 
\Big(
\frac{2}{\pi}
\Big)^{1/2}
\
\sum\limits_{n = 0}^\infty \ 
\frac{(-1)^n x^{2 \hsy n}}{(2 n)!}
\ 
\int\limits_0^1 \
(1 - t^2)^{\nu - \frac{1}{2}}
\hsx 
t^{2 \hsy n} 
\ \td t
\\[15pt] 
&=\ 
\Big(
\frac{2}{\pi}
\Big)^{1/2}
\
\sum\limits_{n = 0}^\infty \ 
\frac{(-1)^n x^{2 \hsy n}}{(2 n)!}
\ 
\frac{1}{2} \ 
\int\limits_0^1 \
u^{n - \frac{1}{2}} 
\hsx 
(1 - u)^{\nu - \frac{1}{2}} 
\ \td u
\\[15pt] 
&=\ 
\frac{1}{\sqrt{2 \hsy \pi}} 
\
\sum\limits_{n = 0}^\infty \ 
\frac{(-1)^n x^{2 \hsy n}}{(2 n)!}
\ 
B \Big(n + \frac{1}{2}, \nu + \frac{1}{2}\Big)
\\[15pt] 
&=\ 
\frac{1}{\sqrt{2 \hsy \pi}} 
\
\sum\limits_{n = 0}^\infty \ 
\frac{(-1)^n x^{2 \hsy n}}{(2 n)!}
\ 
\frac{\Gamma \Big(n + \frac{1}{2}\Big) \hsx \Gamma \Big(\nu + \frac{1}{2}\Big)}{\Gamma (n + \nu + 1)}
\\[15pt] 
&=\ 
\frac{1}{\sqrt{2 \hsy \pi}} 
\hsx
\Gamma \Big(\nu + \frac{1}{2}\Big) 
\
\sum\limits_{n = 0}^\infty \ 
\frac{(-1)^n x^{2 \hsy n}}{(2 n)!}
\
\frac{\sqrt{\pi} (2 n)!}{2^{2 \hsy n} (n !)}
\hsx
\frac{1}{\Gamma (n + \nu + 1)}
\\[15pt] 
&=\ 
\frac{1}{\sqrt{2}} 
\hsx
\Gamma \Big(\nu + \frac{1}{2}\Big) 
\ 
\sum\limits_{n = 0}^\infty \ 
\frac
{
(-1)^n \hsx \Big(\ds\frac{x}{2}\Big)^{2 \hsy n}
}
{n ! \hsx \Gamma (n + \nu + 1)}
\\[15pt] 
&=\ 
\frac{1}{\sqrt{2}} 
\hsx
\Gamma \Big(\nu + \frac{1}{2}\Big) 
\hsx 
\Big(
\frac{x}{2}
\Big)^{-\nu}
\tJ_{\nu} (x) 
\qquad (\tcf. \ 2.29).
\end{align*}
\\[-.25cm]

\qquad
{\small\bf EXAMPLE} \ 
Take 
$
\ds
\nu = \frac{1}{2}
$
$-$then
\[
\tJ_{1/2} (x) \ 
\ = \ 
\Big(
\frac{2}{\pi}
\Big)^{1/2}
\hsx
\frac{\sin x}{\sqrt{x}},
\]
so
\allowdisplaybreaks\begin{align*}
\widehat{f}_{1/2} (x) \ 
&=\ 
\frac{1}{\sqrt{2}} 
\hsx
\Gamma (1) 
\hsx 
\Big(
\frac{x}{2}
\Big)^{-1/2}
\hsx
\tJ_{1/2} (x) 
\\[15pt]
&=\ 
\Big(
\frac{2}{\pi}
\Big)^{1/2}
\hsx
\frac{\sin x}{x},
\end{align*}
in agreement with 21.13.
\\

\qquad
{\small\bf LEMMA} \ 
If $\nu > 0$, then $f_\nu \in \Lp^2 (-\infty, \infty)$.
\\[-.25cm]

\qquad
{\small\bf \un{N.B.}} \ 
\[
f_0 \not\in \Lp^2 (-\infty, \infty).
\]


\chapter{
$\boldsymbol{\S}$\textbf{22}.\quad  PALEY-WIENER}
\setlength\parindent{2em}
\setcounter{theoremn}{0}
\renewcommand{\thepage}{\S22-\arabic{page}}

\qquad 
Let 
\[
\sE_0 (A)
\ = \ 
\{f \in \sE_0 : T (f) \leq A\},
\]
where $0 < A < \infty$.
\\[-.25cm]

\begin{x}{\small\bf NOTATION} \ 
$\PW (A)$ is the subset of $\sE_0 (A)$ consisting of those $f$ such that 
$\restr{f}{\R} \in \Lp^2 (-\infty, \infty)$.
\\[-.5cm]

[Note: \ 
The elements of $\PW (A)$ are called \un{Paley-Wiener functions}.]
\\[-.25cm]
\end{x}

\qquad
{\small\bf \un{N.B.}} \ 
The elements of $\PW (A)$  are bounded on the real axis (cf. 17.29) and 
\[
f (x) \ra 0 
\quad 
\text{as} 
\quad 
\abs{x} \ra \infty
\qquad (\tcf. \ 17.34).
\]
\\[-1.5cm]

\begin{x}{\small\bf LEMMA} \ 
$\PW (A)$ is a vector space.
\\[-.25cm]
\end{x}

\begin{x}{\small\bf LEMMA} \ 
$\PW (A)$ is an inner product space:
\[
\langle
f, g 
\rangle
\ = \ 
\int\limits_{-\infty}^\infty \ 
f (x) \hsx \ov{g (x)}
\ \td x.
\]
\\[-1.25cm]
\end{x}

\begin{x}{\small\bf LEMMA} \ 
$\PW (A)$ is closed under differentiation (cf. 17.8 and 17.31).  
\\[-.5cm]

[Note: \ 
If $f \in \PW (A)$, then
\[
\normx{f^\prime \hsy}_2 
\ \leq \ 
\norm{f}_2 \hsy T (f)
\ \leq \ 
\norm{f}_2 \hsy A.
\]
Therefore 
\[
\frac{\td}{\td z} : \PW (A) \ra \PW (A)
\]

\noindent
is a bounded linear transformation (but it is not surjective).]
\\[-.25cm]
\end{x}

\begin{x}{\small\bf CONSTRUCTION} \ 
Given $\phi \in \Lp^2 [-A, A]$ $(0 < A < \infty)$, put
\[
f (z) 
\ = \ 
\frac{1}{\sqrt{2 \hsy \pi}} \ 
\int\limits_{-A}^A \ 
\phi (t) 
\hsx
e^{\sqrt{-1} \hsx z \hsy t}
\ \td t.
\]
Then $f \in \sE_0(A)$ (cf. 17.19).  
Taking $z$ to be real and $\phi$ to be zero for $\abs{t} > A$, 
it follows that 
$\restr{f}{\R}= \widehat{\phi}$, 
thus by Plancherel
$\norm{\restr{f}{\R}}_2 = \norm{\phi}_2$, 
so $f \in \PW (A)$.  
Therefore this procedure determines an isometric injection
\[
\Lp^2 [-A, A] 
\hsx \ra \hsx  
\PW (A)
\qquad (\tcf. \ 21.11).
\]
\end{x}

\begin{x}{\small\bf EXAMPLE} \ 
Take
\[
\phi (t) 
\ = \ 
\frac{1}{\sqrt{1 - t^2}}
\qquad (-1 < t < 1).
\]
Then $\phi  \in \Lp^1 [-1, 1]$ but $\phi \notin \Lp^2 [-1, 1]$ .  
Moreover, 
\[
\int\limits_{-1}^1 \ 
\frac{e^{\sqrt{-1} \hsx x \hsy t}}{\sqrt{1 - t^2}}
\ \td t
\]
is not square integrable on the real axis.
\\[-.25cm]
\end{x}

\begin{x}{\small\bf THEOREM} \ 
The arrow 
\[
\Lp^2 [-A, A] 
\hsx \ra \hsx   
\PW (A)
\]
that sends $\phi$ to 
\[
f (z)
\ = \ 
\frac{1}{\sqrt{2 \hsy \pi}} \ 
\int\limits_{-A}^A \ 
\phi (t) 
\hsx
e^{\sqrt{-1} \hsx z \hsy t}
\ \td t
\] 
is an isometric isomorphism.
\\[-.5cm]

PROOF \ 
On the basis of what has been said above, it remains to establish surjectivity.  
If $T (f) = 0$, then $f = 0$ (cf. 17.30), so in this case we can take
$\phi = 0$.  
Assume now that $T (f) > 0$ $-$then
\[
\normx{f^\prime \hsy }_2 
\ \leq \
\norm{f}_2 \hsx T(f) 
\qquad (\tcf. \ 17.31),
\]
thus by iteration
\[
\normx{f^{(n)}}_2 
\ \leq \
\norm{f}_2 \hsx T(f)^n
\]
or still, passing to the Fourier transform (cf. 21.29), 
\[
\int\limits_{-\infty}^\infty \ 
x^{2 n} 
\hsx 
\big|\widehat{f} (x)\big|^2 
\ \td x 
\ \leq \
\normx{\widehat{f} \hsx }_2^2 \ T(f)^{2 n}
\qquad (n = 1, 2, \ldots).
\]
Fix $\varepsilon > 0$: 
\begin{align*}
(T (f) + \varepsilon)^{2 n} 
\hspace{-.25cm}
\int\limits_{\abs{x} \geq T(f) + \varepsilon} 
\big|\widehat{f} (x)\big|^2 
\ \td x  \
&\leq  
\int\limits_{\abs{x} \hsy \geq \hsy T(f) + \varepsilon}  
x^{2 n} 
\hsx 
\big|\widehat{f} (x)\big|^2 
\ \td x 
\\[15pt]
&\leq  \
\normx{\widehat{f} \hsx }_2^2 \ T(f)^{2 n}
\end{align*}
\qquad
$\implies$
\\[-1cm]

\[
\bigg[
\frac{T (f) + \varepsilon}{T (f)}
\bigg]^{2 n}
\ \times \ 
\int\limits_{\abs{x} \hsy \geq \hsy T(f) + \varepsilon} 
\big|\widehat{f} (x)\big|^2 
\ \td x
\ \leq \ 
\normx{\widehat{f} \hsx }_2^2 
\]
\qquad
$\implies$
\\[-1cm]

\[
\bigg[
1 + 
\frac{\varepsilon}{T (f)}
\bigg]^{2 n}
\ \times \ 
\int\limits_{\abs{x} \hsy \geq \hsy T(f) + \varepsilon} 
\big|\widehat{f} (x)\big|^2 
\ \td x
\ \leq \ 
\normx{\widehat{f} \hsx }_2^2 
\]
\qquad
$\implies$
\\[-1cm]

\[
\int\limits_{\abs{x} \hsy \geq \hsy T(f) + \varepsilon} 
\big|\widehat{f} (x)\big|^2 
\ \td x
\ = \ 
0
\qquad \text{(send $n$ to $\infty$)}.
\]
\\[-.75cm]

\noindent
Therefore $\widehat{f} (x) = 0$ almost everywhere if 
$\abs {x} \geq T (x) + \varepsilon$, hence $\widehat{f} (x) = 0$ almost 
everywhere if $\abs{x} \geq T (f)$.
Consequently, 
\[
\widehat{f} \in 
\Lp^2 [-T (f), T (f)] 
\subsetx
\Lp^2 [-A, A] .
\]
And for almost all $x$ (cf. 21.26),
\allowdisplaybreaks
\begin{align*}
f (x) \ 
&=\ 
\lim\limits_{R \ra \infty} \ 
\frac{1}{\sqrt{2 \hsy \pi}} \ 
\int\limits_{-R}^R \ 
\widehat{f} (t)
\hsx
\Big(
1 - \frac{\abs{t}}{R}
\Big)
\hsx
e^{-\sqrt{-1} \hsx x \hsy t}
\ \td t
\\[15pt]
&=\ 
\lim\limits_{R \ra \infty} \ 
\frac{1}{\sqrt{2 \hsy \pi}} \ 
\int\limits_{-A}^A \ 
\widehat{f} (t)
\hsx
\Big(
1 - \frac{\abs{t}}{R}
\Big)
\hsx 
e^{-\sqrt{-1} \hsx x \hsy t}
\ \td t
\\[15pt]
&=\ 
\frac{1}{\sqrt{2 \hsy \pi}} \ 
\int\limits_{-A}^A \ 
\widehat{f} (t)
\hsx
e^{-\sqrt{-1} \hsx x \hsy t}
\ \td t
\\[15pt]
&=\ 
\frac{1}{\sqrt{2 \hsy \pi}} \ 
\int\limits_{-A}^A \ 
\widehat{f} (-t)
\hsx
e^{\sqrt{-1} \hsx x \hsy t}
\ \td t
\\[15pt]
&=\ 
\frac{1}{\sqrt{2 \hsy \pi}} \ 
\int\limits_{-A}^A \ 
\phi (t) 
\hsx
e^{\sqrt{-1} \hsx x \hsy t}
\ \td t,
\end{align*}
where $\phi (t) = \widehat{f} (-t)$.  
But $f (z)$ is entire as is 

\[
\frac{1}{\sqrt{2 \hsy \pi}} \ 
\int\limits_{-A}^A \ 
\phi (t) 
\hsx
e^{\sqrt{-1} \hsx z \hsy t}
\ \td t.
\]
Since they agree almost everywhere on the real line, 
they must agree everywhere in the complex plane.
\\[-.25cm]
\end{x}

\begin{x}{\small\bf EXAMPLE} \ 
Let $f \in \sE_0 (A)$.  
Assume:  \ 
$\forall$ real $x$, 
\[
\abs{f (x)} 
\ \leq \ 
M.
\]
Then the function 
\[
\begin{cases}
\ \ds
\frac{f (z) - f (0)}{z}
\hspace{0.5cm} (z \neq 0)
\\[8pt]
\hspace{0.75cm}
f^\prime (0) 
\hspace{1.15cm}  (z = 0),
\end{cases}
,
\]
belongs to $\sE_0(A)$ and its restriction to the real axis is square integrable.  
Therefore
\[
f (z) 
\ = \ 
f (0) \hsx + \hsx 
\frac{z}{\sqrt{2 \hsy \pi}} \ 
\int\limits_{-A}^A \ 
\phi (t) 
\hsx
e^{\sqrt{-1} \hsx z \hsy t}
\ \td t
\]
for some $\phi \in \Lp^2[-A, A]$.
\\[-.25cm]
\end{x}

\begin{x}{\small\bf ADDENDUM} \ 
Assume that $\phi (t)$ does not vanish almost everywhere in any neighborhood of $A$ (or $-A$) 
$-$then $T(f) = A$ (hence $f$ is of order 1 (cf. 17.3)).  
\\[-.25cm]

[Suppose that $T (f) < A$, so $f \in \sE_0 (B)$ with $B < A$ $-$then
\[
f (z) 
\ = \ 
\frac{1}{\sqrt{2 \hsy \pi}} \ 
\int\limits_{-B}^B \ 
\psi (t) 
\hsx
e^{\sqrt{-1} \hsx z \hsy t}
\ \td t, 
\]
where 
$\psi \in \Lp^2[-B, B]$.  
Extend $\psi$ to $[-A, A]$ by taking it to be zero in 
\[
\begin{cases}
\ 
[-A, -B[ \hspace{0.5cm} (-A \leq t < -B)
\\[4pt]
\ 
]B, A] \hspace{1.25cm} (B < t \leq A)
\end{cases}
.
\]
Then still
\[
f (z) 
\ = \ 
\frac{1}{\sqrt{2 \hsy \pi}} \ 
\int\limits_{-A}^A \ 
\psi (t) 
\hsx
e^{\sqrt{-1} \hsx z \hsy t}
\ \td t.
\]
Accordingly, by the uniqueness of Fourier transforms (cf. 21.12), $\phi (t) = \psi (t)$ 
almost everywhere in $[-A, A]$.  
In particular: \ 
$\phi (t) = 0$ almost everywhere in 
\[
\begin{cases}
\ 
[-A, -B[ \hspace{0.5cm} (-A \leq t < -B)
\\[4pt]
\ 
]B, A] \hspace{1.15cm} (B < t \leq A)
\end{cases}
,
\]
a contradiction.]
\\[-.25cm]
\end{x}

\begin{x}{\small\bf THEOREM} \ 
Let $f \in \sE_0$ $(f \not\equiv 0$).  
Assume: \ 
$\restr{f}{\R} \in \Lp^2 (-\infty, \infty)$.  
Put
\[
\begin{cases}
\  \ds
\hspace{0.4cm}
b
\ = \ 
\underset{r \ra \infty}{\limsupx} \ 
\frac{\log \abs{f (- \sqrt{-1} \hsx r)}}{r}
\ \equiv \ 
h_f (- \sqrt{-1})
\\[18pt]
\  \ds
-a
\ = \ 
\underset{r \ra \infty}{\limsupx} \ 
\frac{\log \abs{f (\sqrt{-1} \hsx r)}}{r}
\ \equiv \ 
h_f (\sqrt{-1})
\end{cases}
.
\]
Then $b \geq a$ and 
\[
f (z) 
\ = \ 
\frac{1}{\sqrt{2 \hsy \pi}} \ 
\int\limits_a^b \ 
\phi (t) 
\hsx
e^{\sqrt{-1} \hsx z \hsy t}
\ \td t
\]
for some $\phi \in \Lp^2[a,b]$.
\\[-.5cm]

[Note: \ 
Since $f \not\equiv 0$, both $a$ and $b$ are finite (cf. 19.4).]
\\[-.25cm]
\end{x}

As will be seen below, this result is a consequence of 22.7 once the preliminaries are out of the way.
\\[-.25cm]

\begin{x}{\small\bf RAPPEL} \ 
If $A_1$, $A_2$ are nonempty sets of real numbers which are bounded above and if 
\[
A_1 + A_2 
\ = \ 
\{a_1 + a_2 : a_1 \in A_1, \hsx a_2 \in A_2\},
\]
then 
\[
\sup (A_1 + A_2) 
\ = \ 
\sup A_1 + \sup A_2.
\]
\\[-1.25cm]
\end{x}

\begin{x}{\small\bf LEMMA} \ 
Let $f \not\equiv 0$ be an entire function of exponential type $-$then
\[
h_f \big(\sqrt{-1} \hsx e^{\sqrt{-1} \hsx \theta} \big)
\hsx + \hsx 
h_f \big(- \sqrt{-1} \hsx e^{\sqrt{-1} \hsx \theta} \big)
\ \geq \ 
0.
\]

PROOF \ 
Work instead with $H_f$ (cf. 19.7).  
Put
\[
\begin{cases}
\ 
A_1 
\ = \ 
\{\Reg \big(\sqrt{-1} \hsx e^{\sqrt{-1} \hsx \theta} w_1\big)  : w_1 \in K_f\}
\\[4pt]
\ 
A_2 
\ = \ 
\{\Reg \big(-\sqrt{-1} \hsx e^{\sqrt{-1} \hsx \theta} w_2\big)  : w_2 \in K_f\}
\end{cases}
,
\]
so that by definition 
\[
\begin{cases}
\ 
H_f \big(\sqrt{-1} \hsx e^{\sqrt{-1} \hsx \theta}\big)
\ = \ 
\sup A_1
\\[11pt]
\ 
H_f \big(-\sqrt{-1} \hsx e^{\sqrt{-1} \hsx \theta}\big)
\ = \ 
\sup A_2
\end{cases}
.
\]
Consider now $A_1 + A_2$, a generic element of which has the form
\[
\Reg \big(\sqrt{-1} \hsx e^{\sqrt{-1} \hsx \theta} w_1\big)
\hsx +\hsx
\Reg \big(-\sqrt{-1} \hsx e^{\sqrt{-1} \hsx \theta} w_2\big).
\]
In particular: \ 
$\forall \ w \in K_f$, 
\[
\Reg \big(\sqrt{-1} \hsx e^{\sqrt{-1} \hsx \theta} w\big)
\hsx +\hsx
\Reg \big(-\sqrt{-1} \hsx e^{\sqrt{-1} \hsx \theta} w\big)
\ = \ 
0
\in A_1 + A_2.
\]
Therefore
\[
\sup (A_1 + A_2) 
\ \geq \ 
0
\]
\qquad \qquad
$\implies$
\[
\sup A_1 + \sup A_2
\ = \ 
\sup (A_1 + A_2) 
\ \geq \ 
0
\]
\qquad \qquad
$\implies$
\[
H_f \big(\sqrt{-1} \hsx e^{\sqrt{-1} \hsx \theta}\big)
\hsx +\hsx
H_f \big(-\sqrt{-1} \hsx e^{\sqrt{-1} \hsx \theta}\big)
\ \geq \ 
0.
\]
\\[-1.cm]
\end{x}

\begin{x}{\small\bf APPLICATION} \ 
Take $\theta = 0$ $-$then 
\[
h_f (\sqrt{-1}) + h_f (-\sqrt{-1}) 
\ \geq \ 
0,
\]
i.e., 
\[
h_f (-\sqrt{-1}) 
\ \geq \ 
- h_f (\sqrt{-1}) 
\]
or still, $b \geq a$.
\\[-.25cm]
\end{x}


\begin{x}{\small\bf P-L-P} \ 
Let $F$ be holomorphic in $\Img z > 0$ and continuous in $\Img z \geq 0$.  
Assume: \ 
\[
\log \abs{F (z)} 
\ = \ 
\tO (\abs{z}) 
\qquad (\abs{z} \gg 0)
\]
and
\[
\abs{F (x)} 
\ \leq \ 
M 
\qquad (-\infty < x < \infty)
\]
and
\[
\underset{r \ra \infty}{\limsupx} \ 
\frac{\log \abs{F (\sqrt{-1}\hsx r)}}{r}
\ = \ 
K.
\]
Then for $\Img z \geq 0$, 
\[
\abs{F (z)}
\ \leq \ 
M \hsy e^{K \hsy \Img z}.
\]
\\[-1cm]
\end{x}

Turning to the proof of 22.10, we have
\[
\begin{cases}
\ 
\abs{f (z)} 
\ \leq \
M \hsy e^{-a \hsy \Img z}
\hspace{0.5cm} (\Img z \geq 0)
\\[4pt]
\ 
\abs{f (z)} 
\ \leq \
M \hsy e^{b \hsy \abs{\Img z}}
\hspace{0.55cm} (\Img z \leq 0)
\end{cases}
.
\]
Put
\[
g (z) 
\ = \ 
e^{- \sqrt{-1} \hsx c \hsy z} f (z) 
\qquad \Big(c = \frac{a + b}{2}\Big).
\]
Then
\[
\abs{g (z) }
\ \leq \ 
M 
\hsx
\exp((1/2) (b - a) \abs{\Img z})
\]
\qquad \qquad
$\implies$
\[
g \in \sE_0 ((1/2) (b - a))
\]
if $b > a$ (cf. infra).  
Setting 

\[
C 
\ = \ 
(1/2) (b - a),
\]
it then follows from 22.7 that $\exists \ \psi \in \Lp^2 [-C, C]$: 

\[
g (z)
\ = \ 
\frac{1}{\sqrt{2 \hsy \pi}} \ 
\int\limits_{-C}^C \ 
\psi (t) 
\hsx
e^{\sqrt{-1} \hsx z \hsy t}
\ \td t
\]

\qquad 
$\implies$

\[
f (z)
\ = \ 
\frac{1}{\sqrt{2 \hsy \pi}} \ 
\int\limits_{-C}^C \ 
\psi (t) 
\hsx
e^{\sqrt{-1} \hsx z \hsy (t + c)}
\ \td t
\]
\qquad\qquad 
$\implies$
\[
f (z)
\ = \ 
\frac{1}{\sqrt{2 \hsy \pi}} \ 
\int\limits_a^b \ 
\phi (t) 
\hsx
e^{\sqrt{-1} \hsx z \hsy t}
\ \td t,
\]
where $\phi (t) = \psi (t - c)$.
\\[-.25cm]

[Note: \ 
If $a = b$, then $g$ is bounded, hence is a constant, call it $X$: 
\[
X 
\ = \ 
e^{-\sqrt{-1} \hsx c \hsy z}
\hsx
f (z)
\]
\qquad\qquad 
$\implies$
\[
f (x)
\ = \ 
X
\hsx 
e^{\sqrt{-1} \hsx c \hsy x}
\qquad (z = x + \sqrt{-1} \ 0)
\]
\qquad\qquad 
$\implies$
\[
\abs{f (x)} 
\ = \ 
X,
\]
an impossibility ($f \not\equiv 0$ and $\restr{f}{\R} \in \Lp^2 (-\infty, \infty)$).]
\\

\begin{x}{\small\bf REMARK} \ 
The indicator diagram $K_f$ of $f$ is a subset of 
$[\sqrt{-1} \hsx a, \sqrt{-1} \hsx b]$.
\\[-.5cm]

[Let $w \in K_f$ $-$then 
\[
-H_f (-1) 
\ \leq \ 
\Reg w 
\ \leq \ 
H_f (1)
\qquad (\tcf. \ 18.7)
\]
or still, 
\[
- h_f(-1) 
\ \leq \ 
\Reg w 
\ \leq \ 
h_f (1)
\qquad (\tcf. \ 19.7).
\]
But

\[
\begin{cases}
\  \ds
h_f (1)
\ = \ 
\underset{r \ra \infty}{\limsupx} \ 
\frac
{\log 
\big|
f \big(r \hsy e^{\sqrt{-1} \ 0}\big)
\big|
}
{r}
\\[18pt]
\  \ds
h_f (-1)
\ = \ 
\underset{r \ra \infty}{\limsupx} \ 
\frac
{\log 
\big|
f \big(r \hsy e^{\sqrt{-1} \hsx \pi}\big)
\big|
}
{r}
\end{cases}
.
\]
And

\[
\begin{cases}
\  \ds
\big|
f \big(r \hsy e^{\sqrt{-1} \ 0}\big)
\big|
\ = \ 
\abs{ f (r)} 
\ \leq \ 
M
\\[18pt]
\  \ds
\big|
f \big(r \hsy e^{\sqrt{-1} \hsx \pi}\big)
\big|
\ = \ 
\abs{ f (-r)} 
\ \leq \ 
M
\end{cases}
\]
\qquad\qquad 
$\implies$
\[
\begin{cases}
\ 
h_f (1)
\ \leq \ 
0
\\[4pt]
\ 
h_f (-1)
\ \leq \ 
0
\end{cases}
\]
\qquad\qquad 
$\implies$
\[
0 
\ \leq \ 
- h_f(-1) 
\ \leq \ 
\Reg w 
\ \leq \ 
h_f (1) 
\ \leq \ 
0
\qquad (\tcf. \ 18.9).
\]
Therefore $w$ is necessarily pure imaginary.  
Finally
\[
- H_f (\sqrt{-1}) 
\ \leq \ 
\Img w 
\ \leq \ 
H_f (- \sqrt{-1}) 
\qquad (\tcf. \ 18.8)
\]
or still, 
\[
- h_f (\sqrt{-1}) 
\ \leq \ 
\Img w 
\ \leq \ 
h_f (- \sqrt{-1}) 
\qquad (\tcf. \ 19.7.)
\]
\qquad\qquad 
$\implies$
\[
a 
\ \leq \ 
\Img w 
\ \leq \ 
b.]
\]

[Note: \ 
If $\phi (t)$ does not vanish in any neighborhood of $a$ and does not vanish
in any neighborhood of $b$, then 
\[
K_f 
\ = \ 
[\sqrt{-1}\hsx a, \sqrt{-1}\hsx b]\hsx.]
\]

The functions 
\[
\frac{1}{\sqrt{2\hsx A}} \hsx 
\exp \Big(- \frac{\sqrt{-1}\hsx t \hsy n \hsy \pi}{A}\Big)
\qquad (n = 0, \pm1, \ldots)
\]
constitute an orthonormal basis for $\Lp^2 [-A, A]$.  
Therefore the functions
\[
\frac{1}{\sqrt{2 \hsy \pi}} \hsx 
\frac{1}{\sqrt{2\hsx A}} \ 
\int\limits_{-A}^A \ 
\exp \Big(- \frac{\sqrt{-1}\hsx t \hsy n \hsy \pi}{A}\Big)
\hsx 
e^{\sqrt{-1} \hsx z \hsy t} 
\ \td t
\]
constitute an orthonormal basis for $\PW(A)$, 
i.e., the functions 
\[
\Big(
\frac{A}{\pi} \ 
\Big)^{1/2}
\ 
\frac{\sin (A z - n \pi)}{A z - n \pi}
\]
constitute an orthonormal basis for $\PW(A)$.
\\[-.5cm]

[Note: \ 
Matters simplify when $A = \pi$: \ 
The functions 
\[
\frac{\sin \pi (z - n)}{\pi (z - n)}
\]
constitute an orthonormal basis for $\PW(\pi)$.  
In this connection, observe that if $f (z)$ belongs to $\PW (A)$, then 
$
\ds
f \Big(\frac{\pi \hsy z}{A}\Big)
$ 
belongs to $\PW(\pi)$.]
\\[-.25cm]
\end{x}

\begin{x}{\small\bf THEOREM} \ 
Let $f \in \PW(A)$ $-$then there is an expansion
\[
f (z) 
\ = \ 
\sum\limits_{n \hsy = \hsy -\infty}^\infty \  
c_n 
\hsx 
\Big(
\frac{A}{\pi} \ 
\Big)^{1/2}
\
\frac{\sin (A z - n \pi)}{A z - n \pi}
\]
in $\PW(A)$, where
\[
c_n 
\ = \ 
\Big(
\frac{\pi}{A} \ 
\Big)^{1/2}
\hsx
f\Big(\frac{n \pi}{A}\Big), 
\]
so
\[
\norm{f}^2 
\ = \ 
\sum\limits_{n \hsy = \hsy -\infty}^\infty \  
\abs{c_n}^2 
\ = \ 
\frac{\pi}{A} \ 
\sum\limits_{n \hsy = \hsy -\infty}^\infty \  
\Big|
f\Big(\frac{n \pi}{A}\Big)
\Big|^2.
\]
\\[-1.cm]
\end{x}

\qquad
{\small\bf \un{N.B.}} \ 
Therefore
\[
f (z)
\ = \ 
\sum\limits_{n \hsy = \hsy -\infty}^\infty \  
f\Big(\frac{n \pi}{A}\Big)
\hsx
\frac{\sin (A z - n \pi)}{A z - n \pi}.
\]
\\[-1cm]

\begin{x}{\small\bf LEMMA} \ 
The series 
\[
\sum\limits_{n \hsy = \hsy -\infty}^\infty \  
f\Big(\frac{n \pi}{A}\Big)
\hsx
\frac{\sin (A z - n \pi)}{A z - n \pi}
\]
converges uniformly on every horizontal strip $\abs{\Img z} \leq h$. 
\\[-.25cm]
\end{x}

\begin{x}{\small\bf EXAMPLE} \ 
Take $A = \pi$ $-$then
\[
f (z)
\ = \ 
\sum\limits_{n \hsy = \hsy -\infty}^\infty \  
f (n) 
\hsx
\frac{\sin \pi (z - n)}{\pi (z - n)}.
\]
Accordingly, if $f (n) = 0$ for $n = 0, \pm 1,  \pm 2, \ldots$, then $f \equiv 0$ (cf. 19.14).
\\[-.25cm]
\end{x}

\begin{x}{\small\bf NOTATION} \ 
$\ell^2$ is the set of sequences 
$c_0, c_{\pm 1} c_{\pm 2}, \ldots$ of complex numbers such that 
\[
\sum\limits_{n \hsy = \hsy -\infty}^\infty \  
\abs{c_n}^2 
\ < \ 
\infty.
\]
\\[-1.25cm]
\end{x}

\begin{x}{\small\bf LEMMA} \ 
The arrow 
\[
\ell^2 \ra \PW (\pi)
\]
that sends $\{c_n\}$ to 
\[
f (z)
\ = \ 
\sum\limits_{n \hsy = \hsy -\infty}^\infty \  
c_n
\hsx
\frac{\sin \pi (z - n)}{\pi (z - n)}
\]
is an isometric isomorphism.
\\[-.25cm]
\end{x}


\begin{x}{\small\bf EXAMPLE} \ 
Put
\[
\begin{cases}
\ 
c_n = 0 
\hspace{1.5cm} (n \leq 0)
\\[4pt]
\ \ds
c_n = \frac{(-1)^n}{n}
\hspace{0.6cm} (n > 0)
\end{cases}
\]
and let 

\[
f (z) 
\ = \ 
\sum\limits_{n = 1}^\infty \ 
\frac{(-1)^n}{n} 
\hsx 
\frac{\sin \pi (z - n)}{\pi (z - n)}.
\]
\\[-.75cm]

\noindent
Then $f \in \PW(\pi)$, yet the product $z \hsy f (z)$ does not belong to $\PW(\pi)$ 
(but, of course, it does belong to $\sE_0 (\pi)$ (cf. 17.15)).
\\[-.5cm]

[If $z \hsy f (z)$ was a Paley-Wiener function, then it would be bounded on the real axis (cf. 17.29), 
thus the same would be true of its derivative $z f^\prime (z) + f (z)$ (cf. 17.24 (or quote 22.4)).  
But

\[
f^\prime (z) 
\ = \ 
\sum\limits_{n = 1}^\infty \ 
\frac{1}{n} 
\hsx
\frac
{
\pi^2 (z - n)
\hsx 
\cos \pi z - \pi \sin \pi z
}
{
\pi^2 (z - n)^2
}
\]
\qquad \qquad
$\implies$
\[
k\hsx f^\prime (k) 
\ = \ 
(-1)^k \ 
\sum\limits_{\substack{n = 1 \\ n \neq k}}^\infty \ 
\Big(
\frac{1}{n}  - \frac{1}{n-k} 
\Big)
\]
\qquad \qquad
$\implies$
\[
\abs{k \hsx f^\prime (k)}
\ = \ 
\Big(
1 + \frac{1}{2} + \cdots + \frac{1}{k} 
\Big)
- \frac{2}{k}
\]
\qquad \qquad
$\implies$
\[
\abs{k \hsx f^\prime (k)}
\ra \infty 
\qquad 
\text{as $k \ra \infty$.}
\]
However
\[
f (k) \ra 0
\qquad
\text{as $k \ra \infty$.}
\]
Therefore 
\[
\{k\hsx f^\prime (k) + f(k) : k = 1, 2, \ldots\}
\]
is not bounded.]
\\[-.25cm]
\end{x}

Moving on: \ 
\\[-.25cm]

\begin{x}{\small\bf LEMMA} \ 
$\forall$ real $x$, $y$: \ 
\[
\frac{\sin A (x - y)}{A (x - y)}
\ = \ 
\sum\limits_{n \hsy = \hsy -\infty}^\infty \ 
\frac{\sin (A x - n \pi)}{A x - n \pi}
\ \cdot \ 
\frac{\sin (A y - n \pi)}{A y - n \pi}.
\]
\\[-1cm]
\end{x}

\begin{x}{\small\bf APPLICATION} \ 
Let $f \in \PW(A)$ $-$then 
\[
f (x) 
\ = \ 
\frac{A}{\pi} \ 
\int\limits_{-\infty}^\infty \ 
f (y) 
\hsx 
\frac{\sin A (x - y)}{A (x - y)}
\ \td y.
\]

[Start with the RHS: 
\allowdisplaybreaks
\begin{align*}
\frac{A}{\pi}
\ 
\int\limits_{-\infty}^\infty \
& 
f (y) 
\frac{\sin A (x - y)}{A (x - y)}
\ \td y
\\[15pt]
&=\ 
\frac{A}{\pi}
\ 
\int\limits_{-\infty}^\infty \ 
f (y) 
\ 
\sum\limits_{n \hsy = \hsy -\infty}^\infty \ 
\frac{\sin (A x - n \pi)}{A x - n \pi}
\ \cdot \ 
\frac{\sin (A y- n \pi)}{A y - n \pi}
\ \td y
\\[15pt]
&=\ 
\sum\limits_{n \hsy = \hsy -\infty}^\infty \ 
\frac{A}{\pi}
\hsx
\Big(
\int\limits_{-\infty}^\infty \ 
f (y) 
\hsx 
\frac{\sin (A y- n \pi)}{A y - n \pi}
\ \td y
\Big)
\
\frac{\sin (A x - n \pi)}{A x - n \pi}
\\[15pt]
&=\ 
\sum\limits_{n \hsy = \hsy -\infty}^\infty \ 
\frac{A}{\pi}
\hsx
\Big(
\Big(\frac{\pi}{A}\Big)^{1/2}
\ 
\int\limits_{-\infty}^\infty \ 
f (y) 
\hsx 
\Big(\frac{A}{\pi}\Big)^{1/2}
\frac{\sin (A y- n \pi)}{A y - n \pi}
\ \td y
\Big)
\
\frac{\sin (A x - n \pi)}{A x - n \pi}
\\[15pt]
&=\ 
\sum\limits_{n \hsy = \hsy -\infty}^\infty \ 
\frac{A}{\pi}
\hsx
\Big(
\Big(\frac{\pi}{A}\Big)^{1/2}
\hsx
c_n
\Big)
\
\frac{\sin (A x - n \pi)}{A x - n \pi}
\\[15pt]
&=\ 
\sum\limits_{n \hsy = \hsy -\infty}^\infty \ 
\Big(\frac{A}{\pi}\Big)^{1/2}
\hsx
c_n
\
\frac{\sin (A x - n \pi)}{A x - n \pi}
\\[15pt]
&=\ 
\sum\limits_{n \hsy = \hsy -\infty}^\infty \ 
f \Big(\frac{n \pi}{A}\Big)
\hsx
\frac{\sin (A x - n \pi)}{A x - n \pi}
\\[15pt]
&=\ 
f(x).]
\end{align*}
\\[-1cm]

[Note: \ 
Consequently, 
\allowdisplaybreaks
\begin{align*}
\abs{f(x)} \ 
&\leq \ 
\frac{A}{\pi}
\ 
\int\limits_{-\infty}^\infty \ 
\abs{f (y)} 
\hsx 
\Big|\frac{\sin A (x - y)}{A (x - y)}\Big|
\ \td y
\\[15pt]
&\leq \ 
\frac{A}{\pi}
\ 
\Big(
\int\limits_{-\infty}^\infty \ 
\abs{f (y)}^2
\hsx
\ \td  y
\Big)^{1/2}
\Big(
\int\limits_{-\infty}^\infty \ 
\hsx
\Big|\frac{\sin A (x - y)}{A (x - y)}\Big|^2
\ \td  y
\Big)^{1/2}
\\[15pt]
&=\ 
\frac{A}{\pi}
\hsx 
\norm{f}_2
\hsx
\frac{1}{\sqrt{A}}
\Big(
\int\limits_{-\infty}^\infty \ 
\frac{\sin^2 y}{y^2}
\ \td  y
\Big)^{1/2}
\\[15pt]
&=\ 
\frac{A}{\pi}
\hsx 
\norm{f}_2
\hsx
\frac{1}{\sqrt{A}}
\hsx
\sqrt{\pi}
\qquad (\tcf. \ 21.18)
\\[15pt]
&=\ 
\Big(
\frac{A}{\pi}
\Big)^{1/2}
\hsx
\norm{f}_2.
\end{align*}
Moreover, this estimate is sharp: \ 
Take 
$\ A = \pi\ $, 
$n = 0$, 
$\ \ds f (z) = \frac{\sin \pi \hsy z}{\pi \hsy z}\ $ 
$-$then for real $x$, 
\[
\abs{f(x)}
\ \leq \ 
1 
\ = \ 
\norm{f}_2,
\]
and $f (0) = 1$.]
\\[-.25cm]
\end{x}

\begin{x}{\small\bf REMARK} \ 
the following result is of importance in sampling theory:
\[
\sum\limits_{n \hsy = \hsy -\infty}^\infty \ 
\abs{
\frac{\sin \pi (x - n)}{\pi (x - n)}
}^2
\ < \ 
2.
\]

[There is no loss of generality in imposing the restriction 
$
\ds
-\frac{1}{2} 
< x \leq 
\frac{1}{2} 
$,
hence
\allowdisplaybreaks
\begin{align*}
\sum\limits_{n \hsy = \hsy -\infty}^\infty \ 
\abs{\frac{\sin \pi (x - n)}{\pi (x - n)}}^2
&\leq\ 
1 + 
\sum\limits_{n \neq 0} \ 
\frac{1}{\pi^2 \hsx \abs{x - n}^2}
\\[15pt]
&\leq\  
1 + \frac{1}{\pi^2}
\
\sum\limits_{n = 1}^\infty \ 
\bigg[
\frac{1}{(n-x)^2}
\hsx + \hsx
\frac{1}{(n+x)^2}
\bigg]
\\[15pt]
&\leq\  
1 + \frac{1}{\pi^2}
\
\sum\limits_{n = 1}^\infty \ 
\raisebox{-.15cm}
{$
\Bigg[
$}
\hsx
\frac{1}{\Big(n - \ds\frac{1}{2}\Big)^2}
\hsx + \hsx
\frac{1}{\Big(n + \ds\frac{1}{2}\Big)^2}
\hsx
\raisebox{-.15cm}
{$
\Bigg]
$}
\\[15pt]
&=\  
1 + \frac{1}{\pi^2}
\hsx
\bigg[
\hsx
\sum\limits_{n = 1}^\infty \ 
\frac{1}{\Big(n +\ds\frac{1}{2}\Big)^2}
\hsx + \hsx
\frac{1}{\Big(\ds\frac{1}{2}\Big)^2}
\\[15pt]
&\hspace{3.5cm}
\hsx - \hsx
\sum\limits_{n = 2}^\infty \ 
\frac{1}{\Big(n - \ds\frac{1}{2}\Big)^2}
\hsx + \hsx
2 \ 
\sum\limits_{n = 2}^\infty \ 
\frac{1}{\Big(n - \ds\frac{1}{2}\Big)^2}
\hsx
\raisebox{-.15cm}
{$
\Bigg]
$}
\\[15pt]
&=\  
1 + \frac{1}{\pi^2}
\hsx
\bigg[
\hsx
2^2 + 2\ 
\sum\limits_{n = 2}^\infty \ 
\frac{1}{\Big(n - \ds\frac{1}{2}\Big)^2}
\hsx
\raisebox{-.15cm}
{$
\Bigg]
$}
\\[15pt]
&<\ 
1 + \frac{1}{\pi^2}
\hsx
\bigg[
\hsx
2^2 + 2\ 
\int\limits_1^\infty \ 
\frac{1}{\Big(t - \ds\frac{1}{2}\Big)^2}
\ \td t
\bigg]
\\[15pt]
&=\ 
1 + \frac{1}{\pi^2}
\
[2^2 + 2^2]
\\[15pt]
&=\ 
1 + 
2 \hsy
\Big(
\frac{2}{\pi}
\Big)^2
\\[15pt]
&<\
1 + 1
\\[15pt]
&=\ 
2.]
\end{align*}
\\[-1cm]
\end{x}

\begin{x}{\small\bf THEOREM} \ 
Let $f \in \sE_0 (A)$.  
Assume:  \ 
$\forall$ real $x$, 
\[
\abs{f (x)}
\ \leq \ 
M.
\]
Then 
\[
f (z) 
\ = \ 
f^\prime (0) 
\hsx 
\frac{\sin A \hsy z}{A} 
+ 
f(0) 
\hsx 
\frac{\sin A \hsy z}{A} 
+
\sum\limits_{n \neq 0} \ 
f \Big(\frac{n \hsy \pi}{A}\Big)
\hsx
\Big(\frac{A \hsy z}{n \hsy \pi}\Big)
\hsx
\frac{\sin (A z - n \pi)}{A z - n \pi}.
\]

PROOF \ 
Apply 22.16 to the function figuring in 22.7, hence
\[
\frac{f(z) - f(0)}{z} 
\ = \ 
f^\prime(0) 
\hsx 
\frac{\sin A \hsy z}{A \hsy z} 
+
\sum\limits_{n \neq 0} \ 
\frac
{
f \Big(\ds\frac{n \hsy \pi}{A}\Big) - f(0)
}
{\ds\frac{n \hsy \pi {\color{white}|}}{A}} 
\
\frac{\sin (A z - n \pi)}{A z - n \pi}
\]
\qquad
$\implies$
\begin{align*}
f (z) 
\ = \ 
f^\prime (0) 
\ 
\frac{\sin A \hsy z}{A} 
&+ 
f(0) 
+ 
\sum\limits_{n \neq 0} \ 
f \Big(\frac{n \hsy \pi}{A}\Big)
\hsx
\Big(\frac{A \hsy z}{n \hsy \pi}\Big)
\
\frac{\sin (A z - n \pi)}{A z - n \pi}
\\[15pt]
&+
(-f(0)) (\sin A \hsy z) \ 
\sum\limits_{n \neq 0} \ 
(-1)^n
\hsx
\Big(\frac{A \hsy z}{n \hsy \pi}\Big)
\hsx
\frac{1}{A z - n \pi}.
\end{align*}
But for $w$ nonintegral, 
\[
\frac{\pi}{\sin \pi \hsy w}
\ = \ 
\sum\limits_{n \hsy = \hsy -\infty}^\infty \ 
\frac{(-1)^n}{n + w} 
\ = \ 
\frac{1}{w} + 2 \hsy w \ 
\sum\limits_{n = 1}^\infty \ 
\frac{(-1)^n}{w^2 - n^2}.
\]
Therefore
\allowdisplaybreaks
\begin{align*}
\sum\limits_{n \neq 0} \ 
(-1)^n
\hsx
\Big(\frac{A \hsy z}{n \hsy \pi}\Big)
\hsx
\frac{1}{A z - n \pi} \ 
&=\ 
2 \hsy A \hsy z \ 
\sum\limits_{n = 1}^\infty \ 
\frac{(-1)^n}{A^2 z^2 - n^2 \pi^2}
\hspace{3cm}
\\[15pt]
&=\ 
2 \hsy A \hsy z \ 
\sum\limits_{n = 1}^\infty \ 
\frac{(-1)^n}{\pi^2 \hsy ((A \hsy z/\pi)^2 - n^2}
\\[15pt]
&=\ 
\frac{2 \hsy A \hsy z}{\pi^2} \ 
\sum\limits_{n = 1}^\infty \ 
\frac{(-1)^n}{(A \hsy z/\pi)^2 - n^2}
\\[15pt]
&=\ 
\frac{1}{\pi} \ 
2 \Big(\frac{A z}{\pi}\Big)
\ 
\sum\limits_{n = 1}^\infty \ 
\ds\frac{(-1)^n}
{\Big(\ds\frac{A z}{\pi}\Big)^2- n^2}
\\[15pt]
&=\ 
\frac{1}{\pi} \ 
\raisebox{-.25cm}
{$
\Bigg[
$}
\hsx
\ds\frac{\pi}{\sin \pi\Big(\ds\frac{A z}{\pi}\Big)} 
\hsx - \hsx 
\ds\frac{1}{\ds\frac{A z}{\pi}}
\hsx
\raisebox{-.25cm}
{$
\Bigg]
$}
\\[15pt]
&=\ 
\frac{1}{\pi} \ 
\bigg[
\frac{\pi}{\sin A z}
\hsx - \hsx 
\ds\frac{\pi}{A z}
\bigg]
\\[15pt]
&=\ 
\frac{1}{\sin A z} - \frac{1}{A z}.
\end{align*}
And so
\begin{align*}
f (0) + (- f(0)) 
&(\sin A z) \ 
\sum\limits_{n \neq 0} \ 
(-1)^n
\hsx
\Big(\frac{A \hsy z}{n \hsy \pi}\Big)
\hsx
\frac{1}{A z - n \pi} \ 
\\[15pt]
&=\ 
f (0) + (- f(0)) (\sin A z) \ 
\bigg[
\frac{1}{\sin A z} - \frac{1}{A z}
\bigg]
\\[15pt]
&=\ 
f (0) - f(0) + f(0) \ 
\frac{\sin A z}{A z}
\\[15pt]
&=\ 
f(0) \ 
\frac{\sin A z}{A z}.
\end{align*}
\\[-1cm]

Take $A = 1$ $-$then the functions 
\[
\frac{1}{\sqrt{\pi}} 
\ 
\frac{\sin (z - n \pi)}{z - n \pi} 
\]
constitute an orthonormal basis for $\PW (1)$ (the canonical choice \ldots) .
\\[-.25cm]
\end{x}

\begin{x}{\small\bf RAPPEL} \ 
Let 
\[
P_n (t) 
\ = \ 
\frac{1}{2^n \hsy n !} \ 
\frac{\td^n}{\td t^n} 
\hsx 
(t^2 - 1)^n
\]
be the $n^\nth$ Legendre polynomial (cf. 8.17) $-$then the functions
\[
\sqrt{n + \frac{1}{2}} \ P_n (t) 
\qquad (n = 0, 1, \ldots)
\]
constitute an orthonormal basis for $\Lp^2 [-1, 1]$.
\\[-.25cm]
\end{x}

\begin{x}{\small\bf LEMMA} \ 
We have 
\[
\frac{1}{\sqrt{2 \hsy \pi}} \ 
\int\limits_{-1}^1 \ 
P_n (t) 
\hsx
e^{\sqrt{-1} \hsx x \hsx t} 
\ \td t 
\ = \ 
(\sqrt{-1})^n 
\hsx
\frac{\tJ_{n + \frac{1}{2}} (x)}{\sqrt{x}}.
\]
\\[-1cm]
\end{x}

\begin{x}{\small\bf EXAMPLE} \ 
Take $n = 0$ $-$then $P_0 (t) = 1$ and 
\[
\frac{1}{\sqrt{2 \hsy \pi}} \ 
\int\limits_{-1}^1 \ 
P_0 (t) 
\hsx
e^{\sqrt{-1} \hsx x \hsx t} 
\ \td t 
\ = \ 
\Big(
\frac{2}{\pi}
\Big)^{1/2}
\
\frac{\sin x}{x}
\ = \ 
\frac{\tJ_{\frac{1}{2}} (x)}{\sqrt{x}}.
\]
\\[-1cm]
\end{x}

\begin{x}{\small\bf SCHOLIUM} \ 
The functions 
\[
\sqrt{n + \frac{1}{2}} 
\hsx
(\sqrt{-1})^n
\hsx
\frac{\tJ_{n + \frac{1}{2}} (z)}{\sqrt{z}}
\]
constitute an orthonormal basis for $\PW (1)$.
\\[-.25cm]
\end{x}

\begin{x}{\small\bf APPLICATION} \ 
Let 
\[
\phi_n (t) 
\ = \ 
\sqrt{n + \frac{1}{2}} 
\hsx
P_n (t).
\]
Then in $\Lp^2 [-1, 1]$, 
\[
\begin{cases}
\ \ds
\langle
e^{\sqrt{-1} \hsx x \hsx -}, \phi_n
\rangle
\ = \ 
\int\limits_{-1}^1 \ 
e^{\sqrt{-1} \hsx x \hsy t} \hsy \phi_n (t) 
\ \td t 
\ = \ 
\sqrt{2 \hsy \pi} \ \widehat{\phi}_n (x)
\\[26pt]
\ \ds
\langle
e^{\sqrt{-1} \hsx y \hsx -}, \phi_n
\rangle
\ = \ 
\int\limits_{-1}^1 \ 
e^{\sqrt{-1} \hsx y \hsy t} \hsy \phi_n (t) 
\ \td t 
\ = \ 
\sqrt{2 \hsy \pi} \ \widehat{\phi}_n (y)
\end{cases}
.
\]
Thus, by Parseval, 
\begin{align*}
\langle
e^{\sqrt{-1} \hsx x \hsx -}, e^{\sqrt{-1} \hsx y \hsx -}
\rangle \ 
&=\ 
\sum\limits_{n = 0}^\infty \ 
\langle
e^{\sqrt{-1} \hsx x -}, \phi_n
\rangle
\ 
\ov
{
\langle
e^{\sqrt{-1} \hsx y -}, \phi_n
\rangle
}
\\[15pt]
&=\ 
2 \hsy \pi \ 
\sum\limits_{n = 0}^\infty \ 
\widehat{\phi}_n (x) 
\hsx
\widehat{\phi}_n (-y).
\end{align*}
But
\allowdisplaybreaks
\begin{align*}
\langle
e^{\sqrt{-1} \hsx x \hsx -}, e^{\sqrt{-1} \hsx y \hsx -}
\rangle \ 
&=\ 
\int\limits_{-1}^1 \ 
e^{\sqrt{-1} \hsx (x - y) \hsy t} 
\ \td t 
\\[15pt]
&=\ 
2 \hsx
\frac{\sin (x - y)}{x - y}.
\end{align*}
On the other hand, 
\begin{align*}
&
2 \hsy \pi \ 
\sum\limits_{n = 0}^\infty \ 
\hsx 
\widehat{\phi}_n (x) 
\hsx 
\widehat{\phi}_n (-y) \ 
\\[15pt]
&=\ 
2 \hsy \pi \ 
\sum\limits_{n = 0}^\infty \ 
\sqrt{n + \frac{1}{2}} 
\
(\sqrt{-1})^n 
\hsx
\frac{\tJ_{n + \frac{1}{2}} (x)}{\sqrt{x}}
\hsx
\sqrt{n + \frac{1}{2}}
(\sqrt{-1})^n 
\hsx
\frac{\tJ_{n + \frac{1}{2}} (-y)}{\sqrt{-y}}
\qquad (\tcf. \ 22.27)
\\[15pt]
&=\ 
2 \hsy \pi \ 
\sum\limits_{n = 0}^\infty \ 
\Big(n + \frac{1}{2}\Big)
\hsx
(\sqrt{-1})^{2 \hsy n} 
\hsx
\frac{\tJ_{n + \frac{1}{2}} (x) }{\sqrt{x}}
\hsx
\frac{\tJ_{n + \frac{1}{2}} (-y) }{\sqrt{-y}}
\\[15pt]
&=\ 
2 \hsy \pi \ 
\sum\limits_{n = 0}^\infty \ 
\Big(n + \frac{1}{2}\Big) 
\hsx
(-\sqrt{-1})^{2 \hsy n} 
\hsx
(-1)^n
\hsx
\frac{\tJ_{n + \frac{1}{2}} (x) }{\sqrt{x}}
\hsx
\frac{\tJ_{n + \frac{1}{2}} (y) }{\sqrt{y}}.
\end{align*}
And
\begin{align*}
(\sqrt{-1})^{2 \hsy n} 
\hsx
(-1)^n
&=\ 
((\sqrt{-1})^2)^n
\hsx
(-1)^n
\\[11pt]
&=\ 
(-1)^n
\hsx
(-1)^n
\\[11pt]
&=\ 
(-1)^{2 \hsy n} 
\\[11pt]
&=\ 
1.
\end{align*}
Therefore

\[
\frac{\sin (x - y)}{x - y}
\ = \ 
\pi \ 
\sum\limits_{n = 0}^\infty \ 
\Big(
n + \frac{1}{2}
\Big) 
\hsx
\frac{\tJ_{n + \frac{1}{2}} (x) }{\sqrt{x}}
\hsx
\frac{\tJ_{n + \frac{1}{2}} (y) }{\sqrt{y}}.
\]
\\[-.25cm]
\end{x}

\chapter{
$\boldsymbol{\S}$\textbf{23}.\quad  DISTRIBUTION FUNCTIONS}
\setlength\parindent{2em}
\setcounter{theoremn}{0}
\renewcommand{\thepage}{\S23-\arabic{page}}

\qquad 
Suppose given a function 
$F : \R \ra \R$.
\\[-.25cm]

\begin{x}{\small\bf DEFINITION} \ 
$F$ is \un{increasing} if $F (x) \leq F(y)$ whenever $x \leq y$ and $F$ is 
\un{strictly increasing} if $F(x) < F(y)$ whenever $x < y$.
\\[-.25cm]
\end{x}

Suppose given an increasing function $F : \R \ra \R$.
\\[-.25cm]

\begin{x}{\small\bf NOTATION} \ 
Write
\[
\begin{cases}
\ 
F(x^+) 
\ = \ 
\lim\limits_{h \hsy \ra \hsy 0} \ 
F (x + h)
\\[18pt]
F(x^-) 
\ = \ 
\lim\limits_{h \hsy \ra \hsy 0} \ 
F (x - h)
\end{cases}
\qquad (h > 0)
\]
or still, 
\[
\begin{cases}
\ 
F(x^+) 
\ = \ 
\inf\limits_{y \hsy > \hsy x} \ 
F (y)
\\[18pt]
F(x^-) 
\ = \ 
\sup\limits_{y \hsy < \hsy x} \ 
F (y)
\end{cases}
\]
and put
\[
\begin{cases}
\ 
F(\infty) 
\ = \ 
\sup\limits_{x \hsy \in \hsy \R} \ 
F (x)
\\[18pt]
\ 
F(-\infty) 
\ = \ 
\inf\limits_{x \hsy \in \hsy \R} \ 
F (x)
\end{cases}
.
\]
\\[-.25cm]
\end{x}

\begin{x}{\small\bf DEFINITION} \ 
$F$ is \un{continuous from the right} 
if $\forall \ x$, 
\[
F (x^+) 
\ = \ 
F(x).
\]
\\[-1cm]

A 
\un{distribution function} 
is an increasing function 
$F : \R \ra \R$
which is continuous from the right subject to 
\[
F(\infty) 
\ = \ 
1, 
\qquad 
F(-\infty) 
\ = \ 
0.
\]
\\[-1.25cm]
\end{x}

\begin{x}{\small\bf EXAMPLE} \ 
The function 
\[
I (x) \ = \ 
\begin{cases}
\ 
0 
\hspace{0.5cm} 
(x < 0)
\\[4pt]
\ 
1
\hspace{0.5cm} 
(x \geq 1)
\end{cases}
\]
is a distribution function,  the \un{unit step function}.
\\[-.25cm]
\end{x}

\begin{x}{\small\bf DEFINITION} \ 
Suppose that $F$ is a distribution function.
\\[-.25cm]

\qquad \textbullet \quad
A point $x$ such that 
$F(x)$ $(= F(x^+)) = F(x^-)$ is called a \un{continuity point} of $F$.
\\[-.25cm]

\qquad \textbullet \quad
A point $x$ such that  
$F(x)$ $(= F(x^+)) \neq F(x^-)$ is called a \un{discontinuity point} of $F$.
\\[-.25cm]
\end{x}

\begin{x}{\small\bf DEFINITION} \ 
Suppose that $F$ is a distribution function $-$then the quantity
\[
j_x 
\ = \ 
F(x^+) - F(x^-)
\]
is called the \un{jump} of $F$ at $x$.
\\[-.5cm]

[Note: \ 
$j_x$ is positive at a discontinuity point and zero at a continuity point.]
\\[-.25cm]
\end{x}

\begin{x}{\small\bf LEMMA} \ 
The set
\[
\{x : j_x > 0\}
\]
is at most countable.
\\[-.25cm]
\end{x}

Therefore the set of continuity points of a distribution function is dense in $\R$.
\\[-.25cm]


\begin{x}{\small\bf REMARK} \ 
There exist distribution functions whose set of discontinuity points is dense in $\R$.
\\[-.5cm]

[Let 
$\{q_n: n = 1, 2, \ldots\}$ be an enumeration of $\Q$ and consider
\[
F(x) 
\ = \ 
\sum\limits_{q_n \leq x} \ 
2^{-n}, 
\]
noting that 
$
\ds
\sum\limits_{n = 1}^\infty \ 
2^{-n} = 1
$.]
\\
\end{x}

\begin{x}{\small\bf NOTATION} \ 
$\Bo (\R)$ is the $\sigma$-algebra of Borel subsets of $\R$.
\\[-.25cm]
\end{x}

\begin{x}{\small\bf LEMMA} \ 
If $f$ is a Lebesgue measurable function, then there exists a Borel measurable function $g$ such that $f = g$ almost everywhere.
\\[-.25cm]
\end{x}

\begin{x}{\small\bf CONSTRUCTION} \ 
Let $F$ be a distribution function $-$then there exists a unique Borel measure $\mu_F$ on $\R$ characterized by the condition
\[
\mu_F(\hsy ]a, b]) 
\ = \ 
F (b) - F (a)
\]
for all $a$, $b \in \R$.  
Here
\[
F (x) 
\ = \ 
\mu_F (\hsy ]-\infty, x])
\]
and 
\[
j_x 
\ = \ 
\mu_F(\{x\}).
\]
Moreover, 
\[
1 
\ = \ 
F (\infty) 
\ = \ 
\mu_F (\R),
\]
so $\mu_F$ is a probability measure on the line.
\\[-.5cm]

[Note: \
We have
\[
\begin{cases}
\ 
\mu_F ([a,b[\hsy) 
\ = \ 
F(b^-) - F(a^-)
\\[4pt]
\ 
\mu_F ([a,b]) 
\ = \ 
F(b) - F(a^-)
\\[4pt]
\mu_F (\hsy]a,b[\hsy) 
\ = \ 
F(b^-) - F(a)
\end{cases}
.]
\]
\\[-1.25cm]
\end{x}

\begin{x}{\small\bf EXAMPLE} \ 
Take $F = I$ $-$then $\mu_I = \delta_0$.
\\[-.25cm]
\end{x}

\begin{x}{\small\bf LEMMA} \ 
Any bounded Borel measurable function on $\R$ is $\mu_F$-integrable.
\\[-.25cm]
\end{x}

\begin{x}{\small\bf REMARK} \ 
The considerations in 23.11 can be reversed.  
For suppose that $\mu$ is a probability measure on the line.  
Put
\[
F_\mu (x) 
\ = \ 
\mu(\hsy ]-\infty, x]).
\]
Then $F_\mu$ is a distribution function and 
\[
\mu_{F_\mu} 
\ = \ 
\mu.
\]
In fact, 
\[
]a,b]
\ = \ 
]-\infty, b] \ - \ ]-\infty, a],
\]
thus
\begin{align*}
\mu_{F_\mu} (\hsy ]a,b]) \ 
&=\ 
F_\mu (b) - F_\mu (a)
\\[11pt]
&=\ 
\mu(\hsy]-\infty, b]) \hsx - \hsx  \mu(]-\infty, a])
\\[11pt]
&=\ 
\mu(\hsy]-\infty, b] \ - \  ]-\infty, a]) 
\\[11pt]
&=\ 
\mu(\hsy]a, b]).
\end{align*}

[Note: \ 
In the other direction, 
\[
F_{\mu_F} \ = \ F.]
\]
\\[-1.25cm]
\end{x}

There are three kinds of ``pure'' distribution functions, viz.: \ 
discrete, absolutely continuous, and singular.
\\[-.25cm]

\begin{x}{\small\bf DEFINITION} \ 
A distribution function $F$ is said to be \un{discrete} if there 
\\[.25cm]
is a sequence $\{x_n\} \subset \R$ (possibly finite) 
and positive numbers $j_n$ such that 
$
\ds
\sum\limits_n \ j_n = 1
$ 
and 
\[
F (x) 
\ = \ 
\sum\limits_n \ j_n \hsy I (x - x_n).
\]

[Note: \ 
Accordingly, 
\[
\mu_F 
\ = \ 
\sum\limits_n \ j_n \hsy \delta_{x_n}.]
\]
\\[-1.25cm]
\end{x}

\begin{x}{\small\bf LEMMA} \ 
Suppose that $F$ is a discrete distribution function $-$then a Borel measurable function $f$ is integrable with respect to 
$\mu_F$ 
iff
\[
\sum\limits_n \ j_n \hsy \abs{f(x_n)}
\ < \ 
\infty, 
\]
in which case
\[
\int \ 
f 
\ \td \mu_F 
\ = \ 
\sum\limits_n \ j_n \hsy f(x_n).
\]
\\[-1.cm]
\end{x}

\begin{x}{\small\bf RAPPEL} \ 
An increasing function 
$\phi : \R \ra \R$ 
is differentiable almost everywhere and its derivative $\phi^\prime$ is Lebesgue measurable, nonnegative, and
\[
\int\limits_a^b \ 
\phi^\prime (t) 
\ \td t
\ \leq \ 
\phi (b) - \phi (a)
\]
for all $a$ and $b$.
\\[-.25cm]
\end{x}

\begin{x}{\small\bf APPLICATION} \ 
Suppose that $F$ is a distribution function $-$then $F$ is differentiable almost everywhere 
and its derivative $F^\prime$ is Lebesgue measurable, nonnegative, and integrable: 
\[
\normx{F^\prime}_1 
\ = \ 
\int\limits_{-\infty}^\infty \ 
F^\prime (t) 
\ \td t
\ \leq \ 
F (\infty) - F (-\infty)
\ = \ 
1.
\]
\\[-1cm]
\end{x}

\begin{x}{\small\bf DEFINITION} \ 
A function 
$F : \R \ra \R$ 
is 
\un{absolutely continuous} 
if 
$\forall \ \varepsilon > 0$, $\exists \ \delta > 0$ 
such that for any finite set of disjoint intervals 
$]a_1, b_1[, \ldots, ]a_N, b_N[$, 
\[
\sum\limits_{j = 1}^N \ 
(b_j - a_j)
\ < \ 
\delta 
\ \implies \ 
\sum\limits_{j = 1}^N \ 
\abs{F (b_j) - F(a_j)} 
\ < \ 
\varepsilon. 
\]

[Note: \ 
An absolutely continuous function is necessarily uniformly continuous, the converse being false.]
\\[-.25cm]
\end{x}

\begin{x}{\small\bf EXAMPLE} \ 
If $F$ is everywhere differentiable and if $F^\prime$ is bounded, 
then $F$ is absolutely continuous (use the mean value theorem).
\\[-.25cm]
\end{x}

\begin{x}{\small\bf RAPPEL} \ 
If 
$f \in \Lp^1 (-\infty, \infty)$ and if 
$
\ds 
F(x) = \int\limits_{-\infty}^x \ f(t) \ \td t
$, 
then $F$ is absolutely continuous and $F^\prime = f$ almost everywhere.
\\[-.25cm]
\end{x}

\begin{x}{\small\bf EXAMPLE} \ 
The prescription 
\[
F (x) 
\ = \ 
\int\limits_{-\infty}^x \ 
\frac{1}{\sqrt{2 \hsy \pi}} \hsx e^{-y^2 / 2}
\ \td y
\
\]
defines an  absolutely continuous distribution function.
\\[-.25cm]
\end{x}

\begin{x}{\small\bf CRITERION} \ 
Suppose that $F$ is a distribution function $-$then $F$ is absolutely continuous iff 
$\mu_F$ is absolutely continuous with respect to the restriction of Lebesgue measure to $\Bo(\R)$.
\\[-.25cm]
\end{x}

So, under the assumption that $F$ is absolutely continuous, the Radon-Nikodym theorem 
implies that $\mu_F$ admits a density 
$f \in \Lp^1 (-\infty, \infty)\hsy :$
\[
\forall \ S \in \Bo (\R), 
\quad
\mu_F (S) \ = \ \int\limits_S \ f.
\]
\\[-.75cm]

Matters can then be made  precise.
\\[-.25cm]


\begin{x}{\small\bf THEOREM} \ 
If $F$ is an absolutely continuous distribution function, 
\\[.5cm]
\noindent
then
$\forall \ x$, 
$
\ds
F(x) = 
\int\limits_{-\infty}^x \ 
F^\prime (t)
\ \td t.
$

PROOF \ 
For $h > 0$, 
\[
\mu_F(\hsy ]x, x + h]) \ = \ 
\begin{cases}
\ 
F (x + h) - F(x)
\\[18pt]
\ 
\ds
\int\limits_x^{x + h} \ f
\end{cases}
\]
and 
\[
\mu_F(\hsy ]x - h, x]) \ = \ 
\begin{cases}
\ 
F (x) - F(x - h)
\\[18pt]
\ 
\ds
\int\limits_{x - h}^x \ f
\end{cases}
.
\]
But on general grounds, 
\[
\begin{cases}
\ 
\ds
\lim\limits_{h \ra 0} \ 
\frac{1}{h} \ 
\int\limits_x^{x + h}
f 
\ = \ 
f(x)
\\[26pt]
\ 
\ds
\lim\limits_{h \ra 0} \ 
\frac{1}{h} \ 
\int\limits_{x - h}^x \ 
f 
\ = \ 
f(x)
\end{cases}
\]
almost everywhere.  
Therefore
\[
\begin{cases}
\ \ds
\lim\limits_{h \ra 0} \ 
\frac{F (x + h) - F(x)}{h} 
\ = \ 
f(x)
\\[26pt]
\ \ds
\lim\limits_{h \ra 0} \ 
\frac{F (x) - F(x - h)}{h} 
\ = \ 
f(x)
\end{cases}
\]
almost everywhere, hence $F^\prime (x) = f(x)$ almost everywhere.  
Finally, $\forall \ x$, 
\[
F (x) 
\ = \ 
\mu_F(\hsy ]-\infty, x]) 
\ = \ 
\int\limits_{-\infty}^x \ 
f 
\ = \ 
\int\limits_{-\infty}^x \ 
F^\prime.
\]
\\[-1.25cm]
\end{x}


\begin{x}{\small\bf DEFINITION} \ 
An increasing continuous function $F : \R \ra \R$ is said to be \un{singular} if $F^\prime = 0$ almost everywhere.
\\[-.25cm]
\end{x}

Trivially, a constant function is singular.
\\[-.25cm]

\begin{x}{\small\bf EXAMPLE} \ 
There exist singular distribution functions.
\\[-.25cm]

[Let $\Theta$ denote the Cantor function on $[0,1]$ and put 
$
\begin{cases}
\ 
\Theta (x) = 0 \hspace{0.5cm} (x < 0)
\\[4pt]
\ 
\Theta (x) = 1 \hspace{0.5cm} (x > 1)
\end{cases}
$

\begin{tikzpicture}[scale=2.0]
 
\draw[<->] (-1,0) -- (5,0) node[right] {$$};
\node[label={{}}] at (-.13,-.25) {$0$};

\node[label={{}}] at (4,-.25) {$1$};
\draw[] (4,0) -- (4,0) node[] {$|$};     
\draw[] (0,2) -- (0,2) node[] {$-$}; 

\draw[->]  (0,-.5) -- (0,2.5) node[right] {$$};

\node[label={{}}] at (-.13,2) {$1$};

\draw[-, line width=0.3mm, harvardcrimson] (1.333,1) -- (2.667,1) node[right] {$$};
\node[label={{}}] at (1.333,-.25) {\textcolor{harvardcrimson}{\scriptsize $\frac{1}{3}$}};  
\draw[] (1.333,0) -- (1.333,0) node[] {\scriptsize $|$};   
\node[label={{}}] at (2.667,-.25) {\textcolor{harvardcrimson}{\scriptsize $\frac{2}{3}$}}; 
\draw[] (2.667,0) -- (2.667,0) node[] {\scriptsize $|$};  

\draw[-, line width=0.30mm,, darkviolet] (.444,0.5) -- (.889,0.5) node[right] {$$};
\node[label={{}}] at (.444,-.25) {\textcolor{darkviolet}{\scriptsize $\frac{1}{9}$}};
\draw[] (.444,0) -- (.444,0) node[] {\scriptsize $|$};       
\node[label={{}}] at (.889,-.25) {\textcolor{darkviolet}{\scriptsize $\frac{2}{9}$}}; 
\draw[] (.889,0) -- (.889,0) node[] {\scriptsize $|$};  
\draw[-, line width=0.30mm,, darkviolet] (3.111,1.5) -- (3.555,1.5) node[right] {$$};
\node[label={{}}] at (3.111,-.25) {\textcolor{darkviolet}{\scriptsize $\frac{7}{9}$}};    
\draw[] (3.111,0) -- (3.111,0) node[] {\scriptsize $|$};   
\node[label={{}}] at (3.555,-.25) {\textcolor{darkviolet}{\scriptsize $\frac{8}{9}$}}; 
\draw[] (3.555,0) -- (3.555,0) node[] {\scriptsize $|$};

\draw[-, line width=0.30mm,, officegreen] (.148,0.25) -- (.300,0.25) node[right] {$$};
\node[label={{}}] at (.148,-.25) {\textcolor{officegreen}{\scriptsize $\frac{1}{27}$}};  
\draw[] (.148,0) -- (.148,0) node[] {\scriptsize $|$};     
\node[label={{}}] at (.300,-.25) {\textcolor{officegreen}{\scriptsize $\frac{2}{27}$}}; 
\draw[] (.300,0) -- (.300,0) node[] {\scriptsize $|$};  
\draw[-, line width=0.30mm,, officegreen] (1.037,.75) -- (1.185,.75) node[right] {$$};
\node[label={{}}] at (1.037,-.25) {\textcolor{officegreen}{\scriptsize $\frac{7}{27}$}};   
\draw[] (1.037,0) -- (1.037,0) node[] {\scriptsize $|$};   
\node[label={{}}] at (1.185,-.25) {\textcolor{officegreen}{\scriptsize $\frac{8}{27}$}}; 
\draw[] (1.185,0) -- (1.185,0) node[] {\scriptsize $|$}; 

\draw[-, line width=0.30mm,, officegreen] (2.815,1.25) -- (2.963,1.25) node[right] {$$};
\node[label={{}}] at (2.815,-.25) {\textcolor{officegreen}{\scriptsize $\frac{19}{27}$}};   
\draw[] (2.815,0) -- (2.815,0) node[] {\scriptsize $|$};    
\node[label={{}}] at (2.963,-.25) {\textcolor{officegreen}{\scriptsize $\frac{20}{27}$}}; 
\draw[] (2.963,0) -- (2.963,0) node[] {\scriptsize $|$};  
\draw[-, line width=0.30mm,, officegreen] (3.704,1.75) -- (3.852,1.75) node[right] {$$};
\node[label={{}}] at (3.704,-.25) {\textcolor{officegreen}{\scriptsize $\frac{25}{27}$}};  
\draw[] (3.704,0) -- (3.704,0) node[] {\scriptsize $|$};     
\node[label={{}}] at (3.852,-.25) {\textcolor{officegreen}{\scriptsize $\frac{26}{27}$}}; 
\draw[] (3.852,0) -- (3.852,0) node[] {\scriptsize $|$};  

\draw[-, line width=0.30mm,, debianred] (.049,.125) -- (.099,.125) node[right] {$$}; 
\draw[-, line width=0.30mm,, debianred] (.346,.375) -- (.395,.375) node[right] {$$}; 
\draw[-, line width=0.30mm,, debianred] (.938,.625) -- (.988,.625) node[right] {$$}; 
\draw[-, line width=0.30mm,, debianred] (1.235,.875) -- (1.284,.875) node[right] {$$}; 

\draw[-, line width=0.30mm,, debianred] (2.716,1.125) -- (2.765,1.125) node[right] {$$}; 
\draw[-, line width=0.30mm,, debianred] (3.012,1.375) -- (3.062,1.375) node[right] {$$}; 
\draw[-, line width=0.30mm,, debianred] (3.605,1.625) -- (3.654,1.625) node[right] {$$}; 
\draw[-, line width=0.30mm,, debianred] (3.901,1.875) -- (3.951,1.875) node[right] {$$}; 

\draw[<-, line width=0.30mm,  blue] (-.5,0) -- (0,0) node[right] {$$};
\draw[->, line width=0.30mm, blue] (4,2) -- (4.5,2.) node[right] {$$};
\end{tikzpicture}
\\[-.25cm]

\noindent
$-$then 
$\Theta$ 
is a singular distribution function.  
Therefore
\[
\int\limits_0^1 \ 
\Theta^\prime (t) 
\ \td t 
\ = \ 0 
\ < \ 1 
\ = \ 
\Theta (1) - \Theta (0) 
\qquad 
(\tcf. \ 23.17).]
\]

[Note: \ 
The Cantor function is increasing on $[0,1]$ but there are refined versions of 
$\Theta$ 
that are strictly increasing on $[0,1]$.
The existence of strictly increasing singular functions lies deeper and such functions are more difficult to construct and the proofs that they are singular 
and strictly increasing more involved.   
A famous example, the ?$(x)$ function, was given by 
Minkowski\footnote[3]{\vspace{.11 cm} “Zur Geometrie der Zahlen,”\textit{Gesammelte Abhandlungen} \textbf{2} (1911), 50–51.}.  
Minkowski's objective was to establish a one-one correspondence between the
rational numbers of (0, 1) and the quadratic irrationals of (0, 1).   
It was 
Denjoy\footnote[4]{\vspace{.11 cm}''Sur une fonction r´eelle de Minkowski,'' \textit{J. Math. Pures Appl.} \textbf{17} (1938), 105–151.}
who first proved that ?$(x)$ is totally singular.  
A survey article by
Salem\footnote[5]{\vspace{.11 cm}"On some singular monotonic functions which are strictly increasing," 
\textit{Proc. Amer. Math. Soc.} \textbf{53} (1943), 427–439.}
discusses singular monotonic functions more generally.] 
\\[-.25cm]
\end{x}

\begin{x}{\small\bf LEMMA} \ 
And absolutely continuous distribution function $F$ cannot be singular.
\\[-.5cm]

PROOF \ 
For suppose $F$ was singular $-$then in view of 23.24, $\forall \ x$, 
\[
F (x)
\ = \ 
\int\limits_{-\infty}^x \ 
F^\prime (t) \ \td t
\ = \ 
0,
\]
an impossibility.
\\[-.25cm]
\end{x}

Given a distribution function $F$, 
let $\{x_n\}$ be its set of discontinuity points (which for this discussion we shall assume is not empty).  
Define
$\Phi : \R \ra \R$ 
by the prescription 
\[
\Phi (x) 
\ = \
\sum\limits_n \ 
j_{x_n} \hsy I (x - x_n).
\]
Then
$\Phi$ 
is increasing, continuous from the right, and 
\[
\Phi (-\infty) \ = \ 0, 
\quad
\Phi (\infty) \ \equiv \ a \leq 1.
\]
If $F \neq \Phi$, put
\[
\Psi (x) 
\ = \ 
F(x) - \Phi (x).
\]
Then $\Psi$ is increasing, continuous, and 
\[
\Psi (-\infty) \ = \ 0, 
\quad
\Psi (\infty) \ \equiv \ b \leq 1.
\]
\\[-1cm]

\begin{x}{\small\bf NOTATION} \ 
Let
\[
\begin{cases}
\ \ds
F_d (x) \ = \ \frac{1}{a} \hsy \Phi (x)
\\[8pt]
\ \ds
F_c (x) \ = \ \frac{1}{b} \hsy \Psi (x)
\end{cases}
.
\]
Therefore
$
\begin{cases}
\ F_d
\\[4pt]
\ F_c
\end{cases}
$
are distribution functions and 
\[
F 
\ = \ 
a \hsy F_d + b \hsy F_c 
\qquad (a + b = 1).
\]

[Note: \ 
$F_d$ is referred to as the discrete part of $F$ while $F_c$ is referred to as the continuous part of $F$.  
Here 
$0 \leq a \leq 1$, 
$0 \leq b \leq 1$, 
with the understanding that
\[
\begin{cases}
\ 
a = 1 \iff F = F_d
\\[4pt]
\ 
b = 1 \iff F = F_c
\end{cases}
.
\]
\\[-1.25cm]
\end{x}

\qquad
{\small\bf \un{N.B.}} \ 
More can be said about $F_c$ (cf. infra).
\\[-.25cm]

Given a continuous distribution function $F$, there are two possibilities: \ 
Either $F^\prime = 0$ almost everywhere (in which case $F$ is singular) 
or else $F^\prime \neq 0$ almost everywhere.  
Assuming that the second possibility is in force, define 
$\Phi : \R \ra \R$ 
by the prescription
\[
\Phi (x)
\ = \ 
\int\limits_{-\infty}^x \ 
F^\prime (t)
\ \td t.
\]
Then 
$\Phi$ 
is increasing, absolutely continuous, and 
\[
\Phi (-\infty) \ = \ 0, 
\quad
\Phi (\infty) \ \equiv \ u \leq 1.
\]
If $F \neq \Phi$, put
\[
\Psi (x) 
\ = \ 
F (x) - \Phi (x).
\]
Then $\Psi$ is increasing, continuous, and 
\[
\Psi (-\infty) \ = \ 0, 
\quad
\Psi (\infty) \ \equiv \ v \leq 1.
\]
In addition $\Phi^\prime = F^\prime$ almost everywhere, hence $\Psi^\prime = 0$ almost everywhere, 
hence $\Psi$ is singular.
\\[-.25cm]

\begin{x}{\small\bf NOTATION} \ 
Let 
\[
\begin{cases}
\ \ds
F_{a c} (x) \ = \ \frac{1}{u} \hsy \Phi (x)
\\[8pt]
\ \ds
F_s (x) \ = \ \frac{1}{v} \hsy \Psi (x)
\end{cases}
.
\]
Therefore
$
\begin{cases}
\ 
F_{a  c}
\\[4pt]
\ 
F_s
\end{cases}
$
are distribution functions and 
\[
F 
\ = \ 
u \hsy F_{a  c}
+ 
v \hsy F_s
\qquad 
(u + v = 1).
\]

[Note: \ 
$F_{a c}$ is  referred to as the absolutely continuous part of $F$ while $F_s$ is referred to as the singular part of $F$.  
Here 
$0 \leq u \leq 1$, 
$0 \leq v \leq 1$, 
with the understanding that
\[
\begin{cases}
\ 
u = 1 \iff F = F_{a c}
\\[4pt]
\ 
v = 1 \iff F = F_s
\end{cases}
.]
\]
\\[-1.25cm]
\end{x}

Now let $F$ be an arbitrary distribution function, thus
\[
F 
\ = \ 
a \hsy F_d
+ 
b \hsy F_c.
\]
Since $F_c$ is a continuous distribution function, the preceding discussion is
applicable to it.
Write
\[
\begin{cases}
\ 
\text{$F_{a c}$ in place of $(F_c)_{a c}$}
\\[4pt]
\ 
\text{$F_s$ in place of $(F_c)_s$}
\end{cases}
.
\]
Then
\[
F_c 
\ = \ 
u \hsy F_{a c} + v \hsy F_s
\]
\qquad \qquad
$\implies$
\[
F \ = \ 
a F_d + b(u \hsy F_{a c} + v \hsy F_s).
\]
And
\[
a + b\hsy u + b \hsy v 
\ = \ 
a + b 
\ = \ 1.
\]
\\[-1cm]

\begin{x}{\small\bf SCHOLIUM} \ 
Every distribution function $F$ admits a (unique) decomposition
\[
F 
\ = \ 
A \hsy F_d + B \hsy F_{a c} + C \hsy F_s, 
\]
where
\[
A + B + C 
\ = \ 
1
\qquad 
(A \geq 0, \ B \geq 0, \ C \geq 0),
\]
and $F_d$ is a discrete distribution function, $F_{a c}$  is an absolutely continuous distribution function, 
and $F_s$ is a singular distribution function.
\\[-.25cm]
\end{x}

\begin{x}{\small\bf DEFINITION} \ 
Let $F_1$, $F_2$ be distribution functions $-$then their \un{convolution} is the function 
\[
F_1 * F_2  (x) 
\ = \ 
\int\limits_{-\infty}^\infty \ 
F_1 (x - y) 
\ \td \mu_{F_2} (y).
\]
\\[-1.25cm]
\end{x}

\qquad
{\small\bf \un{N.B.}} \ 
The integral defining $F_1 * F_2$ exists (cf. 23.13).
\\[-.25cm]

\begin{x}{\small\bf LEMMA} \ 
The convolution $F_1 * F_2$ is a distribution function. 
\\[-.25cm]
\end{x}


\begin{x}{\small\bf FORMALITIES} \ 
We have 
\[
F_1 * F_2
\ = \ 
F_2 * F_1
\]
and 
\[
F_1 * (F_2 * F_3)
\ = \ 
(F_1 * F_2) * F_3.
\]
Furthermore, 
\[
F 
\ = \ F * I 
\ = \ I * F.
\]
\\[-1.25cm]
\end{x}

\begin{x}{\small\bf THEOREM} \ 
Suppose that $F = F_1 * F_2$.  
\\[-.25cm]

\qquad \textbullet \quad
If $F_1$, $F_2$ are discrete, then $F$ is discrete.
\\[-.25cm]

\qquad \textbullet \quad
If either $F_1$ or $F_2$ is continuous, then $F$ is continuous.
\\[-.25cm]

\qquad \textbullet \quad
If either $F_1$ or $F_2$ is absolutely continuous, then $F$ is absolutely continuous.
\\[-.25cm]

\qquad \textbullet \quad
If $F_1$ is discrete and $F_2$ is singular, then $F$ is singular.
\\[-.25cm]

\qquad \textbullet \quad
If $F_1$, $F_2$ are singular, then $F$ is continuous.
\\[-.25cm]

[Note: \ 
$F$ might be singular, or $F$ might be absolutely continuous, or $F$ might be a mixture of both.]
\\[-.25cm]
\end{x}

\begin{center}
APPENDIX
\\
\end{center}

An \un{integrator} 
is an increasing function 
$F : \R \ra \R$ 
which is continuous from the right.  
A distribution function is therefore an integrator but not conversely.  
\\[-.5cm]

Every integrator $F$ gives rise to a unique Borel measure $\mu_F$ characterized by the condition 
\[
\mu_F (\hsy ]a, b] \hsy)
\ = \ 
F(b) - F(a).
\]
\\[-.75cm]

\qquad
{\small\bf \un{N.B.}} \ 
Given integrators $F$ and $G$, 
$\mu_F = \mu_G$ iff $F - G$ is a constant.
\\[-.25cm]

\qquad
{\small\bf LEMMA} \ 
If $F$ is a continuously differentiable integrator, then 
$\td \mu_F (x) = F^\prime (x) \td x$.
\\[-.25cm]

\qquad
{\small\bf DEFINITION} \ 
The completion $\thickbar{\mu}_F$ of $\mu_F$ is called the 
\un{Lebesgue-Stieltjes measure} 
associated with $F$.
\\[-.25cm]

\qquad
{\small\bf EXAMPLE} \ 
Take $F (x) = x$ $-$then $\thickbar{\mu}_F$ is Lebesgue measure.
\\[-.25cm]

Denote by $\gA_F \supset \Bo (\R)$ the domain of $\thickbar{\mu}_F$.
\\[-.25cm]

\qquad
{\small\bf LEMMA} \ If $X \in \gA_F$, 
then there is a Borel set $S$ and a $Z \in \gA_F$ of Lebesgue-Stieltjes measure 0 such that 
$X = S \cup Z$.
\\[-.25cm]

Technically, one should distinguish between 
$
\ds
\int \hsx f \hsx \td \mu_F
$
and 
$
\ds
\int \hsx f \hsx \td \thickbar{\mu}_F
$
but this is unnecessary if $f$ is Borel measurable.
\\[-.25cm]

\qquad
{\small\bf NOTATION} \ 
Write 
$
\ds
\int\limits_a^b \ 
$
in place of 
$
\ds
\int\limits_{[a, b]}
$.
\\[.5cm]

\qquad
{\small\bf INTEGRATION BY PARTS} \ 
If $F$, $G$ are integrators, then 
\[
\int\limits_a^b \ 
G(x^+) 
\ \td \mu_F (x) 
\hsx + \hsx 
\int\limits_a^b \ 
F(x^-) 
\ \td \mu_G (x)
\ = \ 
F (b^+) \hsy G(b^+) - F(a^-) \hsy G(a^-).
\]
\\[-.5cm]

[Note: \ 
$G$ is continuous from the right so 
$G(x^+) = G(x)$ 
and 
$G(b^+) = G(b)$.]


\chapter{
$\boldsymbol{\S}$\textbf{24}.\quad  CHARACTERISTIC FUNCTIONS}
\setlength\parindent{2em}
\setcounter{theoremn}{0}
\renewcommand{\thepage}{\S24-\arabic{page}}

\qquad 
Let $F : \R \ra \R$ be a distribution function.
\\[-.5cm]

\begin{x}{\small\bf DEFINITION} \ 
The \un{characteristic function} $\gf$ of $F$ is the Fourier transform of $\mu_F$, i.e., 
\[
\gf (x) 
\ = \ 
\int\limits_{-\infty}^\infty \ 
e^{\sqrt{-1} \hsx x \hsy t} 
\ \td \mu_F (t).
\]

[Note: \ 
The integral defining $\gf$ exists (cf. 23.13).]
\\[-.75cm]
\end{x}

Obviously, 
\[
\gf (0) = 1, 
\quad 
\abs{\gf (x)} \leq 1, 
\quad
\ov{\gf (x)} =  \gf(-x).
\]
\\[-1cm]

\qquad
{\small\bf \un{N.B.}} \ 
We have
\[
\begin{cases}
\ \ds
\Reg \hsx \gf (x) 
\ = \ 
\int\limits_{-\infty}^\infty \ 
\cos (x \hsy t) 
\ \td \mu_F (t)
\\[26pt]
\ \ds
\Img \gf (x) 
\ = \ 
\int\limits_{-\infty}^\infty \ 
\sin (x \hsy t) 
\ \td \mu_F (t)
\end{cases}
.
\]
\\[-1cm]

\begin{x}{\small\bf LEMMA} \ 
$\gf (x)$ is a uniformly continuous function of $x$ (cf. 21.1).
\\[-.25cm]
\end{x}

\begin{x}{\small\bf DEFINITION} \ 
A distribution function $F : \R \ra \R$ is \un{symmetric} if $\forall \ x$, 
\[
\mu_F (]-\infty, x]) 
\ = \ 
\mu_F ([-x, \infty[ \hsx).
\]
Therefore
\[
\mu_F (S)
\ = \ 
\mu_F (-S)
\]
for all $S \in \Bo (\R)$.
\\[-.5cm]

[Note: \ 
Write
\[
]-\infty, -x[
\ \cup \
[-x, \infty[
\ \hsx = \hsx \ 
]-\infty, \infty[
\]
or still, 
\[
\big(
]-\infty, -x] 
- \{-x\}
\big)
\hsx \cup \hsx
[-x, \infty[
\ \hsx = \hsx \ 
]-\infty, \infty[ \hsx .
\]
Then
\[
\mu_F \big(\hsx ]-\infty, -x] - \{-x\}\big)
\hsx + \hsx 
\mu_F \big([-x, \infty[\hsx\big)
\ = \
\mu_F \big(\hsx ]-\infty, \infty[\hsx\big)
\]
\qquad \qquad
$\implies$
\[
\mu_F \big(\hsx ]-\infty, -x] \big)
\hsx - \hsx
\mu_F \big(\{-x\}\big)
\hsx + \hsx 
\mu_F \big([-x, \infty[\hsx\big)
\ = \
1
\]
\qquad \qquad
$\implies$
\[
F (-x) 
\hsx - \hsx
(F (-x) - F(-x^-)) 
\hsx + \hsx 
\mu_F \big([-x, \infty[\hsx\big)
\ = \
1
\]
\qquad \qquad
$\implies$
\[
 F(-x^-)
\hsx + \hsx 
\mu_F \big([-x, \infty[\hsx\big)
\ = \
1
\]
\qquad \qquad
$\implies$
\[
\mu_F \big([-x, \infty[\hsx\big)
\ = \
1 -  F(-x^-).
\]
Accordingly, $F$ is symmetric iff $\forall \ x$, 
\[
F (x) 
\ = \
1 -  F(-x^-).]
\]
\\[-1.5cm]
\end{x}

Given any distribution function $F$, 
the assignment $x \ra 1 - F (-x^-)$ is a distribution function, call it $(-1) \hsy F$, thus
\[
\td \mu_{(-1)F} (t)
\ = \ 
\td \mu_F (-t)
\]
and the characteristic function $(-1) \hsy \gf$ of $(-1) \hsy F$ is $\gf (-x)$ ($= \ov{\gf(x)}$).
\\[-.5cm]

[Note: \ 
$F$ is symmetric iff $F = (-1) \hsy F$.]
\\[-.25cm]

\begin{x}{\small\bf REMARK} \ 
$\Reg \hsx \gf (x)$ is a characteristic function.  
Proof: \ 

\[
\Reg \hsx \gf (x)
\ = \ 
\frac{1}{2} (\gf (x) + \ov{\gf (x)})
\]
and
\[
\frac{1}{2} \hsy F + \frac{1}{2} (-1) \hsy F
\]
is a distribution function.
\\[-.25cm]
\end{x}

\begin{x}{\small\bf LEMMA} \ 
$F$ is symmetric iff $\gf$ is real.
\\[-.5cm]

PROOF \ 
If $F$ is symmetric, then $\mu_F = \mu_{(-1)F}$, so
\allowdisplaybreaks
\begin{align*}
\gf (x) \ 
&=\ 
\int\limits_{-\infty}^\infty \ 
e^{\sqrt{-1} \hsx x \hsy t} 
\ \td \mu_F (t)
\\[15pt]
&=\ 
\int\limits_{-\infty}^\infty \ 
e^{-\sqrt{-1} \hsx x \hsy t} 
\ \td \mu_F (-t)
\\[15pt]
&=\ 
\int\limits_{-\infty}^\infty \ 
e^{-\sqrt{-1} \hsx x \hsy t} 
\ \td \mu_{(-1)F} (t)
\\[15pt]
&=\ 
\int\limits_{-\infty}^\infty \ 
e^{-\sqrt{-1} \hsx x \hsy t} 
\ \td \mu_F (t)
\\[11pt]
&=\ 
\gf (-x) 
\\[11pt]
&=\ 
\ov{\gf (x)}.
\end{align*}
I.e.: \ 
$\gf$ is real.  
\ 
Conversely, if $\gf$ is real, then $F$ and $(-1)F$ have the same characteristic function, 
hence 
$F = (-1)F$ (cf. 24.16).
\\[-.25cm]
\end{x}

\begin{x}{\small\bf LEMMA} \ 
We have 
\[
1 - \Reg \hsx \gf (2 x) 
\ \leq \ 
4 (1 - \Reg \hsx \gf (x))
\]
and
\[
\abs{\Img \gf (x)}
\ \leq \ 
\Big(
\frac{1}{2} 
\hsx 
(1 - \Reg \hsx \gf (2 x) )
\Big)^{1/2}.
\]

PROOF \ 
Write
\allowdisplaybreaks
\begin{align*}
1 - \Reg \hsx \gf (2 x)  \ 
&=\ 
\int\limits_{-\infty}^\infty \ 
(1 - \cos (2 \hsy x \hsy t)) 
\ \td \mu_F (t)
\\[15pt]
&=\ 
\int\limits_{-\infty}^\infty \ 
2 \hsy
(1 - (\cos (x \hsy t))^2 )
\ \td \mu_F (t)
\\[15pt]
&\leq \ 
\int\limits_{-\infty}^\infty \ 
4 (1 - \cos (x \hsy t))
\ \td \mu_F (t)
\\[15pt]
&=\ 
4 (1 - \Reg \hsx \gf (x))
\end{align*}
and
\begin{align*}
\abs{\Img \gf (x)}\ 
&=\ 
\bigg|
\int\limits_{-\infty}^\infty \ 
\sin (x \hsy t) 
\ \td \mu_F (t)
\bigg|
\\[15pt]
&\leq \ 
\bigg(
\int\limits_{-\infty}^\infty \ 
(\sin (x \hsy t) )^2
\ \td \mu_F (t)
\bigg)^{1/2}
\\[15pt]
&=\ 
\bigg(
\int\limits_{-\infty}^\infty \ 
\frac{1}{2}
(1 - \cos (2 \hsy x \hsy t))
\ \td \mu_F (t)
\bigg)^{1/2}
\\[15pt]
&=\ 
\Big(
\frac{1}{2}(1 - \Reg \hsx \gf (2 \hsy x))
\Big)^{1/2}.
\end{align*}
\\[-1cm]
\end{x}

\begin{x}{\small\bf REMARK} \ 
Elementary inequalities of this type (of which there are a number \ldots) 
can be used to preclude a function from being a characteristic function.  
E.g.: \ 
The function 
\[
\exp \big( - \abs{x}^\alpha \big)
\qquad (\alpha > 2)
\]
is not a characteristic function since the first inequality above is violated for small $x$.  
\\[-.5cm]

[Note: \ 
On the other hand, the function 
\[
\exp \big( - \abs{x}^\alpha \big)
\qquad (0 < \alpha \leq 2)
\]
is a characterisric function: 
\\[-.5cm]

\qquad \textbullet \quad
$0 < \alpha \leq 1)$ \ (apply 24.24)
\\[-.25cm]

\qquad \textbullet \quad
$\alpha = 2$ \ (immediate)
\\[-.25cm]

\qquad \textbullet \quad
$1 < \alpha < 2$ \ (trickier).]
\\[-.25cm]
\end{x}

\begin{x}{\small\bf ASYMTOTICS} \ 
Let $F$ be a distribution function, 
$\gf$ its characteristic function.
\\[-.5cm]

\qquad \textbullet \quad
Suppose that $F$ is discrete $-$then

\[
F (x) 
\ = \ 
\sum\limits_n \ 
j_n \hsy I (x - x_n)
\]
\qquad \qquad
$\implies$
\[
\mu_F 
\ = \ 
\sum\limits_n \ 
j_n \hsy \delta_{x_n}
\]
\qquad \qquad
$\implies$
\[
\gf (x) 
\ = \ 
\sum\limits_n \ j_n \hsy e^{\sqrt{-1} \hsx x \hsy x_n}
\]
\qquad \qquad
$\implies$
\[
\underset{\abs{x} \ra \infty}{\limsupx} \ 
\abs{\gf (x)}
\ = \ 1.
\]
\qquad\qquad \textbullet \quad
Suppose that $F$ is absolutely continuous $-$then 
$F^\prime \in \Lp^1 (-\infty, \infty)$ (cf. 23.18) and 
\[
F (x) 
\ = \ 
\int\limits_{-\infty}^x \ 
F^\prime (t) 
\ \td t
\qquad (\tcf. \ 23.24)
\]
\qquad \qquad
$\implies$
\begin{align*}
\gf (x) \ 
&=\ 
\int\limits_{-\infty}^\infty \
e^{\sqrt{-1} \hsx x \hsy t} 
\hsx
F^\prime (t) 
\ \td t
\\[15pt]
&\equiv \ 
\sqrt{2 \hsy \pi} 
\hsx (F^\prime)^{\widehat{\ }}
\end{align*}
\qquad \qquad
$\implies$
\[
\gf (x) \in C_0 (-\infty, \infty)
\qquad (\tcf. \ 21.6)
\]
\qquad \qquad
$\implies$
\[
\lim\limits_{\abs{x} \ra \infty} \ 
\abs{\gf (x)} 
\ = \ 
0.
\]

\qquad \textbullet \quad
Suppose that $F$ is singular $-$then as can be seen by example, 
\[
\underset{\abs{x} \ra \infty}{\limsupx} \ 
\abs{\gf (x)}
\]
might be 0 or it might be 1 or it might be between 0 and 1.
\\[-.25cm]

Put
\[
S (A) 
\ = \ 
\int\limits_{0}^A \ 
\frac{\sin t}{t} 
\ \td t
\qquad (A \geq 0).
\]
Then $S (A)$ is bounded and 
\[
\int\limits_{0}^A \ 
\frac{\sin t \hsy \theta}{t} 
\ \td t
\ = \ 
\sgn \theta \cdot S(A \abs{\theta}).
\]

[Note: \ 
Recall that

\[
\int\limits_0^\infty \ 
\frac{\sin t}{t} 
\ \td t
\ = \ 
\frac{\pi}{2}.]
\]
\\[-.25cm]
\end{x}

\begin{x}{\small\bf INVERSION FORMULA} \ 
Let $F$ be a distribution function, 
$\gf$ its characteristic function $-$then at any two continuity points 
$a < b$ of $F$,  
\[
F (b) - F(a) 
\ = \ 
\lim\limits_{A \ra \infty} \ 
\frac{1}{2 \hsy \pi}\ 
\int\limits_{-A}^A \ 
\frac
{
e^{-\sqrt{-1} \hsx a \hsy x} 
- 
e^{-\sqrt{-1} \hsx b \hsy x} 
}
{
\sqrt{-1} \hsx x
}
\
\gf (x)
\ \td x.
\]

PROOF \ 
Denoting by $I_A$ the entity inside the limit, insert
\[
\gf (x) \ 
\ = \ 
\int\limits_{-\infty}^\infty \
e^{\sqrt{-1} \hsx x \hsy t} 
\ \td \mu_F (t)
\]
and write
\[
I_A 
\ = \ 
\frac{1}{2 \hsy \pi}\ 
\int\limits_{-\infty}^\infty \ 
\bigg(
\int\limits_{-A}^A \ 
\frac
{
e^{\sqrt{-1} \hsx x \hsy (t - a)} 
- 
e^{\sqrt{-1} \hsx x \hsy (t - b)} 
}
{
\sqrt{-1} \hsx x
}
\ \td x
\bigg)
\ \td \mu_F (t)
\]
or still, 
\[
I_A 
\ = \ 
\int\limits_{-\infty}^\infty \ 
\bigg[
\frac{\sgn (t - a)}{\pi}
S(A \abs{t-a})
\hsx - \hsx
\frac{\sgn (t - b)}{\pi}
S(A \abs{t-b})
\bigg]
\ \td \mu_F (t).
\]
The integrand is bounded and converges as $A \ra \infty$ to the function 
\[
\phi_{a, b} (t) \ = \ 
\begin{cases}
\ \
0
\hspace{1cm}
(t < a)
\\[4pt]
\ \ds
1/2
\hspace{0.75cm}
(t = a)
\\[4pt]
\ \
1
\hspace{1.0cm}
(a < t < b)
\\[4pt]
\ \ds
1/2
\hspace{0.8cm}
(t = b)
\\[4pt]
\ \
0
\hspace{1.05cm}
(b < t)
\end{cases}
.
\]
Therefore
\begin{align*}
\lim\limits_{A \ra \infty} \ I_A \ 
&=\ 
\int\limits_{-\infty}^\infty \ 
\phi_{a, b} (t) 
\ \td \mu_F (t)
\\[15pt]
&=\ 
\frac{1}{2} \hsx \mu_F (\{a\}) + \mu_F (]a,b[\hsx) + \frac{1}{2} \hsx \mu_F (\{b\})
\\[15pt]
&=\ 
\frac{1}{2} \hsx (F (a) - F(a^-)) 
\hsx + \hsx 
(F (b^-) - F(a))
\hsx + \hsx 
\frac{1}{2} \hsx (F (b) - F(b^-)) 
\\[15pt]
&=\ 
F (b) - F(a).
\end{align*}
\\[-1.5cm]
\end{x}

\begin{x}{\small\bf REMARK} \ 
Using similar methods, $\forall \ a$, 
\[
j_a 
\ = \ 
\mu_F (\{a\}) 
\ = \ 
\lim\limits_{A \ra \infty} \
\frac{1}{2 A} \ 
\int\limits_{-A}^A \ 
e^{-\sqrt{-1} \hsx a \hsy x} 
\ 
\gf (x)
\ \td x.
\]
\\[-1.5cm]
\end{x}

\begin{x}{\small\bf THEOREM} \ 
If $\gf \in \Lp^1 (-\infty, \infty)$, then $F$ is continuous and its derivative $F^\prime$ exists.  
Moreover, 
\[
F^\prime (t) 
\ = \ 
\frac{1}{2 \pi} \ 
\int\limits_{-\infty}^\infty \ 
e^{-\sqrt{-1} \hsx t \hsy x} 
\ 
\gf (x)
\ \td x,
\]
hence is continuous.
\\[-.5cm]

PROOF \ 
Since $\gf \in \Lp^1 (-\infty, \infty)$, the same is true of 
\[
\frac
{
e^{-\sqrt{-1} \hsx a \hsy x} 
- 
e^{-\sqrt{-1} \hsx b \hsy x} 
}
{
\sqrt{-1} \hsx x
}
\ 
\gf (x), 
\]
so per 24.9, 
\[
F (b) - F(a) 
\ = \ 
\frac{1}{2 \hsy \pi}\ 
\int\limits_{-\infty}^\infty \ 
\frac
{
e^{-\sqrt{-1} \hsx a \hsy x} 
- 
e^{-\sqrt{-1} \hsx b \hsy x} 
}
{
\sqrt{-1} \hsx x
}
\ 
\gf (x)
\td x.
\]
To confirm that $F$ is continuous, fix $t$ and let $\delta$ be a positive parameter such that 
$a = t - \delta$, $b = t + \delta$ are continuity points of $F$ $-$then 
\[
F (t + \delta) - F (t - \delta)  
\ =  \ 
\frac{\delta}{\pi} \ 
\int\limits_{-\infty}^\infty \ 
\frac{\sin \delta x}{\delta x} 
\hsx
e^{-\sqrt{-1} \hsx t \hsy x} 
\ 
\gf (x) 
\ \td x
\]
\qquad \qquad
$\implies$
\allowdisplaybreaks
\begin{align*}
\abs{F (t + \delta) - F (t - \delta)}
&\leq\ 
\frac{\delta}{\pi} \ 
\int\limits_{-\infty}^\infty \ 
\Big|
\frac{\sin \delta x}{\delta x}
\Big| 
\hsx 
\abs{\gf (x)} 
\ \td x
\\[15pt]
&\leq \ 
\frac{\delta}{\pi} \ 
\int\limits_{-\infty}^\infty \ 
\abs{\gf (x)} 
\ \td x.
\end{align*}
Now let $\delta \ra 0$, thus
\[
F (t^+) - F (t^-) 
\ = \ 
0,
\]
so $F$ is continuous at $t$.  
Next, for any $h$ (positive or negative), 
\[
\frac{F (t + h) - F (t)}{h} 
\ = \ 
\frac{1}{2 \hsy \pi} \ 
\int\limits_{-\infty}^\infty \ 
\frac
{
e^{-\sqrt{-1} \hsx t \hsy x} 
- 
e^{-\sqrt{-1} \hsx (t + h) \hsy x} 
}
{
\sqrt{-1} \hsx h \hsy x
}
\ 
\gf (x)
\ \td x
\]
\qquad \qquad
$\implies$
\allowdisplaybreaks
\begin{align*}
F^\prime (t) \ 
&=\ 
\lim\limits_{h \ra 0} \ 
\frac{F (t + h) - F (t)}{h} 
\\[15pt]
&=\ 
\frac{1}{2 \hsy \pi} \ 
\int\limits_{-\infty}^\infty \ 
\lim\limits_{h \ra 0} \ 
\frac
{
e^{-\sqrt{-1} \hsx t \hsy x} 
- 
e^{-\sqrt{-1} \hsx (t + h) \hsy x} 
}
{
\sqrt{-1} \hsx h \hsy x
}
\ 
\gf (x)
\ \td x
\\[15pt]
&=\ 
\frac{1}{2 \hsy \pi} \ 
\int\limits_{-\infty}^\infty \ 
e^{-\sqrt{-1} \hsx t \hsy x} 
\ 
\gf (x)
\ \td x.
\end{align*}

[Note: \ 
$\forall \ t$,
\[
\abs{F^\prime (t) }
\ \leq \ 
\frac{1}{2 \hsy \pi} \ 
\norm{\hsy \gf \hsy}_1 
\ < \ 
\infty.
\] 
Therefore $F$ is absolutely continuous (cf. 23.20).]
\\[-.25cm]
\end{x}

\begin{x}{\small\bf THEOREM} \ 
Suppose that $F_1$, $F_2$ are distribution functions.  
Put
$F = F_1 * F_2$ $-$then
\[
\gf
\ = \ 
\gf_1 \cdot \gf_2.
\]
[$\forall \ x$, 
\\[-1cm]
\begin{align*}
\gf (x) \ 
&=\ 
\int\limits_{-\infty}^\infty \ 
e^{\sqrt{-1} \hsx t \hsy x} 
\ \td \mu_F (t)
\\[15pt]
&=\ 
\int\limits_{-\infty}^\infty \ 
\int\limits_{-\infty}^\infty \ 
e^{\sqrt{-1} \hsx  x \hsy (t_1 + t_2)} 
\ \td \mu_{F_1} (t_1)
\ \td \mu_{F_2} (t_2)
\\[15pt]
&=\ 
\int\limits_{-\infty}^\infty \ 
e^{\sqrt{-1} \hsx x \hsy t_1} 
\ \td \mu_{F_1} (t_1)
\hsx \cdot \hsx
\int\limits_{-\infty}^\infty \ 
e^{\sqrt{-1} \hsx x \hsy t_2} 
\ \td \mu_{F_2} (t_2)
\\[15pt]
&=\ 
\gf_1 (x) \cdot \gf_2 (x).]
\end{align*}
\\[-1cm]
\end{x}

\begin{x}{\small\bf EXAMPLE} \ 
Given a distribution function $F$, consider the convolution 
\[
F * (-1) \hsy F.
\]
Then its characteristic function is
\[
\gf (x) \hsy \gf (-x) 
\ = \ 
\gf (x) \hsy \ovs{\gf (x)} 
\ = \ 
\abs{\gf (x)}^2. 
\]
\\[-1.25cm]
\end{x}

\begin{x}{\small\bf RAPPEL} \ 
$\forall \ t$, $\forall \ \sigma > 0$, 
\[
\int\limits_{-\infty}^\infty \ 
\exp \Big(-\sqrt{-1} \hsx t \hsy x
\hsx - \hsx 
\frac{\sigma^2 \hsy x^2}{2}\Big)
\td x
\ = \ 
\frac{\sqrt{2 \hsy \pi}}{\sigma} 
\hsx
\exp \Big(-\frac{t^2}{2 \sigma^2}\Big).
\]
\\[-1cm]
\end{x}

\qquad
{\small\bf \un{N.B.}} \ 
Given real variables $u$, $v$, let
\[
\phi (v) 
\ = \ 
\frac{1}{\sqrt{2 \hsy \pi}}
\hsx 
\exp \Big(-\frac{v^2}{2}\Big).
\]
Then
\[
\Phi (u) 
\ = \ 
\int\limits_{-\infty}^u \ 
\phi (v) 
\ \td v
\]
is an absolutely continuous distribution function with density $\phi (v)$ and characteristic function 
\[
\exp \Big(-\frac{x^2}{2}\Big).
\]
So, 
$\forall \ \sigma > 0$, 
$\ds\Phi_\sigma  (u) \equiv \Phi\Big(\frac{u}{\sigma} \Big)$ 
is an absolutely continuous distribution function with 
\\[.25cm]
density 
$
\ 
\ds\phi_\sigma  (v) \equiv \frac{1}{\sigma} \ \phi\Big(\frac{v}{\sigma} \Big)
\ 
$  
and characteristic function 
\[
\exp \Big(-\frac{1}{2} \hsy \sigma^2 \hsy x^2\Big).
\]
\\[-1.25cm]

\begin{x}{\small\bf LEMMA} \ 
Two distribution functions
$
\ 
\begin{cases}
\ 
F
\\
\ 
G
\end{cases}
$
that agree at all continuity points common to both agree everywhere.
\\[-.25cm]

PROOF \ 
Let
$
\begin{cases}
\ 
S
\\
\ 
T
\end{cases}
$
be the set of discontinuity points of 
$
\begin{cases}
\ 
F
\\
\ 
G
\end{cases}
$
$-$then 
$S \cup T$ 
\\

\noindent
is at most countable, hence its complement $D$ is dense.  
And on $D$, $F = G$.  
If $x_0$
is arbitrary and if $x_n \in D$ approaches $x_0$ from the right, then

\[
F (x_0) 
\ = \ 
\lim F (x_n) 
\ = \ 
\lim G (x_n) 
\ = \ 
G (x_0).
\]
\\[-1.5cm]
\end{x}

\begin{x}{\small\bf THEOREM} \ 
Suppose that $F_1$, $F_2$ are distribution functions.  
Assume:  \ 
$\gf_1 = \gf_2$ $-$then $F_1 = F_2$.
\\[-.5cm]

PROOF \ 
Write
\[
\begin{cases}
\ \ds
\gf_1 (x)
\ = \ 
\int\limits_{-\infty}^\infty \
e^{\sqrt{-1} \hsx x \hsy s}
\ \td \mu_{F_1} (s)
\\[26pt]
\ \ds
\gf_2 (x)
\ = \ 
\int\limits_{-\infty}^\infty \
e^{\sqrt{-1} \hsx x \hsy s}
\ \td \mu_{F_2} (s)
\end{cases}
.
\]
Then $\forall \ t$, $\forall \ \sigma > 0$, 
\[
\int\limits_{-\infty}^\infty \
\gf_1 (x) 
\hsx
\exp 
\Big(- \sqrt{-1} \hsx x \hsy t - \frac{\sigma^2 \hsy x^2}{2}
\Big)
\td x
\ = \ 
\int\limits_{-\infty}^\infty \
\gf_2 (x) 
\hsx
\exp
\Big(
- \sqrt{-1} \hsx x \hsy t- \frac{\sigma^2 \hsy x^2}{2}
\Big)
\td  x
\]
or still, 
\begin{align*}
\int\limits_{-\infty}^\infty \
\bigg[
\int\limits_{-\infty}^\infty \
\exp
&
\bigg(
- \sqrt{-1} \hsx x\hsy (t - s) - \frac{\sigma^2 \hsy x^2}{2}
\bigg)
\td  x
\bigg]
\ \td \mu_{F_1} (s)
\\[15pt]
&= \ 
\int\limits_{-\infty}^\infty \
\bigg[
\int\limits_{-\infty}^\infty \
\
\exp
\bigg(
- \sqrt{-1} \hsx x\hsy (t - s) - \frac{\sigma^2 \hsy x^2}{2}
\bigg)
\td  x
\bigg]
\ \td \mu_{F_2} (s)
\end{align*}
or still, 
\[
\frac{\sqrt{2 \hsy \pi}}{\sigma} \ 
\int\limits_{-\infty}^\infty \
\exp
\bigg(
- \frac{(t - s)^2}{2 \hsy \sigma^2}
\bigg)
\ \td \mu_{F_1} (s)
\ = \ 
\frac{\sqrt{2 \hsy \pi}}{\sigma} \ 
\int\limits_{-\infty}^\infty \
\exp
\bigg(
- \frac{(t - s)^2}{2 \hsy \sigma^2}
\bigg)
\ \td \mu_{F_2} (s)
\]
or still, 
\[
2 \hsy \pi
\int\limits_{-\infty}^\infty \
\frac{1}{\sigma \hsx \sqrt{2 \hsy \pi}} \ 
\exp
\bigg(
- \frac{(t - s)^2}{2 \hsy \sigma^2}
\bigg)
\ \td \mu_{F_1} (s)
\ = \ 
2 \hsy \pi
\int\limits_{-\infty}^\infty \
\frac{1}{\sigma \hsx \sqrt{2 \hsy \pi}} \ 
\exp
\bigg(
- \frac{(t - s)^2}{2 \hsy \sigma^2}
\bigg)
\ \td \mu_{F_2} (s)
\]
or still, 
\[
2 \hsy \pi
\int\limits_{-\infty}^\infty \
\phi_\sigma (t - s)
\ \td \mu_{F_1} (s)
\ = \ 
2 \hsy \pi
\int\limits_{-\infty}^\infty \
\phi_\sigma (t - s)
\ \td \mu_{F_2} (s)
\]
or still, 
\[
2 \hsx \pi \ 
\big(
\Phi_\sigma * F_1
\big)
\ = \ 
2 \hsx \pi \ 
\big(
\Phi_\sigma * F_2
\big)
\]
\qquad 
$\implies$
\[
\Phi_\sigma * F_1
\ = \ 
\Phi_\sigma * F_2
\]
\qquad 
$\implies$
\[
 F_1 * \Phi_\sigma
\ = \ 
F_2 * \Phi_\sigma
\]
\qquad 
$\implies$
\[
\int\limits_{-\infty}^\infty \ 
F_1 (t - s) 
\ \td \mu_{\Phi_\sigma} (s)
\ = \ 
\int\limits_{-\infty}^\infty \ 
F_2 (t - s) 
\ \td \mu_{\Phi_\sigma} (s)
\]
\qquad 
$\implies$
\[
\int\limits_{-\infty}^\infty \ 
F_1 (t - s) 
\hsx
\exp \Big(- \frac{s^2}{2 \sigma^2}\Big)
\ \td s
\ = \ 
\int\limits_{-\infty}^\infty \ 
F_2 (t - s) 
\hsx
\exp \Big(- \frac{s^2}{2 \sigma^2}\Big) 
\ \td s
\]
\qquad 
$\implies$
\[
\int\limits_{-\infty}^\infty \ 
F_1 (t - \sigma u) 
\hsx
\exp \Big(- \frac{u^2}{2}\Big)
\ \td u
\ = \ 
\int\limits_{-\infty}^\infty \ 
F_2 (t - \sigma u) 
\hsx
\exp \Big(- \frac{u^2}{2}\Big)
\ \td u.
\]
Now let $\sigma \ra 0$ and use dominated convergence to see that 
$F_1 (t) = F_2 (t)$ at all continuity points $t$ common to both, 
so $F_1 = F_2$ period (cf. 24.15).
\\[-.5cm]
\end{x}

\begin{x}{\small\bf REMARK} \ 
The demand is that 
$\gf_1 = \gf_2$ 
everywhere and this cannot be weakened to equality on some finite interval (cf. 24.26).
\\[-.5cm]
\end{x}

\begin{x}{\small\bf LEMMA} \ 
If $\gf_1, \gf_2, \ldots$ is a sequence of characteristic functions that converges uniformly 
on compact subsets of $\R$ to a function $f$, then $\gf \equiv f$ is a characteristic function.
\\[-.5cm]
\end{x}

\begin{x}{\small\bf EXAMPLE} \ 
Let
\[
F_n (t) 
\ = \ 
\begin{cases}
\hspace{.5cm}
0 
\hspace{2.5cm}
(t < -n)
\\[11pt]
\ \ds
\frac{n+t}{2 n}
\hspace{2.1cm}
(-n \leq t < n)
\\[11pt]
\hspace{.5cm}
1
\hspace{2.5cm}
(n \leq t)
\end{cases}
.
\]
Then $F_n$ is a distribution function whose characteristic function $\gf_n$ is given by
\[
\gf_n (x) 
\ = \ 
\frac{\sin x \hsy n}{x \hsy n}
\qquad (n = 1, 2, \ldots).
\]
Therefore
\[
\lim\limits_{n \ra \infty} \ 
\gf_n (x) 
\ = \ 
\begin{cases}
\ 
1 \quad \text{if} \quad x = 0
\\[4pt]
\ 
0 \quad \text{if} \quad x \neq 0
\end{cases}
,
\]
which shows that 24.18 can fail under the weaker assumption of mere pointwise convergence.
\\[-.5cm]
\end{x}

\begin{x}{\small\bf DEFINITION} \ 
A continuous function 
$f: \R \ra \Cx$ is said to be 
\un{positive definite} 
if for any finite sequence 
$x_1, x_2, \ldots, x_n$ 
of real numbers and for any finite sequence 
$\xi_1, \xi_2, \ldots, \xi_n$ 
of complex numbers, 
\[
\sum\limits_{k = 1}^n \ 
\sum\limits_{\ell = 1}^n \ 
f (x_k - x_\ell) 
\hsx
\xi_k
\hsx
\bar{\xi}_\ell
\ \geq \ 
0.
\]

E.g.: \ 
Every characteristic function $\gf$ is positive definite.  
Proof: \ 
\allowdisplaybreaks
\begin{align*}
\sum\limits_{k = 1}^n \ 
\sum\limits_{\ell = 1}^n \ 
\gf (x_k - x_\ell) 
\hsx
\xi_k
\hsx
\bar{\xi}_\ell\ 
&=\ 
\sum\limits_{k = 1}^n \ 
\sum\limits_{\ell = 1}^n \ 
\bigg(
\int\limits_{-\infty}^\infty \ 
e^{\sqrt{-1} \hsx (x_k - x_\ell) \hsy t}
\ \td \mu_F (t)
\bigg)
\xi_k
\hsx
\bar{\xi}_\ell\ 
\\[15pt]
&=\ 
\int\limits_{-\infty}^\infty \ 
\sum\limits_{k = 1}^n \ 
\sum\limits_{\ell = 1}^n \ 
e^{\sqrt{-1} \hsx (x_k - x_\ell) \hsy t}
\hsx
\xi_k
\hsx
\bar{\xi}_\ell\
\ \td \mu_F (t)
\\[15pt]
&=\ 
\int\limits_{-\infty}^\infty \ 
\bigg(
\sum\limits_{k = 1}^n \ 
e^{\sqrt{-1} \hsx  x_k  \hsy t}
\hsx
\xi_k
\bigg)
\ 
\bigg(
\sum\limits_{\ell = 1}^n \ 
e^{-\sqrt{-1} \hsx x_\ell \hsy t}
\hsx
\bar{\xi}_\ell\
\bigg)
\ \td \mu_F (t)
\\[15pt]
&=\ 
\int\limits_{-\infty}^\infty \ 
\bigg|
\sum\limits_{k = 1}^n \ 
e^{\sqrt{-1} \hsx  x_k  \hsy t}
\hsx
\xi_k
\bigg|^2
\ \td \mu_F (t)
\\[15pt]
&\geq\ 
0.
\end{align*}
\\[-1.5cm]
\end{x}

Conversely: 
\\[-.25cm]

\begin{x}{\small\bf THEOREM} \ 
A positive definite function 
$f : \R \ra \Cx$ 
such that $f (0) = 1$ 
is a characteristic function.  
\\[-.25cm]
\end{x}

We shall preface the proof with a lemma.
\\[-.25cm]

\begin{x}{\small\bf LEMMA} \ 
Suppose that 
$\phi \in \Lp^1 [-A, A]$.  
Assume: \ 
$\phi$ is bounded, say 
$\sup \abs{\phi} \leq M$, 
and
\[
\Phi (x) 
\ = \ 
\int\limits_{-A}^A \ 
e^{\sqrt{-1} \hsx x \hsy t} 
\phi (t) 
\ \td t
\ \geq \ 
0.
\]
Then 
$\Phi \in \Lp^1 [-\infty, \infty]$.
\\[-.5cm]

PROOF \ 
Put
\[
G(X) 
\ = \ 
\int\limits_{-X}^X \ 
\Phi.
\]
Then $G$ is increasing, thus it need only be shown that $G$ is bounded.  
To this end, introduce
\[
F (X) 
\ = \ 
\frac{1}{X} \ 
\int\limits_X^{2 \hsy X} \ 
G.
\]
Then
\[
F (X) 
\ \geq \ 
\frac{G(X)}{X} \ 
\int\limits_X^{2 X} \ 
1
\ = \ 
G(X),
\]
so it will be enough to prove that $F$ is bounded.
\allowdisplaybreaks
\begin{align*}
\text{\textbullet} \quad
G(X) \ 
&=\ 
\int\limits_{-X}^X \ 
\Phi
\hspace{8.95cm}
\\[15pt]
&=\ 
\int\limits_{-X}^X \ 
\Big(
\int\limits_{-A}^A \ 
e^{\sqrt{-1} \hsx x \hsy t} 
\phi (t) 
\ \td t
\Big)
\ \td x
\\[15pt]
&=\ 
\int\limits_{-A}^A \ 
\Big(
\int\limits_{-X}^X \ 
e^{\sqrt{-1} \hsx x \hsy t} 
\ \td x
\Big)
\hsx 
\phi (t) 
\ \td t
\\[15pt]
&=\ 
\int\limits_{-A}^A \ 
\bigg(
\frac{e^{\sqrt{-1} \hsx x \hsy t} }{\sqrt{-1} \hsx t} 
\ 
\bigg|_{x = -X}^{x = X}
\bigg)
\hsx 
\phi (t) 
\ \td t
\\[15pt]
&=\ 
\int\limits_{-A}^A \ 
\frac{e^{\sqrt{-1} \hsx X \hsy t} \hsx - \hsx e^{-\sqrt{-1} \hsx X \hsy t}}{\sqrt{-1} \hsx t} 
\hsx 
\phi (t) 
\ \td t
\\[15pt]
&=\ 
2 \ 
\int\limits_{-A}^A \ 
\frac{\sin X  t}{t} 
\ 
\phi (t) 
\ \td t.
\end{align*}
\allowdisplaybreaks
\begin{align*}
\text{\textbullet} \quad
F(X) \ 
&=\ 
\frac{1}{X} \ 
\int\limits_X^{2 X} \ 
G
\hspace{6cm}
\\[15pt]
&=\ 
\frac{2}{X} \ 
\int\limits_X^{2 X} \ 
\Big(
\int\limits_{-A}^A \ 
\frac{\sin Y t}{t}
\hsx
\phi (t) 
\ \td t
\Big)
\ \td Y
\\[15pt]
&=\ 
\frac{2}{X} \ 
\int\limits_{-A}^A \ 
\Big(
\int\limits_X^{2 X} \ 
\frac{\sin Y t}{t}
\ \td Y
\Big)
\phi (t) 
\ \td t
\\[15pt]
&=\ 
\frac{2}{X} \ 
\int\limits_{-A}^A \ 
\bigg(
\frac{- \cos Y \hsy t }{t^2} 
\ 
\bigg|_{Y = X}^{Y = 2 X}
\bigg)
\hsx 
\phi (t) 
\ \td t
\\[15pt]
&=\ 
\frac{2}{X} \ 
\int\limits_{-A}^A \ 
\frac{\cos X \hsy t -  \cos 2 \hsy X \hsy t  }{t^2} 
\ 
\phi (t) 
\ \td t
\\[15pt]
&=\ 
\frac{2}{X} \ 
\int\limits_{-A}^A \ 
\frac{1 - 2 \hsy\ds\sin^2 \frac{X \hsy t}{2} \hsx - \hsx (1 - 2 \hsy \sin^2 X \hsy t)}{t^2} 
\ 
\phi (t) 
\ \td t
\\[15pt]
&=\ 
\frac{4}{X} \ 
\int\limits_{-A}^A \ 
\frac{\sin^2 X \hsy t}{t^2} 
\ 
\phi (t) 
\ \td t
\ - \ 
\frac{4}{X} \ 
\int\limits_{-A}^A \ 
\frac{\ds\sin^2 \frac{X \hsy t}{2}}{t^2} 
\ 
\phi (t) 
\ \td t.
\end{align*}
To bound the first term, write
\begin{align*}
\abs
{
\frac{4}{X} \ 
\int\limits_{-A}^A \ 
\frac{\sin^2 X \hsy t}{t^2} 
\ 
\phi (t) 
\ \td t
} \ 
&\leq \ 
\frac{4 \hsy M}{X} \ 
\int\limits_{-A}^A \ 
\frac{\sin^2 X \hsy t}{t^2} 
\ \td t
\\[15pt]
&\leq \ 
4 \hsy M \ 
\int\limits_{-\infty}^\infty \ 
\frac{\sin^2 t}{t^2} 
\ \td t
\\[15pt]
&< \ 
\infty.
\end{align*}
Ditto for the second term.
\\[-.25cm]
\end{x}

Passing to the proof of 24.21, let
\[
f_A (x) 
\ = \ 
\frac{1}{\sqrt{2 \hsy \pi} \hsx A} \ 
\int\limits_0^A \ 
\int\limits_0^A \ 
f (u - v) 
\hsx
e^{\sqrt{-1} \hsx x \hsy u} 
\hsx
e^{-\sqrt{-1} \hsx x \hsy v} 
\ \td u \hsx \td v
\qquad (A > 0).
\]
The fact that $f$ is positive definite then implies by approximation that 
$f_A (x) \geq 0$.  
Now make the change of variable 
$u = u$, $v = u - t$ to get
\[
f_A (x) 
\ = \ 
\frac{1}{\sqrt{2 \hsy \pi}} \ 
\int\limits_{-A}^A \ 
e^{\sqrt{-1} \hsx x \hsy t} 
\hsx
\bigg(
1 - \frac{\abs{t}}{A}
\bigg)
\hsy 
f (t) 
\ \td t.
\]
This done, in 24.22 take

\[
\phi( t) 
\ = \ 
\bigg(
1 - \frac{\abs{t}}{A}
\bigg)
f (t),
\]
the conclusion being that 
$f_A \in \Lp^1 [-\infty, \infty]$.  
But then 21.17 is applicable, so 
\[
\bigg(
1 - \frac{\abs{t}}{A}
\bigg)
f (t) 
\ = \ 
\frac{1}{\sqrt{2 \hsy \pi}} \ 
\int\limits_{-\infty}^\infty \ 
f_A (x)  
\hsx 
e^{-\sqrt{-1} \hsx t \hsy x} 
\ \td x, 
\]
i.e., 
\[
\bigg(
1 - \frac{\abs{t}}{A}
\bigg)
f (t) 
\ = \ 
\frac{1}{\sqrt{2 \hsy \pi}} \ 
\int\limits_{-\infty}^\infty \ 
f_A (-x)  
\hsx 
e^{\sqrt{-1} \hsx t \hsy x} 
\ \td x
\]
if $\abs{t} \leq A$.  
In particular: 
\[
1 
\ = \ 
f (0) 
\ = \ 
\frac{1}{\sqrt{2 \hsy \pi}} \ 
\int\limits_{-\infty}^\infty \ 
f_A (-x)  
\ \td x.
\]
Therefore
\[
F_A (x)
\ = \ 
\frac{1}{\sqrt{2 \hsy \pi}} \ 
\int\limits_{-\infty}^x \ 
f_A (-y)  
\ \td y
\]
is a distribution function whose characteristic function is
\[
\chisubBminusAAB (t)
\hsx
\bigg(
1 - \frac{\abs{t}}{A}
\bigg)
f (t).
\]
Finally, put
\[
\gf_n (t) 
\ = \ 
\chisubBminusnnB (t)
\bigg(
1 - \frac{\abs{t}}{n}
\bigg)
f (t) 
\qquad (n = 1, 2, \ldots).
\]
Then $\gf_n \ra f$ uniformly on compact subsets of $\R$, thus, as the $\gf_n$ are characteristic functions, 
the same is true of $\gf \equiv f$ (cf. 24.18).
\\[-.5cm]

\begin{x}{\small\bf EXAMPLE} \ 
If $\gf$ is a characteristic function, then $e^{\gf - 1}$ is a characteristic function.
\\[-.5cm]
\end{x}

\begin{x}{\small\bf POLYA CRITERION} \ 
Suppose that $f:\R \ra \R$ is continuous.  
\\
Assume: \ 
$f(0) = 1$, $f (-x) = f(x)$, 
\[
f\Big(\frac{x_1 + x_2}{2}\Big)
\ \leq \ 
\frac{f(x_1) + f(x_2)}{2}
\qquad (x_1, x_2 > 0),
\]
and 
$\lim\limits_{x \ra \infty} \ f(x) = 0$ $-$then $f$ is the characteristic function of an 
absolutely continuous distribution function $F$.
\\[-.5cm]

PROOF \ 
Because $f$ is a continuous, convex function, its derivative $D_+ f$ from the right exists for $x > 0$.  
As such, it is increasing and here
\[
D_+ f(x) 
\ \leq \ 
0 
\ (x > 0), 
\quad 
\lim\limits_{x \ra \infty} \ D_+ f(x) \ = \ 0.
\]
In addition, 
\[
f (x) 
\ = \ 
f(0) 
\hsx + \hsx 
\int\limits_0^x \ 
D_+ f(y)
\ \td y
\]
\qquad \qquad
$\implies$
\[
0 
\ = \ 
f (\infty)
\ = \ 
f(0) 
\hsx + \hsx 
\lim\limits_{x \ra \infty} \ 
\int\limits_0^x \ 
D_+ f(y)
\ \td y
\]
\qquad \qquad
$\implies$
\[
1
\ = \ 
f (0)
\ = \ 
-
\lim\limits_{x \ra \infty} \ 
\int\limits_0^x \ 
D_+ f(y)
\ \td y.
\]
Therefore
$D_+ f$
is integrable on 0 to $\infty$.  
Put
\[
\phi_X (t) 
\ = \ 
\frac{1}{2 \hsy \pi}
\int\limits_{-X}^X \ 
f (x) 
\hsx
e^{-\sqrt{-1} \hsx t \hsy x}
\ \td x.
\]
Then
\allowdisplaybreaks
\begin{align*}
\phi_X (t) \ 
&=\ 
\frac{1}{\pi}
\int\limits_0^X \ 
f (x) \hsx \cos t x
\ \td x
\\[15pt]
&=\ 
\Big(
\frac{\sin X \hsy t}{\pi \hsy t} 
\Big)
\hsx 
f (X)
\hsx -\hsx
\frac{1}{\pi \hsy t}
\int\limits_0^X \ 
D_+ f (x) \hsx \sin t x
\ \td x.
\end{align*}
So for $t \neq 0$, 
\allowdisplaybreaks
\begin{align*}
\phi (t) \ 
&\equiv \ 
\lim\limits_{X \ra \infty} \ 
\phi_X (t) 
\\[15pt]
&=\ 
-
\frac{1}{\pi \hsy t} \ 
\int\limits_0^\infty \ 
D_+ f(x) 
\hsy 
\sin t \hsy x 
\ \td x
\\[15pt]
&=\ 
-
\frac{1}{\pi \hsy t} \ 
\sum\limits_{k = 0}^\infty \ 
\int\limits_{k \hsy \pi / t}^{(k+1) \pi / t} \ 
D_+ f(x) 
\hsy 
\sin t \hsy x 
\ \td x
\\[15pt]
&=\ 
-
\frac{1}{\pi \hsy t} \ 
\sum\limits_{k = 0}^\infty \ 
\int\limits_0^{\pi / t} \ 
(-1)^k 
\hsy 
D_+ f(x + (k \hsy \pi /t)) 
\hsy 
\sin t \hsy x 
\ \td x.
\end{align*}
Since
\[
\sum\limits_{k = 0}^\infty \ 
(-1)^k 
\hsy
D_+ f(x + (k \hsy \pi /t)) 
\]
is an alternating series whose terms are decreasing in absolute value with 
\[
\lim\limits_{k \ra \infty} \ 
D_+ f(x + (k \hsy \pi /t)) 
\ = \ 
0,
\]
it is boundedly convergent and since the first term is
\[
D_+ f (x) 
\ \leq \ 
0,
\]
it follows that
\allowdisplaybreaks
\begin{align*}
\phi (t) \ 
&=\ 
-
\frac{1}{\pi \hsy t} \ 
\int\limits_0^{\pi / t} \ 
\bigg(
\sum\limits_{k = 0}^\infty \ 
(-1)^k 
\hsy 
D_+ f(x + (k \hsy \pi /t)
\bigg) 
\hsx 
\sin t \hsy x 
\ \td x
\\[15pt]
&\geq \ 
0.
\end{align*}

\noindent
Now multiply $\phi(t)$ by $\cos xt$ and integrate with respect to $t$ from 0 to $T$: 
\[
\int\limits_0^T \ 
\phi (t) 
\hsx 
\cos x \hsy t 
\ \td t \ 
\ = \ 
-
\frac{1}{\pi} \ 
\int\limits_0^\infty \ 
D_+ f(y)
\ \td y
\ 
\int\limits_0^T \ 
\frac{\cos x \hsy t \hsx \sin y \hsy t}{t}
\ \td t.
\]
Next, let $T \ra \infty$: 
\[
\lim\limits_{T \ra \infty} \ 
\int\limits_0^T \ 
\frac{\cos x \hsy t \hsx \sin y \hsy t}{t}
\ \td t 
\quad = \quad
\begin{cases}
\ \ds
0 \hspace{0.5cm} (\abs{x} > y)
\\[8pt]
\ \ds
\frac{\pi}{4} \hspace{0.5cm} (\abs{x} = y)
\\[11pt]
\ \ds
\frac{\pi}{2} \hspace{0.5cm} (\abs{x} < y)
\end{cases}
\]
$\implies$
\allowdisplaybreaks
\begin{align*}
\lim\limits_{T \ra \infty} \ 
\int\limits_0^T \ 
\phi (t) 
\hsx 
\cos x \hsy t 
\ \td t \ 
&=\ 
-\frac{1}{2} \ 
\int\limits_x^\infty \ 
D_+ f(y)
\ \td y \ 
\\[15pt]
&=\ 
-\frac{1}{2} \ 
\Big(
\int\limits_0^\infty \ 
D_+ f(y)
\ \td y \ 
\hsx - \hsx
\int\limits_0^x \ 
D_+ f(y)
\ \td y \ 
\Big)
\\[15pt]
&=\ 
-\frac{1}{2} \ 
(1 - (f (x) - 1))
\\[15pt]
&=\ 
\frac{1}{2} \ 
f (x).
\end{align*}
In particular: \ 
\[
\lim\limits_{T \ra \infty} \ 
\int\limits_0^T \ 
\phi (t) 
\ \td t \ 
\ = \ 
\frac{1}{2} \hsy f (0)
\ = \ 
\frac{1}{2},
\]
so, being nonnegative, $\phi$ is integrable on 0 to $\infty$, 
or still, being even, $\phi$ is integrable on $-\infty$ to $\infty$.  
And
\[
f (x) 
\ = \ 
\int\limits_{-\infty}^\infty \ 
\phi (t) 
\hsx
e^{\sqrt{-1} \hsx x \hsy t} 
\ \td t,
\]
thus to finish, let
\[
F (x) 
\ = \ 
\int\limits_{-\infty}^x \ 
\phi (t) 
\hsx
\ \td t.
\]
\\[-1.cm]
\end{x}

\begin{x}{\small\bf EXAMPLE} \ 
The function $e^{-\abs{x}}$ satisfies the assumptions of 24.24 but the function $e^{-\abs{x}^2}$ 
does not satisfy the assumptions of 24.24 (even though it is a characteristic function).
\\[-.25cm]
\end{x}

\begin{x}{\small\bf EXAMPLE} \ 
The functions
\[
\begin{cases}
\ \ds
1 - \abs{x}
\hspace{0.5cm} \Big(0 \leq x \leq \frac{1}{2}\Big)
\\[11pt]
\ \ds
\frac{1}{ 4 \hsy \abs{x}}
\hspace{.95cm} \Big(\abs{x} \geq \frac{1}{2}\Big)
\end{cases}
,\quad
\begin{cases}
\ \ds
1 - \abs{x}
\hspace{0.5cm} (\abs{x} \leq 1)
\\[11pt]
\ \ds
\ \ 
0
\hspace{1.3cm} (\abs{x} \geq 1)
\end{cases}
\]
satisfy the assumptions of 24.24.
\\[-.5cm]

[Note: \ 
This shows that distinct characteristic functions can coincide on a finite interval.]
\\[-.25cm]
\end{x}


\chapter{
$\boldsymbol{\S}$\textbf{25}.\quad  
\normalsize HOLOMORPHIC CHARACTERISTIC FUNCTIONS}
\setlength\parindent{2em}
\setcounter{theoremn}{0}
\renewcommand{\thepage}{\S25-\arabic{page}}

\qquad 
Let $F : \R \ra \R$ be a distribution function.
\\

\begin{x}{\small\bf DEFINITION} \ 
Let 
$k = 0, 1, 2, \ldots \hsx .$
\\[-.25cm]

\qquad \textbullet \quad
$
\ds
\alpha_k 
\ = \ 
\int\limits_{-\infty}^\infty \ 
t^k 
\ \td \mu_F (t)
$
\\

\noindent
is the \un{moment} of order $k$ of $F$.
\\[-.25cm]

\qquad \textbullet \quad
$
\ds
\beta_k 
\ = \ 
\int\limits_{-\infty}^\infty \ 
\abs{t}^k 
\ \td \mu_F (t)
$
\\

\noindent
is the \un{absolute moment} of order $k$ of $F$.
\\[-.5cm]

[Note: \ 
$\alpha_k$ exists iff $\beta_k$ exists.]
\\[-.25cm]
\end{x}

\begin{x}{\small\bf INEQUALITIES} \ 

\[
\begin{cases}
\ 
\alpha_{2 \hsy k} 
\hspace{0.5cm} 
= \ \beta_{2 \hsy k}
\qquad (\alpha_0 = \beta_0 = 1)
\\[4pt]
\ 
\alpha_{2 \hsy k - 1} \ \leq \  \abs{\alpha_{2 \hsy k - 1}} \ \leq \ \beta_{2 \hsy k - 1}
\\[4pt]
\ 
\beta_{k - 1}^2 
\hspace{0.35cm}  
\leq \ \beta_{k - 2}  \hsy \beta_k
\\[4pt]
\ 
\beta_1 
\hspace{0.75cm} 
\leq \ \beta_2^{1/2} 
\leq \cdots \leq 
\beta_2^{1/k} 
\end{cases}
.
\]
\\[-.75cm]
\end{x}

\begin{x}{\small\bf LEMMA} \ 
If $\gf$ has a derivative of order $n$ at $x = 0$, 
then all the moments of $F$ up to order $n$ or up to order $n -1$ exist according to whether $n$ is even or odd. 
\\[-.25cm]
\end{x}

\begin{x}{\small\bf EXAMPLE} \ 
Take $n = 1$ (odd) $-$then it can happen that $\gf^\prime (x)$ exists and is continuous for all values of $x$, 
yet the first moment of $F$ does not exist.  
\\[-.5cm]

[Put
\[
C 
\ = \ 
\sum\limits_{j= 2}^\infty \ 
\frac{1}{j^2 \log j}.
\]
Then
\[
F (t) 
\ = \ 
C^{-1} \ 
\sum\limits_{j= 2}^\infty \ 
\frac{1}{2 \hsy j^2 \log j}
\ 
\big[
I (t - j) + I (t + j)
\big]
\]
is a distribution function whose characteristic function is 
\[
\gf (x) 
\ = \ 
C^{-1} \ 
\sum\limits_{j= 2}^\infty \ 
\frac{\cos j x}{j^2 \log j}.
\] 
To see the claim per $\gf^\prime (x)$, note that
\[
C^{-1} \ 
\sum\limits_{j= 2}^\infty \ 
\frac{\cos j x}{\log j}
\]
is the Fourier series of an integrable function, hence on general grounds, the series
\[
C^{-1} \ 
\sum\limits_{j= 2}^\infty \ 
\frac{- \sin j x}{j \log j}
\]
is uniformly convergent (or proceed directly via the uniform Dirichlet test).  
On the other hand, 
\[
\int\limits_{-\infty}^\infty \ 
\abs{t} 
\ \td \mu_F (t) 
\ = \ 
C^{-1} \ 
\sum\limits_{j= 2}^\infty \ 
\frac{1}{j \log j}
\ = \ 
\infty.]
\]
\\[-1.25cm]
\end{x}

\begin{x}{\small\bf REMARK} \ 
A characteristic function may be nowhere differentiable.
\\[-.5cm]

[The function 
\[
\gf (x) 
\ = \ 
\sum\limits_{j= 0}^\infty \ 
\frac{1}{2^{j+1}} 
\
e^{\sqrt{-1} \hsx x \hsy 5^{^j}}
\]
is the characteristic function of 
\[
F (t) 
\ = \ 
\sum\limits_{j= 0}^\infty \ 
\frac{1}{2^{j+1}} 
\
I (t - 5^j).]
\]
\\[-1.25cm]
\end{x}

\begin{x}{\small\bf LEMMA} \ 
If the moment $\alpha_k$ of order $k$ of $F$ exists, then $\gf$ is $k$-times differentiable and 
\[
\gf^{(k)} (x)
\ = \ 
\big(\sqrt{-1}\big)^k \ 
\int\limits_{-\infty}^\infty \ 
t^k 
\hsx 
e^{\sqrt{-1} \hsx x \hsy t}
\ \td \mu_F (t)
\]
is a continuous function of $x$.  
\\[-.5cm]

[Note: \ 
In particular, 
\[
\gf^{(k)} (0) 
\ = \ 
\big(\sqrt{-1}\big)^k 
\hsx 
\alpha_k.]
\]
\\[-1.25cm]
\end{x}

\begin{x}{\small\bf SCHOLIUM} \ 
The existence of the derivatives of all orders at the origin for $\gf$ is equivalent to the existence of the moments 
of all orders for $F$. 
\\[-.25cm]
\end{x}

\begin{x}{\small\bf DEFINITION} \ 
A characteristic function $\gf$ is said to be a 
\un{holomorphic} \un{characteristic function} 
if for some $\delta > 0$ it coincides with a  function $\fg$ which is holomorphic in the disk 
$\abs{z} < \delta$.
\\[-.25cm]
\end{x}

\begin{x}{\small\bf THEOREM} \ 
If $\gf$ is a holomorphic characteristic function, then $\gf$ is holomorphic in the strip containing the origin of the form 
$-\alpha < \Img z < \beta$ $(\alpha > 0, \hsx \beta > 0$ 
(either $\alpha$ or $\beta$ or both might be $\infty$)) and in that strip, 
\[
\gf (z)
\ = \ 
\int\limits_{-\infty}^\infty \ 
\hsx 
e^{\sqrt{-1} \hsx z \hsy t}
\ \td \mu_F (t).
\]

PROOF  \ 
It is clear that $\gf$ has derivatives of all orders at the origin 
($\forall \ n$, $\gf^{(n)} (0) = \fg^{(n)} (0)$), 
hence $F$ has moments of all orders (cf. 25.7).  
Moreover, 
\[
\abs{\gf^{(2 k)} (0)}
\ = \ 
\alpha_{2 k} 
\ = \ 
\beta_{2 k}, 
\quad
\abs{\gf^{(2 k -1)} (0)}
\ = \ 
\abs{\alpha_{2 k}-1}.
\]
Thus the series
\[
\sum\limits_{k = 0}^\infty \ 
\frac{\abs{\alpha_k}}{k !}
\hsx 
r^k
\]
is convergent if $0 \leq r < \delta$, thus the series
\[
\sum\limits_{k = 0}^\infty \ 
\frac{\beta_{2 k}}{(2 k) !}
\hsx 
r^{2 k}
\]
is convergent if $0 \leq r < \delta$.  
It is also true that the series

\[
\sum\limits_{k = 1}^\infty \ 
\frac{\beta_{2 k - 1}}{(2 k - 1) !}
\hsx 
r^{2 k - 1}
\]
is convergent if $0 \leq r < \delta$.
In fact, its radius of convergence $R$ is
\[
\underset{k \ra \infty}{\liminfx} \ 
\bigg[
\frac{\beta_{2 k -1}}{(2 k -1) !}
\bigg]^{- 1/(2k -1)}.
\]
But
\[
(\beta_{2 k -1})^{1/ (2k -1)} 
\ \leq \ 
(\beta_{2 k})^{1/2k}
\qquad (\tcf. \ 25.2).
\]
So
\allowdisplaybreaks
\begin{align*}
R \ 
&\geq\ 
\underset{k \ra \infty}{\liminfx} \ 
(\beta_{2 \hsy k})^{-1/2 \hsy k}
\
\big[
(2 \hsy k -1) !
\big]^{1/(2 \hsy k - 1)}
\\[15pt]
&=\ 
\underset{k \ra \infty}{\liminfx} \ \ 
(\beta_{2 \hsy k})^{-1/2 \hsy k}
\
\big[
(2 \hsy k) !
\big]^{1/(2 \hsy k - 1)}
\qquad 
\Big(
\lim\limits_{k \ra \infty} \ 
( 2 \hsy k)^{1/(2 \hsy k - 1)}
\ = \ 1
\Big)
\\[15pt]
&\geq\ 
\underset{k \ra \infty}{\liminfx} \ 
\bigg[
\frac{\beta_{2 \hsy k}}{(2 \hsy k) !}
\bigg]^{- 1/2 \hsy k}.
\end{align*}
Applying now the monotone convergence theorem, we have
\allowdisplaybreaks
\begin{align*}
\int\limits_{-\infty}^\infty \ 
e^{r \hsy \abs{t}}
\ \td \mu_F (t)
&=\ 
\int\limits_{-\infty}^\infty \ 
\sum\limits_{n = 0}^\infty \ 
\frac{r^n \hsy \abs{t}^n}{n !}
\ \td \mu_F (t)
\\[15pt]
&=\ 
\sum\limits_{n = 0}^\infty \ 
\bigg(
\int\limits_{-\infty}^\infty \ 
\abs{t}^n
\ \td \mu_F (t)
\bigg)
\
\frac{r^n}{n !}
\\[15pt]
&=\ 
\sum\limits_{n = 0}^\infty \ 
\frac{\beta_n}{n !} \hsy r^n
\\[15pt]
&<\ 
\infty
\qquad (0 \leq r < \delta).
\end{align*}
And this implies that 
\[
\int\limits_{-\infty}^\infty \ 
e^{r \hsy t}
\ \td \mu_F (t)
\]
exists when $-\delta < r < \delta$.  
Put 
\[
\begin{cases}
\ \ds
\alpha 
\ = \ 
\sup
\Big\{
r \geq 0 : 
\int\limits_{-\infty}^\infty \ 
e^{r \hsy t}
\ \td \mu_F (t)
\ < \ 
\infty
\Big\}
\\[26pt]
\ \ds
\beta 
\ = \ 
\sup
\Big\{
r \geq 0 : 
\int\limits_{-\infty}^\infty \ 
e^{-r \hsy t}
\ \td \mu_F (t)
\ < \ 
\infty
\Big\}
\end{cases}
\implies
\begin{cases}
\ 
\alpha 
\ \geq \ 
\delta
\\[4pt]
\ 
\beta 
\ \geq \ 
\delta
\end{cases}
.
\]
Then the integral 
\[
\int\limits_{-\infty}^\infty \ 
e^{\sqrt{-1} \hsx z \hsy t}
\ \td \mu_F (t)
\]
is defined if 
$-\alpha < \Img z < \beta$, is a holomorphic function of $z$ in this strip, 
and agrees with $\gf$ on the real axis.
\\[-.25cm]
\end{x}

\begin{x}{\small\bf RAPPEL} \ 
\ 
Suppose that the power series 
$
\ 
\ds
f (z) 
\ \equiv \ 
\sum\limits_{n = 0}^\infty \ 
a_n \hsy z^n
\ 
$ 
has a 
\\[.25cm]
positive radius of convergence $R$.  
Assume: \ 
$\forall \ n \geq 0$, $a_n \geq 0$ $-$then the point $z = R$ is a singularity for $f (z)$.
\\[-.25cm]
\end{x}

\begin{x}{\small\bf DEFINITION} \ 
Let $\gf$ be a holomorphic characteristic function and take 
$\alpha$, $\beta$ as in 25.9 $-$then the strip 
$-\alpha < \Img z < \beta$ 
is called the 
\un{strip of analyticity} of $\gf$.
\\[-.25cm]
\end{x}

\begin{x}{\small\bf ADDENDUM} \ 
$-\sqrt{-1} \hsx \alpha$
(if $\alpha$ is finite) 
and
$\sqrt{-1} \hsx \beta$
(if $\beta$ is finite) 
are
singularities for $\gf$, hence 
$-\alpha < \Img z < \beta$ 
is the largest strip in which $\gf$ is holomorphic.
\\[-.5cm]

[Put
\[
\begin{cases}
\ \ds
f_- (z) 
\ = \ 
\int\limits_{-\infty}^0 \ 
e^{z \hsy t} 
\ \td \mu_F (t)
\\[26pt]
\ \ds
f_+ (z) 
\ = \ 
\int\limits_0^\infty \ 
e^{z \hsy t}  
\ \td \mu_F (t)
\end{cases}
.
\]
Then
\[
\int\limits_{-\infty}^\infty \ 
e^{r \hsy t} 
\ \td \mu_F (t)
\ < \ 
\infty 
\qquad 
( - \beta < r < \alpha)
\]
\qquad\qquad 
$\implies$

\[
\begin{cases}
\ \ds
\int\limits_0^\infty \ 
e^{r \hsy t} 
\ \td \mu_F (t)
\ < \ 
\infty 
\qquad (r < 0)
\\[26pt]
\ \ds
\int\limits_{-\infty}^0 \ 
e^{r \hsy t} 
\ \td \mu_F (t)
\ < \ 
\infty 
\qquad (r > 0)
\end{cases}
.
\]

\noindent
Therefore
\[
\begin{cases}
\ 
\text{$f_-$ is holomorphic in $\Reg z > - \beta$}
\\[4pt]
\ 
\text{$f_+$ is holomorphic in $\Reg z < \ \alpha$}
\end{cases}
.
\]
And 
\[
\gf (- \sqrt{-1} \hsx z)
\ = \ 
f_+ (z) \hsx + \hsx  f_- (z) 
\qquad (-\beta < \Reg z < \alpha).
\]
Working now with $f_+$, we have
\[
f_+^{(n)} (0)
\ = \ 
\int\limits_0^\infty \ 
t^n
\ \td \mu_F (t)
\ \geq \ 
0.
\]

Consider the power series
\[
f_+ (z)
\ = \ 
\sum\limits_{n = 0}^\infty \ 
\frac{f_+^{(n)} (0)}{n !}
\hsx 
z^n.
\]
Its radius of convergence is $\geq \alpha$ but it cannot be $> \alpha$ since otherwise 
$\exists \ \varepsilon > 0$: 
\[
\int\limits_0^\infty \ 
e^{(\alpha + \varepsilon)  t} 
\ \td \mu_F (t)
\ = \ 
\sum\limits_{n = 0}^\infty \ 
\frac{f_+^{(n)} (0)}{n !}
\hsx 
(\alpha + \varepsilon)^n
\ < \ 
\infty,
\]
contradicting the definition of $\alpha$.  
But its coefficients are $\geq 0$, hence $z = \alpha$ is a singularity for $f_+ (z)$ (cf. 25.10).  
Since
\[
\gf (- \sqrt{-1} \hsx z)
\ = \ 
f_+ (z) \hsx + \hsx f_- (z) 
\qquad (-\beta < \Reg z < \alpha)
\]
and since $f_-$ is holomorphic in $\Reg z > -\beta$, 
it follows that $\alpha$ is a singularity for 
$\gf (- \sqrt{-1} \hsx z)$ 
or still, 
$- \sqrt{-1} \hsx \alpha$
is a singularity for $\gf (z)$.]
\\[-.5cm]

[Note: \ 
To establish that $\sqrt{-1} \hsx \beta$ is a singularity for $\gf$, 
consider the characteristic function 
$(-1) \gf$ of $(-1) F$.]
\\[-.25cm]
\end{x}

\begin{x}{\small\bf REMARK} \ 
There are characteristic functions which are not holomorphic characteristic functions, 
yet can be continued into regions other than strips.
\\[-.5cm]

[Consider 
$
\gf (x) = e^{-\abs{x}}
$ 
$-$then it can be continued into the half-planes $\Reg z > 0$ and $\Reg z \leq 0$, 
yet there is no continuation into a disk centered at the origin.]
\\[-.25cm]
\end{x}

Given a characteristic function $\gf$, put

\[
I (r) 
\ = \ 
\int\limits_{-\infty}^\infty \ 
e^{r \hsy t} 
\ \td \mu_F (t)
\qquad (-\infty < r < \infty)
\]
and let
\[
\begin{cases}
\ \ds
\un{\alpha}
\ = \ 
\underset{t \ra \infty}{\liminfx} \ 
-
\frac{\log (1 - F (t))}{t}
\\[26pt]
\ \ds
\un{\beta}
\ = \ 
\underset{t \ra \infty}{\liminfx} \ 
-
\frac{\log F (-t)}{t}
\end{cases}
.
\]
\\[-.25cm] 

\qquad
{\small\bf \un{N.B.}} \ 
Equivalently, 
\[
\begin{cases}
\ \ds
\un{\alpha}
\ = \ 
-
\underset{t \ra \infty}{\limsupx} \ 
\frac{\log (1 - F (t))}{t}
\\[26pt]
\ \ds
\un{\beta}
\ = \ 
-
\underset{t \ra \infty}{\limsupx} \ 
\frac{\log F (-t)}{t}
\end{cases}
.
\]
\\[-.75cm]

\begin{x}{\small\bf LEMMA} \ 
$I (r)$ is defined for all points 
$r \in \ ]-\un{\beta}, \un{\alpha} \hsy[$, 
where it is understood that 
$\un{\beta}$ 
(respectively
$\un{\alpha}$)
is to be taken as infinite if $F (-t) = 0$ (respectively $1 - F(t) = 0$) for some $t > 0$.  
\\[-.5cm]

PROOF \ 
Noting that 
$\un{\alpha} \geq 0$, 
$\un{\beta} \geq 0$, 
consider the interval 
$[0, \un{\alpha}[$.  
Since $I(0) = 1$, take 
$\un{\alpha} > 0$ and 
$0 < r < \un{\alpha}$.  
Choose 
$r_0 : r < r_0 < \un{\alpha}$ 
and then choose 
$T = T (r_0) > 0$: 
\[
t 
\ \geq \ 
T 
\implies 
-
\frac{\log (1 - F (t))}{t}
\ \geq \ 
r_0
\]
or still, 

\[
t 
\ \geq \ 
T 
\implies 
1 - F (t) 
\ \leq \ 
e^{-t \hsy r_0}.
\]
There is no loss of generality in assuming that $T$ is a continuity point of $F$ 
$(\implies F(T^-) = F(T))$, so if $A > T$, 
\allowdisplaybreaks
\begin{align*}
\int\limits_T^A \ 
e^{r \hsy t} 
\ \td \mu_{F-1} (t)\ 
&=\ 
e^{r \hsy A} 
\hsy 
(F (A^+) - 1) - 
e^{r \hsy T} 
\hsy 
(F (T^-) - 1) 
\hsx - \hsx
r \hsx
\int\limits_T^A \ 
\big(F(t^+) - 1\big) 
\hsy
e^{r \hsy t} 
\td t
\\[15pt]
&=\ 
e^{r \hsy A} 
\hsy 
(F (A) - 1) - 
e^{r \hsy T} 
\hsy 
(F (T) - 1) 
\hsx - \hsx
r
\int\limits_T^A \ 
\big(F(t) - 1\big) 
\hsy
e^{r \hsy t} 
\ \td t
\\[15pt]
&\leq\ 
e^{r \hsy T} 
\hsy 
(1 - F(T)) 
\ + \ 
r 
\int\limits_T^A \ 
e^{r \hsy t} 
(1 - F(t)) 
\ \td t
\\[15pt]
&\leq\ 
e^{r \hsy T} 
\hsy 
(1 - F(T)) 
\ + \ 
r 
\int\limits_T^A \ 
e^{r \hsy t} 
\hsx
e^{-t \hsy r_0}
\ \td t, 
\end{align*}
hence sending $A$ to $\infty$, 
\allowdisplaybreaks
\begin{align*}
\int\limits_T^\infty \ 
e^{r \hsy t} 
\ \td \mu_F (t)\ 
&=\ 
\int\limits_T^\infty \ 
e^{r \hsy t} 
\ \td \mu_{F-1} (t)\ 
\\[15pt]
&\leq\ 
e^{r \hsy T} 
(1 - F(T))
\hsx + \hsx
r 
\hsy
\int\limits_T^\infty \ 
e^{(r - r_0) t}
\ \td t
\\[15pt]
&<\ 
\infty.
\end{align*}
Meanwhile
\[
\int\limits_{-\infty}^T \ 
e^{r \hsy t} 
\ \td \mu_F (t)\ 
\ \leq \ 
e^{r \hsy T} 
\hsx 
F (T) 
\ < \ 
\infty.
\]
Consequently, $I (r)$ is defined for all 
$r \in [0, \un{\alpha} \hsy [$.  
And, analogously, $I(r)$ is defined for all 
$r \in \ ]-\un{\beta}, 0]$.
\\[-.5cm]

[Note: \ 
$I (r)$ is defined for all $r > 0$ if $1 - F (t) = 0$ for some $t > 0$ 
and for all $r < 0$ if $F (-t) = 0$ for some $t > 0$.]
\\[-.25cm]
\end{x}

\begin{x}{\small\bf REMARK} \ 
$I (r)$ does not exist if $r > \un{\alpha}$ ($\un{\alpha}$ finite) 
or if 
$r < -\un{\beta}$ ($\un{\beta}$ finite).  
\\[-.5cm]

\noindent
E.g.: \ 
Suppose that for some $r > 0$, 
$
\ 
\ds
\int\limits_{-\infty}^\infty \ 
e^{r \hsy s} 
\ \td \mu_F (s) = C < \infty
\ 
$
$-$then $\forall \ t > 0$, 

\[
e^{r \hsy T}  (1 - F(T)) 
\ \leq \ 
\int\limits_t^\infty \ 
e^{r \hsy s} 
\ \td \mu_F (s)
\ \leq \ 
C
\]
\qquad \qquad
$\implies$

\[
\underset{t \ra \infty}{\liminfx} \ 
-
\frac{\log (1 - F (t))}{t}
\ \geq \ 
r,
\]
i.e., $r \leq \un{\alpha}$.
\\[-.5cm]

[Note: \ 
In general, 
nothing can be said about the existence of $I (r)$ when 
$r = \un{\alpha}$ 
or when 
$r = -\un{\beta}$.]
\\[-.25cm]
\end{x}

\begin{x}{\small\bf THEOREM} \ 
If 
$\un{\alpha} > 0$, 
$\un{\beta} > 0$, 
then $\gf$ is a holomorphic characteristic function.  
\\[-.5cm]

PROOF \ 
On the basis of 25.14, the integral
\[
\int\limits_{-\infty}^\infty \ 
e^{\sqrt{-1} \hsx z \hsy t} 
\ \td \mu_F (t)
\]
is defined and holomorphic in the region 
$-\un{\alpha} < \Img z < \un{\beta}$ 
and coincides with $\gf(z)$ on the real axis.
\\[-.25cm]
\end{x}

\begin{x}{\small\bf REMARK} \ 
If $\gf$ is a holomorphic characteristic function, then 
\[
\begin{cases}
\ 
\alpha
\ = \ 
\un{\alpha}
\\[4pt]
\ 
\beta
\ = \ 
\un{\beta}
\end{cases}
,
\]
where, by definition (cf. 25.9), 
\[
\begin{cases}
\ \ds
\alpha 
\ = \ 
\sup
\Big\{
r \geq 0 : 
\int\limits_{-\infty}^\infty \ 
e^{r \hsy t}
\ \td \mu_F (t)
\hsx < \hsx \infty
\Big\}
\\[26pt]
\ \ds
\beta 
\ = \ 
\sup
\Big\{
r \geq 0 : 
\int\limits_{-\infty}^\infty \ 
e^{-r \hsy t}
\ \td \mu_F (t)
\hsx < \hsx \infty
\Big\}
\end{cases}
.
\]
\\[-.25cm]
\end{x}

\begin{x}{\small\bf RAIKOV CRITERION} \ 
Suppose there exists a positive constant $R$ such that 
$\forall \ 0 < r < R$: 
\[
\begin{cases}
\ 
1 - F (t) 
\ = \ 
\tO \big(e^{-r \hsy t}\big)
\\[8pt]
\quad
F (-t) 
\ = \ 
\tO \big(e^{-r \hsy t}\big)
\end{cases}
\qquad (t \ra \infty)
.
\]
Then $f$ is a holomorphic characteristic function and its strip of analyticity (cf. 25.11) contains the strip 
$\abs{\Img z} < R$.
\\[-.5cm]

[In view of the foregoing, this is immediate.]
\\[-.25cm]
\end{x}

\begin{x}{\small\bf LEMMA} \ 
Let $\gf$ be a holomorphic characteristic function $-$then 
\[
\abs{\gf (z)} 
\ \leq \ 
\gf (\sqrt{-1} \hsx \Img z) 
\qquad (-\alpha < \Img z < \beta).
\]
[In the strip 
$-\alpha < \Img z < \beta$, 
\[
\gf (z) 
\ = \ 
\int\limits_{-\infty}^\infty \ 
e^{\sqrt{-1} \hsx z \hsy t} 
\ \td \mu_F (t).]
\]
\\[-1.25cm]
\end{x}

\begin{x}{\small\bf APPLICATION} \ 
A a holomorphic characteristic function $\gf$ has no zeros on the segment of the imaginary axis inside its strip of analyticity.
\\[-.5cm]

[For such a zero would force $\gf$ to vanish on a horizontal line within its strip of analyticity 
which in turn would imply that $\gf \equiv 0$.]
\\[-.25cm]
\end{x}

\begin{x}{\small\bf LEMMA} \ 
Let $\gf$ be a holomorphic characteristic function $-$then 
$\log \hsx \gf (\sqrt{-1} \hsx r)$ 
is convex as a function of the real variable 
$-\alpha < r < \beta$.
\\[-.5cm]

PROOF \ 
Bearing in mind that $\gf (\sqrt{-1} \hsx r) > 0$, consider the second derivative of 
$\log \gf (\sqrt{-1} \hsx r)$:
\[
\frac
{\gf (\sqrt{-1} \hsx r) \hsx \cdot \hsx \gf^{\prime\prime}  (\sqrt{-1} \hsx r)
\hsx - \hsx
\big(\gf^\prime (\sqrt{-1} \hsx r)\big)^2}
{\gf (\sqrt{-1} \hsx r)^2}.
\]
Then
\allowdisplaybreaks
\begin{align*}
\gf (\sqrt{-1} \hsx r) \hsx &\cdot \hsx \gf^{\prime\prime}  (\sqrt{-1} \hsx r)
\hsx - \hsx
\big(\gf^\prime (\sqrt{-1} \hsx r)\big)^2
\\[15pt]
&=\ 
\int\limits_{-\infty}^\infty \ 
e^{- r \hsy t} 
\ \td \mu_F (t)
\hsx \cdot  \hsx
\int\limits_{-\infty}^\infty \ 
t^2 
\hsx 
e^{- r \hsy t} 
\ \td \mu_F (t)
\ - \ 
\bigg(
\int\limits_{-\infty}^\infty \ 
t
\hsx 
e^{- r \hsy t} 
\ \td \mu_F (t)
\bigg)^2,
\end{align*}
which is nonnegative 
(Schwarz inequality applied to the measure 
$e^{- r \hsy t} \hsx \td \mu_F (t)$).
\\[-.25cm]
\end{x}

\begin{x}{\small\bf APPLICATION} \ 
For any holomorphic characteristic function $\gf$, the function 
\[
\frac{\log \gf (\sqrt{-1} \hsx r)}{r}
\]
in an increasing function of the real variable 
$0 < r < \beta$.
\\[-.5cm]

[In fact, 
$\log \hsx \gf (\sqrt{-1} \hsx r)$ is convex in $[0,\beta \hsy[$ 
and 
$\log \hsx \gf (\sqrt{-1} \ 0) = \log \gf (0) = \log 1 = 0$.]
\\[-.25cm]
\end{x}


\chapter{
$\boldsymbol{\S}$\textbf{26}.\quad  ENTIRE CHARACTERISTIC FUNCTIONS}
\setlength\parindent{2em}
\setcounter{theoremn}{0}
\renewcommand{\thepage}{\S26-\arabic{page}}

\qquad 
A holomorphic characteristic function $\gf$ is said to be \un{entire} if its strip of analyticity is the complex plane, 
i.e., if $\alpha = \infty$, $\beta = \infty$.
\\[-.5cm]

\begin{x}{\small\bf RAPPEL} \ 
\\[-.5cm]
\[
\quad
\begin{cases}
\ \ds
\un{\alpha} 
\ = \ 
\underset{t \ra \infty}{\liminfx} \ 
-\hsx
\frac{\log (1 - F (t))}{t}
\\[18pt]
\ \ds
\un{\beta} 
\ = \ 
\underset{t \ra \infty}{\liminfx} \ 
-\hsx
\frac{\log F (-t)}{t}
\end{cases}
.
\]
\\[-1.cm]
\end{x}

\begin{x}{\small\bf SCHOLIUM} \ 
A characteristic function $\gf$ is entire iff 
$\un{\alpha} = \infty$, 
$\un{\beta} = \infty$ 
(cf. 25.17).
\\[-.5cm]
\end{x}

\begin{x}{\small\bf SUBLEMMA} \ 
Suppose that $\gf$ is an entire characteristic function $-$then 
\[
M(r; \gf) 
\ = \ 
\max ( \gf (\sqrt{-1} \hsx r), \gf ( -\sqrt{-1} \hsx r)).
\]

PROOF \ 
For all real $x$ and $y$, 
\[
\abs{\gf (x + \sqrt{-1} \hsx y}
\ \leq \ 
\gf  (\sqrt{-1} \hsx y)
\qquad (\tcf. \ 25.19).
\]
\\[-1.5cm]
\end{x}

\begin{x}{\small\bf LEMMA} \ 
Suppose that $\gf$ is an entire characteristic function \text{$-$then $\forall \ t > 0$,} 
\[
M(r; \gf) 
\ \geq \
\frac{1}{2} \hsx 
e^{r \hsy t} \hsy (1 - F(t) + F (-t)).
\]

PROOF \ 
\allowdisplaybreaks
\begin{align*}
M(r; \gf) \ 
&=\ 
\max ( \gf (\sqrt{-1} \hsx r), \gf ( -\sqrt{-1} \hsx r))
\\[15pt]
&\geq\ 
(\gf (\sqrt{-1} \hsx r)
\hsx + \hsx  
\gf ( -\sqrt{-1} \hsx r)) / 2
\\[15pt]
&=\ 
\frac{1}{2} \ 
\bigg(
\int\limits_{-\infty}^\infty \ 
e^{- r \hsy s} 
\ \td \mu_F (s) 
\hsx + \hsx 
\int\limits_{-\infty}^\infty \ 
e^{r \hsy s} 
\ \td \mu_F (s) 
\bigg)
\\[15pt]
&=\ 
\int\limits_{-\infty}^\infty \ 
\cosh (r s) \ 
\ \td \mu_F (s)
\\[15pt]
&\geq\ 
\int\limits_{\abs{s} \geq t} \ 
\cosh (r s) \ 
\td \mu_F (s)
\\[15pt]
&\geq\ 
(\cosh r t) \ 
\int\limits_{\abs{s} \geq t} \ 
\td \mu_F (s) 
\\[15pt]
&\geq\ 
\frac{1}{2} \hsx
e^{r \hsy t} \ 
\int\limits_{\abs{s} \geq t} \ 
\td \mu_F (s).
\end{align*}
But
\allowdisplaybreaks
\begin{align*}
\int\limits_{\abs{s} \geq t} \ 
\td \mu_F (s) \ 
&=\ 
\mu_F ([t, \infty[ \hsy ) +\mu_F( \hsy ]-\infty, -t])
\\[15pt]
&=\ 
\mu_F ([t, \infty[ \hsy ) + F(-t).
\end{align*}
And
\[
[t, \infty[ 
\ \ = \ 
\R 
\ - \  
]-\infty, t[
\]
\qquad 
$\implies$
\allowdisplaybreaks
\begin{align*}
\mu_F ([t, \infty[ ) \ 
&=\ 
1 - \mu_F ( \hsy ]-\infty, t[ \hsy )
\\[11pt]
&\geq \ 
1 - \mu_F ( \hsy ]-\infty, t])
\\[11pt]
&= \ 
1 - F(t).
\end{align*}
\\[-1.25cm]
\end{x}

\begin{x}{\small\bf THEOREM} \ 
The order of an entire charcteristic function $\gf$ cannot be less than one except for the case when $\gf \equiv 1$ 
(i.e., when $F = I$ (cf. 23.4)).  
\\[-.5cm]

PROOF \ 
If $F \neq I$, then 

\[
1 - F (a) + F (-a) 
\ > \ 
0
\]
for some $a > 0$.  
Now take $t = a$ in 26.4.
\\[-.5cm]

[Note: \ 
It can be shown that there exist entire characteristic functions of any order $\geq 1$ (including $\infty$).]
\\[-.25cm]
\end{x}

\begin{x}{\small\bf TERMINOLOGY} \ 
Let $F$ be a distribution function.
\\[-.25cm]

\qquad \textbullet \quad
$F$ is \un{bounded to the left} if $F(a) = 0$ for some real $a$.   
When this is so, one puts

\[
\lext [F] 
\ = \ 
\sup \{a : F(a) = 0\}
\]
and calls $\lext [F]$ the \un{left extremity} of $F$.
\\[-.25cm]

\qquad \textbullet \quad
$F$ is \un{bounded to the right} if $F(b) = 1$ for some real $b$.   
When this is so, one puts

\[
\rext [F] 
\ = \ 
\inf \{b : F(b) = 1\}
\]
and calls $\rext [F]$ the \un{right extremity} of $F$.
\\[-.25cm]
\end{x}

\begin{x}{\small\bf DEFINITION} \ 
A distribution function $F$ such that $F (a) = 0$ and $F (b) = 1$ for some real $a$ and $b$ is said to be \un{finite}.
\\[-.25cm]
\end{x}

\begin{x}{\small\bf THEOREM} \ 
Let $\gf$ be an entire characteristic function.  
Assume: \ 
$\gf$ is of exponential type $-$then its distribution function $F$ is finite.  
Moreover, 
\[
\begin{cases}
\ \ds
\rext [F]
\ = \ 
\underset{r \ra \infty}{\limsupx} \ 
\frac{\log \abs{\gf (- \sqrt{-1} \hsx r)}}{r}
\\[18pt]
\ \ds
\lext [F]
\ = \ 
-\hsx
\underset{r \ra \infty}{\limsupx} \ 
\frac{\log \abs{\gf ( \sqrt{-1} \hsx r)}}{r}
\end{cases}
.
\]

PROOF  \ 
It will be enough to deal with $\lext[F]$.  
So choose $M > 0$, $K > 0$: 

\[
\abs{\gf (z)}
\ \leq \ 
M \hsy e^{K \hsy \abs{z}}.
\]
Then
\[
\log \abs{\gf ( \sqrt{-1} \hsx r)}
\ \leq \ 
\log M + K \hsy r
\]
\qquad \qquad
$\implies$
\[
\underset{r \ra \infty}{\limsupx} \ 
\frac{\log \abs{\gf ( \sqrt{-1} \hsx r)}}{r}
\ \leq \ 
K
\]
or still, 
\[
\underset{r \ra \infty}{\limsupx} \ 
\frac{\log {\gf ( \sqrt{-1} \hsx r)}}{r}
\ \leq \ 
K
\qquad (\tcf. \ 25.19)
\]
or still, 
\[
\lim\limits_{r \ra \infty} \ 
\frac{\log {\gf ( \sqrt{-1} \hsx r)}}{r}
\ \leq \ 
K
\qquad (\tcf. \ 25.22). 
\]
Denote this limit by $-a$, hence
\[
\frac{\log {\gf ( \sqrt{-1} \hsx r)}}{r}
\ \leq \ 
-a
\]
for all $r > 0$.  
Given an arbitrary $\varepsilon > 0$, let $t_1 < t_2 = a - \varepsilon$, thus
\allowdisplaybreaks
\begin{align*}
e^{- r \hsy t_2} \hsy 
(F (t_2) - F(t_1)) \ 
&=\ 
e^{- r \hsy t_2} \hsy 
\mu_F ( \hsy ]t_1, t_2] \hsy )
\hspace{4cm}
\\[15pt]
&\leq\ 
e^{- r \hsy t_2} \hsy 
\mu_F ([t_1, t_2])
\\[15pt]
&=\ 
e^{- r \hsy t_2} \hsy 
\int\limits_{t_1}^{t_2} \ 
\ \td \mu_F (t)
\\[15pt]
&=\ 
\int\limits_{t_1}^{t_2} \ 
e^{- r \hsy t_2}
\ \td \mu_F (t)
\\[15pt]
&\leq\ 
\int\limits_{t_1}^{t_2} \ 
e^{- r \hsy t}
\ \td \mu_F (t)
\\[15pt]
&\leq\ 
f (\sqrt{-1} \hsx r) 
\\[15pt]
&\leq\ 
e^{- a \hsy r}
\end{align*}
\\[-1.75cm]
\allowdisplaybreaks
\begin{align*}
\implies \hspace{0.5cm}&
\\[11pt]
&
F(t_2) - F(t_1) 
\ \leq \ 
e^{- \varepsilon\hsy r}
\\[11pt]
\implies \hspace{0.5cm}&
\\[11pt]
&
F(t_2) - F(t_1) 
\ = \ 
0
\qquad (\text{let $r \ra \infty$})
\\[11pt]
\implies \hspace{0.5cm}&
\\[11pt]
&
F(t_2) 
\ = \ 
0
\qquad (\text{let $t_1 \ra -\infty$})
\\[11pt]
\implies \hspace{0.5cm}&
\\[11pt]
&
F(a - \varepsilon) 
\ = \ 
0
\\[11pt]
\implies \hspace{0.5cm}&
\\[11pt]
&
\lext [F] 
\ \geq \ 
a.
\end{align*}
To reverse this, put
\[
\lambda_F 
\ = \ 
\lext [F] .
\]
Then 
\allowdisplaybreaks
\begin{align*}
\gf (\sqrt{-1} \hsx r) \ 
&=\ 
\int\limits_{\lambda_F}^\infty \ 
e^{- r \hsy t} 
\ \td \mu_F (t) 
\\[15pt]
&\leq \  
e^{\lambda_F \hsy r}
\end{align*}
\qquad \qquad
$\implies$

\[
a 
\ = \ 
-
\lim\limits_{r \ra \infty} \ 
\frac{\log \gf (\sqrt{-1} \hsx r)}{r} 
\ \geq \ 
\lambda_F.
\]

\noindent
Therefore
\[
a 
\ = \ 
\lambda_F
\ = \ 
\lext [F],
\]
the contention.
\\[-.25cm]
\end{x}

\qquad
{\small\bf \un{N.B.}} \ 
It is a corollary that the distribution function of an entire characteristic function of order 1 and of maximal type is not finite.
\\[-.25cm]

\begin{x}{\small\bf REMARK} \ 
Compare the above result with that of 22.10.
\\[-.25cm]
\end{x}

A \un{degenerate} distribution function is, by definition, of the form 
\[
F (t) 
\ = \ 
I (t - C),
\]
$C$ a real constant.
\\[-.25cm]

\qquad
{\small\bf \un{N.B.}} \ 
The associated characteristic function is 
\[
\gf (x) 
\ = \ 
e^{\sqrt{-1} \hsx C \hsy x}, 
\]
hence is entire of exponential type, hence further is of order 1 and type $\abs{C}$ provided $C \neq 0$.
\\[-.25cm]

\begin{x}{\small\bf LEMMA} \ 
If $F$ is degenerate, then $F$ is finite and 
\[
\rext [F]
\ = \ 
\lext [F].
\]

PROOF \ 
\[
\begin{cases}
\ \ds
\rext [F]
\ = \ 
\lim\limits_{r \ra \infty} \ 
\frac{\log e^{C \hsy r}}{r} 
\ = \ 
C
\\[18pt]
\ \ds
\lext [F]
\ = \ 
-
\lim\limits_{r \ra \infty} \ 
\frac{\log e^{-C \hsy r}}{r} 
\ = \ 
- (- C)
\ = \ 
C
\end{cases}
.
\]
\\[-.25cm]
\end{x}

\begin{x}{\small\bf CONSTRUCTION} \ 
Suppose that $F \neq I$ is a finite distribution function.  
Let
\[
\begin{cases}
\ 
a \ = \ \lext [F]
\\[4pt]
\ 
b \ = \ \rext [F]
\end{cases}
.
\]
Then
\allowdisplaybreaks
\begin{align*}
\gf (x) \ 
&=\ 
\int\limits_{-\infty}^\infty \ 
e^{\sqrt{-1} \hsx x \hsy t} 
\ \td \mu_F (t)
\\[15pt]
&=\ 
\int\limits_a^b \ 
e^{\sqrt{-1} \hsx x \hsy t} 
\ \td \mu_F (t).
\end{align*}
But the integral 
\[
\int\limits_a^b \ 
e^{\sqrt{-1} \hsx z \hsy t} 
\ \td \mu_F (t)
\]
represents an entire function, thus $f$ is an entire function of exponential type (cf. 17.19), 
thus is of order 1 (cf. 26.5).
\\[-.25cm]
\end{x}

\qquad
{\small\bf \un{N.B.}} \ 
\[
T (\gf) 
\ = \ 
\max (-a, b).
\]
For, by definition, 
\[
T (\gf) 
\ = \ 
\underset{r \ra \infty}{\limsupx} \ 
\frac{M (r; \gf)}{r}.
\]
On the other hand, 
\[
a 
\ = \ 
-
\lim\limits_{r \ra \infty} \ 
\frac{\log \gf (\sqrt{-1} \hsx r)}{r}
\]
and
\[
b
\ = \ 
\lim\limits_{r \ra \infty} \ 
\frac{\log \gf (-\sqrt{-1} \hsx r)}{r}.
\]
And
\[
M (r; \gf)
\ = \ 
\max(\gf ( \sqrt{-1} \hsx r),\gf ( -\sqrt{-1} \hsx r))
\qquad (\tcf. \ 26.3)
\]
\qquad\qquad 
$\implies$
\[
T (\gf) 
\ \geq \ 
\max (-a, b).
\]
In the other direction, 

\[
\gf ( \sqrt{-1} \hsx r)
\ \leq \ 
e^{- a \hsy r} 
\quad \text{and} \quad 
\gf ( -\sqrt{-1} \hsx r)
\ \leq \ 
e^{b \hsy r} 
\]
\qquad\qquad 
$\implies$
\[
M(r; \gf) 
\ \leq \ 
\max (e^{- a \hsy r} , e^{b \hsy r} )
\]
\qquad\qquad 
$\implies$
\[
T (\gf) 
\ \leq \ 
\max (- a, b).]
\]
\\[-,75cm]

\begin{x}{\small\bf EXAMPLE} \ 
If 
\[
F (t) 
\ = \ 
I (t - C) 
\qquad (C \neq 0), 
\]
then 
\[
a 
\ = \ 
b
\ = \ 
C.
\]

\qquad \textbullet \quad
$a > 0 \implies \max (-a, a) = a = C$
\\[-.25cm]

\qquad \textbullet \quad
$a < 0 \implies \max (-a, a) = -a = -C = \abs{C}.$
\\[-.25cm]

I.e.: \ 
$T (f) = \abs{C}$ in agreement with what has been said earlier.
\\[-.25cm]
\end{x}

\begin{x}{\small\bf REMARK} \ 
There is no entire characteristic function of order 1 and of minimal type (apply 17.18).
\\[-.25cm]
\end{x}

\begin{x}{\small\bf LEMMA} \ 
If $F$ is a finite distribution function and if $F$ is nondegenerate, then its characteristic function $\gf$ 
has an infinity of zeros (they need not be real).  
\\[-.5cm]

PROOF \ 
Since $\gf$ is bounded on the real axis, the conclusion that $\gf$ has finitely many zeros is untenable (cf. \S7).
\\[-.25cm]
\end{x}

\begin{x}{\small\bf REMARK} \ 
An infinitely divisible entire characteristic function has no 
zeros.\footnote[2]{\vspace{.11 cm}
E. Lukacs,  
\textit{Characteristic Functions},  Griffin, 1970, pp. 258-259.}
\\[-.25cm]
\end{x}


\begin{x}{\small\bf NOTATION} \ 
Given a distribution function $F$, let
\[
T (t) 
\ = \ 
1 - F(t) + F(-t) 
\qquad (t > 0).
\]
Let $K$ and $\alpha$ be positive constants.
\\[-.25cm]
\end{x}

\begin{x}{\small\bf SUBLEMMA} \ 
The integral
\[
I (z) 
\ = \ 
\int\limits_0^\infty \ 
\exp \big(\sqrt{-1} \hsx z \hsy t \hsy - \hsy K \hsy t^{1 + \alpha}\big) 
\ \td t
\]
defines an entire function of order 
$
\ds
1 + \frac{1}{\alpha}
$.
\\[-.5cm]

[Consider the expansion
\[
I (z) 
\ = \ 
\sum\limits_{n = 0}^\infty \ 
c_n \hsy z^n,
\]
where
\[
c_n 
\ = \ 
\frac{(\sqrt{-1})^n}{n !} 
\ 
\Gamma \Big(\frac{n + 1}{1 + \alpha} \Big)
\hsx
\frac{1}{(1 + \alpha) K^{(n + 1)/(1 + \alpha)}}.]
\]

[Note: \ 
To within a constant factor, $I(z)$ is an entire characteristic function.  
Accordingly, 
\begin{align*}
M(r; I) \ 
&=\ 
\max(I (\sqrt{-1} \hsx r), I(-\sqrt{-1} \hsy r)) 
\qquad (\tcf. \ 26.3)
\\[15pt]
&=\ 
\int\limits_0^\infty \ 
\exp \big(r \hsy t \hsy - \hsy K \hsy t^{1 + \alpha}\big) 
\ \td t.]
\end{align*}
\\[-.75cm]
\end{x}

\begin{x}{\small\bf LEMMA} \ 
Let $F$ be a distribution function.  
Assume:  \ 
$\exists \ A > 0$ such that
\[
t 
\ \geq \ 
A 
\ \implies \ 
T (t) 
\ \leq \ 
\exp \big(- K \hsy  t^{1 + \alpha}\big).
\]
Then the associated characteristic function $\gf$ is entire (cf. 25.18) and its order is 
$
\ds
\leq 
1 + \frac{1}{\alpha}
$.
\\[-.25cm]

PROOF \ 
Take $A > 0$ to be a continuity point of $F$ and let $R > A$ $-$then for $r > 0$: 
\allowdisplaybreaks
\begin{align*}
\int\limits_A^R \ 
e^{r \hsy t} 
\ \td \mu_F (t) 
&=\ 
\int\limits_A^R \ 
e^{r \hsy t} 
\ \td \mu_{F- 1} (t) 
\\[15pt]
&=\ 
e^{r \hsy R} 
\big(F (R^+) - 1\big)
\hsx - \hsx 
e^{r \hsy A} 
\big(F (A^-) - 1 \big)
\hsx - \hsx 
r
\int\limits_A^R \ 
(F(t^+) - 1) \hsx e^{r \hsy t} 
\ \td t
\\[15pt]
&=\ 
e^{r \hsy R} 
\big(F (R) - 1\big)
\hsx - \hsx 
e^{r \hsy A} 
\big(F (A) - 1\big)
\hsx - \hsx 
r
\int\limits_A^R \ 
(F(t) - 1) \hsx e^{r \hsy t} 
\ \td t
\\[15pt]
&\leq\ 
e^{r \hsy A} \hsy
(1 - F(A)) 
\hsx + \hsx 
r \hsx 
\int\limits_A^R \ 
e^{r \hsy t}  
\hsy 
(1 - F(t)) 
\ \td t
\end{align*}
$\implies$
\allowdisplaybreaks
\begin{align*}
\int\limits_A^\infty \ 
e^{r \hsy t} 
\ \td \mu_F (t) 
&\leq\ 
e^{r \hsy A} \hsy
(1 - F(A)) 
\hsx + \hsx 
r \hsx 
\int\limits_A^\infty \ 
e^{r \hsy t}  
\hsy 
(1 - F(t)) 
\ \td t
\\[15pt]
&\leq\ 
e^{r \hsy A} \hsy
(1 - F(A)) 
\hsx + \hsx 
r \hsx 
\int\limits_A^\infty \ 
\exp (r \hsy t - K \hsy t^{1 + \alpha})  
\ \td t
\\[15pt]
&\leq\ 
e^{r \hsy A} \hsy
(1 - F(A)) 
\hsx + \hsx 
r \hsx 
\int\limits_0^\infty \ 
\exp (r \hsy t - K \hsy t^{1 + \alpha})  
\ \td t.
\end{align*}
But
\[
\int\limits_{-\infty}^A \ 
e^{r \hsy t} \hsy
\ \td \mu_F (t)
\ \leq \ 
e^{r \hsy A} \hsy F (A).
\]
Therefore
\[
\int\limits_{-\infty}^\infty \ 
e^{r \hsy t} \hsy
\ \td \mu_F (t)
\ \leq \ 
e^{r \hsy A} \hsy
\hsx + \hsx 
r \hsx 
\int\limits_0^\infty \ 
\exp (r \hsy t - K \hsy t^{1 + \alpha})  
\ \td t.
\]
And analogously, 
\[
\int\limits_{-\infty}^\infty \ 
e^{- r \hsy t} \hsy
\ \td \mu_F (t)
\ \leq \ 
e^{r \hsy A} \hsy
\hsx + \hsx 
r \hsx 
\int\limits_0^\infty \ 
\exp (r \hsy t - K \hsy t^{1 + \alpha})  
\ \td t.
\]
These estimates then enable one to estimate $M(r; \gf)$: 
\allowdisplaybreaks
\begin{align*}
M(r; \gf) \ 
&=\ 
\max (\gf (\sqrt{-1} \hsx r), \gf (- \sqrt{-1} \hsx r)) 
\qquad (\tcf. \ 26.3)
\\[15pt]
&\leq\ 
e^{r \hsy A} \hsy
\hsx + \hsx 
r \hsx 
\int\limits_0^\infty \ 
\exp (r \hsy t - K \hsy t^{1 + \alpha})  
\ \td t
\\[15pt]
&=\ 
M(r; e^{z \hsy A}) 
\hsx + \hsx 
M(r; z \hsy I(z)).
\end{align*}
The order of $e^{z \hsy A}$ is 1 whereas the order of $I(z)$ is 
$
\ds
1 + \frac{1}{\alpha}
$ 
(cf. 26.17), 
hence the order of $z \hsy I(z)$ is also 
$
\ds
1 + \frac{1}{\alpha}
$ 
(cf. 2.36), thus for any $\varepsilon > 0$, 
\[
M(r; e^{z \hsy A}) 
\hsx + \hsx 
M(r; z \hsy I(z)) 
\ < \ 
\exp\Big(r^{1 + \frac{1}{\alpha} + \varepsilon}\Big)
\qquad (r \gg 0),
\]
which implies that the order of $\gf$ is 
$
\ds
\leq \ 
1 + \frac{1}{\alpha}
$.
\\[-.25cm]
\end{x}

\begin{x}{\small\bf THEOREM} \ 
The characteristic function $\gf$ of a distribution function $F$ is entire of order 1 and of maximal type iff
\[
t > 0 
\implies
T (t) > 0
\]
and
\[
\lim\limits_{t \ra \infty} \ 
\frac{\log \log \frac{1}{T (t)}}{\log t}  
\ = \ 
\infty.
\]

PROOF \ 
\\[-.25cm]

\qquad \textbullet \quad
\un{Necessity} \ 
It is clear that the first condition 
\[
t > 0 
\implies
T (t) > 0
\]
holds (simply note that $F$ is not finite).  
To see that the second condition holds, 
let $\varepsilon > 0$ be given and choose $R$: 
\[
r 
\ \geq \ 
R 
\implies 
\exp \big( r^{1 + \varepsilon}\big)
\ \geq \ 
M(r; \gf).
\]
But $\forall \ t > 0$, 
\[
M(r; \gf) 
\ \geq \ 
\frac{1}{2} \hsx e^{r \hsy t} \hsy T(t)
\qquad (\tcf. \ 26.4).
\]
Therefore
\[
T (t) 
\ \leq \  
2 \hsy 
\exp \big(-r \hsy t + r^{1 + \varepsilon}\big).
\]
Choose $t \geq 2 \hsy R^\varepsilon$ and taking 
$
\ds
r = \Big(\frac{t}{2}\Big)^{1 / \varepsilon}
$, 
we have
\[
T (t) 
\ \leq \  
2 \hsy 
\exp \Big(-\Big(\frac{t}{2}\Big)^{(1 + 1 / \varepsilon)}\Big)
\]
\qquad\qquad 
$\implies$
\[
\underset{t \ra \infty}{\liminfx} \ 
\frac{\log \log \frac{1}{T (t)}}{\log t} 
\ \geq \ 
1 + (1/\varepsilon)
\]
\qquad\qquad 
$\implies$
\[
\lim\limits_{t \ra \infty} \ 
\frac{\log \log \frac{1}{T (t)}}{\log t} 
\ = \ 
\infty,
\]
$\varepsilon$ being arbitrary.
\\[-.25cm]

\qquad \textbullet \quad
\un{Sufficiency} \ 
Given $\varepsilon > 0$, 
\[
\frac{\log \log \frac{1}{T (t)}}{\log t} 
\ \geq \  
1 + \frac{1}{\varepsilon}
\qquad (t \gg 0)
\]
\qquad\qquad 
$\implies$
\[
T (t) 
\ \leq \ 
\exp \Big( - t^{1 + \frac{1}{\varepsilon}}\Big)
\qquad (t \gg 0).
\]
Therefore $\gf$ is entire of order
\[
\leq \ 
1 + \frac{1}{\frac{1}{\varepsilon}} 
\ = \ 
1 + \varepsilon 
\qquad (\tcf. \ 26.18).
\]
But $F \neq I$, hence $\rho (\gf) = 1$ (cf. 26.5).  
Now $\gf$ cannot be of minimal type (cf. 26.13) nor can $\gf$ be of intermediate type (cf. 26.8 
($F$ is not finite due to the assumption on $T$)), 
thus $\gf$ must be of maximal type.
\\[-.25cm]
\end{x}

While a discussion of entire characteristic functions of order $> 1$ will be omitted, 
there is an important result of a negative nature.
\\[-.25cm]

\begin{x}{\small\bf THEOREM} \ 
If $p$ is a polynomial of degree $> 2$, then $e^p$ is not a characteristic function.
\\[-.25cm]
\end{x}
\newpage

\begin{center}
APPENDIX
\vspace{.5cm}
\end{center}

Let $F : \R \ra \R$ $-$then $F$ is an \un{NBV function} if $F$ is of bounded variation, 
if $F$ is continuous from the right, and if $F(-\infty) = 0$.
\\[-.25cm]

\qquad
{\small\bf NOTATION} \ 
$T_F$ is the total variation function associated with an NBV function $F$.  
So: \ 
\\[-.25cm]

\qquad \textbullet \quad
$T_F$ is increasing. 
\\[-.25cm]

\qquad \textbullet \quad
$T_F$ is continuous from the right. 
\\[-.25cm]

\qquad \textbullet \quad
$T_F(-\infty) = 0$, $T_F(\infty) < \infty$.
\\

\qquad
{\small\bf RAPPEL} \ 
The distribution functions $F$ are in a one-to-one correspondence with the probability measures on the line: $F \ra \mu_F$.
\\[-.25cm]

This can be generalized: \ 
The NBV functions $F$ are in a one-to-one correspondence with the finite signed measures on the line: $F \ra \mu_F$.
\\[-.25cm]


\qquad
{\small\bf NOTATION} \ 
$\abs{\mu_F}$ is the total variation measure associated with an NBV function $F$.  
So
\\[-.5cm]

\qquad \textbullet \quad
$\abs{\mu_F} (\R)  < \infty$.
\\[-.25cm]

\qquad \textbullet \quad
$\abs{\mu_F} = \mu_{T_F}$. 
\\

\qquad
{\small\bf \un{N.B.}} \ 
For the record, 
\[
F(t) 
\ = \ 
\mu_F (\hsy ]-\infty, t])
\]
and 
\[
T_F (t) 
\ = \ 
\mu_{T_F} (\hsy ]-\infty, t])
\ = \ 
\abs{\mu_F} (\hsy ]-\infty, t]).
\]
\\[-1cm]

\qquad
{\small\bf EXAMPLE} \ 
\[
\mu_{T_F} / \mu_{T_F} (\R)
\]
is a probability measure on the line.  
\\[-.25cm]

\qquad
{\small\bf LEMMA} \ 
Any bounded Borel measurable function on $\R$ is $\mu_F$-integrable (cf. 23.13).
\\

{\small\bf DEFINITION} \ 
Given an NBV function $F$, put
\[
\gf (x) 
\ = \ 
\int\limits_{-\infty}^\infty \ 
e^{\sqrt{-1} \hsx x \hsy t} 
\ \td \mu_F (t),
\]
the Fourier transform of $\mu_F$.
\\

Obviously, 
\[\abs{\gf (x)} 
\ \leq \ 
\abs{\mu_F} (\R) 
\ < \ 
\infty.
\]
\\[-1.cm]

\qquad
{\small\bf DEFINITION} \ An NBV function $F$ is \un{constant outside a finite interval} $[T^\prime, T^{\prime\prime}]$
if
\[
\begin{cases}
\ 
F (t) = 0 
\hspace{0.6cm} 
(t < T^\prime)
\\[4pt]
\ 
F (t) = C 
\hspace{0.5cm} 
(t > T^{\prime\prime})
\end{cases}
\]
for some real number $C$.
\\

\qquad
{\small\bf \un{N.B.}} \ 
Under these circumstances, 
\[
\int\limits_{-\infty}^\infty \ 
e^{\sqrt{-1} \hsx z \hsy t} 
\ \td \mu_F (t)
\ = \ 
\int\limits_{T^\prime}^{T^{\prime\prime}} \ 
e^{\sqrt{-1} \hsx z \hsy t} 
\ \td \mu_F (t)
\]
and the integral on the right is defined for all complex $z$, 
thus $\gf$ admits a continuation as an entire function and, as such, is of exponential type.
\\[-.5cm]

[Put
\[
\gT_\gf (z) 
\ = \ 
\int\limits_{-\infty}^\infty \ 
e^{\sqrt{-1} \hsx z \hsy t} 
\ \td \mu_{T_F} (t), 
\]
the ``characteristic function'' of $T_F$ $-$then 
\[
M(r; \gT_\gf) 
\ = \ 
\max (\gT_\gf (\sqrt{-1} \hsx r), \gT_\gf (-\sqrt{-1} \hsx r))
\qquad 
(\tcf. \ 26.3).
\]
On the other hand,
\allowdisplaybreaks
\begin{align*}
\abs{\gf (x + \sqrt{-1} \hsx y)} \ 
&=\ 
\bigg| 
\hsx
\int\limits_{-\infty}^\infty \ 
e^{\sqrt{-1} \hsx z \hsy t}  
\ \td \mu_F (t)
\hsx
\bigg|\ 
\\[15pt]
&\leq \ 
\int\limits_{-\infty}^\infty \ 
e^{- y \hsy t}
\ \td \mu_{T_F} (t)
\\[15pt]
&= \ 
\gT_\gf (\sqrt{-1} \hsx y)
\end{align*}
\qquad\qquad 
$\implies$
\[
M(r; \gf) 
\ \leq \ 
M(r; \gT_\gf ).
\]
But 
\[
\gT_\gf (\sqrt{-1} \hsx r) 
\ \leq \ 
e^{-T^{\prime} \hsy r} \hsx \mu_{T_F} (\R)
\]
and
\[
\gT_\gf (- \sqrt{-1} \hsx r) 
\ \leq \ 
e^{T^{\prime\prime}  \hsy r} \hsx \mu_{T_F} (\R).
\]
Therefore
\[
M (r; f) 
\ \leq \ 
\exp \big(\max \big(\abs{T^\prime}, \abs{T^{\prime\prime}} \big) \hsy r\big), 
\]
so $\gf$ is of exponential type.]
\\

\qquad
{\small\bf THEOREM} \ 
Suppose that $F$ is an NBV function.  
Assume:  \ 
$\gf$ can be extended into the complex plane as an entire function of exponential type.  
Let

\[
\begin{cases}
\ \ds
a \ = \ 
-\ 
\underset{r \ra \infty}{\limsupx} \ 
\frac{\abs{\gf (\sqrt{-1} \hsx r)}}{r}
\\[18pt]
\ \ds
b \ = \ 
\ 
\underset{r \ra \infty}{\limsupx} \ 
\frac{\abs{\gf (- \sqrt{-1} \hsx r)}}{r}
\end{cases}
.
\]

\noindent
Then $a$ and $b$ are finite (sic).  
Moreover, $F$ is constant outside a finite interval and in fact $[a, b]$ is the smallest finite interval outside of which $F$ is constant.
\\[-.25cm]

PROOF \ 
We shall work initially with $b$ to show that $F$ is constant to the right of $b$.  
To this end, note that for any pair $t_1 < t_2$ of continuity points of $F$: 

$
\ds
F (t_2) - F(t_1) 
\ = \ 
\lim\limits_{r \ra \infty} \ 
\int\limits_{-r}^r \ 
\frac
{e^{- \sqrt{-1} \hsx t_1 \hsy x} \hsx - \hsx e^{- \sqrt{-1} \hsx t_2 \hsy x}}
{2 \hsy \pi \hsx \sqrt{-1} \hsx x}
\ 
\gf (x) 
\ \td x
\qquad (\tcf. \ 24.9).
$
\\

\noindent
Now specialize and take 
$\ b \hsx  < \hsx t_1 \hsx < \hsx t_2\ $ ($t_2$ arbitrary) and let 
$\ 2 \hsy \varepsilon \hsx = \hsx  t_1 - b > 0\ $ 

\noindent
$(\implies b < b + \varepsilon = t_1 - \varepsilon < t_1)$.  
Put
\[
f (z) 
\ = \ 
\big(
1 - e^{- \sqrt{-1} \hsx (t_2 - t_1) \hsy z}
\big)
\hsx
\gf (z) 
\hsx
e^{- \sqrt{-1} \hsx (b + \varepsilon) \hsy z}.
\]
Then 
\\[-.25cm]

\qquad \textbullet \quad
$f$ is entire of exponential type.  
\\[-.25cm]

\qquad \textbullet \quad
$f$ is bounded on the real axis.
\\[-.25cm]

\qquad \textbullet \quad
$f (- \sqrt{-1} \hsx r)$ $(0 \leq r < \infty)$ is bounded.
\\[-.25cm]

\noindent
Therefore $(\ldots)$ $f$ is bounded in the lower half-plane: \ 
$\abs{f} < M$.  
And
\[
2 \hsy \pi \hsy \sqrt{-1} \hsx (F (t_2) - F(t_1)) 
\ = \ 
\lim\limits_{r \ra \infty} \ 
\int\limits_{-r}^r \ 
\frac{f(x)}{x} 
\hsx \cdot \hsx 
e^{- \sqrt{-1} \hsx \varepsilon \hsy x} 
\ \td x.
\]
Since the integrand is entire $(f (0) = 0)$, the integration interval can be replaced by a semi-circular arc of radius 
$r$ centered at the origin and situated in the lower half-plane, hence
\allowdisplaybreaks
\begin{align*}
\bigg| 
\hsx
\int\limits_{-r}^r \ 
\frac{f(x)}{x} 
\hsx \cdot \hsx 
e^{- \sqrt{-1} \hsx \varepsilon \hsy x}  
\ \td x
\hsx
\bigg|\ 
&\leq \ 
\int\limits_\pi^{2 \hsy \pi}  \ 
\abs{f \big(r \hsy e^{\sqrt{-1} \hsx \theta}\big)} 
\hsx 
e^{\varepsilon \hsy r \hsy \sin \theta} 
\ \td \theta
\\[12pt]
&\leq \ 
M \ 
\int\limits_0^\pi  \ 
e^{-\varepsilon \hsy r \hsy \sin \theta} 
\ \td \theta
\\[12pt]
&\leq \ 
2 \hsy M \ 
\int\limits_0^{\pi / 2}  \ 
e^{-\varepsilon \hsy r \hsy \sin \theta} 
\ \td \theta
\\[12pt]
&\leq \ 
2 \hsy M \ 
\int\limits_0^{\pi / 2}  \ 
e^{- (2 \hsy \varepsilon \hsy r \hsy  \theta) / \pi} 
\ \td \theta
\\[12pt]
&\ra \ 
0 
\qquad (r \ra \infty)
\end{align*}
\qquad\qquad 
$\implies$
\[
\lim\limits_{r \ra \infty} \ 
\int\limits_{-r}^r \ 
\frac{f(x)}{x} 
\hsx \cdot \hsx 
e^{- \sqrt{-1} \hsx \varepsilon \hsy x} 
\ \td x
\ = \ 
0
\]
\qquad\qquad 
$\implies$
\[
F (t_2) - F(t_1)
\ = \ 
0
\]
\qquad\qquad 
$\implies$
\[
F (t_2) 
\ = \  
F(t_1)
\ = \  
F (b + 2 \hsy \varepsilon), 
\]
proving that $F$ is constant to the right of $b$.  
By a similar argument, one finds that $F$ is constant to the left of $a$, 
thus equals $F(-\infty) = 0$ there.  
Finally, if 
$[T^\prime, T^{\prime\prime}]$ 
is a finite interval outside of which $F$ is constant, then 
$T^\prime \leq a$, 
$b \leq T^{\prime\prime}$.  
E.g.: \ 
\allowdisplaybreaks
\begin{align*}
\abs{\gf (\sqrt{-1} \hsx r)} \ 
&\leq \ 
\gT_\gf (\sqrt{-1} \hsx r) 
\\[11pt]
&\leq 
e^{-T^\prime \hsy r} \hsx \mu_{T_F} (\R)
\end{align*}
\qquad\qquad 
$\implies$
\[
a 
\ = \ 
- 
\underset{r \ra \infty}{\limsupx} \ 
\frac{\log \abs{f (\sqrt{-1} \hsx r)}}{r}
\ \geq \ 
T^\prime.
\]


\chapter{
$\boldsymbol{\S}$\textbf{27}.\quad  ZERO THEORY: \ BERNSTEIN FUNCTIONS}
\setlength\parindent{2em}
\setcounter{theoremn}{0}
\renewcommand{\thepage}{\S27-\arabic{page}}

\qquad 
Let $\sB_0 (A)$ be the subset of $\sE_0 (A)$ consisting of those $f$ which are bounded on the real axis.
\\[-.5cm]

[Note: \ 
The elements of 
$\sB_0 (A)$ 
are called \un{Bernstein functions}.]
\\[-.25cm]

\qquad
{\small\bf \un{N.B.}} \ 
If 
$f \in \sB_0(A)$ 
and if $T (f) = 0$, then $f$ is a constant (cf. 17.18).
\\[-.5cm]

[Note: \ 
Accordingly, if 
$f \in \sB_0(A)$ 
is not a constant, then $T (f) > 0$ and 
$\rho (f) = 1$ (with $T (f) = \tau (f)$) (cf. 17.3).]
\\[-.25cm]

\begin{x}{\small\bf EXAMPLE} \ 
Take $A = 1$ $-$then 
$e^{\sqrt{-1} \hsx z} \in \sB_0(1)$.
\\[-.25cm]
\end{x}

\begin{x}{\small\bf EXAMPLE} \ 
Suppose that $F \neq I$ is a finite distribution function $-$then its characteristic function 
$\gf \in \sB_0(A)$, where $A = \max (-a, b)$ (cf. 26.11).  
\\[-.5cm]

[Note: \ 
Take
\[
F (t) 
\ = \ 
I (t - 1).
\]
Then \ 
$\gf (z) \hsx = \hsx e^{\sqrt{-1} \hsx z}$.]
\\[-.25cm]
\end{x}

\begin{x}{\small\bf LEMMA} \ 
$\PW (A)$ is a subset of 
$\sB_0 (A)$ 
(cf. 17.29).
\\[-.25cm]
\end{x}

\begin{x}{\small\bf LEMMA} \ 
$\sB_0 (A)$  
is a vector space (under pointwise addition and scalar multiplication) and,  
when equipped with the supremum norm, is a Banach space (cf. 17.17).
\\[-.25cm]
\end{x}

\begin{x}{\small\bf LEMMA} \ 
$\sB_0 (A)$ 
is closed under differentiation (cf. 17.24).
\\[-.25cm]
\end{x}

\begin{x}{\small\bf LEMMA} \ 
If $f \in \sB_0 (A)$ is not a constant, then $n (r) = \tO (r)$, i.e.
$
\ds 
\frac{n (r)}{r}
$ 
remains bounded as $r \ra \infty$ (cf. 4.31).
\\[-.25cm]
\end{x}

\begin{x}{\small\bf NOTATION} \ 
Given 
$f \in \sB_0(A)$, 
let 
$z_n = r_n \hsy e^{\sqrt{-1} \hsx \theta_n}$ $(n = 1, 2, \ldots)$ be the nonzero zeros of $f$ 
repeated according to multiplicity with 
\[
0 
\ < \ 
\abs{z_1}
\ \leq \ 
\abs{z_2}
\ \leq \ 
\cdots \hsx .
\]

[Note 
\[
\frac{1}{z_n}
\ = \ 
\frac{e^{-\sqrt{-1} \hsx \theta_n}}{r_n}
\ = \ 
\frac{\cos \theta_n}{r_n}
\ - \ 
\sqrt{-1} \ 
\frac{\sin \theta_n}{r_n}.]
\]
\\[-1.cm]
\end{x}

\begin{x}{\small\bf LEMMA} \ 
If $f \in \sB_0 (A)$ is not a constant, then 
\[
S (r) 
\ = \ 
\sum\limits_{\abs{z_n} \hsy \leq \hsy r} \ 
\frac{1}{z_n}
\]
remains bounded as $r \ra \infty$.
\\[-.5cm]

[One can extract a proof from the material in \S6.  
To proceed directly, assume for convenience that 
$\abs{f(0)} = 1$ and choose $K > 0$ : $n (r) \leq K \hsy r$ (cf. 27.6) $-$then 
\[
\abs{S (r) - S (R)}
\ \leq \ 
2 \hsy K 
\qquad (R \leq r \leq 2 \hsy R)
\]
\qquad \qquad
$\implies$
\[
\int\limits_R^{2 \hsy R} \ 
S (r) 
\hsx 
r 
\ \td r
\ = \ 
\frac{3}{2} 
\hsx
R^2 \hsy S (R) + \tO(R^2).
\]
Under the supposition that $f (z)$ is zero free on $\abs{z} = r$, write
\allowdisplaybreaks
\begin{align*}
S (r) \ 
&=\ 
\frac{1}{2 \hsy \pi \sqrt{-1}} 
\ 
\int\limits_C \ 
\frac{f^\prime (z)}{f (z)} 
\hsx \cdot \hsx 
\frac{1}{z}
\ \td z
\ - \ 
\frac{f^\prime (0)}{f (0)} 
\\[15pt]
&=\ 
\frac{1}{2 \hsy \pi} \ 
\int\limits_0^{2 \hsy \pi}  \ 
\Big(
\frac{\partial}{\partial x} - \sqrt{-1} \hsx \frac{\partial}{\partial y} 
\Big)
\hsx 
\log \big | f \big(r \hsy e^{\sqrt{-1} \hsx \theta}\big)\big |
\ \td \theta 
\ - \ 
\frac{f^\prime (0)}{f (0)} 
\end{align*}
\qquad 
$\implies$
\allowdisplaybreaks
\begin{align*}
\frac{3}{2} 
\hsx
R^2 \hsy S (R) 
&=\ 
\int\limits_R^{2 \hsy R} \ 
S (r) 
\hsx 
r 
\ \td r
\hsx + \hsx
\tO(R^2) \ 
\\[15pt]
&=\ 
\frac{1}{2 \hsy \pi} \  
\iint\limits_{R \hsy \leq \hsy \abs{z} \hsy \leq \hsy 2 R} \ 
\Big(
\frac{\partial}{\partial x} - \sqrt{-1} \hsx \frac{\partial}{\partial y} 
\Big)
\hsx 
\log \abs{f (z)} 
\ \td x \hsy \td y
\hsx + \hsx
\tO (R^2)
\\[15pt]
&=\ 
\frac{1}{2 \hsy \pi}   
\int\limits_0^{2 \hsy \pi}  \hsy 
\big(
2 R \log 
\big|
f (2 \hsy R \hsy e^{\sqrt{-1} \hsx \theta})
\big|
\hsy - \hsy 
R 
\log 
\big|
f (R \hsy e^{\sqrt{-1} \hsx \theta})
\big|
\big)
\hsx 
e^{- \sqrt{-1} \hsx \theta}
\hsy \td \theta
\hsx + \hsx
\tO (R^2)
\end{align*}
$\implies$ 
\allowdisplaybreaks
\begin{align*}
\frac{3}{2} 
R^2 &\hsy \abs{S (R)} \ 
\\[15pt]
&\leq \
\frac{R}{2 \hsy \pi}  
\int\limits_0^{2 \hsy \pi}  
\big(
2  
\hsx 
\big| 
\log \big|f (2 \hsy R \hsy e^{\sqrt{-1} \hsx \theta})
\big| 
\hsy
\big| 
\hsx + \hsx 
\big|
\log \big|f (R \hsy e^{\sqrt{-1} \hsx \theta})
\big|
\hsy
\big|
\big)
\ \td \theta
\ + \ 
\tO (R^2).
\end{align*}
Estimating the integral in the usual way gives rise to another $\tO (R^2)$, 
so in the end
\[
\frac{3}{2} 
R^2 \hsy \abs{S (R)}
\ \leq \ 
\tO(R^2)
\]
\qquad \qquad
$\implies$
\[
\abs{S (R)}
\ \leq \ 
\tO(1)
\qquad (R \ra \infty).]
\]
\\[-1.cm]
\end{x}

\begin{x}{\small\bf CARLEMAN FORMULA} \ 
Suppose that $f (z)$ is holomorphic for $\Img z \geq 0$ and let 
$z_k = r_k \hsy e^{\sqrt{-1} \hsx \theta_k}$ $(k = 1, \ldots, n)$ be its zeros in the region
\[
\{z : \Img z \geq 0, \ 
1 \leq \abs{z} \leq R\}.
\]
Then
\allowdisplaybreaks
\begin{align*}
\sum\limits_{k = 1}^n \ 
\Big(
\frac{1}{r_k} - \frac{r_k}{R^2}
\Big)
\hsx 
\sin \theta_k 
&=\ 
\frac{1}{\pi \hsy R} \ 
\int\limits_0^\pi \ 
\log 
\big|
f (R \hsy e^{\sqrt{-1} \hsx \theta})
\big|
\hsx
\sin \theta
\ \td \theta
\\[15pt]
&\hspace{1cm}
\hsx + \hsx 
\frac{1}{2 \hsy \pi} \ 
\int\limits_1^R \ 
\Big(
\frac{1}{x^2} - \frac{1}{R^2}
\Big)
\hsx
\log \abs{f (x) \hsy f(-x)} 
\ \td x
\hsx + \hsx 
A (R),
\end{align*}
where $A (R)$ is a bounded function of $R$.
\\[-.5cm]

[Note: \ 
Replace 1 by $\rho > 0$ $-$then $A (R)$ depends on $\rho$ and 
\[
A(\rho, R) 
\ = \ 
- 
\Img 
\frac{1}{2 \hsy \pi} \ 
\int\limits_0^\pi \ 
\log f \big(\rho \hsy e^{\sqrt{-1} \hsx \theta}\big)
\hsx
\Big(
\frac{\rho \hsy e^{\sqrt{-1} \hsx \theta}}{R^2}
\hsx - \hsx 
\frac{e^{-\sqrt{-1} \hsx \theta}}{\rho}
\Big)
\ \td \theta,
\]
thus if $f (0) = 1$, 
\[
\lim\limits_{\rho \ra 0}
A(\rho, R) 
\ = \ 
\frac{1}{2} \ 
\Img 
f^\prime (0),
\]
so
\allowdisplaybreaks
\begin{align*}
\sum\limits_{r_k \hsy \leq \hsy R} \ 
&
\Big(
\frac{1}{r_k} - \frac{r_k}{R^2}
\Big)
\hsx 
\sin \theta_k \ 
\\[15pt]
&=\ 
\frac{1}{\pi \hsy R} \ 
\int\limits_0^\pi \ 
\log 
\big|
f (R \hsy e^{\sqrt{-1} \hsx \theta})
\big|
\hsx
\sin \theta
\ \td \theta
\\[15pt]
&
\hspace{1cm}
\hsx + \hsx 
\frac{1}{2 \hsy \pi} \ 
\int\limits_0^R \ 
\Big(
\frac{1}{x^2} - \frac{1}{R^2}
\Big)
\hsx
\log \abs{f (x) \hsy f(-x)} 
\ \td x
\hsx + \hsx 
\frac{1}{2} \ 
\Img 
f^\prime (0).]
\end{align*}
\\[-1cm]
\end{x}

\begin{x}{\small\bf THEOREM} \ 
If $f \in \sB_0 (A)$ is not a constant, then the series
\[
\sum\limits_{n = 1}^\infty \ 
\frac{\sin \theta_n}{r_n}
\]
\\[-1cm]
is absolutely convergent.
\\[-.5cm]

PROOF \ 
Apply 27.9 to $f (z)$, $f (-z)$ and add the results.  
In this way we are led to
\allowdisplaybreaks 
\begin{align*}
\sum\limits_{k = 1}^n \ 
\Big(
\frac{1}{r_k} - \frac{r_k}{R^2}
\Big)
\hsx 
&
\sin \theta_k \qquad
(0 \leq \theta_k \leq \pi)
\\[11pt]
&
\ + \ 
\sum\limits_{\ell = 1}^m \ 
\Big(
\frac{1}{r_\ell} - \frac{r_\ell}{R^2}
\Big)
\hsx 
\sin (\theta_\ell + \pi) \qquad
(- \pi \leq \theta_\ell \leq 0).
\end{align*}
But 
$\sin \theta_k = \abs{\sin \theta_k}$, \
$\sin (\theta_\ell + \pi) = -\sin \theta_\ell = \abs{\sin \theta_\ell }$, 
hence 
\[
\sum\limits_{r_n \hsy \leq \hsy R} \ 
\Big(
1 - \frac{r_n^2}{R^2}
\Big)
\hsx 
\frac{\abs{\sin \theta_n}}{r_n}
\ < \ 
C
\qquad (R \gg 0)
\]
for some constant $C > 0$.  
And this implies that
\[
\sum\limits_{r_n \hsy \leq \hsy R/2} \ 
\Big(
1 - \frac{1}{4}
\Big)
\hsx 
\frac{\abs{\sin \theta_n}}{r_n}
\ < \ 
C.
\]
Now send $R$ to $\infty$.
\\[-.5cm]

[Note: \ 
The zeros on the real axis do not figure in the calculation.]
\\[-.25cm]
\end{x}

\qquad
{\small\bf \un{N.B.}} \ 
Restated, 27.10 says that 
\[
\sum\limits_{n = 1}^\infty \ 
\Big|
\Img \frac{1}{z_n}
\Big|
\ < \ 
\infty.
\]

[Note: \ 
In traditional terminology, an entire function $f$ of exponential type is said to be 
\un{class A} 
if
\[
\sum\limits_{n = 1}^\infty \ 
\Big|\Img \frac{1}{z_n}\Big|
\ < \ 
\infty.
\]
\\[-.75cm]

\noindent
Characterization: \ 
$f$ is class A iff
\[
\sup\limits_{R > 1} \ 
\int\limits_1^R \ 
\frac{\log \abs{f (x) \hsy f(-x)}}{x^2}
\ \td x 
\ < \ 
\infty.]
\]
\\[-.75cm]

\begin{x}{\small\bf APPLICATION} \ 
Given $\varepsilon > 0$, let $\Omega (\varepsilon)$ be the sector 

\[
\abs{\arg z} 
\ < \ 
\varepsilon
\ \cup \ 
\abs{\arg z - \pi } 
\ < \ 
\varepsilon.
\]
Then
\[
\sum\limits_{k = 1}^\infty \ 
\frac{1}{\abs{z_{n_k}}}
\ < \ 
\infty,
\]
where $z_{n_k}$ runs through the zeros of $f$ which are not in $\Omega (\varepsilon)$.
\\[-.25cm]
\end{x}

\begin{x}{\small\bf THEOREM} \ 
If $f \in \sB_0 (A)$ is not a constant, then 
\[
\lim\limits_{r \ra \infty} \ 
\frac{n (r)}{r}
\ = \ 
\frac{h_f (\sqrt{-1}) + h_f (-\sqrt{-1})}{\pi}.
\]

[This is a substantial reinforcement of 27.6.  
For a proof, consult 
B. 
Levin\footnote[2]{\vspace{.11 cm}
\textit{Lectures on Entire Functions}, A.M.S., 1996, pp. 127-130.} 
(see also 
P. 
Koosis
\footnote[3]{\vspace{.11 cm}
\textit{The Logarithmic Integral $\tI$}, Cambridge University Press, 1988, pp. 69-76.}).]
\\[-.25cm]
\end{x}

\begin{x}{\small\bf REMARK} \ 
One can say more.  
Thus let $n_+ (r)$ be the number of zeros of $f$ with real part $\geq 0$ and modulus $\leq r$ 
and let $n_- (r)$ be the number of zeros of $f$ with real part $< 0$ and modulus $\leq r$ $-$then 

\[
n (r) 
\ = 
\ n_+ (r) \hsx + \hsx n_- (r).
\]
Moreover, it can be shown that
\[
\lim\limits_{r \ra \infty} \ 
\frac{n_+ (r)}{r}
\ = \ 
\frac{h_f (\sqrt{-1}) + h_f (-\sqrt{-1})}{2 \hsy \pi}
\]
and
\[
\lim\limits_{r \ra \infty} \ 
\frac{n_- (r)}{r}
\ = \ 
\frac{h_f (\sqrt{-1}) + h_f (-\sqrt{-1})}{2 \hsy \pi}.
\]
\\[-.75cm]
\end{x}

\begin{x}{\small\bf EXAMPLE} \ 
Take 
$f (z) = e^{\sqrt{-1} \hsx z}$ $-$then $n (r) \equiv 0$.  
On the other hand, 
\[
h_f (\sqrt{-1})
\ = \ 
\underset{r \ra \infty}{\limsupx} \ 
\frac
{\log \big|e^{\sqrt{-1} \hsx (\sqrt{-1} \hsx r)}\big|}
{r}
\ = \ 
\underset{r \ra \infty}{\limsupx} \ 
\frac{\log e^{-r}}{r}
\ = \ -1
\]
and
\[
h_f (-\sqrt{-1})
\ = \ 
\underset{r \ra \infty}{\limsupx} \ 
\frac
{\log \big|e^{\sqrt{-1} \hsx (-\sqrt{-1} \hsx r)}\big|}
{r}
\ = \ 
\underset{r \ra \infty}{\limsupx} \ 
\frac{\log e^r}{r}
\ = \ 1.
\]
Therefore
\[
h_f (\sqrt{-1}) + h_f (-\sqrt{-1})
\ = \ 
- 1 + 1 
\ = \ 
0.
\]
\\[-1cm]
\end{x}

\begin{x}{\small\bf LEMMA}\footnote[4]{\vspace{.11 cm}
R. Boas \textit{Entire Functions}, Academic Press, 1954, p. 116.}  \ 
If $f \in \sB_0 (A)$ is not a constant, then 
\[
H_f (1) = 0
\quad \text{and} \quad 
H_f (-1) = 0
\]
or still, 
\[
h_f (1)
\ = \ 
\underset{r \ra \infty}{\limsupx} \ 
\frac{\log \abs{f (r)}}{r}
\ = \ 
0
\]
and 
\[
h_f (-1)
\ = \ 
\underset{r \ra \infty}{\limsupx} \ 
\frac{\log \abs{f (-r)}}{r}
\ = \ 
0.
\]

[Note: \ 
This result is a consequence of ``Ahlfors-Heins theory'' and is valid for any entire function $f$ 
of exponential type in the Cartwright class, i.e., such that
\[
\int\limits_{-\infty}^\infty \ 
\frac{\log^+ \abs{f (x)}}{1 + x^2}
\ \td x
\ < \ 
\infty.]
\]
\\[-1.cm]
\end{x}

\begin{x}{\small\bf COROLLARY} \ 
The indicator diagram $K_f$ of $f$ is a segment of the imaginary axis (or a point) (cf. 18.9).
\\[-.25cm]
\end{x}


\begin{x}{\small\bf LEMMA} \ 
Let 
$K = [\sqrt{-1} \hsx A, \sqrt{-1} \hsx B]$ $(A \leq B)$ $-$then
\[
H_K \big(e^{\sqrt{-1} \hsx \theta}\big) 
\ = \ 
a \hsy \abs{\sin \theta} + b \sin \theta, 
\]
where
\[
a 
\ = \ 
\frac{B - A}{2}, 
\quad
b = \frac{- B - A}{2}.
\]
\\[-1cm]
\end{x}

\begin{x}{\small\bf EXAMPLE} \ 
Take $A = B$, call it $C$ $-$then 
\[
a \ = \ \frac{C - C}{2} \ = \ 0, 
\quad 
b \ = \ \frac{-C - C}{2} \ = \ -C
\]
and
\[
H_K \big(e^{\sqrt{-1} \hsx \theta}\big) 
\ = \ 
-C \sin \theta
\qquad (\tcf. \  18.2).
\]
\\[-1.cm]
\end{x}

\begin{x}{\small\bf EXAMPLE} \ 
Take $A = -c$, $B = c$ with $c > 0$, then 

\[
a 
\ = \ 
\frac{c - (-c)}{2} \ = \ c, 
\quad
b = \frac{- c + c}{2} \ = \ 0
\]
and
\[
H_K \big(e^{\sqrt{-1} \hsx \theta}\big) 
\ = \ 
a \hsy \abs{\sin \theta}
\qquad (\tcf. \  18.5).
\]
\\[-.25cm]
\end{x}

\begin{x}{\small\bf RAPPEL} \ 
If $f \in \sB_0 (A)$ is not a constant, then 
\[
T (f) 
\ = \ 
\tau (f) 
\ = \ 
\sup\limits_{0 \hsy \leq \hsy \theta \leq 2 \hsy \pi} \ 
h_f \big(e^{\sqrt{-1} \hsx \theta}\big)
\qquad (\tcf. \ 19.10).
\]

Recalling that $H_f \ (= H_{K_f}$ (cf. 18.17)) $= h_f$ (cf. 19.7), we have
\allowdisplaybreaks
\begin{align*}
\sup\limits_{0 \hsy \leq \hsy \theta \leq 2 \hsy \pi} \ 
h_f \big(e^{\sqrt{-1} \hsx \theta}\big) \ 
&=\ 
\sup\limits_{0 \hsy \leq \hsy  \theta \leq 2 \hsy \pi} \ 
(a \hsy \abs{\sin \theta} + b \sin \theta)
\\[11pt]
&=\ 
\max (a + b, a - b)
\\[11pt]
&=\ 
a + \abs{b}.
\end{align*}
But
\[
\begin{cases}
\ 
a + b 
\ = \ 
h_f \big(\sqrt{-1}\hsx\big)
\\[8pt]
\ 
a - b 
\ = \ 
h_f \big(-\sqrt{-1}\hsx\big)
\end{cases}
.
\]
Therefore
\[
T (f) 
\ = \ 
\max \big(h_f (\sqrt{-1}), h_f (-\sqrt{-1})\big).
\]
\\[-1.cm]
\end{x}

\begin{x}{\small\bf SCHOLIUM} \ 
If 
$h_f \big(\sqrt{-1}\hsx\big) = h_f \big(-\sqrt{-1}\hsx\big)$, 
then
\allowdisplaybreaks
\begin{align*}
\lim\limits_{r \ra \infty} \ 
\frac{n (r)}{r} \ 
&=\ 
\frac{h_f \big(\sqrt{-1}\hsx\big) + h_f \big(-\sqrt{-1}\hsx\big)}{\pi}
\qquad (\tcf. \ 27.12)
\\[15pt]
&=\ 
2 \hsx \frac{T (f)}{\pi}.
\end{align*}
\\[-1.25cm]
\end{x}

\begin{x}{\small\bf LEMMA} \ 
\[
K_f 
\ = \ 
\big[
\sqrt{-1} \hsx 
\big(
-
h_f \big(\sqrt{-1}\hsx\big), 
\sqrt{-1} \hsx 
h_f \big(-\sqrt{-1}\hsx\big)
.
\big]
\]

PROOF \ 
Writing 
$K_f = [\sqrt{-1} \hsx A, \sqrt{-1} \hsx B]$, 
it is a question of explicating $A$ and $B$.  
But
\[
\begin{cases}
\ 
a + b 
\ = \ 
h_f \big(\sqrt{-1}\hsx\big)
\\[8pt]
\ 
a - b 
\ = \ 
h_f \big(-\sqrt{-1}\hsx\big)
\end{cases}
.
\]
And
\[
\begin{cases}
\ \ds
a 
\ = \ 
\frac{B - A}{2}
\\[18pt]
\ \ds
b 
\ = \ 
\frac{-B - A}{2}
\end{cases}
\quad \implies \quad
\begin{cases}
\ \ds
\frac{B - A}{2} + \frac{-B - A}{2}
\ = \ 
- A
\\[18pt]
\ \ds
\frac{B - A}{2} - \frac{-B - A}{2}
\ = \ 
B
\end{cases}
\]

\[
\hspace{3cm}
\implies
\quad
\begin{cases}
\ \ds
- A 
\ = \ 
h_f \big(\sqrt{-1}\hsx\big)
\\[18pt]
\ \ds
\ B
\ = \ 
h_f \big(-\sqrt{-1}\hsx\big)
\end{cases}
\]
\qquad \qquad
$\implies$
\[
K_f 
\ = \ 
\big[
\sqrt{-1} \hsx 
\big(
-
h_f \big(\sqrt{-1}\hsx\big), 
\sqrt{-1} \hsx 
h_f \big(-\sqrt{-1}\hsx\big)
.
\big]
\]
\\[-1.cm]
\end{x}

\begin{x}{\small\bf APPLICATION} \ 
$K_f$ reduces to a point iff 
\[
h_f \big(\sqrt{-1}\hsx\big)
\hsx + \hsx 
h_f \big(-\sqrt{-1}\hsx\big)
\ = \ 
0,
\]
hence $K_f$ reduces to a point iff
\[
\lim\limits_{r \ra \infty} \ 
\frac{n (r)}{r} \ 
\ =\ 
0.
\]
\\[-1.cm]
\end{x}

\begin{x}{\small\bf EXAMPLE} \ 
Suppose that $c \neq 0$ is real and let 
$
\ds
f (z) = e^{\sqrt{-1} \hsx c \hsy z}
$
$-$then
\[
h_f \big(e^{\sqrt{-1} \hsx \theta}\big) 
\ = \ 
- c \hsy \sin \theta
\qquad (\tcf. \ 19.2)
\]
\qquad \qquad
$\implies$
\[
\qquad
\begin{cases}
\ 
h_f \big(\sqrt{-1}\hsx\big)
\hspace{0.5cm}
\ = \ 
-c
\\[11pt]
\ 
h_f \big(-\sqrt{-1}\hsx\big)
\ = \ 
c
\end{cases}
\quad \implies \quad
K_f 
\ = \ 
\{\sqrt{-1}\hsx c\}.
\]
\\[-.5cm]

\noindent
And $T (f) = \abs{c}$.
\\[-.25cm]
\end{x}

\begin{x}{\small\bf EXAMPLE} \ 
Suppose that $F \neq I$ is a finite distribution function, $\gf$ its characteristic function (cf. 27.2) $-$then

\[
\begin{cases}
\ 
\rext [F] 
\ = \ 
h_\gf \big(-\sqrt{-1}\big) 
\\[11pt]
\ 
\lext [F] 
\ = \ 
-
h_\gf \big(\sqrt{-1}\big) 
\end{cases}
\qquad 
(\tcf. \ 26.8)
\]
and

\[
-h_\gf \big(\sqrt{-1}\big) 
\ \leq \ 
h_\gf \big(-\sqrt{-1}\big) 
\]
in agreement with 27.22 (cf. 22.13).
\\[-.5cm]

[Note: \ 
Recall too that
\[
T (f) 
\ = \ 
\max (-\lext [F] , \rext [F])
\qquad 
(\tcf. \ 26.11).]
\]
\\[-1.25cm]
\end{x}

\begin{x}{\small\bf EXAMPLE} \ 
Given $\phi \in \Lp^1 [-A, A]$ $(0 < A < \infty)$, put
\[
f (z) 
\ = \ 
\frac{1}{\sqrt{2 \hsy \pi}} \ 
\int\limits_{-A}^A \ 
\phi (t) 
\hsx 
e^{\sqrt{-1} \hsx z \hsy t}
\ \td t.
\]
Then $f \in \sB_0(A)$ (cf. 17.19).  
Assume further that $\phi (t)$ does not vanish almost everywhere in any neighborhood of $A$ (or $-A$) $-$then
\[
\begin{cases}
\ \
A
\ \ = \ 
h_f \big(-\sqrt{-1}\big) 
\\[11pt]
\ 
- A
\ = \ 
- h_f \big(\sqrt{-1}\big) 
\end{cases}
\implies
T (f) = A
\]
\\[-.5cm]

\qquad 
$\implies$
\allowdisplaybreaks
\begin{align*}
\lim\limits_{r \ra \infty} \ 
\frac{n (r)}{r} \ 
&=\ 
2 \hsx \frac{T (f)}{\pi} 
\qquad (\tcf. \ 27.21)
\\[11pt]
&=\ 2 \hsx \frac{A}{\pi}.
\end{align*}
\\[-1cm]
\end{x}

\begin{x}{\small\bf NOTATION} \ 
Put
\[
D 
\ = \ 
\frac{h_f \big(\sqrt{-1}\big)  \hsx + \hsx  h_f \big(-\sqrt{-1}\big) }{\pi}.
\]
\\[-1.25cm]
\end{x}


\begin{x}{\small\bf DEFINITION} \ 
The zeros of $f$ \un{have a density} if $D > 0$.
\\[-.25cm]
\end{x}

\begin{x}{\small\bf RAPPEL} \ 
Take $\alpha > 0$ $-$then the series
\[
\sum\limits_{n = 1}^\infty \ 
\frac{1}{r_n^\alpha}
\]
converges iff the integral
\[
\int\limits_0^\infty \ 
\frac{n(t)}{t^{\alpha + 1}}
\ \td t
\]
converges.
\\[-.25cm]
\end{x}

\begin{x}{\small\bf LEMMA} \ 
If the zeros of $f$ have a density, then the series
\[
\sum\limits_{n = 1}^\infty \ 
\frac{1}{r_n}
\]
is divergent.
\\[-.5cm]

[In 27.29, take $\alpha = 1$: 
\allowdisplaybreaks
\begin{align*}
\int\limits_0^\infty \ 
\frac{n(t)}{t^{2}}
\ \td t \ 
&=\ 
\int\limits_0^\infty \ 
\frac{n(t)}{t}
\
\cdot
\ \frac{\td t}{t}
\\[15pt]
&=\ 
\int\limits_0^\infty \ 
\frac{(n(t)/ t)}{D}
\cdot
D 
\hsx
\ \frac{\td t}{t}
\end{align*}
is divergent (cf. 27.12).]
\\[-.5cm]

[Note: \ 
The convergence exponent is equal to 1 (cf. 4.10).  
Therefore $f$ is of divergence class (cf. 4.24).]
\\[-.25cm]
\end{x}

\begin{x}{\small\bf THEOREM} \ 
If $f \in \sB_0(A)$ is not a constant and if the zeros of $f$ have a density, then the series
\[
\sum\limits_{n = 1}^\infty \ 
\frac{\cos \theta_n}{r_n}
\]
is convergent.
\\[-.25cm]
\end{x}


\begin{x}{\small\bf REMARK} \ 
According to 27.10, the series
\[
\sum\limits_{n = 1}^\infty \ 
\frac{\sin \theta_n}{r_n}
\]
is absolutely convergent.  
On the other hand, in view of 27.30, the series

\[
\sum\limits_{n = 1}^\infty \ 
\frac{\cos \theta_n}{r_n}
\]
is not absolutely convergent.  
\\[-.25cm]
\end{x}

Before tackling the proof, we shall first set up the relevant generalities.
\\[-.25cm]

\begin{x}{\small\bf RAPPEL} \ 
Given a sequence $a_1, a_2, \ldots$, put
\[
\sigma_n
\ = \ 
\frac{a_1 + a_2 + \cdots + a_n}{n}.
\]
Assume: \ 
$
\ds
\lim\limits_{n \ra \infty} \ 
a_n = 0
\ 
$
$-$then 
$
\ 
\ds
\lim\limits_{n \ra \infty} \ 
\sigma_n = 0.
$
\\[-.25cm]
\end{x}

\begin{x}{\small\bf APPLICATION} \ 
If $a_n \ra L$, then $\sigma_n \ra L$.  
\\[-.5cm]

[In fact, $a_n - L \ra 0$, so 
\[
\frac{(a_1 - L) + (a_2 - L) + \cdots + (a_n - L)}{n}
\ra 
0
\]
or still, 
$\sigma_n - L \ra 0$.]
\\[-.25cm]
\end{x}

\begin{x}{\small\bf RAPPEL} \ 
Given an infinite series 
$
\ds
\sum\limits_1^\infty \ 
a_n
$, 
let $s_n$ denote its $n^\nth$ partial sum and put
\[
\sigma_n
\ = \ 
\frac{s_1 + s_2 + \cdots + s_n}{n}.
\]
Assume: \ 
$\{\sigma_n\}$ converges to $S$ and 
$
\ds
a_n = \tO \Big(\frac{1}{n}\Big)
$
$-$then 
$\{s_n\}$ 
converges to $S$ .
\\[-.25cm]

[Note: \ 
In other words , if \ 
$
\ds
\sum\limits_1^\infty \ 
a_n
$ 
is $(C,1)$ summable to $S$ and if 
$
\ds
a_n = \tO \Big(\frac{1}{n}\Big)
$,  
then \ 
$
\ds
\sum\limits_1^\infty \ 
a_n
$ 
is convergent to $S$.]
\\
\end{x}

\qquad
{\small\bf \un{N.B.}} \ 
\[
\frac{\cos \theta_n}{r_n}
\ = \ 
\tO \Big(\frac{1}{n}\Big).
\]

[For 
\[
\frac{n (r_n)}{r_n} 
\ = \ 
\frac{n}{r_n}
\ra D.]
\]
\\[-.75cm]

\begin{x}{\small\bf JENSEN FORMULA} \ 
Suppose that $f (z)$ is holomorphic in $\abs{z} < R$ with $f (0) = 1$ $-$then 

\[
\int\limits_0^r \ 
\frac{n(t)}{t}
\ \td t
\ = \ 
\frac{1}{2 \hsy \pi} \ 
\int\limits_0^{2 \hsy \pi} \ 
\log \big| f \big(r \hsy e^{\sqrt{-1} \hsx \theta} \big)\big|
\ \td \theta
\qquad (0 < r < R).
\]
\\[-1cm]
\end{x}

\begin{x}{\small\bf CARLEMAN FORMULA} \ 
(bis) \ 
Suppose that $f (z)$ is holomorphic for $\Reg z \geq 0$ and let 
$
\ds
z_k = r_k \hsy e^{\sqrt{-1} \hsx \theta_k}
$ 
$(k = 1, \ldots, n)$ 
be its zeros in the region
\[
\{z : \Reg z \geq 0, 1 \leq \abs{z} \leq R\}.
\]
Then
\allowdisplaybreaks
\begin{align*}
\sum\limits_{k = 1}^n
&
\Big(
\frac{1}{r_k} - \frac{r_k}{R^2}
\Big) 
\hsx 
\cos \theta_k
\ 
\\[15pt]
&=\ 
\frac{1}{\pi \hsy R} \ 
\int\limits_{-\frac{\pi}{2}}^{\frac{\pi}{2}} \ 
\log 
\big |
f \big(R e^{\sqrt{-1} \hsx \theta}\big)
\big | 
\hsx 
\cos \theta 
\ \td \theta
\\[15pt]
&
\hspace{1cm}
+\ 
\frac{1}{2 \hsy \pi} \ 
\int\limits_1^R \ 
\Big(
\frac{1}{x^2} - \frac{1}{R^2}
\Big) 
\hsy
\log \abs{f (\sqrt{-1} \hsx x) \hsy f(-\sqrt{-1} \hsx x)}
\ \td x 
\ + \ 
A (R), 
\end{align*}
where $A (R)$ is a bounded function of $R$.  
\\[-.375cm]

[Note: \ 
If $f (0) = 1$, then 
\allowdisplaybreaks
\begin{align*}
\sum\limits_{r_k \hsy \leq \hsy R}
&
\Big(
\frac{1}{r_k} - \frac{r_k}{R^2}
\Big) 
\hsx 
\cos \theta_k
\ 
\\[15pt]
&=\ 
\frac{1}{\pi \hsy R} \ 
\int\limits_{-\frac{\pi}{2}}^{\frac{\pi}{2}} \ 
\log 
\big |
f \big(R e^{\sqrt{-1} \hsx \theta}\big)
\big | 
\hsx 
\cos \theta 
\ \td \theta
\\[15pt]
&
\hspace{1cm}
+\
\frac{1}{2 \hsy \pi} \ 
\int\limits_0^R \ 
\Big(
\frac{1}{x^2} - \frac{1}{R^2}
\Big) 
\hsy
\log \abs{f (\sqrt{-1} \hsx x) \hsy f(-\sqrt{-1} \hsx x)}
\ \td x 
\hsx - \hsx 
\frac{1}{2} \ 
\Reg f^\prime (0).]
\end{align*}
\\[-1cm]
\end{x}

Proceeding to the proof of 27.31, it will be assumed that $f (0) = 1$.
\\[-.25cm]

[Note: \ 
Zeros of $f (z)$ on the imaginary axis do not participate 
(
$
\ds
\cos \Big(\pm \hsy\frac{\pi}{2}\Big)
= 0
$
).]
\\[-.25cm]

\un{Step 1:}\ 
In the formula
\[
\sum\limits_{r_k \hsy \leq \hsy R}
\Big(
\frac{1}{r_k} - \frac{r_k}{R^2}
\Big) 
\hsx 
\cos \theta_k 
\hsx + \hsx 
\frac{1}{2} \ 
\Reg f^\prime (0)
\ = \ 
\ldots, 
\]
replace $f (z)$ by $f (-z)$ to get
\allowdisplaybreaks
\begin{align*}
\sum\limits_{r_\ell \hsy \leq \hsy R}
&
\Big(
\frac{1}{r_\ell} - \frac{r_\ell}{R^2}
\Big) 
\hsx 
\cos (\theta_\ell + \pi) 
\hsx - \hsx 
\frac{1}{2} \ 
\Reg f^\prime (0)
\\[15pt]
&=\ 
\frac{1}{\pi \hsy R} \ 
\int\limits_{-\frac{\pi}{2}}^{\frac{\pi}{2}} \ 
\log 
\big |
f \big(- R e^{\sqrt{-1} \hsx \theta}\big)
\big |
\hsx 
\cos \theta 
\ \td \theta
\\[15pt]
&
\hspace{1cm}
\ + \ 
\frac{1}{2 \hsy \pi} \ 
\int\limits_0^R \ 
\Big(
\frac{1}{x^2} - \frac{1}{R^2}
\Big) 
\hsy
\log 
\big |
f (-\sqrt{-1} \hsx x) \hsy f(\sqrt{-1} \hsx x)
\big |
\ \td x 
\end{align*}
or still, 
\[
- 
\sum\limits_{r_\ell \hsy \leq \hsy R}
\Big(
\frac{1}{r_\ell} - \frac{r_\ell}{R^2}
\Big) 
\hsx 
\cos \theta_\ell
\hsx - \hsx 
\frac{1}{2} \ 
\Reg f^\prime (0)
\ = \ 
\cdots \hsx .
\]
Therefore
\allowdisplaybreaks
\begin{align*}
\sum\limits_{r_k \hsy \leq \hsy R}
&
\Big(
\frac{1}{r_k} - \frac{r_k}{R^2}
\Big) 
\hsx 
\cos \theta_k 
\hsx + \hsx 
\frac{1}{2} \ 
\Reg f^\prime (0)
\hsx + \hsx
\sum\limits_{r_\ell \hsy \leq \hsy R}
\Big(
\frac{1}{r_\ell} - \frac{r_\ell}{R^2}
\Big) 
\hsx 
\cos \theta_\ell
\hsx + \hsx 
\frac{1}{2} \ 
\Reg f^\prime (0)
\\[15pt]
&=\ 
\frac{1}{\pi \hsy R} \ 
\int\limits_{-\frac{\pi}{2}}^{\frac{\pi}{2}} \ 
\log 
\big |
f \big(R e^{\sqrt{-1} \hsx \theta}\big)
\big |
\hsx 
\cos \theta 
\ \td \theta
\ - \ 
\frac{1}{\pi \hsy R} \ 
\int\limits_{-\frac{\pi}{2}}^{\frac{\pi}{2}} \ 
\log 
\big |
f \big(- R e^{\sqrt{-1} \hsx \theta}\big)
\big |
\hsx 
\cos \theta 
\ \td \theta.
\end{align*}
\\[-.75cm]

\un{Step 2:}\ 
\allowdisplaybreaks
\begin{align*}
-
\frac{1}{\pi \hsy R} \ 
\int\limits_{-\frac{\pi}{2}}^{\frac{\pi}{2}} \ 
\log 
&
\big |
f \big(- R e^{\sqrt{-1} \hsx \theta}\big)
\big |
\hsx 
\cos \theta 
\ \td \theta \ 
\\[15pt]
&= \  
-
\frac{1}{\pi \hsy R} \ 
\int\limits_{-\frac{\pi}{2}}^{\frac{\pi}{2}} \ 
\log 
\big |
f \big(R e^{\sqrt{-1} \hsx (\theta + \pi)}\big)
\big |
\hsx 
\cos \theta 
\ \td \theta
\\[15pt]
&= \  
-
\frac{1}{\pi \hsy R} \ 
\int\limits_{\frac{\pi}{2}}^{\frac{3 \hsy \pi}{2}} \ 
\log 
\big |
f \big(R e^{\sqrt{-1} \hsx \theta}\big)
\big |
\hsx 
\cos (\theta - \pi)
\ \td \theta
\\[15pt]
&= \  
\frac{1}{\pi \hsy R} \ 
\int\limits_{\frac{\pi}{2}}^{\frac{3 \hsy \pi}{2}} \ 
\log 
\big |
f \big(R e^{\sqrt{-1} \hsx \theta}\big)
\big |
\hsx 
\cos \theta 
\ \td \theta.
\end{align*}

\un{Step 3:}\quad
Therefore
\allowdisplaybreaks
\begin{align*}
\sum\limits_{r_k \hsy \leq \hsy R}
&
\Big(
\frac{1}{r_k} - \frac{r_k}{R^2}
\Big) 
\hsx 
\cos \theta_k 
\hsx + \hsx 
\frac{1}{2} \ 
\Reg f^\prime (0)
\\[15pt]
&
\hspace{2cm}
\hsx + \hsx
\sum\limits_{r_\ell \hsy \leq \hsy R}
\Big(
\frac{1}{r_\ell} - \frac{r_\ell}{R^2}
\Big) 
\hsx 
\cos \theta_\ell
\hsx + \hsx 
\frac{1}{2} \ 
\Reg f^\prime (0)
\\[15pt]
&=\  
\frac{1}{\pi \hsy R} \ 
\int\limits_{-\frac{\pi}{2}}^{\frac{\pi}{2}} \ 
\log 
\big |
f \big(R e^{\sqrt{-1} \hsx \theta}\big)
\big |
\hsx 
\cos \theta 
\ \td \theta
\\[15pt]
&
\hspace{2cm}
\ + \ 
\frac{1}{\pi \hsy R} \ 
\int\limits_{\frac{\pi}{2}}^{\frac{3 \hsy \pi}{2}} \ 
\log 
\big |
f \big(R e^{\sqrt{-1} \hsx \theta}\big)
\big |
\hsx 
\cos \theta 
\ \td \theta
\\[15pt]
&=\  
\frac{1}{\pi \hsy R} \ 
\int\limits_{-\frac{\pi}{2}}^0 \ 
\log 
\big |
f \big(R e^{\sqrt{-1} \hsx \theta}\big)
\big |
\hsx 
\cos \theta 
\ \td \theta
\\[15pt]
&
\hspace{2cm}
\ + \ 
\frac{1}{\pi \hsy R} \ 
\int\limits_0^{\frac{\pi}{2}} \ 
\log 
\big |
f \big(R e^{\sqrt{-1} \hsx \theta}\big)
\big |
\hsx 
\cos \theta 
\ \td \theta
\\[15pt]
&
\hspace{3cm}
\ + \ 
\frac{1}{\pi \hsy R} \ 
\int\limits_{\frac{\pi}{2}}^{\frac{3 \hsy \pi}{2}} \ 
\log 
\big |
f \big(R e^{\sqrt{-1} \hsx \theta}\big)
\big | 
\hsx 
\cos \theta 
\ \td \theta
\\[15pt]
&=\ 
\frac{1}{\pi \hsy R} \ 
\int\limits_{\frac{3 \hsy\pi}{2}}^{2 \hsy \pi} \ 
\log 
\big |
f \big(R e^{\sqrt{-1} \hsx (\theta - 2 \hsy \pi)}\big)
\big |
\hsx 
\cos (\theta - 2 \hsy \pi)
\ \td \theta
\\[15pt]
&
\hspace{2cm}
\ + \ 
\frac{1}{\pi \hsy R} \ 
\int\limits_0^{\frac{3 \hsy \pi}{2}} \ 
\log 
\big |
f \big(R e^{\sqrt{-1} \hsx \theta}\big)
\big | 
\hsx 
\cos \theta 
\ \td \theta
\\[15pt]
&=\ 
\frac{1}{\pi \hsy R} \ 
\int\limits_0^{2 \hsy \pi} \ 
\log 
\big |
f \big(R e^{\sqrt{-1} \hsx \theta}\big)
\big |
\hsx 
\cos \theta 
\ \td \theta.
\end{align*}
Summary: 
\[
\sum\limits_{r_n \leq r}
\Big(
\frac{1}{r_n} - \frac{r_n}{r^2}
\Big) 
\hsx 
\cos \theta_n 
\hsx + \hsx 
\Reg f^\prime (0)
\ = \ 
\frac{1}{\pi \hsy r} \ 
\int\limits_0^{2 \hsy \pi} \ 
\log 
\big |
f \big(r e^{\sqrt{-1} \hsx \theta}\big)
\big |
\hsx 
\cos \theta 
\ \td \theta.
\]
\\[-.5cm]

\un{Step 4:}\ 

\[
\int\limits_0^r \ 
\frac{n (t)}{t} 
\ \td t
\ = \ 
\frac{1}{2 \hsy \pi} \ 
\int\limits_0^{2 \hsy \pi} \ 
\log \big | f \big(r \hsy e^{\sqrt{-1} \hsx \theta}\big)\big | 
\td \theta
\]
\qquad 
$\implies$

\[
\frac{1}{r} \ 
\int\limits_0^r \ 
\frac{n (t)}{t} 
\ \td t
\ = \ 
\frac{1}{2 \hsy \pi \hsy r} \ 
\int\limits_0^{2 \hsy \pi} \ 
\log \big | f \big(r \hsy e^{\sqrt{-1} \hsx \theta}\big)\big | 
\td \theta
\]
\qquad 
$\implies$

\[
\lim\limits_{r \ra \infty} \ 
\frac{1}{r} \ 
\int\limits_0^r \ 
\frac{n (t)}{t} 
\ \td t
\ = \ 
D 
\ = \ 
\lim\limits_{r \ra \infty} \ 
\frac{1}{2 \hsy \pi \hsy r} \ 
\int\limits_0^{2 \hsy \pi} \ 
\log \big | f \big(r \hsy e^{\sqrt{-1} \hsx \theta}\big)\big | 
\td \theta.
\]
[Given $\varepsilon > 0$, choose $t_0$: 
\[
t > t_0 
\implies 
D - \varepsilon 
< 
\frac{n (t)} {t} 
< 
D + \varepsilon.
\]
Write
\[
\frac{1}{r} \ 
\int\limits_0^r \ 
\frac{n (t)}{t} 
\ \td t
\ = \ 
\frac{1}{r} \ 
\int\limits_0^{t_0} \ 
\frac{n (t)}{t} 
\ \td t
\hsx + \hsx
\frac{1}{r} \ 
\int\limits_{t_0}^r \ 
\frac{n (t)}{t} 
\ \td t
\qquad (r > t_0).
\]
Then 
\[
\frac{(r - t_0) \hsy (D - \varepsilon)}{r}
\ < \ 
\frac{1}{r} \ 
\int\limits_{t_0}^r \ 
\frac{n (t)}{t} 
\ \td t
\ < \ 
\frac{(r - t_0) \hsy (D + \varepsilon)}{r}
\]

\qquad 
$\implies$ $(r \ra \infty)$

\[
D - \varepsilon
\ \leq \ 
\lim\limits_{r \ra \infty} \ 
\frac{1}{r} \
\int\limits_{t_0}^r \ 
\frac{n (t)}{t} 
\ \td t
\ \leq \ 
D + \varepsilon.]
\]

\un{Step 5:}\quad
We have

\allowdisplaybreaks
\begin{align*}
h_f \big(e^{\sqrt{-1} \hsx \theta}\big)\
&=\ 
a \hsy \abs{\sin \theta} + b \sin \theta
\\[15pt]
&=\ 
\frac{h_f \big(\sqrt{-1}\big) + h_f \big(-\sqrt{-1}\big)}{2}
\hsx 
\abs{\sin \theta} 
\ + \ 
\frac{h_f \big(\sqrt{-1}\big) - h_f \big(-\sqrt{-1}\big)}{2}
\hsx 
\sin \theta 
\\[15pt]
&=\ 
\frac{\pi \hsy D}{2} 
\hsx 
\abs{\sin \theta} + b \sin \theta
\end{align*}
\qquad
$\implies$
\allowdisplaybreaks
\begin{align*}
\frac{1}{2 \hsy \pi} \ 
\int\limits_0^{2 \hsy \pi} \ 
h_f \big(e^{\sqrt{-1} \hsx \theta}\big)
\ \td \theta \ 
&=\ 
\frac{1}{2 \hsy \pi} \ 
\int\limits_0^{2 \hsy \pi} \ 
\frac{\pi \hsy D}{2} 
\abs{\sin \theta}
\ \td \theta \
\\[15pt]
&=\ 
\frac{D}{4}  \ 
\int\limits_0^{2 \hsy \pi} \ 
\abs{\sin \theta}
\ \td \theta 
\\[15pt]
&=\ 
D.
\end{align*}
\\[-.5cm]

\un{Step 6:}\quad 
Given $\varepsilon > 0$, choose $r_0$:
\\[-.25cm]
 
$r > r_0$
$\implies$
\[
- 2 \varepsilon 
\ < \ 
\int\limits_0^{2 \hsy \pi} \ 
\big(
h_f \big(e^{\sqrt{-1} \hsx \theta}\big)
\hsx + \hsx
\varepsilon
\hsx - \hsx 
\frac{1}{r} 
\hsx
\log \big | f \big(r e^{\sqrt{-1} \hsx \theta}\big)\big |
\big)
\ \td \theta
\ < \ 
2 \varepsilon .
\]
But for $r_0 \gg 0$, 
\[
\frac{1}{r} 
\hsx
\log \big |f \big(r e^{\sqrt{-1} \hsx \theta}\big)\big |
\ < \ 
h_f \big(e^{\sqrt{-1} \hsx \theta}\big) + \varepsilon
\]
uniformly in $\theta$ (inspect the first part of the proof of 19.7), thus

\[
- 2 \varepsilon 
\ < \ 
\int\limits_0^{2 \hsy \pi} \ 
\Big(
h_f \big(e^{\sqrt{-1} \hsx \theta}\big)
\hsx + \hsx
\varepsilon
\hsx - \hsx 
\frac{1}{r} 
\hsx
\log \big | f \big(r e^{\sqrt{-1} \hsx \theta}\big) \big |
\Big)
\cos \theta
\ \td \theta
\ < \ 
2 \varepsilon
\]
and so
\[
\lim\limits_{r \ra \infty} \
\frac{1}{r} \ 
\int\limits_0^{2 \hsy \pi} \ 
\log \big | f \big(re^{\sqrt{-1} \hsx \theta}\big) \big |
\hsx
\cos \theta
\ \td \theta
\ = \ 
\int\limits_0^{2 \hsy \pi} \ 
h_f \big(e^{\sqrt{-1} \hsx \theta}\big)
\cos \theta
\ \td \theta.
\]
\\[-1cm]

\un{Step 7:}\ 
\allowdisplaybreaks
\begin{align*}
\qquad \text{\textbullet} \quad
\int\limits_0^\pi \ 
\abs{\sin \theta} \cos \theta
\ \td \theta \ 
&=\ 
\int\limits_0^\pi \ 
\sin \theta
\hsx
\cos \theta
\ \td \theta \ 
\hspace{5.75cm}
\\[15pt]
&=\ 
\frac{1}{2} \ 
\int\limits_0^\pi \ 
\sin 2 \hsy \theta
\ \td \theta \ 
\\[15pt]
&=\ 
\frac{1}{2} 
\ - \ 
\frac{\cos 2 \theta}{2}
\bigg|_0^\pi 
\\[15pt]
&=\ 
\frac{1}{4} 
\hsx 
(- \cos 2 \pi + \cos 0)
\\[15pt]
&=\ 
0.
\end{align*}
\allowdisplaybreaks
\begin{align*}
\qquad \text{\textbullet} \quad
\int\limits_\pi^{2 \hsy \pi} \ 
\sin \theta
\hsx
\cos \theta
\ \td \theta \ 
&=\ 
\frac{1}{2} \
\int\limits_\pi^{2 \hsy \pi} \
\sin 2 \theta
\ \td \theta 
\hspace{5.75cm}
\\[15pt]
&=\ 
\frac{1}{2} 
\ - \ 
\frac{\cos 2 \theta}{2}
\bigg|_\pi^{2 \hsy \pi}
\\[15pt]
&=\ 
\frac{1}{4} 
\hsx 
(- \cos 4 \pi + \cos 2 \pi)
\\[15pt]
&=\ 
0.
\end{align*}
Consequently, 
\[
\frac{1}{\pi} \ 
\int\limits_0^{2 \hsy \pi} \
h_f \big(e^{\sqrt{-1} \hsx \theta}\big) 
\cos \theta
\ \td \theta
\ = \ 
0,
\]
which implies that
\[
\lim\limits_{r \ra \infty} \
\frac{1}{\pi \hsy r} \ 
\int\limits_0^{2 \hsy \pi} \
\log \big | f \big(r e^{\sqrt{-1} \hsx \theta}\big)\big |
\cos \theta 
\ \td \theta
\ = \ 
0.
\]
\\[-.5cm]

Summary: 

\[
\lim\limits_{r \ra \infty} \
\sum\limits_{r _n \hsy \leq \hsy r} \
\Big(
\frac{1}{r_n} - \frac{r_n}{r^2}
\Big)
\hsx 
\cos \theta_n 
\ = \ 
- 
\Reg f^\prime (0).
\]
\\[-1cm]

\un{Step 8:}\ 
Let $r$ take the values $m / D$, where $m$ is an integer $-$then
\[
\abs
{
m - n \Big(\frac{m}{D}\Big)
}
\ = \ 
o (m) 
\qquad (m \ra \infty)
\]
\qquad \qquad
$\implies$

\[
\lim\limits_{m \ra \infty} \
\sum\limits_{n = 1}^m \ 
\frac{\cos \theta_n}{r_n}
\hsx
\Big(
1 - \frac{r_n^2 \hsy D^2}{m^2}
\Big)
\ = \ 
- 
\Reg f^\prime (0).
\]
\\[-.75cm]

\un{Step 9:}\ 
Let
\[
\gamma_m 
\ = \ 
\sum\limits_{n = 1}^m \
\frac{\cos \theta_n}{r_n}
\hsx
\Big(
1 - \frac{r_n^2 \hsy D^2}{m^2}
\Big).
\]
Then 
\allowdisplaybreaks
\begin{align*}
(m+1)^2 \hsy \gamma_{m+1} - m^2 \hsy \gamma_m \ 
&=\ 
(2m + 1) \ 
\sum\limits_{n = 1}^m \
\frac{\cos \theta_n}{r_n}
\hsx + \hsx 
\frac{\cos \theta_{m+1}}{r_{m+1}}
\hsx
((m+1)^2 - D^2 \hsy r_{m+1}^2).
\end{align*}

[Starting from the LHS, 
\allowdisplaybreaks
\begin{align*}
(m+1)^2 \hsy 
&
\gamma_{m+1} 
- m^2 \hsy \gamma_m \
\\[15pt]
&=\ 
\sum\limits_{n = 1}^{m+1} \  
\frac{\cos \theta_n}{r_n}
\hsx
(m^2 \hsx + \hsx 2m + 1 \hsx - \hsx D^2 r_n^2)
\ - \ 
\sum\limits_{n = 1}^m \
\frac{\cos \theta_n}{r_n}
\hsx 
(m^2 \hsx \hsx - \hsx D^2 r_n^2)
\\[15pt]
&=\ 
\sum\limits_{n = 1}^m \
\frac{\cos \theta_n}{r_n}
\hsx 
m^2 
\hsx - \hsx 
\sum\limits_{n = 1}^m \
\frac{\cos \theta_n}{r_n}
\hsx 
m^2 
\hsx + \hsx 
\frac{\cos \theta_{m+1} }{r_{m+1} }
\hsx 
m^2 
\\[15pt]
&
\hspace{1.5cm}
\hsx + \hsx 
\sum\limits_{n = 1}^{m+1} \ 
\frac{\cos \theta_n}{r_n}
(2 m+1 \hsx - \hsx D^2 r_n^2)
\hsx + \hsx 
\sum\limits_{n = 1}^m \
\frac{\cos \theta_n}{r_n}
\hsx 
D^2 r_n^2
\\[15pt]
&=\ 
(2m + 1) \ 
\sum\limits_{n = 1}^m \
\frac{\cos \theta_n}{r_n}
\hsx + \hsx 
\frac{\cos \theta_{m+1}}{r_{m+1}}
(2m + 1)
\hsx + \hsx 
\frac{\cos \theta_{m+1}}{r_{m+1}}
\hsx
m^2
\\[15pt]
&
\hspace{1.5cm}
\hsx - \hsx
\sum\limits_{n = 1}^m \
\frac{\cos \theta_n}{r_n}
\hsx D^2 \hsy r_n^2
\hsx + \hsx 
\sum\limits_{n = 1}^m \
\frac{\cos \theta_n}{r_n}
\hsx D^2 \hsy r_n^2
\hsx - \hsx
\frac{\cos \theta_{m+1}}{r_{m+1}}
\hsx D^2 \hsy r_n^2
\\[15pt]
&=\ 
(2m + 1) \ 
\sum\limits_{n = 1}^m \
\frac{\cos \theta_n}{r_n}
\hsx + \hsx 
\frac{\cos \theta_{m+1} }{r_{m+1} }
\hsx 
(m^2 + 2m+ 1 \hsx - \hsx D^2 \hsy r_n^2).]
\end{align*}
\\[-.5cm]

\un{Step 10:}\quad
Write
\[
\sum\limits_{n = 1}^m \
\frac{\cos \theta_n}{r_n}
\ = \ 
\frac{(m + 1)^2 \gamma_{m + 1} - m^2 \gamma_m}{2 m + 1}
\ + \ 
A_m, 
\]
where

\[
A_m 
\ = \ 
- \ 
\frac
{
\ds
\frac{\cos \theta_{m+1}}{r_{m+1}} 
\hsx 
((m + 1)^2 - D^2 \hsy r_{m+1}^2)
}
{
2 m + 1
}
\hsx.
\]
Claim: 

\[
\lim\limits_{m \ra \infty} \ 
A_m 
\ = \ 
0.
\]

[Take absolute values: 
\allowdisplaybreaks
\begin{align*}
\abs{A_m} \ 
&=\ 
\abs
{
\frac{\cos \theta_{m+1}}{r_{m+1}} 
\ \cdot \ 
\frac{1}{2m + 1}
\ \cdot \ 
((m + 1)^2 - D^2 \hsy r_{m+1}^2)
}
\\[15pt]
&\leq 
\frac{1}{r_{m+1}} \ 
\abs
{
\frac{1}{2m+1}
\
(m^2 + 2 m + 1 - D^2 \hsy r_{m+1}^2)
}
\\[15pt]
&=\ 
\abs
{
\frac{m^2}{2 m + 1}
\hsx
\frac{1}{r_{m+1}}
\ + \ 
\frac{1}{r_{m+1}}
\ - \ 
\frac{D^2 r_{m+1}}{2 m + 1}
}\hsy .
\end{align*}
\\[-1cm]

\qquad \textbullet \quad

\[
\frac{m^2}{2 m + 1}
\hsx
\frac{1}{r_{m+1}}
\ = \ 
\frac{m^2}{2 m + 1}
\hsx
\frac{1}{m+1}
\frac{m+1}{r_{m+1}}
\ \ra \ 
\frac{D}{2}
\qquad (m \ra \infty).
\]
\\[-.5cm]

\qquad \textbullet \quad
\allowdisplaybreaks
\begin{align*}
\frac{1}{r_{m+1}} \ 
&= \ 
\frac{1}{m+1}
\frac{m+1}{r_{m+1}}
\\[15pt]
&\ra
0 \hsy D 
\qquad (m \ra \infty)
\\[15pt]
&= 
0
.
\end{align*}
\\[-.5cm]

\qquad \textbullet \quad
\allowdisplaybreaks
\begin{align*}
- \frac{D^2 r_{m+1}}{2 m + 1} \ 
&=\ 
- \hsy D^2 \ 
\frac{r_{m+1}}{m+1}
\hsx 
\frac{m+1}{2 m + 1}
\\[15pt]
&\ra
- D^2 
\hsx
\frac{1}{D}
\hsx
\frac{1}{2} 
\qquad (m \ra \infty)
\\[15pt]
&= 
-
\frac{D}{2}.
\end{align*}
\\[-1cm]

\un{Step 11:}\quad
Form
\allowdisplaybreaks
\begin{align*}
\frac{1}{p} \ 
\sum\limits_{m = 1}^p \
&
\Big(
\sum\limits_{n = 1}^m \
\frac{\cos \theta_n}{r_n}
\Big)
\\[15pt]
&=\ 
\frac{1}{p} \ 
\sum\limits_{m = 1}^p \
\Big(
\frac{(m + 1)^2 \hsy \gamma_{m + 1} - m^2 \hsy \gamma_m}{2 m + 1}
\ + \
A_m
\Big)
\\[15pt]
&=\
\frac{1}{p} \ 
\Big(
- \frac{\gamma_1}{3}
\ + \ 
\sum\limits_{m = 2}^p \
\frac{2 m^2}{4 m^2 - 1} 
\hsx 
\gamma_m 
\ + \ 
\frac{(p+1)^2}{2 p + 1}
\hsx 
\gamma_{p + 1}
\ + \ 
\sum\limits_{m = 1}^p \
A_m
\Big)
\\[15pt]
&=\
\frac{1}{p} \ 
\Big(
- \gamma_1
\ + \ 
\sum\limits_{m = 1}^p \
\frac{2 m^2}{4 m^2 - 1} 
\hsx 
\gamma_m 
\ + \ 
\frac{(p+1)^2}{2 p + 1}
\hsx 
\gamma_{p + 1}
\ + \ 
\sum\limits_{m = 1}^p \
A_m
\Big).
\end{align*}
\\[-.5cm]

\un{Step 12:}\ 
The series
\[
\sum\limits_{n = 1}^\infty \
\frac{\cos \theta_n}{r_n}
\]
is $(C, 1)$ summable to $-\Reg f^\prime (0)$, hence the series 
\[
\sum\limits_{n = 1}^\infty \
\frac{\cos \theta_n}{r_n}
\]
is convergent to $- \Reg f^\prime (0)$ (cf. 27.35).
\\[-.5cm]

[Let $p \ra \infty$ in the expression above to see what happens.  
First, 
$
\ds
- \frac{\gamma_1}{p}
\ra 0$ 
$(p \ra \infty)$.  
Second, 

\[
\begin{cases}
\hspace{0.5cm}
\gamma_m 
\hspace{0.75cm}
\ra 
\ 
- \Reg f^\prime (0)
\hspace{0.5cm} (m \ra \infty)
\\[18pt]
\  \ds
\frac{2 \hsy m^2}{4 \hsy m^2 - 1}
\ra 
\hspace{0.5cm}
\frac{1}{2}
\hspace{1.65cm} (m \ra \infty)
\end{cases}
\]

$\implies$
\[
\frac{1}{p}\ 
\sum\limits_{m = 1}^p \ 
\frac{2 \hsy m^2}{4 \hsy m^2 - 1}
\hsx 
\gamma_m 
\ra 
- \frac{1}{2}
\hsx
\Reg f^\prime (0)
\quad (p \ra \infty)
\qquad (\tcf. \ 27.34).
\]
Third, 

\[
\frac{1}{p}\ 
\frac{(p + 1)^2}{2 p + 1}
\ 
\gamma_{p + 1} 
\ra 
- \frac{1}{2}
\hsx
\Reg f^\prime (0)
\quad (p \ra \infty).
\]
Fourth, 

\[
\frac{1}{p}\ 
\sum\limits_{m = 1}^p \ 
A_m 
\ra 0
\quad (p \ra \infty)
\qquad (\tcf. \ 27.33).]
\]
\\[-.75cm]

This completes the proof of 27.31 which, as a bonus, serves to establish that 
\[
\sum\limits_{n = 1}^\infty \
\frac{\cos \theta_n}{r_n}
\ = \ 
-
\Reg f^\prime (0)
\qquad (f (0) = 1).
\]
On the other hand, the series

\[
\sum\limits_{n = 1}^\infty \
\frac{\sin \theta_n}{r_n}
\]
is absolutely convergent (cf. 27.10), thus is convergent, the only new wrinkle being that

\[
\frac{1}{\pi} \ 
\int\limits_0^{2 \hsy \pi} \ 
h_f (e^{\sqrt{-1} \hsx \theta})
\hsx 
\sin \theta
\ \td \theta
\ = \ 
\frac{1}{\pi} \ 
\int\limits_0^{2 \hsy \pi} \ 
(a \hsy \abs{\sin \theta} + b \sin \theta) 
\hsx 
\sin \theta
\ \td \theta
\]
is equal to 

\[
b 
\ = \ 
\frac{h_f (\sqrt{-1}) - h_f (-\sqrt{-1})}{2}
\ \equiv \ 
b_f
\]
and this might not vanish (cf. 27.25).  
The upshot, therefore, is that 

\[
\sum\limits_{n = 1}^\infty \
\frac{\sin \theta_n}{r_n}
\ = \ 
\Img f^\prime (0) + b_f 
\qquad (f (0) = 1).
\]
\\[-1.cm]

\begin{x}{\small\bf SCHOLIUM} \ 
If $f (0) = 1$ and $b_f = 0$, then 
\allowdisplaybreaks
\begin{align*}
\sum\limits_{n = 1}^\infty \
\frac{1}{z_n}\ 
&= \ 
\sum\limits_{n = 1}^\infty \
\frac{\cos \theta_n}{r_n}
\ - \ 
\sqrt{-1} \ 
\sum\limits_{n = 1}^\infty \
\frac{\sin \theta_n}{r_n}
\\[15pt]
&= \ 
\Reg f^\prime (0) - \sqrt{-1} \hsx f^\prime (0)
\\[15pt]
&= \ 
- f^\prime (0).
\end{align*}

[Note: \ 
When $f (0) \neq 1$ (but $f (0) \neq 0$), the formula becomes 
\[
\sum\limits_{n = 1}^\infty \
\frac{1}{z_n}
\ = \ 
-\hsx
\frac{f^\prime (0)}{f (0)}.]
\]
\\[-1.cm]
\end{x}

\begin{x}{\small\bf REMARK} \ 
Write
\[
f (z) 
\ = \ 
f (0) \hsx e^{c \hsy z}\
\prod\limits_{n = 1}^\infty \ 
\Big(
1 - \frac{z}{z_n}
\Big)
\hsx 
e^{z / z_n}.
\]
Then
\[
c 
\ = \ 
-\hsx
\frac{f^\prime (0)}{f (0)}
\]
and

\[
f (z) 
\ = \ 
f (0) \ 
\lim\limits_{R \ra \infty} \ 
\prod\limits_{\abs{z_n} \hsy < \hsy R} \ 
\Big(
1 - \frac{z}{z_n}
\Big),
\]
the convergence of the product being conditional.
\\[-.25cm]
\end{x}

\begin{x}{\small\bf EXAMPLE} \ 
Take

\[
f (z) 
\ = \ 
\frac
{(e^{\sqrt{-1} \hsx z} - 1) \hsx (e^{-\sqrt{-1} \hsx z} + \sqrt{-1}\hsx )}
{\sqrt{-1} \hsx z}.
\]
Then
\[
\begin{cases}
\ \ds
f (0) 
\ = \ 
\sqrt{-1} \hsx + 1
\\[18pt]
\ \ds
f^\prime (0) = \frac{(\sqrt{-1} \hsx - 1)}{2}
\hsy \sqrt{-1}
\end{cases}
\implies \ 
\frac{f^\prime (0)}{f (0)}
\ = \ 
-\frac{1}{2}
\]
\\[-.25cm]

\noindent
and the theory predicts that 

\[
\sum\limits_{n = 1}^\infty \ 
\frac{\cos \theta_n}{r_n}
\ = \ 
\frac{1}{2}.
\]
To establish this, note that the zeros of $f (z)$ are at
\[
\pm 2 \hsy \pi, \hsx
\pm 4 \hsy \pi, \hsx
\ldots
\]
and at 
\[
\frac{\pi}{2}, \hsx
-\frac{3 \hsy \pi}{2}, \hsx
\frac{5 \hsy\pi}{2}, \hsx
-\frac{7 \hsy\pi}{2}, \hsx
\ldots \hsx .
\] 
Those of the first kind make no contribution 
(since the corresponding terms of the series cancel in pairs) 
but there is a contribution from those of the second kind, viz.

\[
\frac{2}{\pi} \hsx
\Big(
1 - \frac{1}{3} + \frac{1}{5} - \frac{1}{7} + \cdots
\Big)
\ = \ 
\frac{1}{2}.
\]
\\[-1cm]


[Note: \ 
As regards 

\[
\sum\limits_{n = 1}^\infty \
\frac{\sin \theta_n}{r_n}, 
\]
it is clear that $\sin \theta_n = 0$ $\forall \ n$.  
To see that here $b_f = 0$, work on $[-1, 1]$ and let
\[
\phi (t)\ = \ 
\begin{cases}
\ 
1 
\hspace{1.5cm} (-1 \leq t \leq 0)
\\[4pt]
\ 
\sqrt{-1}
\hspace{0.8cm} (0 < t \leq 1)
\end{cases}
.
\]
Then
\[
f (z) 
\ = \ 
\int\limits_{-1}^1 \ 
\phi (t) 
\hsx 
e^{\sqrt{-1} \hsx z \hsy t}
\ \td t,
\]
hence
\[
\begin{cases}
\ \ \ 
1 = h_f (- \sqrt{-1} \hsx )
\\[11pt]
\ 
- 1 = - h_f ( \sqrt{-1} \hsx )
\end{cases}
\qquad (\tcf. \ 27.26)
\]
\\[-.75cm]

\qquad \qquad
$\implies$
\[
b_f 
\ = \ 
\frac{1-1}{2} 
\ = \ 
0.]
\]
\\[-.25cm]
\end{x}

\chapter{
$\boldsymbol{\S}$\textbf{28}.\quad  ZERO THEORY: \ PALEY-WIENER FUNCTIONS}
\setlength\parindent{2em}
\setcounter{theoremn}{0}
\renewcommand{\thepage}{\S28-\arabic{page}}

\qquad 
Recall that $\PW (A)$ is the subset of $\sE_0 (A)$ consisting of those $f$ such that 
$\restr{f}{\R} \in \Lp^2 (-\infty, \infty)$ (cf. 22.1).
\\[-.25cm]

\begin{x}{\small\bf EXAMPLE} \ 
Take $A = \pi$ $-$then
\[
\Big(
1 - \frac{\sin \pi \hsy z}{\pi \hsy z}
\Big) 
/ (\pi \hsy z)^2
\in \PW (\pi)
\]
has no real zeros.
\\[-.25cm]
\end{x}

\begin{x}{\small\bf EXAMPLE} \ 
Take $A = \pi$ $-$then
\[
\Big(
1 - \frac{\sin \pi \hsy z}{\pi \hsy z}
\Big) 
/ \pi \hsy z
\in \PW (\pi)
\]
has exactly one real zero.
\\[-.25cm]
\end{x}

\begin{x}{\small\bf EXAMPLE} \ 
Take $A = 1$ $-$then
\[
\frac{e^{\sqrt{-1} \hsx z} - 1}{z} \in \PW (1)
\]
has infinitely many real zeros.
\\[-.25cm]
\end{x}

\begin{x}{\small\bf RAPPEL} \ 
The elements $f \in \PW (A)$ have the form
\[
f (z) 
\ = \ 
\frac{1}{\sqrt{2 \hsy \pi}} \ 
\int\limits_{-A}^A \ 
\phi (t) 
\hsx 
e^{\sqrt{-1} \hsx z \hsy t}
\ \td t
\qquad (0 < A < \infty)
\]
for some $\phi \in \Lp^2[-A, A]$ (cf. 22.7).
\\[-.5cm]

[Note: \ 
The prescription 
\[
\phi (t) 
\ = \ 
\lim\limits_{R \ra \infty} \ 
\frac{1}{\sqrt{2 \hsy \pi}} \ 
\int\limits_{-R}^R \ 
f (x)
\hsx
e^{-\sqrt{-1} \hsx t \hsy x}
\ \td x
\quad (\Lp^2)
\]
computes $\phi$ in terms of $f$.]
\\[-.25cm]
\end{x}


\begin{x}{\small\bf DEFINITION} \ 
Suppose that $f \in \PW (A)$ $-$then $f$ is called a \un{band-pass} function 
if there exists an interval $[-B,B]$ $(0 < B < A)$ in which $\phi = 0$ almost everywhere.
\\[-.25cm]
\end{x}

\begin{x}{\small\bf LEMMA} \ 
If $f \not\equiv 0$ is a real integrable band-pass function, then $f$ has at least one real zero.
\\[-.5cm]

PROOF \ 
Take  $\phi \equiv 0$ in $[-B,B]$, hence
$
\ds
\int\limits_{-\infty}^\infty \ 
f (x) 
\ \td x 
= 0
$, 
so $f$ must change sign somewhere in $\R$.
\\[-.25cm]
\end{x}

More is true.
\\[-.25cm]

\begin{x}{\small\bf THEOREM} \ 
If $f \not\equiv 0$ is a real band-pass function, then $f$ has infinitely many real zeros.
\\[-.5cm]

[The point of departure is the following observation: \ 
$\forall \ g \in \PW (B)$ 
$(\subset \PW (A))$, 
\[
\langle g, f \rangle
\ = \ 
\langle \psi, \phi\rangle,
\]
where 
\[
g (z) 
\ = \ 
\frac{1}{\sqrt{2 \hsy \pi}} \ 
\int\limits_{-B}^B \ 
\psi (t) 
\hsx
e^{\sqrt{-1} \hsx z \hsy t}
\ \td t.
\]
With this in mind, assume that $f$ has but finitely many real zeros.  
One then arrives at a contradiction by exhibiting a real $g \in \PW (B)$ such that 
$\langle g, f \rangle \neq 0$.
\\[-.25cm]

\qquad \textbullet \quad
$f (x)$ is of constant sign: \ 
Take 
\[
g (z) 
\ = \ 
\Big(
\frac{1}{z}
\hsx
\sin 
\Big(
\frac{B}{2} \hsy z
\Big)
\Big)^2.
\]
\\[-1.25cm]

\qquad \textbullet \quad
$f (x)$ is not of constant sign, thus has zeros of odd order, say $x_1, \ldots, x_n$ 
(these are the zeros at which $f$ changes sign).  
Now construct a real $g \in \PW (B)$ whose real zeros are precisely the 
$x_k$ $(k = 1, \ldots, n)$, each $x_k$ being of 
order 1 (per $g$).  
Therefore 
$g (x) \hsy f (x) \geq 0$ 
$\forall \ x$ 
or 
$g (x) \hsy f (x) \leq 0$ 
$\forall \ x$, 
so 
$\langle g, f \rangle \neq 0$.]
\\[-.25cm]
\end{x}

\begin{x}{\small\bf RAPPEL} \ 
Let $f$ be a continuously differentiable complex valued function on $[a,b]$. 
Assume: \ 
$f (a) = f (b) = 0$ $-$then 
\[
\int\limits_a^b \ 
\abs{f(x)}^2 
\ \td x \
\ \leq \ 
\Big(
\frac{b - a}{\pi}
\Big)^2 \ 
\int\limits_a^b \ 
\abs{f^\prime(x)}^2 
\ \td x \
\]
with equality iff 
\[
f (x) 
\ = \ 
C 
\hsx 
\sin 
\Big(
\pi
\hsx
\frac{x - a}{b - a}
\Big).
\]

[This is known as 
\un{Wirtinger's inequality}
\footnote[2]{\vspace{.11 cm}
G. Folland, 
\textit{Real Analysis}, Wiley-Interscience, 1984, p. 247.}
.]
\\[-.25cm]
\end{x}

\begin{x}{\small\bf THEOREM} \ 
Let $f \in \PW (A)$ be nonzero $-$then $\abs{f} > 0$ on at least one open interval 
of the real axis of length 
$
\ds
> 
\frac{\pi}{A}
$.
\\[-.5cm]

PROOF \ 
One need only consider the situation when $f$ has infinitely many real zeros. 
So suppose that 
$a < b$ are two consecutive zeros of $f$ and that, moreover, 
$
\ds 
b - a
\leq 
\frac{\pi}{A}
$.  
Since $f$ is not a sine function on any interval, 
\begin{align*}
\int\limits_a^b \ 
\abs{f(x)}^2 
\ \td x \
&<\ 
\Big(
\frac{b - a}{\pi}
\Big)^2 \ 
\int\limits_a^b \ 
\abs{f^\prime(x)}^2 
\ \td x \
\\[15pt]
&\leq 
\Big(
\frac{1}{A}
\Big)^2 
\ 
\int\limits_a^b \ 
\abs{f^\prime(x)}^2 
\ \td x ,
\end{align*}
which implies by addition that

\[
\norm{f}_2 
\ < \ 
\frac{1}{A}
\hsx
\norm{f^\prime}_2 .
\]
But

\[
\norm{f^\prime}_2 
\ \leq \ 
\norm{f}_2 
\hsx
T (f)
\qquad (\tcf. \ 17.31).
\]
Therefore

\[
\norm{f}_2 
\ < \ 
\frac{T (f)}{A}
\hsx
\norm{f}_2 
\]
\qquad 
$\implies$
\[
A 
\ < \
T (f),
\]
a contradiction
\\[-.25cm]
\end{x}

\begin{x}{\small\bf EXAMPLE} \ 
The Paley-Wiener function 

\[
\frac{\sin A \hsy x}{A \hsy x}
\]
has just one zero free open interval of length 
$
\ds
> 
\frac{\pi}{A}
$, 
namely
$
\ds
\Big]
- 
\frac{\pi}{A}, 
\frac{\pi}{A}
\Big[\hsx.
$
\\[-.25cm]
\end{x}


\chapter{
$\boldsymbol{\S}$\textbf{29}.\quad  INTERMEZZO}
\setlength\parindent{2em}
\setcounter{theoremn}{0}
\renewcommand{\thepage}{\S29-\arabic{page}}

\qquad
Given $\phi \in \Lp^1 [a,b]$, let
\[
f (z) 
\ = \ 
\int\limits_a^b \ 
\phi (t) 
\hsx 
e^{\sqrt{-1} \hsx z \hsy t}
\ \td t.
\]
Then $f (z)$ is a Bernoulli function and subject to suitable restrictions on $\phi$, 
the overall program is to study the position of the zeros of $f (z)$.
\\[-.25cm]

\qquad
{\small\bf \un{N.B.}} \ 
It is sometimes convenient to ``normalize'' the interval and take 
$[a,b] = [0,1]$ 
or 
$[a,b] = [-1,1]$.
\\[-.5cm]

\qquad \textbullet \quad
Thus
\begin{align*}
\int\limits_a^b \ 
\phi (t) 
&
\hsx 
e^{\sqrt{-1} \hsx z \hsy t}
\ \td t
\ = \ 
(b - a) \hsy e^{\sqrt{-1} \hsx a \hsy z}
\int\limits_0^1 \ 
\phi (a + (b - a) \hsy t)
\hsx
e^{\sqrt{-1} \hsx (b - a) \hsy z \hsy t}
\ \td t.
\end{align*}
\\[-1.25cm]

\qquad \textbullet \quad
Thus
\\[-1cm]
\begin{align*}
\int\limits_a^b \ 
\phi (t) 
&
\hsx 
e^{\sqrt{-1} \hsx z \hsy t}
\ \td t
\\[-18pt]
&
\hspace{-0.5cm}
= \  
\frac{1}{2} \hsx (b - a) \hsy e^{\frac{1}{2} (a + b)\sqrt{-1} \hsx z} \ 
\int\limits_{-1}^1 \ 
\phi 
\Big(
\frac{1}{2}
(b + a) 
+ 
\frac{1}{2}
(b - a)
\hsy t
\Big)
\hsx
e^{\frac{1}{2} (b - a) \sqrt{-1} \hsx z \hsy t}
\ \td t.
\end{align*}
\\[-.5cm]
The theory developed in \S 27 is applicable under the following conditions.
\\[-.25cm]

\qquad \textbullet \quad
Assume: \ 
$f (0) \neq 0$.  
\\[-.5cm]

[Note: \ 
Nothing of substance is lost in so doing.  
For if $f (0) = 0$, then 
\[
\frac{f (z)}{z} 
\ = \ 
- \sqrt{-1} \ 
\int\limits_a^b \ 
\psi (t) 
e^{\sqrt{-1} \hsx z \hsy t}
\ \td t, 
\]
where
\[
\psi (t) 
\ = \ 
\int\limits_a^t \ 
f (s)
\ \td s.]
\]

\qquad \textbullet \quad
Assume: \ 
There is no $\alpha > a$ such that 
\[
\int\limits_a^\alpha \ 
\abs{\phi (t) }
\ \td t
\ = \ 
0
\]
and there is no $\beta < b$ such that 
\[
\int\limits_\beta^b \ 
\abs{\phi (t) }
\ \td t
\ = \ 
0.
\]

[Note: \ 
Accordingly,
\[
\begin{cases}
\ 
a 
\ = \ 
- 
h_f \big(\sqrt{-1} \hsx\big)
\\[4pt]
\ 
b 
\ = \ 
h_f \big(-\sqrt{-1} \hsx\big)
\end{cases}
,
\]
and
\[
T (f) 
\ = \ 
\max\big(h_f \big(\sqrt{-1} \hsx\big), h_f \big(-\sqrt{-1} \hsx\big)\big).]
\]

Therefore in review: 
\\

1. \quad
$
\ds
\lim\limits_{r \ra \infty} \ 
\frac{n (r)}{r}
\ = \ 
\frac{b - a}{\pi} 
\ \equiv \ 
D 
\ > \ 
0.
$
\\[.25cm]

2. \quad
$
\ds
\sum\limits_{n = 1}^\infty \
\frac{\sin \theta_n}{r_n}
\ 
$ 
is absolutely convergent and has sum
\[
\Img \frac{f^\prime (0)}{f (0)} 
- 
\frac{(a + b)}{2}.
\]

3. \quad
$
\ds
\sum\limits_{n = 1}^\infty \
\frac{\cos \theta_n}{r_n}
\ 
$ 
is conditionally convergent and has sum
\[
- \Reg \frac{f^\prime (0)}{f (0)}.
\]
\\[-1cm]

\qquad
{\small\bf \un{N.B.}} \ 
Matters simplify if $a = -A$, $b = A$.
\\[-.25cm]

\begin{x}{\small\bf EXAMPLE} \ 
The zeros of $f (z)$ which lie on the imaginary axis constitute a ``thin'' set 
(if there are any at all) (cf. 27.11).  
Still, their number may be infinite.
\\[-.5cm]

[Working on \ $[0,1]$, \ choose constants 
$
\ 
\ds
0 < \mu < \frac{1}{2}
$, 
\ 
$\nu > 2$, and put 
$\alpha = \nu/\mu$.
Define $\phi \in \Lp^1 [0,1]$ by letting 

\[
\phi (t) 
\ = \ 
(-\alpha)^k \hsy e^{-\nu^k} 
\quad (\mu^k - \alpha^{-k} < t \leq \mu^k) 
\qquad (k = 1, 2, \ldots)
\]
and taking $\phi (t) = 0$ elsewhere on $[0,1]$.  
Given any positive integer $n$, we have
\allowdisplaybreaks
\begin{align*}
\bigg|
\hsx
\int\limits_0^{\mu^{n+1}} \ 
\phi (t) 
\hsy
e^{-\alpha^n \hsy t}
\ \td t
\hsx
\bigg|
&\leq \
\int\limits_0^{\mu^{n+1}} \ 
\abs{\phi (t) }
\ \td t
\\[11pt]
&= \
\sum\limits_{k = n+1}^\infty \ 
e^{-\nu^k} 
\\[11pt]
&< \
e^{-\nu^{n+1}} 
\ 
\sum\limits_{j = 0}^\infty \ 
e^{-\nu^j} 
\\[11pt]
&= \
e^{-\nu^{n+1}} 
\int\limits_0^1 \ 
\abs{\phi (t) }
\ \td t
\end{align*}
and
\allowdisplaybreaks
\begin{align*}
\bigg|
\hsx
\int\limits_{\mu^{n-1} - \alpha^{-n + 1}}^1\ 
\phi (t) 
\hsy
e^{-\alpha^n \hsy t}
\ \td t
\hsx
\bigg|
\ 
&\leq \
e^{- \alpha^n} 
\hsy
\big(
\mu^{n-1} - \alpha^{-n+1}
\big)
\int\limits_0^1
\abs{\phi (t) }
\ \td t
\\[11pt]
&= \
e^{- \nu^n / \mu + \alpha} \ 
\int\limits_0^1 \ 
\abs{\phi (t) }
\ \td t
\end{align*}
and 
\[
\int\limits_{\mu^n - \alpha^{-n}}^{\mu^n} \ 
\phi (t) 
\hsy
e^{- \alpha^{n}\hsy t} 
\ \td t
\ = \ 
(-1)^n
\hsy
(e - 1)
\hsy
e^{-2 \nu^n}.
\]
Therefore
\allowdisplaybreaks
\begin{align*}
\bigg|
\hsx
e^{2 \nu^n} \ 
\int\limits_0^1 \ 
\phi (t) \hsx e^{- \alpha^n \hsy t}
\ \td t
&
\ - \ 
(e - 1) \hsy (-1)^n
\hsx
\bigg|
\\[15pt]
&<\ 
\big(
e^{\nu^n (2 - \nu)} 
\hsx + \hsx 
e^{\nu^n  (2 - 1/\mu) + \alpha}
\big)
\ 
\int\limits_0^1 \ 
\abs{\phi (t)}
\ \td t.
\end{align*}
So for $n \gg 0$,
\[
\sgn \ 
\int\limits_0^1 \ 
\phi (t) \hsx e^{- \alpha^n \hsy t}
\ \td t
\ = \ 
\sgn (-1)^n,
\]
thus at some $x_0$: 
$-\alpha^{n+1} \leq x_0 \leq - \alpha^n$, 
\[
\int\limits_0^1 \ 
\phi (t) \hsx e^{x_0 \hsy t}
\ \td t
\ = \ 
0
\]
or still, 
\[
f \Big(\frac{x_0}{\sqrt{-1}}\Big)
\ = \ 
0.]
\]
\\[-1.cm]
\end{x}

\begin{x}{\small\bf NOTATION} \ 
Let
\[
F (z) 
\ = \ 
\int\limits_a^b \ 
\phi (t) \hsx e^{z \hsy t}
\ \td t.
\]
Then
\[
f (z) 
\ = \ 
F (\sqrt{-1} \hsx z).
\]
\\[-1.cm]
\end{x}

\begin{x}{\small\bf LEMMA} \ 
Take $[a,b] = [-1, 1]$ $-$then
\[
F \big(r \hsy e^{\sqrt{-1}\hsx \theta}\big) 
\ = \ 
o(e^{r \hsy \abs{\cos \theta}})
\qquad (r \ra \infty)
\]
uniformly with respect to $\theta$.
\\[-.5cm]

PROOF \ 
Assume first that $\theta = 0$ and write
\allowdisplaybreaks
\begin{align*}
\abs{F (r)} \ 
&=\ 
\bigg|
\hsx
\int\limits_{-1}^1 \ 
\phi (t)
\hsy
e^{r \hsy t}
\ \td t
\hsx
\bigg|
\\[15pt]
&=\ 
\bigg|
\hsx
\int\limits_{-1}^{1 - \delta} \ 
\phi (t)
\hsy
e^{r \hsy t}
\ \td t
\ + \ 
\int\limits_{1 - \delta}^1 
\phi (t)
\hsy
e^{r \hsy t}
\ \td t
\hsx
\bigg|
\\[15pt]
&\leq\ 
e^{(1 - \delta) \hsy r}
\ 
\int\limits_{-1}^{1 - \delta} \ 
\abs{\phi (t)} 
\ \td t
\ + \ 
e^r
\ 
\int\limits_{1 - \delta}^1 
\abs{\phi (t)} 
\ \td t.
\end{align*}
Given $\varepsilon > 0$, choose $\delta > 0$: 
\[
\int\limits_{1 - \delta}^1 \ 
\abs{\phi (t)} 
\ \td t
\ < \ 
\frac{\varepsilon}{2}
\]
and then choose $r_0 \gg 0$: 
\[
e^{- \delta \hsy r} 
\ 
\int\limits_{-1}^{1 - \delta} \ 
\abs{\phi (t)} 
\ \td t
\ < \ 
\frac{\varepsilon}{2}
\qquad (r > r_0).
\]
Therefore
\[
\abs{F (r)} 
\ < \ 
\varepsilon \hsy e^r
\qquad (r > r_0).
\]
I.e.: \ 
$F (r) = \txo(e^r)$ $(\cos 0 = 1)$.  
\ 
Next 
\[
F (\sqrt{-1} \hsx x) 
\ = \ 
\int\limits_{-1}^1 \ 
\phi (t) 
\hsy
\cos x t 
\ \td t
\hsx + \hsx 
\sqrt{-1} \ 
\int\limits_{-1}^1 \ 
\phi (t) 
\hsy
\sin x t 
\ \td t
\]
and the two integrals on the right approach 0 as $x \ra \infty$ (Riemann-Lebesgue lemma).  
These facts, in conjunction with Phragm\'en-Lindel\"of, then imply that the function 
$e^{-z} F (z)$ tends uniformly to zero in the sector 
$
\ds
0 \leq \theta \leq \frac{\pi}{2}
$
which gives the result in this range.
And so on \ldots\hsx.
\\[-.25cm]
\end{x}

\begin{x}{\small\bf RAPPEL} \ 
If $\phi$ is absolutely continuous on $[a,b]$, then its derivative $\phi^\prime$ exists almost everywhere.  
Moreover, $\phi^\prime \in \L^1 [a,b]$ and 
\[
\phi (t) 
\ = \ 
\phi(a) 
\ + \ 
\int\limits_a^t \ 
\phi^\prime (s) 
\ \td s
\qquad (a \leq t \leq b).
\]
\\[-1.cm]
\end{x}

\begin{x}{\small\bf THEOREM} \ 
Take $[a,b] = [-1, 1]$ and assume that $\phi$ is absolutely continuous with $\phi (1) = \phi(-1) = 1$ $-$then 
the zeros of $f (z)$ are determined asymptotically by the formula 
\[
z 
\ = \ 
\pm m \hsy \pi + \varepsilon_m,
\]
where $m$ is a positive integer and $\varepsilon_m \ra 0$ $(m \ra \infty)$.
\\[-.5cm]

PROOF \ \
We shall work instead with $F (z)$, thereby shifting the claim to 
$
\ 
\pm \hsx m \hsy \pi \sqrt{-1} \  + \ \varepsilon_m$.  
So $\forall \ z \neq 0$, integrate by parts and write
\[
F (z) 
\ = \ 
\frac{e^z - e^{-z}}{z}
\ - \ 
\frac{1}{z}
\ 
\int\limits_{-1}^1 \
\phi^\prime (t) 
\hsx 
e^{z \hsy t}
\ \td t
\]
or still, 
\[
z \hsy F (z) 
\ = \ 
e^z - e^{-z}
\ - \ 
\int\limits_{-1}^1 \
\phi^\prime (t) 
\hsx
e^{z \hsy t}
\ \td t,
\]
a relation that is valid $\forall \ z$.
Since $\phi^\prime$ is integrable, 29.3 is applicable (replace the $\phi$ there by $\phi^\prime$), hence
\[
\int\limits_{-1}^1 \
\phi^\prime (t) 
\hsx
e^{z \hsy t}
\ \td t
\ = \ 
\txo\big(e^{r \abs{\cos \theta}}\big)
\qquad (r \ra \infty)
\]
uniformly with respect to $\theta$.  
If generically, $\varepsilon_r$ is a function of $r$ and $\theta$ which tends to 0 uniformly in $\theta$ 
as $r \ra \infty$, then at a zero of $F (z)$, 
\[
e^z (1 + \varepsilon_r) 
\ = \ 
e^{-z} (1 + \varepsilon_r) 
\]
\qquad \qquad 
$\implies$
\[
e^{2 \hsy z}
\ = \ 
1 + \varepsilon_r
\]
\qquad \qquad
$\implies$
\[
2 \hsy z 
\ = \ 
\pm 2 \hsy m \hsy \pi \sqrt{-1} \ + \  \varepsilon_m
\]
\qquad \qquad 
$\implies$
\[
z
\ = \ 
\pm m \hsy \pi \sqrt{-1} \ + \  \varepsilon_m.
\]
To reverse this, note that $\sinh z$ has exactly one zero at each point $\pm m \hsy \pi \hsy \sqrt{-1}$.  
Choosing $\delta > 0$ small, surround each of these points by a circle of radius $\delta$, 
thus on the circle 
\[
\abs{\sinh z} 
\ > 
K (\delta) 
\ > 
 0
\]
and
\[
z \hsy F (z)
\ = \ 
\sinh z \hsx (1 + \varepsilon_m),
\]
where $\varepsilon_m > 0$ $(m \ra \infty)$.  
So for large $m$, $z \hsy F (z)$ has the same number of zeros inside the circle as $\sinh z$, i.e., one.
\\[-.25cm]
\end{x}

\begin{x}{\small\bf REMARK} \ 
The supposition that $\phi (1) = \phi(-1) = 1$ is not unduly restrictive at least if 
$\phi (1)$, $\phi(-1)$ are real and positive: \ 
Consider
\[\psi (t) 
\ = \ 
\bigg[
\frac{\phi (-1)}{\phi (1)}
\bigg]^{t/2}
\hsx
\frac{\phi (t) }{\sqrt{\phi (1) \hsy \phi (-1)}}
\]
and define $w$ by the relation
\[
z 
\ = \ 
w + \frac{1}{2} 
\hsx
\log \frac{\phi (-1)}{\phi (1)}.
\]
Then
\allowdisplaybreaks
\begin{align*}
f (z) \ 
&= \ 
\sqrt{\phi (1) \hsy \phi (-1)} \ 
\int\limits_{-1}^1 \ 
\psi (t) 
\hsx 
e^{w \hsy t} 
\ \td t
\\[15pt]
& \equiv \ 
\sqrt{\phi (1) \hsy \phi (-1)} \ 
g (w)
\end{align*}
and $\psi$ is absolutely continuous with $\psi (1) = \psi(-1) = 1$.
\\[-.25cm]
\end{x}

\begin{x}{\small\bf EXAMPLE} \ 
The situation can be different if 
$\phi (-1) = 0$ and $\phi(1) = 0$.  
To see this, let 

\[
\phi (t) \ = \ 
\begin{cases}
\ 
1 - t
\qquad (0 < t \leq 1)
\\[4pt]
\ 
1 + t 
\qquad (-1 \leq t \leq 0)
\end{cases}
.
\]
\\[-.5cm]

\noindent
Then
\[
\phi (t) 
\ = \ 
\int\limits_{-1}^t \ 
\phi^\prime (s) 
\ \td s
\]
is absolutely continuous and 
\[
F (z) 
\ = \ 
\frac{\ds 4 \sinh^2 \Big(\frac{z}{2}\Big)}{z^2}.
\]
However, the zeros are at the points $\pm 2 \hsy m \hsy \pi \sqrt{-1} \hsx$, hence the pattern has changed.
\\[-.25cm]
\end{x}

\begin{x}{\small\bf THEOREM} \ 
Take $[a,b] = [-1, 1]$ and assume that $\phi$ is of bounded variation and continuous at 1 and $-1$  
with $\phi (1) = \phi(-1) = 1$ $-$then 
the zeros of $f (z)$ lie within a horizontal strip $\abs{\Img z} \leq C$.
\\[-.5cm]

PROOF \ 
An equivalent assertion is that the zeros of $F (z)$ lie within a vertical strip $\abs{\Reg z} \leq C$.  
Thus let $\Reg z = x > 0$, and for $\delta > 0$ small, write
\[
z \hsy F (z) 
\ = \ 
e^z - e^{-z} 
\ - \ 
\int\limits_{-1}^{1 - \delta} \ 
e^{z \hsy t}
\ \td \phi
\ - \ 
\int\limits_{1 - \delta}^1 \ 
e^{z \hsy t}
\ \td \phi.
\]
Then
\allowdisplaybreaks
\begin{align*}
\bigg|
\hsx
\int\limits_{-1}^{1 - \delta} \ 
e^{z \hsy t}
\ \td \phi
\hsx
\bigg|
\ 
&\leq \ 
e^{x \hsy (1 - \delta)} \ 
\int\limits_{-1}^{1 - \delta} \ 
\abs{ \td \phi}
\\[15pt]
&< \ 
K \hsy e^{x \hsy (1 - \delta)}
\end{align*}
and
\allowdisplaybreaks
\begin{align*}
\bigg|
\hsx
\int\limits_{1 - \delta}^1 \ 
e^{z \hsy t}
\ \td \phi
\bigg|
\ 
&\leq \ 
e^x \hsy \max\limits_{1 - \delta \hsy < \hsy t_1 \hsy < \hsy t_2 \hsy \leq \hsy 1} \abs{\phi (t_2) - \phi (t_1)}
\\[15pt]
&=\ 
e^x \hsy M (\delta).
\end{align*}
Therefore
\[
\abs{z \hsy F (z) } 
\ \geq \ 
e^x (1 - e^{- 2 x} - K \hsy e^{-\delta x} - M (\delta)).
\]
Bearing in mind that $\phi (t)$ is continuous at $t = 1$, choose $\delta$ so small that 
$
\ds
M (\delta) < \frac{1}{4}
$. 
This done, choose $x$ so large that 
\[
e^{- 2 x} + K \hsy e^{-\delta x} 
\ < \ 
\frac{1}{4}.
\]
Then 
\allowdisplaybreaks
\begin{align*}
e^x (1 - e^{- 2 x} - K \hsy e^{-\delta x} - M (\delta)) \ 
&>\ 
e^x \Big(1 - \frac{1}{2}\Big)
\\[11pt]
&=\ 
\frac{e^x}{2}
\\[11pt]
&>\  
0.
\end{align*}
Consequently, for $x \gg 0$, $F (z)$ has no zeros.  
And, analogously, for $x \ll 0$, $F (z)$ has no zeros.  
\\[-.25cm]
\end{x}

\begin{x}{\small\bf REMARK} \ 
The result goes through if the assumption on $\phi$ at the endpoints is weakened to 
$\phi (1^-) \neq 0$, $\phi (-1^+) \neq 0$.
\\[-.25cm]
\end{x}

\begin{x}{\small\bf EXAMPLE} \ 
Let $\phi$ be defined on $]0,1[$.  
Suppose that $\phi$ is positive and 
increasing and
\[
\begin{cases}
\ 
\phi (1^-) 
\ < \ \infty
\\[4pt]
\ 
\phi (0^+) 
\ > \ 0
\end{cases}
.
\]
Then $\phi$ can be extended to a function of bounded variation on $[0,1]$.  
Taking $[a,b] = [0, 1]$, write
\[
\int\limits_0^1 \ 
\phi (t) 
\hsy 
e^{\sqrt{-1} \hsx z \hsy t}
\ \td t
\ = \ 
\frac{1}{2} 
\hsx 
e^{ \frac{1}{2}   \sqrt{-1} \hsx z}
\ \cdot \ 
\int\limits_{-1}^1 \ 
\phi \Big(\frac{1 + t}{2}\Big)
\hsx
e^{ \frac{1}{2}   \sqrt{-1} \hsx z \hsy t}
\ \td t
\]
to conclude that the zeros of $f (z)$ lie within a horizontal strip $\abs{\Img z} \leq C$. 
\\[-.25cm]
\end{x}

\begin{spacing}{1.55}
\begin{x}{\small\bf RAPPEL} \ 
Suppose that $\phi \in C[a,b]$.  
Given $\delta > 0$, let $\omega (\delta)$ be the supremum of 
$\abs{\phi (t_2) - \phi (t_1)}$ 
computed over all points $t_1$, $t_2$ in $[a,b]$ such that $\abs{t_2 - t_1} < \delta$   
$-$then $\omega (\delta)$ is called the 
\un{modulus of continuity} of $\phi$.  
As a function of $\delta$, $\omega$ is continuous and increasing and 
$
\lim\limits_{\delta \ra 0} \  
\omega (\delta)
= 0
$.  
In addition, 
$\omega (\delta) \geq A \hsy \delta$ for some $A > 0$ provided $\phi$ is not a constant.
\end{x}
\end{spacing}

\begin{x}{\small\bf THEOREM} \ 
Take $[a,b] = [-1, 1]$ and let $\phi \in C[-1,1]$, where 
$\phi (\pm 1) = 1$ $-$then all the zeros of 
\[
F (z) 
\ = \ 
\int\limits_{-1}^1 \ 
\phi (t) 
\hsy 
e^{ z \hsy t}
\ \td t
\]
which are sufficiently large in modulus lie in the set
\[
\abs{x}
\ \leq \ 
K \hsy r \hsy \omega \Big(\frac{1}{r}\Big)
\qquad (x = \Reg z, \ r = \abs{z}).
\]

\begin{spacing}{1.55}
PROOF \ 
It can be assumed that $\phi$ is not a constant (since otherwise $F (z)$ is
proportional to 
$
\ds
\frac{\sinh z}{z}
$ 
and there is nothing to prove).  
Proceeding, subdivide $[-1, 1]$ into $2 m$ equal parts and write
\end{spacing}

\[
\phi (t) 
\ = \ 
\phi \Big(\frac{j}{m}\Big) - \psi_j (t) 
\qquad 
\Big(\frac{j-1}{m} \leq t \leq \frac{j}{m}\Big).
\]
Then
\[
\abs{ \psi_j (t)} 
\ \leq \ 
\omega \Big(\frac{1}{m}\Big).
\]
\\[-.75cm]

\noindent
There are now two cases: \ 
$x > 0$ or $x < 0$, and it will be enough to consider the first of these.  
To begin with, 

\allowdisplaybreaks
\begin{align*}
F (z) \ 
&=\ 
\sum\limits_{j = -m + 1}^m \
\int\limits_{(j-1)/m}^{j / m} \ 
\Big(
\phi\Big(\frac{j}{m}\Big) 
- 
\psi_j (t) 
\Big)
e^{z \hsy t} 
\ \td t
\\[15pt]
&=\ 
\sum\limits_{j = -m + 1}^m \
\phi\Big(\frac{j}{m}\Big) 
\int\limits_{(j-1)/m}^{j / m} \ 
e^{z \hsy t} 
\ \td t
\ - \ 
\sum\limits_{j = -m + 1}^m \
\int\limits_{(j-1)/m}^{j / m} \ 
\psi_j (t) 
\hsx
e^{z \hsy t} 
\ \td t
\\[15pt]
&=\ 
I_1 + I_2.
\end{align*}
\\[-.75cm]

\qquad \textbullet \quad
\vspace{-.5cm}
\allowdisplaybreaks
\begin{align*}
\abs{I_2} \ 
&\leq \ 
\sum\limits_{j = -m + 1}^m \
\int\limits_{(j-1)/m}^{j / m} \ 
e^{x \hsy t} 
\hsx
\omega \Big(\frac{1}{m}\Big)
\ \td t
\\[15pt]
&=\ 
\omega \Big(\frac{1}{m}\Big) \ 
\int\limits_{-1}^1 \ 
e^{x \hsy t}
\ \td t
\\[15pt]
&=\ 
\omega \Big(\frac{1}{m}\Big) \ 
\frac{e^x - e^{-x}}{x}.
\end{align*}
\\[-.75cm]

\qquad \textbullet \quad
\allowdisplaybreaks
\begin{align*}
I_1 \ 
&=\ 
\sum\limits_{j = 0}^{2 m - 1} \
\phi \Big(1 - \frac{j}{m}\Big) 
\hsx
\frac
{
e^{z (1 - j/m)} 
- 
e^{z(1 - (j+1)/m)}
}
{z}
\\[15pt]
&=\ 
\frac{e^z}{z} 
+ 
\frac{e^z}{z}  \ 
\sum\limits_{j = 1}^{2 m - 1} \
\phi \Big(1 - \frac{j}{m}\Big) 
\hsx
\big(
e^{-z\hsy j/m}
- 
e^{-z\hsy (j+1)/m}
\big)
-
\frac{e^z}{z} 
\hsx
e^{-z/m}
\\[15pt]
&=\ 
\frac{e^z}{z} 
+ 
\frac{e^z}{z}  \ 
\sum\limits_{j = 1}^{2 m - 1} \
\Big(
\phi \Big(1 - \frac{j}{m}\Big) 
\hsx - \hsx 
\phi \Big(1 - \frac{j-1}{m}\Big) 
\Big) 
e^{-z\hsy j/m}
- 
\phi
\Big(-1 + \frac{1}{m}
\Big)
\hsx 
\frac{e^{-z}}{z}
\\[15pt]
&=\ 
\frac{e^z}{z} 
+ 
\frac{e^z}{z} \hsy I_3
-
\phi \Big(-1 + \frac{1}{m}\Big) 
\hsy
\frac{e^{-z}}{z}.
\end{align*}
\\[-.75cm]

\qquad \textbullet \quad
\vspace{-.5cm}
\allowdisplaybreaks
\begin{align*}
\abs{I_3} \ 
&\leq \ 
\sum\limits_{j = 1}^\infty \ 
\omega \Big(\frac{1}{m}\Big)
e^{- j \hsy x/m}
\\[15pt]
&=\ 
\omega \Big(\frac{1}{m}\Big) 
\ 
\frac{e^{-x/m}}{1 - e^{-x/m}}
\\[15pt]
&\leq\ 
\omega \Big(\frac{1}{m}\Big) 
\ 
\frac{m}{x}.
\end{align*}
\\[-.75cm]

[Note: \ 
For $\alpha > 0$,
\begin{align*}
1 + \alpha 
\hsx \leq \hsx 
e^\alpha
&\implies
\alpha 
\hsx \leq \hsx 
e^\alpha - 1
\\[11pt]
&\implies
\alpha
\hsx \leq \hsx
\frac{1 - e^{-\alpha}}{e^{-\alpha}}
\\[11pt]
&\implies
\alpha \hsy e^{-\alpha} 
\hsx \leq \hsx 
1 - e^{-\alpha}
\\[11pt]
&\implies
\frac{e^{-\alpha}}{1 - e^{-\alpha}}
\hsx \leq \hsx
\frac{1}{\alpha}.]
\end{align*}
Setting $m = [r]$, we have
\[
\omega \Big(\frac{1}{[r]}\Big)
\ \leq \ 
2 \hsy 
\omega \Big(\frac{1}{r}\Big)
\qquad (r \gg 0).
\]
\\[-.25cm]
Therefore
\allowdisplaybreaks
\begin{align*}
z \hsy F (z) \ 
&=\ 
z \hsy I_1 + z \hsy I_2
\\[15pt]
&=\ 
z 
\Big(
\frac{e^z}{z} + \frac{e^z}{z} \hsy I_3 
\hsx - \hsx
\phi
\Big(
-1 + \frac{1}{[r]}
\Big)
\hsy
\frac{e^{-z}}{z} 
\Big)
 + z \hsy I_2
 \\[15pt]
&=\ 
e^z
\Big(
1 + I_3 - \phi (-1 + \frac{1}{[r]}
\Big)
\hsy
e^{-2 z}
\Big)
 + z \hsy I_2
\\[15pt]
&=\ 
e^z
\Big(
1 + \tO \Big(\frac{r \omega(1/r)}{x}\Big)
\ - \ 
(1 + \txo (1)
)
\hsy 
e^{-2z}
\Big)
 + z \hsy I_2,
\end{align*}
where $o (1) \ra 0$ $(r \ra \infty)$.  
Next
\[
z I_2 
\ = \ 
e^z \hsy e^{-z} z I_2.
\]
And
\allowdisplaybreaks
\begin{align*}
\abs{e^{-z} z I_2} \ 
&\leq\ 
e^{-x} \hsy r \abs{I_2}
\\[15pt]
&\leq\ 
e^{-x} \hsy r \hsy \omega \bigg(\frac{1}{[r]}\bigg)
\frac{e^x - e^{-x}}{x}
\\[15pt]
&\leq\ 
2 \hsy r \hsy \omega \Big(\frac{1}{r}\Big) 
\hsx 
\frac{1- e^{-2x}}{x}
\\[15pt]
&=\ 
\tO \Big(\frac{r \omega(1/r)}{x}\Big).
\end{align*}
So in summary: \ 
$\forall \ r \gg 0$, 
\[
z \hsy F (z) 
\ = \ 
e^z 
\Big(
1 + 
\tO \Big(\frac{r \omega(1/r)}{x}\Big)
\ - \ 
(1 + \txo (1)) \hsy e^{-2z}
\Big).
\]
If $K > 0$ and if 
$
\ds
x > 
K \hsy r \hsy \omega \Big(\frac{1}{r}\Big)
$, 
then $x > A \hsy K$ (cf. 29.11), thus if $K$ is sufficiently large
\[
\Big |
\tO \Big(\frac{r \omega(1/r)}{x}\Big)
\ - \ 
(1 + \txo (1)) \hsy e^{-2z}
\Big |
\ \leq \ 
\frac{1}{2} 
\qquad (r \gg 0).
\]
\\[-1cm]

\noindent
But this implies that 
\[
1 
+ 
\tO \Big(\frac{r \omega(1/r)}{x}\Big) 
-
(1 + \txo (1)) \hsy e^{-2z}
\]
is bounded away from 0, hence $F (z)$ does not vanish in the region 
$
\ds
x > 
K \hsy r \hsy \omega \Big(\frac{1}{r}\Big)
$.
\\[-.25cm]
\end{x}

\begin{x}{\small\bf REMARK} \ 
The condition 
$\phi(\pm1) = 1$ 
can be replaced by the condition 
$\phi(\pm1) \neq 0$.
\\[-.25cm]
\end{x}

\begin{x}{\small\bf DEFINITION} \ 
A step function $\phi$ on $[0,1]$ of the form 
\[
\phi (t) 
\ = \ 
c_j 
\qquad (t_j < t < t_{j+1}),
\]
where
\[
0 = t_0 < t_1 < \cdots < t_n < t_{n+1} = 1
\]
and
\[
0 < c_0 < c_1 < \cdots < c_n, 
\]
is said to be \un{exceptional} if the $t_j$ are rational numbers.
\\[-.25cm]
\end{x}

\begin{x}{\small\bf NOTATION} \ 
Write $E (1, 0)$ for the set of exceptional step functions on $[0,1]$.
\\[-.25cm]
\end{x}

\begin{x}{\small\bf THEOREM} \ 
If $\phi \in \Lp^1[0,1]$ is positive and increasing on $]0,1[$ and if 
$\phi \notin E (1, 0)$, 
then the zeros of $f (z)$ lie in the open upper half-plane.
\\[-.5cm]

[We shall postpone the proof until later (cf. 34.2).]
\\[-.5cm]

[Note: \ 
In terms of $F (z)$, the conclusion is that its zeros  lie in the open left half-plane.]
\\[-.25cm]
\end{x}

\begin{x}{\small\bf EXAMPLE} \ 
The zeros of the real entire function 
\[
z \ra 
\int\limits_0^z \ 
e^{-t^2} 
\ \td t
\]
with the exception of $z = 0$ lie inside the region 
$\Reg z^2 < 0$ 
(a spiral in the complex plane).
\\[-.5cm]

[Write
\allowdisplaybreaks
\begin{align*}
\int\limits_0^z \ 
e^{-t^2} 
\ \td t \ 
&=\ 
\frac{z}{2} \ 
\int\limits_0^1 \ 
\frac{1}{\sqrt{t}} \hsx 
e^{-z^2} 
\ \td t \
\\[15pt]
&=\ 
\frac{z}{2} \ 
\int\limits_0^1 \ 
\frac{1}{\sqrt{1-t}} \hsx 
e^{-z^2 (1 - t)} 
\ \td t 
\\[15pt]
&=\ 
\frac{z}{2} 
\hsx
e^{-z^2} \ 
\int\limits_0^1 \ 
\frac{1}{\sqrt{1-t}} \hsx 
e^{-z^2 \hsy t} 
\ \td t.]
\end{align*}

[Note: \ 
The \un{error function} is defined by 
\[
\erf z 
\ = \ 
\frac{2}{\sqrt{\pi}} \ 
\int\limits_0^z \ 
e^{-t^2} 
\ \td t 
\]
and the \un{complementary error function} is defined by
\\[-.75cm]

\[
\erf_c \hsx z 
\ = \ 
\frac{2}{\sqrt{\pi}} \ 
\int\limits_z^\infty \ 
e^{-t^2} 
\ \td t.
\]
Therefore
\[
\erf z  + \erf_c \hsx z 
\ = \ 
1.
\]
\\[-.75cm]

The \un{Fresnel integrals} are defined by
\[
\begin{cases}
\ \ds
C (z) 
\ = \ 
\int\limits_0^z \ 
\cos \Big(\frac{\pi}{2} \hsx t^2\Big) 
\ \td t 
\\[26pt]
\ \ds
S (z) 
\ = \ 
\int\limits_0^z \ 
\sin \Big(\frac{\pi}{2} \hsx t^2\Big) 
\ \td t 
\end{cases}
.
\]
Accordingly, in terms of the error function,  
\[
C (z) + \sqrt{-1}\hsx S (z) 
\ = \ 
\frac{1 + \sqrt{-1}}{2} 
\hsx 
\erf \Big(\frac{\sqrt{\pi}}{2} (1 - \sqrt{-1}\hsy )\hsy z\Big).]
\]
Consider a step function $\phi$ per 29.14 $-$then
\[
f (z) 
\ = \ 
\sum\limits_{j = 0}^n \ 
c_j 
\ 
\int\limits_{t_j}^{t_j + 1} \ 
e^{\sqrt{-1} \hsx z \hsy t}
\ \td t 
\qquad (\implies f (0) > 0)
\]
$\implies$
\allowdisplaybreaks
\begin{align*}
\sqrt{-1} \hsx z \hsy f (z) \ 
&=\ 
c_0 
\Big(
e^{\sqrt{-1} \hsx z \hsy t_1}
\hsx - \hsx 
e^{\sqrt{-1} \hsx z \hsy t_0}
\Big)
\hsx + \hsx
c_1
\Big(
e^{\sqrt{-1} \hsx z \hsy t_2}
\hsx - \hsx 
e^{\sqrt{-1} \hsx z \hsy t_1}
\Big)
\\[15pt]
&
\hspace{2.5cm}
\ + \cdots + \
c_n 
\Big(
e^{\sqrt{-1} \hsx z \hsy t_{n+1}}
\hsx - \hsx 
e^{\sqrt{-1} \hsx z \hsy t_n}
\Big)
\\[15pt]
&=\ 
c_n \hsy e^{\sqrt{-1} \hsx z} 
- 
c_0 
- 
e^{\sqrt{-1} \hsx z \hsy t_1} 
\hsy 
(c_1 - c_0) 
- \cdots -
e^{\sqrt{-1} \hsx z \hsy t_n}
\hsy
(c_n - c_{n-1})
\end{align*}
\qquad 
$\implies$
\[
\abs{
\sqrt{-1} \hsx x \hsy f (x)
}
\ \geq \ 
c_n - c_0 - (c_1 - c_0) - \cdots - (c_n - c_{n-1})
\ = \ 
0.
\]
\\[-1.cm]
\end{x}

\begin{x}{\small\bf LEMMA} \ 
If for some $x \neq 0$, 
\[
\abs{\sqrt{-1} \hsx x \hsy f (x)}
\ = \ 
0,
\]
then 
$\phi \in E(1, 0)$.
\\[-.5cm]

PROOF \ 
The assumption implies that 
\[
e^{\sqrt{-1} \hsx x}
\ = \ 
1, 
\quad
e^{\sqrt{-1} \hsx x \hsy t_1}
\ = \ 
1, 
\ \ldots, \ 
e^{\sqrt{-1} \hsx x \hsy t_n}
\ = \ 
1, 
\]
from which the existence of integers 
$q, p_1, \ldots, p_n$ such that
\[
x 
\ = \ 
2 \hsy \pi \hsy q, 
\quad
x  \hsy t_1
\ = \ 
2 \hsy \pi \hsy p_1, 
\ \ldots, 
x \hsy t_n
\ = \ 
2 \hsy \pi \hsy p_n, 
\quad
\]
so 
\[
t_j 
\ = \ 
\frac{p_j}{q}.
\]
And this shows that 
$\phi \in E(1, 0)$.
\\[-.5cm]

[Note: \ 
If $x$ is positive, then $q$ and the $p_j$ are positive but if $x$ is negative, 
then $q$ and the $p_j$ are negative and we write
\[
t_j 
\ = \ 
\frac{-p_j}{-q}.]
\]
\\[-1.25cm]
\end{x}

If $\phi$ is a step function and if 
$\phi \notin E(1, 0)$, 
then
\[
x \neq 0 
\implies 
\abs{\sqrt{-1} \hsx x \hsy f (x)}
> 0,
\]
thus $f (z)$ has no real zeros.  
Now fix $y < 0$ and consider 
\allowdisplaybreaks
\begin{align*}
f (z) \ 
&= \ 
f (x + \sqrt{-1} \hsx y) 
\\[15pt]
&= \ 
\int\limits_0^1 \ 
\phi (t) 
\hsx
e^{\sqrt{-1} \hsx (x + \sqrt{-1} \hsx y)}
\ \td t
\\[15pt]
&= \ 
\int\limits_0^1 \ 
\big(
\phi (t) 
\hsy 
e^{-y \hsy t}
\big)
\hsx 
e^{\sqrt{-1} \hsx x}
\ \td t.
\end{align*}
Since $y$ is negative, the function 
$\phi (t) \hsy e^{-y \hsy t} $
is positive and increasing on $]0,1[$ and it is obviously not in $E(0,1)$.  
Therefore, on the basis of 29.16, 
\[
\int\limits_0^1 \ 
\big(
\phi (t) 
\hsy 
e^{-y \hsy t}
\big)
\hsx 
e^{\sqrt{-1} \hsx x}
\ \td t
\]
does not vanish on the real axis, so $f (z)$ does not vanish on the line 
$\Img z = y$.
\\[-.25cm]

\begin{x}{\small\bf SCHOLIUM} \ 
If $\phi$ is a step function and if 
$\phi \notin E(1, 0)$, 
then the zeros of $f (z)$ lie in the open upper half-plane.
\\[-.5cm]

[Note: \ 
This is an important point of principle: \ 
If $\phi$ is a step function, then it either is in $E(1, 0)$ or it isn't and if it isn't, then the truth of 29.16 for those $\phi$ 
which are not step functions implies the truth of 29.16 for those step functions $\phi \notin E(1, 0)$.]
\\[-.25cm]
\end{x}

\begin{x}{\small\bf LEMMA} \ 
If $\phi \in E(1, 0)$, 
then $f (z)$ has a real zero.
\\[-.25cm]

PROOF \ 
Let 
\[
t_1 
\ = \ 
\frac{p_1}{q_1} 
\  (q_1 > 0), 
\quad
t_2 
\ = \ 
\frac{p_2}{q_2} 
\  (q_2 > 0), 
\quad
\cdots, 
\quad
t_n 
\ = \ 
\frac{p_n}{q_n} 
\  (q_n > 0). 
\quad
\]
Put
\[
q 
\ = \ 
q_1 \cdots q_n, 
\quad 
a_j 
\ = \ 
\frac{p_j q}{q_j}
\quad
\big(
\implies t_j = \frac{a_j}{q} 
\quad 
(j = 1, \ldots, n)
\big) 
\]
and set $x = 2 \pi q$ $-$then
\[
e^{\sqrt{-1} \hsx x}
\ = \ 
e^{\sqrt{-1} \hsx 2 \hsy \pi \hsy q}
\ = \ 
1
\]
and
\[
e^{\sqrt{-1} \hsx x \hsy t_j}
\ = \ 
e^{\sqrt{-1} \hsx 2 \hsy \pi \hsy q \hsy t_j}
\ = \ 
e^{\sqrt{-1} \hsx 2 \hsy \pi \hsy a_j}
\ = \ 
1
\qquad (j = 1, \ldots, n).
\]
Therefore
\allowdisplaybreaks
\begin{align*}
\sqrt{-1} \hsx &(2 \hsy \pi \hsy q) f (2 \hsy \pi \hsy q) 
\\[11pt]
&=\ 
c_n \hsy e^{\sqrt{-1} \hsx 2 \hsy \pi \hsy q} 
- 
c_0 
- 
e^{\sqrt{-1} \hsx 2 \hsy \pi \hsy q \hsy t_1} \hsy (c_1 - c_0) 
\hsx - \cdots - \hsx 
e^{\sqrt{-1} \hsx 2 \hsy \pi \hsy q \hsy t_n} \hsy (c_n - c_{n-1}) 
\\[11pt]
&=\ 
c_n - c_0 - (c_1 - c_0) 
\hsx - \cdots - \hsx
(c_n - c_{n-1})
\\[11pt]
&=\
0
\end{align*}
\qquad 
$\implies$
\[
f (x) 
\ = \ 
f(2 \hsy \pi \hsy q) 
\ = \ 
0.
\]
\\[-1.cm]
\end{x}

\begin{x}{\small\bf THEOREM} \ 
If $\phi \in E(1, 0)$, 
then $f (z)$ has an infinity of real zeros.
\\[-.5cm]

PROOF \ 
Write
\[
\sqrt{-1} \hsx z f (z) 
\ = \ 
P 
\big(
e^{\sqrt{-1} \hsx z/q} 
\big),
\]
where $P$ is a polynomial of degree $q$ $-$then 
$P (1) = 0$ (set $z = 0$), hence
\[
\sqrt{-1} \hsx z f (z) 
\ = \ 
\big(
e^{\sqrt{-1} \hsx z/q} - 1
\big)
\hsx
P_1 
\big(
e^{\sqrt{-1} \hsx z/q} 
\big).
\]
Therefore
\[
\pm 2 \hsy \pi \hsy q, 
\hsx
\pm 4 \hsy \pi \hsy q, 
\hsx
\ldots
\]
are zeros of $f (z)$.
\\[-.25cm]
\end{x}

Let $u = e^{\sqrt{-1} \hsx z/q}$ $-$then 
\allowdisplaybreaks
\begin{align*}
\sqrt{-1} \hsx z f (z) \ 
&=\ 
c_0 
\big(
u^{a_1} - 1
\big)
\hsx + \hsx
c_1 
\big(
u^{a_2} - u^{a_1}
\big)
\hsx + \cdots + \hsx 
c_n 
\big(
u^q - u^{a_n}
\big)
\\[11pt]
&=\ 
(u - 1) 
\hsy 
\big(
c_0 + c_0 u + \cdots + c_0 u^{a_1 - 1} + 
c_1 u^{a_1} + \cdots + c_n u^{q-1}
\big)
\\[11pt]
&=\ 
(u - 1) 
\hsy 
P_1 (u).
\end{align*}
Thanks to wellknown generalities (explicated in \S30 (cf. 30.13)), 
the structure of the  coefficients of $P_1$ confines the zeros of $P_1$ to the closed unit disk $\abs{u} \leq 1$, 
thus, in terms of $z$:
\allowdisplaybreaks
\begin{align*}
\big|
e^{\sqrt{-1} \hsx z / q}
\big|
\ \leq \ 
1
&\implies
\big|
e^{\sqrt{-1} \hsx (x + \sqrt{-1} \hsx y) / q}
\big|
\ \leq \ 
1
\\[11pt]
&\implies
\big|
e^{(\sqrt{-1} \hsx x - y) / q}
\big|
\ \leq \ 
1
\\[11pt]
&\implies
e^{-y / q} 
\ \leq \ 
1
\\[11pt]
&\implies
-y / q 
\ \leq \  
0
\\[11pt]
&\implies
y 
\ \leq \  
0.
\end{align*}

[Note: \ 
Any zero of $P_1$ on the unit circle $\abs{u} = 1$ is necessarily simple, so the real zeros of $f (z)$ are simple.]
\\[-.25cm]

\begin{x}{\small\bf LEMMA} \ 
If $\phi \in E(1, 0)$, 
then the zeros of $f (z)$ lie on a finite set of horizontal straight lines 
$\Img z = b_k$ 
$(b_k \geq 0, 1 \leq k \leq s, s \leq q)$.  
\\[-.5cm]

[In terms of distinct roots $w_1 = 1$, $w_2, \ldots, w_s$ of $P$, 
\[
b_k 
\ = \ 
- q \log \abs{w_k}.]
\]

[Note: \ 
These lines are not necessarily distinct.  
E.g., if 
$w_k = \sqrt{-1} \hsx$, the associated horizontal straight line is the real axis and the zeros are situated at
\[
q \hsy \frac{\pi}{2}, \ 
q \hsy \Big(\frac{\pi}{2} \pm 2 \hsy \pi\Big), \ 
q \hsy \Big(\frac{\pi}{2} \pm 4 \hsy \pi\Big), \ 
\ldots \ .]
\]
\\[-1.25cm]
\end{x}

Here is an application of 29.16.
\\[-.25cm]

\begin{x}{\small\bf THEOREM} \ 
If $\phi \in \Lp^1 [0,1]$ is positive and differentiable on $]0,1[$ with 
\[
\alpha 
\ \leq \ 
- \frac{\phi^\prime (t) }{\phi (t)}
\ \leq \ 
\beta 
\qquad (0 < t < 1)
\]
and if 
\[
\phi (t) 
\ \neq \ 
C e^{- \alpha \hsy t}, \ 
C e^{- \beta \hsy t}, 
\]
then the zeros of 
\[
F (z) 
\ = \ 
\int\limits_0^1 \ 
\phi (t) 
\hsy 
e^{z \hsy t} 
\ \td t
\]
are confined to the open strip 
$\alpha < \Reg z < \beta$.
\\[-.5cm]

PROOF \ 
Write
\[
F (z) 
\ = \ 
\int\limits_0^1 \ 
e^{\beta \hsy t}
\hsy
\phi (t) 
\hsy 
e^{(z - \beta) \hsy t} 
\ \td t.
\]
Then
\[
\frac{\td}{\td t} 
\hsy
\big(
e^{\beta \hsy t}
\hsy
\phi (t) 
\big)
\ = \ 
e^{\beta \hsy t}
\hsy
\phi (t) 
\Big(
\frac{\phi^\prime (t)}{\phi (t)} 
+
\beta
\Big)
\ \geq \ 
0.
\]
Therefore the zeros of $F (z)$ are restricted by the relation
\[
\Reg (z - \beta) 
\ < \ 
0
\qquad (\tcf. \ 29.16).
\]
Write
\[
F (z) 
\ = \ 
e^z 
\ 
\int\limits_0^1 \ 
e^{-\alpha \hsy t}
\hsy
\phi (1 - t) 
\hsy 
e^{(\alpha - z) \hsy t} 
\ \td t.
\]
Then
\[
\frac{\td}{\td t} 
\hsy
\big(
e^{-\alpha \hsy t}
\hsy
\phi (1 - t) 
\big)
\ = \ 
e^{-\alpha \hsy t}
\hsy
\phi (1 - t) 
\Big(
- \frac{\phi^\prime (1-t) }{\phi (1-t)} 
- 
\alpha
\Big)
\ \geq \ 
0.
\]
Therefore the zeros of $F (z)$ are restricted by the relation
\[
\Reg (\alpha - z) 
\ < \ 
0
\qquad (\tcf. \ 29.16).
\]
But 
\[
\begin{cases}
\ 
\Reg (z - \beta) 
\ < \ 0
\\[11pt]
\ 
\Reg (\alpha - z) 
\ < \ 
0
\end{cases}
\implies 
\alpha 
\ < \ 
\Reg z 
\ < \ 
\beta.
\]
\\[-.75cm]
\end{x}

\begin{x}{\small\bf EXAMPLE} \ 
Take $\phi (t) = \exp (-e^t)$ \ $-$then 
\[
- \frac{\phi^\prime (t) }{\phi (t)}
\ = \ 
e^t
\]
and
\[
1 
\ \leq \ 
e^t
\ \leq \ 
e
\qquad (0 < t < 1).
\]
Consequently, $\forall \ \varepsilon > 0$, the zeros of 
\[
F (z) 
\ = \ 
\int\limits_0^1 \
\exp (- e^t) 
\hsy
e^{z \hsy t} 
\ \td t
\]
are confined to the open strip
\[
1 - \varepsilon
\ < \ 
\Reg z
\ < \ 
e + \varepsilon
\]
or still, to the closed strip
\[
1 
\ \leq \ 
\Reg z
\ \leq \ 
e.
\]
\\[-1cm]
\end{x}

\begin{x}{\small\bf EXAMPLE} \ 
Given a complex parameter $\mu$, let 
\[
E (z; \mu)
\ = \ 
\sum\limits_{n = 0}^\infty \ 
\frac{z^n}{\Gamma (\mu + n)},
\]
an entire function of $z$.  
In particular: 
\[
\begin{cases}
\ 
e^z
\ = \ 
E (z; 1)
\\[4pt]
\ 
z e^z
\ = \ 
E (z; 0)
\end{cases}
\]
and 
\[
z^{1  - \mu} \hsy e^z 
\ = \ 
E (z; \mu)
\qquad (\mu = -1, -2, \ldots).
\]
\\[-1.25cm]

\un{Differential Equations} : \ 
\\[-.25cm]

\qquad \textbullet \quad
$
(\mu - 1) \hsy E (z; \mu)
+ 
z \hsy E^\prime (z; \mu) 
\ = \ 
E (z; \mu - 1)
$
\\[-.25cm]

\qquad \textbullet \quad
$
E (z; \mu)
-
E^\prime (z; \mu) 
\ = \ 
(\mu - 1) \hsy E (z; \mu +1)
$
\\[-.25cm]

Suppose now that $\mu > 1$ $-$then 
\[
E (z; \mu)
\ = \ 
\int\limits_0^1 \ 
\phi (t) 
\hsx
e^{z \hsy t}
\ \td t, 
\]
where 
\[
\phi (t) 
\ = \
\frac{(1 - t)^{\mu - 2}}{\Gamma (\mu - 1)}, 
\]
thus
\[- 
\frac{\phi^\prime (t) }{\phi (t)}
\ = \ 
\frac{\mu - 2}{1 - t}
\qquad (0 < t < 1)
\]

\[
\implies
\hspace{2cm}
\begin{cases}
\ \ds
-
\frac{\phi^\prime (t) }{\phi (t)}
\ \leq \ \mu - 2
\qquad (1 < \mu < 2)
\\[18pt]
\ \ds
-
\frac{\phi^\prime (t) }{\phi (t)}
\ \geq \ \mu - 2
\qquad (\mu > 2)
\end{cases}
. \hspace{2cm}
\]
So, the zeros of $E (z; \mu)$ lie in the region 
$\Reg z < \mu - 2$ 
if 
$1 < \mu < 2$ 
and in the region 
$\Reg z > \mu - 2$ 
if 
$\mu > 2$.
\\[-.25cm]

\un{$1 < \mu < 2$}: \quad 
The zeros of $E (z; \mu)$ are simple.  
In fact, if $E (z; \mu)$ had a multiple zero $z_0$, then 
\[
E (z_0; \mu + 1)
\ = \ 
0.
\]
But
\[
\mu + 1 
\ > \ 
2 
\hsx \implies \hsx 
\Reg z_0 
\ > \ 
(\mu + 1) - 2
\ = \ 
\mu - 1
\ > \ 
0
\]
in contradiction to 
\[
\Reg z_0 
\ < \ 
\mu - 2
\ < \ 
0.
\]
\\[-.75cm]

\un{$2 \leq \mu \leq 3$}: \quad 
First
\[
E (z; 2) 
\ = \ 
\frac{e^z - 1}{z}
\]
and its zeros are simple and lie on the imaginary axis.  
Assume, therefore, that 
$2 < \mu \leq 3$ $-$then the zeros of 
$E(z; \mu)$ are also simple.  
For at a multiple zero $z_0$, we would have
\[
E (z_0; \mu - 1) 
\ = \ 
0
\]
from which 
\[
\Reg z_0 
\ \leq \ 
\mu - 1 - 2 
\ \leq \ 
3 - 3 
\ = \ 
0,
\]
contradicting
\[
\Reg z_0 
\ > \ 
\mu - 2 
\ > \ 
0.
\]
\\[-1.25cm]
\end{x}

\begin{x}{\small\bf EXAMPLE} \ 
The 
\un{incomplete gamma function} 
is defined by the rule
\[
\gamma (\alpha, z)
\ = \ 
\int\limits_0^z \ 
e^{-t} 
\hsy 
t^{\alpha - 1}
\ \td t
\qquad (\Reg \alpha > 0).
\]
As a function of $z$, $\gamma (\alpha, z)$ is holomorphic with the potential exception of a branch point at  the origin, 
the principal branch being determined by introducing a cut along the negative real $t$ axis and requiring 
$t^{\alpha - 1}$ to have its principal value.  
Expanding $e^{-t}$ and integrating gives
\[
\gamma (\alpha, z)
\ = \ 
z^\alpha 
\ 
\sum\limits_{n = 0}^\infty \
(-1)^n 
\hsx 
\frac{z^n}{n ! (n + \alpha)}\hsy ,
\]
the right hand side providing an extension of the left hand side to all $\alpha \neq 0$, 
$-1, -2, \ldots \ .$  
Put
\[
\gamma^* (\alpha, z) 
\ = \ 
\frac{\gamma (\alpha, z) }{z^\alpha \Gamma (\alpha)}\hsy .
\]
Then 
$\gamma^* (\alpha, z)$ 
is entire and 
\[
\gamma^* (\alpha, z) 
\ = \ 
e^{-z} \ 
\sum\limits_{n = 0}^\infty \ 
\frac{z^n}{\Gamma (\alpha + n + 1)}
\]
or still, 
\[
\gamma^* (\alpha, z) 
\ = \ 
e^{-z} \ 
E (z; 1 + \alpha).
\]
Specializing what has been said in 29.25, we can thus say the following.
\\[-.25cm]

\qquad \textbullet \quad
For 
$0 < \alpha < 1$, 
all the zeros of 
$\gamma^* (\alpha, z) $
lie in the region 
$\Reg z < \alpha - 1$.
\\[-.25cm]

\qquad \textbullet \quad
For 
$\alpha > 1$, 
all the zeros of 
$\gamma^* (\alpha, z) $
lie in the region 
$\Reg z > \alpha - 1$.
\\[-.25cm]

\qquad \textbullet \quad
For 
$0 < \alpha \leq 2$, 
all the zeros of 
$\gamma^* (\alpha, z) $
are simple.
\\[-.5cm]

[Note: \ 
\[
\gamma^* (0, z) 
\ \equiv \ 
1 
\quad \text{and} \quad 
\gamma^* (-n, z) 
\ = \ 
z^n
\qquad 
(n = 1, 2, \ldots).]
\]
\\[-1.25cm]
\end{x}

\begin{x}{\small\bf EXAMPLE} \ 
Consider the error function 
\[
\erf z 
\ = \ 
\frac{2}{\sqrt{\pi}} 
\int\limits_0^z \ 
e^{-t^2}
\ \td t
\qquad (\tcf. \ 29.17).
\]
Then $\erf z$ has a simple zero at $z = 0$ and no other real zeros.  
Since
\[
\erf z 
\ = \ 
\frac{1}{\sqrt{\pi}} 
\hsx
\gamma \Big(\frac{1}{2}, t^2\Big),
\]
the nonreal zeros of $\erf z$ coincide with the zeros of 
$
\ds
\gamma^* \Big(\frac{1}{2}, z^2\Big)
$, 
these lying in the region 
$
\ds
\Reg z^2 < - \frac{1}{2}
$ 
(which, when explicated, is seen to consist of two curvilinear sectors placed symmetrically with respect to the real axis and bounded by
the components of the hyperbola 
\ 
$
\ds
y^2 - x^2 
\ = \ 
\frac{1}{2}
$ 
\ 
$(z = x + \sqrt{-1} \hsx y)$.
\\[-.5cm]

[Note: \ 
It can be shown that the zeros of $\erf z $ are simple.  
In addition, the nonreal zeros of 
$\erf z$ are comprised of two sequences 
$z_n^+$, $z_n^-$, 
$(n = \pm 1, \pm 2, \ldots)$ 
which are symmetric with respect to the real axis and contained in the region 
\ 
$
\ds
y^2 - x^2 
\ > \ 
\frac{1}{2}.
$
And asymptotically, 
\[
\big(z_n^\pm\big)^2 
\ = \ 
2 \hsy \pi \hsx n \sqrt{-1} \hsx 
\hsx - \hsx 
\frac{1}{2}
\log \abs{n} 
\hsx - \hsx 
\sqrt{-1} 
\  
\frac{\pi}{4}
\hsx
\sgn n
\hsx - \hsx 
\log (\pi \sqrt{2})
\hsx + \hsx 
\tO \Big(\frac{\log \abs{n}}{\abs{n}}\Big)
\quad 
(n \ra \infty).]
\]
\\[-.25cm]
\end{x}

\chapter{
$\boldsymbol{\S}$\textbf{30}.\quad  TRANSFORM THEORY: \ JUNIOR GRADE}
\setlength\parindent{2em}
\setcounter{theoremn}{0}
\renewcommand{\thepage}{\S30-\arabic{page}}

\qquad 
If $\phi \in \Lp^1 [0,1]$, then by definition 
\[
f (z) 
\ = \ 
\int\limits_0^1 \ 
\phi (t) 
\hsx 
e^{\sqrt{-1} \hsx z \hsy t}
\ \td t
\]
or still, 
\[
f (z) 
\ = \ 
C (z) 
\hsx + \hsx 
\sqrt{-1} \hsx S (z),
\]
where
\[
\begin{cases}
\ \ds
C (z) 
\ = \ 
\int\limits_0^1 \ 
\phi (t) 
\hsx 
\cos z \hsy t
\ \td t
\\[26pt]
\ \ds
S (z) 
\ = \ 
\int\limits_0^1 \ 
\phi (t) 
\hsx 
\sin z \hsy t
\ \td t
\end{cases}
.
\]
\\[-.25cm]

\begin{x}{\small\bf EXAMPLE} \ 
Take 
$
\ds
\phi (t) = 
\frac{1}{\sqrt{1 - t^2}} \ \
(0 \leq t < 1)$ 
$-$then
\[
\frac{2}{\pi} \ 
\int\limits_0^1 \ 
\frac{\cos z \hsy t}{\sqrt{1 - t^2}}
\ \td t
\ = \ 
\tJ_0 (z).
\]
\\[-1.25cm]
\end{x}

Extend $\phi$ to an even function 
$\widetilde{\phi}$ 
on 
$[-1, 1]$ 
and let 
\[
\widetilde{C} (z) 
\ = \ 
\int\limits_{-1}^1 \ 
\widetilde{\phi} (t)
\hsx 
\cos z \hsy t 
\ \td t,
\]
thus 
\[
\widetilde{C} (z) 
\ = \ 
\sum\limits_{n = 0}^\infty \ 
\frac{(-1)^n z^{2 n}}{2 n !}
\ 
\int\limits_{-1}^1 \ 
\tilde{\phi} (t)
\hsx 
t^{2 n} 
\ \td t.
\]
\\[-.5cm]

\begin{x}{\small\bf RAPPEL} \ 
The $n^\nth$ Appell polynomial $\tJ_n^*$ associated with a real entire function $f$ is defined by 

\[
\tJ_n^* (f; z) 
\ = \ 
\sum\limits_{k = 0}^\infty \ 
\binom{n}{k}
\hsx
\gamma_k 
\hsy 
z^{n - k}
\qquad (\tcf. \ 12.4).
\]
\\[-1.25cm]
\end{x}

\begin{x}{\small\bf LEMMA} \ 
We have
\[
\tJ_n^* (\widetilde{C}; z) 
\ = \ 
\int\limits_{-1}^1 \ 
\widetilde{\phi} (t) 
\hsy 
\hsx 
(z + \sqrt{-1} \hsx t)^n
\ \td t.
\]

PROOF \ 
Expand the RHS: 
\allowdisplaybreaks
\begin{align*}
\int\limits_{-1}^1 \ 
\widetilde{\phi} (t)
\hsx 
(z + \sqrt{-1} \hsx t)^n
\ \td t \ 
&=\ 
\int\limits_{-1}^1 \ 
\widetilde{\phi} (t)
\hsx 
(\sqrt{-1} \hsx t + z)^n
\ \td t \ 
\\[15pt]
&=\ 
\sum\limits_{k = 0}^n \ 
\binom{n}{k}
\hsx
(\sqrt{-1} \hsx)^k
\ 
\Big(
\int\limits_{-1}^1 \ 
\widetilde{\phi} (t)
\hsx 
t^k
\ \td t
\Big)
z^{n - k}
\\[15pt]
&=\ 
\sum\limits_{k = 0}^{[n/2]} \ 
\binom{n}{2 k}
\hsx
(-1)^k
\ 
\Big(
\int\limits_{-1}^1 \ 
\widetilde{\phi} (t)
\hsx 
t^{2 k}
\ \td t
\Big)
z^{n - 2 k}.
\end{align*}
On the other hand, from the definitions, 
\allowdisplaybreaks
\begin{align*}
&
\gamma_0
\ = \ 
\int\limits_{-1}^1 \ 
\widetilde{\phi} (t)
\ \td t,
\hspace{1.35cm}
\gamma_1
\ = \ 
0,
\\[15pt]
&
\gamma_2
\ = \ 
-
\int\limits_{-1}^1 \ 
\widetilde{\phi} (t)
\hsx
t^2
\ \td t,
\hspace{.6cm}
\gamma_3
\ = \ 
0,
\\[15pt]
&
\gamma_4
\ = \ 
\int\limits_{-1}^1 \ 
\widetilde{\phi} (t)
\hsx 
t^4
\ \td t,
\hspace{1cm}
\gamma_5
\ = \ 
0,
\\[15pt]
&
\hspace{2cm}
\vdots
\end{align*}
\\[-1.25cm]
\end{x}

\begin{x}{\small\bf RAPPEL} \ 
The $n^\nth$ Jensen polynomial $\tJ_n$ associated with a real entire function $f$ is defined by 
\[
\tJ_n (f; z) 
\ = \ 
\sum\limits_{k = 0}^n \ 
\binom{n}{k}
\hsx 
\gamma_k 
\hsx z^k
\qquad (\tcf. \ 12.1).
\]
\\[-1cm]
\end{x}

\begin{x}{\small\bf LEMMA} \ 
We have
\[
\tJ_n (\widetilde{C}; z) 
\ = \ 
\int\limits_{-1}^1 \ 
\widetilde{\phi} (t)
\hsx
(1 + \sqrt{-1} \hsx z \hsy t)^n 
\ \td t.
\]

PROOF \ 
In fact, 
\allowdisplaybreaks
\begin{align*}
\tJ_n (\widetilde{C}; z) \ 
&=\ 
z^n 
\hsy 
\tJ_n^* \Big(\widetilde{C}; \frac{1}{z}\Big)
\\[15pt]
&=\ 
z^n 
\int\limits_{-1}^1 \ 
\widetilde{\phi} 
 (t)
\hsx
\Big(
\frac{1}{z} + \sqrt{-1} \hsx t
\Big)^n
\ \td t
\\[15pt]
&=\ 
z^n 
\int\limits_{-1}^1 \ 
\widetilde{\phi} 
 (t)
\hsx
\Big(
\frac{1 + \sqrt{-1} \hsx z \hsy t}{z}
\Big)^n
\ \td t
\\[15pt]
&=\ 
\int\limits_{-1}^1 \ 
\widetilde{\phi} 
 (t)
\hsx
(1 + \sqrt{-1} \hsx z \hsy t)^n 
\ \td t.
\end{align*}
\\[-1cm]
\end{x}

\begin{x}{\small\bf EXAMPLE} \ 
Take 
$\phi (t) = (1 - t^{2 \hsy p})^\lambda$, 
where 
$p = 1, 2, \ldots$, 
and 
$\lambda > -1$ 
$-$then the real polynomial
\[
\int\limits_{-1}^1 \ 
(1 - t^{2 \hsy p})^\lambda 
\hsx 
(1 + \sqrt{-1} \hsx z \hsy t)^n
\ \td t
\qquad (n > 1)
\]
has real zeros only, hence the real entire function 
\[
\int\limits_0^1 \ 
(1 - t^{2 \hsy p})^\lambda 
\hsx 
\cos z \hsy t
\ \td t
\]
has real zeros only (being in $\sL - \sP$ (cf. 12.14)).
\\[-.5cm]

[Note: \ 
It is known that for 
$
\ds
\nu > -\frac{1}{2}
$, 
\[
\tJ_\nu (z) 
\ = \ 
\frac{2}{\sqrt{\pi} \ \Gamma(\nu + \frac{1}{2})}
\
\Big(\frac{z}{2}\Big)^\nu
\
\int\limits_0^1 \ 
(1 - t^2)^{\nu - \frac{1}{2}}
\hsx 
\cos z \hsy t
\ \td t.
\]
But then 
$
\ds
\nu  -\frac{1}{2} > -1
$, 
so the zeros of $\tJ_\nu (z)$ are real (cf. 12.33) (matters there require only that $\nu > -1$).]
\\[-.25cm]
\end{x}

\begin{x}{\small\bf REMARK} \ 
Let $\lambda = k = 1, 2, \ldots,$ 
and replace $z$ by $z \hsy  k^{1/2 p}$: 
\[
\int\limits_0^1 \ 
(1 - t^{2 \hsy p})^k 
\hsx 
\cos z k^{1/2 p} 
\hsx 
t
\ \td t.
\]
Then make the change of variable 
$t = x \hsy k^{-1/2 p}$:
\[
k^{-1/2 p} \ 
\int\limits_0^{k^{1/ 2 p}} \ 
\Big(
1 - \frac{x^{2 \hsy p}}{k}
\Big)^k
\hsx
\cos z x 
\ \td x.
\]
Now replace $x$ by $t$ and form 
\[
\lim\limits_{k \ra \infty} \ 
\int\limits_0^{k^{1/ 2 p}} \ 
\Big(
1 - \frac{t^{2 p}}{k}
\Big)^k
\hsx
\cos z \hsy t 
\ \td t
\]
to see that the real entire function 
\[
\Phi_{2 p} (z) 
\ = \ 
\int\limits_0^\infty \ 
\exp (- t^{2 p}) 
\hsy 
\cos z \hsy t 
\ \td t
\]
has real zeros only (cf. 12.34).
\\[-.25cm]
\end{x}

\begin{x}{\small\bf THEOREM} \ 
Suppose that $\phi (t)$ is positive, strictly increasing, and continuous on $[0,1[$ and 
\[
\int\limits_0^1 \ 
\phi (t) 
\ \td t
\ = \ 
\lim\limits_{\varepsilon \ra 0} \ 
\int\limits_0^{1 - \varepsilon} \ 
\phi (t) 
\ \td t
\]
exists $-$then the real entire function 
\[
C (z)
\ = \ 
\int\limits_0^1 \ 
\phi (t) 
\cos z \hsy t
\ \td t
\]
has real zeros only.
\\[-.25cm]
\end{x}

\qquad
{\small\bf \un{N.B.}} \ 
Accordingly, 
\[
\lim\limits_{n \ra \infty} \ 
\frac
{
\phi
\Big(
\frac{1}{n}
\Big)
\hsx + \hsx 
\phi
\Big(
\frac{2}{n}
\Big)
\hsx + \hsx 
\cdots
\hsx + \hsx 
\phi
\Big(
\frac{n-1}{n}
\Big)
}
{n}
\ = \ 
\int\limits_0^1 \ 
\phi (t) 
\ \td t.
\]

[The expression on the left (sans the limit) is bounded from below by
\[
\int\limits_0^{1 - \frac{1}{n}} \ 
\phi (t) 
\ \td t
\]
and from above by
\[
\int\limits_{\frac{1}{n}}^1 \ 
\phi (t) 
\ \td t.]
\]
\\[-1cm]

\begin{x}{\small\bf REMARK} \ 
The assumptions on $\phi$ can be weakened (cf 31.1) but the methods utilized in arriving at 30.8 
are instructive and can be employed in other situations as well. 
\\[-.25cm]
\end{x}

\begin{x}{\small\bf LEMMA} \ 
Suppose given polynomials
\[
\begin{cases} 
\ 
P (z) 
\ = \ 
a_n (z - z_1) (z - z_2) \cdots (z - z_n)
\\[4pt]
\ 
Q (z) 
\ = \ 
\bar{a}_n (1 - \bar{z}_1 z) (1 - \bar{z}_2 z) \cdots (1 - \bar{z}_n z)
\end{cases}
.
\]
Assume: \ 
The zeros of $P (z)$ lie in the region $\abs{z} \geq 1$ $-$then the zeros of 
\[
P (z)
\hsx + \hsx 
\gamma \hsy z^k \hsy Q (z) 
\qquad (\abs{\gamma} = 1, \ k = 1, 2, \ldots)
\]
lie on the unit circle $\abs{z} = 1$.
\\[-.5cm]

PROOF \ 
There are two points.
\\[-.25cm]

\qquad \textbullet \quad
If $\abs{w} > 1$, then 
\[
\abs{\frac{z - w}{1 - \bar{w} z}}
\ =_{\raisebox{-.2cm}{$\hspace{-.3cm}<$}}^{\raisebox{.12cm}{$\hspace{-.3cm}>$}} \ 
1
\quad \text{for} 
\abs{z} 
\ =_{\raisebox{-.2cm}{$\hspace{-.3cm}>$}}^{\raisebox{.12cm}{$\hspace{-.3cm}<$}} \ 
1.
\]
\\[-1cm]

\qquad \textbullet \quad
If $\abs{\omega} = 1$, then 
\[
\abs{\frac{z - \omega}{1 - \bar{\omega} z}}
\ = \ 
\abs{\frac{z - \omega}{\omega - z}}
\quad \text{for} 
\abs{z} 
\ =_{\raisebox{-.2cm}{$\hspace{-.3cm}>$}}^{\raisebox{.12cm}{$\hspace{-.3cm}<$}} \ 
1.
\]
Therefore the equality is possible only when $\abs{z} = 1$.
\\[-.25cm]
\end{x}


\begin{x}{\small\bf REMARK} \ 
If $\abs{z_i} > 1$ $(i = 1, \ldots, n)$, then the zeros of 

\[
P (z)
\hsx + \hsx 
\gamma \hsy z^k \hsy Q (z) 
\]
are simple.
\\[-.5cm]

[Let 
$p (z) = P (z)$, 
$q (z) = - \gamma \hsy z^k \hsy Q (z)$, 
and suppose that $z_0$ is a multiple zero of 
$p (z) - q (z)$ $-$then
\[
\begin{cases}
\ 
p (z_0) \ = \ q (z_0)
\\[4pt]
\ 
p^\prime (z_0) \ = \ q^\prime (z_0)
\end{cases}
.
\]
Since $p (z)$ and $q (z)$ do not vanish on $\abs{z} = 1$, it follows that 
\[
\frac{p^\prime}{p} (z_0) 
\ = \ 
\frac{q^\prime}{q} (z_0)
\]
or still, 
\[
\sum\limits_{i = 1}^n \ 
\frac{1}{z_0 - z_i}
\ = \ 
\sum\limits_{i = 1}^n \ 
\frac{1}{\ds z_0  - \frac{1}{\bar{z}_i} }+ \frac{k}{z_0}
\]
or still, 
\[
\sum\limits_{i = 1}^n \ 
\frac{1}{\ds 1 - \frac{z_i}{z_0}}
\ = \ 
\sum\limits_{i = 1}^n \ 
\frac{1}{\ds 1 - \frac{1}{\bar{z}_i \hsz z_0}} + k.
\]
But
\[
\begin{cases}
\ \ds
\abs{w} < 1 
\implies 
\Reg \frac{1}{1 - w} > \frac{1}{2}
\\[18pt]
\ \ds
\abs{w} > 1 
\implies 
\Reg \frac{1}{1 - w} < \frac{1}{2}
\end{cases}
.
\]
Therefore
\[
\Reg
\Big(
\sum\limits_{i = 1}^n \ 
\frac{1}{\ds 1 - \frac{z_i}{z_0}}
\Big)
\ < \ 
\frac{n}{2}
\]
while
\[
\Reg
\Big(
\sum\limits_{i = 1}^n \ 
\frac{1}{\ds 1 - \frac{1}{\bar{z}_i \hsz z_0}}
\Big)
\ > \ 
\frac{n}{2},
\]
from which the evident contradiction.]
\\[-.25cm]
\end{x}

Let 
\[
P (z) 
\ = \ 
a_0 + a_1 z+ \cdots + a_n z^n
\]
be a real polynomial whose zeros lie in the region 
$\abs{z} \geq 1$.  
Put 
$\zeta = e^{\sqrt{-1} \hsx z}$
$-$then 
\[
\begin{cases}
\ 
P (\zeta) 
\ = \ 
a_0 + a_1 \zeta+ \cdots + a_n \zeta^n
\\[4pt]
\ 
Q (\zeta) 
\ = \ 
a_0\zeta^n + a_1 \zeta^{n-1}+ \cdots + a_n 
\end{cases}
\]
and 
\allowdisplaybreaks
\begin{align*}
&
P (\zeta) + \zeta^n Q (\zeta)  \ = \ 0
\\[11pt]
&
\hspace{1cm}
\implies
\abs{z} = 1 
\qquad (\tcf. \ 30.10)
\\[11pt]
&
\hspace{1cm}
\implies
z \in \R.
\end{align*}
\\[-1cm]

\begin{x}{\small\bf LEMMA} \ 
The trigonometric polynomial 
\[
\sum\limits_{k = 0}^n \ 
a_{n-k} 
\hsy
\cos k z
\]
has real zeros only.
\\[-.5cm]

PROOF \ 
Write
\allowdisplaybreaks
\begin{align*}
\zeta^{-n} 
\hsy 
(P (\zeta) - \zeta^n Q (\zeta) ) \ 
&=\ 
2 a_n + a_{n-1} (\zeta + \zeta^{-1}) + \cdots + a_0 (\zeta^n + \zeta^{-n})
\\[11pt]
&=\ 
2 (a_n + a_{n-1} \cos z + \cdots + a_0 \cos n z)
\\[11pt]
&=\ 
2 \ 
\sum\limits_{k = 0}^n \ 
a_{n-k} 
\hsy
\cos k z.
\end{align*}
\\[-.75cm]
\end{x}


\begin{x}{\small\bf ENESTR\"OM-KAKEYA CRITERION} \ 
Let 
\[
p (z) 
\ = \ 
a_0 + a_1 z + \cdots + a_n z^n, 
\]
where
\[
a_0 \hsx > \hsx a_1 \hsx >  \hsx \cdots \hsx > \hsx a_n > 0.
\]
Then the zeros of $p$ lie in the region $\abs{z} > 1$.
\\[-.5cm]

PROOF \ 
Assuming that $\abs{z} \leq 1$ $(z \neq 1)$, we have
\allowdisplaybreaks
\begin{align*}
&
\big|
\hsx
(1 - z) \hsy 
(a_0 + a_1 z + \cdots + a_n z^n)
\hsx
\big|
\\[11pt]
&
\hspace{1cm}
=\ 
\abs{a_0 - (a_0 - a_1) z - \cdots - (a_{n-1} - a_n) z^n - a_n z^{n+1}}
\\[11pt]
&
\hspace{1cm}
\geq\ 
a_0 - 
\abs{(a_0 - a_1) z + \cdots + (a_{n-1} - a_n) z^n + a_n z^{n+1}}
\\[11pt]
&
\hspace{1cm}
>\ 
a_0 - 
((a_0 - a_1)  + \cdots + (a_{n-1} - a_n) + a_n )
\\[11pt]
&
\hspace{1cm}
=\ 
0.
\end{align*}

[Note: \ 
If instead
\[
a_0 \hsx \geq \hsx a_1 \hsx \geq  \hsx \cdots \hsx \geq \hsx a_n > 0, 
\]
then the zeros of $p$ lie in the region  $\abs{z} \geq 1$.]
\\[-.25cm]
\end{x}

\begin{x}{\small\bf APPLICATION} \ 
If 
\[
0 \hsx < \hsx a_0 \hsx < \hsx a_1 \hsx < \hsx \cdots \hsx < \hsx a_n
\]
and if 
\[
P (z) 
\ = \ 
\sum\limits_{k = 0}^n \ 
a_{n-k} z^k, 
\]
then the zeros of $P$ lie in the region $\abs{z} > 1$, thus the zeros of the trigonometric
polynomial
\[
\sum\limits_{k = 0}^n \ 
a_k \hsy \cos k z
\]
are real (and simple (cf. 30.11)).
\\[-.25cm]
\end{x}

\begin{x}{\small\bf FACT} \ 
For any continuous function $f (t)$ on $[0,1]$, 
\\[-1cm]

\allowdisplaybreaks
\begin{align*}
&
\lim\limits_{n \ra \infty} \ 
\frac{
\ds
\phi
\Big(
\frac{1}{n}
\Big)_{\textcolor{white}{|}}
\hsy
f
\Big(
\frac{1}{n}
\Big)
\hsx + \hsx
\phi
\Big(
\frac{2}{n}
\Big)
\hsy
f
\Big(
\frac{2}{n}
\Big)
\hsx + \hsx
\cdots
\hsx + \hsx
\phi
\Big(
\frac{n-1}{n}
\Big)
\hsy
f
\Big(
\frac{n-1}{n}
\Big)
}
{n}
\ = \ 
\int\limits_0^1 \ 
\phi (t) 
\hsx 
f(t) 
\ \td t.
\end{align*}

PROOF \ 
Given $\varepsilon > 0$, choose $\delta > 0$: 
\[
\int\limits_{1 - \delta}^1 \ 
\phi (t)  
\ \td t
\ < \ 
\varepsilon.
\]
Then
\[
\lim\limits_{n \ra \infty} \ 
\frac{1}{n}
\sum\limits_{k = 1}^{[(1 - \delta) n]} \ 
\phi
\Big(
\frac{k}{n}
\Big)
\hsy
f
\Big(
\frac{k}{n}
\Big)
\ = \ 
\int\limits_0^{1 - \delta} \ 
\phi (t) 
\hsx 
f(t) 
\ \td t.
\]
On the other hand, with 
$
M = \sup\limits_{[0,1]} \ \abs{f}
$, 
we have 
\allowdisplaybreaks
\begin{align*}
\Big|
\hsy
\frac{1}{n} \ 
\sum\limits_{k = [(1 - \delta) n] + 1}^{n-1} \ 
\phi
\Big(
\frac{k}{n}
\Big)
\hsx
f
\Big(
\frac{k}{n}
\Big)
\hsy
\Big|\ 
&\leq 
\frac{M}{n} \ 
\sum\limits_{k = [(1 - \delta) n] + 1}^{n-1} \ 
\phi
\Big(
\frac{k}{n}
\Big)
\\[15pt]
&\leq \ 
M \ 
\int\limits_{1 - \delta}^1 \ 
\phi (t)  
\ \td t
\\[15pt]
&\leq \ 
M \hsx \varepsilon.
\end{align*}
\\[-1.25cm]
\end{x}

With these preliminaries established, the proof of 30.8 is straightforward.
Indeed, for $n = 1, 2, \ldots$, 
\[
0 
\hsx < \hsx
\phi
\Big(
0
\Big)
\hsx < \hsx
\phi
\Big(
\frac{1}{n}
\Big)
\hsx < \hsx
\hsx \cdots \hsx
\hsx < \hsx
\phi
\Big(
\frac{n-1}{n}
\Big),
\]
so a specialization of the preceding generalities implies that the zeros of the trigonometric polynomial
\[
\phi
(
0
)
\hsx + \hsx
\phi
\Big(
\frac{1}{n}
\Big)
\hsy
\cos z
\hsx + \hsx
\cdots
\hsx + \hsx
\phi
\Big(
\frac{n-1}{n}
\Big)
\hsy
\cos (n - 1) z
\]
are real, as are the zeros of the trigonometric polynomial
\[
\phi
(
0
)
\hsx + \hsx
\phi
\Big(
\frac{1}{n}
\Big)
\hsy
\cos \frac{z}{n}
\hsx + \hsx
\cdots
\hsx + \hsx
\phi
\Big(
\frac{n-1}{n}
\Big)
\hsy
\cos
\Big(
\frac{n-1}{n}
\Big) z.
\]
But (cf. 30.15)
\[
\lim\limits_{n \ra \infty} \ 
\frac{1}{n} \ 
\sum\limits_{k = 0}^{n-1} \ 
\phi
\Big(
\frac{k}{n}
\Big)
\hsx
\cos \frac{k}{n} z
\ = \ 
\int\limits_0^1 \ 
\phi (t) 
\hsx 
\cos z \hsy t 
\ \td t, 
\]
the convergence being uniform on compact subsets of $\Cx$, thereby terminating the proof of 30.8.
\\[-.5cm]

[Note: \ 
The zeros of 
\[
\sum\limits_{k = 0}^{n-1} \ 
\phi
\Big(
\frac{k}{n}
\Big)
\hsx
\cos \Big(\frac{k}{n} z\Big)
\]
are not only real but they are also simple (cf. 30.14).  
Still, additional argument is needed in order to conclude that the zeros of 
\[
C (z) 
\ = \ 
\int\limits_0^1 \ 
\phi (t) 
\hsx 
\cos z \hsy t 
\ \td t
\]
are simple (cf. 31.1).]
\\[-.5cm]

\begin{x}{\small\bf REMARK} \ 
Work instead with 
\[
\zeta^{-n} 
\hsy 
(P (\zeta) - 
\zeta^n
Q (\zeta)  )
\]
to see that the trigonometric polynomial
\[
2 \hsy \sqrt{-1} \ 
\sum\limits_{k = 0}^n \ 
a_{n-k} 
\hsy
\sin k z
\]
has real zeros only.  
Pass now to 
\[
\phi
\Big(
\frac{1}{n}
\Big)
\hsy
\sin z
\hsx + \hsx
\cdots
\hsx + \hsx
\phi
\Big(
\frac{n-1}{n}
\Big)
\hsy
\sin (n-1) z
\]
and proceed as  above, the bottom line being that the zeros of the real entire function 
\[
S (z) 
\ = \ 
\int\limits_0^1 \ 
\phi (t) 
\hsx 
\sin z \hsy t 
\ \td t
\]
are real.
\\[-.25cm]
\end{x}

\begin{x}{\small\bf EXAMPLE} \ 
The zeros of 
\[
\frac{\cos z}{z^2}
\hsx 
(\tan z - z)
\ = \ 
\int\limits_0^1 \ 
t
\hsx 
\sin z \hsy t 
\ \td t
\]
are real.
\\[-.5cm]

[Note: \ 
Consequently, 
$\tan z - z$ 
has real zeros only.]
\\[-.25cm]
\end{x}

\begin{x}{\small\bf EXAMPLE} \ 
The zeros of 
\[
\tJ_1 (z) 
\ = \ 
-\tJ_0^\prime (z) 
\ = \ 
\frac{2}{\pi} \ 
\int\limits_0^1 \ 
\frac{t}{\sqrt{1 - t^2}}
\hsx
\sin z \hsy t 
\ \td t
\]
are real (cf 12.33).
\\[-.25cm]
\end{x}

\begin{x}{\small\bf EXAMPLE} \ 
Consider 
\[
\int\limits_0^1 \ 
(1 - t^2)
\hsx
\cos z \hsy t 
\ \td t.
\]
Then its zeros are real (cf. 30.6).
\\[-.5cm]

[Since $1 - t^2$ is decreasing, this is not a special case of 30.8.  
But 
\[
\int\limits_0^1 \ 
(1 - t^2)
\hsx
\cos z \hsy t 
\ \td t
\ = \ 
\frac{2}{z} \ 
\int\limits_0^1 \ 
t
\hsx
\sin z \hsy t 
\ \td t, 
\]
so it is a special case of 30.16.]
\\[-.5cm]

[Note: \ 
In detail, 
\allowdisplaybreaks
\begin{align*}
\int\limits_0^1 \ 
t
\hsx
\sin z \hsy t 
\ \td t \ 
&=\ 
- \frac{1}{2} \ 
\int\limits_0^1 \ 
\sin z \hsy t 
\ \td t (1 - t^2) 
\\[15pt]
&=\ 
- \frac{1}{2} \
(\sin z \hsy t) \hsy (1 - t^2) 
\hsy 
\bigg|_0^1
\ + \ 
\frac{z}{2} \ 
\int\limits_0^1 \ 
\cos z \hsy t (1 - t^2)
\ \td t
\\[15pt]
&=\ 
\frac{z}{2} \ 
\int\limits_0^1 \ 
\cos z \hsy t (1 - t^2)
\ \td t.]
\end{align*}
\\[-1cm]
\end{x}

\begin{x}{\small\bf REMARK} \ 
If in 30.8, the assumption that $\phi (t)$ is 
positive, strictly increasing, and continuous on $[0,1[$ 
is replaced by the assumption that $\phi (t)$ is 
positive, strictly decreasing, and continuous on $[0,1]$ , 
then $C (z)$ may have nonreal zeros.  
\\[-.5cm]

[Consider 
\[
\int\limits_0^1 \ 
e^{- t} 
\hsx 
\cos z \hsy t 
\ \td t
\ = \ 
\frac{(z \sin z - \cos z) + 1}{e \hsy (z^2 + 1)}
\hsx.]
\]
\\[-.25cm]
\end{x}

\chapter{
$\boldsymbol{\S}$\textbf{31}.\quad  TRANSFORM THEORY: \ SENIOR GRADE}
\setlength\parindent{2em}
\setcounter{theoremn}{0}
\renewcommand{\thepage}{\S31-\arabic{page}}

\qquad 
The following result supercedes 30.8.
\\[-.25cm]

\begin{x}{\small\bf THEOREM} \ 
If $\phi \in \Lp^1 [0,1]$ 
is positive and increasing on 
$]0,1[$, 
then the zeros of 
\[
C (z) 
\ = \ 
\int\limits_0^1 \ 
\phi (t) 
\hsx
\cos z \hsy t 
\ \td t
\]
are real and simple.  
Furthermore, the positive zeros of 
$C(z)$ 
lie in the intervals
\[
\Big]
\frac{\pi}{2}
, 
\frac{3 \pi}{2}
\Big[
\hsx
,
\quad 
\Big]
\frac{3 \pi}{2}
, 
\frac{5 \pi}{2}
\Big[
\hsx
,
\quad
\Big]
\frac{5 \pi}{2}
, 
\frac{7 \pi}{2}
\Big[
\hsx
,
\quad
\ldots
\]
and only in these intervals.  
Finally, each of these intervals contains exactly one zero of 
$C(z)$.
\\[-.5cm]

[Note: \ 
$C (z)$ is even, hence $C (z_0) = 0$ iff $C (-z_0) = 0$.]
\\[-.5cm]

The proof is spelled out in the lines below.
\\[-.25cm]

\un{Step 1:}
\[
C\Big(\frac{\pi}{2}\Big) 
\ = \ 
\int\limits_0^1 \ 
\phi (t) 
\hsx
\cos \frac{\pi}{2} \hsx t 
\ \td t
\ > \ 
0.
\]

\un{Step 2:}
\\[-.25cm]

\qquad \textbullet \quad
$
\ds
C\Big(\frac{\pi}{2} + 2 \hsy \pi \hsy n\Big) 
\ > \  
0
\qquad
(n = 1, 2, \ldots).
$
\\[-.25cm]

[We have
\allowdisplaybreaks
\begin{align*}
\int\limits_0^1 \ 
\phi (t) 
&
\cos (2 \hsy \pi \hsy n + \frac{\pi}{2}) \hsy t
\ \td t \ 
\\[15pt]
&=\ 
\int\limits_0^{1/(4n + 1)} \ 
\phi (t) 
\cos (4 \hsy  n + 1) \hsy \frac{\pi}{2}\hsy t
\ \td t \ 
\ + \ 
\sum\limits_{k = 0}^n \ 
\int\limits_{\frac{4k + 1}{4n + 1}}^{\frac{4k + 5}{4n + 1}} \ 
\phi (t) 
\cos (4 n + 1) \frac{\pi}{2} \hsy t
\ \td t \ 
\\[15pt]
&\geq \ 
\int\limits_0^{1/(4n + 1)} \ 
\phi (t) 
\hsx
\cos (4 n + 1) \hsy \frac{\pi}{2} \hsy t
\ \td t 
\\[15pt]
&> \ 
0.]
\end{align*}
\\[-.25cm]

\qquad \textbullet \quad
$
\ds
C\Big(\frac{3 \pi}{2} + 2 \hsy \pi \hsy n\Big) 
\ < \  
0
\qquad
(n = 0, 1, 2, \ldots).
$
\\

[We have
\allowdisplaybreaks
\begin{align*}
\int\limits_0^1 \ 
\phi (t) 
&
\cos (4 n + 3) \hsy  \frac{\pi}{2} \hsy t
\ \td t \ 
\\[15pt]
&=\ 
\int\limits_0^{2/(4n + 3)} \ 
\phi (t) 
\hsx
\cos (4 n + 3) \hsy \frac{\pi}{2} \hsy t
\ \td t \ 
\ + \ 
\int\limits_{2/(4n + 3)}^{3/(4n + 3)} \ 
\phi (t) 
\hsx
\cos (4 n + 3) \hsy \frac{\pi}{2} \hsy t
\ \td t \ 
\\[15pt]
&\hspace{3cm}
\ + \ 
\sum\limits_{k = 0}^n \ 
\int\limits_{\frac{4k + 3}{4n + 3}}^{\frac{4k + 7}{4n + 3}} \ 
\phi (t) 
\cos (4 n + 3) \hsy \frac{\pi}{2} \hsy t
\ \td t \ 
\\[15pt]
&\leq \ 
\int\limits_{2/(4n + 3)}^{3/(4n + 3)} \ 
\phi (t) 
\hsx
\cos (4 n + 3) \hsy \frac{\pi}{2} \hsy t
\ \td t 
\\[15pt]
&< \ 
0.]
\end{align*}

So far then
\[
C\Big(\frac{\pi}{2}\Big) \ > 0, 
\quad 
C\Big(\frac{3 \pi}{2}\Big) \ < 0, 
\quad 
C\Big(\frac{5 \pi}{2}\Big) \ > 0, 
\quad 
C\Big(\frac{7 \pi}{2}\Big) \ < 0, 
\quad 
\ldots, 
\]
which implies that each of the intervals
\[
\Big]
\frac{\pi}{2}
, 
\hsx
\frac{3 \pi}{2}
\Big[
\hsx
,
\quad 
\Big]
\frac{3 \pi}{2}
, 
\hsx
\frac{5 \pi}{2}
\Big[
\hsx
,
\quad
\Big]
\frac{5 \pi}{2}
, 
\hsx
\frac{7 \pi}{2}
\Big[
\hsx
,
\quad
\ldots
\]
contains at least one zero of 
$C(z)$, 
as do the intervals symmetric to them.  
The objective now is to show that any such interval contains but one zero of 
$C(z)$, 
that said zero is simple, and that there are no other zeros.
\\[-.5cm]
\end{x}

To move forward, assume without loss of generality that 
$C (0) = 1$.
\\[-.25cm]

\begin{x}{\small\bf RAPPEL} \ 
\[
\int\limits_0^r \ 
\frac{n (t)}{t}
\ \td t
\ = \ 
\frac{1}{2 \pi}\ 
\int\limits_0^{2 \pi} \ 
\log 
\big|
C (r \hsy e^{\sqrt{-1} \hsx \theta})
\big|
\ \td \theta
\qquad (\tcf. \ 27.36).
\]

Let $n^* (t)$ denote the number of points 
$
\ 
\ds 
\pm 
\Big(
\frac{\pi}{2} + \pi n
\Big)
\ 
$
$(n = 1, 2, \ldots)$ in the interval
$]-t, t[$ $(t > 0)$, thus 
$n^* (t) = 0$ 
for 
$
\ds 
\abs{t} < 
\frac{3 \pi}{2}
$ 
and 
\[
n^* (t) 
\ = \ 2 k 
\quad \text{if} \quad 
\frac{\pi}{2} + \pi k 
\ < \ 
t 
\ < \ 
\frac{\pi}{2} + \pi (k + 1)
\qquad (k = 1, 2, \ldots).
\]
\\[-1.5cm]
\end{x}

To derive a contradiction, suppose that 
$C (z_0) = 0$ ($\implies C (-z_0) = 0$), 
where $z_0$ is either not in one of the intervals above or is a multiple zero of one thereof.  
Choose $K > 0$: 
\[
n (t) 
\ \geq \ 
n^* (t) 
\ \ (0 < t < K), 
\quad
n (t) 
\ \geq \ 
n^* (t)  + 2
\ \ (t >  K).
\]
\\[-.75cm]

\un{Step 3:} \ 
Take 
$
\ds 
r = \pi n + 
\frac{3 \pi}{2}
$ 
$-$then
\allowdisplaybreaks
\begin{align*}
\int\limits_0^r \ 
\frac{n (t)}{t}
\ \td t \ 
&\geq 
\sum\limits_{k = 1}^n \ 
(2 k + 2) 
\ 
\int\limits_{\frac{\pi}{2} + \pi k}^{\frac{\pi}{2} + \pi (k+1)} \ 
\ \frac{\td t}{t} \ 
+ \ 
\tO (1)
\\[15pt]
&=\ 
2\ 
\sum\limits_{k = 1}^n \ 
(k + 1) 
\hsx 
\log \Big(1 + \frac{1}{k + \frac{1}{2}}\Big)
\ + \ 
\tO (1)
\\[15pt]
&=\ 
2\ 
\sum\limits_{k = 1}^n \ 
(k + 1) 
\hsx
\Big(
1 
+
\frac{1}{k + \frac{1}{2}}
-
\frac{1}{2 \hsy(k + \frac{1}{2})^2}
\Big)
\ + \ 
\tO (1)
\\[15pt]
&=\ 
2 \ 
\sum\limits_{k = 1}^n \ 
1
\ + \ 
\sum\limits_{k = 1}^n \ 
\frac{1}{k + \frac{1}{2}}
\ - \ 
\sum\limits_{k = 1}^n \ 
\frac{k+1}{(k + \frac{1}{2})^2}
\ + \ 
\tO (1)
\\[15pt]
&=\ 
2 n 
\ + \ 
\tO (1)
\\[15pt]
&=\ 
2 \hsy \frac{r}{\pi}
\ + \ 
\tO (1).
\end{align*}

\un{Step 4:} \ 
Since 
\[
C (x) \ra 0 
\quad \text{as} \quad 
x \ra \pm \infty
\]
and since the exponential type of $C (z)$ is $\leq 1$, 
\[
\frac
{
\big|
C \big(r \hsy e^{\sqrt{-1} \hsx \theta}\big)
\big|
}
{e^{\abs{r \sin \theta}}}
\ \ra \ 0
\qquad (r \ra \infty)
\]
uniformly in $\theta$.  
Therefore
\allowdisplaybreaks
\begin{align*}
\frac{1}{2 \hsy \pi} \ 
\int\limits_0^{2 \hsy \pi} \ 
\log 
\big|
C \big(r \hsy e^{\sqrt{-1} \hsx \theta}\big)
\big|
\ \td \theta
&=\ 
\frac{1}{2 \hsy \pi} \ 
\int\limits_0^{2 \hsy \pi} \ 
\log 
\Big|
\frac
{
C \big(r \hsy e^{\sqrt{-1} \hsx \theta}\big)
}
{e^{\abs{r \sin \theta}}}
\ \cdot \ 
e^{\abs{r \sin \theta}}
\Big|
\ \td \theta
\\[15pt]
&=\ 
\frac{1}{2 \hsy \pi} \ 
\int\limits_0^{2 \hsy \pi} \ 
\log 
\Big|
\frac
{
C \big(r \hsy e^{\sqrt{-1} \hsx \theta}\big)
}
{e^{\abs{r \sin \theta}}}
\Big|
\ \td \theta
\ + \ 
\frac{1}{2 \hsy \pi} \ 
\int\limits_0^{2 \hsy \pi} \ 
\abs{r \sin \theta}
\ \td \theta
\\[15pt]
&\leq\ 
\log \txo (1) + 2 \hsx \frac{r}{\pi}.
\end{align*}

\un{Step 5:} \ 
Combine the data: 
\allowdisplaybreaks
\begin{align*}
\log \txo (1) + 2 \hsx \frac{r}{\pi} \ 
&\geq \ 
\frac{1}{2 \hsy \pi} \ 
\int\limits_0^{2 \hsy \pi} \ 
\log 
\big|
C \big(r \hsy e^{\sqrt{-1} \hsx \theta}\big)
\big|
\ \td \theta
\\[15pt]
&=\ 
\int\limits_0^r \ 
\frac{n (t)}{t}
\ \td t \ 
\\[15pt]
&\geq \
2 \hsx \frac{r}{\pi}
\ + \ 
\tO (1)
\end{align*}
\qquad 
$\implies$
\[
\log \txo (1) 
\ \geq \ 
\tO (1), 
\]
an impossibility.
\\

\begin{x}{\small\bf THEOREM} \ 
If $\phi \in \Lp^1 [0,1]$ 
is positive and increasing on 
$]0,1[$ 
and is not exceptional (cf. 29.14), then the zeros of 
\[
S (z) 
\ = \ 
\int\limits_0^1 \ 
\phi (t) 
\hsx
\sin z \hsy t 
\ \td t
\]
are real and simple.  
Furthermore, the positive zeros of $S (z)$ lie in the intervals
\[
]\pi, \hsx 2 \pi[
\hsx ,
\quad 
]2 \pi, \hsx 3 \pi[
\hsx ,
\quad
]3 \pi, \hsx 4 \pi[
\hsx ,
\quad
\ldots
\]
and only in these intervals.  
Finally, each of these intervals contains exactly one zero of $S (z)$.
\\[-.5cm]

[Note: \ 
$S (z)$ is odd, hence $S (z_0) = 0$ iff $S (-z_0) = 0$.]
\\[-.25cm]
\end{x}

The proof is spelled out in the lines below.
\\

\un{Step 1:}
\[
S (0) 
\ = \ 
\int\limits_0^1 \ 
\phi (t) \hsy \sin 0 \hsx t
\ \td t
\ = \ 
0.
\]
And
\[
S^\prime (z) 
\ = \ 
\int\limits_0^1 \ 
\phi (t) \hsy t \hsy \cos z \hsy t
\ \td t
\]
\qquad
\qquad 
$\implies$
\allowdisplaybreaks
\begin{align*}
S^\prime (0) \ 
&= \ 
\int\limits_0^1 \ 
\phi (t) \hsy t \hsy \cos 0 \hsx t
\ \td t 
\\[15pt]
&= \ 
\int\limits_0^1 \ 
\phi (t) \hsy t 
\ \td t 
\\[15pt]
&> \  
0.
\end{align*}
Therefore 0 is  a simple zero of $S (z)$.
\\[-.75cm]

\un{Step 2:}
\[
S (\pi) 
\ = \ 
\int\limits_0^1 \ 
\phi (t) \hsy \sin \pi t
\ \td t
\ > \ 
0.
\]

\un{Step 3:}
\\[-.25cm]

\qquad \textbullet \quad
$S (\pi + 2 \hsy \pi \hsy n) > 0
\qquad (n = 1, 2, \ldots).$
\\[-.25cm]

[We have
\allowdisplaybreaks
\begin{align*}
\int\limits_0^1 \ 
&
\phi (t) \hsy \sin (2 n + 1) \pi \hsy t
\ \td t \ 
\\[15pt]
&=\ 
\int\limits_0^{1 / (2 n + 1)} \ 
\phi (t) \hsy \sin (2 n + 1) \pi \hsy t
\ \td t \ 
\ + \ 
\sum\limits_{k = 0}^{n-1} \ 
\int\limits_{\frac{2 k + 1}{2 n + 1}}^{\frac{2 k + 3}{2 n + 1}}
\phi (t) \hsy \sin (2 n + 1) \pi \hsy t
\ \td t \ 
\\[15pt]
&\geq \ 
\int\limits_0^{1 / (2 n + 1)} \ 
\phi (t) \hsy \sin (2 n + 1) \pi \hsy t
\ \td t \ 
\\[15pt]
&> \ 
0.]
\end{align*}

\qquad \textbullet \quad
$S (2 \hsy \pi \hsy n) \ < \  0
\qquad (n = 1, 2, \ldots).$
\\[-.25cm]

[We have
\allowdisplaybreaks
\begin{align*}
\int\limits_0^1 \ 
&
\phi (t) \hsy \sin 2 \hsy \pi \hsy n \hsy t
\ \td t \ 
\\[15pt]
&=\ 
\sum\limits_{k = 0}^{n-1} \ 
\int\limits_{\frac{k}{n}}^{\frac{k + 1}{n}}
\phi (t) \hsy \sin 2 \hsy \pi \hsy n \hsy t
\ \td t \ 
\\[15pt]
&=\ 
\sum\limits_{k = 0}^{n-1} \ 
\int\limits_0^{1/n}
\phi \Big(t + \frac{k}{n}\Big) 
\hsy 
\sin 2 \hsy \pi \hsy n \hsy t
\ \td t \ 
\\[15pt]
&=\ 
\sum\limits_{k = 0}^{n-1} \ 
\int\limits_0^{1/2 \hsy n}
\Big(
\phi \Big(t + \frac{k}{n}\Big) 
\hsx - \hsx 
\phi \Big(\frac{k+1}{n} - t\Big) 
\Big)
\hsy 
\sin 2 \hsy \pi \hsy n \hsy t
\ \td t \ 
\\[15pt]
&< \ 
0.]
\end{align*}

[Note: \ 
The function 
$\sin 2 \pi n t$ is positive on 
$
\ds
\Big]0, \hsx \frac{1}{2 n}\Big[
$ 
and
\[
\phi \Big(t + \frac{k}{n}\Big) 
\hsx - \hsx 
\phi \Big(\frac{k+1}{n} - t\Big) 
\qquad 
\Big(0 < t < \frac{1}{2 n}\Big) 
\]
is nonpositive and increasing, thus a priori
\[
\sum\limits_{k = 0}^{n-1} \ 
\int\limits_0^{1/2 \hsy n}
\Big(
\phi \Big(t + \frac{k}{n}\Big) 
\hsx - \hsx 
\phi \Big(\frac{k+1}{n} - t\Big) 
\Big)
\hsy 
\sin 2 \hsy \pi \hsy n \hsy t
\ \td t \ 
\ \leq \ 
0, 
\]
with equality only if $\forall \ k$
\[
\phi \Big(t + \frac{k}{n}\Big) 
\hsx - \hsx 
\phi \Big(\frac{k+1}{n} - t\Big) 
\ = \  
0
\]
\begin{spacing}{1.75}
\noindent
almost everywhere and this means zero on 
$
\ds
\Big]0, \frac{1}{2 n}\Big[
$ 
(if negative anywhere on 
$
\ds
\Big]0, \frac{1}{2 n}\Big[
$, 
then it is negative from there to the left giving a negative integral), 
hence $\phi (t)$ would be a constant in each of the intervals 
$
\ 
\ds
\frac{k}{n} < t < \frac{k+1}{n}
\ 
$ 
$(k = 0, \ldots, n-1)$, 
a scenario excluded by the assumption 
$\phi \notin E (1, 0)$.]
\\[-1cm]
\end{spacing}

So far then 
\[
S (\pi) > 0, 
\quad
S (2 \pi) < 0, 
\quad
S (3 \pi) > 0, 
\quad
S (4 \pi) < 0, 
\ldots
\]
which implies that each of the intervals 
\[
]\pi, \hsx 2 \pi[
\hsx ,
\quad 
]2 \pi, \hsx 3 \pi[
\hsx ,
\quad
]3 \pi, \hsx 4 \pi[
\hsx ,
\quad
\ldots
\]
contains at least one zero from 
$S(z)$, 
as do the intervals symmetric to them (recall too that 0 is a simple zero of $S (z)$).  
The remaining details are similar to those figuring in 31.1 and will be omitted.
\\

\begin{x}{\small\bf LEMMA} \ 
If $\phi \in \Lp^1 [0,1]$ 
is positive and increasing on 
$]0,1[$ 
and if $\phi \notin E (1, 0)$, then $C (z)$ and $S (z)$ have no common zeros.
\\[-.5cm]

PROOF \ 
The zeros of 
\allowdisplaybreaks
\begin{align*}
f (z) \ 
&=\ 
\int\limits_0^1 \ 
\phi (t) 
\hsy 
e^{\sqrt{-1} \hsx z \hsy t}
\ \td t
\\[15pt]
&=\ 
C (z) + \sqrt{-1} \hsx S (z)
\end{align*}
lie in the open upper half-plane (cf. 29.16).  
On the other hand, as has been seen above, 
the zeros of $C (z)$ and $S (z)$ are real, so 
\[
\begin{cases}
\ 
C (x_0) = 0
\\[8pt]
\ 
S (x_0) = 0
\end{cases}
\implies f (x_0) = 0,
\]
which cannot be.]
\\[-.25cm]
\end{x}

\chapter{
$\boldsymbol{\S}$\textbf{32}.\quad  APPLICATION OF INTERPOLATION}
\setlength\parindent{2em}
\setcounter{theoremn}{0}
\renewcommand{\thepage}{\S32-\arabic{page}}

\qquad 
Let $f \in \sB_0(A)$ and assume that $f$ is not a constant, hence $T (f) > 0$.
\\[-.25cm]

\begin{x}{\small\bf RAPPEL} \ 
(cf. 17.22) $\forall$ real $x$, 
\[
f^\prime (x) 
\ = \ 
\frac{4 T (f) }{\pi^2} 
\ 
\sum\limits_{k = -\infty}^\infty \ 
(-1)^k
\hsx
\frac{1}{(2 k + 1)^2} 
\hsx
f\Big(x + \frac{2k + 1}{2 T (f)} \hsx \pi\Big),
\]
the convergence being uniform on compact subsets of $\R$.
\\[-.25cm]
\end{x}

\begin{x}{\small\bf THEOREM} \ 
$\forall \ x, \alpha \in \R$, there is an expansion
\[
\sin \alpha
\hsx \cdot \hsx f^\prime (x) 
\hsx - \hsx
A \cos \alpha
\hsx \cdot \hsx f (x)
\ = \ 
A \sin^2 \alpha
\ 
\sum\limits_{k = -\infty}^\infty \ 
\frac{(-1)^{k-1}}{(\alpha - k \pi)^2} 
\
f\Big(x + \frac{k \pi - \alpha}{A}\Big),
\]
\\[-.25cm]
the convergence being uniform on compact subsets of $\R$.
\\[-.5cm]

[Note: \ 
Replace $k$ by $k + 1$ and take 
$
\ds 
\alpha = \frac{\pi}{2}
$, 
$A = T (f)$ 
to recover 31.1.]
\\[-.5cm]

PROOF \ 
Write
\[
f (z) 
\ = \ 
f (0) 
\hsx + \hsx 
\frac{z}{\sqrt{2 \hsy \pi}} \ 
\int\limits_{-A}^A \ 
\phi (t) 
\hsx 
e^{\sqrt{-1} \hsx z \hsy t}
\ \td t
\]
for some $\phi \in \Lp^2 [-A, A]$ (cf. 22.8), so
\allowdisplaybreaks
\begin{align*}
\sin \alpha
\hsx 
&
\cdot \hsx f^\prime (x) 
\hsx 
- 
\hsx
A \cos \alpha
\hsx \cdot \hsx f (x)
\\[11pt]
&=\ 
- A \cos \alpha 
\hsx \cdot \hsx 
f^\prime (0) 
\\[15pt]
&
\hspace{1cm}
+ \ 
\frac{1}{\sqrt{2 \hsy \pi}} \ 
\int\limits_{-A}^A \ 
\phi (t) 
\hsx 
\frac{\partial}{\partial t}
\Big(
e^{\sqrt{-1} \hsx x \hsy t}
(t \sin \alpha + \sqrt{-1} \hsx A \cos \alpha)
\Big)
\ \td t.
\end{align*}
Now develop
\[
-
\sqrt{-1} \hsx
e^{\sqrt{-1} \hsx\frac{\alpha}{A} \hsy t}_{\textcolor{white}{\bigg|}}
\hsy
(t \sin \alpha + \sqrt{-1} \hsx A \cos \alpha)
\]
into a Fourier series: 
\[
A \sin^2 \alpha
\ 
\sum\limits_{k = -\infty}^\infty \ 
(-1)^k
\ 
\frac{1}{(\alpha - k \pi)^2}
\ 
e^{\frac{\sqrt{-1} \hsx k \hsy \pi}{A} \hsy t}
\] 
$\implies$
\allowdisplaybreaks
\begin{align*}
\sin \alpha
\hsx 
&
\cdot \hsx
f^\prime (x) 
\hsx - \hsx
A \cos \alpha
\hsx \cdot \hsx
f (x) \ 
\\[15pt]
&=\ 
- 
A \cos \alpha
\hsx \cdot \hsx
f (0)
\\[15pt]
&
\quad
+ 
\frac{\sqrt{-1} \hsx A \hsy \sin^2 \alpha}{\sqrt{2 \hsy \pi}} 
\int\limits_{-A}^A 
\phi (t)  
\frac{\partial}{\partial t}
\Big(
\sum\limits_{k = -\infty}^\infty \ 
\frac{(-1)^k}{(\alpha - k \pi)^2}
\exp \Big(\sqrt{-1} \hsx t (x + \frac{k \hsy \pi - \alpha}{A}
\Big)
\Big)
\Big)
\ \td t
\\[15pt]
&=\ 
- 
A \cos \alpha
\hsx \cdot \hsx
f (0)
\\[15pt]
&
\hspace{1cm}
\hsx - \hsx
A \sin^2 \alpha
\ 
\sum\limits_{k = -\infty}^\infty \ 
\frac{(-1)^k}{(\alpha - k \pi)^2}
\ 
\big(
f
\Big(
x + \frac{k \pi - \alpha}{A}
\Big)
- f (0)
\big)
\\[15pt]
&=\ 
A \sin^2 \alpha
\ 
\sum\limits_{k = -\infty}^\infty \ 
\frac{(-1)^{k-1}}{(\alpha - k \pi)^2}
\ 
f
\Big(
x + \frac{k \pi - \alpha}{A}
\Big),
\end{align*}
since
\[
\sum\limits_{k = -\infty}^\infty \ 
\frac{(-1)^k}{(\alpha - k \pi)^2}
\ = \ 
- 
\frac{\td}{\td \alpha} 
\hsx
\frac{1}{\sin \alpha} 
\ = \ 
\frac{\cos \alpha}{\sin^2 \alpha}.
\]
\\[-.75cm]
\end{x}

\begin{x}{\small\bf APPLICATION} \ 
$\forall B \in \R$,
\allowdisplaybreaks
\begin{align*}
\sin A (x - B) 
&
\hsx \cdot \hsx
f^\prime (x) 
\hsx - \hsx
A \cos A (x - B) 
\hsx \cdot \hsx
f (x) 
\\[15pt]
&=\ 
A \sin^2 A (x - B) 
\ 
\sum\limits_{k = -\infty}^\infty \ 
\frac{(-1)^{k-1}}{(A (x - B) - k \pi)^2}
\
f
\Big(
\frac{k \hsy \pi}{A} + B
\Big).
\end{align*}


[Replace $\alpha$ by $A (x - B)$ in 32.2.]
\\[-.25cm]
\end{x}

\qquad
{\small\bf \un{N.B.}} \ 
If 
$
\ 
\ds
f \Big(\frac{k \hsy \pi}{A} + B\Big) = 0
$ 
$\forall \ k$, 
then
\[
f (x) 
\ = \ 
C \sin A (x - B) 
\qquad (C \neq 0)
\]
and its zeros are at the points 
$
\ds
\frac{k \hsy \pi}{A} + B
$.
\\[-.25cm]

\begin{x}{\small\bf NOTATION} \ 
$\R \sB_0 (A)$ is the subset of $\sB_0 (A)$ consisting of those nonconstant $f$ which are real on the real axis.
\\[-.25cm]
\end{x}

\begin{x}{\small\bf DEFINITION} \ 
Let $f \in \R \sB_0 (A)$ $-$then $f$ is \un{standard} of level $B$ if $\exists \ n = 0$ 
or 1 and $B \in \R$ such that $\forall \ k \in \Z$, 
\[
(-1)^{n+k}
\hsx
f
\Big(
\frac{k \hsy  \pi}{A} + B
\Big)
\ \geq \ 
0.
\]

[Note: \ 
If $f$ is standard of level $B$, then $-f$ is standard of level $B$.]
\\[-.25cm]
\end{x}

\begin{x}{\small\bf EXAMPLE} \ 
Take $A = 1$, $B = 0$ $-$then if $n = 0$, 

\[
\ldots, 
f(- 2 \pi) \geq 0, \
f(- \pi) \leq 0, \
f(0) \geq 0, \
f(\pi) \leq 0, \
f( 2 \pi) \geq 0, \
\ldots,
\]
with a reversal of signs if $n = 1$.
\\[-.25cm]
\end{x}

\begin{x}{\small\bf EXAMPLE} \ 
Take $A = 1$, 
$
\ds
B = \frac{\pi}{2}$ $-$then if $n = 0$, 
\allowdisplaybreaks
\begin{align*}
\ldots, 
f
\Big(
-
\frac{5 \hsy \pi}{2}
\Big)
\ \leq \ 
0, 
\
f
\Big(
-
\frac{3 \hsy \pi}{2}
\Big)
\ \geq \ 
0, 
\
&
f
\Big(
-
\frac{\pi}{2}
\Big)
\ \leq \ 
0, 
\
f
\Big(
\frac{\pi}{2}
\Big)
\ \geq \ 
0, 
\
\\[11pt]
&
f
\Big(
\frac{3 \hsy \pi}{2}
\Big)
\ \leq \ 
0, 
\
f
\Big(
\frac{5 \hsy \pi}{2}
\Big)
\ \geq \ 
0, 
\hsx
\ldots, 
\end{align*}
with a reversal of signs if $n = 1$.
\\[-.25cm]
\end{x}

\begin{x}{\small\bf LEMMA} \ 
If $f \in \R \sB_0 (A)$ is standard of level $B$, then $\forall \ x \in \R$, 
\allowdisplaybreaks
\begin{align*}
\sin A (x 
&
- B) 
\hsx \cdot \hsx 
f^\prime (x) 
\hsx - \hsx A
\cos A (x - B)
\hsx \cdot \hsx 
f (x) 
\\[15pt]
&=\ 
(-1)^{n - 1} \ 
A \sin^2 A (x - B) 
\ 
\sum\limits_{k = -\infty}^\infty \ 
\frac{1}{(A (x - B) - k \pi)^2}
\
\Big|
f
\Big(
\frac{k \hsy \pi}{A} + B
\Big)
\Big|.
\end{align*}
\\[-1.25cm]
\end{x}

\begin{x}{\small\bf THEOREM} \ 
If $f \in \R \sB_0 (A)$ is standard of level $B$, then $\forall \ p \in \Z$, 
the \un{ambient interval}
\[
I_p
\ = \ 
\Big]
\frac{(p-1)\hsy \pi}{A} + B, 
\frac{p \hsy \pi}{A} + B
\Big[
\]
contains at most one zero of $f$ and if there is one, then it must be simple.
\\[-.5cm]

\begin{spacing}{1.75}
PROOF \ 
Suppose that for some $p \in \Z$, 
$f (x_0) = 0$ $(x_0 \in I_p)$ $-$then $\exists \ k \in \Z$ such that 
$
\ds
f \Big(\frac{k \hsy \pi}{A} + B\Big) \neq 0
$, 
hence
\\[-1cm]
\end{spacing}
\[
\sin A (x_0 - B) 
\hsx \cdot \hsx 
f^\prime (x_0) 
\ = \ 
(-1)^{n-1} 
A \sin^2 A (x_0 - B) 
\hsy
M (x_0) 
\qquad (M (x_0) > 0)
\]
\qquad 
$\implies$
\allowdisplaybreaks
\begin{align*}
f^\prime (x_0) \ 
&=\ 
(-1)^{n-1} 
A \sin A (x_0 - B) 
\hsy
M (x_0) 
\\[11pt]
&=\ 
(-1)^{n-1} 
\hsy
(-1)^{p-1} 
A
\hsy
\abs{\sin A (x_0 - B) }
\hsy
M (x_0) 
\end{align*}

\qquad 
$\implies$
\[
(-1)^{n+p} 
\hsx
f^\prime (x_0) 
\ > \ 
0,
\]
\begin{spacing}{1.5}
\noindent
which implies that $x_0$ is simple.  
If now 
$f(x_1) = 0$, 
$f(x_2) = 0$ 
with $x_1 < x_2$ 
and 
$f (x) \neq 0$ $(x_1 < x < x_2)$, 
then we shall arrive at a contradiction by showing that there would be another zero of $f$ between $x_1$ and $x_2$.  
To see this, choose a small $h > 0$ with the property that 
$f (x)$ and $f^\prime (x)$ 
have the same sign in 
$]x_1, x_1 + h[$ 
and opposite signs in 
$]x_2 - h, x_2[$ 
$(\implies x_1 + h < x_2 - h)$.
\\[-.25cm]
\end{spacing}

\qquad \textbullet \quad
\un{$n + p$ \ even}: \ 
Therefore 
$f^\prime (x_1) > 0$, 
$f^\prime (x_2) > 0$ 
and it can be assumed that 
$f^\prime (x)$ 
is positive in 
$]x_1, x_1 + h[$ 
and 
$]x_2 - h, x_2[$.  
But then
\[
\begin{cases}
\ 
x_1 < x < x_1 + h \implies f (x) > 0
\\[4pt]
\ 
x_2 - h < x < x_2 \implies f (x) < 0
\end{cases}
.
\]
\\[-.25cm]

\qquad \textbullet \quad
\un{$n + p$ \ odd}: \ 
Therefore 
$f^\prime (x_1) < 0$, 
$f^\prime (x_2) < 0$ 
and it can be assumed that 
$f^\prime (x)$ 
is negative in 
$]x_1, x_1 + h[$ 
\ and \ 
$]x_2 - h, x_2[$.  
But then
\[
\begin{cases}
\ 
x_1 < x < x_1 + h \implies f (x) < 0
\\[4pt]
\ 
x_2 - h < x < x_2 \implies f (x) > 0
\end{cases}
.
\]
\\[-1.cm]
\end{x}

\begin{x}{\small\bf LEMMA} \ 
If $f \in \R \sB_0 (A)$ is standard of level $B$, then 
\[
\sup\limits_{x \in \R} x^2 \abs{f (x)}
\ = \ 
\infty.
\]

PROOF  \ 
Assuming this is false, let
\[
g (z) 
\ = \ 
f (z) 
\hsy 
(z - x_0)^2 
\qquad 
\Big(x_0 \in I_1 
\hsx = \ 
\Big]
B, \frac{\pi}{A} + B
\Big[\hsx\Big).
\]
Then $g \in \R \sB_0 (A)$ is standard of level $B$.  
But $x_0$ is a zero of $g$ of multiplicity $\geq 2$, an impossibility (cf. 32.9).
\\[-.25cm]
\end{x}

\begin{x}{\small\bf THEOREM} \ 
If $f \in \R \sB_0 (A)$ is standard of level $B$, then all the zeros of $f$ are real.
\\[-.5cm]

PROOF  \ 
Suppose that $f(z_0) = 0$ for some $z_0 \in \Cx - \R$.  
Since $f$ is real, $f (\bar{z}_0) = 0$ and the function 
\[
g (z) 
\ = \ 
\frac{f (z)}{(z - z_0) (z - \bar{z}_0)}
\]
belongs to $\R \sB_0 (A)$.  
As such, it is standard of level $B$ and 
\[
\sup\limits_{x \in \R} x^2 \abs{g (x)}
\ < \ 
\infty,
\]
which contradicts 32.10.
\\[-.25cm]
\end{x}

\begin{x}{\small\bf EXAMPLE} \ 
Given $\phi \in \Lp^1 [0,1]$ real $\not\equiv 0$, let

\[
C (z) 
\ = \ 
\int\limits_0^1 \
\phi (t) 
\hsx 
\cos z \hsy t 
\ \td t.
\]
Then $C \in \R \sB_0 (1)$.  
Assume: \ 
$\forall \ k \in \Z$, 
\[
(-1)^k 
\hsy
C( k \hsy \pi) 
\ > \ 
0.
\]
Then all the zeros of $C$ are real and each ambient interval $I_p$ contains a single zero and it is simple.
\\[-.25cm]
\end{x}

We have yet to examine what happens at the endpoints of an $I_p$.
\\[-.25cm]

\begin{x}{\small\bf THEOREM} \ 
If $f \in \R \sB_0 (A)$ is standard of level $B$, and if for some $p \in \Z$, 
\[
f \Big(\frac{p\hsy \pi}{A} + B \Big)
\ = \ 
0,
\]
then 
\[
x_p 
\ \equiv \ 
\frac{p\hsy \pi}{A} + B
\]
is a zero of multiplicity $\leq 2$ and $f$ cannot have zeros in both ambient intervals 
$I_p$ 
and 
$I_{p+1}$.  
Moreover, if $x_p$ is a zero of multiplicity 2, then 
\[
(-1)^{n+p} f^{\prime\prime} (x_p)
\ < \ 
0
\]
\\[-1cm]
and
\[
(-1)^{n+p} f (x) 
\ < \ 
0
\qquad 
(x \in I_p \cup I_{p+1}),
\]
while if 
$x_{p-1}$ 
(or $x_{p+1}$) 
is a zero, then 
$x_{p-1}$ 
(or $x_{p+1}$) 
must be simple.
\\[-.5cm]

PROOF\ 
This is elementary, albeit detailed.
\\[-.25cm]

\qquad \textbullet \quad
If 
$f (x_p) = 0$, 
$f^\prime (x_p) = 0$, 
then
\[
(-1)^{n+p} f^{\prime\prime} (x_p)
\ < \ 
0, 
\]
hence in particular, $x_p$ is a zero of multiplicity $\leq 2$.  
Thus let
\[
g (z) 
\ = \ 
\frac{f (z)}{(z - x_p)^2}\hsx .
\]
Then 
$g \in \R \sB_0 (A)$ 
and we claim that $g$ is standard of level $B$ if 
\[
(-1)^{n+p} f^{\prime\prime} (x_p)
\ \geq \ 
0.
\]
For it is clear that 
\[
(-1)^{n+k} 
g \Big(\frac{k\hsy \pi}{A} + B \Big)
\ \geq \ 
0
\]
$\forall \ k \neq p$, so take $k = p$ and consider
\[
(-1)^{n+p}
g \Big(\frac{p\hsy \pi}{A} + B \Big)
\]
or still, 
\[
(-1)^{n+p}
g (x_p)
\]
or still, 
\[
\lim\limits_{h \ra 0}  \ 
(-1)^{n+p} 
\hsy 
g(x_p + h)
\]
or still, 
\[
\lim\limits_{h \ra 0}  \ 
(-1)^{n+p} 
\hsy 
\frac{f (x_p + h)}{(x_p + h - x_p)^2}
\]
or still, 
\[
\lim\limits_{h \ra 0}  \ 
(-1)^{n+p} 
\hsy 
\frac{f (x_p + h)}{h^2}
\]
or still, 
\[
\lim\limits_{h \ra 0}  \ 
(-1)^{n+p} 
\hsy 
\frac{f^\prime (x_p + h)}{2 \hsy h}
\]
or still, 
\[
\lim\limits_{h \ra 0}  \ 
(-1)^{n+p} 
\hsy 
\frac{f^{\prime\prime} (x_p + h)}{2}
\]
or still, 
\[
\frac{1}{2}
\hsx
(-1)^{n+p} 
\hsy 
f^{\prime\prime} (x_p) 
\ \geq \ 
0.
\]
Therefore $g$ is standard of level $B$.  
But
\[
\sup\limits_{x \in \R} x^2 \abs{g (x)}
\ < \ 
\infty,
\]
contradicting 32.10.  
Accordingly, the supposition 
\[
(-1)^{n+p} 
\hsy 
f^{\prime\prime} (x_p) 
\ \geq \ 
0
\]
is untenable, leaving 
\[
(-1)^{n+p} 
\hsy 
f^{\prime\prime} (x_p) 
\ < \ 
0.
\]
\\[-1cm]

\qquad \textbullet \quad
To see that $f$ cannot have zeros in both intervals
$I_p$ 
and 
$I_{p+1}$, 
assume the opposite: 
\[
\begin{cases}
\ 
f (x_1) = 0 \qquad (x_1 \in I_p)
\\[4pt]
\ 
f (x_2) = 0 \qquad (x_2 \in I_{p+1})
\end{cases}
.
\]

\noindent
Then $x_1$ is the only zero of $f$ in $I_p$ and it is simple, 
whereas $x_2$ is the only zero of $f$ in $I_{p+1}$ and it is simple (cf. 32.9).
Now form

\[
g (z) 
\ = \ 
\frac{f (z) \hsy (z - x_p)^2}{(z - x_1)\hsy (z - x_2)}.
\]
Then 
$g \in \R \sB_0 (A)$ 
and $g$ is standard of level $B$: $\forall \ k \in \Z$, 

\[
(-1)^{n+k} 
g \Big(\frac{k\hsy \pi}{A} + B \Big).
\]
Here the point is slightly subtle and explains the presence of two factors in the denominator rather than just one factor.  
For
\[
\frac{(p-1)\hsy \pi}{A} + B 
\ < \ 
x_1
\ < \ 
x_2,
\]
so
\[
\frac{k \pi}{A} + B 
\ \leq \ 
\frac{(p-1)\hsy \pi}{A} + B 
\]
\qquad \qquad
$\implies$

\[
\frac{k \pi}{A} + B - x_1 
\ < \ 
0, 
\quad 
\frac{k \pi}{A} + B - x_2 
\ < \ 
0
\]
\qquad \qquad
$\implies$
\[
\Big(
\frac{k \pi}{A} + B - x_1 
\Big)
\hsx
\Big(
\frac{k \pi}{A} + B - x_2 
\Big)
\ > \ 
0.
\]
What remains is obvious and one then comes to a contradiction, $x_p$ being a zero of $g$ of multiplicity $> 2$.
\\[-.25cm]

\qquad \textbullet \quad
Suppose that $x_p$ is a zero of multiplicity 2 $-$then $f$ has no zeros in 
$I_p \cup I_{p+1}$.  
E.g.: \ 
Let $x_1 \in I_p$ be a zero of $f$ and put
\[
g (z) 
\ = \ 
\frac{f (z)}{(z - x_1) (z - x_p)}.
\]
Then 
$g \in \R \sB_0 (A)$ 
is standard of level $B$.  
On the other hand, 
\[
\sup\limits_{x \in \R} x^2 \abs{g (x)}
\ < \ 
\infty,
\]
which is incompatible with 32.10.  
Bearing in mind that
\[
(-1)^{n+p} f^{\prime\prime} (x_p)
\ < \ 
0,
\]
it then follows that
\[
(-1)^{n+p} f (x)
\ < \ 
0
\qquad (x \in I_p \cup I_{p+1}).
\]
Thus choose a small $h > 0$ with the property that
\[
(-1)^{n+p} \ 
\begin{cases}
\ 
f (x)
\\[11pt]
\ 
f^\prime (x)
\end{cases}
\quad \text{and} \ (-1)^{n+p} \quad
\begin{cases}
\ 
f^\prime (x)
\\[11pt]
\ 
f^{\prime\prime} (x)
\end{cases}
\]
have the same sign in 
$]x_p, x_p + h[$ 
and opposite signs in 
$]x_p - h, x_p[$.  
Working first with 
$]x_p, x_p + h[$ 
and assuming, as we may, that
\[
x \in \  ]x_p, x_p + h[\ 
\implies
(-1)^{n+p} f^{\prime\prime} (x)
\\[11pt]
\ < \ 
0,
\]
thence
\\[-1cm]
\allowdisplaybreaks
\begin{align*}
x \in \ ]x_p, x_p + h[\ 
&\implies
(-1)^{n+p} f^\prime (x)
\ < \ 
0
\\[11pt]
&\implies
(-1)^{n+p} f (x)
\ < \ 
0.
\end{align*}
But $f$ has no zeros in $I_{p+1}$, so 
\[
(-1)^{n+p} f (x)
\ < \ 
0 
\qquad (x \in I_{p+1}).
\]
As for 
$]x_p - h, x_p[$, 
it can be assumed that

\[
x \in \  ]x_p - h, x_p[\ 
\implies
(-1)^{n+p} f^{\prime\prime} (x)
\\[11pt]
\ < \ 
0,
\]
thence
\allowdisplaybreaks
\begin{align*}
x \in \  ]x_p - h, x_p[\ 
&\implies
(-1)^{n+p} f^\prime (x)
\ > \ 
0
\\[11pt]
&\implies
(-1)^{n+p} f (x)
\ < \ 
0.
\end{align*}
But $f$ has no zeros in $I_p$, so 
\[
(-1)^{n+p} f (x)
\ < \ 
0 
\qquad (x \in I_p).
\]

\qquad \textbullet \quad
That $x_p$ and $x_{p-1}$ cannot both be zeros of multiplicity 2 is ruled out by consideration of
\[
g (z) 
\ = \ 
\frac{f (z)}{(z - x_{p-1}) (z - x_p)}.
\]

The zero theory for 
$f^\prime$ 
can be reduced to that for $f$.  
To begin with, matters are trivial if 
\[
f (x) 
\ = \ 
C \hsy \sin A (x - B)
\qquad (C \neq 0),
\]
so this case can be ignored.  
Suppose, therefore, that
$
\ds
f \Big(\frac{k \hsy \pi}{A} + B\Big) \neq 0
$
for some $k$ and in 32.8 take
\[
x 
\ = \ 
\frac{p \hsy \pi}{A}
+ 
\frac{\pi}{2 \hsy A}
+ B
\qquad (p \in \Z).
\]
Then
\allowdisplaybreaks
\begin{align*}
\cos A \hsy \Big(\frac{p \hsy \pi}{A} +  \frac{\pi}{2 \hsy A} + B - B\Big) \ 
&=\ 
\cos \Big( p \pi + \frac{\pi}{2}\Big)
\\[11pt]
&=\ 
\cos p \hsy \pi 
\hsy
\cos \Big(\frac{\pi}{2}\Big) 
\hsx - \hsx 
\sin p \hsy \pi 
\hsy
\sin  \Big(\frac{\pi}{2}\Big)
\\[11pt]
&=\ 
0
\end{align*}
\\[-1cm]
and
\\[-1cm]
\allowdisplaybreaks
\begin{align*}
\sin A \hsy \Big(\frac{p \hsy \pi}{A} +  \frac{\pi}{2 \hsy A} + B - B\Big) \ 
&=\ 
\sin \Big(p \hsy \pi + \frac{\pi}{2}\Big)
\\[11pt]
&=\ 
\sin p \hsy \pi 
\hsy
\cos \Big(\frac{\pi}{2}\Big) 
\hsx + \hsx 
\sin  \Big(\frac{\pi}{2}\Big)
\hsy
\cos p \hsy \pi 
\\[11pt]
&=\ 
(-1)^p
\end{align*}
\\[-1cm]
\qquad 
$\implies$
\[
(-1)^p
\hsx
f^\prime \Big(\frac{p\hsy \pi}{A} + \frac{\pi}{2 \hsy A}  + B \Big)
\ = \ 
(-1)^{n-1} 
\hsx 
M (p) 
\qquad (M (p) > 0)
\]
 \qquad
$\implies$
\[
(-1)^{n-1} \hsy (-1)^p
\hsx
f^\prime \Big(\frac{p\hsy \pi}{A} + \frac{\pi}{2 \hsy A}  + B \Big)
\ > \ 
0
\]
 \qquad
$\implies$

\[
(-1)^{n^\prime} \hsy (-1)^p
f^\prime \Big(\frac{p\hsy \pi}{A} + \frac{\pi}{2 \hsy A}  + B \Big)
\ > \ 
0,
\]
where
\[
\begin{cases}
\ 
n^\prime = 0 \hspace{0.25cm} \text{if} \hspace{0.25cm} n = 1
\\[4pt]
\ 
n^\prime = 1 \hspace{0.25cm} \text{if} \hspace{0.25cm} n = 0
\end{cases}
.
\]
\\[-.5cm]

\noindent
I.e.: \ 
$f^\prime$ is standard of level 
$
\ds
\frac{\pi}{2 \hsy A}  + B
$.
\\[-.25cm]
\end{x}

\qquad
{\small\bf \un{N.B.}} \ 
The ambient interval per 
$f^\prime$ 
is 
\[
I_p^\prime 
\ = \ 
\Big] 
\frac{(p-1) \hsy \pi}{A} 
\hsx + \hsx 
\frac{\pi}{2 \hsy A} 
\hsx + \hsx 
B
, \ 
\frac{p \hsy \pi}{A} 
\hsx + \hsx 
\frac{\pi}{2 \hsy A} 
\hsx + \hsx 
B
\Big[
\hsx 
.
\]
\\[-1cm]

\begin{x}{\small\bf LEMMA} \ 
The zeros of 
$f^\prime$
are real 
(cf. 32.11).
\\[-.25cm]
\end{x}

\begin{x}{\small\bf LEMMA} \ 
The zeros of 
$f^\prime$
are simple.
\\[-.5cm]

PROOF \ 
The only possibility for a nonsimple zero is at an endpoint of an ambient interval (cf. 32.9) 
and at such an endpoint, 
$f^\prime$ 
does not vanish.
\\[-.25cm]
\end{x}

\begin{x}{\small\bf LEMMA} \ 
$\forall \ p \in \Z$, 
$f^\prime$ 
has a zero in the ambient interval $I_p^\prime$ (it being
necessarily unique).
\\[-.5cm]

PROOF \ 
We have 
\[
(-1)^{n^\prime} \hsy (-1)^{p-1}
\hsx
f^\prime \Big(\frac{(p-1)\hsy \pi}{A} + \frac{\pi}{2 \hsy A}  + B \Big)
\ > \ 
0
\]
and 
\[
(-1)^{n^\prime} \hsy (-1)^p
\hsx
f^\prime \Big(\frac{p\hsy \pi}{A} + \frac{\pi}{2 \hsy A}  + B \Big)
\ > \ 
0.
\]
\\[-.5cm]

\qquad \textbullet \quad
\un{$p$ \hsx even}:\ 
Then
\[
(-1)^{n^\prime}
\hsx
f^\prime \Big(\frac{(p-1)\hsy \pi}{A} + \frac{\pi}{2 \hsy A}  + B \Big)
\ < \ 
0
\]
while 
\[
(-1)^{n^\prime} 
\hsx
f^\prime \Big(\frac{p\hsy \pi}{A} + \frac{\pi}{2 \hsy A}  + B \Big)
\ > \ 
0.
\]
\\[-.5cm]

\qquad \textbullet \quad
\un{$p$ \hsx odd}:\ 
Then
\[
(-1)^{n^\prime}
\hsx
f^\prime \Big(\frac{(p-1)\hsy \pi}{A} + \frac{\pi}{2 \hsy A}  + B \Big)
\ > \ 
0
\]
while 
\[
(-1)^{n^\prime}
\hsx
f^\prime \Big(\frac{p\hsy \pi}{A} + \frac{\pi}{2 \hsy A}  + B \Big)
\ < \ 
0.
\]
But this means that $f^\prime$ has a zero in $I_p^\prime$.
\\[-.25cm]
\end{x}

\begin{x}{\small\bf LEMMA} \ 
Take $C$ per 32.12 
$(\implies A = 1, B = 0)$ 
$-$then $C^\prime$ is standard 
\\[-.5cm]

\noindent
of level 
$
\ds 
\frac{\pi}{2}
$ 
and $n = 0 \implies n^\prime = 1$
\\[-.5cm]

\qquad 
$\implies$
\[
(-1)^1 \hsy (-1)^k
\hsx 
C^\prime \Big(k \hsy \pi + \frac{\pi}{2}\Big)
\ > \ 
0.
\]
And all the zeros of $C^\prime$ are real, each ambient interval $I_p^\prime$ contains a single zero and this zero is simple.
\\[-.25cm]
\end{x}


There is another situation which arises in the applications.
\\[-.25cm]

\begin{x}{\small\bf DEFINITION} \ 
Let 
$f \in \R \sB_0 (A)$ 
$-$then $f$ is \un{semi-standard} of level $B$ if 
$\exists \ n = 0$ or 1 and $B \in \R$ such that $\forall \ k \in \Z$, 
\[
\begin{cases}
\ \ds
(-1)^{n + k}
\hsx 
f \Big(\frac{k \hsy \pi}{A} + B\Big) 
\ \leq \ 
0
\qquad (k \geq 1)
\\[18pt]
\ \ds
(-1)^{n + k}
\hsx 
f \Big(\frac{k \hsy \pi}{A} + B\Big) 
\ \geq \ 
0
\qquad (k \leq 0)
\end{cases}
.
\]

[Note: \ 
A fundamental class of examples is dealt with in the next \S.]
\\[-.25cm]
\end{x}

Suppose that $f$ is semi-standard of level $B$.  
Fix 
$
\ds 
x _0 \in I_1 = \ 
\big]B, \frac{\pi}{A} + B \big[
$
and let
\[
g (z) 
\ = \ 
(x_0 - z) \hsy f (z).
\]
Impose the condition 
\[
\sup\limits_{x \in \R} \abs{x \hsy f (x)}
\ < \ 
\infty.
\]
Then $g$ is standard of level $B$.  
But $g (x_0) = 0$, thus $g$ has a unique zero in $I_1$, viz. $x_0$.  
Therefore 
\[
x \in I_1 
\implies 
f (x) \neq 0.
\]
In addition, however, 
\[
(-1)^{n + 1}
\hsy 
g^\prime (x_0)
\ > \ 
0 
\qquad (\tcf. \ 32.9).
\]
So
\[
f (x_0) 
\ = \ 
g^\prime (x_0)
\]
\qquad 
$\implies$
\allowdisplaybreaks
\begin{align*}
(-1)^n 
\hsx 
f (x_0) 
&=\ 
(-1)^n \hsy (-1)^1
\hsy
g^\prime (x_0)
\\[11pt]
&=\ 
(-1)^{n + 1}
\hsy 
g^\prime (x_0)
\\[11pt]
&> \ 
0.
\end{align*}
Therefore
\[
x \in I_1 
\implies 
(-1)^n 
\hsy 
f (x) 
\ > \ 
0.
\]
\\[-1.25cm]

\begin{x}{\small\bf THEOREM} \ 
Suppose that $f$ is semi-standard of level $B$ and 
\[
\sup\limits_{x \in \R} \abs{x \hsy f (x)}
\ < \ 
\infty.
\]
Then all the zeros of $f$ are real (cf. 32.11).  
Furthermore, the ambient interval
\[
I_p
\ = \ 
\Big]
\frac{(p-1) \hsy \pi}{A} + B, 
\frac{p \hsy \pi}{A} + B
\Big[
\qquad (p \in \Z, \hsx p \neq 1)
\]
contains at most one zero of $f$ and if there is one, then it must be simple.  
Finally, 
\[
x \in I_1 
\implies 
(-1)^n 
\hsy 
f (x) 
\ > \ 
0.
\] 

Picture: 
\[
\begin{tikzpicture}[scale=1.75]
      

      
      \draw[-] (-3,0) -- (3,0) node[right] {$$};

      \draw[] (-3,-.015) node[] {$\text{\textbullet}$};
      \draw[] (-1,-.015) node[] {$\text{\textbullet}$};
      \draw[] (1,-.015) node[] {$\text{\textbullet}$};
      \draw[] (3,-.015) node[] {$\text{\textbullet}$};

      \draw[] (-3,-.1) node[below] {$\ds\frac{-\pi}{A} + B$};
      \draw[] (-1,-.1) node[below] {$\ds\frac{0 \hsy \pi}{A} + B$};
      \draw[] (1,-.1) node[below] {$\ds\frac{\pi}{A} + B$};
      \draw[] (3,-.1) node[below] {$\ds\frac{2 \hsy \pi}{A} + B$};
      
      \draw[] (-2,.1) node[above] {$I_0$};
      \draw[] (0,.1) node[above] {$I_1$};
      \draw[] ( 2,.1) node[above] {$I_2$};
      

\end{tikzpicture}
\]
\\[-.75cm]
\end{x}

\begin{x}{\small\bf THEOREM} \ 
Suppose that $f$ is semi-standard of level $B$ and 
\[
\sup\limits_{x \in \R} \abs{x \hsy f (x)}\ < \ 
\infty.
\]

\qquad \textbullet \quad
If 
$f (B) = 0$, 
then its multiplicity is equal to 1 and there are no zeros of $f$ in $I_0 \cup I_1$.
\\[-.5cm]

[Apply 32.13 to 
\[
g (z) 
\ = \ 
(B - z) \hsy f (z).
\]
Then per $g$, $B$ is a zero of multiplicity 2, hence $(p = 0)$
\[
(-1)^n
\hsy
g (x)
\ < \ 
0
\qquad (x \in I_0 \cup I_1)
\]
\qquad \qquad
$\implies$
\[
(-1)^n
\hsy
(B - x)
\hsy
f (x)
\ < \ 
0
\qquad (x \in I_0)
\]
\qquad \qquad
$\implies$
\[
(-1)^n
\hsy
f (x)
\ < \ 
0
\qquad (x \in I_0).
\]
On the other hand, a priori, 
\[
(-1)^n
\hsy
f (x)
\ > \ 
0
\qquad (x \in I_1).
\]

\qquad \textbullet \quad
If 
$
\ds
f \Big(\frac{\pi}{A} + B\Big) = 0
$, 
then its multiplicity is equal to 1 and there are no 
\\[-.5cm]

\noindent
zeros of $f$ in $I_1 \cup I_2$.
\\[-.5cm]

[Apply 32.13 to 
\[
g (z) 
\ = \ 
\Big(\frac{\pi}{A} + B - z \Big) \hsy f (z).
\]
Then per $g$, 
$
\ds
\frac{\pi}{A} + B
$  
is a zero of multiplicity 2, hence $(p = 1)$
\[
(-1)^{n+1} 
\hsy
g (x)
\ < \ 
0
\qquad (x \in I_1 \cup I_2)
\]
\qquad \qquad
$\implies$
\[
(-1)^{n+1} 
\hsy
\Big(\frac{\pi}{A} + B - x \Big) \hsy f (x)
\ < \ 
0
\qquad (x \in I_2)
\]
\qquad \qquad
$\implies$
\[
(-1)^n
\hsy
\Big(x - \frac{\pi}{A} - B\Big) \hsy f (x)
\ < \ 
0
\qquad (x \in I_2)
\]
\qquad \qquad
$\implies$
\[
(-1)^n
\hsy
f (x)
\ < \ 
0
\qquad (x \in I_2).
\]
On the other hand, a priori, 
\[
(-1)^n 
\hsy 
f (x) 
\ > \ 
0
\qquad 
(x \in I_1).]
\]
\\[-1.25cm]
\end{x}


\begin{x}{\small\bf REMARK} \ 
The condition 
\[
\sup\limits_{x \in \R} \abs{x \hsy f (x)}
\ < \ 
\infty
\]
is not automatic (consider $\sin A (x - B)$).
\\[-.25cm]
\end{x}


\chapter{
$\boldsymbol{\S}$\textbf{33}.\quad  ZEROS OF $W_{A, \alpha}$}
\setlength\parindent{2em}
\setcounter{theoremn}{0}
\renewcommand{\thepage}{\S33-\arabic{page}}

\qquad 
Working on $]0,A[$ $(A > 0)$, 
suppose that $\phi$ is defined on $]0,A[$ and is integrable on $[0,A]$.    
Assume further that $\phi$ is positive and increasing on $]0,A[$.
\\[-.25cm]

\begin{x}{\small\bf NOTATION} \ 
Given $\alpha \in [0,\pi[$, let
\[
W_{A, \alpha} (z) 
\ = \ 
\int\limits_0^A \ 
\phi (t) 
\hsx
\sin (z \hsy t +\alpha) 
\ \td t,
\]
thus 
\[
W_{A, \alpha} (z)
\ = \ 
(\sin \alpha) 
\hsx 
C_A (z) 
\hsx + \hsx 
(\cos \alpha) 
\hsx 
S_A (z),
\]
where
\[
\begin{cases}
\ \ds
C_A (z) 
\ = \ 
\int\limits_0^A \ 
\phi (t) 
\hsx
\cos z \hsy t 
\ \td t
\\[18pt]
\ \ds
S_A (z) 
\ = \ 
\int\limits_0^A \ 
\phi (t) 
\hsx
\sin z \hsy t 
\ \td t
\end{cases}
.
\]
\\[-1cm]
\end{x}

It is clear that 
$W_{A, \alpha} \in \R\sB_0 (A)$.
\\[-.25cm]

\begin{x}{\small\bf LEMMA} \ 
$W_{A, \alpha}$ is semi-standard of level 
$
\ds
-\frac{\alpha}{A}
$.
\\[-.25cm]

PROOF \ 
In 32.18, take $n = 0$, the issue being $\forall \ k \in \Z$ the inequalities
\[
\begin{cases}
\ \ds
(-1)^k 
\hsx
W_{A, \alpha} 
\Big(
\frac{k \hsy \pi - \alpha}{A}
\Big)
\ \leq \ 
0
\qquad (k \geq 1)
\\[18pt]
\ \ds
(-1)^k 
\hsx
W_{A, \alpha} 
\Big(
\frac{k \hsy \pi - \alpha}{A}
\Big)
\ \geq \ 
0
\qquad (k \leq 0)
\end{cases}
.
\]

\qquad \textbullet \quad
\un{$k = 0$}: \ 
Here
\[
W_{A, \alpha} 
\Big(
-
\frac{\alpha}{A}
\Big)
\ = \ 
\int\limits_0^A \ 
\phi (t) 
\hsy
\sin
\Big(
\frac{\alpha (A - t)}{A}
\Big)
\ \td t 
\ \geq \ 
0
\]
and
\[
W_{A, \alpha} 
\Big(
-
\frac{\alpha}{A}
\Big)
\ = \ 
0
\]
iff $\alpha = 0$.
\\[-.25cm]

\qquad \textbullet \quad
$\un{k = 1, 2, \ldots} :$ \ 
Here

\[
W_{A, \alpha} 
\Big(
\frac{k \hsy \pi - \alpha}{A}
\Big)
\ = \ 
\frac{A}{k \hsy \pi - \alpha}
\
\int\limits_0^{k \hsy \pi} \ 
\phi
\Big(
\frac{A (s - \alpha)}{k \hsy \pi - \alpha}
\Big)
\hsy 
\sin s 
\ \td s
\]
and
\[
\frac{A}{k \hsy \pi - \alpha}
\ > \
0.
\]
\\[-1cm]

\qquad \textbullet \quad
\un{$\lra: k$ odd}: \ 
Split the interval of integration 
$[\alpha, k \hsy \pi]$ 
into the closed subintervals
$[\alpha, \pi], [\pi, 3 \pi], 
\ldots, 
[k \hsy \pi - 2 \pi,  k \hsy \pi]$ 
$-$then the integral over each of these subintervals is nonnegative, hence
\[
(-1)^k 
\hsx
W_{A, \alpha} 
\Big(
\frac{k \hsy \pi - \alpha}{A}
\Big)
\ \leq \ 
0.
\]
\\[-1cm]

\qquad \textbullet \quad
\un{$\lra: k$ even}: \ 
Split the interval of integration 
$[\alpha, k \hsy \pi]$ 
into the closed subintervals
$[\alpha, 2 \pi], [2 \pi, 4 \pi], 
\ldots, 
[k \hsy \pi - 2 \pi,  k \hsy \pi]$ 
$-$then the integral over each of these subintervals is nonpositive, hence
\[
(-1)^k 
\hsx
W_{A, \alpha} 
\Big(
\frac{k \hsy \pi - \alpha}{A}
\Big)
\ \leq \ 
0.
\]
\\[-1cm]

\qquad \textbullet \quad
$\un{k = -1, -2, \ldots}:$ \ 
Here
\[
W_{A, \alpha} 
\Big(
\frac{k \hsy \pi - \alpha}{A}
\Big)
\ = \ 
\frac{A}{k \hsy \pi - \alpha}
\
\int\limits_{- \alpha}^{- k \hsy \pi} \ 
\phi
\Big(
\frac{A (s + \alpha)}{\alpha - k \hsy \pi}
\Big)
\hsy 
\sin s 
\ \td s
\]
and
\[
\frac{A}{k \hsy \pi - \alpha}
\ < \
0.
\]
\\[-1cm]

\qquad \textbullet \quad
\un{$\lra: k$ odd}: \ 
Split the interval of integration 
$[- \alpha,- k \hsy \pi]$ 
into the closed subintervals
$[-\alpha, \pi], [\pi, 3 \pi], 
\ldots, 
[-k \hsy \pi - 2 \pi, - k \hsy \pi]$ 
$-$then the integral over each of these subintervals is nonpositive, hence
\[
(-1)^k 
\hsx
W_{A, \alpha} 
\Big(
\frac{k \hsy \pi - \alpha}{A}
\Big)
\ \geq \ 
0.
\]
\\[-1cm]

\qquad \textbullet \quad
\un{$\lra: k$ even}: \ 
Split the interval of integration 
$[- \alpha,- k \hsy \pi]$ 
into the closed subintervals
$[-\alpha, 0], [0, 2 \pi], 
\ldots, 
[-k \hsy \pi - 2 \pi, - k \hsy \pi]$ 
$-$then the integral over each of these subintervals is nonpositive, hence
\[
(-1)^k 
\hsx
W_{A, \alpha} 
\Big(
\frac{k \hsy \pi - \alpha}{A}
\Big)
\ \geq \ 
0.
\]
\\[-1cm]
\end{x}

\begin{x}{\small\bf APPLICATION} \ 
If $\phi$ is bounded on 
$]0, A[$, 
then all the zeros of 
$W_{A, \alpha}$ 
are real.  
Furthermore, the ambient interval
\[
I_p
\ = \ 
\Big]
\frac{(p-1) \pi - \alpha}{A}, 
\frac{p \pi - \alpha}{A}
\Big[
\qquad 
(p \in \Z, \hsx p \neq 1)
\]
contains at most one zero of 
$W_{A, \alpha}$ 
and if there exists one, then it must be simple.  
Finally, 
\allowdisplaybreaks
\begin{align*}
x \in I_1 
&\implies
(-1)^n \hsx W_{A, \alpha}  (x) > 0
\\[11pt]
&\implies
\hspace{1.15cm}
W_{A, \alpha}  (x) > 0
\qquad (n = 0).
\end{align*}

[In fact, 
\[
\sup\limits_{x \in \R} \ 
\abs{x \hsy W_{A, \alpha}  (x)} 
\ \leq \ 
2 \ 
\lim\limits_{t \uparrow A} \ 
\phi (t) 
\ < \ 
\infty, 
\]
so one can quote 32.19.]
\\[-.25cm]
\end{x}

A finer analysis will lead to more precise results.
\\

\qquad \textbullet \quad
\un{$k \geq 1$ ($k$ odd)}: \ 
Suppose that 
\[
W_{A, \alpha} \Big(\frac{k \hsy \pi - \alpha}{A}\Big)
\ = \ 
0.
\]
Then there exist constants
\[
0 < c_0 \leq c_1 \leq \cdots \leq c_{(k-1)/2}
\]
and points
\[
t_{-1} \ = \ 0, 
\quad 
t_j 
\ = \ 
A 
\hsx
\frac{(2j + 1) \pi - \alpha}{k \hsy \pi - \alpha}
\]
such that 
\[
\phi (t) 
\ = \ 
c_j 
\qquad 
(t_{j-1} < t < t_j) \ 
\Big(0 \leq j \leq \frac{k-1}{2}\Big).
\]
Therefore
\[
W_{A, \alpha} (x) 
\ = \ 
\frac{2}{x} 
\hsx 
\sin 
\Big(
\frac{A \hsy \pi \hsy x}{k \hsy \pi - \alpha}
\Big)
\ 
\sum\limits_{j = 0}^{(k-1)/2}  \ 
c_j 
\hsx 
\sin 
\Big(
\frac{2 \hsy j \hsy \pi - \alpha}{k \hsy \pi - \alpha}
\hsx 
A x + \alpha
\Big).
\]
\\[-.75cm]

\qquad \textbullet \quad
\un{$k \geq 1$ ($k$ even)}: \ 
Suppose that 
\[
W_{A, \alpha} \Big(\frac{k \hsy \pi - \alpha}{A}\Big)
\ = \ 
0.
\]
Then there exist constants
\[
0 < c_0 \leq c_1 \leq \cdots \leq c_{(k-2)/2}
\]
and points
\[
t_{-1} \ = \ 0, 
\quad 
t_j 
\ = \ 
A 
\hsx
\frac{(2j + 2) \pi - \alpha}{k \hsy \pi - \alpha}
\]
such that 
\[
\phi (t) 
\ = \ 
c_j 
\qquad 
(t_{j-1} < t < t_j) \ 
\Big(0 \leq j \leq \frac{k-2}{2}\Big).
\]
Therefore
\[
W_{A, \alpha} (x) 
\ = \ 
\frac{2}{x} 
\hsx 
\sin 
\Big(
\frac{A \hsy \pi \hsy x}{k \hsy \pi - \alpha}
\Big)
\ 
\sum\limits_{j = 0}^{(k-2)/2}  \ 
c_j 
\hsx 
\sin 
\Big(
\frac{(2 \hsy j + 1) \hsy \pi - \alpha}{k \hsy \pi - \alpha}
\hsx A x + \alpha
\Big).
\]
\\[-.75cm]

\qquad \textbullet \quad
\un{$k \leq -1$ ($k$ odd)}: \ 
Suppose that 
\[
W_{A, \alpha} \Big(\frac{k \hsy \pi - \alpha}{A}\Big)
\ = \ 
0.
\]
Then there exist constants
\[
0 < c_0 \leq c_1 \leq \cdots \leq c_{(-k-1)/2}
\]
and points
\[
t_{-1} \ = \ 0, 
\quad 
t_j 
\ = \ 
A 
\hsx
\frac{(2j + 1) \pi + \alpha}{\alpha - k \hsy \pi}
\]
such that 
\[
\phi (t) 
\ = \ 
c_j 
\qquad 
(t_{j-1} < t < t_j) \ 
\Big(0 \leq j \leq \frac{- k - 1}{2}\Big).
\]
Therefore
\[
W_{A, \alpha} (x) 
\ = \ 
\frac{2}{x} 
\hsx 
\sin 
\Big(
\frac{A \hsy \pi \hsy x}{\alpha - k \hsy \pi}
\Big)
\ 
\sum\limits_{j = 0}^{(-k - 1)/2}  \ 
c_j 
\hsx 
\sin 
\Big(
\frac{2 \hsy j  \hsy \pi + \alpha}{\alpha - k \hsy \pi}
\hsx A x + \alpha
\Big).
\]
\\[-.5cm]

\qquad \textbullet \quad
\un{$k \leq -1$ ($k$ even)}: \ 
Suppose that 
\[
W_{A, \alpha} \Big(\frac{k \hsy \pi - \alpha}{A}\Big)
\ = \ 
0.
\]
Then there exist constants
\[
0 < c_0 \leq c_1 \leq \cdots \leq c_{-k/2}
\]
and points
\[
t_{-1} \ = \ 0, 
\quad 
t_j 
\ = \ 
A 
\hsx
\frac{2 j \pi + \alpha}{\alpha - k \hsy \pi}
\]
such that 
\[
\phi (t) 
\ = \ 
c_j 
\qquad 
(t_{j-1} < t < t_j) 
\qquad 
\Big(0 \leq j \leq -\frac{k}{2}\Big).
\]
Therefore
\[
W_{A, \alpha} (x) 
\ = \ 
\frac{2}{x} 
\hsx 
\sin 
\Big(
\frac{A \hsy \pi \hsy x}{\alpha - k \hsy \pi}
\Big)
\ 
\sum\limits_{j = 0}^{-k/2}  \ 
c_j 
\hsx 
\sin 
\Big(
\frac{(2 \hsy j - 1) \hsy \pi + \alpha}{\alpha - k \hsy \pi}
\hsx A x + \alpha
\Big).
\]
\\[-.75cm]

\begin{x}{\small\bf NOTATION} \ 
Write
\[
E (A, \alpha, k)
\]
for the set of those $\phi$ such that 
\[
W_{A, \alpha} \Big(\frac{k \hsy \pi - \alpha}{A}\Big)
\ = \ 
0
\]
for some $k \in \Z - \{0\}$ and put
\[
E (A, \alpha) 
\ = \ 
\bigcup\limits_k \ 
E (A, \alpha, k).
\]

[Note: \ 
In general, 
\[
E (A, \alpha, k_1)
\cap
E (A, \alpha, k_2)
\ \neq \ 
\emptyset.]
\]
\\[-1.25cm]
\end{x}

\begin{x}{\small\bf RECONCILIATION} \ 
Take $A = 1$, $\alpha = 0$, hence
\[
W_{1, 0} (z) 
\ = \ 
\int\limits_0^1 \ 
\phi (t) 
\hsx 
\sin z \hsy t 
\ \td t.
\]
Recall now the definition of ``exceptional'' from 29.14 and the notation $E(1, 0)$ from 29.15 
$-$then the claim is that the two possible meanings of $E(1, 0)$ are one and the same.  
To see this, consider
\[
W_{1, 0}
\Big(
\frac{k \hsy \pi - \alpha}{A}
\Big)
\ \equiv \ 
W_{1, 0} (k \hsy \pi) 
\qquad (k = \pm 1, \pm 2, \ldots),
\]
there being no loss of generality in assuming that $k =1, 2, \ldots \hsx .$
\\[-.25cm]

\qquad \textbullet \quad
\un{$k$ odd}: \ 
Here
\[
W_{1, 0} (k \hsy \pi) 
\ > \ 
0
\qquad 
(k = 1, 3, \ldots) 
\quad 
(\tcf. \ 31.3).
\]
Therefore
\[
E (1, 0, k \ \text{odd}) 
\ = \ 
\emptyset. 
\]
\\[-.5cm]

\qquad \textbullet \quad
\un{$k$ even}: \ 
Suppose that
\[
W_{1, 0} (2 n \pi)
\ = \ 
0
\quad 
\text{for some $n = 1, 2, \ldots$}\ .
\]
I.e.: 
\[
\int\limits_0^1 \ 
\phi (t) 
\hsx 
\sin 2 \hsy n  \pi  t 
\ \td t 
\ = \ 
0.
\]
But this implies that $\phi$ is exceptional (look at the proof of 31.3).  
Therefore
\[
E (1, 0, k \ \text{even})
\]
is comprised of exceptional $\phi$, so
\[
\bigcup\limits_{n = 1}^\infty \ 
E (1, 0, 2 \hsy  n)
\]
is contained in the 
$E(1, 0)$ 
per 29.15.  
To turn matters around, take an exceptional $\phi$ and write
\allowdisplaybreaks
\begin{align*}
f (z) \ 
&=\ 
\int\limits_0^1 \ 
\phi (t) 
\hsx 
e^{\sqrt{-1} \hsx z \hsy t}
\ \td t 
\\[15pt]
&=\ 
C (z) 
\hsx + \hsx 
\sqrt{-1} \hsx S (z)
\end{align*}
where, of course, 
\[
S (z) 
\ \equiv \ 
W_{1, 0} (z).
\]
Then in the notation of 29.20, 
\[
f (2 \pi  q) 
\ = \ 
0
\]
\qquad 
$\implies$
\[
C (2 \pi  q) 
\hsx + \hsx 
\sqrt{-1} \hsx S(2 \pi  q) 
\ = \ 
0
\]
\qquad 
$\implies$
\[
S (2 \pi  q) 
\ = \ 
0
\hsx \implies \hsx 
W_{1, 0} (2 \pi  q)
\ = \ 
0
\]
\qquad 
$\implies$
\[
\phi \in E (1, 0, 2 q).
\]
Conclusion: 
\[
E (1, 0) 
\subsetx 
\bigcup\limits_{n = 1}^\infty \ 
E (1, 0, 2 n)
\subsetx
E (1, 0) .
\]
\\[-1.5cm]
\end{x}


\begin{x}{\small\bf REMARK} \ 
If $\phi \in E (A, \alpha)$, then 
\[
\sup\limits_{x \in \R} \ 
\abs{x \hsy W_{A, \alpha}  (x)} 
\ < \ 
\infty.
\]

[Note: \ 
Accordingly, all the particulars of the semi-standard theory developed at the end of \S32 are in force but the detailed 
explication thereof will be left to the reader.]
\\[-.25cm]
\end{x}

\begin{x}{\small\bf LEMMA} \ 
If $\phi \notin E (A, \alpha)$, then 
\[
\begin{cases}
\ \ds
(-1)^k 
\hsx
W_{A, \alpha} 
\Big(
\frac{k \hsy \pi - \alpha}{A}
\Big)
\ < \ 
0
\qquad (k \geq 1)
\\[18pt]
\ \ds
(-1)^k 
\hsx
W_{A, \alpha} 
\Big(
\frac{k \hsy \pi - \alpha}{A}
\Big)
\ > \ 
0
\qquad (k \leq -1)
\end{cases}
\]
and at $k = 0$, 
\[
W_{A, \alpha} \Big(- \frac{\alpha}{A}\Big)
\ > \ 
0
\qquad (0 < \alpha < \pi).
\]
\\[-1cm]
\end{x}

\begin{x}{\small\bf LEMMA} \ 
If $\phi \notin E (A, \alpha)$ and if 
\[
\sup\limits_{x \in \R} \ 
\abs{x \hsy W_{A, \alpha}  (x)} 
\ < \ 
\infty, 
\]
then all the zeros of 
$W_{A, \alpha}$ 
are real (cf. 32.11) and simple (cf. infra).
\\[-.5cm]

PROOF \ 
The ambient interval
\[
I_p
\ = \ 
\Big]
\frac{(p-1) \pi }{A}
 -  
\frac{\alpha}{A}, \
\frac{p \pi}{A}
-
\frac{\alpha}{A}
\Big[
\qquad 
(p \in \Z, \hsx p \neq 0, 1)
\]
contains exactly one zero of 
$W_{A, \alpha}$ 
and it is simple (cf. 32.19).
\\

\qquad \textbullet \quad
\un{$p = 0$}:\quad
$
\ds
I_0
\ = \ 
\Big]
-
\frac{\pi}{A} - \frac{\alpha}{A}, 
\hsx
-\frac{\alpha}{A} 
\Big[
\hsx .$   
\
If $0  < \alpha < \pi$, then 
\[
(-1)^1 
\hsx 
W_{A, \alpha}
\Big(
-\frac{\pi}{A} - \frac{\alpha}{A}
\Big)
\ > \ 
0
\]
\qquad 
$\implies$
\[
W_{A, \alpha}
\Big(
-\frac{\pi}{A} - \frac{\alpha}{A}
\Big)
\ < \ 
0.
\]
Meanwhile, 
\[
W_{A, \alpha}
\Big(- \frac{\alpha}{A}
\Big)
\ > \ 
0.
\]
\begin{spacing}{1.75}
\noindent 
So 
$W_{A, \alpha}$ 
has a (unique) zero in $I_0$ and it is simple (cf. 32.19).  
If $\alpha = 0$, then 
$
\ds
W_{A, 0}
\Big(-\frac{0}{A}
\Big)
 = 0 
$
and its multiplicity is equal to 1 and there are no zeros of 
$W_{A, 0}$ 
in 
$I_0 \cup I_1$ 
(cf. 32.20).
\\[-.25cm]
\end{spacing}

\qquad \textbullet \quad
\un{$p = 1$}:\quad
$
\ds
I_1 
\ = \ 
\Big]
-
\frac{\alpha}{A}, 
\
\frac{\pi}{A} - \frac{\alpha}{A}
\Big[
$
.  
\ 
In this situation, 
\[
x \in I_1 
\implies 
W_{A, \alpha} (x) 
\ > \ 
0 
\qquad (n = 0), 
\]
thus in $I_1$, 
$W_{A, \alpha}$ 
is zero free.
\\[-.25cm]

[Note: \ 
$
\ds
\frac{k \hsy \pi - \alpha}{A}
$ 
is a zero of 
$W_{A, \alpha}$ 
only when $k = 0$, $\alpha = 0$.]
\\[-.25cm]
\end{x}

\begin{x}{\small\bf THEOREM} \ 
If $\phi \notin E (A, \alpha)$, then all the zeros of 
$W_{A, \alpha}$ 
are real and simple.
\\[-.25cm]

PROOF \ 
The idea is to reduce things to the bounded case, i.e., to 33.8.  
To this end, for $n > 1$, let
\[
\phi_n (t) 
\ = \ 
\phi (t) 
\qquad 
\Big(0 < t \leq A - \frac{1}{n}\Big)
\]
and 
\[
\phi_n (t) 
\ = \ 
\phi
\Big(A - \frac{1}{n}\Big)
\hsx + \hsx 
t - A + \frac{1}{n} 
\qquad 
\Big(A - \frac{1}{n} \leq t < A\Big).
\]
Then 
$\phi_n \notin E (A, \alpha)$ 
and
\allowdisplaybreaks
\begin{align*}
\int\limits_0^A \ 
\abs{\phi (t)  - \phi_n (t)}
\ \td t
&=\ 
\int\limits_{A - \frac{1}{n}}^A \
\abs{\phi (t)  - \phi_n (t)}
\ \td t
\\[15pt]
&\leq 
\int\limits_{A - \frac{1}{n}}^A \
\abs{\phi (t)}
\hsx \td t
\hsx + \hsx 
\frac{1}{2 n^2}
\\[15pt]
&\ra 
0
\qquad (n \ra \infty).
\end{align*}
Put
\[
W_{A, \alpha, n} (z) 
\ = \ 
\int\limits_0^A \ 
\phi_n (t) 
\hsx 
\sin (z \hsy t + \alpha) 
\ \td t.
\]
Then 
$W_{A, \alpha, n} \ra W_{A, \alpha}$ 
uniformly on compact subsets of $\Cx$.  
On the other hand, $\phi_n$ is bounded on 
$]0, A[$, hence
\[
\sup\limits_{x \in \R} \ 
\abs{x \hsy W_{A, \alpha, n}  (x)} 
\ < \ 
\infty
\qquad (\tcf. \ 33.3).
\]
Therefore all the zeros of 
$W_{A, \alpha, n}$ are real and simple (cf. 33.8), so all the zeros of 
$W_{A, \alpha}$ 
are real and it remains to establish their simplicity.
\\[-.25cm]

\qquad \textbullet \quad
\un{$0 \hsx < \hsx \alpha \hsx < \hsx \pi$}: \ 
Given $p \in \Z$, let $D_p$ be the rectangle 
\[
\Big\{z : \abs{\Img z} \leq 1, 
\ 
\frac{(p-1)\pi}{A} - \frac{\alpha}{A} 
\ \leq \  \Reg z \ \leq \ 
\frac{p \pi}{A} - \frac{\alpha}{A} 
\Big\}.
\]
Then for $z \in \partial D_p$ and $n \gg 0$, 
\[
\abs{W_{A, \alpha, n} (z) - W_{A, \alpha} (z)} 
\ < \ 
\min\limits_{\partial D_p} \ 
\abs{W_{A, \alpha}}
\ \leq \ 
\abs{W_{A, \alpha} (z)}.
\]
But this implies by Rouch\'e that 
$W_{A, \alpha}$ 
and 
$W_{A, \alpha, n}$ 
have the same number of zeros inside $D_p$.  
\\[-.25cm]

\qquad \textbullet \quad
\un{$0 \hsx = \hsx  \alpha$}: \ 
At level 0, 1, work with $D_0 \cup D_1$ rather then $D_0$ and $D_1$ separately. 
\\[-.25cm]
\end{x}

Implicit in the foregoing is a description of the position of the zeros of 
$W_{A, \alpha}$ 
(what was said in the proof of 33.8 is valid in general).  

\begin{x}{\small\bf EXAMPLE} \ 
By definition, 
\[
W_{1, \frac{\pi}{2}} (z)
\ = \ 
\int\limits_0^1 \ 
\phi (t) 
\hsx 
\cos (z \hsy t) 
\ \td t.
\]
Assuming that $\phi \notin E(0,1)$ 
(a restriction that is actually unnecessary \ldots), 
the theory predicts that all the zeros of 
$W_{1, \frac{\pi}{2}}$ 
are real.  
As for their position, 
$W_{1, \frac{\pi}{2}}$ 
has a zero in each of the ambient intervals
\[
I_2 
\ = \ 
\Big]
\frac{\pi}{2}, \hsx \frac{3 \pi}{2}
\Big[, 
\quad
I_3 
\ = \ 
\Big]
\frac{3 \pi}{2}, \hsx \frac{5 \pi}{2}
\Big[, 
\quad
I_4 
\ = \ 
\Big]
\frac{5 \pi}{2}, \hsx \frac{7 \pi}{2}
\Big[, 
\quad
\ldots
\]
and this zero is unique and simple.  
Moreover, 
\[
C \Big(\frac{\pi}{2}\Big) \ > \ 0, 
\quad
C \Big(\frac{3 \pi}{2}\Big) \ < \ 0, 
\quad
C \Big(\frac{5 \pi}{2}\Big) \ > \ 0, 
\quad
C \Big(\frac{7 \pi}{2}\Big) \ < \ 0, 
\quad
\ldots
\]
and 
$
\ds
I_1 = 
\Big]- \frac{\pi}{2}, \hsx \frac{\pi}{2}\Big[
$ 
is zero free.  
All the positive zeros of 
$W_{1, \frac{\pi}{2}}$ 
are thereby accounted 
\\[-.5cm]

\noindent
for so 31.1 has been recovered.
\\[-.25cm]
\end{x}

\begin{x}{\small\bf LEMMA} \ 
We have
\[
\int\limits_0^A \ 
\phi (t) 
\hsx 
\cos (z \hsy t + \alpha) 
\ \td t 
\ = \quad
\begin{cases}
\ \ds
W_{A, \alpha \hsy + \hsy \frac{\pi}{2}}
\hspace{0.75cm} 
\Big(0 \leq \alpha < \frac{\pi}{2}\Big)
\\[18pt]
\ \ds
- W_{A, \alpha \hsy - \hsy \frac{\pi}{2}}
\hspace{0.5cm} 
\Big(\frac{\pi}{2} \leq \alpha < \pi\Big)
\end{cases}
.
\]
\\[-.25cm]
\end{x}


\chapter{
$\boldsymbol{\S}$\textbf{34}.\quad  ZEROS OF $f_A$}
\setlength\parindent{2em}
\setcounter{theoremn}{0}
\renewcommand{\thepage}{\S34-\arabic{page}}

\begin{x}{\small\bf NOTATION} \ 
Given $\phi \in \Lp^1 [0,A]$, put
\[
f_A (z)
\ = \ 
\int\limits_0^A \ 
\phi (t) 
\hsx 
e^{\sqrt{-1} \hsx z \hsy t}
\ \td t,
\]
thus
\[
f_A (z)
\ = \ 
C_A (z)+ \sqrt{-1} \hsx S_A (z),
\]
\[
\text{where} \hspace{2cm}
\begin{cases}
\ \ds
C_A (z)
\ = \ 
\int\limits_0^A \ 
\phi (t) 
\hsx 
\cos z \hsy t
\ \td t
\\[18pt]
\ \ds
S_A (z)
\ = \ 
\int\limits_0^A \ 
\phi (t) 
\hsx 
\sin z \hsy t
\ \td t
\end{cases}
.
\]

[Note: \ 
To be in agreement with \S30, drop the ``A'' if $A =  1$.]
\\[-.5cm]
\end{x}

\begin{x}{\small\bf THEOREM} \ 
If $\phi \in \Lp^1 [0,A]$ is positive and increasing on $]0,A[$ and if $\phi$ is not a step function, 
then the zeros of $f_A (z)$ lie in the open upper half-plane.  
\\[-.5cm]
\end{x}

\qquad
{\small\bf \un{N.B.}} \ 
Since $\phi$ is not a step function, it follows that $\forall \ \alpha$,
\[
\phi \notin E (A, \alpha).
\]
\\[-1.5cm]

\noindent
Therefore all the zeros of $W_{A, \alpha}$ are real and simple (cf. 33.9) and this persists to all 
$\alpha \in \R$ (elementary verification).
\\[-.5cm]

\begin{x}{\small\bf REMARK} \ 
Take $A = 1$ $-$then this result implies 29.16 (granted 29.19).
\\[-.5cm]
\end{x}

Let $P$ and $Q$ be nonconstant real entire functions.
\\[-.5cm]

\begin{x}{\small\bf CHEBOTAREV CRITERION} \ 
Assume: \ 
\\[-.5cm]

\qquad \textbullet \quad
$P$ and $Q$ have no common zeros.
\\[-.5cm]

\qquad \textbullet \quad
$\forall \ \mu, \nu \in \R$, 
$\mu^2 + \nu^2 \neq 0$, 
the combination  
$\mu \hsy P + \nu \hsy Q$ has no zeros in 
$\Cx - \R$.  
\\[-.5cm]

\qquad \textbullet \quad
$\exists \ x_0 \in \R$ such that 
\[
P (x_0) \hsy Q^\prime (x_0) 
\hsx - \hsx 
Q (x_0) 
\hsy
P^\prime (x_0)  
\ > \ 
0.
\]
\\[-.25cm]
Then
\[
F (z) 
\ = \ 
P (z) + \sqrt{-1} \hsz Q (z)
\]
has all its zeros in the open upper half-plane.   
\\[-.5cm]

[Note: \ 
It is an a posteriori conclusion that $\forall \ x \in \R$,
\[
P (x) \hsy Q^\prime (x) 
\hsx - \hsx 
Q (x)
\hsy
P^\prime (x)   
\ > \ 
0.]
\]
\\[-1.25cm]
\end{x}

\begin{x}{\small\bf REMARK} \ 
Compare the above with what has been said in \S16: \ 
There it was a question of nonconstant real polynomials in the open lower half-plane, 
hence the sign switch to 
\[
Q (x_0) \hsy P^\prime (x_0)
\hsx - \hsx 
P (x_0) \hsy Q^\prime (x_0) 
\ > \ 
0.
\]
\\[-1.25cm]
\end{x}

\qquad
{\small\bf \un{N.B.}} \ 
It is clear that $F (z)$ has no zeros on the real axis: 
\[
F (x_0) 
\ = \ 
P (x_0) + \sqrt{-1} \hsz Q (x_0) 
\ = \ 
0
\quad 
\implies 
\quad
\begin{cases}
\ 
P (x_0)
\ = \ 
0
\\[8pt]
\ 
Q (x_0)
\ = \ 
0
\end{cases}
.
\]
\\[-1.25cm]

Proceeding to the proof, begin by noting that the meromorphic function 
\[
\theta (z)
\ = \ 
\frac{Q (z)}{P (z)}
\]
does not take on real values for $\Img z \neq 0$, thus it maps the open upper half-plane 
either onto itself or onto the open lower half-plane.  
But
\[
P (x_0) \hsy Q^\prime (x_0) 
\hsx - \hsx 
Q (x_0) \hsy P^\prime (x_0)
\ > \ 
0
\]
\qquad 
$\implies$
\[
\theta^\prime (x_0) 
\ > \ 
0,
\]
so $\theta (z)$  maps the open upper half-plane onto itself.  
Since
\[
\frac{P + \sqrt{-1} \hsz Q}{P - \sqrt{-1} \hsz Q}
\ = \ 
\frac{1 + \sqrt{-1} \hsx \theta}{1 - \sqrt{-1} \hsx \theta}
\
, 
\]
it then follows that 
\[
\Img z > 0 
\ \implies \ 
\bigg| 
\hsx
\frac{P (z) + \sqrt{-1} \hsz Q (z)}{P (z) - \sqrt{-1} \hsz Q (z)}
\hsx
\bigg|
\ < \ 
1.
\] 
Next
\[
\begin{cases}
\ 
P (\bar{z}) 
\ = \ 
\ovs{P (z)}
\\[4pt]
\ 
Q (\bar{z}) 
\ = \ 
\ovs{Q (z)}
\end{cases}
,
\]
hence
\[
P (z_0) + \sqrt{-1} \hsz Q (z_0) 
\ = \ 
0
\]
\qquad \qquad
$\implies$
\[
P (\bar{z}_0) - \sqrt{-1} \hsz Q (\bar{z}_0) 
\ = \ 
0.
\]
Accordingly, it need only be shown that 
$P - \sqrt{-1} \hsz Q$ 
has no zeros in the open upper half-plane. 
However

\[
\frac{P (z) + \sqrt{-1} \hsz Q (z)}{P (z) - \sqrt{-1} \hsz Q (z)}
\] 
is unbounded near any zero of 
$P - \sqrt{-1} \hsz Q$ 
which is not a zero of 
$P + \sqrt{-1} \hsz Q$.  
And this means that any zero of 
$P - \sqrt{-1} \hsz Q$ 
in the open upper half-plane
must be a zero of 
$P + \sqrt{-1} \hsz Q$.  
But

\[
\begin{cases}
\ 
P (z_0) - \sqrt{-1} \hsz Q (z_0) 
\ = \ 
0
\\[4pt]
\ 
P (z_0) + \sqrt{-1} \hsz Q (z_0) 
\ = \ 
0
\end{cases}
\qquad (\Img z_0 > 0)
\]
\[
\implies
\qquad
\begin{cases}
\ 
2 P (z_0) \ = \ 0 
\hspace{1.75cm}
\implies 
P(z_0) \ = \ 0
\\[4pt]
\ 
- 2 \hsy \sqrt{-1} \hsx  Q (z_0) \ = \ 0 
\hspace{0.43cm}
\implies 
Q(z_0) \ = \ 0
\end{cases}
,
\hspace{2cm}
\]
contradicting the assumption that $P$ and $Q$ have no common zeros. 
\\[-.25cm]

Having dispensed with the preparation, we are now in a position to give the proof of 34.2.  
Bearing in mind that 
\[
f_A (z) 
\ = \ 
C_A (z) + \sqrt{-1} \hsx  S_A (z), 
\]
start by writing 
\[
W_{A, \alpha}  (z)
\ = \ 
(\sin \alpha) \hsy C_A (z) 
\hsx + \hsx 
(\cos \alpha) \hsy S_A (z) .
\]
Then there are three items to be checked.
\\[-.25cm]

1. \quad
$C_A$ and $S_A$ have no common zeros.  
To see this, observe that 
\[
\begin{cases}
\ \ds
W_{A, \frac{\pi}{2}} (z)
\ = \ 
C_A (z)
\\[4pt]
\ \ds
\ 
W_{A, 0} (z)
\ = \ 
S_A (z)
\end{cases}
,
\]
so the zeros of $C_A (z)$ and $S_A (z)$ are real and simple.  
If 
$C_A (x_0) = 0$, 
$S_A (x_0) = 0$ 
for some $x_0 \in \R$, then 
$C_A^\prime (x_0) \neq 0$, 
$S_A^\prime (x_0) \neq0$ 
and taking 
\[
\alpha 
\ = \ 
\arctan 
\Big(- \frac{S_A^\prime (x_0) }{C_A^\prime (x_0) }\Big),
\]
we have 
\allowdisplaybreaks
\begin{align*}
W_{A, \alpha}^\prime (x_0) \ 
&=\
(\sin \alpha) \hsy C_A^\prime (x_0) 
\hsx + \hsx 
(\cos \alpha) \hsy S_A^\prime (x_0) 
\\[11pt]
&=\ 
0
\end{align*}
for a suitable choice of arctan.  
But this implies that $x_0$ is  a zero of $W_{A, \alpha}$ of multiplicity $\geq 2$ which cannot be.
\\[-.25cm]

2. \quad
$\forall \ \mu, \nu \in \R$, 
$\mu^2 + \nu^2 \neq 0$, 
the combination  
$\mu \hsy C_A + \nu \hsy S_A$ has no zeros in 
$\Cx - \R$.  
The cases 
$\mu \neq 0$, 
$\nu = 0$, 
and 
$\mu = 0$, 
$\nu \neq 0$ 
being obvious, consider the remaining four possibilities.
\\[-.25cm]

\qquad \textbullet \quad 
$\un{\mu > 0, \hsx \nu > 0}:$ \quad 
Write
\[
\mu \hsy C_A \hsx + \hsx \nu \hsy S_A 
\ = \ 
\sqrt{\mu^2 + \nu^2} 
\
\bigg(
\frac{\mu}{\sqrt{\mu^2 + \nu^2}} 
\ 
C_A 
\ + \ 
\frac{\nu}{\sqrt{\mu^2 + \nu^2}} 
\ 
S_A 
\bigg)
\]
and determine $\alpha$ by
\[
\sin \alpha
\ = \ 
\frac{\mu}{\sqrt{\mu^2 + \nu^2}}, 
\quad 
\cos \alpha
\ = \ 
\frac{\nu}{\sqrt{\mu^2 + \nu^2}}
\hsx
.
\]
\\[-.25cm]

\qquad \textbullet \quad 
$\un{\mu < 0, \hsx \nu < 0}:$ \quad 
Write
\[
\mu \hsy C_A \hsx + \hsx \nu \hsy S_A 
\ = \ 
-\sqrt{\mu^2 + \nu^2} 
\
\bigg(
\frac{-\mu}{\sqrt{\mu^2 + \nu^2}} 
\ 
C_A 
\ + \ 
\frac{-\nu}{\sqrt{\mu^2 + \nu^2}} 
\ 
S_A 
\bigg)
\]
and determine $\alpha$ by
\[
\sin \alpha
\ = \ 
\frac{-\mu}{\sqrt{\mu^2 + \nu^2}}, 
\quad 
\cos \alpha
\ = \ 
\frac{-\nu}{\sqrt{\mu^2 + \nu^2}}
\hsx
.
\]
\\[-.25cm]

\qquad \textbullet \quad 
$\un{\mu < 0, \hsx \nu > 0}:$ \quad 
Write
\allowdisplaybreaks
\begin{align*}
\mu \hsy C_A \hsx + \hsx \nu \hsy S_A  \ 
&= \ 
\sqrt{\mu^2 + \nu^2} 
\
\bigg(
-
\frac{-\mu}{\sqrt{\mu^2 + \nu^2}} 
\ 
C_A 
\ + \ 
\frac{\nu}{\sqrt{\mu^2 + \nu^2}} 
\ 
S_A 
\bigg)
\\[15pt]
&= \ 
\sqrt{\mu^2 + \nu^2} 
\
((- \sin \alpha ) C_A 
\hsx + \hsx 
(\cos \alpha) S_A)
\\[15pt]
&= \ 
\sqrt{\mu^2 + \nu^2} 
\
((\sin -\alpha ) C_A 
\hsx + \hsx 
(\cos -\alpha) S_A).
\end{align*}
\\[-.75cm]

\qquad \textbullet \quad 
$\un{\mu > 0, \hsx \nu < 0}:$ \quad 
Write
\allowdisplaybreaks
\begin{align*}
\mu \hsy C_A \hsx + \hsx \nu \hsy S_A  \ 
&= \ 
\sqrt{\mu^2 + \nu^2} 
\
\bigg(
\frac{\mu}{\sqrt{\mu^2 + \nu^2}} 
\ 
C_A 
\ - \ 
\frac{-\nu}{\sqrt{\mu^2 + \nu^2}} 
\ 
S_A 
\bigg)
\\[15pt]
&= \ 
\sqrt{\mu^2 + \nu^2} 
\
((\sin \alpha ) C_A 
\hsx - \hsx 
(\cos \alpha) S_A)
\\[15pt]
&= \ 
\sqrt{\mu^2 + \nu^2} 
\
(-(\sin \hsy -\alpha ) C_A 
\hsx - \hsx 
(\cos \hsy -\alpha) S_A)
\\[15pt]
&= \ 
-\sqrt{\mu^2 + \nu^2} 
\
((\sin \hsy -\alpha ) C_A 
\hsx + \hsx 
(\cos ( \hsy-\alpha)) S_A)
\hsx
.
\end{align*}
\\[-.25cm]

3. \quad
$\exists \ x_0 \in \R$ such that

\[
C_A  (x_0) \hsy S_A^\prime (x_0) 
\hsx - \hsx 
S_A (x_0) \hsy C_A ^\prime (x_0)
\ \neq \ 
0.
\]
In fact, 
\begin{align*}
\allowdisplaybreaks
C_A  (0) \hsy S_A^\prime (0) 
\hsx - \hsx 
S_A (0) \hsy   C_A ^\prime (0)
\ 
&=\ 
C_A  (0) \hsy S_A^\prime (0) 
\\[11pt]
&=\ 
\bigg(
\int\limits_0^1 \ 
\phi (t) 
\ \td t
\bigg)
\ 
\bigg(
\int\limits_0^1 \ 
\phi (t) 
\hsx
t
\ \td t
\bigg)
\\[11pt]
&>\ 
0.
\end{align*}

\begin{x}{\small\bf REMARK} \ 
If $\phi$ is a step function and if $\phi \in E(A, \alpha)$, 
then $f_A (z)$ has an infinity of real zeros (cf. 29.21) (all of which are simple) 
and there is an analog of 29.22.
\\[-.25cm]
\end{x}

\begin{x}{\small\bf NOTATION} \ 
Given $\phi \in \Lp^1 [0,A]$, let
\[
\begin{cases}
\ \ds
\gC_A (z)
\ = \ 
\int\limits_0^A \ 
\phi (A - t) 
\cos z \hsy t 
\ \td t
\\[26pt]
\ \ds
\ggS_A (z)
\ = \ 
\int\limits_0^A \ 
\phi (A - t) 
\sin z \hsy t 
\ \td t
\end{cases}
.
\]
\\[-.75cm]
\end{x}


\begin{x}{\small\bf IDENTITIES} \ 
\[
f_A (z) \hsx e^{-\sqrt{-1} \hsx A \hsy z}
\ = \ 
\gC_A (z) - \sqrt{-1} \hsx \ggS_A (z)
\]
and
\[
\begin{cases}
\ 
C_A (z) 
\ = \ 
\gC_A (z) \hsy \cos A z \hsx + \hsx \ggS_A (z) \hsy \sin A z
\\[8pt]
\ 
S_A (z) 
\ = \ 
\gC_A (z) \hsy \sin A z \hsx - \hsx \ggS_A (z) \hsy \cos A z
\end{cases}
.
\]
\\[-.75cm]
\end{x}

\begin{x}{\small\bf RAPPEL} \ 
If 0 and $A$ are the effective limits of integration 
(thus excluding the possibility that $\phi = 0$ almost everywhere), 
then $f_A (z)$ has an infinity of zeros (see the initial comments in \S29).
\\[-.25cm]
\end{x}

\begin{x}{\small\bf LEMMA} \ 
Put
\[
H (s) 
\ = \ 
-
\frac{y}{\pi (y^2 + s^2)}
\qquad (y \in \R).
\]
Then
\[
\int\limits_{-\infty}^\infty \ 
e^{\sqrt{-1} \hsx s \hsy t}
\hsx 
H(s) 
\ \td s
\ = \ 
e^{y \hsy \abs{t}}.
\]
\\[-.25cm]
\end{x}

\begin{x}{\small\bf THEOREM} \ 
If $\phi \in \Lp^1 [0,A]$ is real and if 
\[
\gC_A (x) 
\ \geq \ 
0 
\qquad (x \in \R),
\]
then $f_A (z)$ has no zeros in the open lower half-plane.
\\[-.5cm]

PROOF \ 
Let $z = x + \sqrt{-1} \hsx y$  $(y < 0)$ and write
\allowdisplaybreaks
\begin{align*}
f_A (z)
\hsx 
e^{-\sqrt{-1} \hsx A \hsy z}\ 
&=\ 
\int\limits_0^A \ 
\phi (t) 
\hsx 
e^{\sqrt{-1} \hsx z \hsy t}
\hsx 
e^{-\sqrt{-1} \hsx A \hsy z}
\ \td t
\\[11pt]
&=\ 
\int\limits_0^A \ 
\phi (t) 
\hsx 
e^{\sqrt{-1} \hsx z (t - A)}
\ \td t
\\[11pt]
&=\ 
\int\limits_0^A \ 
\phi (t) 
\hsx 
e^{-\sqrt{-1} \hsx z (A -t)}
\ \td t
\\[11pt]
&=\ 
\int\limits_0^A \ 
\phi (A - t) 
\hsx 
e^{-\sqrt{-1} \hsx x \hsy t}
\hsx
e^{y \hsy t}
\ \td t
\\[11pt]
&=\ 
\int\limits_0^A \ 
\phi (A - t) 
\hsx 
e^{-\sqrt{-1} \hsx x \hsy t}
\ 
\bigg(
\int\limits_{-\infty}^\infty \ 
e^{\sqrt{-1} \hsx s \hsy t} 
\hsx 
H (s) 
\ \td s
\bigg)
\\[11pt]
&=\ 
\int\limits_{-\infty}^\infty \ 
H (s) 
\ 
\bigg(
\int\limits_0^A \ 
e^{\sqrt{-1} \hsx (s - x)t }\ 
\hsx
\phi (A - t) 
\ \td t
\bigg)
\ \td s
\\[11pt]
&=\ 
\int\limits_{-\infty}^\infty \ 
H (s + x) 
\big(
\gC_A (s) 
\hsx + \hsx
\sqrt{-1} \hsx 
\ggS_A (s)
\big)
\ \td s.
\end{align*}
But 
$\gC_A  \not\equiv 0$ 
(consult the Appendix below), hence
\allowdisplaybreaks
\begin{align*}
\Reg \big(f_A (z) \hsx e^{-\sqrt{-1} \hsx A \hsy z}\big) \ 
&=\ 
-
\int\limits_{-\infty}^\infty \ 
\frac{1}{\pi (y^2 + (s + x)^2)}
\
\gC_A (s) 
\ \td s
\\[15pt]
&> \ 
0.
\end{align*}
\\[-1.25cm]
\end{x}

\begin{x}{\small\bf REMARK} \ 
Any real zero of $f_A (z)$ (if there is one) is necessarily simple.
\\[-.25cm]
\end{x}

\begin{x}{\small\bf EXAMPLE} \ 
If $\phi \in C [0,A]$ is real, 
$\phi (0) = 0$, 
$\phi (A) > 0$, 
and the function 
\[
t \ra \phi ((A - \abs{t})_+)
\]
is positive definite on $\R$, then 
\[
\gC_A (x) 
\ \geq \ 
0 
\qquad (x \in \R),
\]
so 34.11 is applicable.
\\[-.25cm]
\end{x}


\[
\text{APPENDIX}
\]

\qquad
{\small\bf M\"UNTZ CRITERION} \ 
If $\lambda_1, \lambda_2, \ldots$ is a strictly increasing sequence of real numbers such that 
\[
\sum\limits_{n = 1}^\infty \
\frac{1}{\lambda_n}
\ = \ 
\infty,
\]
then the set
\[
\{1, t^{\lambda_1}, t^{\lambda_2}, \ldots \}
\]
is total in $C [0,1]$.
\\

\qquad
{\small\bf EXAMPLE} \ 
The set
\[
\{t^0, t^2, t^4, \ldots \}
\]
is total in $C [0,1]$.
\\

\qquad
{\small\bf APPLICATION} \ 
If $\psi \in \Lp^1 [0,1]$ and if 
\[
\begin{cases}
\ \ \ds
\int\limits_0^1 \ 
\psi (t) 
\ \td t
\ = \ 0
\\[26pt]
\ \ \ds
\int\limits_0^1 \ 
t^{2 k}
\hsx
\psi (t) 
\ \td t
\ = \ 0
\end{cases}
\qquad 
(k = 1, 2, \ldots),
\]
then $\psi = 0$ almost everywhere.
\\[-.5cm]

[Let
\[
\Psi (t) 
\ = \ 
\int\limits_0^t \ 
\psi (s) 
\ \td s.
\]
Then $\Psi$ is absolutely continuous and 
$\Psi (0) = 0$, 
$\Psi (1) = 0$. 
Now integrate by parts to get
\allowdisplaybreaks
\begin{align*}
0 \ 
&=\ 
\int\limits_0^1 \ 
t^{2 k}
\hsx
\psi (t) 
\ \td t
\\[15pt]
&=\ 
- 2 k \ 
\int\limits_0^1 \ 
t^{2 k-1}
\hsx
\Psi (t) 
\ \td t
\qquad (k = 1, 2, \ldots).
\end{align*}
Therefore
\allowdisplaybreaks
\begin{align*}
&
\int\limits_0^1 \ 
t^0 
\hsx
(t \hsy \Psi (t)) 
\ \td t 
\ = \ 0
\qquad (k = 1)
\\[15pt]
&
\int\limits_0^1 \ 
t^2 
\hsx
(t \hsy \Psi (t)) 
\ \td t 
\ = \ 0
\qquad (k = 2)
\\[15pt]
&
\int\limits_0^1 \ 
t^4 
\hsx
(t \hsy \Psi (t)) 
\ \td t 
\ = \ 0
\qquad (k = 3)
\\[8pt]
&
\hspace{2cm}
\vdots
\end{align*}
Define a bounded linear functional $\mu$ on $C [0,1]$ by the rule 
\[
\mu (g) 
\ = \ 
\int\limits_0^1 \ 
g (t)
\hsx
(t \hsy \Psi (t)) 
\ \td t.
\]
Then 
\[
\mu (t^{2 \hsy k}) 
\ = \ 
0 
\qquad (k = 0, 1, 2, \ldots)
\]
\qquad 
$\implies$
\[
\mu \ \equiv \ 0
\]
\qquad 
$\implies$
\[
t \hsy \Psi (t) \ = \  0
\quad (0 \leq t \leq 1)
\]
\qquad 
$\implies$
\[
\Psi (t) \ = \  0
\quad (0 \leq t \leq 1).
\]
But this implies that $\psi = 0$ almost everywhere.]
\\

\qquad
{\small\bf THEOREM}  
If $C_A (z) \equiv 0$, then $\phi = 0$ almost everywhere  
($\implies f_A( z) \equiv 0$).  
\\[-.5cm]

PROOF \ 
Consider the expansion 
\allowdisplaybreaks
\begin{align*}
\int\limits_0^A \ 
\phi (t) 
\hsx 
\cos z \hsy t
\ \td t\ 
&=\ 
\int\limits_0^A \ 
\phi (t) 
\  
\sum\limits_{k = 0}^\infty \ 
\frac{(-1)^k \hsy (z \hsy t)^{2 \hsy k}}{(2 k)!}
\ \td t
\\[15pt]
&=\ 
\sum\limits_{k = 0}^\infty \ 
\frac{(-1)^k \hsy }{(2 k)!}
\ 
\bigg(
\int\limits_0^A \ 
t^{2 \hsy k} 
\hsx
\phi (t) 
\ \td t\ 
\bigg)
\ 
z^{2 \hsy k}, 
\end{align*}
\\[-.25cm]
hence
\[
\int\limits_0^A \ 
t^{2 \hsy k} 
\hsx
\phi (t) 
\ \td t
\ = \ 
0
\qquad (k = 0, 1, 2, \ldots)
\]
or still (letting $t = s \hsy A$), 
\[
A^{2 k + 1} \ 
\int\limits_0^1 \ 
s^{2 k} 
\hsx 
\phi (s \hsy A)
\ \td s 
\ = \ 
0
\qquad (k = 0, 1, 2, \ldots).
\]
Consequently, $\phi (s \hsy A)$ vanishes almost everywhere 
$(0 \leq s \leq 1)$, so $\phi (t)$ vanishes almost everywhere 
$(0 \leq t \leq A)$.
\\[-.25cm]

\qquad
{\small\bf \un{N.B.}} \ 
If $\gC_A (z) \equiv 0$, then $\phi = 0$ almost everywhere 
($\implies f_A( z) \equiv 0$) 
(argue analogously). 
\\[-.25cm]

\qquad
{\small\bf REMARK} \ 
If $f_A( z) \equiv 0$, then $\phi = 0$ almost everywhere.
\\[-.5cm]

[In fact, 
\allowdisplaybreaks
\begin{align*}
C_A (z) \ 
&=\ 
\int\limits_0^A \ 
\phi (t) 
\hsx 
\cos z \hsy t
\ \td t
\\[15pt]
&=\ 
\int\limits_0^A \ 
\phi (t) 
\hsx 
\frac{e^{\sqrt{-1} \hsx z \hsy t} \hsx + \hsx e^{-\sqrt{-1} \hsx z \hsy t}}{2}
\ \td t
\\[15pt]
&=\ 
\frac{f_A (z) + f_A (-z)}{2}
\\[15pt]
&\equiv\ 
0.]
\end{align*}
\\[-.25cm]


\chapter{
$\boldsymbol{\S}$\textbf{35}.\quad  MISCELLANEA}
\setlength\parindent{2em}
\setcounter{theoremn}{0}
\renewcommand{\thepage}{\S35-\arabic{page}}

\qquad
Here there will be found a number of complements, some theoretical, others disguised as ``examples''.
\\[-.25cm]

\begin{x}{\small\bf LEMMA} \ 
If $\phi \in \Lp^1 [0,A]$ 
is real valued and continuously differentiable and if $\phi (A) \neq 0$, then
\[
C_A (z) 
\ = \ 
\int\limits_0^A \ 
\phi (t) 
\hsx 
\cos z t 
\ \td t
\]
has an infinite number of real zeros.  
\\[-.5cm]

PROOF \ 
In fact, 
\allowdisplaybreaks
\begin{align*}
x \hsy C_A (x) \ 
&= \ 
\phi (A) 
\hsy
\sin (x A) 
\ - \ 
\int\limits_0^A \ 
\phi^\prime (t) 
\hsx 
\sin (x \hsy t) 
\ \td t
\\[15pt]
&=\ 
\phi (A) 
\hsy
\sin (x \hsy A) 
\ + \ 
\txo (1) 
\qquad (\abs{x} \ra \infty).
\end{align*}
\\[-1.25cm]
\end{x}

\begin{x}{\small\bf CHAKALOV CRITERION}\footnote[2]{\vspace{.11 cm}
{\fontencoding{OT2}\selectfont
Spianie p1AN}
\textbf{36} (1927), pp. 51-92.} \ 
Suppose given a sequence
\[
\cdots
\ < \ 
a_{-2}
\ < \ 
a_{-1}
\ < \ 
a_0
\ < \ 
a_1
\ < \ 
a_2
\ < \ 
\cdots
\]
and real numbers
\[
\ldots, 
\hsx
A_{-2}, 
\hsx
A_{-1}, 
\hsx
A_0, 
\hsx
A_1, 
\hsx
A_2, 
\hsx
\ldots, 
\]
where
\[
A_k 
\ \neq \ 
0, 
\qquad 
k = 0, \pm1, \pm2, \ldots \hsx.
\]
Assume: \ 
$\exists$ integers $p$ and $q$ with $p < q$ such that $A_k$ and $A_{k+1}$ 
have the same sign for $k < p$ and for $k \geq q$.  
Put 
\[
R_n (z) 
\ = \ 
\sum\limits_{k = -n + 1}^n \ 
\frac{A_k}{z - a_k}
\]
and impose the condition that
\[
R (z) 
\ = \ 
\lim\limits_{n \ra \infty} \ 
R_n (z)
\]
uniformly on compact subsets of 
$\Cx - \{a_k\}_{-\infty}^\infty$ 
$-$then $R (z)$ has no more than $q - p$ nonreal zeros.  
\\[-.5cm]
\end{x}

Maintaining the setup of 35.1, introduce the meromorphic function 
\[
R (z) 
\ = \ 
\frac{C_A (z)}{\cos (z A)}
\]
and put
\[
R_n (z) 
\ = \ 
\sum\limits_{k = -n + 1}^n \ 
(-1)^k
\ 
\frac
{\ds C_A \hsx \Big(\frac{k - \frac{1}{2} \pi}{A}\Big)}
{
\ds z - \frac{\Big(k - \frac{1}{2}\Big) \pi}{A}
}.
\]
Abbreviate
\[
\frac{\big(k - \frac{1}{2}\big) \pi}{A}
\quad \text{to} \quad 
a_k.
\]
\\[-1cm]

\begin{x}{\small\bf LEMMA} \ 
We have
\\[-.25cm]
\[
R (z) 
\ = \ 
\lim\limits_{n \ra \infty} \ 
R_n (z)
\]
uniformly on compact subsets of 
$\Cx - \{a_k\}_{-\infty}^\infty$.
\\[-.5cm]
\end{x}

Next
\allowdisplaybreaks
\begin{align*}
\lim\limits_{k \hsy \ra \hsy \pm \infty} \ 
(-1)^k \ 
a_k 
\hsx 
C_A (a_k)\ 
&=\ 
\phi (A) \ 
\lim\limits_{k \hsy \ra \hsy \pm \infty} \ 
(-1)^k \ 
\sin (a_k A)
\\[11pt]
&=\ 
\phi (A) \ 
\lim\limits_{k \hsy \ra \hsy \pm \infty} \ 
(-1)^k 
\
\sin
\bigg(
\frac{\ds \big(k - \frac{1}{2}\big) \pi}{A}
\hsx 
A
\bigg)
\\[11pt]
&=\ 
\phi (A) \ 
\lim\limits_{k \hsy \ra \hsy \pm \infty} \ 
(-1)^k 
\hsy
(-1)
\hsy
(-1)^k
\\[15pt]
&=\ 
- \phi (A) \ 
\\[11pt]
&\neq
0.
\end{align*}
If now
\[
A_k 
\ \equiv \ 
(-1)^k
\hsx 
C_A (a_k), 
\]
then the sequence
\[
\ldots, \hsx A_{-2}, \hsx A_{-1}, \hsx A_0,\hsx  A_1, \hsx A_2, \hsx \ldots
\]
has but a finite number of sign changes. 
\\[-.5cm]

[E.g.: \ 
Suppose that $L \equiv -\phi (A)$ is positive and send $k$ to $+\infty$ $-$then from some point on, 
$A_k$ is also positive: 
\allowdisplaybreaks
\begin{align*}
k \gg 0 
&\implies
\abs{a_k \hsy A_k - L} \ < \ \frac{L}{2}
\\[11pt]
&\implies
\frac{L}{2} \ < \ a_k \hsy A_k \ < \ \frac{3 L}{2}
\\[11pt]
&\implies
0 \ < \ \frac{L}{2 \hsy a_k} \ < \ A_k.]
\end{align*}

[Note: \ 
These considerations also serve to show that the number of $k$ for which $A_k = 0$ is finite.]
\\[-.25cm]

\begin{x}{\small\bf LEMMA} \ 
If $\phi \in \Lp^1 [0,A]$ 
is real valued and continuously differentiable and if $\phi (A) \neq 0$, then
\[
C_A (z) 
\ = \ 
\int\limits_0^A \ 
\phi (t) 
\hsy 
\cos z t 
\ \td t
\]
has at most a finite number of nonreal zeros.
\\[-.5cm]

[Thanks to what has been said above, one has only to invoke 35.2.]
\\[-.5cm]
\end{x}

\qquad
{\small\bf \un{N.B.}} \ 
Therefore
\[
C_A \in * - \sL - \sP
\qquad (\tcf. \ 10.36).
\]
\\[-1.25cm]

\begin{x}{\small\bf EXAMPLE} \ 
Take $\phi (t) = e^{-t}$ $-$then the zeros of 
\allowdisplaybreaks
\begin{align*}
C_A (z) \ 
&= \ 
\int\limits_0^A \ 
e^{- t}
\hsy 
\cos z t 
\ \td t
\\[15pt]
&=\ 
\frac{e^{-A} \hsy (z \sin A z \hsx - \hsx \cos A z) + 1}{z^2 + 1}
\\[15pt]
&=\ 
\frac{\sqrt{-1}}{2} \ 
\bigg[
\frac{e^{A ( - 1 \hsx - \sqrt{-1} \hsx z)} - 1}{z - \sqrt{-1}}
\ - \ 
\frac{e^{A ( - 1  \hsx + \sqrt{-1} \hsx z)} - 1}{z + \sqrt{-1}}
\bigg]
\end{align*}
\\
lie in the horizontal strip
\[
- 1 \ < \ y \ < \ 1
\qquad 
\Big(\tcf. \ 29.23 \ 
\Big(
\hsx
\Big|
\frac{\phi^\prime (t)}{\phi (t)}
\Big|
\ = \ 
1
\Big)\Big).
\]
The number of real zeros is infinite (cf. 35.1) while the number of nonreal zeros is finite (cf. 35.4).  
And the estimate 
$-1 < y < 1$ cannot be improved provided $A$ is allowed to vary, i.e., given $\varepsilon > 0$, in 
\[
-1 
\ < \ 
y
\ < \ 
-1 + \varepsilon 
\ \cup \ 
1 - \varepsilon
\ < \ 
y
\ < \ 
1
\]
there is a zero if $A \gg 0$.  
Finally, any compact subset $S$ of 
$-1 < y < 1$ is zero free for $A \gg 0$.
\ 
Proof:   
In $S$, 
\\
\[
\lim\limits_{A \ra \infty} \ 
\int\limits_0^A \ 
e^{-t}
\hsx 
\cos z \hsy t 
\ \td t
\ = \ 
\frac{1}{z^2 + 1}
\]
\\
and the function on the right has no zeros there.
\\[-.5cm]

[Note: \ 
As a function of $A$, the number of nonreal zeros is unbounded.]
\\[-.25cm]
\end{x}

\begin{x}{\small\bf NOTATION} \ 
(cf. 34.1) \ 
Given $\phi \in \L^1 (-\infty, \infty)$ put
\[
f_\infty (z) 
\ = \ 
\int\limits_{-\infty}^\infty \ 
\phi (t) 
\hsx 
e^{\sqrt{-1} \hsx z \hsy t}
\ \td t,
\]
thus
\[
f_\infty (z) 
\ = \ 
C_\infty (z) + \sqrt{-1} \hsx S_\infty (z),
\]
where

\[
\begin{cases}
\quad \ds
C_\infty (z) \ = \ \int\limits_{-\infty}^\infty \ \phi (t) \hsx \cos z t \ \td t
\\[18pt]
\quad \ds
S_\infty (z) \ = \ \int\limits_{-\infty}^\infty \ \phi (t) \hsx \sin z t \ \td t
\end{cases}
.
\]
\\[-1.25cm]
\end{x}

\qquad
{\small\bf \un{N.B.}} \ 
If $\phi$ is real and even (odd), then one can work instead with 
\[
\begin{cases}
\quad \ds
C_\infty (z) \ \equiv \ \int\limits_0^\infty \ \phi (t) \hsx \cos z t \ \td t
\\[18pt]
\quad \ds
S_\infty (z) \ \equiv \ \int\limits_0^\infty \ \phi (t) \hsx \sin z t \ \td t
\end{cases}
.
\]
\\[-1.25cm]

\begin{x}{\small\bf EXAMPLE} \ 
Suppose that $2 n$ is an even positive integer and take 
\[
\phi (t) 
\ = \ 
\exp (- t^{2 n})
\qquad (n = 1, 2, \ldots).
\] 
Then 
\[
\int\limits_{-\infty}^\infty \ 
\exp (- t^{2 n}) 
\hsx 
e^{\sqrt{-1} \hsz z \hsy t}
\ \td t
\ = \ 
\sqrt{\pi} \ 
\exp 
\Big(
- \frac{z^2}{4}
\Big) 
\]
has no zeros but
\[
\int\limits_{-\infty}^\infty \ 
\exp (-t^{4, \hsx 6, \hsx \ldots}) 
\hsx 
e^{\sqrt{-1} \hsz z \hsy t}
\ \td t
\]
has an infinity of real zeros though it has no complex zeros (cf. 12.34).
\\[-.5cm]

[Note:  \ 
Put 
\[
f_n (z) 
\ = \ 
\int\limits_{-\infty}^\infty \ 
\exp (- t^{2 n}) 
\hsx 
e^{\sqrt{-1} \hsz z \hsy t}
\ \td t
\qquad (n = 1, 2, \ldots).
\]
Then $f_n \in \sL - \sP$ is transcendental and satisfies the differential equation
\[
f_n^{(2n - 1)} (z) 
\ = \ 
\frac{(-1)^n}{2 n} 
\hsx 
z 
\hsy 
f_n (z).
\]
Therefore all the zeros of $f_n$ are simple (see the Appendix \hsy to \S13).]
\\[-.25cm]
\end{x}

\begin{x}{\small\bf REMARK} \ 
Consider
\[
\int\limits_0^A
\exp (- t^2) 
\hsx 
\cos z t 
\ \td t.
\]
Then 35.1 and 35.4 are applicable  and there is an $A$ with the property that 

\[
\int\limits_0^A
\exp (- t^2) 
\hsx 
\cos z t 
\ \td t
\]
has a nonreal zero (but no characterization is known of those $A$ for which this happens) 
(the situation in 35.5 is simpler although a complete explication is lacking there too).
\\[-.25cm]
\end{x}

\begin{x}{\small\bf EXAMPLE} \ 
The zeros of 
\[
\int\limits_{-\infty}^\infty \ 
\exp (-t^{4, \hsx 6, \hsx \ldots}) 
\hsx 
e^t
\hsy
e^{\sqrt{-1} \hsz z \hsy t}
\ \td t
\]
lie on the line $\Img z = 1$.
\\[-.5cm]

[If $z = a + \sqrt{-1} \hsz b$ is a zero, write
\[
e^t
\hsy
e^{\sqrt{-1} \hsz z \hsy t}
\ = \ 
e^{\sqrt{-1} \hsx (-  \sqrt{-1} \hsx + z) \hsy t},
\]
hence $-\sqrt{-1} + z$ is real, so $b = 1$.]
\\[-.25cm]
\end{x}

\begin{x}{\small\bf EXAMPLE} \ 
Fix $\alpha > 1$, $\alpha \neq 2 n$ $(n = 1, 2, \ldots)$, 
take 
$\phi (t) = \exp (-t^\alpha)$, 
and put

\[
\Phi_\alpha (z) 
\ = \ 
\int\limits_0^\infty \ 
\exp (- t^\alpha) 
\hsx 
\cos z t 
\ \td t.
\]
Then $\Phi_\alpha$ has an infinite number of nonreal zeros and a finite number of real zeros,
there being at least 
$
\ds
2 \hsx \Big[\frac{\alpha}{2}\Big]
$
of the latter if $\alpha > 2$.
\\[-.25cm]
\end{x}

\begin{x}{\small\bf LEMMA} \ 
We have
\[
\lim\limits_{x \ra \infty} \ 
x^{\alpha + 1} \hsx \Phi_\alpha (x) 
\ = \ 
\Gamma (\alpha + 1) \hsx \sin \Big(\frac{\pi \hsy \alpha}{2}\Big).
\]

PROOF \ 
There are seven steps.
\\[-.25cm]

\un{Step 1:} \ 
Integrate by parts to get
\[
x^{\alpha + 1} \hsx \Phi_\alpha (x) 
\ = \ 
x^\alpha \ 
\int\limits_0^\infty \
\sin x \hsy t 
\hsx \cdot \hsx
\alpha \hsy t^{\alpha - 1}
e^{- t^\alpha}
\ \td t.
\]

\un{Step 2:} \ 
Make the change of variable $u = x^\alpha  \hsy t^\alpha$, hence
\[
x^{\alpha + 1} \hsx \Phi_\alpha (x) 
\ = \ 
\int\limits_0^\infty \ 
\sin u^{1/\alpha}
\hsx \cdot \hsx 
e^{-x^{-\alpha} u}
\ \td u,
\]
a.k.a. the Laplace transform of $\sin u^{1/\alpha}$ at $x^{-\alpha}$.
\\[-.25cm]

\un{Step 3:} \ 
Rewrite the right hand side in terms of a complex exponential,  so
\[
x^{\alpha + 1} \hsx \Phi_\alpha (x) 
\ = \ 
\Img \ 
\int\limits_0^\infty \ 
\exp (\sqrt{-1} \hsx u^{1/\alpha} \hsx - x^{-\alpha} u)
\ \td u.
\]

\un{Step 4:} \ 
Move the contour of integration up to a straight line going from 0 to $\infty$ 
placed at a ``small'' angle $\theta$ to the positive real axis, call it $\ell_\theta$.
\\[-.25cm]

\un{Step 5:} \ 
By Jordan's lemma, the integral around the curved part is small when 
$s = x^{-\alpha} > 0$ 
is small and on 
$\ell_\theta$ 
the integrand is bounded by an absolutely integrable function, 
thus the result is continuous as a function of $s$ all the way to 0 (dominated convergence).  
Therefore
\[
\lim\limits_{x \ra \infty} \ 
x^{\alpha + 1} \hsx \Phi_\alpha (x) 
\ = \ 
\Img \ 
\int\limits_0^{\infty, \theta}
\exp (\sqrt{-1} \hsx u^{1/\alpha})
\ \td u,
\]
the symbol
$
\ds
\int\limits_0^{\infty, \theta}
\ldots
$
being an abbreviation for the integral along $\ell_\theta$.
\\[.5cm]

\un{Step 6:} \ 
Now change the variable and let 
$\ds u = v \exp\Big(\frac{\sqrt{-1} \hsx \pi \hsy \alpha}{2}\Big)$: 
\allowdisplaybreaks
\begin{align*}
\Img \ 
\int\limits_0^\infty \ 
&
\exp \Big(\sqrt{-1} \hsx v^{1/a} \hsx \exp \Big(\frac{\sqrt{-1} \hsx \pi}{2}\Big)\Big)
\ \cdot \ 
\exp \Big(\frac{\sqrt{-1} \hsx \pi \hsy a}{2} \Big)
\ \td v\ 
\\[15pt]
&\hspace{1cm}=\ 
\Img \ 
\Big(
\exp \Big(\frac{\sqrt{-1} \hsx \pi \hsy a}{2}\Big)
\ 
\int\limits_0^\infty \ 
\exp \big(-v^{1/a} 
\big)
\ \td v
\Big)
\\[15pt]
&\hspace{1cm}=\ 
\sin \Big(\frac{\pi \hsy a}{2}\Big) \ 
\int\limits_0^\infty \ 
\exp \big(- v^{1/a}\big)
\ \td v.
\end{align*}

[Note: \ Strictly speaking, this is a rotation of contours, not a change of variable.]
\\

\un{Step 7:} \ 
In 
\[
\int\limits_0^\infty \ 
\exp \big(-v^{1/a}\big)
\ \td v,
\]
let
$w = v^{1/a}$, so
\allowdisplaybreaks
\begin{align*}
\td w \ 
&=\ 
\frac{1}{a} 
\hsx 
v^\frac{1}{a}
\hsx
v^{-1}
\hsx
\td v
\\[11pt]
&=\ 
\frac{1}{a} 
\hsx 
w \cdot w^{-a}
\hsx
\td v
\\[11pt]
&=\ 
\frac{1}{a} 
\hsx 
w^{1 - a} 
\hsx
\td v
\end{align*}
\qquad 
$\implies$
\allowdisplaybreaks
\begin{align*}
\int\limits_0^\infty \ 
\exp \big(-v^{1/a}\big)
\ \td v\ 
&=\  
a \ 
\int\limits_0^\infty \ 
\exp (- w) w^{a - 1} 
\ \td w
\\[15pt]
&=\ 
a \Gamma (a) 
\\[15pt]
&=\ 
\Gamma (a + 1). 
\end{align*}
Returning to 35.10, the assumption on $\alpha$ implies that 
$
\ds
\sin \Big(\frac{\pi \hsy \alpha}{2}\Big)
\neq 0
$.
Consequently, $\Phi_\alpha$ cannot have an infinite number of real zeros.  
But $\Phi_\alpha$ does have an infinite number of zeros (cf. \S7), 
from which it follows that 
$\Phi_\alpha$ has an infinite number of nonreal zeros.
\\[-.5cm]

There remains the claim that the number (finite) of real zeros of 
$\Phi_\alpha$
is 
$
\ds
\geq 
2 \hsx \Big[\frac{\alpha}{2}\Big]
$
if $\alpha > 2$.  
To this end, choose $m \geq 1$: 
\[
2 m 
\ < \ 
\alpha
\ < \ 
2 m + 2.
\]
Write
\[
\frac{2}{\pi} \ 
\int\limits_0^\infty \ 
\Phi_\alpha (x) 
\hsy
\cos x \hsy t 
\ = \ 
e^{- t^\alpha},
\]
differentiate $2 m$ times with respect to $t$, and then put $t = 0$: 
\\[-.25cm]

\qquad 
$\implies$
\allowdisplaybreaks
\begin{align*}
\int\limits_0^\infty \ 
\Phi_\alpha (x) 
\hsy
x^2 
\ \td x \ 
&=\ 
0
\\[11pt]
&\hspace{-1cm}
\vdots
\\[11pt]
\int\limits_0^\infty \ 
\Phi_\alpha (x) 
\hsy
x^{2 m} 
\ \td x \
&=\ 
0.
\end{align*}
Accordingly, 
\[
\int\limits_0^\infty \ 
\Phi_\alpha (x) 
\hsy
x^2
\hsy P(x^2)
\ \td x
\ = \ 0,
\]
where $P$ is any polynomial of degree $\leq m - 1$.
\\[-.5cm]

For sake of argument, suppose now that 
$\Phi_\alpha (x)$ 
changes sign at most $k \leq m - 1$ times $(x > 0),$ e.g., at
\[
0 
\ < \ 
x_1
\ < \ 
x_2
\ < \ 
\cdots
\ < \ 
x_k.
\]
Introduce
\[
P (x^2) 
\ = \ 
(x_1^2 - x^2) 
\hsx 
(x_2^2 - x^2) 
\hsx \cdots \hsx
(x_k^2 - x^2) .
\]
Then
\[
\Phi_\alpha (x) 
\hsy
x^2
\hsy P(x)
\]
is never negative 
($\Phi_\alpha (0)$ is positive) while
\[
\int\limits_0^\infty \ 
\Phi_\alpha (x) 
\hsy
x^2
\hsy P(x^2)
\ \td x
\ = \ 0,
\]
a contradiction.
\\[-.5cm]

So in conclusion, 
$\Phi_\alpha (x)$ 
changes sign at least 
$
\ds
m = \Big[\frac{\alpha}{2}\Big]
$ 
times $(x > 0)$, 
thus being even, the number of real zeros of 
$\Phi_\alpha$ is 
$
\ds
\geq 
2 \hsx \Big[\frac{\alpha}{2}\Big]
$
if 
$\alpha > 2$.
\\[-.25cm]
\end{x}

\qquad
{\small\bf \un{N.B.}} \ 
This analysis breaks down if 
$1 < \alpha < 2$.  
However, in this case it can be shown that 
$\Phi_\alpha$  
has no real zeros.\footnote[2]{\vspace{.11 cm}
A. Wintner,  
\textit{American J. Math.} \textbf{58} (1936), pp. 64-66.}
\\[-.5cm]

[Note: \ 
A crucial preliminary to the proof is the fact that
\[
\scalebox{1.25}{$\ds e^{- \abs{t}^\alpha}$}
\]
is the characteristic function of an absolutely continuous distribution function 
(which is definitely not an ``elementary'' function).]
\\[-.25cm]

\begin{x}{\small\bf REMARK} \ 
Take $\phi \in \Lp^1 (0, \infty)$ real valued and twice differentiable $-$then 
under appropriate decay conditions on 
$\phi$, 
$\phi^\prime$, 
$\phi^{\prime\prime}$, 
the assumption that 
$\phi^\prime (0) \neq 0$ 
implies that 
\[
C_\infty (z) 
\ = \ 
\int\limits_0^\infty \ 
\phi (t) \hsy \cos z t 
\ \td t
\]
has an infinite number of nonreal zeros and a finite number of real zeros (if any at all).
\\[-.5cm]

[Supposing that $C_\infty (z)$ is of order $< 2$, consider the formula
\[
x^2 \hsx C_\infty (x) 
\ = \ 
- \phi^\prime (0) 
\hsx + \hsx 
\int\limits_0^\infty \ 
 \phi^{\prime\prime} (t)
\hsy \cos x \hsy t 
\ \td t
\]
that arises upon a double integration by parts.]
\\[-.5cm]

[Note: \ 
Since
\[
\frac{\td}{\td t} \exp \big(- t^\alpha\big)
\ = \ 
\exp (- t^\alpha) \hsx (- \alpha \hsy t^{\alpha - 1})
\]
vanishes at $t = 0$, this fact cannot be used to circumvent the analysis in 35.10.]
\\[-.25cm]
\end{x}

\begin{x}{\small\bf EXAMPLE} \ 
The zeros of the function 
\[
\int\limits_{-\infty}^\infty \ 
\exp (- t^{4 n} + t^{2 n} + t^2) 
\hsx 
e^{\sqrt{-1} \hsx z \hsy t} 
\ \td t
\qquad (n = 1, 2, \ldots)
\]
are real.
\\[-.25cm]
\end{x}

\begin{x}{\small\bf DEFINITION} \ 
Let $\phi \in \Lp^1  (-\infty, \infty)$ subject to 
\[
\phi (-t) 
\ = \ 
\ovs{\phi (t)}.
\]
Then $\phi$ is said to be of \un{regular growth} if 
\[
\phi (t) 
\ = \ 
\tO \big(e^{-\abs{t}^b}\big)
\qquad (\abs{t} \ra \infty)
\]
for some constant $b > 2$.
\\[-.25cm]
\end{x}

\begin{x}{\small\bf LEMMA} \ 
Suppose that $\phi$ is of regular growth $-$then $f_\infty$ is a real entire function of order 
\[
\leq \ 
\frac{b}{b - 1}
\ < \ 
2.
\]

PROOF \ 
The computation
\allowdisplaybreaks
\begin{align*}
\ov{f_\infty  (x)} \ 
&=\ 
\int\limits_{-\infty}^\infty \ 
\ovs{\phi (t)} 
\hsx 
e^{-\sqrt{-1} \hsx x \hsy t} 
\ \td t
\\[15pt]
&=\ 
\int\limits_{-\infty}^\infty \
\phi (-t)
\hsx 
e^{-\sqrt{-1} \hsx x \hsy t} 
\ \td t
\\[15pt]
&=\ 
\int\limits_{-\infty}^\infty \
\phi (t)
\hsx 
e^{\sqrt{-1} \hsx x \hsy t} 
\ \td t
\\[15pt]
&=\ 
f_\infty  (x) 
\end{align*}
shows that $f_\infty$ is real.  
Define now $\beta > 0$ by writing $b = 2 + \beta$, hence
\[
\abs{\phi (t)}
\ \leq \ 
M 
\hsy e^{-\abs{t}^{2 + \beta}}
\qquad (M > 0)
\]
\qquad 
$\implies$
\allowdisplaybreaks
\begin{align*}
\abs{f_\infty  (z)} \ 
&\leq \ 
2 \hsy M \ 
\int\limits_0^\infty \
\hsy e^{-\abs{t}^{2 + \beta}}
\hsx 
e^{\abs{z} \hsy t}
\ \td t
\\[15pt]
&=\ 
2 \hsy M \ 
\int\limits_0^\infty \
\exp 
\big(
\abs{z}  t - \abs{\hsy t \hsy }^{2 + \beta}
\big)
\ \td t.
\end{align*}
But
\[
\abs{z}  t - \abs{\hsy t \hsy }^{2 + \beta}
\ < \ 
\abs{z}  t
\]
if 
\[
0 \ < \  t \ < \  2 \abs{z}^{\frac{1}{1 + \beta}}
\]
and
\allowdisplaybreaks
\begin{align*}
\abs{z}  t - \abs{\hsy t \hsy }^{2 + \beta}
&< \ 
\Big(\frac{t}{2}\Big)^{1 + \beta}
t - t^{2 + \beta}
\\[11pt]
&< \ 
- \frac{1}{2} \hsx \hsy t^{2 + \beta}
\end{align*}
if 
\[
\abs{t} 
\ > \ 
2 \hsy \abs{z}^{\frac{1}{1 + \beta}}.
\]
Therefore
\allowdisplaybreaks
\begin{align*}
\abs{f_\infty  (z)} \ 
&\leq\ 
2 \hsy M \ 
\bigg[
\int\limits_0^{2 \hsy \abs{z}^{\frac{1}{1 + \beta}}} \ 
\ + \ 
\int\limits_{2 \abs{z}^{\frac{1}{1 + \beta} }}^\infty \ 
\bigg]
\ \hsx
\exp 
\big(
\abs{z} \hsy t - \abs{t}^{2 + \beta}
\big)
\ \td t
\\[15pt]
&\leq\ 
2 \hsy M \ 
\Big[
\abs{z}^{-1} \hsx 
\exp
\Big(\ds 2 \hsy \abs{z}^{\frac{2 + \beta}{1 + \beta}}\Big)
\Big]
\ + \ 
\int\limits_0^\infty \ 
\exp \Big(- \frac{1}{2} \hsy t^{2 + \beta}\Big)
\ \td t.
\end{align*}
And so the integral defining $f_\infty (z)$ is an entire function of order
\[
\leq \ 
\frac{2 + \beta}{1 + \beta} 
\ = \ 
\frac{b}{b - 1}
\ < \ 
2.
\]
\\[-1.5cm]
\end{x}

\qquad
{\small\bf \un{N.B.}} \ 
\[
\gen f_\infty 
\ \leq \ 
\rho ( f_\infty ) 
\ < \ 
2
\qquad (\tcf. \  6.2)
\]
\qquad \qquad
$\implies$
\[
\gen f_\infty 
\ = \ 
0
\quad \text{or} \quad 
\gen f_\infty 
\ = \ 
1.
\]
\\[-1.5cm]

\begin{x}{\small\bf RAPPEL} \ 
Suppose that the real polynomial 
\[
P (z) 
\ = \ 
a_0 + a_1 z + \cdots + a_n z^n
\]
has real zeros only $-$then $\forall \ f \in \sL - \sP$, the function 
\[
P
\Big(
\frac{\td}{\td z}
\Big) 
\hsx 
f (z) 
\ \equiv \ 
a_0 f (z) + a_1 f^\prime (z) + \cdots + a_n f^{(n)} (z)
\]
is in $\sL - \sP$ (easy extension of 12.10).
\\[-.5cm]
\end{x}

\begin{x}{\small\bf PROPAGATION PRINCIPLE} \ 
If $\phi$ is of regular growth and  if 
\[
f_\infty (z) 
\ = \ 
\int\limits_{-\infty}^\infty \ 
\phi (t)
\hsx 
e^{\sqrt{-1} \hsx z \hsy t} 
\ \td t
\]
has real zeros only, then $\forall \ f \in \sL - \sP$, the function  
\[
\int\limits_{-\infty}^\infty \ 
\phi (t)
\hsx 
f(\sqrt{-1} \hsx t)
\hsx 
e^{\sqrt{-1} \hsx z \hsy t} 
\ \td t
\]
has real zeros only.
\\[-.5cm]

PROOF \ 
Per \S12, write
\[
f (z) 
\ = \ 
\sum\limits_{n = 0}^\infty \ 
\frac{\gamma_n}{n !} \hsx z^n.
\]
Then on compact subsets of $\Cx$, 
\[
P_n (z) 
\ \equiv \ 
\tJ_n 
\Big(
f ; \frac{z}{n}
\Big)
\ra f (z)
\]
uniformly (cf. 12.9).  
Moreover, $\exists \ K > 0$: $\forall \ n$, 
\[
\Big|
\tJ_n 
\Big(
f ; \frac{z}{n}
\Big)
\Big|
\ < \ 
\exp \big( K ( \abs{z}^2 + 1)\big).
\]
The preliminaries in place, by hypothesis $f_\infty \in \sL - \sP$, thus
\[
P_n
\Big(
\frac{\td}{\td z}
\Big) 
\hsx 
f_\infty 
\in \sL - \sP 
\qquad (\tcf. \ 35.16).
\]
But
\\[-1cm]
\allowdisplaybreaks
\begin{align*}
P_n
\Big(
\frac{\td}{\td z}
\Big) 
\hsx 
f_\infty (z) \
&= \ 
\int\limits_{-\infty}^\infty \  
\phi (t)
\hsx 
P_n(\sqrt{-1} \hsx \hsy t)
\hsx
e^{\sqrt{-1} \hsx z \hsy t} 
\ \td t
\\[11pt]
&
\ra \ 
\int\limits_{-\infty}^\infty \ 
\phi (t)
\hsx 
f (\sqrt{-1} \hsx \hsy t)
\hsx
e^{\sqrt{-1} \hsx z \hsy t} 
\ \td t
\qquad (n \ra \infty).
\end{align*}
\\[-1.25cm]
\end{x}

\begin{x}{\small\bf EXAMPLE} \ 
Take $f (z) = (z + \alpha)^n$ $(n = 1, 2, \ldots)$ ($\alpha$ real) $-$then
\[
f(\sqrt{-1} \hsx \hsy t)
\ = \ 
(\sqrt{-1} \hsx \hsy t + \alpha)^n.
\]
Therefore the zeros of the function
\[
\int\limits_{-\infty}^\infty \ 
\phi (t) 
\hsx 
(\sqrt{-1} \hsx \hsy t + \alpha)^n 
\hsx
e^{\sqrt{-1} \hsx z \hsy t} 
\ \td t
\]
\\[-.75cm]
are real if $f_\infty \in \sL - \sP$.
\\[-.5cm]
\end{x}

\begin{x}{\small\bf EXAMPLE} \ 
Take $f (z) = e^{ b z}$ ($b$ real) $-$then 
\[
f(\sqrt{-1} \hsx \hsy t)
\ = \ 
e^{b \sqrt{-1} \hsx t} 
\ = \ 
\cos bt + \sqrt{-1} \hsx \sin bt.
\]
Therefore the zeros of the function 
\[
\int\limits_{-\infty}^\infty \ 
\phi (t) 
\hsx 
(\cos b t +\sqrt{-1} \hsx \sin b t) 
\hsx 
e^{\sqrt{-1} \hsx z \hsy t} 
\ \td t
\]
are real if $f_\infty \in \sL - \sP$.
\\[-.25cm]
\end{x}

\begin{x}{\small\bf EXAMPLE} \ 
Take $f (z) =e^{a \hsy z^2}$ ($a$ real and $< 0$) $-$then
\[
f (\sqrt{-1} \hsx \hsy t) 
\ = \ 
e^{a (\sqrt{-1} \hsx \hsy t)^2}
\ = \ 
e^{-a t^2}
\ = \ 
e^{\lambda t^2}
\qquad (\lambda = -a).
\]
Therefore the zeros of the function 
\[
\int\limits_{-\infty}^\infty \ 
\phi (t) 
\hsx 
e^{\lambda t^2}
\hsx
e^{\sqrt{-1} \hsx z \hsy t}
\ \td t
\qquad (\lambda > 0)
\]
are real if $f_\infty \in \sL - \sP$.
\\[-.25cm]
\end{x}

\begin{x}{\small\bf RAPPEL} \ 
Suppose that $f$ is a real entire function of genus 0 or 1 and write
\[
\abs{f (x + \sqrt{-1} \hsx y)}^2
\ = \ 
\sum\limits_{n = 0}^\infty \ 
\Lambda_n (f) (x) 
\hsy 
y^{2 n} 
\qquad (\tcf. \ 13.8)
\]
or still, 
\[
\abs{f (x + \sqrt{-1} \hsx y)}^2
\ = \ 
\sum\limits_{n = 0}^\infty \ 
L_n (f) (x) 
\hsy 
y^{2 n}
\qquad (\tcf. \ 13.9).
\]
Then $f \in \sL - \sP$ iff $\forall \ n \geq 0$ and $\forall \ x \in \R$, 

\[
L_n (f) (x) 
\ \geq \ 
0
\qquad (\tcf. \ 13.7).
\]
\\[-1.25cm]
\end{x}

\begin{x}{\small\bf APPLICATION} \ 
$f_\infty \in \sL - \sP$ iff $\forall \ n \geq 0$ and $\forall \ x \in \R$, 
\[
\int\limits_{-\infty}^\infty \ 
\int\limits_{-\infty}^\infty \ 
\phi (s) 
\hsx 
\phi (t) 
\hsx 
e^{\sqrt{-1} \hsx (s + t) \hsy x}
\hsx 
(s - t)^{2 n}
\td s 
\hsy 
\td t
\ \geq \ 
0.
\]

[In fact, 
\allowdisplaybreaks
\begin{align*}
\abs{f_\infty (x + \sqrt{-1} \hsx y)}^2 \ 
&= \ 
f_\infty (x + \sqrt{-1} \hsx y)
\hsx
f_\infty (x - \sqrt{-1} \hsx y)
\\[15pt]
&=\ 
\sum\limits_{n = 0}^\infty \ 
\frac{y^{2 n}}{n !} \ 
\int\limits_{-\infty}^\infty \ 
\int\limits_{-\infty}^\infty \ 
\phi (s) 
\hsx 
\phi (t) 
\hsx 
e^{\sqrt{-1} \hsx (s + t) \hsy x}
\hsx 
(s - t)^{2 n}
\ \td s 
\hsy 
\td t.]
\end{align*}
\\[-1.5cm]
\end{x}


\begin{x}{\small\bf EXAMPLE} \ 
Take 
\[
\phi (t) 
\ = \ 
\exp (- t^{2 \hsy k}) 
\qquad (k \geq 2) 
\quad 
(\tcf. \ 35.7).
\]
Then is it obvious that $\forall \ n \geq 0$ and $\forall \ x \in \R$, the expression 
\[
\int\limits_{-\infty}^\infty \ 
\int\limits_{-\infty}^\infty \ 
\phi (s) 
\hsx 
\phi (t) 
\hsx 
e^{\sqrt{-1} \hsx (s + t) \hsy x}
\hsx 
(s - t)^{2 n}
\td s 
\hsy 
\td t
\]
is nonnegative?
\\[-.25cm]
\end{x}

\begin{x}{\small\bf RAPPEL} \ 
Suppose that $f$ is a real entire function of genus 0 or 1 $-$then 
$f \in \sL - \sP$  iff 
\[
\frac{\partial^2}{\partial y^2}
\hsx 
\abs{f (x + \sqrt{-1} \hsx y)}^2
\ \geq \ 
0.
\]
 
[Examine the proof of 13.12.]
\\[-.25cm]
\end{x}

\begin{x}{\small\bf APPLICATION} \ 
$f_\infty \in \sL - \sP$  iff $\forall \ x, y \in \R$, 
\[
\int\limits_{-\infty}^\infty \ 
\int\limits_{-\infty}^\infty \ 
\phi (s) 
\hsx 
\phi (t) 
\hsx 
e^{\sqrt{-1} \hsx (s + t) \hsy x}
\hsx 
e^{(s - t) \hsy y}
\hsy (s - t)^2 
\ 
\td s 
\hsx 
\td t 
\ \geq \ 
0.
\]

[Differentiate
\[
\abs{f_\infty (x + \sqrt{-1} \hsx y)}^2
\ = \ 
f_\infty (x + \sqrt{-1} \hsx y)
\hsx
f_\infty (x - \sqrt{-1} \hsx y)
\]
twice with respect to $y$.]
\\[-.25cm]
\end{x}

One can employ 35.24 to ascertain that the zeros of certain real entire functions are real. 
\\[-.25cm]

\begin{x}{\small\bf EXAMPLE} \ 
We have
\[
\begin{cases}
\ 
\abs{\sin z}^2 
\ = \ 
\sin^2 x + \sinh^2 y
\\[8pt]
\ 
\abs{\cos z}^2 
\ = \ 
\cos^2 x + \sinh^2 y
\end{cases}
.
\]
And
\[
\begin{cases}
\ \ds
\frac{\partial^2}{\partial y^2}
\hsx 
\abs{\sin (x + \sqrt{-1} \hsx y)}^2 
\ = \ 
2 (\cosh^2 y + \sinh^2 y) 
\ \geq \ 2 
\ > \ 
0
\\[18pt]
\ \ds
\frac{\partial^2}{\partial y^2}
\hsx 
\abs{\cos (x + \sqrt{-1} \hsx y)}^2 
\ = \ 
2 (\cosh^2 y + \sinh^2 y) 
\ \geq \ 2 
\ > \ 
0
\end{cases}
.
\]
Therefore the zeros of $\sin z$ and $\cos z$ are real (\ldots).
\\[-.5cm]

[Note: \ 
It is a corollary that the zeros of 
\[
\begin{cases}
\ \ds
\tJ_{\frac{1}{2}} (z) 
\ = \ 
\Big(
\frac{2}{\pi z}
\Big)^{\frac{1}{2}} 
\hsx 
\sin z
\\[18pt]
\ \ds
\tJ_{-\frac{1}{2}} (z) 
\ = \ 
\Big(
\frac{2}{\pi \hsy z}
\Big)^{\frac{1}{2}} 
\hsx 
\cos z
\end{cases}
\]
are real.]
\\[-.5cm]
\end{x}

\begin{x}{\small\bf EXAMPLE} \ 
Recall from 12.33 that the zeros of the Bessel function $\tJ_\nu (z)$ $(\nu > -1)$ are real.  
This important point can also be established via 35.24.  
Thus put
\[
\tJ_\nu (z) 
\ = \ 
z^{-\nu} 
\hsx
\tJ_\nu (z).
\]
Then it can be shown that
\[
\frac{\partial^2}{\partial y^2}
\hsx 
\abs{\tJ_\nu (x + \sqrt{-1} \hsx y)}^2 
\ \geq \ 
4 
\hsx 
(\nu + 1) 
\hsx \abs{\tJ_{\nu+1} (x)}^2,
\]
from which the contention.
\\[-.75cm]
\end{x}

In terms of the modified Bessel functions, let
\[
K_z (\alpha)
\ = \ 
\frac{\pi}{2} \ 
\frac{I_{-z}(\alpha) - I_z (\alpha)}{\sin \pi z} 
\qquad (\alpha > 0).
\]
Then
\[
K_z (\alpha)
\ = \ 
\int\limits_0^\infty \ 
e^{- \alpha \hsy \cosh t}
\hsx 
\cosh z t 
\ \td t
\]
or still, 
\allowdisplaybreaks
\begin{align*}
K_{\sqrt{-1} \hsx z} (\alpha) \ 
&=\ 
\int\limits_0^\infty \ 
e^{- \alpha \hsy \cosh  \hsy t}
\hsx 
\cosh \sqrt{-1} \hsx z \hsy t
\ \td t
\\[15pt]
&=\ 
\int\limits_0^\infty \ 
e^{- \alpha \hsy \cosh  t}
\hsx 
\cos z t 
\ \td t.
\end{align*}
\\[-.75cm]

\begin{x}{\small\bf EXAMPLE} \ 
Take 
$\phi (t) = e^{- \alpha \hsy \cosh t}$ 
$-$then $\phi$ is of regular growth and the claim is that all the zeros of
\[
C_\infty (z)
\ = \ 
\int\limits_0^\infty \ 
e^{- \alpha \hsy \cosh t}
\hsx 
\cos z t 
\ \td t
\]
are real.
\\[-.5cm]

[A ``special function'' manipulation leads to the relation
\allowdisplaybreaks
\begin{align*}
&\hspace{-.75cm}
\abs{K_{\sqrt{-1} \hsx z} (\alpha) }^2 
\\[15pt]
&= \ 
\abs{K_{\sqrt{-1} \hsx x} (\alpha) }^2 
\hsx + \hsx 
y^2 \ 
\int\limits_0^1 \ 
t^{y-1} 
\hsx 
2^{F_1}
\bigg[
\overset{\ds y+1, y+1}{\underset{\ds 2}{\ }}
; 
1 - t
\bigg]
\ 
\big(
K_{\sqrt{-1} \hsx x} \big(\frac{\alpha}{\sqrt{t}}\big)
\big)^2 
\ \td t.
\end{align*}
Therefore
\[
\frac{\partial^2}{\partial y^2}
\hsx 
\abs{K_{\sqrt{-1} \hsx z} (\alpha) }^2 
\ = \ 
\int\limits_0^1 \ 
\frac{\partial^2}{\partial y^2}
\hsx 
f_t (y) 
\hsx 
\big(
K_{\sqrt{-1} \hsx x} \big(\frac{\alpha}{\sqrt{t}}\big)
\big)^2 
\ \frac{\td t}{t},
\]
where
\[
f_t (y) 
\ = \ 
y^2 
\hsx 
t^y 
\hsx 
2^{F_1}
\bigg[
\overset{\ds y+1, y+1}{\underset{\ds 2}{\ }}
; 
1 - t
\bigg].
\]
But $f_t (y)$ is an (even) absolutely monotonic function of $y$ when $0 < t < 1$, hence
\[
\frac{\partial^2}{\partial y^2}
\hsx
f_t (y) 
\ \geq \ 
0
\qquad (0 < t < 1).]
\]
\\[-1.cm]
\end{x}

\begin{x}{\small\bf RAPPEL} \ 
If $f \in \sL - \sP$, then $\forall \ \lambda \in \R$, either $f_\lambda \in \sL - \sP$ or $f_\lambda \equiv 0$ (cf. 14.9). 
\\[-.25cm]
\end{x}

\begin{x}{\small\bf EXAMPLE} \ 
Take 
\[
f (z) 
\ = \ 
K_{\sqrt{-1} \hsx z} (\alpha) 
\qquad (\alpha > 0).
\]
Then $\forall \ \lambda \in \R$, the real entire function 
\[
K_{\sqrt{-1} \hsx (z + \sqrt{-1} \hsx \lambda)} (\alpha) 
\hsx + \hsx 
K_{\sqrt{-1} \hsx (z - \sqrt{-1} \hsx \lambda)} (\alpha) 
\ = \ 
2 \ 
\int\limits_0^\infty \ 
e^{- \alpha \hsy \cosh t}
\hsx 
\cosh (\lambda t) 
\hsx 
\cos z t 
\ \td t
\]
has real zeros only.
\\[-.5cm]

[Note: \ 
Since
\[\cosh (\lambda t) 
\ = \ 
\cos (\sqrt{-1} \hsx \lambda \hsy t),
\]
one could also quote 35.17.]
\\[-.25cm]
\end{x}


\chapter{
$\boldsymbol{\S}$\textbf{36}.\quad  LOCATION, LOCATION, LOCATION}
\setlength\parindent{2em}
\setcounter{theoremn}{0}
\renewcommand{\thepage}{\S36-\arabic{page}}

\qquad 
Let $f \not\equiv 0$ be a real entire function $-$then for any real number $\lambda$, 
\[
f_\lambda (z) 
\ = \ 
f (z + \sqrt{-1} \hsx \lambda) 
\hsx + \hsx 
f (z - \sqrt{-1} \hsx \lambda) 
\qquad (\tcf. \ 14.1). 
\]
\\[-1.25cm]

\begin{x}{\small\bf NOTATION} \ 
Given $A \geq 0$ $(A < \infty)$, put
\[
A_\lambda 
\ = \ 
\big(
\max (A^2 - \lambda^2, 0)
\big)^{1/2}.
\]
\\[-1.25cm]
\end{x}

\begin{x}{\small\bf RAPPEL} \ 
Let $f \in A - \sL - \sP$ and take $\lambda > 0$ $-$then
\[
f_\lambda \in A - \sL - \sP
\qquad (\tcf. \ 15.8).
\]
\\[-1.5cm]
\end{x}

\begin{x}{\small\bf THEOREM} \ 
Suppose that $\phi$ is of regular growth and 
\[
f_\infty (z) 
\ = \ 
\int\limits_{-\infty}^\infty \ 
\phi (t) 
\hsx 
e^{\sqrt{-1} \hsx z \hsy t}
\ \td t
\]
is in 
$A - \sL - \sP$ 
$-$then for $\lambda > 0$, 
\[
(f_\infty)_\lambda (z)
\ = \ 
\int\limits_{-\infty}^\infty \ 
\phi (t) 
\hsx 
\big(
e^{\lambda\hsy t} + e^{-\lambda \hsy t}
\big) 
\hsy
e^{\sqrt{-1} \hsx z \hsy t}
\ \td t
\]
is in 
$A_\lambda - \sL - \sP$.
\\[-.5cm]

[Note: \ 
Specialize to $A = 0$ and in 35.17, take
\[
f (z) 
\ = \ 
\cos \lambda z.
\]
Then
\[
f (\sqrt{-1} \hsx t) 
\ = \ 
\cos \sqrt{-1} \hsx  \lambda t
\ = \ 
\cosh \lambda t
\ = \ 
\frac{\big(
e^{\lambda  \hsy t} + e^{-\lambda  \hsy t}
\big) }{2}, 
\]
so a priori,  
\[
(f_\infty)_\lambda 
\in 
\sL - \sP.]
\]
\\[-1.5cm]
\end{x}

\begin{x}{\small\bf LEMMA} \ 
Suppose that $\phi$ is of regular growth and 
\[
f_\infty (z) 
\ = \ 
\int\limits_{-\infty}^\infty \ 
\phi (t) 
\hsx 
e^{\sqrt{-1} \hsx z \hsy t}
\ \td t
\]
is in 
$A - \sL - \sP$ 
$-$then for 
$\lambda_1 > 0, \hsx \lambda_2 > 0, \hsx \ldots, \hsx \lambda_N > 0$, 
the zeros of 
\[
\big(
\cdots 
((f_\infty)_{\lambda_1} )_{\lambda_2} 
\cdots
\big)_{\lambda_N}
\ = \ 
\int\limits_{-\infty}^\infty \ 
\phi (t) 
\ 
\prod\limits_{k = 1}^N \
\big(
e^{\lambda_k \hsy t} + e^{-\lambda_k \hsy t}
\big)
\hsy 
e^{\sqrt{-1} \hsx z \hsy t}
\ \td t
\]
are in the strip
\[
\abs{\Img z}
\ \leq \ 
\big(
\max \big(A^2 
\hsx - \hsx  
\sum\limits_{k = 1}^N \
\lambda_k^2, \hsx 0 \big)
\big)^{1/2}.
\]
\\[-1.25cm]
\end{x}

\begin{x}{\small\bf THEOREM} \ 
Suppose that $\phi$ is of regular growth and 
\[
f_\infty (z) 
\ = \ 
\int\limits_{-\infty}^\infty \ 
\phi (t) 
\hsx 
e^{\sqrt{-1} \hsx z \hsy t}
\ \td t
\]
is in 
$A - \sL - \sP$ 
\ 
$-$then the function
\[
\int\limits_{-\infty}^\infty \ 
\phi (t) 
\hsx 
e^{\frac{1}{2} \hsx \lambda^2 \hsy t^2}
\hsx 
e^{\sqrt{-1} \hsx z \hsy t}
\ \td t
\qquad (\lambda > 0)
\]
is in 
$A_\lambda - \sL - \sP$.
\\[-.5cm]

PROOF  \ 
Given a positive integer $N$, the zeros of the function 
\[
\int\limits_{-\infty}^\infty \ 
\phi (t) 
\hsx
\Big(
\cosh \frac{\lambda t}{N}
\Big)^{N^2}
\hsx 
e^{\sqrt{-1} \hsx z \hsy t}
\ \td t
\]
lie in the strip
\allowdisplaybreaks\begin{align*}
\abs{\Img z}\ 
&\leq 
(\max(A^2 
- 
\Big(
\frac{\lambda}{N}
\Big)^2
\hsy 
N^2
, \hsx 0)
\big)^{1/2}
\\[11pt]
&=\ 
(\max(A^2 - \lambda^2, \hsx 0))^{1/2}
\qquad (\tcf. \ 36.4).
\end{align*}
But
\[
\int\limits_{-\infty}^\infty \ 
\phi (t) 
\hsx
\Big(
\cosh \frac{\lambda t}{N}
\Big)^{N^2}
\hsx 
e^{\sqrt{-1} \hsx z \hsy t}
\ \td t
\ \ra \ 
\int\limits_{-\infty}^\infty \ 
\phi (t) 
\hsx 
e^{\frac{1}{2} \hsx \lambda^2 \hsy t^2}
\hsx
e^{\sqrt{-1} \hsx z \hsy t}
\ \td t
\qquad (N \ra \infty)
\]
uniformly on compact subsets of $\Cx$.
\\[-.5cm]

[Note: \ 
To supply the details for this contention, use the inequality
\[
\cosh r 
\ \leq \ 
\exp\Big(\frac{r^2}{2}\Big)
\qquad (-\infty < r < \infty)
\]
to get
\allowdisplaybreaks\begin{align*}
C (N, t) \ 
&\equiv \ 
\Big(
\cosh \frac{\lambda t}{N}
\Big)^{N^2}
\\[11pt]
&\leq \
\exp\Big(\frac{1}{2} \hsy \lambda^2 \hsy t^2\Big).
\end{align*}
We then claim that
\[
\lim\limits_{N \ra \infty} \ 
C (N, t)
\ = \ 
\exp\Big(\frac{1}{2} \hsy \lambda^2 \hsy t^2\Big)
\]
or still, 
\[
N^2 \log \cosh 
\frac{\lambda \hsy t}{N}  
\ \ra \ 
\frac{\lambda^2 \hsy t^2}{2}
\qquad (N \ra \infty)
\]
or still, 

\[
\Big(\frac{N}{\lambda \hsy t}\Big)^2 
\hsx 
\log \cosh \frac{\lambda t}{N}
\ \ra \ 
\frac{1}{2}
\qquad (N \ra \infty).
\]
But letting 
$
\ds 
\ 
s = \frac{\lambda t}{N}
$, 
\[
\lim\limits_{s \ra 0} \ 
\frac{\log \cosh s}{s^2}
\ = \ 
\frac{1}{2}
\]
by L'Hospital.  
Now fix a compact subset $S$ of $\Cx$ and let $K > 0$ be a bound for the $\abs{\Img z}$ $(z \in S)$ $-$then
\allowdisplaybreaks
\begin{align*}
\Big|
\hsx
\phi (t) \hsy\big(C (N, t) 
&
- \exp 
\Big(\frac{1}{2} \hsy \lambda^2 \hsy t^2\Big)\big)
\hsx
e^{\sqrt{-1} \hsx z \hsy t}
\hsx
\Big|
\\[11pt]
&\leq \ 
\abs{\phi (t)}
\hsx
\Big|
C (N, t) - \Big(\frac{1}{2} \hsy \lambda^2 \hsy t^2\Big)
\Big|
\hsx
e^{K \hsy \abs{t}}
\\[11pt]
&\leq \ 
M \hsy e^{- \abs{t}^b} 
\hsx
\big(
\exp \Big(\frac{1}{2} \hsy \lambda^2 \hsy t^2\Big) 
\hsx - \hsx
C (N, t)
\big)
\hsx
e^{K \hsy \abs{t}}
\\[11pt]
&\leq \ 
M \hsy e^{- \abs{t}^b} 
\hsx
\exp \Big(\frac{1}{2} \hsy \lambda^2 \hsy t^2\Big) 
\hsx
e^{K \hsy \abs{t}}
\\[11pt]
&\in \ 
\Lp^1(-\infty, \infty)
\qquad (b > 2),
\end{align*}
so dominated convergence is applicable.]
\\[-.25cm]
\end{x}

\qquad
{\small\bf \un{N.B.}} \ 
For use below, subject the data to a relabeling: \ 
$f_\infty \in A - \sL - \sP$ 
implies that the function 
\[
\int\limits_{-\infty}^\infty \ 
\phi (t) 
\hsx 
e^{\lambda \hsy t^2}
\hsx
e^{\sqrt{-1} \hsx z \hsy t}
\ \td t
\qquad (\lambda > 0)
\]
is in 
\[
A_{\sqrt{2} \hsx \lambda} - \sL - \sP, 
\]
where
\[
A_{\sqrt{2} \hsx \lambda} 
\ = \ 
(\max (A^2 - 2 \lambda, \hsx 0))^{1/2}
\qquad (\tcf. \ 35.20).
\]
\\[-1.25cm]

\begin{x}{\small\bf NOTATION} \ 
Put 
\[
f_\infty (z; \lambda)
\ = \ 
\int\limits_{-\infty}^\infty \ 
\phi (t) 
\hsx 
e^{\lambda \hsy t^2}
\hsx
e^{\sqrt{-1} \hsx z \hsy t}
\ \td t
\qquad (\lambda \in \R),
\]
thus in particular, 
\[
f_\infty (z; 0)
\ = \ 
f_\infty (z).
\]
\\[-1.25cm]
\end{x}

\begin{x}{\small\bf LEMMA} \ 
For every real number $\lambda$, 
\[
\phi (t; \lambda) 
\ \equiv \ 
\phi (t) \hsy e^{\lambda \hsy t^2}
\]
is of regular growth.
\\[-.5cm]

PROOF \ 
By definition, for some $\beta > 0$, 
\[
e^{\abs{t}^{2 + \beta}} 
\hsx
\phi (t)
\]
stays bounded as $\abs{t} \ra \infty$.  
Let 
$
\ds
\beta^\prime = \frac{\beta}{2}
$ and consider 
\[
e^{\abs{t}^{2 + \beta^\prime}} 
\hsx
e^{\lambda \hsy t^2}
\hsx
\abs{\phi (t)}
\ = \ 
e^{t^2 (\lambda +  \abs{t}^{\beta^\prime})}
\hsy
\abs{\phi (t)}
\]
which is eventually
\[
\leq \ 
e^{\abs{t}^{2 + \beta}} 
\hsx
\abs{\phi (t)}
\]
once
\[
\lambda + \abs{t}^{\beta^\prime} 
\ < \ 
\abs{t}^\beta.
\]
\\[-1.25cm]
\end{x}

\begin{x}{\small\bf APPLICATION} \ 
If $\lambda_1 < \lambda_2$ and if the zeros of 
$f_\infty (z; \lambda_1)$
 lie in the strip 
$\{z : \abs{\Img z} \leq A\}$, 
then the zeros of 
$f_\infty (z; \lambda_2)$
lie in the strip
\[
\big\{z : \abs{\Img z} \leq A_{\sqrt{2 (\lambda_2 - \lambda_1)}}\hsy\big\}.
\]

[Simply write
\allowdisplaybreaks
\begin{align*}
f_\infty (z; \lambda_2) \ 
&=\ 
\int\limits_{-\infty}^\infty \ 
\phi (t) 
\hsx 
e^{\lambda_2 \hsy t^2} 
\hsx 
e^{\sqrt{-1} \hsx z \hsy t}
\ \td t
\\[15pt]
&=\ 
\int\limits_{-\infty}^\infty \ 
\phi (t) 
\hsx 
e^{\lambda_1 \hsy t^2} 
\hsx 
e^{(\lambda_2 - \lambda_1) t^2}
\hsx
e^{\sqrt{-1} \hsx z \hsy t}
\ \td t
\\[15pt]
&=\ 
\int\limits_{-\infty}^\infty \ 
\phi (t; \lambda_1)
\hsx 
e^{(\lambda_2 - \lambda_1) t^2}
\hsx 
e^{\sqrt{-1} \hsx z \hsy t}
\ \td t
\end{align*}
and use the assumption that the zeros of 
\[
f_\infty (z; \lambda_1)
\ = \ 
\int\limits_{-\infty}^\infty \ 
\phi (t; \lambda_1)
\hsx 
e^{\sqrt{-1} \hsx z \hsy t}
\ \td t
\]
lie in the strip 
$\{z : \abs{\Img z} \leq A\}$.]
\\[-.5cm]
\end{x}

\begin{x}{\small\bf SCHOLIUM} \ 
If the zeros of 
$f_\infty (z)$ 
lie in the strip
$\{z : \abs{\Img z} \leq A\}$, 
then 
\\[-.25cm]

\noindent
the zeros of 
$f_\infty (z; \lambda)$ $(\lambda > 0)$ 
are real when 
$A^2 - 2 \lambda \leq 0$, 
i.e., provided 
$
\ds
\frac{A^2}{2} \leq \lambda.
$
\\[-.25cm]
\end{x}

\begin{x}{\small\bf SCHOLIUM} \ 
If the zeros of 
$f_\infty (z; \lambda_1)$ 
are real and if 
$\lambda_1 < \lambda_2$, 
then the zeros of 
$f_\infty (z; \lambda_2)$ 
are real. 
\\[-.25cm]
\end{x}

There is more to be said but before so doing we shall install some machinery.
\\[-.25cm]

\begin{x}{\small\bf NOTATION} \ 
Given a complex constant $\gamma$ and an entire function $f$ of order $< 2$, let
\[
e^{\gamma \hsy D^2}\hsy f (z)
\ = \ 
\sum\limits_{n = 0}^\infty \ 
\frac{\gamma^n}{n !} \hsy f^{(2 n)} (z)
\]
or, equivalently
\[
e^{\gamma \hsy D^2}\hsy f (z)
\ = \ 
\sum\limits_{n = 0}^\infty \ 
\frac{f^{(n)} (0)}{n !} \hsx e^{\gamma \hsy D^2} \hsy z^n.
\]
\\[-1.25cm]
\end{x}

\begin{x}{\small\bf EXAMPLE} \ 
Suppose that $\phi$ is of regular growth $-$then 
$f_\infty$ 
is a real entire function of order $< 2$ (cf. 35.15) and
\[
f_\infty (z; \lambda)
\ = \ 
e^{-\lambda \hsy D^2}\hsx f_\infty (z).
\]
\\[-1.25cm]
\end{x}

\begin{x}{\small\bf LEMMA} \ 
Either series defining 
$e^{\gamma \hsy D^2}\hsy f (z)$ 
converges absolutely and uniformly on compact subsets of $\Cx$, 
hence represents an entire function.
\\[-.25cm]
\end{x}

\begin{x}{\small\bf LEMMA} \ 
$\forall$ complex constant $c$, 
\[
e^{c^2 \hsy D^2 / 2}\hsy f (z)
\ = \ 
\frac{1}{\sqrt{2 \hsy \pi}} \ 
\int\limits_{-\infty}^\infty \ 
e^{-t^2 / 2} 
\hsy
f (z + c t)
\ \td t.
\]

PROOF \ 
Bearing in mind that
\[
\int\limits_{-\infty}^\infty \ 
e^{-t^2 / 2}
\hsy 
t^{2 n} 
\ \td t
\ = \ 
\sqrt{2 \hsy \pi} 
\ 
\frac{(2 n)!}{2^n \hsy n !}
\]
and
\[
\int\limits_{-\infty}^\infty \ 
e^{-t^2 / 2}
\hsy 
t^{2 n + 1} 
\ \td t
\ = \ 
0
\]
for $n = 0, 1, 2, \ldots$, we have
\allowdisplaybreaks
\begin{align*}
e^{c^2 \hsy D^2 / 2}\hsy f (z) \ 
&= \ 
\sum\limits_{n = 0}^\infty \ 
\frac{c^{2 n}}{2^n \hsy n !}
\hsy 
f^{(2n)} (z)
\\[15pt]
&= \ 
\frac{1}{\sqrt{2 \hsy \pi}} \ 
\sum\limits_{k = 0}^\infty \ 
\int\limits_{-\infty}^\infty \ 
e^{-t^2 / 2}
\hsx 
\frac{f^{(k)} (z)}{k!}
\hsy (c t)^k
\ \td t
\\[15pt]
&= \ 
\frac{1}{\sqrt{2 \hsy \pi}} \ 
\int\limits_{-\infty}^\infty \ 
e^{-t^2 / 2}
\ 
\bigg(
\sum\limits_{k = 0}^\infty \ 
\frac{f^{(k)} (z)}{k!}
\hsy (c t)^k
\bigg)
\ \td t
\\[15pt]
&= \ 
\frac{1}{\sqrt{2 \hsy \pi}} \ 
\int\limits_{-\infty}^\infty \ 
e^{-t^2 / 2} 
\hsy
f (z + c t)
\ \td t.
\end{align*}

[Note: \ 
The interchange of summation and integration is legal.]
\\[-.25cm]
\end{x}

\begin{x}{\small\bf LEMMA} \ 
The order of 
\[
e^{\gamma \hsy D^2} \hsy f (z)
\]
is $\leq 2$.
\\[-.5cm]

PROOF \ 
For $\varepsilon > 0$ and sufficiently small, 
\[
f (z) 
\ = \ 
\tO (e^{\abs{z}^{\rho + \varepsilon}})
\qquad (\rho = \rho (f)),
\]
where
$\rho + \varepsilon < 2$, 
so $\exists$ a constant $C > 0$: 
\[
\abs{f (z)} 
\ \leq \ 
C \hsy \exp(\abs{z}^{\rho + \varepsilon}).
\]
Choose $c$ such that 
$
\ds
\gamma = \frac{c^2}{2}
$ 
$-$then 
\allowdisplaybreaks
\begin{align*}
e^{\gamma \hsy D^2} \hsy f (z) \ 
&=\ 
e^{c^2 \hsy D^2 / 2} \hsy f (z)
\\[15pt]
&= \ 
\frac{1}{\sqrt{2 \hsy \pi}} \ 
\int\limits_{-\infty}^\infty \
e^{- t^2 / 2} 
\hsy 
f (z + c t) 
\ \td t
\qquad (\tcf. \ 36.14).
\end{align*}
Therefore
\allowdisplaybreaks
\begin{align*}
\Big|\hsx e^{\gamma \hsy D^2} \hsy f (z) \hsx \Big| \ 
&\leq \ 
\frac{1}{\sqrt{2 \hsy \pi}} \ 
\int\limits_{-\infty}^\infty \
e^{- t^2 / 2} 
\hsy
\abs{f (z + c t)}
\ \td t
\\[15pt]
&\leq \ 
\frac{C}{\sqrt{2 \hsy \pi}} \ 
\int\limits_{-\infty}^\infty \
e^{- t^2 / 2} 
\hsy 
\exp\big(\abs{z + c t}^{\rho + \varepsilon}\big)
\ \td t
\\[15pt]
&\leq \ 
\frac{C}{\sqrt{2 \hsy \pi}} \ 
\int\limits_{-\infty}^\infty \
e^{- t^2 / 2} 
\hsy
\exp\big((\abs{z} + \abs{c t})^{\rho + \varepsilon}\big)
\ \td t
\\[15pt]
&\leq \ 
\frac{C}{\sqrt{2 \hsy \pi}} \ 
\int\limits_{-\infty}^\infty \
e^{- t^2 / 2} 
\hsy
\exp\big(2^{\rho + \varepsilon} (\abs{z}^{\rho + \varepsilon} + \abs{c t}^{\rho + \varepsilon})\big)
\ \td t
\\[15pt]
&\leq \ 
\frac{C}{\sqrt{2 \hsy \pi}} \ 
\bigg(
\int\limits_{-\infty}^\infty \
e^{- t^2 / 2} 
\hsy
\exp\big(2^{\rho + \varepsilon} \hsy \abs{c t}^{\rho + \varepsilon})
\ \td t
\bigg)
\exp\big(2^{\rho + \varepsilon} \hsy \abs{z}^{\rho + \varepsilon})
\\[15pt]
&\leq \ 
\frac{C}{\sqrt{2 \hsy \pi}} \ 
\bigg(
\int\limits_{-\infty}^\infty \
\ldots
\bigg)
\hsy
\exp\big(4 \abs{z}^{\rho + \varepsilon}\big),
\end{align*}
from which the assertion.
\\[-.25cm]
\end{x}

\begin{x}{\small\bf LEMMA} \ 
Given complex constants $\mu$ and $\nu$, 

\[
e^{\mu \hsy D^2} \hsy e^{\nu \hsy D^2}  \hsy f (z) 
\ = \ 
e^{(\mu + \nu) D^2}  \hsy f (z)
\ = \ 
 e^{\nu \hsy D^2} \hsy e^{\mu \hsy D^2} \hsy f (z).
\]

[Note: \ 
Thanks to 36.15, it makes sense to apply 
$e^{\mu \hsy D^2}$ to 
$e^{\nu \hsy D^2} \hsy f (z)$ 
and 
$e^{\nu \hsy D^2}$ to 
$e^{\mu \hsy D^2} \hsy f (z)$.]
\\[-.25cm]
\end{x}

\begin{x}{\small\bf RAPPEL} \ 
Define polynomials 
$\widetilde{H}_n (z)$ 
by the rule

\[
\widetilde{H}_n (z)
\ = \ 
(-1)^n 
\hsy
e^{z^2 / 2}
\ 
\frac{\td^n}{\td z^n}
\ 
e^{-z^2 / 2}
\qquad (n = 0, 1, 2, \ldots).
\]
Then the zeros of the 
$\widetilde{H}_n (z)$ 
are real and simple. 
\\[-.5cm]

[Note: \ 
This is but one of several variations on the definition of 
``Hermite polynomial'' 
(cf. 8.17).]
\\[-.25cm]
\end{x}

\begin{x}{\small\bf SUBLEMMA} \ 
Given a nonzero complex constant $c$, 

\[
e^{-c^2 \hsy D^2 /2} \hsy z^n
\ = \ 
c^n \hsy \widetilde{H}_n \Big(\frac{z}{c}\Big)
\qquad (n = 0, 1, 2, \ldots).
\]
\\[-1.25cm]
\end{x}

\begin{x}{\small\bf LEMMA} \ 
Suppose that $f (z)$ has a multiple zero at the origin $-$then there is a positive constant 
$\lambda_1$ 
such that for all 
$\lambda \in \ ]0, \lambda_1[$, \ 
$e^{\lambda \hsy D^2} \hsy f (z)$ 
has a nonreal zero.
\\[-.5cm]

PROOF \ 
Write
\[
f (z) 
\ = \ 
\sum\limits_{n = k}^\infty \ 
c_n \hsy z^n,
\]
where $k \geq 2$ and $c_k \neq 0$.  
Take $c$ positive and consider
\allowdisplaybreaks\begin{align*}
e^{c^2 \hsy D^2 / 2} 
\hsy f (z) \ 
&=\ 
\sum\limits_{n = k}^\infty \ 
c_n 
\hsy 
e^{c^2 \hsy D^2 / 2} 
\hsy 
z^n 
\\[15pt]
&=\ 
\sum\limits_{n = k}^\infty \ 
c_n \big(\sqrt{-1} \hsx c\big)^n
\hsy
\widetilde{H}_n \Big(\frac{- \sqrt{-1} \hsx z}{c}\Big).
\end{align*}
Now replace $z$ by $c w$ and instead consider
\allowdisplaybreaks\begin{align*}
F_c (w) \ 
&=\ 
\big(\sqrt{-1} \hsx c\big)^{-k} 
\hsy
e^{c^2 \hsy D^2 / 2} \hsy f (c w)
\\[11pt]&=\ 
\sum\limits_{n = k}^\infty \ 
c_n \big(\sqrt{-1} \hsx c\big)^{n-k}
\hsy
\widetilde{H}_n (- \sqrt{-1} \hsx w).
\end{align*}
The point then is that 
$\widetilde{H}_k (- \sqrt{-1} \hsx w)$ 
has a nonreal zero, thus if $c > 0$ is sufficiently small, then the same holds for 
$F_c (w)$ 
(quote Rouch\'e).  
And this suffices \ldots \hsy.
\\[-.25cm]
\end{x}

\begin{x}{\small\bf THEOREM} \ 
If the zeros of 
$f_\infty (z)$ 
lie in the strip 
$\{z : \abs{\Img z} \leq A\}$, 
\\[-.5cm]

\noindent
then the zeros of 
$f_\infty (z; \lambda)$ 
$(\lambda > 0)$ 
are real when 
$A^2 - 2 \hsy \lambda \leq 0$, 
i.e., 
provided 
$
\ds
\frac{A^2}{2}
\leq \lambda
$ 
\\[-.5cm]

\noindent
(cf. 36.9), and are simple when 
$A^2 - 2 \hsy \lambda < 0$, 
provided 
$
\ds
\frac{A^2}{2}
< \lambda
$.
\\[-.5cm]

PROOF \ 
The issue is simplicity.  
So suppose that 
\[
f_\infty (z; \lambda) 
\ = \ 
 e^{-\lambda \hsy D^2} \hsy f_\infty (z)	
\qquad (\tcf. \ 36.12)
\]
has a multiple zero at $z = a$.  
Without essential loss of generality, take $a = 0$ and apply 36.19 to 
$f_\infty (z; \lambda)$ 
and secure $\varepsilon > 0$: 
\[
e^{\varepsilon \hsy D^2} 
\hsx
e^{-\lambda \hsy D^2}
\ 
f (z)
\]
has a nonreal zero, imposing simultaneously the restriction
\[
A^2 
\ < \ 
2 (\lambda - \varepsilon).
\]
But
\allowdisplaybreaks\begin{align*}
e^{\varepsilon \hsy D^2} \hsy e^{-\lambda \hsy D^2} \hsy f_\infty (z) \ 
&=\ 
e^{- (\lambda - \varepsilon) \hsy D^2} \hsy f_\infty (z)
\qquad (\tcf. \ 36.16)
\\[11pt]
&=\ 
f_\infty (z; \lambda - \varepsilon),
\end{align*}
a function with real zeros only.  
Contradiction.
\\[-.25cm]
\end{x}

\begin{x}{\small\bf REMARK} \ 
Take $A = 0$, 
thus 
$f_\infty (z)$
is in 
$\sL - \sP$, 
as is 
$f_\infty (z; \lambda)$ 
$(\lambda > 0)$ 
and its zeros are simple.
\\[-.25cm]
\end{x}

\begin{x}{\small\bf LEMMA} \ 
Let $f$ be a real entire function of order $< 2$.  
Assume: \ 
$f \in A - \sL - \sP$ 
\ 
$-$then 
\[
e^{-\lambda \hsy D^2} \hsy f (z)
\qquad (\lambda > 0)
\]
is in 
$A_{\sqrt{2 \hsy \lambda}} - \sL - \sP$
(cf. 36.5).
\\[-.25cm]

PROOF \ 
Let 
$T^{^{\text{\small$\gamma$}}}$ 
be the translation operator
\[
T^{^{\text{\small$\gamma$}}} \hsy f (z) 
\ = \ 
f (z + \gamma).
\]
Then
\allowdisplaybreaks\begin{align*}
e^{-\lambda \hsy D^2} \hsy f (z)
\ 
&=\ 
e^{(\sqrt{-1} \hsx \sqrt{2 \hsy \lambda} \hsx )^2 \hsy D^2/ 2} \hsy f (z)
\\[15pt]
&=\ 
\lim\limits_{N \ra \infty} \ 
2^{-N}
\ 
\big(
T^{\sqrt{-1} \hsx \sqrt{2 \hsy \lambda} \hsx / \sqrt{N}} 
\hsx + \hsx 
T^{-\sqrt{-1} \hsx \sqrt{2 \hsy \lambda} \hsx / \sqrt{N}}  
\big)^N
\hsy 
f (z),
\end{align*}
the convergence being uniform on compact subsets of $\Cx$.  
But $\forall \ N$, the function 
\[
\big(
T^{\sqrt{-1} \hsx \sqrt{2 \hsy \lambda} \hsx / \sqrt{N}} 
\hsx + \hsx 
T^{-\sqrt{-1} \hsx \sqrt{2 \hsy \lambda} \hsx / \sqrt{N}}  
\big)^N
\hsy 
f (z)
\]
is in 
\[
A_{\sqrt{2 \hsy \lambda}}
\ = \ 
(\max(A^2 - 2 \hsy \lambda, \hsx 0))^{1/2}
\qquad (\tcf. \ 36.2).
\]
\\[-1.25cm]
\end{x}

\qquad
{\small\bf \un{N.B.}} \ 
In general, this estimate cannot be improved as can be seen by taking 
$f (z) = z^2 + A^2$: 
\[
e^{-\lambda \hsy D^2} 
\hsy
f (z) 
\ = \ 
z^2 + A^2 - 2 \hsy \lambda.
\]
\\[-1.25cm]

\begin{x}{\small\bf LEMMA} \ 
Let $f$ be a real entire function of order $< 2$.  
Assume: \ 
$f \in A - \sL - \sP$ 
and 
$A^2 < 2 \hsy \lambda$
$-$then 
all the zeros of 
\[
e^{-\lambda \hsy D^2} 
\hsy
f (z)
\]
are real and simple.
\\[-.5cm]

[From the above, the reality is clear and the simplicity can be established as in 36.20.]
\\[-.25cm]
\end{x}

\begin{x}{\small\bf NOTATION} \ 
\\[-.25cm]

\qquad \textbullet \quad
$S - \sL - \sP$ 
denotes the subclass of 
$\sL - \sP$ 
whose zeros are simple. 
\\[-.25cm]

\qquad \textbullet \quad
$* - S - \sL - \sP$ 
denotes the subclass of 
$* -\sL - \sP$ 
consisting of all real entire functions which are the product of a real polynomial and a function in 
$S - \sL - \sP$.
\\[-.25cm]
\end{x}

\begin{x}{\small\bf LEMMA} \ 
$S - \sL - \sP$  
and
$* - S - \sL - \sP$ 
are closed under differentiation.
\\[-.25cm]
\end{x}

\begin{x}{\small\bf NOTATION} \ 
Given complex constants 
$\gamma$, 
$c$ 
and an entire function $F$ of order $< 2$, 
define 
$\Gamma_{\gamma, c} F (z)$
by the prescription 
\[
\Gamma_{\gamma, c} F (z)
\ = \ 
(z - c) \hsy F (z) 
\hsx - \hsx 
2 \hsy \gamma \hsy F^\prime (z).
\]
\\[-1.25cm]
\end{x}

\qquad
{\small\bf \un{N.B.}} \ 
The order of 
$\Gamma_{\gamma, c} F (z)$ 
is $\hsx < 2$ (cf. 2.25 and 2.31).
\\[-.25cm]

\begin{x}{\small\bf LEMMA} \ 
$\forall \ \gamma, \hsx \forall \ c$, 
\[
e^{-\gamma \hsy D^2} 
\hsy
((z - c) \hsy F (z)) 
\ = \ 
\Gamma_{\gamma, c} 
\hsy 
e^{-\gamma \hsy D^2} 
\hsy
F (z).
\]

[Note: \ 
The order of 
\[
e^{-\gamma \hsy D^2} 
\hsy
F (z)
\]
is $\hsx < 2$ (cf. 36.15).]
\\[-.25cm]
\end{x}

\begin{x}{\small\bf LEMMA} \ 
$\forall \ \gamma \neq 0, \hsx \forall \ c$, 
\[
\Gamma_{\lambda, c} F (z) 
\ = \ 
- 2 \hsy \gamma 
\hsy
\exp\Big(\frac{(z - c)^2}{4 \gamma}\Big)
\ \frac{\td}{\td z} \ 
\Big(
\exp\Big(-\frac{(z - c)^2}{4 \gamma}\Big)
\Big) 
\hsy
F (z).
\]
\\[-1.25cm]
\end{x}

\begin{x}{\small\bf APPLICATION} \ 
Given 
$\lambda > 0$ and $a$ real, the class 
$* - S - \sL - \sP$ 
is closed under the operator 
$\Gamma_{\lambda, a}$.
\\[-.5cm]

[If $f (z)$ is in 
$* - S - \sL - \sP$, 
then 
\[
\exp\Big(-\frac{(z - a)^2}{4 \lambda}\Big) 
\hsy 
f (z)
\]
is in 
$* - S - \sL - \sP$ 
($a$ being real), as is its derivative (cf. 36.25), so all but a finite number of zeros of the latter are real and simple.  
The same then holds for 
$\Gamma_{\lambda, a} \hsy f (z)$, 
itself a real entire function of order $< 2$.]
\\[-.25cm]
\end{x}

\begin{x}{\small\bf LEMMA} \ 
Suppose that $\lambda$ is positive and $c$ is nonreal.  
Let $f$ be a real entire function of order $< 2$ and assume that 
\[
e^{- \lambda \hsy D^2}
\hsy 
f (z) 
\in 
* - S - \sL - \sP.
\]
Then
\[
e^{- \lambda \hsy D^2}
\hsy 
((z - c) (z - \bar{c})
\hsy
f (z) )
\in 
* - S - \sL - \sP.
\]

PROOF \ 
Write
\allowdisplaybreaks\begin{align*}
(z - c) (z - \bar{c})
\ 
&=\ 
z^2 - (c + \bar{c}) z + c \hsy \bar{c}
\\[11pt]
&=\ 
z^2 - 2 a z + a^2 + b^2,
\end{align*}
where 
$c = a + \sqrt{-1} \hsx b$. 
With 

\[
P (z) 
\ = \ z^2 + b^2 
\qquad (b \neq 0),
\]
we thus have
\allowdisplaybreaks\begin{align*}
(T^{-a} P) (z) \ 
&= \ 
P (z - a) 
\\[11pt]
&=\ 
(z - a)^2 + b^2 
\\[11pt]
&=\ 
z^2 - 2 a z + a^2 + b^2
\\[11pt]
&=\ 
(z - c) (z - \bar{c}).
\end{align*}
But on the basis of the definitions, 
$e^{- \lambda D^2}$ 
commutes with the translation operators 
$T^{^{\text{\small$\gamma$}}},$ 
hence
\allowdisplaybreaks\begin{align*}
e^{- \lambda \hsy D^2}
\hsy 
((z - c) (z - \bar{c})
\hsy
f (z)) \  
&=\ 
e^{- \lambda \hsy D^2}
\hsy
\big((T^{\hsy -a} P)  (z) \hsx f (z)\big) 
\\[11pt]
&=\ 
e^{- \lambda \hsy D^2}
\hsy
(T^{\hsy -a} P \cdot T^{\hsy -a + a} f)
\\[11pt]
&=\ 
e^{- \lambda \hsy D^2}
\hsy
(T^{\hsy -a}  (P \cdot T^a f))
\\[11pt]
&=\ 
T^{\hsy -a} (e^{- \lambda \hsy D^2} (P \cdot T^a f)).
\end{align*}
Since 
$* - S - \sL - \sP$ 
is closed under translation by a real constant, matters therefore reduce to showing that 
\[
e^{- \lambda \hsy D^2} (P \cdot T^a f)
\in 
* - S - \sL - \sP
\]
or still, to showing that 
\[
e^{- \lambda \hsy D^2}
\hsy
((z - \sqrt{-1} \hsz \abs{b}) \hsy (z + \sqrt{-1} \hsz \abs{b}) 
\hsx
T^a f (z)
\in 
* - S - \sL - \sP
\]
or still, to showing that 
\[
\Gamma_{\lambda, -\sqrt{-1}\hsx \abs{b}}  
\hsx \circ \hsx 
\Gamma_{\lambda, -\sqrt{-1}\hsx \abs{b}}  
\ 
\big(
e^{-\lambda \hsy D^2} 
\hsy 
T^a (f) (z)
\big)
\in 
* - S - \sL - \sP 
\qquad (\tcf. \ 36.27).
\]
And for this, cf. 36.31 and 36.32 infra.
\\[-.25cm]
\end{x}

\begin{x}{\small\bf SUBLEMMA} \ 
Fix positive contants 
$\lambda$ 
and 
$\beta$ 
$-$then 
\[
\Gamma_{\lambda, \sqrt{-1}\hsx \sqrt{\beta}} 
\hsx \circ \hsx 
\Gamma_{\lambda, -\sqrt{-1}\hsx \sqrt{\beta}}  
\ = \ 
\Gamma_{\lambda, 0}^2 
\hsx + \hsx
\beta.
\]

PROOF \ 
\[
\Gamma_{\lambda, -\sqrt{-1}\hsx \sqrt{\beta}}  \hsx F (z) 
\ = \ 
(z + \sqrt{-1} \hsx \sqrt{\beta} \hsx) \hsx F (z) 
\hsx - \hsx  
2 \hsy \lambda \hsy F^\prime (z)
\]

\qquad 
$\implies$
\allowdisplaybreaks\begin{align*}
\Gamma_{\lambda, \sqrt{-1}\hsx \sqrt{\beta}} 
\hsx 
&\circ 
\hsx 
\Gamma_{\lambda, -\sqrt{-1}\hsx \sqrt{\beta}}  
\hsx 
F (z) 
\\[11pt]
&=\ 
(z - \sqrt{-1} \hsx \sqrt{\beta} \hsx) 
\ 
\big(
(z + \sqrt{-1} \hsx \sqrt{\beta}  \hsx F (z) - 2 \hsy \lambda \hsy F^\prime (z) 
\big)
\\[11pt]
&
\hspace{1.5cm}
\ - \ 
2 \hsy \lambda 
\big(
F (z)
\hsx + \hsx  
(z + \sqrt{-1} \hsx \sqrt{\beta} \hsx) \hsx F^\prime (z)
- 
2 \hsy \lambda \hsy F^{\prime\prime} (z) 
\big)
\\[11pt]
&= \ 
(z^2 + \beta) \hsy F (z) 
- 2 \hsy \lambda \hsy 
(z - \sqrt{-1} \hsx \sqrt{\beta} \hsx + z +  \sqrt{-1} \hsx \sqrt{\beta} \hsx) 
\hsx
F^\prime (z) 
\\[11pt]
&
\hspace{1.5cm}
\ - \ 
2 \hsy \lambda \hsy F (z)
\ + \ 
4 \hsy \lambda^2 \hsy F^{\prime\prime} (z) 
\\[11pt]
&= \ 
z^2 \hsy F (z) 
\ - \ 
2 \hsy \lambda \hsy \big(2 \hsy z \hsy F^\prime (z) + F (z)\big)
\ + \ 
4 \hsy \lambda^2 \hsy F^{\prime\prime} (z) 
\ + \ 
\beta \hsy F (z).
\end{align*}
Meahwhile
\allowdisplaybreaks\begin{align*}
\Gamma_{\lambda, 0}^2  F (z) \ 
&= \ 
\Gamma_{\lambda, 0} 
\ \circ \ 
\Gamma_{\lambda, 0} F (z)
\\[11pt]
&= \ 
\Gamma_{\lambda, 0} (z F (z) - 2 \hsy \lambda \hsy F^\prime (z) )
\\[11pt]
&= \ 
z ( z F (z) - 2 \lambda F^\prime (z) )
- 2 \hsy \lambda (z F^\prime (z)  + F (z) 
- 2 \hsy \lambda \hsy  F^{\prime\prime} (z) )
\\[11pt]
&= \ 
z^2 \hsy F (z) 
\ - \ 
2 \hsy \lambda \hsy \big(2 \hsy z \hsy F^\prime (z) + F (z)\big)
\ + \ 
4 \hsy \lambda^2 \hsy F^{\prime\prime} (z).
\end{align*}
\\[-1.25cm]
\end{x}

\begin{x}{\small\bf LEMMA} \ 
Fix positive constants 
$\lambda$ 
and 
$\beta$ 
$-$then 
$* - S - \sL - \sP$ 
is closed under the operator 
\[
\Gamma_{\lambda, 0}^2 
\hsx + \hsx
\beta
\qquad (\lambda > 0, \ \beta > 0).
\]

[We shall relegate the proof of this to the Appendix of this \S.]
\\[-.25cm]
\end{x}

\begin{x}{\small\bf THEOREM} \ 
Suppose that 
$\forall \ \varepsilon > 0$, 
all but a finite number of zeros of 
$f_\infty (z)$
lie in the strip 
$\abs{\Img z} \leq \varepsilon$ 
$-$then 
$\forall \ \lambda > 0$, 
the function 
\[
f_\infty (z; \lambda) 
\ = \ 
\int\limits_{-\infty}^\infty \ 
\phi (t) 
\hsy
e^{\lambda \hsy t^2} 
\hsy
e^{\sqrt{-1} \hsx z \hsy t}
\ \td t
\]
belongs to 
$* - S - \sL - \sP$.
\\[-.5cm]

PROOF \ 
Fix
$\lambda > 0$ 
and choose 
$\varepsilon > 0$ : $\varepsilon^2 < 2 \hsy \lambda$.  
By assumption, there are only a finite number of zeros  of 
$f_\infty (z)$ 
outside the strip 
$\abs{\Img z} \leq \varepsilon$, 
hence
\[
f_\infty (z)
\ = \ 
(z - c_1) \hsy (z - \bar{c}_1) \cdots (z - c_n) \hsy (z - \bar{c}_n) \hsx f (z),
\]
where
\[
\abs{\Img c_k} 
\ > \ 
\varepsilon
\qquad (k = 1, \ldots, n)
\]
and $f (z)$ is a real entire function of order $< 2$ whose zeros lie in the strip 
$\abs{\Img z} \leq \varepsilon$, 
thus the zeros of 
$e^{-\lambda \hsy D^2} f (z)$
lie in the strip
\[
(\max (\varepsilon^2 - 2 \hsy \lambda, 0))^{1/2}
\qquad (\tcf. \ 36.22).
\]
But 
$\varepsilon^2$ 
is less than 
$2 \hsy \lambda$, 
so all the zeros of 
$e^{-\lambda \hsy D^2} f (z)$
are real and simple (cf. 36.23) or still, 
\[
e^{-\lambda \hsy D^2}
\hsy
f (z) 
\in S - \sL - \sP.
\]
Therefore
\allowdisplaybreaks\begin{align*}
f_\infty (z; \lambda) \ 
&=\ 
e^{-\lambda \hsy D^2}
\hsx
f_\infty (z)
\qquad (\tcf. \ 36.12)
\\[11pt]
&=\ 
e^{-\lambda \hsy D^2}
\hsx
((z - c_1) \hsy (z - \bar{c}_1) \cdots (z - c_n) \hsy (z - \bar{c}_n) \hsx f (z))
\\[11pt]
&\in 
* - S - \sL - \sP
\end{align*}
via iteration of 36.30.
\\[-.25cm]
\end{x}

\qquad
{\small\bf \un{N.B.}} \ 
In consequence, all but a finite number of the zeros of 
$f_\infty (z; \lambda)$ 
are real and simple and in particular 
$f_\infty (z; \lambda)$ 
has at most a finite number of nonreal zeros.
\\[-.25cm]

\begin{x}{\small\bf REMARK} \ 
The result remains valid if 
$f_\infty (z) $
is replaced by an arbitrary real entire function $f$ of order $< 2$, 
the role of 
$f_\infty (z; \lambda)$
being played by 
$e^{-\lambda \hsy D^2} \hsy f (z)$.
\\[-.25cm]
\end{x}

\begin{x}{\small\bf THEOREM} \ 
Let $f$ be a real entire function of order $< 2$.  
Assume: \ 
Given any 
$\lambda_0 > 0$, $\forall \ \varepsilon> 0$, 
all but a finite number of zeros of 
$e^{-\lambda_0 \hsy D^2} \hsy f (z)$ 
lie in the strip  
$\abs{\Img z} \leq \varepsilon$ 
$-$then 
$\forall \  \lambda > 0$, 
all but a finite number of zeros of 
$e^{-\lambda \hsy D^2} \hsy f (z)$ 
are real and simple.
\\[-.5cm]

PROOF \ 
Take 
$
\ds
\lambda_0 = \frac{\lambda}{2}
$
and put
\[
f_0 (z) 
\ = \ 
e^{-\lambda_0 \hsy D^2} \hsy f (z), 
\]
a real entire function of order $< 2$ (cf. 36.15).  
Now write
\allowdisplaybreaks\begin{align*}
e^{-\lambda \hsy D^2} \hsy f (z) \ 
&=\ 
e^{-(\lambda_0 + \lambda_0)D^2} \hsy f (z)
\\[11pt]
&=\ 
e^{-\lambda_0 \hsy D^2}
\hsx 
e^{-\lambda_0 \hsy D^2} 
\hsx 
f (z) 
\qquad (\tcf. \ 36.16)
\\[11pt]
&=\ 
e^{-\lambda_0 \hsy D^2}
\hsy 
f_0 (z) 
\end{align*}
and apply 36.34.
\\[-.25cm]
\end{x}

\begin{x}{\small\bf LEMMA} \ 
Let $f$ be a real entire function of order $< 2$.  
Assume: \ 
$f$ has $2 \hsy K$ nonreal zeros $-$then 
$\forall \ \lambda > 0$, 
$e^{-\lambda \hsy D^2}  f$ 
has at most $2 \hsy K$ nonreal zeros. 
\\[-.5cm]

[Work first with $f_\lambda$ (use 16.5).]
\\[-.25cm]
\end{x}

\begin{x}{\small\bf THEOREM} \ 
Let $f$ be a real entire function of order $< 2$.  
Assume: \ 
$f$ has $2 \hsy K$ nonreal zeros and $K$ is $\leq$ the number of real zeros of $f$.  
Fix $A > 0$: $f \in A - \sL - \sP$ \ $-$then 
\[
e^{-\lambda \hsy D^2} \hsy f (z) 
\qquad (0 < 2 \hsy \lambda < A^2)
\]
is in 
$\un{A}- \sL - \sP$ for some 
$\un{A} < \big(A^2 - 2 \hsy \lambda\big)^{1/2}$.  
\\[-.5cm]

PROOF \ 
$e^{-\lambda \hsy D^2} \hsy f$ 
has at most 
$2 \hsy K$ 
nonreal zeros and they lie in the strip
\[
\big\{z : \abs{\Img z} \leq \big(A^2 - 2 \hsy \lambda\big)^{1/2} \big\}
\qquad (\tcf. \ 36.22),
\]
thus it will be enough to show that 
$e^{-\lambda \hsy D^2} \hsy f$ 
does not vanish on the line
\[
\big\{z : \abs{\Img z} = \big(A^2 - 2 \hsy \lambda\big)^{1/2} \big\}
\]
if 
$0 < 2 \hsy \lambda < A^2$.  
Write
\[
f (z) 
\ = \ 
(z - a_1) 
\hsy
\cdots
\hsy
(z - a_K) 
\hsy
g (z), 
\]
where $a_1, \ldots, a_K$ are real zeros of $f$ and $g$ (like $f$) is a real entire function of order $< 2$ 
$-$then $f$ and $g$ have the same nonreal zeros, hence
$e^{-\lambda \hsy D^2} \hsy g$ 
has at most $K$ nonreal zeros in the open upper half-plane, these being subject to the restriction that their 
imaginary parts are positive and 
$\leq \big(A^2 - 2 \hsy \lambda\big)^{1/2}$.  
Set 
$h_0 = e^{-\lambda \hsy D^2} \hsy g$ 
and define $h_1, \ldots,  h_K$ by 
\\[-.5cm]
\[
h_k 
\ = \ 
\Gamma_{\lambda, a_k} 
\hsy 
h_{k-1}
\qquad (k = 1, \ldots, K).
\]
Then 
$h_0, h_1, \ldots, h_K$ 
are real entire functions of order $< 2$.  
And (cf. 36.27)
\allowdisplaybreaks\begin{align*}
h_1 \ 
&= \ 
\Gamma_{\lambda, a_1} \hsy h_0
\\[11pt]
&= \ 
\Gamma_{\lambda, a_1} \hsy 
e^{-\lambda \hsy D^2} \hsy g
\\[11pt]
&=\ 
e^{-\lambda \hsy D^2} \hsy ((z - a_1) g),
\end{align*}
so in the end
\[
h_K 
\ = \ 
e^{-\lambda \hsy D^2} \hsy f.
\]
If now $h_K$ has a zero $z_K$ on the line
\[
\big\{ z : \abs{\Img z} = \big(A^2 - 2 \hsy \lambda\big)^{1/2} \big\}, 
\]
then there are complex numbers $z_0, \ldots, z_{K-1}$ in the open upper half-plane such that 
$h_k (z_k) = 0$ and 
\[
\abs{z_{k+1} - \Reg z_k} 
\ \leq \ 
\Img z_k 
\qquad (k = 0, 1, \ldots, K-1) 
\quad 
\text{(Jensen \ldots)}.
\]
Therefore 
$\Img z_{k+1} \leq \Img z_k$ 
and 
$\Img z_{k+1} = \Img z_k$ 
iff 
$z_{k+1} = z_k$.  
Since 
$h_0 (z_0) = 0$, 
it follows that 
$\Img z_0  \leq \big(A^2 - 2 \hsy \lambda\big)^{1/2}$ 
from which 
\allowdisplaybreaks\begin{align*}
\Img z_K \ 
& = \ 
\big(A^2 - 2 \hsy \lambda\big)^{1/2}
\\[11pt]
&\leq \
\Img z_{K-1}
\ \leq \
\cdots
\ \leq \
\Img z_0
\ \leq \ 
\big(A^2 - 2 \hsy \lambda\big)^{1/2}
\end{align*}
\qquad 
$\implies$
\[
z_0 \ = \ z_1 \ = \ \cdots \ = \ z_K
\]
and we claim that $z_0$ is a zero of $h_0$ of multiplicity $> K$.  
First
\allowdisplaybreaks\begin{align*}
0 \ 
&=\ 
h_1 (z_1) 
\\[11pt]
&=\ 
h_1 (z_0) 
\\[11pt]
&=\ 
(z_0 - a_1) \hsy h_0 (z_0) 
- 
2 \hsy \lambda \hsy h_0^\prime (z_0) 
\\[11pt]
&=\ 
- 
2 \hsy \lambda \hsy  h_0^\prime (z_0) 
\end{align*}
\qquad 
$\implies$
\[
h_0^\prime (z_0)
\ = \ 
0.
\]
Next
\allowdisplaybreaks\begin{align*}
0 \ 
&=\ 
h_2 (z_2) 
\\[11pt]
&=\ 
h_2 (z_1) 
\\[11pt]
&=\ 
(z_0 - a_2) \hsy h_1 (z_1) 
- 
2 \hsy \lambda \hsy h_1^\prime (z_1) 
\\[11pt]
&=\ 
- 
2 \hsy \lambda \hsy  h_1^\prime (z_1) 
\\[11pt]
&=\ 
- 
2 \hsy \lambda \hsy  h_1^\prime (z_0) 
\end{align*}
\qquad \qquad
$\implies$
\[
h_1^\prime (z_0)
\ = \ 
0.
\]
But
\[
h_1 (z) 
\ = \ 
(z - a_1) h_0 (z) 
- 
2 \hsy \lambda \hsy  h_0^\prime (z) 
\]
\qquad \qquad
$\implies$
\[
h_1^\prime (z) 
\ = \ 
h_0 (z) + (z - a_1) \hsy h_0^\prime (z)  
- 
2 \hsy \lambda \hsy  h_0^{\prime\prime} (z) 
\]
\qquad \qquad
$\implies$
\allowdisplaybreaks\begin{align*}
0 \ 
&=\ 
h_1^\prime (z_0) 
\\[11pt]
&=\ 
h_0 (z) + (z_0  -  a_1) \hsy h_0^\prime (z_0) 
- 
2 \hsy \lambda \hsy  h_0^{\prime\prime} (z_0) 
\\[11pt]
&=\ 
- 
2 \hsy \lambda \hsy  h_0^{\prime\prime} (z_0) 
\end{align*}
\qquad \qquad
$\implies$
\[
h_0^{\prime\prime} (z_0) 
\ = \ 0.
\]
ETC.
However the claim leads to a contradiction: \ 
$h_0 = e^{-\lambda \hsy D^2} g$ 
has at most $K$ nonreal zeros in the open upper half-plane.
\\[-.25cm]
\end{x}

\qquad
{\small\bf \un{N.B.}} \ 
The condition on $K$ is obviously fulfilled if the number of real zeros of $f$ is infinite.
\\[-.25cm]

\[
\text{APPENDIX}
\]

Here a proof of 36.32 will be sketched.  
So take an 
$f \in * - S - \sL - \sP$ 
\ 
$-$then the claim is that
\[
(\Gamma_\lambda^2 + \beta) \hsy f
\qquad (\Gamma_\lambda^2 \equiv \Gamma_{\lambda, 0}^2)
\]
remains within  
$* - S - \sL - \sP$ 
and for this, it can be assumed that $f$ has infinitely many real zeros.
\\

SETUP \ 
Write
\[
f (z) 
\ = \ 
e^{a z^2 + b z} 
\hsx 
Q (z) 
\ 
\prod\limits_{n = 1}^\infty \ 
\Big(
1 - \frac{z}{\lambda_n}
\Big)^{z /  \lambda_n}, 
\]
where $a$ is real and $\leq 0$, 
$b$ is real, 
$Q (z)$ is a real polynomial, 
the $\lambda_n$ are real and distinct with 
\[
\sum\limits_{n = 1}^\infty \ 
\frac{1}{\lambda_n^2} 
\ < \ 
\frac{1}{4 \hsy \beta}
\qquad (\tcf. \ 10.19).
\]
Choose a positive constant $B$ such that $\abs{t} \geq B$

\[
\implies 
\hspace{1.cm}
Q (t) \neq 0, 
\quad 
\frac{\td}{\td t} 
\ 
\frac{Q^\prime (t)}{Q(t)}
\ < \ 
0, 
\quad \text{and} \quad 
\Big| \hsx
\frac{b}{t} + \frac{Q^\prime (t)}{t \hsy Q(t)}
\Big|
\ < \ 
\frac{1}{4 \lambda}.
\hspace{1cm}
\]

\noindent
Assume further that the zeros of $f (z)$ that lie in $\abs{z} \geq B$ are real and simple.
\\

\qquad
{\small\bf NOTATION} \ 
For $R > 0$, put
\[
f_R (z) 
\ = \ 
e^{a z^2 + b z} 
\hsx 
Q (z) 
\ 
\prod\limits_{\abs{\lambda_n} \hsy < \hsy R}\ 
\Big(
1 - \frac{z}{\lambda_n}
\Big)^{z /  \lambda_n}.
\]
\\

\qquad
{\small\bf \un{N.B.}} \ 
\[
(\Gamma_\lambda^2 + \beta) \hsy f_R 
\ \in \ 
* - \sL - \sP
\]
and
\[
(\Gamma_\lambda^2 + \beta) \hsy f_R 
\ \ra \ 
(\Gamma_\lambda^2 + \beta) \hsy f
\qquad (R \ra \infty)
\]
uniformly on compact subsets of $\Cx$.
\\

\qquad
{\small\bf LEMMA} \ 

\[
\frac{\Gamma_\lambda \hsy f_R (z) }{f_R (z)}
\ = \ 
(1 - 4 \hsy \lambda \hsy a) z 
\hsx - \hsx 
2 \hsy \lambda \hsy b 
\hsx - \hsx 
2 \hsy \lambda \ \frac{Q^\prime (z)}{Q (z)} 
\hsx - \hsx 
2 \hsy \lambda \ 
\sum\limits_{\abs{\lambda_n} \hsy < \hsy R} \ 
\frac{z}{\lambda_n (z - \lambda_n)}.
\]
\\[-.75cm]

\qquad
{\small\bf APPLICATION} \ 
If 
$\lambda^\prime$, 
$\lambda^{\prime\prime}$ 
are two consecutive real zeros of $f_R (z)$ such that 
$\lambda^\prime < \lambda^{\prime\prime} \leq -B$ 
\ 
or 
\ 
$B \leq \lambda^\prime < \lambda^{\prime\prime}$, 
then 
\[
\frac{\Gamma_\lambda \hsy f_R (z) }{f_R (z)}
\]
has exactly one real zero between 
$\lambda^\prime$ 
and 
$\lambda^{\prime\prime}$.

[In fact, 

\[
\lim\limits_{t \downarrow \lambda^\prime} \ 
\frac{\Gamma_\lambda \hsy f_R (t) }{f_R (t)}
\ = \ 
- \infty,  
\quad 
\lim\limits_{t \uparrow \lambda^{\prime\prime}} \ 
\frac{\Gamma_\lambda \hsy f_R (t) }{f_R (t)}
\ = \ 
\infty
\]
and

\[
\frac{\Gamma_\lambda \hsy f_R (t) }{f_R (t)}
\]
is strictly increasing in the interval 
$]\lambda^\prime, \lambda^{\prime\prime}[$\hsy .]
\\[-.25cm]

\qquad
{\small\bf LEMMA} \ 
Suppose that
\[
\frac{\Gamma_\lambda \hsy f_R (r_0) }{f_R (r_0)}
\ = \ 
0
\qquad (r_0 \in \R,\ \abs{r_0} \geq B).
\]
Then the real numbers 
\[
f_R (r_0) 
\quad \text{and} \quad 
(\Gamma_\lambda^2 + \beta) \hsy f_R (r_0) 
\]
are of opposite sign.
\\[-.5cm]

PROOF \ 
Trivially, 
\[
r_0 
\ = \ 
\frac{2 \hsy \lambda \hsy f_R^\prime (r_0)}{f_R (r_0)}.
\]
Therefore
\[
\frac{r_0 }{2 \hsy \lambda}
\ = \ 
2 \hsy a \hsy r_0 
\hsx + \hsx 
b 
\hsx + \hsx 
\frac{Q^\prime (r_0)}{Q (r_0)}
\ + \ 
\sum\limits_{\abs{\lambda_n} \hsy < \hsy R}
\frac{r_0 }{\lambda_n \hsy \big(r_0 - \lambda_n\big)}
\]
\qquad 
$\implies$
\allowdisplaybreaks
\begin{align*}
\frac{1 }{2 \hsy \lambda} \ 
&=\ 
2 a 
\hsx + \hsx 
\Big(
\frac{b}{r_0} 
\hsx + \hsx 
\frac{Q^\prime (r_0)}{r_0 \hsy Q (r_0)}
\Big)
\ + \ 
\sum\limits_{\abs{\lambda_n} \hsy < \hsy R}
\frac{1}{\lambda_n \hsy \big(r_0 - \lambda_n\big)}
\\[15pt]
&\leq \ 
\Big(
\frac{b}{r_0} 
\hsx + \hsx 
\frac{Q^\prime (r_0)}{r_0 \hsy Q (r_0)}
\Big)
\ + \ 
\sum\limits_{\abs{\lambda_n} \hsy < \hsy R}
\frac{1}{\lambda_n \hsy \big(r_0 - \lambda_n\big)}
\\[15pt]
&\leq \ 
\Big|
\hsx
\frac{b}{r_0} 
\hsx + \hsx 
\frac{Q^\prime (r_0)}{r_0 \hsy Q (r_0)}
\hsx
\Big|
\ + \ 
\Big|
\sum\limits_{\abs{\lambda_n} \hsy < \hsy R}
\frac{1}{\lambda_n \hsy \big(r_0 - \lambda_n\big)}
\Big|
\\[15pt]
&< \ 
\frac{1}{4 \hsy \lambda} 
\ + \ 
\Big|
\hsx
\sum\limits_{\abs{\lambda_n} \hsy < \hsy R}
\frac{1}{\lambda_n \hsy \big(r_0 - \lambda_n\big)}
\hsx
\Big|
\end{align*}
\qquad 
$\implies$
\allowdisplaybreaks
\begin{align*}
\frac{1}{4 \hsy \lambda} \ 
&<\ 
\Big|
\hsx
\sum\limits_{\abs{\lambda_n} \hsy < \hsy R} \ 
\frac{1}{\lambda_n \hsy \big(r_0 - \lambda_n\big)}
\hsx
\Big|
\\[15pt]
&\leq \ 
\sum\limits_{\abs{\lambda_n} \hsy < \hsy R} \ 
\frac{1}{\abs{\lambda_n} \hsy \abs{r_0 - \lambda_n}}
\\[15pt]
&\leq \ 
\Big(
\sum\limits_{\abs{\lambda_n} \hsy < \hsy R} \ 
\frac{1}{\lambda_n^2}
\Big)^{1/2}
\ 
\Big(
\sum\limits_{\abs{\lambda_n} \hsy < \hsy R} \ 
\frac{1}{\big(r_0 - \lambda_n\big)^2}
\Big)^{1/2}
\\[15pt]
&< \ 
\frac{1}{2 \sqrt{\beta}} 
\ 
\Big(
\sum\limits_{\abs{\lambda_n} \hsy < \hsy R} \ 
\frac{1}{\big(r_0 - \lambda_n\big)^2}
\Big)^{1/2}
\end{align*}
\qquad 
$\implies$
\allowdisplaybreaks
\begin{align*}
\sum\limits_{\abs{\lambda_n} \hsy < \hsy R} \ 
\frac{1}{\big(r_0 - \lambda_n\big)^2} \ 
&>\ 
\Big(\frac{1}{4 \hsy \lambda}\Big)^2 
\hsx 
(2 \sqrt{\beta})^2
\\[15pt]
&=\ 
\frac{\beta}{4 \lambda^2}.
\end{align*}
Moving on, 

\allowdisplaybreaks
\begin{align*}
\frac{(\Gamma_\lambda^2 + \beta) \hsy f_R (r_0)}{f_R (r_0)} \ 
&=\ 
\beta 
\hsx - \hsx
2 \hsy \lambda 
\hsx + \hsx
4 \hsy \lambda^2
\ 
\frac{f_R^{\prime\prime} (r_0) \hsy f_R (r_0) - f_R^\prime (r_0)^2}{f_R (r_0)^2} 
\\[15pt]
&=\ 
\beta 
\hsx - \hsx
2 \hsy \lambda 
\hsx + \hsx
4 \hsy \lambda^2
\ 
\frac{\td}{\td t} 
\ 
\Big(
\frac{f_R^\prime (t)}{f_R (t)}
\Big)
\bigg|_{t = r_0}
\\[15pt]
&=\ 
\beta 
\hsx - \hsx
2 \hsy \lambda 
\hsx + \hsx
4 \hsy \lambda^2
\ 
\bigg(
2 a 
\hsx + \hsx
\frac{\td}{\td t} 
\ 
\Big(
\frac{Q^\prime (t)}{Q(t)}
\Big)
\bigg|_{t = r_0}
\hsx - \hsx
\sum\limits_{\abs{\lambda_n} \hsy < \hsy R} \ 
\frac{1}{\big(r_0 - \lambda_n\big)^2} 
\bigg)
\\[15pt]
&<\ 
\beta 
\hsx + \hsx
4 \hsy \lambda^2
\Big(
-
\sum\limits_{\abs{\lambda_n} \hsy < \hsy R} \ 
\frac{1}{\big(r_0 - \lambda_n\big)^2} 
\Big).
\end{align*}
But
\[
\sum\limits_{\abs{\lambda_n} \hsy < \hsy R} \ 
\frac{1}{\big(r_0 - \lambda_n\big)^2} \ 
\ >\ 
\frac{\beta}{4 \lambda^2},
\]
so
\[
\frac{(\Gamma_\lambda^2 + \beta) \hsy f_R (r_0)}{f_R (r_0)} 
\ < \ 
\beta - \beta 
\ = \ 0.
\]
\\[-.5cm]

\begin{spacing}{1.75}
\qquad
{\small\bf APPLICATION} \ 
If 
$\lambda^\prime $, 
$\lambda^{\prime\prime}$, 
$\lambda^{\prime\prime\prime}$ 
are three consecutive real zeros of $f_R (z)$
such that 
$\lambda^\prime 
< 
\lambda^{\prime\prime}
< 
\lambda^{\prime\prime\prime}
\leq 
-B$
\ 
or 
\ 
$
B 
\leq
\lambda^\prime 
< 
\lambda^{\prime\prime}
< 
\lambda^{\prime\prime\prime}
$ 
and if $r_1$ and $r_2$ are real zeros of 
$
\ 
\ds
\frac{\Gamma_\lambda \hsy f_R (z)}{f_R (z)}
\ 
$
such that 
$\lambda^\prime 
< 
r_1 
< 
\lambda^{\prime\prime}
< 
r_2 
< 
\lambda^{\prime\prime\prime}
$, 
then 
$(\Gamma_\lambda^2 + \beta) \hsy f_R (z)$ 
has a real zero between $r_1$ and $r_2$.  
\\[-.5cm]
\end{spacing}

[As a part of the overall setup, the zeros of 
$f_R (z)$ are real and simple.]
\\[-.25cm]

\qquad
{\small\bf NOTATION} \ 
Given an entire function $F (z)$ and a subset $S$ of $\Cx$, let
\[
N (F (z); S)
\]
denote the number (counting multiplicity) of zeros of $F (z)$ that lie in $S$.
\\[-.25cm]

\qquad
{\small\bf EXAMPLE} \ 
\[
N((\Gamma_\lambda^2 + \beta) f_R (z) ; \Cx)
\ = \ 
N(f_R (z) ; \Cx)
\hsx + \hsx 
2.
\]
\\[-1.25cm]

\qquad
{\small\bf EXAMPLE} \ 
\allowdisplaybreaks\begin{align*}
N((\Gamma_\lambda^2 
&
+ \beta) \hsy f_R (z); \hsx ]-\infty, -B] 
\hsx \cup \hsx 
[B, \infty[\hsx )
\\[11pt]
&
\ \geq \ 
N(f_R (z) ; \hsx   ]-\infty, -B] 
\hsx \cup \hsx [B, \infty[\hsx ) 
\hsx - \hsx 
4.
\end{align*}
\\[-1.25cm]

\qquad
{\small\bf LEMMA} \ 
We have 
\allowdisplaybreaks\begin{align*}
N((\Gamma_\lambda^2 
&
+ \beta) \hsy f_R (z); \Img z \neq 0)
\ \leq \ 
N(f (z) ; \Img z \neq 0)
\hsx + \hsx 
N(f (z) ; \hsx ]-B, B[\hsx) 
\hsx + \hsx 
6.
\end{align*}

PROOF \ 
Rewrite the first term as
\[
N((\Gamma_\lambda^2 + \beta) \hsy f_R (z); \Cx)
\hsx - \hsx 
N((\Gamma_\lambda^2 + \beta) \hsy f_R (z); \R)
\]
and then bound it by
\[N (f_R (z); \Cx) 
\hsx + \hsx 
2 
\hsx - \hsx 
N((\Gamma_\lambda^2 + \beta)  (f_R (z); \hsx ]-\infty, -B] \hsx \cup \hsx [B, \infty[ \hsy )
\]
or stilll, by
\[
N (f_R (z); \Cx) 
\hsx - \hsx 
N (f_R (z); \hsx ]-\infty, -B] \hsx \cup \hsx[ B, \infty[ \hsy ) 
\hsx + \hsx 
6
\]
or stilll, by
\[
N (f_R (z); \Img z \neq 0) 
\hsx + \hsx  
N (f_R (z); \hsx ]-B,\hsx B[ \hsy ) 
\hsx + \hsx 
6
\]
or stilll, by
\[
N (f (z); \Img z \neq 0) 
\hsx + \hsx 
N (f_R (z); \hsx ]-B,\hsx B[ \hsy ) 
\hsx + \hsx 
6.
\]
\\[-.25cm]
Accordingly, 
\[
(\Gamma_\lambda^2 + \beta) \hsy f \in * - \sL - \sP
\]
but there remains the possibility that it might have infinitely many multiple zeros. 
However, if this were the case, then we would have

\[
\lim\limits_{A \ra \infty} \ 
\big(
N((\Gamma_\lambda^2 + \beta) \hsy f (z); ]-A,A[ \hsy)
\hsx - \hsx
N (f (z); \hsx ]-A,A[ \hsy)
\big)
\ = \ 
\infty.
\]
And: 
\\[-.25cm]

\qquad
{\small\bf LEMMA} \ 
Take $A > B$ $-$then $\exists \ R_0 > A$ such that
\[
N((\Gamma_\lambda^2 + \beta) \hsy f (z); \hsx \abs{\Reg z} < A)
\ \leq \ 
N((\Gamma_\lambda^2 + \beta) \hsy f_{R_0} (z); \hsx \abs{\Reg z} < A).
\]
On the other hand, 
\allowdisplaybreaks\begin{align*}
N((\Gamma_\lambda^2 
&
+ \beta) \hsy f (z);
\hsx ]-A,A[ \hsy)\
\\[11pt]
&=\ 
N((\Gamma_\lambda^2 + \beta) \hsy f (z); \hsx \abs{\Reg z} < A)
\\[11pt]
&\leq\
N((\Gamma_\lambda^2 + \beta) \hsy f_{R_0} (z); \hsx \abs{\Reg z} < A)
\\[11pt]
&\leq\
N((\Gamma_\lambda^2 + \beta) \hsy f_{R_0} (z); \Cx) 
\hsx - \hsx
N((\Gamma_\lambda^2 + \beta) \hsy f_{R_0} (z); \abs{\Reg z} \geq A)
\\[11pt]
&\leq\
N((\Gamma_\lambda^2 + \beta) \hsy f_{R_0} (z); \Cx) 
\hsx - \hsx
N((\Gamma_\lambda^2 + \beta) \hsy f_{R_0} (z); \hsx
]-\infty, -A] 
\hsx \cup \hsx 
[A, \infty[\hsx) 
\\[11pt]
&\leq\
N( f_{R_0} (z) ; \Cx) 
\hsx + \hsx
2
\hsx - \hsx
N( f_{R_0} (z) ; \hsx]-R_0, -A] 
\hsy \cup \hsy 
[A, R_0 \hsx [
\hsx + \hsx
4
\\[11pt]
&=\ 
N( f_{R_0} (z) ; \Img z \neq 0) 
\hsx + \hsx
N( f_{R_0} (z) ; \hsx ]-A,A[ \hsy )
\hsx + \hsx
6
\\[11pt]
&\leq\
N (f (z); \Img z \neq 0) 
\hsx + \hsx
N (f (z); \hsx ]-A,A[ \hsy)
\hsx + \hsx
6
\end{align*}
\qquad 
$\implies$
\allowdisplaybreaks\begin{align*}
N((\Gamma_\lambda^2 
+ \beta) \hsy f (z); ]-A,A[ \hsy)
\hsx - \hsx
&
N (f (z); ]-A,A[ \hsy)
\ \ \leq\ \
N (f (z); \Img z)
\hsx + \hsx
6,
\end{align*}
from which a contradiction (send $A$ to $\infty$).
\\[-.25cm]


\chapter{
$\boldsymbol{\S}$\textbf{37}.\quad  THE $\sF_0$\ - CLASS}
\setlength\parindent{2em}
\setcounter{theoremn}{0}
\renewcommand{\thepage}{\S37-\arabic{page}}


\qquad Let $F$ be a real entire function such that
\[
\log M (r; F) 
\ = \ 
\tO (r^4)
\qquad (r \ra \infty) 
\]
and
\[
\int\limits_{-\infty}^\infty \ 
\abs{F (\sqrt{-1} \hsx t)} 
\ \td t 
\ < \ 
\infty.
\]
\\[-1cm]

[Note: \ 
Since $F$ is real, 
$\ov{F (z)} = F (\bar{z})$, 
hence if 
\hsx
$G (t) = F (\sqrt{-1} \hsx t)$, 
then 
\allowdisplaybreaks
\begin{align*}
G (-t) \ 
&=\ 
F (\sqrt{-1} \hsx (-t))
\\[11pt]
&=\ 
F ((-\sqrt{-1} \hsx t))
\\[11pt]
&=\ 
F (\ov{\sqrt{-1}} \hsx t)
\\[11pt]
&=\ 
F (\ov{\sqrt{-1} \hsx t})
\\[11pt]
&=\ 
\ov{F (\sqrt{-1} \hsx t)}
\\[11pt]
&=\ 
\ov{G (t)} \hsx.]
\end{align*}
\\[-1cm]

\begin{x}{\small\bf DEFINITION} \ 
$F \in \sF_0$ provided all its zeros are real and 
\[
\sum\limits_n \ 
\frac{1}{\lambda_n^4}
\ < \ 
\infty
\qquad (F (\lambda_n) = 0, \ \lambda_n \neq 0).
\]

[Note: \ 
The sum is finite or infinite.]
\\[-.25cm]
\end{x}

\begin{x}{\small\bf THEOREM} \ 
Suppose that $F \in \sF_0$ and 
\[
f (z) 
\ \equiv \ 
\int\limits_{-\infty}^\infty \ 
F (\sqrt{-1} \hsx t)
\hsx 
e^{\sqrt{-1} \hsx z \hsy t}
\ \td t.
\]
Then $f \in \sL - \sP$.
\\[-.5cm]

[Note: \  
While not quite obvious, the assumptions on $F$ imply that $f$ is entire (see below).  
Moreover $f$ is real: 
\allowdisplaybreaks
\begin{align*}
\ov{f (x)} \ 
&=\ 
\int\limits_{-\infty}^\infty \ 
\ov{F (\sqrt{-1} \hsx t)} 
\hsx 
e^{- \sqrt{-1} \hsx x \hsy t}
\ \td t
\\[15pt]
&=\ 
\int\limits_{-\infty}^\infty \ 
F (-\sqrt{-1} \hsx t) 
\hsx 
e^{- \sqrt{-1} \hsx x \hsy t}
\ \td t
\\[15pt]
&=\ 
\int\limits_{-\infty}^\infty \ 
F (\sqrt{-1} \hsx t)
\hsx 
e^{\sqrt{-1} \hsx x \hsy t}
\ \td t
\\[15pt]
&=\ 
f (x) \hsx .]
\end{align*}
\\[-1.25cm]
\end{x}

\begin{x}{\small\bf RAPPEL} \ 
If $f_n \in \sL - \sP$ $(n = 1, 2, \ldots)$ and if $f_n \ra f$ uniformly on compact subsets of $\Cx$, 
then 
$f \in \sL - \sP$.
\\[-.25cm]
\end{x}

The proof of 37.2 falls into two cases, according to whether the number of zeros of $F$ is finite or infinite.
\\[-.5cm]

So suppose first that $F$ has finitely many zeros $-$then there exists a real polynomial $P$ and real 
constants $\alpha, \beta, \gamma, \delta$ such that $P$ has only real zeros, 
$\alpha$ is nonnegative, $\max(\alpha, \gamma)$ is positive, and 
\[
F (z) 
\ = \ 
P (z) 
\hsy
\exp(- \alpha^2 \hsy z^4 - \beta^3 \hsy z^3 + \gamma \hsy z^2 + \delta \hsy z).
\]
Choose a positive integer $N$: 
\[
2 \hsy n \hsy \alpha + \frac{3}{2} \hsy n \hsy \beta^2 + \gamma 
\ > \ 
0
\qquad (n \geq N).
\]
Then define $F_n(z)$ $(n \geq N)$ by 
\allowdisplaybreaks
\begin{align*}
&
F_n (z) 
\ = \ 
P (z) 
\hsx 
\Big(
\Big(
1 - \frac{\alpha \hsy z^2}{n}
\Big)
\exp
\Big(
\frac{\alpha \hsy z^2}{n}
\Big)
\Big)^{2 \hsy n^2}
\\[15pt]
&
\hspace{2cm}
\hspace{1cm}
\times  \ 
\Big(
\Big(
1 - \frac{\beta \hsy z}{n}
\Big)
\exp
\Big(
\frac{\beta \hsy z}{n} 
+ 
\frac{\beta^2 \hsy z^2}{2 n^2} 
\Big)
\Big)^{3 \hsy n^3}
\hsx
e^{\gamma \hsy z^2 + \delta \hsy z}
\end{align*}
and set
\[
f_n (z) 
\ = \ 
\int\limits_{-\infty}^\infty \ 
F_n (\sqrt{-1} \hsx t)
\hsx 
e^{\sqrt{-1} \hsx z \hsy t}
\ \td t.
\]
\\[-1cm]

\begin{x}{\small\bf LEMMA} \ 
$f_n \ra f$ uniformly on compact subsets of $\Cx$.
\\[-.5cm]

PROOF \ 
In fact, 
\[
\Big(
\Big(
1 - \frac{\alpha \hsy z^2}{n}
\Big)
\exp
\Big(
\frac{\alpha \hsy z^2}{n}
\Big)
\Big)^{2 \hsy n^2}
\ra e^{- \alpha^2 \hsy z^4}
\]
and
\[
\Big(
\Big(
1 - \frac{\beta \hsy z}{n}
\Big)
\hsx
\exp
\Big(
\frac{\beta \hsy z}{n}
\hsx + \hsx 
\frac{\beta^2 \hsy z^2}{2 \hsy n^2}
\Big)
\Big)^{3 \hsy n^3}
\ra 
e^{- \beta^3 \hsy z^3}
\]
uniformly on compact subsets of $\Cx$.  
On the other hand, 
\[
\Big|
\hsx
\Big(
1 - \frac{\beta \hsy \sqrt{-1} \hsx t}{n}
\Big)
\hsx
\exp
\Big(
\frac{\beta \hsy \sqrt{-1} \hsx t}{n}
\hsx + \hsx 
\frac{\beta^2 \hsy (\sqrt{-1} \hsx t)^2}{2 \hsy n^2}
\Big)
\hsx
\Big|
\ \leq \ 
1
\qquad (t \in \R).
\]
In addition, there are positive constants $C$, $t_0$ such that 
\[
\Big(
\Big(
1 + \frac{\alpha\hsy t^2}{n}
\Big)
\hsx
\exp
\Big(
-
\frac{\alpha\hsy t^2}{n}
\Big)
\Big)^{2 \hsy n^2}
\hsx
e^{- \gamma \hsy t^2}
\ \leq \ 
e^{-C \hsy t^2} 
\qquad (n \geq N, \ \abs{t\hsy} \geq t_0).
\]
And this sets the stage for dominated convergence.
\\[-.25cm]
\end{x}

\begin{x}{\small\bf LEMMA} \ 
$\forall \ n \geq N$, 
$f_n \in \sL - \sP$.
\\[-.5cm]

PROOF \ 
We have
\allowdisplaybreaks
\begin{align*}
&
F_n (z) 
\ = \ 
P (z) 
\hsx 
\Big(
1 - \frac{\alpha \hsy z^2}{n}
\Big)^{2 \hsy n^2}
\hsx
\Big(
1 - \frac{\beta \hsy z}{n}
\Big)^{3 \hsy n^3}
\\[15pt]
&
\hspace{2.5cm}
\times  \ 
\exp \Big(\Big(2 \hsy n \hsy \alpha + \frac{3}{2} \hsx n \hsy \beta^2 + \gamma\Big) z^2
\hsx + \hsx 
(3 \hsy n^2 \hsy \beta + \delta)\hsy z\Big).
\end{align*}
But
\[
2 \hsy n \hsy \alpha 
\hsx + \hsx 
\frac{3}{2} \hsx n \hsy \beta^2
\hsx + \hsx 
\gamma
\ > \ 
0
\]
and replacing $z$ by $\sqrt{-1} \hsx t$ leads to 
\[
-
\Big(2 \hsy n \hsy \alpha 
\hsx + \hsx
\frac{3}{2} \hsx n \hsy \beta^2 
\hsx + \hsx
\gamma\Big) \hsy t^2,
\]
thus an application of 12.37 completes the proof. 
\\[-.5cm]
\end{x}

Taking into account 37.3, it then follows from 37.4 and 37.5 that 
$f \in \sL - \sP$.  
\\[-.25cm]

Suppose now that $F$ has infinitely many zeros (by hypothesis real) and write
\allowdisplaybreaks
\begin{align*}
F (z) \ 
&= \ 
M \hsy z^m \hsx
\exp(A_4 \hsy z^4 + A_3 \hsy z^3 + A_2 \hsy z^2 + A_1 \hsy z)
\\[15pt]
&
\hspace{2cm}
\times  \ 
\prod\limits_{n = 1}^\infty \ 
\Big(
1 - \frac{z}{\lambda_n}
\Big)
\hsx
\exp
\Big(
\frac{z}{\lambda_n} 
+ \frac{z^2}{2 \hsy \lambda_n^2}
+ \frac{z^3}{3 \hsy \lambda_n^3}
\Big),
\end{align*}
where $M \neq 0$ is real, $m$ is a nonnegative integer, $A_1, A_2, A_3, A_4$ are real constants, 
the $\lambda_n$ are real with 
$
\ 
\ds
\sum\limits_{n = 1}^\infty \ 
\frac{1}{\lambda_n^4} < \infty
\ 
$ 
$-$then $\forall \ t \in \R$, 
\[
\abs{F (\sqrt{-1} \hsx t)} 
\ = \ 
\abs{M} \hsy \abs{t\hsy}^m 
\hsx
e^{A_4 \hsy t^4 - A_2 \hsy t^2}
\ 
\prod\limits_{n = 1}^\infty \ 
\Big(
1 \hsx + \hsx  \frac{t^2}{\lambda_n^2}
\Big)^{1/2}
\hsx
\exp\Big(-\frac{t^2}{2 \hsy \lambda_n^2}\Big).
\]
\\[-.75cm]

\begin{x}{\small\bf LEMMA} \ 
There exists a positive integer $N$ with the property that 
\[
\max 
\Big(
- A_4, 
\hsx 
A_2 
+ 
\sum\limits_{k = 1}^n \ 
\frac{1}{\lambda_k^2} 
\hsy
\Big)
\ > \ 0
\qquad (n \geq N).
\]

PROOF \ 
Since
\[
\int\limits_{-\infty}^\infty \ 
\abs{F (\sqrt{-1} \hsx t)}
\ \td t
\ < \ 
\infty,
\]
$A_4$ must be $\leq 0$, thus matters are obvious if $A_4$ is $< 0$.  
Assume, therefore, that $A_4 = 0$ $-$then
\allowdisplaybreaks
\begin{align*}
\abs{F (\sqrt{-1} \hsx t)} \ 
&\geq \ 
\abs{M} \hsy \abs{t\hsy}^m 
\hsx
e^{- A_2 \hsy t^2}
\ 
\prod\limits_{n = 1}^\infty \ 
\exp\Big(-\frac{t^2}{2 \hsy \lambda_n^2}\Big).
\\[15pt]
&=\ 
\abs{M} \hsy \abs{t\hsy}^m 
\hsx
e^{- A_2 \hsy t^2}
\hsx
\exp
\Big(
\Big(
- 
\frac{1}{2}
\ 
\sum\limits_{n = 1}^\infty \ 
\frac{1}{\lambda_n^2}
\Big) 
\hsy t^2
\Big),
\end{align*}
so if 
\[
\frac{1}{2}
\
\sum\limits_{n = 1}^\infty \ 
\frac{1}{\lambda_n^2} 
\ < \ 
\infty,
\]
the condition on $A_2$ is that 
\[
- A_2 - 
\frac{1}{2}
\
\sum\limits_{n = 1}^\infty \ 
\frac{1}{\lambda_n^2} 
\ < \ 
0
\]
or still, 
\[
A_2 
+ 
\frac{1}{2}
\
\sum\limits_{n = 1}^\infty \ 
\frac{1}{\lambda_n^2} 
\ > \ 
0
\]
\qquad 
$\implies$
\[
A_2 
+ 
\sum\limits_{n = 1}^\infty \ 
\frac{1}{\lambda_n^2} 
\ > \ 
0
\]
\qquad 
$\implies$
\[
A_2 
+ 
\sum\limits_{k = 1}^n \ 
\frac{1}{\lambda_k^2} 
\ > \ 
0
\qquad (n \gg 0).
\]
However, in the event that
\[
\frac{1}{2}
\
\sum\limits_{n = 1}^\infty \ 
\frac{1}{\lambda_n^2} 
\ = \ 
\infty,
\]
then it is automatic that 
\[
\max 
\Big(
0, 
\hsx 
A_2 
\hsx + \hsx  
\sum\limits_{k = 1}^n \ 
\frac{1}{\lambda_k^2} 
\Big)
\ > \ 0
\]
$\forall \ n \gg 0$, there being in this case no condition on $A_2$.
\\[-.25cm]

Define $F_n (z)$ $(n \geq N)$ by
\allowdisplaybreaks
\begin{align*}
F_n (z) \ 
&= \ 
M \hsy z^m \hsx
\exp(A_4 \hsy z^4 + A_3 \hsy z^3 + A_2 \hsy z^2 + A_1 \hsy z)
\\[15pt]
&
\hspace{2cm}
\times  \ 
\prod\limits_{k = 1}^n \ 
\Big(
1 - \frac{z}{\lambda_k}
\Big)
\hsx
\exp
\Big(
\frac{z}{\lambda_k} 
+ \frac{z^2}{2 \hsy \lambda_k^2}
+ \frac{z^3}{3 \hsy \lambda_k^3}
\Big)
\\[15pt]
&\equiv \ 
P_n (z) 
\hsy
\exp(A_4 \hsy z^4 + A_{3,n} \hsy z^3 + A_{2, n} \hsy z^2 + A_{1, n} \hsy z),
\end{align*}
where
\[
P_n (z) 
\ = \ 
M \hsy z^m \
\prod\limits_{k = 1}^n \ 
\Big(
1 - \frac{z}{\lambda_k}
\Big)
\]
and 
\[
A_{j, n} 
\ = \ 
A_j 
\hsx + \hsx 
\frac{1}{j} 
\ 
\sum\limits_{k = 1}^n \ 
\frac{1}{\lambda_k^j} 
\qquad (j = 1, 2, 3),
\]
and set
\[
f_n (z) 
\ = \ 
\int\limits_{-\infty}^\infty \ 
F_n (\sqrt{-1} \hsx t) 
\hsx 
e^{\sqrt{-1} \hsx z \hsy t}
\ \td t.
\]
\\[-.75cm]
\end{x}

\begin{x}{\small\bf LEMMA} \ 
$\forall \ n \geq N$, 
$f_n \in \sL - \sP$.  
\\[-.5cm]

PROOF \ 
From the definitions, $F_n \in \sF_0$.  
But $F_n$ has finitely many zeros, hence by the earlier work, 
$f_n \in \sL - \sP$. 
\\[-.25cm]
\end{x}

\begin{x}{\small\bf LEMMA} \ 
$F_n \ra F$ uniformly on compact subsets of $\Cx$.
\\[-.25cm]
\end{x}

\begin{x}{\small\bf LEMMA} \ 
$\forall \ n \geq N$, 
\[
\abs{F_n (\sqrt{-1} \hsx t)} 
\ \leq \ 
\abs{F_N (\sqrt{-1} \hsx t} 
\qquad (t \in \R).
\]

PROOF \ 
This is because 
\[
\bigg|
\hsx
\bigg(
1 - \frac{\sqrt{-1} \hsx t)}{\lambda_n}
\bigg)
\hsx 
\exp
\bigg(
\frac{\sqrt{-1} \hsx t)}{\lambda_n}
\hsx + \hsx
\frac{(\sqrt{-1} \hsx t)^2}{2 \hsy \lambda_n^2}
\hsx + \hsx
\frac{(\sqrt{-1} \hsx t)^3}{3 \hsy \lambda_n^3}
\hsx
\bigg)
\hsx
\bigg|
\ \leq \ 
1
\]
for all $n$ and for all $t$.
\\[-.25cm]
\end{x}

Consequently, $f_n \ra f$ uniformly on compact subsets of $\Cx$, thus 37.3 can be invoked to conclude that 
$f \in \sL - \sP$, 
thereby finishing the proof of 37.2
\\[-.25cm]


\begin{x}{\small\bf LEMMA} \ 
If 
$F \in \sF_0$, 
then 
$\forall \ \lambda > 0$, the function 
\[
e^{\lambda \hsy z^2} \hsy F (z)
\]
is in $\sF_0$, hence the function 
\[
\int\limits_{-\infty}^\infty \ 
F (\sqrt{-1} \hsx t)
\hsx
e^{-\lambda \hsy t^2}
\hsx
e^{\sqrt{-1} \hsx z \hsy t}
\ \td t
\]
is in $\sL - \sP$ (cf. 37.2).
\\[-.5cm]

[Note: \ 
\allowdisplaybreaks
\begin{align*}
\Reg (- \lambda \hsy t^2 + \sqrt{-1} \hsx z \hsy t) \ 
&=\ 
- \lambda \hsy t^2 - t \ \Img z
\\[11pt]
&\leq\ 
- \lambda \hsy t^2 + \abs{t\hsy} \hsx \abs{\Img z}
\\[11pt]
&\leq\ 
- \lambda \hsy t^2 + \abs{t\hsy} \hsy \abs{z\hsy}.  
\end{align*}
As a function of $t$, the max of 
\[
- \lambda \hsy t^2 + \abs{t\hsy} \hsy \abs{z\hsy}
\]
is at 
$
\hsx
\ds
\abs{t\hsy} = \frac{\abs{z\hsy}}{2 \hsy \lambda}
$, 
\hsx
and the maximum value is 
\[
- \lambda 
\hsy
\frac{\abs{z\hsy}^2}{4 \hsy \lambda^2}
\hsx + \hsx 
\frac{\abs{z\hsy}}{2 \hsy \lambda}
\hsy \abs{z\hsy}
\ = \ 
\frac{\abs{z\hsy}^2}{4 \hsy \lambda}.
\]
And then
\[
\bigg|
\hsx
\int\limits_{-\infty}^\infty \ 
F (\sqrt{-1} \hsx t)
\hsx
e^{-\lambda \hsy t^2}
\hsx
e^{\sqrt{-1} \hsx z \hsy t}
\ \td t
\hsx
\bigg|
\ \leq \ 
\bigg(
\int\limits_{-\infty}^\infty \ 
\abs{F (\sqrt{-1} \hsx t)}
\ \td t
\bigg)
\hsx \exp
\Big(
\frac{\abs{z\hsy}^2}{4 \hsy \lambda}
\Big)
\hsx
.]
\]
\\[-1cm]
\end{x}

The foregoing considerations can, in a certain sense, be reversed.
\\[-.25cm]


\begin{x}{\small\bf THEOREM}\footnote[2]{\vspace{.11 cm}
C. Newman,  \textit{Proc. Amer. Math. Soc.} \textbf{61} (1976), pp. 245-251.}
\ 
Let $\mu$ be an even, finite, absolutely continuous Borel measure on the real line.  
Suppose that $\forall \ \lambda < 0$, the function 
\[
\int\limits_{-\infty}^\infty \ 
e^{\lambda \hsy t^2}
\hsy
e^{\sqrt{-1} \hsx z \hsy t}
\ \td \mu(t)
\]
has real zeros only $-$then 
\[
\td \mu(t)
\ = \ 
F (\sqrt{-1} \hsx t) \td t
\]
for some $F \in \sF_0$.
\\[-.25cm]
\end{x}

\qquad
{\small\bf \un{N.B.}} \ 
In this situation, 
$F (\sqrt{-1} \hsx t)$ is nonnegative, even, and admits the decomposition 
\[
F (\sqrt{-1} \hsx t)
\ = \ 
M \hsy t^{2 \hsy m} 
\hsx 
\exp
\Big(
-
\alpha \hsy t^4 - \beta \hsy t^2)
\ 
\prod\limits_j \ 
\Big(
1 + 
\frac{t^2}{a_j^2}
\Big)
\hsx
\exp\Big(-\frac{t^2}{a_j^2}\Big),
\]
where $M > 0$, 
$m = 0, 1, \ldots$, 
$a_j > 0$, 
$
\ 
\ds
\sum\limits_j \ 
\frac{1}{a_j^4} < \infty
\hsx 
$, 
$\alpha > 0$ and $\beta$ is real or 
$\alpha = 0$ 
and 
$
\ds
\beta + 
\sum\limits_j \ 
\frac{1}{a_j^2} > 0
$.]
\\

\begin{spacing}{1.75}
[Note: \ 
The product is over a set of $j$ which may be empty, finite, or infinite and the condition 
$
\ds
\beta + 
\sum\limits_j \ 
\frac{1}{a_j^2} > 0
$ 
is considered to be satisfied if 
$
\ds
\ 
\sum\limits_j \ 
\frac{1}{a_j^2} = \infty
\hsx .]$
\end{spacing}

\begin{x}{\small\bf SUBLEMMA} \ 
$\forall \ x \in \R$, 
\[
(1 + x^2) \hsy \exp (-x^2) 
\ \geq \ 
\exp (-x^4 / 2).
\]

PROOF \ 
$\forall \ y \geq 0$, 
\[
\log (1 + y) 
\ \geq \ 
y - \frac{y^2}{2}.
\]
Therefore

\[
1 + y
\ \geq \ 
\exp\Big(
y - \frac{y^2}{2}\Big)
\]
\qquad 
$\implies$

\[
(1 + y) \hsy \exp (-y) 
\ \geq \ 
\exp \Big(- \frac{y^2}{2}\Big).
\]
Now take $y = x^2$. 
\\[-.25cm]
\end{x}

\begin{x}{\small\bf APPLICATION} \ 
We have
\[
F (\sqrt{-1} \hsx t)
\ \geq \ 
M \hsy t^{2 \hsy m} 
\hsx 
\exp
\Big(
-
\Big(
\alpha 
\hsx + \hsx
\sum\limits_j \ 
\frac{1}{2 \hsy a_j^4}
\Big)
\hsy
t^4 
\hsx - \hsx
\beta \hsy t^2
\Big).
\]
\\[-1.25cm]
\end{x}

Let $\Phi \in \Lp^1(-\infty, \infty)$ be real analytic, positive and even.  
Assume: \ 
\[
\Phi (t)
\ = \ 
\tO
\big(
\exp\big(A \hsy \abs{t\hsy}^a - B \hsy e^{C \hsy \abs{t\hsy}^c}\big)\big)
\qquad (\abs{t\hsy} \ra \infty)
\]
for positive constants $A$, $a \geq 1$, $B$, $C$, $c \geq 1$.
\\[-.25cm]

\qquad
{\small\bf \un{N.B.}} \ 
Therefore $\Phi$ is of regular growth (cf. 35.14).  
\\[-.25cm]

Given any real $\lambda$, put 
\[
\Xi_\lambda (z) 
\ = \ 
\int\limits_{-\infty}^\infty \ 
\Phi (t) 
\hsx
e^{\lambda \hsy t^2}
\hsx
e^{\sqrt{-1} \hsx z \hsy t}
\ \td t.
\]
\\[-1cm]

\begin{x}{\small\bf THEOREM} \ 
If the zeros of $\Xi_0$  lie in the strip 
$\{z : \abs{\Img z} \leq \Delta\}$, 
then the zeros of $\Xi_\lambda$ $(\lambda > 0)$ are real provided 
$
\ds
\frac{\Delta^2}{2} \leq \lambda
$ 
and simple provided 
$
\ds
\frac{\Delta^2}{2} < \lambda
$ 
(cf. 36.20).
\\[-.25cm]
\end{x}

\begin{x}{\small\bf LEMMA} \ 
There does not exist an $F \in \sF_0$ such that 
$\Phi (t) = F(\sqrt{-1} \hsx t)$.  
\\[-.5cm]

PROOF \ 
For if this were the case, then 
\[
\Phi (t) 
\ \geq \ 
M \hsy t^{2 \hsy m} \hsx 
\exp
\Big(
-
\Big(
\alpha 
\hsx + \hsx
\sum\limits_j\ 
\frac{1}{2 \hsy a_j^4}
\Big)
\hsy
t^4 - \beta \hsy t^2
\Big)
\qquad (\tcf. \ 37.13),
\]
so 
\[
M \hsy t^{2 \hsy m} \hsx 
\exp
\Big(
-
\Big(
\alpha 
\hsx + \hsx
\sum\limits_j \ 
\frac{1}{2 \hsy a_j^4}
\Big)
\hsy
t^4 - \beta \hsy t^2
\Big)
\ = \ 
\tO
\big(
\exp\big(A \hsy \abs{t\hsy} - B \hsy e^{C \hsy \abs{t\hsy}}\big)\big).
\]
Setting $T = \abs{t\hsy}$, it thus follows that 
\[
\log M + 2 \hsy m \log T 
\hsx - \hsx
\Big(
\alpha 
\hsx + \hsx
\sum\limits_j \ 
\frac{1}{2 \hsy a_j^4}
\Big)
\hsy
T^4 
\hsx - \hsx
\beta \hsy T^2
\hsx - \hsx
A \hsy T 
\hsx + \hsx
B \hsy e^{C \hsy T}
\]
stays bounded as $T \ra \infty$, an absurdity.]
\\[-.5cm]
\end{x}

Supposing still that the zeros of 
$\Xi_0$
lie in the strip 
$\{z : \abs{\Img z} \leq \Delta\}$, 
there must exist a negative 
$\lambda_0$ 
such that 
$\Xi_{\lambda_0}$
has a nonreal zero (otherwise, taking 
$\td \mu (t) = \Phi (t) \td t$ 
in 37.11 forces 
$\Phi (t) = F(\sqrt{-1} \hsx t)$ 
for some 
$F \in \sF_0$
contradicting 37.15).
\\[-.5cm]

\begin{x}{\small\bf LEMMA} \ 
$\forall$ $\lambda < \lambda_0$, $\Xi_\lambda$ has a nonreal zero.
\\[-.5cm]

PROOF \ 
In fact, if all the zeros of 
$\Xi_\lambda$ 
were real, then all the zeros of 
$\Xi_{\lambda_0}$ 
would also be real (cf. 36.8).
\\[-.75cm]
\end{x}

Let $L$ be the set of 
$\lambda$ 
such that 
$\Xi_\lambda$ 
has a nonreal zero and let $R$ be the set of 
$\lambda$ 
such that all the zeros of 
$\Xi_\lambda$ 
are real $-$then
\[
\lambda_1 \in L, 
\  
\lambda_2 \in R
\implies
\lambda_1  \ < \ \lambda_2.
\]
Therefore the pair $(L,R)$ defines a Dedekind cut and we shall denote its cut point by 
$\Lambda_0$, 
hence
\[
\begin{cases}
\ 
\lambda \ < \ \Lambda_0
\implies 
\lambda \in L
\\[4pt]
\ 
\lambda \ > \ \Lambda_0
\implies 
\lambda \in R
\end{cases}
.
\]
\\[-.75cm]


\qquad
{\small\bf \un{N.B.}} \ 
A priori, 
\[
\Lambda_0
\ \leq \ 
\frac{\Delta^2}{2}
\qquad (\tcf. \ 37.14).
\]
\\[-1.5cm]

\begin{x}{\small\bf LEMMA} \ 
\[
\Lambda_0 \in R.
\]

\begin{spacing}{1.5}
PROOF \ 
Put 
$
\ds
\lambda_n = \Lambda_0 + \frac{1}{n}
$ 
$(n = 1, 2, \ldots)$ $-$then 
$\Xi_{\lambda_n} \ra \Xi_{\Lambda_0}$ 
uniformly on compact subsets of $\Cx$ 
(the assumptions serve to ensure that the 
$\Xi_{\lambda_n}$ 
constitute a normal family).  
But the zeros of 
$\Xi_{\lambda_n}$ 
are real and a zero of 
$\Xi_{\Lambda_0}$ 
is either a zero of 
$\Xi_{\lambda_n}$ 
for all sufficiently large values of $n$ or else is a limit point of the set of zeros of the 
$\Xi_{\lambda_n}$.  
And this means that the zeros of 
$\Xi_{\Lambda_0}$ 
are real, i.e., 
$\Lambda_0 \in R$.  
\\[-1.5cm]
\end{spacing}
\end{x}

\qquad
{\small\bf \un{N.B.}} \ 
Therefore $L$ consists of all $\lambda$ such that 
$\lambda < \Lambda_0$ and $R$ consists of all 
$\lambda$ such that $\Lambda_0 \leq \lambda$.
\\[-.5cm]

\begin{x}{\small\bf THEOREM} \ 
If  
$\lambda < \Lambda_0$, 
then 
$\Xi_\lambda$ 
has a nonreal zero and if 
$\Lambda_0 \leq \lambda$, 
then all the zeros of 
$\Xi_\lambda$
are real.
\\[-.5cm]

[This is a statement of recapitulation.]
\\[-.25cm]
\end{x}

\begin{x}{\small\bf THEOREM} \ 
Suppose that
$\Xi_\lambda$ 
has a multiple real zero $x_0$ $-$then 
$\lambda \leq \Lambda_0$.  
\\[-.5cm]

PROOF \ 
Take $x_0 = 0$ and in 36.19, take 
$f (z) = \Xi_\lambda (z)$ $-$then 
$\forall \ \delta > 0$ and sufficiently small, 
$e^{\delta \hsy D^2} \hsy \Xi_\lambda (z)$ 
has a nonreal zero.  
But
\allowdisplaybreaks
\begin{align*}
e^{\delta \hsy D^2} \hsy \Xi_\lambda (z)\ 
&=\ 
e^{\delta \hsy D^2} \hsy  
e^{-\lambda \hsy D^2} \hsy
\Xi_0 (z)
\qquad (\tcf. \ 36.12)
\\[15pt]
&=\ 
e^{(\delta - \lambda) \hsy D^2} \hsy
\Xi_0 (z)
\hspace{1.25cm} (\tcf. \ 36.16)
\\[15pt]
&=\ 
\Xi_{\lambda - \delta} \hsy (z)
\hspace{2.26cm} (\tcf. \ 36.12), 
\end{align*}
so
\[
\lambda - \delta
\ < \ 
\Lambda_0 
\hsx \implies \hsx 
\lim\limits_{\delta \ra 0} \ 
(\lambda - \delta)
\ \leq \ 
\Lambda_0 
\hsx \implies \hsx 
\lambda
\ \leq \ 
\Lambda_0 .
\]
\\[-1.25cm]
\end{x}

\begin{x}{\small\bf SCHOLIUM} \ 
If 
$\lambda > \Lambda_0$, 
then all the zeros of 
$\Xi_\lambda$ 
are real and simple.
\\[-.25cm]
\end{x}

\begin{x}{\small\bf APPLICATION} \ 
If 
$\Xi_0$ 
has a multiple real zero, then 
$0 \leq \Lambda_0$.  
\\[-.5cm]

[Note: \ 
If 
$\Xi_0$ 
has a nonreal zero, then $\Lambda_0 > 0 \hsx.]$
\\[-.25cm]
\end{x}

\begin{x}{\small\bf CRITERION} \ 
Suppose that there exists a 
$\lambda_0 < \Lambda_0$ 
with the property that 
$\forall \ \varepsilon > 0$, 
all but a finite number of zeros of  
$\Xi_{\lambda_0}$ 
lie in the strip 
$\abs{\Img z} \leq \varepsilon$ 
$-$then 
$\forall \ \lambda \in \ ]\lambda_0, \Lambda_0[\hsx$,  
$\Xi_\lambda \in  * - S - \sL - \sP$.
\\[-.5cm]

[By definition, 
\[
\Xi_{\lambda_0} (z)
\ = \ 
\int\limits_{-\infty}^\infty \ 
\Phi (t) 
\hsx 
e^{\lambda_0 \hsy t^2} 
\hsx 
e^{\sqrt{-1} \hsx z \hsy t}
\ \td t.
\]
Put
\[
\phi (t) \ = \ \Phi (t) 
\hsx 
\scalebox{1.25}{$\ds e^{\lambda_0 \hsy t^2} $}
\]
so that 
\allowdisplaybreaks
\begin{align*}
\Xi_{\lambda_0} (z) \ 
&=\ 
\int\limits_{-\infty}^\infty \ 
\phi (t) 
\hsx 
e^{\sqrt{-1} \hsx z \hsy t}
\ \td t
\\[15pt]
&=\ 
f_\infty (z).
\end{align*}
Pass now to 
\[
f_\infty (z; \lambda - \lambda_0)
\ = \ 
\int\limits_{-\infty}^\infty \ 
\phi (t) 
\hsy
e^{(\lambda - \lambda_0) \hsy t^2}
\hsy
e^{\sqrt{-1} \hsx z \hsy t}
\ \td t,
\]
a function in 
$* - S - \sL - \sP$ 
(cf. 36.33).  
But
\allowdisplaybreaks
\begin{align*}
f_\infty (z; \lambda - \lambda_0)
&=\ 
\int\limits_{-\infty}^\infty \ 
\Phi (t) 
\hsy
e^{\lambda \hsy t^2}
\hsy
e^{\sqrt{-1} \hsx z \hsy t}
\ \td t
\\[15pt]
&=\ 
\Xi_\lambda (z)\hsx .]
\end{align*}
\\[-.25cm]
\end{x}


\chapter{
$\boldsymbol{\S}$\textbf{38}.\quad   $\zeta$, $\xi$, AND $\Xi$}
\setlength\parindent{2em}
\setcounter{theoremn}{0}
\renewcommand{\thepage}{\S38-\arabic{page}}

\qquad
\vspace{-.5cm}
If 
$\zeta (s)$ 
is the Riemann zeta function and if 
\\[-.25cm]

\[
\xi (s) 
\ = \ 
\frac{s (s - 1)}{2} \ 
\pi^{^{\text{\small{$- \frac{s}{2}$}}}} \ 
\Gamma\Big(\frac{1}{2}\Big) \hsx \zeta (s)
\]
is the completed Riemann zeta function, then
\[
\xi (s) 
\ = \ 
\xi (1 - s). 
\]

\begin{x}{\small\bf NOTATION} \ 
Put
\[
\Xi (z) 
\ = \ 
\xi \Big(\frac{1}{2} + \sqrt{-1} \hsx z\Big).
\]
Then $\Xi$ is even, i.e., $\Xi (z) = \Xi (-z)$.
\\[-.25cm]
\end{x}

\begin{x}{\small\bf LEMMA} \ 
$\Xi$ is a real entire function of order 1 and of maximal type.
\\[-.5cm]
\end{x}

\begin{x}{\small\bf LEMMA} \ 
The zeros of $\Xi$ lie in the strip 
$
\ds
\Big\{ z: \abs{\Img z} < \frac{1}{2}\Big\}
$.
\\[-.25cm]

[Note: \ 
Recall that 
$\zeta (s)$ 
is zero free on the lines $\Reg s = 1$, $\Reg s = 0$.]
\\[-.5cm]
\end{x}

\begin{x}{\small\bf LEMMA} \ 
If
$\rho = \alpha + \sqrt{-1} \hsx \beta$ is a zero of $\Xi$, then 
\[
\begin{cases}
\ 
\bar{\rho} \ = \ \alpha - \sqrt{-1} \hsx \beta,  
\\[4pt]
- \rho \ = \ -\alpha - \sqrt{-1} \hsx \beta,  
\\[4pt]
-\bar{\rho} \ = \ -\alpha + \sqrt{-1} \hsx \beta,  
\end{cases}
\text{are also zeros of $\Xi$.}
\]
\\[-.5cm]
\end{x}

\begin{x}{\small\bf LEMMA} \ 
$\Xi$ has an infinity of zeros.
\\[-.5cm]
\end{x}

If $\rho_1, \rho_2, \ldots$ are the zeros of $\Xi$ and if $r_n = \abs{\rho_n}$, and if 
\[
0 < r_1 \leq r_2 \leq \cdots 
\qquad (r_n \ra \infty),
\]
then $\forall \ \varepsilon > 0$, 
\[
\sum\limits_{n = 1}^\infty \ 
\frac{1}{r_n^{1 + \varepsilon}}
\ < \ 
\infty
\]
but
\[
\sum\limits_{n = 1}^\infty \ 
\frac{1}{r_n}
\ = \ 
\infty.
\]

[Note: \ 
Therefore the convergence exponent of the zeros of $\Xi$ is equal to 1.]
\\[-.25cm]

\begin{x}{\small\bf LEMMA} \ 
$\gen \hsx \Xi = 1$ and 
\[
\Xi (z) 
\ = \ 
\Xi (0) \ 
\prod\limits_{n = 1}^\infty \ 
\Big(1 - \frac{z}{\rho_n}\Big)\hsx e^{z/\rho_n}.
\]

[Note: \ 
$\forall \ \rho$, 
\[
\Big(1 - \frac{z}{\rho}\Big)\hsx e^{z/\rho}
\hsx \cdot \hsx 
\Big(1 + \frac{z}{\rho}\Big)\hsx e^{-z/\rho}
\ = \ 
\Big(1 - \frac{z^2}{\rho^2}\Big).]
\]
Therefore
\[
\Xi 
\hsx \in \hsx  
\frac{1}{2} - \sL - \sP.
\]
\\[-1.25cm]
\end{x}

\begin{x}{\small\bf DEFINITION} \ 
The \un{Riemann Hypothesis} (RH) is the statement that all the zeros of $\Xi$ are real.
\\[-.25cm]
\end{x}

\begin{x}{\small\bf LEMMA} \ 
RH holds iff
\[
\Xi \in  \sL - \sP.
\]

[Note: \ 
Since $\sL - \sP$ is closed under differentiation, if the Riemann Hypothesis obtains, then $\forall \ n$, 
\[
\Xi^{(n)} (z) 
\ = \ 
\frac{\td^n}{\td z^n} \hsx \Xi 
\hsx \in \hsx  
\sL - \sP \hsx.]
\]
\\[-1.25cm]
\end{x}

\begin{x}{\small\bf THEOREM} \ 
$\Xi$ has an infinity of real zeros.
\\[-.5cm]

[There are a number of proofs of this result, one of which is delineated below.]
\\[-.25cm]
\end{x}


\begin{x}{\small\bf NOTATION} \ 
Put
\[
\Phi (t) 
\ = \ 
\sum\limits_{n = 1}^\infty \ 
\Big(
4 \hsy \pi^2 \hsy n^4 \hsy e^{^{\frac{9}{2} \hsy t}} 
\hsx - \hsx 
6 \hsy \pi \hsy n^2 \hsy e^{^{\frac{5}{2} \hsy t}} 
\big)
\hsx
\exp\big(- \pi \hsy n^2 \hsy e^{2 t}\big).
\]
\\[-1.25cm]
\end{x}

\begin{x}{\small\bf THEOREM} \ 
$\Xi$ and $\Phi$ are connected by the relation
\[
\Xi (z) 
\ = \ 
\int\limits_{-\infty}^\infty \
\Phi (t) 
\hsx 
e^{\sqrt{-1} \hsx z \hsy t}
\ \td t.
\]
\\[-1.25cm]
\end{x}

\begin{x}{\small\bf RAPPEL} \ 
The \un{theta function} is defined by 
\[
\theta (z) 
\ = \ 
\sum\limits_{n = -\infty}^\infty \
e^{\pi \hsy n^2 \hsy z}
\qquad (\Reg z > 0).
\]
\\[-1.25cm]
\end{x}

\begin{x}{\small\bf LEMMA} \ 
$\Phi$ and $\theta$ are connected by the relation
\[
\Phi (t) 
\ = \ 
\frac{1}{2} \ 
\bigg(
\frac{\td^2}{\td t^2} \hsx - \frac{1}{4} 
\bigg)
\Big(
e^{\frac{t}{2}} \hsx \theta (e^{2 t})
\Big).
\]
\\[-1.25cm]
\end{x}

\begin{x}{\small\bf LEMMA} \ 
$\Phi$ 
is an even function of $t$: $\Phi (t) \hsx = \hsx \Phi (-t)$.
\\[-.5cm]

PROOF \ 
In the functional equation
\[
\theta (x) 
\ = \ 
\Big(\frac{1}{x}\Big)^{1/2} 
\hsx 
\theta\Big(\frac{1}{x}\Big),
\]
take $x = e^{2 t}$, hence
\[
e^{\frac{t}{2}}
\hsx 
\theta (e^{2 t})
\ = \ 
e^{-\frac{t}{2}}
\hsx 
\theta (e^{-2 t}).
\]
\\[-1.25cm]
\end{x}

\begin{x}{\small\bf LEMMA} \ 
$\Phi$ 
is a positive function of $t$: $\Phi (t) > 0$.
\\[-.5cm]

[Note: \ 
In particular, 
\allowdisplaybreaks
\begin{align*}
\Xi (0) \ 
&=\ 
\int\limits_{-\infty}^\infty \
\Phi (t) 
\ \td t
\\[15pt]
&=\ 
2 \ 
\int\limits_0^\infty \ 
\Phi (t) 
\ \td t
\\[15pt]
&>\
0.]
\end{align*}
\\[-1cm]
\end{x}


\begin{x}{\small\bf LEMMA} \ 
We have
\[
\Phi (t) 
\ = \ 
\tO
\Big(
\exp
\Big(
\frac{9}{2} \hsx \abs{\hsy t \hsy} \hsx - \hsx \pi \hsy e^{2 \abs{\hsy t \hsy}}
\Big)
\Big)
\qquad \text{as} \ \abs{\hsy t \hsy} \ra \infty.
\]
\\[-1.25cm]
\end{x}

\begin{x}{\small\bf LEMMA} \ 
$\Phi (t)$ 
admits an analytic continuation into the strip 
$
\ds
\abs{\Img z} < \frac{\pi}{4}$ 
and 
$\forall \ n = 0, 1, 2, \ldots$, 
\[
\lim\limits_{t \ra \frac{\pi}{4}} \ 
\Phi^{(n)} (\sqrt{-1} \hsx t) 
\ = \ 0.
\]

[Note: \ 
$\Phi$ 
cannot be extended to an entire function.]
\\[-.25cm]
\end{x}

\qquad
{\small\bf \un{N.B.}} \ 
Therefore 
$\Phi$ 
is real analytic.
\\[-.25cm]

\begin{spacing}{1.75}
\begin{x}{\small\bf REMARK} \ 
The data above thus fits within the framework of \S37, 
viz. 
$\Phi \in \Lp^1 (-\infty, \infty)$ 
is real analytic, positive and even, the growth constants being 
$
\ds
A = \frac{9}{2}$, 
$a = 1$, 
$B = \pi$, 
$C = 2$, 
$c = 1$.
\\[-.5cm]

[Note: \ 
This theme is pursued in \S39.]
\\[-.75cm]
\end{x}
\end{spacing}

Here is Polya's proof of 38.9. 
To begin with, Fourier inversion is clearly possible, hence 
\[
\Phi (t) 
\ = \ 
\frac{1}{\pi} \ 
\int\limits_0^\infty \ 
\Xi (x) 
\hsx
\cos t x 
\ \td x,
\]
from which
\[
\Phi^{(2 n)} (t) 
\ = \ 
\frac{(-1)^n}{\pi} \ 
\int\limits_0^\infty \ 
\Xi (x) 
\hsx 
x^{2 n}
\hsx
\cos t x 
\ \td x.
\]
Write
\[
\Phi (\sqrt{-1} \hsx t) 
\ = \ 
c_0 \hsx + \hsx c_1 t^2 \hsx + \hsx c_2 t^4 \hsx + \hsx \cdots 
\qquad \Big(\abs{\hsy t \hsy} < \frac{\pi}{4}\Big),
\]
so
\[
(2 n)! \hsx c_n 
\ = \ 
(-1)^n 
\hsx \Phi^{(2 n)} (0) 
\ = \ 
\frac{1}{\pi} \ 
\int\limits_0^\infty \ 
\Xi (x) 
\hsx 
x^{2 n} 
\ \td x .
\]
To get a contradiction, suppose now that the sign of 
$\Xi (x)$ 
is eventually constant, say
$\Xi (x) > 0$ for $x > X$ $-$then  
\allowdisplaybreaks
\begin{align*}
\int\limits_0^\infty \ 
\Xi (x) 
\hsx 
x^{2 n} 
\ \td x \ 
&> \ 
\int\limits_{X + 1}^{X + 2} \ 
\Xi (x) 
\hsx 
x^{2 n} 
\ \td x 
\ - \ 
\int\limits_0^X
\abs{\Xi (x) }
\hsx 
x^{ 2 n} 
\ \td x 
\\[15pt]
&> \ 
(X + 1)^{2 \hsy n} \ 
\int\limits_{X + 1}^{X + 2} \ 
\Xi (x) 
\ \td x 
\ - \ 
X^{2 n} 
\int\limits_0^X \ 
\abs{\Xi (x)}
\ \td x 
\\[15pt]
&>\ 
0 
\qquad (n \gg 0)
\end{align*}
\qquad \qquad
$\implies$
\[
c_n > 0 
\qquad (n \gg 0).
\]
Therefore 
$\Phi^{(2 n)} \big(\sqrt{-1} \hsx t\big)$ 
increases monotonically in $t$ for $n \gg 0$, whereas 
\[
\Phi^{(2 n)}  (\sqrt{-1} \hsx t) \ra 0
\]
for 
$t \ra 0$, 
$
\ds 
t \ra \frac{\pi}{4}
$
(cf. 38.17).
\\[-.25cm]

\begin{x}{\small\bf LEMMA} \ 
If $t > 0$, then 
$\Phi^\prime (t) < 0$.
\\[-.5cm]

[This is a brute force computation (see the Appendix to \S42 for the ``how to'').]
\\[-.25cm]
\end{x}

\begin{x}{\small\bf LEMMA} \ 
$\Phi$ is a strictly decreasing function of $t$ on $[0, \infty[$.
\\[-.25cm]
\end{x}


\chapter{
$\boldsymbol{\S}$\textbf{39}.\quad  THE de BRUIJN-NEWMAN CONSTANT}
\setlength\parindent{2em}
\setcounter{theoremn}{0}
\renewcommand{\thepage}{\S39-\arabic{page}}

\qquad
Take $\Xi$ and $\Phi$ as in \S38, hence

\[
\Xi (z) 
\ = \ 
\int\limits_{-\infty}^\infty \ 
\Phi (t) 
\hsx
e^{\sqrt{-1} \hsx z \hsy t}
\ \td t
\qquad (\tcf. \ 38.11),
\]
\begin{spacing}{1.65}
\noindent
and $\Phi$ meets the growth requirements per \S37 (cf. 38.18).  
Since the zeros of $\Xi$ lie in the strip 
$
\ds
\Big\{z : \abs{\Img z} < \frac{1}{2}\Big\}
$ (cf. 38.3), 
\end{spacing}
\[
\Delta
\ = \ 
\frac{1}{2}
\hsx \implies \hsx 
\frac{\Delta^2}{2\ } 
\ = \ 
\frac{1}{8}.
\]
Given a real $\lambda$, set
\[
\Xi_\lambda (z) 
\ = \ 
\int\limits_{-\infty}^\infty \ 
\Phi (t) 
\hsx e^{\lambda \hsy t^2} 
\hsx
e^{\sqrt{-1} \hsx z \hsy t}
\ \td t
\qquad (\Xi_0 \hsx = \hsx  \Xi).
\]
\begin{spacing}{1.5}
\noindent
Then the zeros of $\Xi_\lambda$ ($\lambda > 0$) are real provided 
$
\ds
\frac{1}{8} \leq \lambda$
and simple provided 
$
\ds
\frac{1}{8} < \lambda$
(cf. 37.14).  
Now introduce $\Lambda_0$ and recall: \ 
If 
$\lambda < \Lambda_0$, 
then 
$\Xi_\lambda$ 
has a nonreal zero and if 
$\Lambda_0 \leq \lambda$, 
then all the zeros of
$\Xi_\lambda$ 
are real (cf. 37.18).
\\[-.25cm]
\end{spacing}

\qquad
{\small\bf \un{N.B.}} \ 
It is automatic that 
\[
\Lambda_0 
\ \leq \ 
\frac{1}{8} .
\]
\\[-1.5cm]

\begin{x}{\small\bf DEFINITION} \ 
$\Lambda_0$ 
is called the 
\un{de Bruijn-Newman constant}.  
\\[-.25cm]

[Note: \ 
Some authorities reserve this term for $4 \hsy \Lambda_0$.]
\\[-.25cm]
\end{x}

\begin{x}{\small\bf LEMMA} \ 
RH holds iff $\Lambda_0 \leq 0$.
\\[-.25cm]
\end{x}

\begin{x}{\small\bf THEOREM}\footnote[2]{\vspace{.11 cm}B. Rogers and T. Tau, \url{https://arxiv.org/pdf/1801.05914}\hsy.}\ 
(Newman Conjecture
\footnote[3]{\vspace{.11 cm}C. Newman,  \textit{Proc. Amer. Math. Soc.} \textbf{61} (1976), pp.245-251.}
) \ 
\\[-.5cm]
\[
\Lambda_0 
\ \geq \ 
0.
\]

[Note: \ 
This is 
``a quantitative version of the dictum that the Riemann Hypothesis, if true, is only barely so''.]
\\[-.25cm]

[Note: \ 
It is true but not obvious that 
$
\ds
\Lambda_0 < \frac{1}{8}$ 
(cf. 39.11).]
\\[-.25cm]

We shall pass on the proof of this result and instead refer the reader to the paper by Rogers and Tau.  
\\[-.25cm]
\end{x}




\begin{x}{\small\bf REMARK} \ 
A simpler proof of the Newman conjecture has been subsequently given by 
A. Dobner\footnote[4]{\vspace{.11 cm}A. Dobner,  \url{https://arxiv.org/pdf/2005.05142}\hsy.}.\ 
Dobner's proof does not rely on any information about the zeros of
the zeta function and, moreover, extends the result more generally to a broad class of $L$-functions. 
\\[-.25cm]
\end{x}

\begin{x}{\small\bf LEMMA} \ 
If $f$ is an entire function of order $< 2$, then the order of 
\[
e^{\lambda \hsy D^2} \hsy f (z)
\]
is $< 2$ (cf. 36.15) and, in fact, the orders of $f (z)$ and $e^{\lambda \hsy D^2} \hsy f (z)$ are equal.
\\[-.25cm]
\end{x}

\begin{x}{\small\bf APPLICATION} \ 
$\Xi_\lambda$ is a real entire function of order 1.
\\[-.5cm]

[Thanks to 36.12, 
\[
\Xi_\lambda (z) 
\ = \ 
e^{-\lambda \hsy D^2} 
\hsy 
\Xi (z).
\]
\\[-1.25cm]
\end{x}

\begin{x}{\small\bf LEMMA} \ 
$\Xi_\lambda$ is of maximal type.
\\[-.5cm]

PROOF \ 
If 
$\Xi_\lambda$ were of finite type, then 
$\Xi_\lambda$ 
would be of exponential type but this is ruled out by the Paley-Wiener theorem (cf. 22.7).
\\[-.25cm]
\end{x}

On general grounds, 
$\Xi_\lambda$ 
has an infinity of zeros but more is true: \ 
$\Xi_\lambda$ 
has an infinity of real zeros (argue as in 38.9).
\\[-.25cm]


\begin{x}{\small\bf LEMMA}\footnote[3]{\vspace{.11 cm}H. Ki et al.,  \textit{Advances in Math.} \textbf{222} (2009), pp. 281-306.}
\ 
Take $\lambda > 0$ $-$then 
$\forall \ \varepsilon > 0$, 
all but a finite number of zeros of 
$\Xi_\lambda (z)$ lie in the strip 
$\abs{\Img z} \leq \varepsilon$.
\\[-.25cm]
\end{x}

\begin{x}{\small\bf APPLICATION} \ 
$\forall \ \lambda > 0$,  
all but a finite number of zeros of 
$\Xi_\lambda$ 
are real and simple (cf. 36.35).
\\[-.25cm]
\end{x}

\begin{x}{\small\bf LEMMA} \ 
Suppose that 
$
\ds
0 < \lambda < \frac{1}{8}
$ $-$then the zeros of 
$\Xi_\lambda$  
lie in the strip 
\[
\{z : \abs{\Img z} \leq A_\lambda\}
\]
for some 
$
\ds
A_\lambda < \Big(\frac{1}{4} - 2 \hsy \lambda \Big)^{1/2}
$. 
\\[-.25cm]

\begin{spacing}{1.5}
PROOF \ 
Choose $\lambda_0 : 0 < \lambda_0 < \lambda$ 
and put 
$
\ds
A_0 = 
\Big(\frac{1}{4} - 2 \hsy \lambda_0\Big)^{1/2}$.  
Since the zeros of 
$\Xi_0 \  (= \Xi)$  
are confined to the strip 
$
\ds
\Big\{z : \abs{\Img z} \leq \frac{1}{2}\Big\}
$
and since
$
\hsx
\ds
\Xi_{\lambda_0} = e^{-\lambda_0 \hsy D^2} \hsy \Xi_0$, 
it follows from 36.5 (and subsequent comment) that the zeros of 
$\Xi_{\lambda_0}$ are confined to the strip 
$
\ds
\{z : \abs{\Img z} \leq A_0\}
$
(the $A^2$ there is 
$
\ds
\Big(\frac{1}{2}\Big)^2$ 
here
$
(f_\infty = \Xi_0)$).  
On the other hand, the number of nonreal zeros of 
$\Xi_{\lambda_0}$ 
is finite (cf. 39.9) and 
$\Xi_{\lambda_0}$ 
has an infinity of real zeros.  
Observing now that 
\[
2 \hsy (\lambda - \lambda_0) 
\ < \ 
A_0^2 
\ = \ 
\frac{1}{4} - 2 \hsy \lambda_0,
\]
on the basis of 36.37, the zeros of 
\end{spacing}
\allowdisplaybreaks
\begin{align*}
\Xi_{\lambda} \ 
&= \ 
e^{-\lambda \hsy D^2} \hsy \Xi_0
\\[11pt]
&= \ 
e^{-(\lambda + \lambda_0 - \lambda_0)\hsy D^2} \hsy \Xi_0
\\[11pt]
&= \ 
e^{-(\lambda - \lambda_0 )\hsy D^2} \hsy e^{- \lambda_0 \hsy D^2} \hsx \Xi_0
\qquad (\tcf. \ 36.16)
\\[11pt]
&= \ 
e^{-(\lambda - \lambda_0 )\hsy D^2}  \hsy\Xi_{\lambda_0}
\end{align*}
lie in the strip
\[
\{z : \abs{\Img z} \leq A_\lambda\}
\]
for some
\[
A_\lambda 
\ < \ 
(A_0^2 - 2(\lambda - \lambda_0))^{1/2}
\ = \ 
\Big(\frac{1}{4} - 2 \hsy \lambda\Big)^{1/2}.
\]
\\[-1.25cm]
\end{x}

\begin{x}{\small\bf THEOREM} \ 
The de Bruijn-Newman constant $\Lambda_0$ is $< \frac{1}{8}$.
\\[-.25cm]

PROOF \ 
Fix 
$\lambda : 0 < \lambda < \ds\frac{1}{8}\ $ 
and choose 
$\lambda_0$ subject to 
\[
A_\lambda^2 
\ < \ 
2 \hsy \lambda_0
\ < \ 
\frac{1}{4} - 2 \hsy \lambda, 
\]
hence
\[
2 \hsy \lambda + 2 \hsy \lambda_0 
\ < \ 
\frac{1}{4} 
\hsx \implies \ 
\lambda + \lambda_0
\ < \ 
\frac{1}{8}.
\]
Now take in 36.22 
$f = \Xi_\lambda$, 
$A = A_\lambda$ 
and conclude that the zeros of 
\[
e^{- \lambda_0 \hsy D^2} \hsx \Xi_\lambda 
\]
are real.  
But
\begin{align*}
e^{- \lambda_0 \hsy D^2} \hsx \Xi_\lambda \ 
&=\ 
e^{- \lambda_0 \hsy D^2} 
\hsx
e^{- \lambda \hsy D^2} 
\hsx
\Xi_0
\\[8pt] \
&=\ 
e^{-(\lambda +  \lambda_0) \hsy D^2} 
\hsx
\Xi_0
\qquad (\tcf. \ 36.16)
\\[8pt] \
&=\ 
\Xi_{\lambda + \lambda_0}.
\end{align*}
And this implies that
\[
\Lambda_0
\ \leq \ 
\lambda + \lambda_0 
\ < \ 
\frac{1}{8}.
\]
\\[-1.25cm]
\end{x}

\begin{x}{\small\bf REMARK} \ 
It has been shown by
D. H. J. Polymath\footnote[2]{\vspace{.11 cm}D.H.J.$\hsy$Polymath,$\hsx$\url{https://github.com/km-git-acc/dbn_upper_bound/blob/master/Writeup/debruijn.pdf};$\hsx$also$\hsx$\url{https://terrytao.wordpress.com/2018/05/04/polymath15-ninth-thread-going-below-0-22/}\hsy} 
\ 
that
$4 \hsy \Lambda_0 \ \leq \ 0.22$.  
A subsequent observation to the calculations in the reference above by 
D. Platt and D. Trudgian\footnote[3]{\vspace{.11 cm}D. Platt and D. Trudgian,  \url{https://arxiv.org/pdf/2004.09765f}}
has yielded 
$4 \hsy \Lambda_0 \ \leq \ 0.2$.
\\[-.25cm]
\end{x}

\begin{x}{\small\bf REMARK} \ 
Consider $\Xi_{1/8}$ $-$then its zeros are real and simple (cf. 37.20).
\\[-.75cm]
\end{x}

Per 
\[
\Xi^{(n)}  (z)
\ = \ 
\frac{\td^n}{\td z^n}
\
\Xi,
\]

\noindent
one has the analog of $\Lambda_0$, call it $\Lambda_0^{(n)}$ 
$\big(\Lambda_0 = \Lambda_0^{(0)}\big).$ 
\\

\qquad
{\small\bf \un{N.B.}} \ 
\[
\Xi_\lambda^{(n)} (z)
\ = \ 
e^{- \lambda \hsy D^2} 
\hsx 
\Xi^{(n)}  (z).
\]
\\[-1.25cm]

\begin{x}{\small\bf THEOREM} \ 
The sequence 
$\{\Lambda^{(n)}\}$ 
is decreasing and its limit is $\leq 0$.
\\[-.5cm]

\begin{spacing}{1.5}
PROOF \ 
By definition, 
$\Lambda^{(n)}$ 
is the infimum of the set of $\lambda$ such that 
$\Xi_\lambda^{(n)}$ has real zeros only.  
But if 
$\hsx \Xi_\lambda^{(n)}$ 
has real zeros only, then the same is true of 
$\hsx \Xi_\lambda^{(n+1)}$, 
hence
$\Lambda^{(n+1)} \leq \Lambda^{(n)}$.  
Next, 
$\forall \ \lambda > 0$,  
$\Xi_\lambda$ has at most a finite number of nonreal zeros (cf. 39.9), 
thus 
$\hsx \Xi_\lambda \in * - \sL - \sP$, 
so 
$\exists \ n$ : 
$\Xi_\lambda^{(n)}$ 
is in 
$\sL - \sP$ (cf. 11.9)
from which 
$\Lambda^{(n)} \leq \lambda$.  
Now send $\lambda$ to 0 and conclude that 
\end{spacing}
\[
\lim\limits_{n \ra \infty} \ 
\Lambda^{(n)}
\ \leq \ 
0.
\]
\\[-.25cm]
\end{x}


\chapter{
$\boldsymbol{\S}$\textbf{40}.\quad  TOTAL POSITIVITY}
\setlength\parindent{2em}
\setcounter{theoremn}{0}
\renewcommand{\thepage}{\S40-\arabic{page}}

 
\qquad 
A sequence 
$\{c_n : n \geq 0\}$ 
$(c_0 \neq 0)$ 
of real numbers is said to be \un{totally positive} if all the minors of all the orders of the infinite lower triangular matrix

\[
\gC: \quad
\begin{bmatrix}
\hsx
c_0 
&
0
&
0
&
0
&
0
&\cdots
\\[4pt]
\hsx
c_1 
&
c_0
&
0
&
0
&
0
&\cdots
\\[4pt]
\hsx
c_2 
&
c_1
&
c_0
&
0
&
0
&\cdots
\\[4pt]
\hsx
c_3 
&
c_2
&
c_1
&
c_0
&
0
&\cdots
\\[4pt]
&& \cdots
\end{bmatrix}
\]
are nonnegative.
\\[-.5cm]

[Note: \ 
Therefore the $c_n$ are nonnegative.]
\\[-.25cm]

\begin{x}{\small\bf LEMMA} \ 
If for some $n$, $c_n = 0$, then 
$\forall \ k = 1, 2, \ldots$, $c_{n+k} = 0$. 
\\[-.5cm]

PROOF \ 
The minor

\[
\begin{vmatrix*}[l]
\hsx
c_n
&c_0 \text{\hsy}
\\[4pt]
\hsx
c_{n+k}
& c_k\ 
\end{vmatrix*}
\ = \ 
- c_0 \hsx c_{n+k}
\]
\\[-.75cm]
\end{x}
is nonnegative.  
But $c_0$ is $\hsx > 0 \hsx$ and $c_{n+k}$ is $\geq 0$, hence $\hsx c_{n+k} = 0$.
\\[-.25cm]

With the understanding that $c_n = 0$, if $n < 0$, put
\[
D (n, r) 
\ = \ 
\begin{vmatrix*}[l]
\ c_n 
&c_{n-1} 
&\cdots 
&c_{n- r +1}  
\\[4pt]
\ c_{n+1} 
&c_n 
&\cdots 
&c_{n- r +2}   
\\[4pt]
\ \ \vdots 
&\vdots
&
&\vdots
\\[4pt]
\ c_{n + r - 1} 
&c_{n + r - 2}
&\cdots
&c_n 
\end{vmatrix*}
\hsx .
\]
Here $n = 0, 1, 2, \ldots $, while 
$r = 1, 2, 3, \ldots \ $. 
\\[-.25cm]

\begin{x}{\small\bf EXAMPLE} \ 
Take $r = 1$ $-$then
\[
D(n, 1) \ = \ c_n.
\]
\\[-1cm]
\end{x}

\begin{x}{\small\bf EXAMPLE} \ 
Take $r = 2$, $-$then

\[
D(n, 2) \ = \ 
\begin{vmatrix}
\ 
c_n
&c_{n-1} \text{\ } 
\\[4pt]
\ 
c_{n+1} 
&c_n 
\end{vmatrix}
.
\]
In particular: 

\[
D(0, 2) \ = \ 
\begin{vmatrix}
\ 
c_0
&0 \ 
\\[4pt]
\ 
c_1 
&c_0 \text{\ }  
\end{vmatrix}
.
\]
\\[-.5cm]
\end{x}

\begin{x}{\small\bf EXAMPLE} \ 
Take $r = 3$, $-$then

\[
D(n, 3) \ = \ 
\begin{vmatrix}
\ 
c_n
&c_{n-1} 
&c_{n-2} \text{\ }  
\\[4pt]
\ 
c_{n+1} 
&c_n 
& c_{n-1}
\\[4pt]
\ 
c_{n+2}
&c_{n+1} 
&c_n\ 
\end{vmatrix}
.
\]
In particular: 
\[
D(0, 3) \ = \ 
\begin{bmatrix}
\ 
c_0 
&
0
&
0
\\[4pt]
\
c_1 
&
c_0 
&
0
\\[4pt]
\
c_2
&
c_1
&
c_0  \text{\ }
\end{bmatrix}
, \quad
D(1, 3) \ = \ 
\begin{vmatrix}
\ 
c_1 
&
c_0 
&
0
\\[4pt]
\
c_2
&
c_1 
&
c_0
\\[4pt]
\
c_3
&
c_2
&
c_1  \text{\ }
\end{vmatrix}
.
\]
\\[-.25cm]
\end{x}


\begin{x}{\small\bf FEKETE CRITERION} \ 
A sequence 
$\{c_n : n \geq 0\}$ 
$(c_0 \neq 0)$ 
of nonnegative real numbers is totally positive if

\[
\forall \ n, \forall \ r, \ D(n, r) \ > \ 0.
\]
\\[-1.25cm]
\end{x}

\begin{x}{\small\bf THEOREM}\footnote[2]{\vspace{.11 cm}
M. Aissen et al.,  \textit{Proc. Nat. Acad. Sci. U.S.A.} \textbf{37} (1951), pp. 303-307.}
 \ 
Suppose that

\[
f (z) 
\ = \ 
\sum\limits_{n = 0}^\infty \ 
c_n \hsy z^n
\]
is a real entire function with $f (0) > 0$ $-$then the sequence 
$c_0, c_1, c_2, \ldots$ is totally positive iff $f$ has a representation of the form

\[
f (z) 
\ = \ 
f (0) \hsy e^{a z} \ 
\prod\limits_{n = 1}^\infty \ 
\Big(1 - \frac{z}{\lambda_n}\Big),
\]
where $a$ is real and $\geq 0$, the $\lambda_n$ are real and $< 0$ with 
$
\ds
\ 
\sum\limits_{n = 1}^\infty \ 
\frac{1}{\lambda_n} 
< 
\infty
$.
\\[-.25cm]
\end{x}

\begin{x}{\small\bf EXAMPLE} \ 
Take $f (z) = e^z$ $-$then the sequence 
$\ds \frac{1}{0 !}, \hsx  \frac{1}{1 !}, \hsx  \frac{1}{2 !}, \ldots$ 
is totally positive.  
\\[-.25cm]
\end{x}

\begin{x}{\small\bf EXAMPLE} \ 
Take $f (z) = (1 + z)^n$ $-$then the sequence 
$\ds \binom{n}{0},  \binom{n}{1},  \binom{n}{2}, \ldots$ 
is totally positive. 
\\[-.25cm]
\end{x}

\begin{x}{\small\bf RAPPEL} \ 
(cf. 10.11) Let $f \not\equiv 0$ be a real entire function $-$then 
$f \in \ent(\hsy]-\infty, 0]\hsy)$ 
iff $f$ has a representation of the form

\[
f (z) 
\ = \ 
C \hsy z^m \hsy e^{a z} \ 
\prod\limits_{n = 1}^\infty \ 
\Big(1 - \frac{z}{\lambda_n}\Big),
\]
where $C \neq 0$ is real, $m$ is a nonnegative integer, $a$ is real and $\geq 0$, 
the $\lambda_n$ are real and $< 0$ with 
$
\ds
\ 
\sum\limits_{n = 1}^\infty \ 
\frac{1}{\lambda_n} 
< 
\infty
$.
\\
\end{x}

\begin{x}{\small\bf NOTATION} \ 
Denote by

\[
\ent_+(\hsy]-\infty, 0]\hsy)
\]
the subset of 
$\ent(\hsy]-\infty, 0]\hsy)$ 
(cf. 10.26) consisting of those $f$ such that 

\[
f (z) 
\ = \ 
f (0) \hsy e^{a z} \ 
\prod\limits_{n = 1}^\infty \ 
\Big(1 - \frac{z}{\lambda_n}\Big)
\]
with $f (0) > 0$.
\\[-.25cm]
\end{x}

\begin{x}{\small\bf SCHOLIUM} \ 
If 

\[
f (z) 
\ = \ 
\sum\limits_{n = 0}^\infty \ 
c_n \hsy z^n
\]
is a real entire function with $f (0) > 0$, then the sequence 
$c_0, c_1, c_2, \ldots$ is totally positive iff 

\[
f \in \ent_+(\hsy]-\infty, 0]\hsy).
\]
\\[-1.25cm]
\end{x}

\begin{x}{\small\bf NOTATION} \ 
Write
\[
\gC : \big[c_{i - j}\big]_{i=1, \hsx j= 1}^\infty \hsx .
\]
So, e.g., 
\[
c_{1-1} = c_0, \quad
c_{1-2} = 0, \quad
c_{2-1} = c_1, \quad
c_{2-2} = c_0, \quad
c_{2-3} = 0, \quad
\quad \text{etc}.
\]
\\[-1.25cm]
\end{x}

\begin{x}{\small\bf NOTATION} \ 
Given a positive integer $n$, let

\[
\begin{cases}
\ 
1 \leq i_1 < i_2 < \cdots < i_n
\\[4pt]
\ 
1 \leq j_1 < j_2 < \cdots < j_n
\end{cases}
\]
be positive integers and let

\[
\gC (i_1, i_2, \ldots, i_n \hsx \big| \hsx j_1, j_2, \ldots, j_n)
\]
denote the $n \times n$ minor obtained from $\gC$ by deleting all the rows and columns except those labeled 
$i_1, i_2, \ldots, i_n$ and $j_1, j_2, \ldots, j_n$ 
respectively.
\\[-.5cm]
\end{x}

\begin{x}{\small\bf THEOREM}\footnote[2]{\vspace{.11 cm}
S. Karlin,  \textit{Total Positivity}, Stanford University Press, (1968) pp. 427-432.} 
\ 
Let

\[
f \in \ent_+(\hsy]-\infty, 0]\hsy).
\]
Assume: \ 
$a$ is equal to 0, the $c_n$ are greater than 0, and the product

\[
\prod\limits_{n = 1}^\infty \ 
\Big(1 - \frac{z}{\lambda_n}\Big)
\]
is infinite $-$then the minor

\[
\gC (i_1, i_2, \ldots, i_n \hsx \big| \hsx j_1, j_2, \ldots, j_n)
\]
is positive if 
$j_1 \leq i_1, j_2 \leq i_2, \ldots, j_n \leq i_n$. 
\\[-.25cm]
\end{x}

\begin{x}{\small\bf APPLICATION} \ 
For 
$n = 0, 1, 2, \ldots$ 
and 
$r = 1, 2, 3, \ldots$, 

\[
D (n, r) 
\ = \ 
\gC (n+1, n+2, \ldots, n + r \hsx \big| \hsx  1, 2, \ldots, r),
\]
so $D( n, r)$ is positive.
\\[-.25cm]
\end{x}

\begin{x}{\small\bf EXAMPLE} \ 
\allowdisplaybreaks
\begin{align*}
D(n, 2) \ 
&= \ 
\begin{vmatrix}
\ 
c_n
&c_{n-1} \ 
\\[4pt]
\ 
c_{n+1} 
&c_n\ 
\end{vmatrix}
\\[18pt]
&= \
c_n^2 - c_{n-1} \hsy c_{n+1}
\\[18pt]
&= \
\gC (n + 1, n + 2 \hsx | \hsx 1, 2) 
\\[18pt]
&> \ 
0.
\end{align*}

[Note: \ 
\[
D (n, 1) 
\ = \ 
c_n
\ = \ 
\gC (n + 1 \hsx | \hsx 1) 
\ > \ 0.]
\]
\\[-1cm]
\end{x}

\begin{x}{\small\bf LEMMA} \ 
Suppose that 
\[
f (z) 
\ = \ 
\sum\limits_{n = 0}^\infty \ 
c_n \hsy z^n
\]
is a real entire function with $f (0) > 0$ and $\forall \ n$, $c_n \geq 0$. 
Assume: \ 
$f \in \sL - \sP$ $-$then
\[
f \in \ent_+(\hsy]-\infty, 0]\hsy).
\]
\\[-1cm]
\end{x}

\begin{x}{\small\bf EXAMPLE} \ 
Take
\[
f (z) 
\ = \ 
\sum\limits_{n = 0}^\infty \ 
\frac{1}{e^{n^2}} \hsx z^n.
\]
Then
\[
f \in \ent_+(\hsy]-\infty, 0]\hsy).
\]

[The Jensen polynomials

\[
\tJ_n (f; z) 
\ = \ 
\sum\limits_{k = 0}^n \ 
\binom{n}{k} \hsx \frac{k !}{e^{k^2}} \hsx z^k
\]
associated with $f$ have real zeros only, thus 
$f \in \sL - \sP$ (cf. 12.14).]
\\[-.25cm]
\end{x}


\chapter{
$\boldsymbol{\S}$\textbf{41}.\quad  CHANGE OF VARIABLE}
\setlength\parindent{2em}
\setcounter{theoremn}{0}
\renewcommand{\thepage}{\S41-\arabic{page}}


\qquad 
Continuing the discussion initiated in \S38, from the definitions

\allowdisplaybreaks
\begin{align*}
\Xi\Big(\frac{z}{2}\Big) \ 
&=\ 
\int\limits_{-\infty}^\infty \ 
\Phi (t) 
\hsy 
\scalebox{1.25}{$\ds e^{\sqrt{-1} \hsx \frac{z}{2} \hsy t}$}
\ \td t
\\[15pt]
&=\ 
2 \ 
\int\limits_0^\infty \ 
\Phi (t) 
\hsy 
\cos z \hsx \frac{t}{2} 
\ \td t
\\[15pt]
&=\ 
4 \ 
\int\limits_0^\infty \ 
\Phi (2 t) 
\hsy 
\cos z \hsy  t 
\ \td t
\\[15pt]
&=\ 
8 \ 
\int\limits_0^\infty \ 
\Phi (t) 
\hsy 
\cos z \hsy  t 
\ \td t,
\end{align*}

\noindent
where, in a flagrant abuse of notation, the ``new'' $\Phi (t)$ is
\[
\Phi (t) 
\ = \ 
\sum\limits_{n = 1}^\infty \ 
\big(
2 \hsy \pi^2 \hsy n^4 \hsy e^{9 t} 
\hsy - \hsy 
3 \hsy \pi \hsy n^2 \hsy e^{5 t}
\big)
\hsx 
\exp\big(- \pi \hsy  n^2 \hsy e^{4 \hsy t}\big).
\]
Expand now the cosine and integrate term by term to get the representation
\allowdisplaybreaks
\begin{align*}
\Sh (z) \ 
&\equiv \ 
\frac{1}{8} \ \Xi\Big(\frac{z}{2}\Big)
\\[15pt]
&=\ 
\sum\limits_{k = 0}^\infty \ 
\frac{(-1)^k}{(2 k)!}
\hsx
b_k \hsy z^{2 k}.
\end{align*}
Here
\[
b_k 
\ = \ 
\int\limits_0^\infty \ 
t^{2 k} 
\hsy 
\Phi (t) 
\ \td t.
\]
\\[-1cm]

\begin{x}{\small\bf NOTATION} \ 
Put

\[
F_\zeta (z)
\ = \ 
\sum\limits_{k = 0}^\infty \ 
\frac{b_k}{(2 k)!} \hsy z^k
\]
and set
\[
C_k 
\ = \ 
\frac{b_k}{(2 k)!}.
\]
\\[-1.25cm]
\end{x}

Accordingly, 
\[
\Sh (z) 
\ = \ 
F_\zeta (- z^2).
\]
Therefore if $z_0$ is a zero of $\Sh (z)$, then $- z_0^2$ is a zero of $F_\zeta (z)$.
\\[-.25cm]

\begin{x}{\small\bf LEMMA} \ 
$F_\zeta$ is a real entire function of order 
$\ds
\frac{1}{2}
$ and of maximal type.
\\[-.5cm]
\end{x}

\begin{x}{\small\bf LEMMA} \ 
$\forall \ k \geq 0$, $C_k$ is positive (cf. 38.15).
\\[-.5cm]
\end{x}

\qquad
{\small\bf \un{N.B.}} \ 
In particular: 
\[
F_\zeta (0) 
\ = \ 
C_0 
\ > \ 
0.
\]
\\[-1.75cm]

\begin{x}{\small\bf SCHOLIUM} \ 
RH is equivalent to the statement that all the zeros of $F_\zeta$ are real and negative.
\\[-.25cm]
\end{x}

\begin{x}{\small\bf SCHOLIUM} \ 
RH is equivalent to the statement that 
\[
F_\zeta \in \ent_+ \big(\hsy]-\infty, 0]\hsy\big).
\]
\\[-1.5cm]
\end{x}

\begin{x}{\small\bf THEOREM} \ 
If RH obtains, then
\[
\forall \ n, \ \forall  \ r, \ 
D (n, r) 
\ > \ 0.
\]
PROOF \ 
In fact, 
\[
\tRH
\implies
F_\zeta 
\hsx \in \hsx 
\ent_+ \big(\hsy]-\infty, 0]\hsy\big).
\]
But if 
\[
F_\zeta 
\hsx \in \hsx  
\ent_+ \big(\hsy]-\infty, 0]\hsy\big),
\]
then
\[
F_\zeta (z)
\ = \ 
F_\zeta (0) \ 
\prod\limits_{n = 1}^\infty \ 
\Big(
1 - \frac{z}{\lambda_n}
\Big)
\]
and, as there is no exponential term, in view of 40.15, 
\[
\forall \ n, \ \forall  \ r, \ 
D (n, r) 
\ > \ 0.
\]
\\[-1.5cm]
\end{x}

\begin{x}{\small\bf THEOREM} \ 
If 
\[
\forall \ n, \ \forall  \ r, \ 
D (n, r) 
\ > \ 0, 
\]
then RH obtains.
\\[-.5cm]

PROOF \ 
The assumption implies that the sequence 
$C_0, C_1, C_2, \ldots$ is totally positive (cf. 40.5), hence
\[ 
F_\zeta 
\hsx \in \hsx  
\ent_+ \big(\hsy]-\infty, 0]\hsy\big)
\qquad (\tcf. \ 40.11),
\]
from which RH.
\\[-.25cm]
\end{x}

\begin{x}{\small\bf SCHOLIUM} \ 
RH is equivalent to the statement that 
\[
\forall \ n, \ \forall  \ r, \ 
D (n, r) 
\ > \ 0.
\]
\\[-1.25cm]
\end{x}

\qquad
{\small\bf \un{N.B.}} \ 
Trivially, 
\[
D (n, 1) 
\ = \ 
C_n 
\ > \ 
0.
\]
\\[-.25cm]


\chapter{
$\boldsymbol{\S}$\textbf{42}.\quad  $D (n, 2)$}
\setlength\parindent{2em}
\setcounter{theoremn}{0}
\renewcommand{\thepage}{\S42-\arabic{page}}


\qquad Here it will be shown that $D (n, 2)$ is positive (cf. 41.8).
\\[-.25cm]

\qquad
{\small\bf \un{N.B.}} \ 
We have
\[
D (0, 2)
\quad = \quad
\begin{vmatrix*}[l]
\ 
C_0  
&\ 
\ 0 
\\[4pt]
\ C_1
&\ 
\ 
C_0 \text{\hsy} 
\end{vmatrix*}
\quad = \quad
C_0^2 
\ > \ 
0,
\]
so it can be assumed that $n \geq 1$.
\\[-.25cm]

\begin{x}{\small\bf LEMMA}
\footnote[2]{\vspace{.11 cm}
G. Csordas and R. Varga, 
\textit{Constr. Approx.} \textbf{4} (1988), pp. 175-198.}
 \ 
$\forall \ t > 0$, 
\[
\frac{\td}{\td t} \hsx \Big(\frac{\Phi^\prime  (t)}{t \hsx \Phi (t)}\Big)
\ < \ 
0.
\]
\\[-1.25cm]
\end{x}

\begin{x}{\small\bf THEOREM} \ 
$\forall \ n \geq 1$,
\[
C_n^2 \hsx - \hsx \Big(1 + \frac{1}{n}\Big) 
\hsx 
C_{n-1} \hsy C_{n+1}
\ \geq \ 
0.
\]

PROOF \ 
Write
\allowdisplaybreaks
\begin{align*}
C_n^2 \hsx - \hsx \Big(1 + \frac{1}{n}\Big) &C_{n-1} \hsy C_{n+1}\ 
\\[15pt]
&= \
\frac{b_n^2}{(2 n !)^2}
\hsx - \hsx
\frac{n+1}{n}
\hsx
\frac{1}{(2n - 2) !}
\hsx
\frac{1}{(2n + 2) !}
\
b_{n-1} \hsy b_{n+1}
\\[15pt]
&= \
\frac{1}{(2 n !)^2}
\Big(
b_n^2 
\hsx - \hsx
\frac{n+1}{n}
\hsx
\frac{(2 n )!}{(2n - 2) !}
\hsx
\frac{(2 n )!}{(2n + 2) !}
\
b_{n-1} \hsy b_{n+1}
\Big)
\\[15pt]
&= \
\frac{1}{(2 n !)^2}
\Big(
b_n^2 
\hsx - \hsx
\frac{n+1}{n}
\hsx
\frac{2n (2n - 1)}{1}
\hsx
\frac{1}{2 (n + 1)(2n + 1)}
\
b_{n-1} \hsy b_{n+1} 
\Big)
\\[15pt]
&= \
\frac{1}{(2 n !)^2}
\Big(
b_n^2 
\hsx - \hsx
\frac{2 n - 1}{2n + 1}
\hsx
b_{n-1} \hsy b_{n+1} 
\Big).
\end{align*}
Put
\[
\Delta_n
\ = \ 
b_n^2 
\hsx - \hsx
\frac{2 n - 1}{2n + 1}
\hsx
b_{n-1} \hsy b_{n+1}
\]
and then make the claim that 
$\Delta_n \geq 0$.  
First
\[
b_n 
\ = \ 
\int\limits_0^\infty \ 
t^{2 n} 
\hsy
\Phi (t) 
\ \td t
\]
\qquad \qquad
$\implies$
\[
b_n 
\ = \ 
- 
\frac{1}{2n + 1}
\ 
\int\limits_0^\infty \ 
t^{2 n + 1} 
\hsy
\Phi^\prime (t) 
\ \td t.
\]
Therefore
\allowdisplaybreaks\begin{align*}
&
\int\limits_0^\infty \ 
\int\limits_0^\infty \ 
u^{2 n} 
\hsy
v^{2 n} 
\hsy
\Phi (u) 
\hsy
\Phi (v) 
\hsy (v^2 - u^2)
\ 
\bigg(
\int\limits_u^v \ 
\hsx - \hsx
\frac{\td}{\td t} 
\ 
\Big(\frac{\Phi^\prime  (t)}{t \hsx\Phi (t)}\Big)
\ \td t
\bigg)
\ \td u \hsy \td v
\\[13pt]
&
\hspace{1cm}
=\ 
\int\limits_0^\infty \ 
\int\limits_0^\infty \ 
u^{2 n - 1} 
\hsy
v^{2 n - 1} 
\hsy
(v^2 - u^2)
\hsy
\big(
v \hsy \Phi (v) \Phi^\prime  (u) 
\hsx - \hsx 
u \hsy \Phi (u) \Phi^\prime  (v) 
\big)
\ \td u \hsy \td v
\\[13pt]
&
\hspace{1cm}
=\ 
- 
(2n - 1) b_{n-1}
\ 
\int\limits_0^\infty \ 
v^{2 n + 2} 
\hsy
\Phi (v)
\ \td v
\ + \ 
(2n + 1) \hsy b_n
\ 
\int\limits_0^\infty \ 
v^{2 n} 
\hsy
\Phi (v)
\ \td v
\\[13pt]
&
\hspace{1.5cm}
\ + \ 
(2n + 1) b_n
\ 
\int\limits_0^\infty \ 
u^{2 n} 
\hsy
\Phi (u)
\ \td u
\ - \ 
(2n - 1) \hsy b_{n-1}
\ 
\int\limits_0^\infty \ 
u^{2 n + 2} 
\hsy
\Phi (u)
\ \td u
\\[13pt]
&
\hspace{1cm}
=\ 
- 
(2n - 1)\hsy  b_{n-1} \hsy b_{n+1}
\hsx + \hsx 
(2n + 1) \hsy b_n^2 
\hsx + \hsx 
(2n + 1) \hsy b_n^2 
\hsx - \hsx 
(2n - 1)\hsy  b_{n-1} \hsy b_{n+1}
\\[13pt]
&
\hspace{1cm}
=\ 
2 \hsy (2n + 1) \hsy b_n^2
\hsx - \hsx 
2  (2n - 1)\hsy  b_{n-1} \hsy b_{n+1}
\\[13pt]
&
\hspace{1cm}
=\ 
2 \hsy (2n + 1) \hsy
\Big(
b_n^2 \hsx - \hsx \frac{2 (2n - 1)}{2 (2n + 1)}
\hsx
 b_{n-1} \hsy b_{n+1}
\Big)
\\[13pt]
&
\hspace{1cm}
=\ 
2 \hsy (2n + 1) \hsy \Delta_n.
\end{align*}
But $\forall \ t > 0$, 
\[
-
\frac{\td}{\td t} \hsx \Big(\frac{\Phi^\prime  (t)}{t \hsx \Phi (t)}\Big)
\ > \ 
0
\qquad (\tcf. \ 41.9).
\]
Consequently, 
\[
(v^2 - u^2) 
\ 
\bigg(
\int\limits_u^v \ 
-
\frac{\td}{\td t} \hsx \Big(\frac{\Phi^\prime  (t)}{t \hsx \Phi (t)}\Big)
\ \td t
\bigg)
\ \td u \hsy \td v
\]
is nonnegative for all $0 \leq u$, $v < \infty$, hence 
$\Delta_n$ is  $\geq 0$, 
as claimed.
\\[-.25cm]
\end{x}

\begin{x}{\small\bf APPLICATION} \ 
$\forall \ n \geq 1$,
\[
C_n^2
\ \geq 
\Big(1 + \frac{1}{n}\Big) C_{n-1} \hsy C_{n+1}
\ > \ 
C_{n-1} \hsy C_{n+1}
\]
\qquad
\qquad 
$\implies$
\[
C_n^2
\ > \ 
C_{n-1} \hsy C_{n+1}
\]
\qquad
\qquad  
$\implies$
\allowdisplaybreaks
\begin{align*}
D (n, 2) \ 
&=\ 
\begin{vmatrix*}[l]
\ C_n 
&\ 
\ C_{n-1} \text{\hsx} 
\\[4pt]
\ C_{n+1}
&\ 
\ C_n  
\end{vmatrix*}
\\[15pt]
&=\ 
C_n^2 - C_{n-1} \hsy C_{n+1}
\\[15pt]
&> \
0.
\end{align*}
\\[-1.5cm]
\end{x}

\begin{x}{\small\bf REMARK} \ 
Put
\[
\Gamma_n 
\ = \ 
F_\zeta^{(n)} (0) 
\quad 
\Big(\implies C_n = \frac{\Gamma_n }{n !}\Big).
\]
Then
\[
\Gamma_n^2 \hsx - \hsx \Gamma_{n-1} \hsy \Gamma_{n + 1} 
\ \geq \ 
0.
\]
I.e.: 
\[
\big(F_\zeta^{(n)} (0) \big)^2 
\hsx - \hsx
F_\zeta^{(n-1)} (0) 
\hsx
F_\zeta^{(n+1)} (0) 
\ \geq \ 
0.
\]
Take now $n = 1$ and, in the notation of 13.6, ask: \ 
Is it true that for ALL real $t$, 
\[
L_1 (F_\zeta) (t) 
\ = \ 
\big(F_\zeta^\prime (t) \big)^2 
\hsx - \hsx 
F_\zeta (t) \hsy F_\zeta^{\prime\prime} (t) 
\ \geq \ 
0\hsy ?
\]
The answer is unknown (although the inequality does hold in  a finite interval containing the origin \ldots).

[Note: \ 
If $\forall \ t$, 
\[
L_1 (F_\zeta) (t) 
\ > \ 
0,
\]
then it would follow that all the real zeros of $F_\zeta$ are simple.]
\\[-.5cm]
\end{x}

There is another proof of the positivity of $D (n,2)$ that is based on a different set of ideas, 
these being important for their associated methodology.
\\[-.25cm]

\begin{x}{\small\bf LEMMA} \ 
$\forall \ t > 0$, 
\[
-
\
\begin{vmatrix*}
\ \Phi (t)  
&\ 
\ \Phi^\prime  (t) \ 
\\[4pt]
\ \Phi^\prime (t) 
&\ 
\ \Phi^{\prime\prime}  (t)   \ 
\end{vmatrix*}
\quad > \quad
0.
\]

PROOF \ 
Owing to 42.1, 
$\forall \ t > 0$, 
\[
\frac{\td}{\td t} \hsx \Big(\frac{\Phi^\prime  (t)}{t \hsx \Phi (t)}\Big)
\ < \ 
0
\]
which, when written out, is equivalent to the inequality
\[
t \hsy ((\Phi^\prime  (t))^2 
\hsx - \hsx 
\Phi (t) \hsy \Phi^{\prime\prime}  (t) )
\hsx + \hsx
\Phi (t) \hsy \Phi^\prime  (t)
\ > \ 
0
\]
or still, 
\[
t \hsy ((\Phi^\prime  (t))^2 
\hsx - \hsx 
\Phi (t) \hsy \Phi^{\prime\prime}  (t) )
\ > 
-
\Phi (t) \hsy \Phi^\prime  (t).
\]
But $\Phi (t)$ is positive (cf. 38.15) and $\Phi^\prime  (t)$ is negative (cf. 38.19).  
Therefore 
\[
-
\Phi (t)
\hsy
\Phi^\prime  (t)
\ > \ 
0
\]
\qquad
\qquad 
$\implies$
\allowdisplaybreaks
\begin{align*}
(\Phi^\prime  (t))^2 
\hsx - \hsx 
\Phi (t) \hsy \Phi^{\prime\prime}  (t) 
\quad
&= \quad
-
\
\begin{vmatrix*}
\ \Phi (t)  
&\ 
\ \Phi^\prime  (t) \ 
\\[4pt]
\ \Phi^\prime (t) 
&\ 
\ \Phi^{\prime\prime}  (t)   \ 
\end{vmatrix*}
\\[15pt]
\quad 
&> \quad
0.
\end{align*}

[Note: \ 
\allowdisplaybreaks
\begin{align*}
\frac{\td^2}{\td t^2} \hsx \log \Phi (t) \ 
&=\ 
\frac{\td}{\td t} \hsx \Big(\frac{\Phi^\prime  (t)}{\Phi (t)}\Big)
\\[15pt]
&=\ 
\frac{\Phi (t) \hsy \Phi^{\prime\prime}  (t)  - (\Phi^\prime  (t))^2 }{\Phi (t)^2}
\\[15pt]
&<\ 
0.]
\end{align*}
\\[-1.25cm]
\end{x}

\qquad
{\small\bf \un{N.B.}} \ 
It is to be emphasized that it is possible to give a proof of 42.5 which is independent of 42.1 
(see the Appendix to this \S).]
\\[-.5cm]

[Note: \ 
It is shown there that the inequality persists to $t = 0$ (or directly: 
\[
\big((\Phi^\prime  (t))^2 \hsx - \hsx \Phi (t) \hsy \Phi^{\prime\prime}  (t)\big)
\hsx 
\Big|_{t = 0}
\quad = \quad 
0^2 - \Phi (0) \hsy \Phi^{\prime\prime} (0)
\ > \ 
0, 
 \]
 $\Phi (0)$ being positive and 
 $\Phi^{\prime\prime}  (0)$ being negative.]
\\[-.25cm]

\begin{x}{\small\bf SUBLEMMA} \ 
Let 
$f_1 (t)$, 
$f_2 (t)$, 
$g_1 (t)$, 
$g_2 (t)$ 
be continuous and absolutely integrable on 
$[0, \infty[\hsy$. 
Assume: \ 
$f_i (t) \hsy g_j (t)$ $(1 \leq i, \hsy j \leq 2)$ 
and 
$f_1 (t) \hsy f_2 (t) \hsy g_1 (t) \hsy g_2 (t)$ 
are also absolutely integrable on 
$[0, \infty[$ \hsy $-$then

\[
\det
\begin{bmatrix*}
\ \ds
\int\limits_0^\infty \ 
f_1 (t) \hsy g_1 (t)
\ \td t
&
\ds
\int\limits_0^\infty \ 
f_1 (t) \hsy g_2 (t)
\ \td t
\\[26pt]
\ \ds
\int\limits_0^\infty \ 
f_2 (t) \hsy g_1 (t)
\ \td t
&
\ds
\int\limits_0^\infty \ 
f_2 (t) \hsy g_2 (t)
\ \td t
\end{bmatrix*}
\hspace{4.4cm}
\]

\[
\hspace{1.8cm}
= 
\quad 
\iint\limits_{0 < u < v < \infty} \ 
\det \ 
\begin{bmatrix*}
\ \ds
f_1 (u) 
&
f_1 (v) \ 
\\[4pt]
\ \ds
f_2 (u) 
&
f_2 (v) \ 
\end{bmatrix*}
\ \bcdot \ 
\det \ 
\begin{bmatrix*}
\ \ds
g_1 (u) 
&
g_1 (v) \ 
\\[4pt]
\ \ds
g_2 (u) 
&
g_2 (v) \ 
\end{bmatrix*}
\ \td u \hsy \td v.
\]
\\[-.5cm]
\end{x}

\begin{x}{\small\bf NOTATION} \ 
Given nonempty subsets $X$ and $Y$ of $\R$ and a real valued function 
$f$ on $X \times Y$, put
\\
\[
f \ 
\begin{bmatrix*}
\ 
x_1
&
\ x_2 \ 
\\[4pt]
\ 
y_1
&
\ y_2 \ 
\end{bmatrix*}
\quad = \quad 
\det
\begin{bmatrix*}
\ 
f(x_1, y_1)
&
\ f(x_1, y_2) \ 
\\[4pt]
\ 
f(x_2, y_1)
&
\ f(x_2, y_2) \ 
\end{bmatrix*}
.
\]
\\[-.5cm]

\noindent
Put

\[
\phi (v, t) 
\ = \ 
\frac{v^{t-1}}{\Gamma (t)}
\qquad (v > 0, t > 0).
\]
\\[-1cm]
\end{x}

\begin{x}{\small\bf LEMMA} \ 
$\forall \ t > 0$,
$\forall \ s> 0$,
\[
\phi (v, t + s) 
\ = \ 
\int\limits_0^v \ 
\phi (u, t)
\hsy 
\phi (v-u, s)
\ \td u.
\]

PROOF \ 
Start with the RHS:
\allowdisplaybreaks\begin{align*}
&
\int\limits_0^v \ 
\frac{u^{t-1}}{\Gamma (t)}
\hsx
\frac{(v - u)^{s-1}}{\Gamma (s)}
\ \td u
\\[15pt]
&
\hspace{1cm}
=\ 
\frac{1}{\Gamma (t)}
\hsx
\frac{1}{\Gamma (s)}
\
\int\limits_0^v \ 
u^{t-1} 
\hsy (v - u)^{s - 1}
\ \td u
\\[15pt]
&
\hspace{1cm}
=\ 
\frac{1}{\Gamma (t)}
\hsx
\frac{1}{\Gamma (s)}
\hsz
v^{s-1}
\
\int\limits_0^v \ 
u^{t-1} 
\hsy
\Big(1 - \frac{u}{v}\Big)^{s-1}
\ \td u
\\[15pt]
&
\hspace{1cm}
=\ 
\frac{1}{\Gamma (t)}
\hsx
\frac{1}{\Gamma (s)}
\hsz
v^{s-1}
\
\int\limits_0^1\ 
(v \hsy w)^{t-1} 
\hsy
(1 - w)^{s-1}
\hsy 
v 
\ \td w
\\[15pt]
&
\hspace{1cm}
=\ 
\frac{1}{\Gamma (t)}
\hsx
\frac{1}{\Gamma (s)}
\hsz
v^{t+s-1}
\hsx
B (t, s)
\\[15pt]
&
\hspace{1cm}
=\ 
\frac{1}{\Gamma (t)}
\hsx
\frac{1}{\Gamma (s)}
\hsz
v^{t+s-1}
\hsx
\frac{\Gamma (t) \hsy \Gamma (s)}{\Gamma (t+s)}
\\[15pt]
&
\hspace{1cm}
=\ 
\frac{v^{t+s-1}}{\Gamma (t+s)}
\\[15pt]
&
\hspace{1cm}
=\
\phi(v, t+s).
\end{align*}

Put
\[
\lambda (t) 
\ = \ 
\int\limits_0^\infty \ 
\Phi (v) 
\hsy 
\phi (v, t) 
\ \td v
\qquad (t > 0).
\]
Then
\allowdisplaybreaks
\begin{align*}
\lambda (2n + 1) \ 
&=\ 
\int\limits_0^\infty \ 
\Phi (v) 
\hsy 
\phi (v, 2n + 1) 
\ \td v
\\[15pt]
&= \ 
\int\limits_0^\infty \ 
\Phi (v) 
\ 
\frac{v^{2n + 1 - 1}}{\Gamma(2n + 1)}
\ \td v
\\[15pt]
&= \ 
\int\limits_0^\infty \ 
\Phi (v) 
\ 
\frac{v^{2 n}}{(2 n)!}
\ \td v
\\[15pt]
&= \ 
\frac{1}{(2 n)!}
\int\limits_0^\infty \ 
\Phi (v) 
\ 
v^{2 n}
\ \td v
\\[15pt]
&= \ 
\frac{b_n}{(2 n)!}
\\[15pt]
&= \ 
C_n.
\end{align*}
\\[-1.25cm]
\end{x}

\begin{x}{\small\bf LEMMA} \ 
$\forall \ t > 0$,
$\forall \ s> 0$,
\allowdisplaybreaks\begin{align*}
\Lambda (s, t) \ 
&\equiv \ 
\lambda (s + t)
\\[15pt]
&= \ 
\int\limits_0^\infty \ 
\Phi (v) 
\hsy
\phi (v, s + t) 
\ \td v
\\[15pt]
&= \ 
\int\limits_0^\infty \ 
\phi (u, s)
\bigg(
\int\limits_0^\infty \ 
\Phi (u + v) 
\hsy 
\phi (v, t) 
\ \td v
\bigg)
\ \td u.
\end{align*}

PROOF \ 
In the double integral, let

\[
\begin{cases}
\ 
x \ = \ u
\\[4pt]
\ 
y \ = \ u + v
\end{cases}
.
\]
Then the Jacobian equals 1, so there is no $\tJ (x,y)$ factor and since $u$ and $v$ are nonnegative, 
if $x$ is varied first, it goes from 0 to $y$.  
This said, upon inverting, thus

\[
\begin{cases}
\ 
u \ = \ x
\\[4pt]
\ 
v \ = \ y - x
\end{cases}
,
\]
we arrive at

\[
\int\limits_{y = 0}^\infty \ 
\int\limits_{x = 0}^y \ 
\phi (x, s) 
\hsy
\phi (y - x, t) 
\hsy
\Phi (y) 
\ \td x 
\hsx 
\td y
\]
or still, 

\[
\int\limits_{y = 0}^\infty \ 
\Phi (y) 
\bigg(
\int\limits_{x = 0}^y \ 
\phi (x, s) 
\hsy
\phi (y - x, t) 
\ \td x 
\bigg)
\ \td y
\]
or still, 

\[
\int\limits_{y = 0}^\infty \ 
\Phi (y) 
\hsy
\phi (y, s + t) 
\ \td y
\qquad (\tcf. \ 42.8)
\]
or still, 

\[
\int\limits_{y = 0}^\infty \ 
\Phi (v) 
\hsy
\phi (v, s + t) 
\ \td v.
\]
\\[-1.25cm]
\end{x}

\begin{x}{\small\bf LEMMA} \ 
If 
$0 < v_1 < v_2$ 
and if 
$0 < t_1 < t_2$, 
then 

\[
\phi 
\begin{bmatrix*}
\ \ds
v_1
&
v_2
\\[4pt]
\ \ds
t_1
&
t_2
\end{bmatrix*}
\ > \ 
0.
\]

PROOF \ 
In fact, 
\allowdisplaybreaks\begin{align*}
\det
\begin{bmatrix*}
\ \ds
\phi (v_1, t_1)
&
\phi (v_1, t_2)
\\[4pt]
\ \ds
\phi (v_2, t_1)
&
\phi (v_2, t_2)
\end{bmatrix*}
&=\ 
\phi(v_1, t_1) \hsy \phi (v_2, t_2) \hsx - \hsx \phi (v_1, t_2) \hsy \phi (v_2, t_2)
\\[26pt]
&=\ 
\frac{v_1^{t_1}}{v_1 \hsy \Gamma(t_1)}
\ 
\frac{v_2^{t_2}}{v_2 \hsy \Gamma(t_2)}
\ - \ 
\frac{v_1^{t_2}}{v_1 \hsy \Gamma(t_2)}
\ 
\frac{v_2^{t_1}}{v_2 \hsy \Gamma(t_1)}
\\[26pt]
&=\ 
\frac{1}{\Gamma(t_1) \hsy \Gamma(t_2)}
\ 
\bigg[
\frac{v_1^{t_1} \hsy v_2^{t_2}}{v_1 \hsy v_2}
\ - \ 
\frac{v_1^{t_2} \hsy v_2^{t_1}}{v_1 \hsy v_2}
\bigg]
\\[26pt]
&=\ 
\frac{1}{\Gamma(t_1) \hsy \Gamma(t_2)}
\ 
\bigg[
v_1^{t_1 - 1} \hsy v_2^{t_1 - 1 + t_2 - t_1}
\ - \ 
v_1^{t_1 - 1 + t_2 - t_1} \hsy v_2^{t_1 - 1}
\bigg]
\\[26pt]
&=\ 
\frac{v_1^{t_1 - 1}\hsy v_2^{t_1 - 1}}{\Gamma(t_1) \hsy \Gamma(t_2)} 
\ 
\bigg[
v_2^{ t_2 - t_1}
\ - \ 
v_1^{t_2 - t_1} 
\bigg]
\\[26pt]
&>\ 
0.
\end{align*}
\\[-1.25cm]
\end{x}


\begin{x}{\small\bf SUBLEMMA} \ 
Let $I$ be an open interval (bounded or unbounded).  
Suppose that $f$ is twice continuously differentiable on $I$ and 

\[
\frac{\td^2}{\td t^2} \hsx f (t)
\ < \ 
0
\qquad (t \in I).
\]
Then for any four points $a, b, c, d$ in $I$ with 
$a < c < d < b$, 

\[
\frac{f (c) - f (a)}{c - a}
\ > \ 
\frac{f (b) - f (d)}{b - d}.
\]

PROOF \ 
By the mean value theorem, 

\[
\begin{cases}
\ \ds
\frac{f (c) - f (a)}{c - a}
\ = \ 
f^\prime (x) 
\qquad (\exists \  x \in \ ]a, c[\hsx)
\\[11pt]
\ \ds
\frac{f (b) - f (d)}{b - d}
\ = \ 
f^\prime (y) 
\qquad (\exists \ y \in \ ]d ,b[\hsx)
\end{cases}
.
\]
But the assumption on $f$ implies that $f^\prime$ is strictly decreasing on $I$, 
hence 
\[
x 
\ < \ 
y 
\implies 
f^\prime (x) 
\ > \ 
f^\prime (y).
\]

[Note: \ 
If $c - a = b - d$, then 

\[
f (c) + f (d) 
\ > \ 
f (a) + f (b).]
\]
\\[-1.5cm]
\end{x}

\qquad
{\small\bf \un{N.B.}} \ 
In the applications (as below), it can happen that during the course of a ``labeling procedure'', 
one has ``$c = d$'', so
\[
\begin{cases}
\ \ds
\frac{f (c) - f (a)}{c - a}
\ = \ 
f^\prime (x) 
\qquad (\exists \ x \in \ ]a, c[\hsx)
\\[11pt]
\ \ds
\frac{f (b) - f (c)}{b - c}
\ = \ 
f^\prime (y) 
\qquad (\exists\  y \in \ ]c, b[\hsx)
\end{cases}
,
\]
thus if $c - a = b - c$, then 
\[
f (c) + f (c) 
\ > \ 
f (a) + f (b).]
\]
\\[-1cm]
Put
\[
K (u, v) 
\ = \ 
\Phi (u + v) 
\qquad (u > 0, v > 0).
\]
\\[-1.5cm]

\begin{x}{\small\bf LEMMA} \ 
If 
$0 < u_1 < u_2$ 
and if 
$0 < v_1 < v_2$, 
then 
\[
K \  
\begin{bmatrix*}
\ \ds
u_1
&
u_2
\\[4pt]
\ \ds
v_1
&
v_2
\end{bmatrix*}
\ < \ 
0.
\]

PROOF \ 
In 42.11, take 

\[
f (t) 
\ = \ 
\log \Phi (t) 
\qquad (\tcf. \ 42.5).
\]
Define $a, b, c, d$, as follows: 
\hspace{1.5cm}
$
\begin{cases}
a \ = \ u_1 + v_1
\\
b \ = \ u_2 + v_2
\\
c \ = \ u_2 + v_1
\\
d \ = \ u_1 + v_2
\end{cases}
$
.

\noindent
Therefore
\hspace{.75cm}
$
\begin{cases}
a \hsx < \hsx c \hsx < \hsx  b
\\
a \hsx < \hsx d \hsx < \hsx b
\end{cases}\,
$
and 
\quad
$c - a \hsx = \hsx b - d$. 
\\

\noindent
Now, while the setup in 42.11 called for $c < d$, if $d < c$, then their roles can be interchanged and the possibility that 
$c = d$ is not excluded (cf. supra).  
Consequently, 
\[
\log \Phi (c) \hsx + \hsx \log \Phi (d)  
\ > \ 
\log \Phi (a) \hsx + \hsx \log \Phi (b)  
\]
\qquad 
$\implies$
\[
\Phi (c) \hsy \Phi (d)
\ > \ 
\Phi (a) \hsy \Phi (b)
\]
\qquad 
$\implies$
\[
\Phi (u_2 + v_1) \hsy \Phi (u_1 + v_2)
\ > \ 
\Phi (u_1 + v_1) \hsy \Phi (u_2 + v_2)
\]
or still, 
\[
\Phi (u_1 + v_1) \hsy \Phi (u_2 + v_2)
\hsx - \hsx 
\Phi (u_1 + v_2) \hsy \Phi (u_2 + v_1)
\ < \ 
0.
\]
And
\allowdisplaybreaks\begin{align*}
K\  
\begin{bmatrix*}
\ \ds
u_1
&
u_2
\\[4pt]
\ \ds
v_1
&
v_2
\end{bmatrix*}
&=\ 
\det
\begin{bmatrix*}
\ \ds
K (u_1, v_1)
&
K (u_1, v_2)
\\[4pt]
\ \ds
K (u_2, v_1)
&
K (u_2, v_2)
\end{bmatrix*}
\\[26pt]
&=\ 
\det
\begin{bmatrix*}
\ \ds
\Phi (u_1 + v_1)
&
\Phi (u_1 + v_2)
\\[4pt]
\ \ds
\Phi (u_2 + v_1)
&
 \hsy \Phi (u_2 + v_2)
\end{bmatrix*}
\\[11pt]
&<\ 
0.
\end{align*}

Put
\[
L (u, t) 
\ = \ 
\int\limits_0^\infty \ 
K (u, v) 
\hsy 
\phi (v, t) 
\ \td v.
\]
\\[-.75cm]
\end{x}

\begin{x}{\small\bf LEMMA} \ 
If 
$0 < u_1 < u_2$ 
and if 
$0 < t_1 < t_2$, 
then 

\[
L\  
\begin{bmatrix*}
\ \ds
u_1
&
u_2
\\[4pt]
\ \ds
t_1
&
t_2
\end{bmatrix*}
\ < \ 
0.
\]

PROOF \ 
Using 42.6, write
\[
L\  
\begin{bmatrix*}
\ \ds
u_1
&
u_2
\\[4pt]
\ \ds
t_1
&
t_2
\end{bmatrix*}
\quad = \quad 
\iint\limits_{0 \hsy < \hsy u \hsy < \hsy v \hsy < \hsy \infty} \ K
\
\begin{bmatrix*}
\ u_1
&
u_2
\\[4pt]
\ u \
&
v
\end{bmatrix*}
\ \phi \ 
\begin{bmatrix*}
\ u
&
v
\\[4pt]
\ t_1 \
&
t_2
\end{bmatrix*}
\ \td u \hsy \td v.
\]
In this connection, it is necessary to observe that
\\[-.75cm]

\allowdisplaybreaks\begin{align*}
\det
\begin{bmatrix*}
\ \ds
\phi (u, t_1)
&
\phi (v, t_1)
\\[4pt]
\ \ds
\phi (u, t_2)
&
\phi (v, t_2)
\end{bmatrix*}
&=\ 
\det
\begin{bmatrix*}
\ \ds
\phi (u, t_1)
&
\phi (u, t_2)
\\[4pt]
\ \ds
\phi (v, t_1)
&
\phi (v, t_2)
\end{bmatrix*}
\\[26pt]
&=\ 
\phi \
\begin{bmatrix*}
\ u
&
v \ 
\\[4pt]
\ t_1
&
t_2
\end{bmatrix*}
.
\end{align*}
But

\[
K \
\begin{bmatrix*}
\ u_1
&
u_2 \ 
\\[4pt]
\ u
&
v
\end{bmatrix*}
\ < \  
0
\qquad (\tcf. \ 42.12)
\]
and

\[
\phi \
\begin{bmatrix*}
\ u
&
v \ 
\\[4pt]
\ t_1
&
t_2
\end{bmatrix*}
\ > \  
0
\qquad (\tcf. \ 42.10).
\]

Therefore

\[
L \
\begin{bmatrix*}
\ u_1
&
u_2 \ 
\\[4pt]
\ t_1
&
t_2
\end{bmatrix*}
\ < \
0.
\]
\\[-1.25cm]
\end{x}

Using the notation of 42.9, we have 
\allowdisplaybreaks\begin{align*}
\Lambda (s, t) \ 
&\equiv \ 
\lambda (s + t) 
\\[15pt]
&=\ 
\int\limits_0^\infty \ 
\phi (u, s) 
\ 
\bigg(
\int\limits_0^\infty \ 
\Phi (u + v)
\hsy
\phi (v, t)
\ \td v
\bigg)
\ \td u
\\[15pt]
&=\ 
\int\limits_0^\infty \ 
\phi (u, s) 
\ 
\bigg(
\int\limits_0^\infty \ 
K (u, v)
\hsy
\phi (v, t)
\ \td v
\bigg)
\ \td u
\\[15pt]
&=\ 
\int\limits_0^\infty \ 
\phi (u, s) 
\hsy
L(u, t)
\ \td u.
\end{align*}
\\[-1cm]

\begin{x}{\small\bf LEMMA} \ 
If 
$0 < s_1 < s_2$ 
and if 
$0 < t_1 < t_2$, then
\[
\Lambda 
\
\begin{vmatrix*}[l]
\ 
s_1
&
s_2
\\[4pt]
\ 
t_1
&
t_2 \text{\hsx}
\end{vmatrix*}
\quad < \quad
0.
\]

PROOF \ 
Appealing once again to 42.6, write

\[
\Lambda \ 
\begin{bmatrix*}[l]
\ s_1
&
s_2 \ 
\\[4pt]
\ t_1
&
t_2
\end{bmatrix*}
\quad = \quad 
\iint\limits_{0 \hsy < \hsy u \hsy < \hsy v \hsy < \hsy \infty} \ \phi
\hsx
\begin{bmatrix*}[l]
\ u
&
v
\\[4pt]
\ s_1 \
&
s_2
\end{bmatrix*}
\ L \hsx 
\begin{bmatrix*}[l]
\ u
&
v
\\[4pt]
\ t_1 \
&
t_2
\end{bmatrix*}
\ \td u \hsy \td v
\]
\\

\noindent
and then apply 42.10 and 42.13.
\\[-.25cm]
\end{x}


\begin{x}{\small\bf SCHOLIUM} \ 
If 
$0 < s_1 < s_2$ 
and if 
$0 < t_1 < t_2$, then
\[
\begin{vmatrix*}[l]
\ 
\lambda (s_1 + t_1)  
&\ 
\lambda (s_1 + t_2) \ 
\\[8pt]
\ 
\lambda (s_2 + t_1)
&\ 
\lambda (s_2 + t_2)\text{\hsy}
\end{vmatrix*}
\quad < \quad
0.
\]

Consider now the determinant
\[
\begin{vmatrix*}[l]
\ 
C_{n-1} 
&\ 
C_n \
\\[4pt]
\ 
C_n
&\ 
C_{n+1} \text{\hsy} 
\end{vmatrix*}
\qquad (n \geq 1),
\]
hence

\[
C_{n-1} 
\ = \ 
\lambda \hsy (2n - 1), 
\quad 
C_n
\ = \ 
\lambda \hsy (2n+ 1), 
\quad
C_{n+1} 
\ = \ 
\lambda \hsy (2n +3).
\]
In 42.15, let
\hspace{0.25cm}
$
\begin{cases}
\ \ds
s_1 \hsx = \hsx t_1 \hsx = \hsx n - \frac{1}{2}
\\[15pt]
\ \ds
s_2 \hsx = \hsx t_2 \hsx = \hsx n + \frac{3}{2}
\end{cases}
.
$
\hspace{.5cm}
%
Then
\hspace{0.5cm}
$
\begin{cases}
s_1 + t_1 \ = \ 2n - 1
\\
s_1 + t_2 \ = \ 2n + 1
\\
s_2 + t_1 \ = \ 2n + 1
\\
s_2 + t_2\ = \ 2n +3
\end{cases}
$
.
\\

\noindent
Therefore

\[
\begin{vmatrix*}[l]
\ 
\lambda (2n - 1)  
&\ 
\lambda (2n + 1) \ 
\\[8pt]
\ 
\lambda (2n + 1))
&\ 
\lambda (2n + 3)\text{\hsy}
\end{vmatrix*}
\ < \ 
0.
\]
I.e.: 
\[
\begin{vmatrix*}[l]
\ 
C_{n-1} 
&\ 
C_n \
\\[4pt]
\ 
C_n
&\ 
C_{n+1} \text{\hspace{0.08cm}}
\end{vmatrix*}
\ < \ 
0
\]
\\[-1cm]
\end{x}
or still, 

\[
C_{n-1} \hsy C_{n+1} - C_n^2 
\ < \ 
0
\]
or still, 

\[
D (n, 2) 
\ = \ 
C_n^2 - C_{n-1} \hsy C_{n+1} 
\ > \ 
0.
\]
\\[-1.cm]

\begin{x}{\small\bf REMARK} \ 
The condition

\[
C_n^2 \hsx - \hsx  C_{n-1} \hsy C_{n+1} 
\ > \ 
0
\]
is weaker than the condition 

\[
C_n^2 -\Big(1 + \frac{1}{n}\Big) \hsx  C_{n-1} \hsy C_{n+1} 
\ \geq \ 
0
\]
and this is because less was used in its derivation (viz. 42.5 as opposed to 42.1).
\\[-.25cm]
\end{x}

A similar but more complicated analysis serves to establish that $D (n, 3)$ is positive 
(for this and addtional information, see
Nuttall\footnote[2]{\vspace{.11 cm}
\url{https://arxiv.org/pdf/1111.1128}\hsy; 
also 
\textit{Constr. Approx.} \textbf{38} (2013), pp. 193-212.}
).

\begin{center}
APPENDIX
\end{center}

\qquad
{\small\bf THEOREM} \ 
$\forall \ t \geq 0$, 
\[
(\Phi^\prime (t))^2 
\hsx - \hsx 
\Phi (t) \hsy \Phi^{\prime\prime} (t)
\ > \ 
0.
\]
\\[-1cm]

We shall proceed via a list of lemmas.
\\[-.25cm]

Write
\[
\Phi (t) 
\ = \ 
\sum\limits_{n = 1}^\infty \ 
a_n (t),
\]
where
\[
\scalebox{1.04}
{$
a_n (t)
\ = \ 
\big(2 \hsy \pi^2 \hsy n^4 \hsy e^{9 t} \hsx - \hsx  3 \hsy \pi \hsy n^2 \hsy e^{5 t}\big) \hsy\exp (- \pi n^2 e^{4 t}),
$}
\]
and put
\[
a (t) \ = \ a_1 (t), 
\quad
\Psi (t) 
\ = \ 
\sum\limits_{n = 2}^\infty \ 
a_n (t),
\]
thus
\[
\Phi (t) 
\ = \ 
a (t) \hsx + \hsx \Psi (t) 
\]
and so
\allowdisplaybreaks\begin{align*}
(\Phi^\prime (t))^2 
\hsx - \hsx 
\Phi (t) \hsy \Phi^{\prime\prime} (t)\ 
&=\ 
(a^\prime (t) \hsx + \hsx  \Psi^\prime (t))^2
\ - \ 
(a (t) + \Psi (t))
\hsy
(a^{\prime\prime} (t) \hsx + \hsx \Psi^{\prime\prime} (t))
\\[11pt]
&=\ 
V (t) \hsx + \hsx  U (t) \hsx + \hsx \big(\Psi^{\prime} (t)\big)^2.
\end{align*}
Here, by definition, 
\[
V (t)
\ = \ 
(a^\prime (t))^2 \hsx - \hsx a(t) \hsy a^{\prime\prime} (t) 
\]
and 
\[
U (t)
\ = \ 
2 \hsy a^\prime (t) \hsy \Psi^\prime (t)
\hsx - \hsx
 a^{\prime\prime} (t) \hsy \Psi (t)
 \hsx - \hsx
 \Phi (t) \hsy \Psi^{\prime\prime} (t).
\]
\\[-1.cm]

\qquad
{\small\bf NOTATION} \ 
Let 

\[
y 
\ = \ 
\pi \hsy e^{4 \hsy t} 
\quad (t \geq 0) 
\ \implies \ 
y \geq \pi.
\]
\\[-1.cm]

\qquad
{\small\bf LEMMA 1} \ 
$\forall \ t \geq 0$, 
\[
0 
\ < \ 
\Psi (t)
\ \leq \ 
64 \hsy e^t \hsy y^2 \hsy e^{-4 \hsy y}.
\]

PROOF \ 
\allowdisplaybreaks\begin{align*}
0 \ 
&< \ 
\Psi (t)
\\[15pt]
&=\ 
\sum\limits_{n = 2}^\infty \ 
\big(
2 \hsy \pi^2 \hsy n^4 \hsy e^{9 t} \hsx - \hsx  3 \hsy \pi \hsy n^2 \hsy e^{5 t}
\big) 
\hsx 
\exp (- \pi \hsy n^2 \hsy e^{4 t})
\\[15pt]
&\leq\ 
2 \hsy e^t 
\ 
\sum\limits_{n = 2}^\infty \ 
n^4 \hsy \pi^2 \hsy e^{8 t}
\hsx 
\exp (- \pi \hsy n^2 \hsy e^{4 t})
\\[15pt]
&=\ 
2 \hsy e^t 
\hsx
\Big(
16 y^2 \hsy e^{- 4 y} 
\ + \ 
\sum\limits_{n = 1}^\infty \ 
y^2 \hsy n^4 \hsy e^{-n^2 \hsy y}
\Big).
\end{align*}
And
\allowdisplaybreaks\begin{align*}
\sum\limits_{n = 3}^\infty \ 
y^2 \hsy n^4 \hsy e^{-n^2 \hsy y} \ 
&\leq\
\int\limits_2^\infty \ 
y^2 \hsy x^4 \hsy e^{- y \hsy x^2}
\ \td x
\\[15pt]
&< \
\int\limits_2^\infty \ 
y^2 \hsy x^5\hsy e^{- t \hsy x^2}
\ \td x
\\[15pt]
&= \
\frac{1}{y} \hsx e^{-4 y} \hsx \big(1 + 4 y + 8 y^2\big)
\\[15pt]
&< \
16 y^2 \hsx e^{-4 y}.
\end{align*}
Therefore
\allowdisplaybreaks\begin{align*}
\Psi (t) \ 
&\leq \
2 \hsy e^t \big(16 y^2 \hsx e^{-4 y} \hsx + \hsx 16 y^2 \hsx e^{-4 y}\big)
\\[11pt]
&= \
64  \hsy e^t  \hsy y^2 \hsx e^{-4 y}.
\end{align*}
\\[-1.25cm]

\qquad
{\small\bf LEMMA 2} \ 
$\forall \ t \geq 0$, 
\[
\abs{\Psi^\prime (t)} 
\ \leq \ 
565 \hsy e^t \hsy y^3 \hsy e^{-4 y}.
\]

PROOF \ 
\[
\abs{\Psi^\prime (t)} 
\ \leq \ 
\Big|
\hsx
\sum\limits_{n = 2}^\infty \ 
\pi \hsy n^2 
\big(
8 \hsy \pi^2 \hsy n^4 \hsy e^{8 t} 
\hsy - \hsy 
30 \hsy \pi \hsy n^2 \hsy e^{4 t} 
\hsy + \hsy 
15
\big)
\hsx
\exp\big(5 \hsy t \hsy - \hsy \pi \hsy n^2 \hsy e^{4 t}\big)
\hsx
\Big|
\]
or still, if $x = e^t$, 
\[
\abs{\Psi^\prime (t)} 
\ = \ 
8 \pi^3 x^5 
\ 
\bigg|
\hsx
\sum\limits_{n = 2}^\infty \ 
n^6 
\Big(
x^8 - \frac{15}{4 \hsy \pi \hsy n^2} \hsx x^4 + \frac{15}{8 \hsy \pi^2 \hsy n^4}
\Big)
\exp( - \pi \hsy n^2 x^4)
\hsx
\bigg|.
\]
To examine 
$
\ 
\ds
\Big|\hsx \sum\limits_{n = 2} \ldots \hsx \Big|
$, 
\ 
first pull out $x^8$: 
\[
x^8 
\ 
\bigg|
\hsx
\sum\limits_{n = 2}^\infty \ 
n^6 
\Big(
1- \frac{15}{4 \hsy \pi \hsy n^2} \hsx \frac{1}{x^4} + \frac{15}{8 \hsy \pi^2 \hsy n^4}\hsx \frac{1}{x^8} 
\Big)
\exp( - \pi \hsy n^2 x^4)
\hsx
\bigg|
\]
and consider
\[
- \frac{15}{4 \hsy \pi \hsy n^2} \frac{1}{x^4} 
\hsx + \hsx  
\frac{15}{8 \hsy \pi^2 \hsy n^4}\hsx \frac{1}{x^8},
\]
\\[-.5cm]

\noindent
which we claim is strictly trapped between -1 and 0.
\\[-.25cm]

\qquad \textbullet \quad
\[
\frac{1}{2 \hsy \pi \hsy n^2} 
\ < \ 
x^4 
\ \implies \ 
\frac{1}{2 \hsy \pi \hsy n^2} \hsx \frac{1}{x^4}
\ < \ 
1
\]
\\[-1.75cm]

\hspace{1cm}  
$\implies$

\[
- 1 
\ + \ 
\frac{1}{2 \hsy \pi \hsy n^2} \hsx \frac{1}{x^4}
\ < \ 
0
\]
\\[-1.75cm]

\hspace{1cm}  
$\implies$

\[
- 15
\ + \ 
\frac{15}{2 \hsy \pi \hsy n^2} \hsx \frac{1}{x^4}
\ < \ 
0
\]
\\[-1.75cm]

\hspace{1cm}  
$\implies$

\[
-
\frac{15}{4 \hsy \pi \hsy n^2} \hsx \frac{1}{x^4}
\ + \ 
\frac{15}{8 \hsy \pi^2 \hsy n^4} \hsx \frac{1}{x^8}
\ < \ 
0.
\]
\\[-1.75cm]

\qquad \textbullet \quad
\[
\frac{4 \hsy \pi \hsy n^2} {15}
\ > \ 
\frac{1}{x^4}
\]
\hspace{2cm} 
$\implies$
\[
\frac{1}{2 \hsy \pi \hsy n^2} \hsx \frac{1}{x^8}
\hsx + \hsx 
\frac{4 \hsy \pi \hsy n^2} {15}
\ > \ 
\frac{1}{x^4}
\]
\hspace{2cm} 
$\implies$
\[
- \frac{1}{x^4}
\hsx + \hsx 
\frac{1}{2 \hsy \pi \hsy n^2} \hsx \frac{1}{x^8}
\ > \ 
-
\frac{4 \hsy \pi \hsy n^2} {15}
\]
\hspace{2cm}  
$\implies$
\[
- \frac{1}{4 \hsy \pi \hsy n^2} \hsx \frac{1}{x^4}
\hsx + \hsx 
\frac{1}{8 \hsy \pi^2 \hsy n^4} \hsx \frac{1}{x^8}
\ > \ 
-
\frac{1} {15}
\]
\hspace{2cm}  
$\implies$
\[
- \frac{15}{4 \hsy \pi \hsy n^2} \hsx \frac{1}{x^4}
\hsx + \hsx 
\frac{15}{8 \hsy \pi^2 \hsy n^4} \hsx \frac{1}{x^8}
\ > \ 
-
1.
\]
\noindent
Accordingly, if 
\[
C_{x, n} 
\ = \ 
- \frac{15}{4 \hsy \pi \hsy n^2} \hsx \frac{1}{x^4}
\hsx + \hsx 
\frac{15}{8 \hsy \pi^2 \hsy n^4} \hsx \frac{1}{x^8}
\hsx , 
\]
then 
\[
-1
\ < \ 
C_{x, n} 
\ < \ 
0
\]
\hspace{1cm}  
$\implies$
\[
0
\ < \ 
1 \hsx + \hsx
C_{x, n} 
\ < \ 
1
\]
\hspace{1cm}  
$\implies$
\[
\abs{1 \hsx + \hsx C_{x, n} }
\ = \ 
1 + C_{x, n}
\ < \ 
1
\]
\hspace{1cm}  
$\implies$
\allowdisplaybreaks\begin{align*}
\bigg|
\hsx
\sum\limits_{n = 2}^\infty \
& 
n^6 
\Big(
1- \frac{15}{4 \hsy \pi \hsy n^2} \frac{1}{x^4} 
\hsx + \hsx
\frac{15}{8 \hsy \pi^2 \hsy n^4}\hsx \frac{1}{x^8} 
\Big)
\exp( - \pi \hsy n^2 x^4)
\hsx
\bigg|
\\[15pt]
&=\ 
\bigg|
\hsx
\sum\limits_{n = 2}^\infty \ 
n^6 
\big(
1 \hsx + \hsx C_{x, n}
\big)
\exp( - \pi \hsy n^2 x^4)
\hsx
\bigg|
\\[15pt]
&\leq \ 
\sum\limits_{n = 2}^\infty \ 
n^6 
\abs{1 \hsx + \hsx C_{x, n}}
\exp( - \pi \hsy n^2 x^4)
\\[15pt]
&< \ 
\sum\limits_{n = 2}^\infty \ 
n^6 \hsx \exp( - \pi \hsy n^2 x^4)
\end{align*}
\hspace{1cm}  
$\implies$
\[
\abs{\Psi^\prime (t)} 
\ < \ 
\frac{8 \hsy y^{13/4}}{\pi^{1/4}}
\ 
\sum\limits_{n = 2}^\infty \ 
n^6 \hsx e^{-n^2 y} \ 
\qquad (y = \pi \hsy x^4 \geq \pi).
\]
And
\allowdisplaybreaks
\begin{align*}
\sum\limits_{n = 2}^\infty \ 
n^6 \hsx e^{-n^2 y} \ 
&<\ 
64 \hsy e^{- 4 y} 
\ + \ 
\int\limits_2^\infty \ 
s^6 \hsx e^{-s^2 y} 
\ \td s
\\[15pt]
&<\ 
64 \hsy e^{- 4 y} 
\ + \ 
\frac{e^{- 4 y}}{2 \hsy y^{7/2}}
\Big(
(4 y)^{5/2}
\ + \ 
\frac{5}{2} \hsy (4 y)^{3/2}
\\[15pt]
&
\hspace{3cm}
\ + \ 
\frac{15}{4} \hsy (4 y)^{1/2}
\ + \ 
\frac{15\hsy e^{4 y} }{8}
\ 
\int\limits_{4 \hsy y}^\infty \ 
\frac{e^{-u}}{\sqrt{u}}
\ \td u
\Big).
\end{align*}
But 
$\ds 
\frac{1}{\sqrt{u}} < 1$ 
for $u \geq 4 y \geq 4 \pi$, hence
\[
e^{4 y} \ 
\int\limits_{4 \hsy y}^\infty \ 
\frac{e^{-u}}{\sqrt{u}}
\ \td u
\ < \ 
1,
\]
so
\[
\sum\limits_{n = 2}^\infty \ 
n^6 \hsx e^{-n^2 y} \ 
\]
is bounded above by 
\[
64 \hsy e^{- 4 y} 
\Big(
1 
\hsx + \hsx 
\frac{1}{4 y}
\hsx + \hsx 
\frac{1}{32 \hsy y^2}
\hsx + \hsx 
\frac{15}{256 \hsy y^3}
\hsx + \hsx 
\frac{15}{1024 \hsy y^{7/2}}
\Big)
\qquad (y \geq \pi).
\]
The expression in parentheses is strictly decreasing, thus is majorized by its value at 
$y = \pi$ 
and it follows that 
\[
\sum\limits_{n = 2}^\infty \ 
n^6 \hsx e^{-n^2 y} \ 
\ < \ 
64 \hsy e^{- 4 y} 
\hsx 
\Big(
1 
\hsx + \hsx 
\frac{13}{40\hsy \pi}
\Big).
\]
Therefore
\allowdisplaybreaks\begin{align*}
\abs{\Psi^\prime (t)} \ 
&< \ 
\frac{8 \hsy y^{13/4}}{\pi^{1/4}}\ 
\Big(
64 \hsy e^{- 4 y} 
\hsx 
\Big(
1 
\hsx + \hsx 
\frac{13}{40\hsy \pi}
\Big)
\Big)
\\[15pt]
&=\ 
512 \hsy 
\hsx 
\Big(
1 
\hsx + \hsx 
\frac{13}{40\hsy \pi}
\Big)
\hsy 
\pi^3 
\hsy 
\exp(13 t - 4 \hsy \pi \hsy e^{4 t})
\\[15pt]
&< \
565 \hsx  \pi^3 \hsy \exp(13 t - 4 \hsy \pi \hsy e^{4 t})
\\[15pt]
&=\
565 \hsx e^t \hsy y^3 \hsy e^{-4 \hsy y}. 
\end{align*}

\qquad
{\small\bf LEMMA 3} \ 
$\forall \ t \geq 0$, 
\[
\abs{\Psi^{\prime\prime} (t)} 
\ \leq \ 
(1.031) \hsy 2^{13} \hsy e^t \hsy y^4 \hsy e^{-4 y}.
\]

PROOF \
Let
\[
p (x) 
\ = \ 
32 x^3 - 224 x^2 + 330 x - 75.
\]
Then $p (x)$ has three distinct positive roots
\[
0 
\hsx < \hsx 
x_1 
\hsx < \hsx 
x_2 
\hsx < \hsx 
x_3 
\ = \ 
5.049720\ldots \ .
\]
Therefore
\[
x \hsx > \hsx x_3 
\implies
p (x) > 0.
\]
On the other hand,
\[
x \hsx > \hsx x_3 
\implies 
0 \hsx < \hsx  p (x) \hsx < \hsx  32 x^3.
\]
These points made, from the definitions, 
\[
\Psi^{\prime\prime} (t)
\ = \ 
\sum\limits_{n = 2}^\infty \ 
\pi \hsy n^2 \hsy 
p (\pi \hsy n^2 \hsy e^{4t})
\hsx
\exp\big(5 t - \pi \hsy n^2 \hsy e^{4t}\big).
\]
But
\[
\pi \hsy n^2 \hsy e^{4t}
\ \geq \ 
4 \hsy \pi 
\ > \ 
x_3
\]
\qquad 
$\implies$
\allowdisplaybreaks\begin{align*}
\abs{\Psi^{\prime\prime} (t)} \ 
&\leq \ 
32
\ 
\sum\limits_{n = 2}^\infty \ 
\pi \hsy n^2 \hsy 
(\pi \hsy n^2 \hsy e^{4t})^3
\hsx
\exp\big(5 t - \pi \hsy n^2 \hsy e^{4t}\big)
\\[15pt]
&=\ 
32 \hsy \pi^4 \hsy e^{17 \hsy t} 
\ 
\sum\limits_{n = 2}^\infty \ 
\frac{n^8}{\exp \big(\pi \hsy n^2 \hsy e^{4t}\big)}
\\[15pt]
&=\ 
32 \hsy \pi^4 \hsy e^{17 \hsy t} 
\ 
\sum\limits_{n = 2}^\infty \ 
\frac{n^8}{\exp \big(n^2 \hsy y\big)}
\\[15pt]
&=\ 
32 \hsy \pi^4 \hsy e^{17 \hsy t} 
\ 
\sum\limits_{n = 2}^\infty \ 
\frac{1}{\exp \big(n^2 \hsy y - 8 \log n\big)}
\\[15pt]
&\leq\ 
32 \hsy \pi^4 \hsy e^{17 \hsy t} 
\ 
\sum\limits_{n = 2}^\infty \ 
\frac{1}{K (y)^n}
\\[15pt]
&=\ 
32 \hsy \pi^4 \hsy e^{17 \hsy t} 
\ 
\frac{1}{K (y)^{\raisebox{.1cm}{\scriptsize 2}} \hsy\Big(\ds 1 - \frac{1}{K (y)}\Big)}
\end{align*}
if
\[
K (y)
\ = \ 
\frac{e^{2 y}}{16}
\]
as then
\[
n^2 \hsy y - 8 \log n
\ \geq \ 
n \log K (y).
\]
But
\allowdisplaybreaks\begin{align*}
\frac{1}{K (y)^{\raisebox{.1cm}{\scriptsize 2}} \hsy\Big(\ds 1 - \frac{1}{K (y)}\Big)^{\textcolor{white}{|}}} \ 
&=\ 
\frac{\ds  2^8 \hsy e^{-4 \hsy y}}{\ds 1 - \frac{16^{\textcolor{white}{|}}}{e^{2 y}}}
\\[15pt]
&\leq \ 
\ds\frac{\ds  2^8 \hsy e^{-4 y}}{\ds 1 - \frac{16^{\textcolor{white}{|}}}{e^{2 \pi}}}
\qquad (y \geq \pi).
\end{align*}
And
\[
\frac{1}{\ds 1 - \frac{16^{\textcolor{white}{|}}}{e^{2 \hsy \pi}}}
\ < \ 
1.031,
\]
leaving 
\[
< \ 
(1.031) \hsy 2^8 \hsy  e^{-4 y}.
\]
Finally
\allowdisplaybreaks
\begin{align*}
\pi^4 \hsy e^{17 t} \ 
&=\ 
e^t \hsy \pi^4 \hsy  e^{16 t}
\\[11pt] 
&=\ 
e^t \hsy y^4.
\end{align*}
\\[-1cm]

\qquad
{\small\bf LEMMA 4} \ 
$\forall \ t \geq 0$, 
\[
0 
\ < \ 
\Phi (t) 
\ < \ 
\frac{203}{202} \hsx a (t).
\]
\\[-1.5cm]

PROOF \
\allowdisplaybreaks\begin{align*}
\Psi (t) \ 
&<\ 
64 \hsy \pi^2 \hsy \exp (9 t - 4 \hsy \pi \hsy e^{4 t})
\\[11pt]
&<\ 
\frac{1}{202} \hsx a (t)
\end{align*}
\qquad\qquad 
$\implies$
\\[-1cm]
\allowdisplaybreaks\begin{align*}
\Phi (t) 
&= \ 
a (t) + \Psi (t)
\\[11pt]
&<\ 
a (t) + \frac{1}{202} \hsx a (t)
\\[11pt]
&=\ 
\frac{203}{202} \hsx a (t).
\end{align*}
\\[-1.25cm]

\qquad
{\small\bf NOTATION} \ 
Put
\[
E (y) 
\ = \ 
e^{2 t} \hsy e^{-2 y} \hsy y^3.
\]
\\[-1cm]

\qquad
{\small\bf LEMMA 5} \ 
$\forall \ t \geq 0$, 
\[
V (t) 
\ \geq \ 
256 \hsy e^{2 t} \hsy e^{-2 y} \hsy y^3
\ 
\equiv \ 
256 \hsy E (y) .	
\]

PROOF \
\allowdisplaybreaks\begin{align*}
V (t) \
&=\ 
16 
\hsy
\exp\big(- 2 \hsy \pi \hsy e^{4 t} \hsy + \hsy  14 \hsy  t\big)
\hsx
\pi^3
\hsx
(15 \hsy - \hsy 12 \hsy \pi \hsy e^{4 t} \hsy + \hsy 4 \hsy \pi^2 \hsy e^{8 t})
\\[11pt]
&=\ 
16 
\hsy
e^{14 t} \hsy e^{-2 y} \hsy \pi^3 \hsy (15 \hsy - \hsy 12 y \hsy + \hsy 4 y^2)
\\[11pt]
&=\ 
16 
\hsy
e^{2 t} \hsy e^{-2 y} \hsy y^3 \hsy (15 \hsy - \hsy 12 y \hsy + \hsy 4 y^2).
\end{align*}
But
\[
15 - 12 y + 4 y^2
\ = \ 
4 \Big(y - \frac{3}{2}\Big)^2 + 6
\]
is an increasing function of $y \geq \pi$, so
\allowdisplaybreaks\begin{align*}
4 \Big(y - \frac{3}{2}\Big)^2 + 6\ 
&\geq \ 
4 \Big(\pi - \frac{3}{2}\Big)^2 + 6
\\[11pt]
&\leq \ 
16.
\end{align*}
Therefore
\[
V (t) 
\ \geq \ 
256 \hsy 
e^{2 t} \hsy e^{-2 y} \hsy y^3 
\ \equiv \ 
256 \hsy E (y).
\]
\\[-1cm]


\qquad
{\small\bf NOTATION} \ 
Write
\[
\begin{cases}
\ 
a (t) 
\ = \ 
e^t \hsy e^{-y} \hsy y \hsy (2 y - 3)
\\[4pt]
\ 
a^\prime (t) 
\ = \ 
- e^t \hsy e^{-y} \hsy y \hsy (15 - 30y + 8 y^2)
\\[4pt]
\ 
a^{\prime\prime} (t)
\ = \ 
e^t \hsy e^{-y} \hsy y \hsy (-75 + 330 y - 224 y^2 + 32 y^3)
\end{cases}
.
\]
\\[-.5cm]

\qquad
{\small\bf LEMMA 6} \ 
$\forall \ t \geq 0$, 
\[
\abs{U (t)} 
\ \leq \ 
56,424 
\hsy
E (y) 
\hsy 
e^{-3y}
\hsy
y^3.
\]

PROOF \
Start from the inequality
\[
\abs{U (t)} 
\ \leq \ 
\abs{2 \hsy a^\prime (t)  \hsy \Psi^\prime (t)} 
\hsx + \hsx 
\abs{a^{\prime\prime} (t) \hsy \Psi (t)}
\hsx + \hsx 
\abs{\Phi (t) \hsy \Psi^{\prime\prime} (t)}
\]
and estimate separately each of the three summands.
\\[-.25cm]

\qquad \textbullet \quad
\allowdisplaybreaks\begin{align*}
\abs{2 \hsy a^\prime (t)  \hsy \Psi^\prime (t)} \ 
&\leq 
\abs{2 \hsy (- e^t \hsy e^{-y} \hsy y \hsy (15 - 30y + 8 y^2)}
\hsx \cdot \hsx 
\abs{565 \hsy e^t \hsy y^3 \hsy e^{-4 \hsy y}}
\\[11pt]
&\leq 
E (y) \hsy A (y),
\end{align*}
where
\[
A (y)
\ = \ 
1,130 \hsy e^{-3y}
\hsy
(15 y + 30 y ^2+ 8 y^3).
\]
\\[-1.25cm]

\qquad \textbullet \quad
\allowdisplaybreaks\begin{align*}
&
\abs{a^{\prime\prime} (t) \hsy \Psi (t)} \ 
\\[11pt]
&
\hspace{1cm}
\leq \
\abs{e^t \hsy e^{-y} \hsy y(-75 + 330 y - 224 y^2 + 32 y^3)}
\hsx \cdot \hsx 
\abs{64 \hsy e^t \hsy y^2 \hsy e^{-4 \hsy y}}
\\[11pt]
&
\hspace{1cm}
\leq \
E (y) \hsy B (y),
\end{align*}
where
\[
B (y)
\ = \ 
64 \hsy e^{-3y}
\hsy
(75 + 330 y + 224 y^2 + 32 y^3).
\]

\qquad \textbullet \quad
\allowdisplaybreaks\begin{align*}
\abs{\Phi (t) \hsy \Psi^{\prime\prime} (t)} \ 
&\leq \
\abs
{\frac{203}{202} \hsx e^t \hsy e^{-y} \hsy y (2y - 3)}
\hsx \cdot \hsx 
\abs{(1.031) \hsy 2^{13} \hsy e^t \hsy y^4 \hsy e^{-4 y}}
\\[11pt]
&\leq \
E (y) \hsy C (y),
\end{align*}
where
\[
C (y)
\ = \ 
8, 562 \hsy e^{-3 y} \hsy (2 y^3 + 3 y^2).
\]
Combining these estimates then gives
\allowdisplaybreaks\begin{align*}
\abs{U (t)} \ 
&\leq \
E (y) \hsy (A (y) \hsx + \hsx  B (y) \hsx + \hsx  C (y))
\\[11pt]
&\leq \
E (y) \hsy 2 e^{- 3 y}
\hsy
(2,400 + 19,035 y + 36,961 y^2 + 14,206 y^3)
\\[11pt]
&\leq \
E (y) \hsy 2 e^{- 3 y}\hsy (14,206 y^3)
\hsx \cdot \hsx 
\frac{2,400 + 19,035 y + 36,961 y^2 + 14,206 y^3}{14,206 y^3}
\\[11pt]
&\leq \
E (y) \hsy 2 e^{- 3 y} \hsy (14,206 y^3) \hsy (1.97)
\\[11pt]
&\leq \
56,424 \hsy E (y)\hsy e^{- 3 y} \hsy y^3.
\end{align*}

Recall now the statement of the theorem: \ 
$\forall \ t \geq 0$, 

\[
(\Phi^\prime (t))^2 
\hsx - \hsx 
\Phi (t) \hsy \Phi^{\prime\prime} (t)
\ > \ 
0.
\]

Proof: 
In fact, 
\allowdisplaybreaks\begin{align*}
V (t) + U (t) \ 
&\geq V (t) - \abs{U (t)}
\\[11pt]
&\geq\ 
256 \hsy E (y) - 56,424 \hsy E (y) \hsy e^{- 3 y} \hsy y^3
\\[11pt]
&\geq\ 
E (y) \hsy (256  - 56,424 \hsy e^{- 3 \pi} \hsy \pi^3)
\\[11pt]
&>\ 
114 \hsy E (y)
\\[11pt]
&>\ 
0.
\end{align*}


\chapter{
$\boldsymbol{\S}$\textbf{43}.\quad  POSITIVE QUADRATIC FORMS}
\setlength\parindent{2em}
\setcounter{theoremn}{0}
\renewcommand{\thepage}{\S43-\arabic{page}}

\qquad
Let $p \not\equiv 0$ be a real polynomial of degree $n \geq 1$: 
\\[-.5cm]
\[
p (z) 
\ = \ 
a_0 + a_1 z + \cdots + a_n z^n
\qquad (a_0 \neq 0).
\]
Let $z_1, \ldots, z_n$ be its zeros and put
\\[-.5cm]
\[
S_0 \ = \ n, 
\quad 
S_k \ = \ z_1^k + z_2^k + \cdots + z_n^k
\qquad (k = 1, 2, \ldots).
\]
\\[-1.5cm]

\begin{x}{\small\bf LEMMA} \ 
There is an expansion 
\[
z \hsx \frac{p^\prime (z)}{p (z)} 
\ = \ 
\sum\limits_{k = 0}^\infty \ 
S_k \hsy z^{-k}
\ = \ 
S_0 +\frac{S_1}{z}  + \cdots \hsy.
\]
In addition, 
\[
\sum\limits_{k = 0}^m \ 
a_{n - k} \hsy S_{m - k} 
\ = \ 
(n - m) \hsy a_{n - m}
\]
if $m < n$ but vanishes if $m \geq n$.
\\[-.25cm]
\end{x}

\begin{x}{\small\bf BORCHARDT-HERMITE CRITERION} \ 
The zeros of $p$ are real iff the determinants
\[
\Delta_k 
\ = \ 
\begin{vmatrix*}[l]
\ 
S_0
&
S_1
&
\cdots
&
S_{k-1} \ 
\\[4pt]
\
S_1
&
S_2
&
\cdots
&
S_k \ 
\\[4pt]
\ \
\vdots
&
\
\vdots
&
&
\ \hsx
\vdots
\\[4pt]
\
S_{k-1} 
&
S_k 
&
\cdots
&
S_{2k-2} \text{\hsy}
\end{vmatrix*}
\qquad (k = 1,2, \ldots, n)
\]
are nonnegative.  
Moreover, the number of distinct zeros of $p$ is equal to the index $k$ of the last 
$\Delta_k \neq 0$ 
in the above sequence.
\\[-.5cm]

[Note: \ 
Spelled out
\[
\Delta_1 
\ = \ 
S_0 \hsy,
\quad
\Delta_2 
\ = \ 
\begin{vmatrix*}[l]
\ 
S_0
&
S_1\ 
\\[4pt]
\
S_1
&
S_2 \text{\hsy}
\end{vmatrix*}
, 
\ldots \hsy .]
\]
\\[-1.25cm]
\end{x}


\qquad
{\small\bf \un{N.B.}} \ 
If 
$\Delta_{k+1} = 0$, 
then 
$\Delta_{k+2} =  \cdots = \Delta_n = 0$.
\\[-.25cm]

\begin{x}{\small\bf EXAMPLE} \ 
Take $n = 2$ and consider 
$p (z) = z^2 - 1$ 
$-$then 
$S_0 = 2$,  
$S_1 = 1 + (-1) = 0$, 
$S_2 = 1^2 + (-1)^2 = 2$, 
hence
\[
\Delta_2 
\ = \ 
\begin{vmatrix*}[l]
\ 
2
&&
0\ 
\\[4pt]
\
0
&&
2\text{\hsy}
\end{vmatrix*}
\ = \ 4.
\]
\\[-1.cm]
\end{x}

\begin{x}{\small\bf EXAMPLE} \ 
Take $n = 2$ and consider 
$p (z) = z^2 + 1$ 
$-$then 
$S_0 = 2$,  
$S_1 = \sqrt{-1} \hsx + (-\sqrt{-1}) = 0$, 
$S_2 = (\sqrt{-1})^2 \hsx + (-\sqrt{-1})^2 = - 1 - 1 = - 2$, 
hence
\[
\Delta_2 
\ = \ 
\begin{vmatrix*}[l]
\ 
2
&
\ \
0\ 
\\[4pt]
\
0
&
-2\text{\hsy}
\end{vmatrix*}
\ = \ 
-4. 
\]
\\[-1.25cm]
\end{x}

\begin{x}{\small\bf EXAMPLE} \ 
Take $n = 2$ and consider 
$p (z) = (z - 1)^2$ 
$-$then 
$S_0 = 2$,  
$S_1 = 1 + 1$, 
$S_2 = 1^2 + 1^2 = 2$, 
hence
\[
\Delta_2 
\ = \ 
\begin{vmatrix*}[l]
\ 
2
&&
2\ 
\\[4pt]
\
2
&&
2\text{\hsy}
\end{vmatrix*}
\ = \ 0.
\]
\\[-1.25cm]
\end{x}

\begin{x}{\small\bf RAPPEL} \ 
Let 
$A = \{a_{i j}\}$ 
be a real symmetric matrix of degree $n$ $-$then the quadratic form $\un{A}$ associated with $A$ 
is the function of $n$ real variables $x_1, \ldots, x_n$ defined by 
\[
\un{A} (\un{x})
\ = \ 
\sum\limits_{i = 1}^n \ 
\sum\limits_{j =1}^n \ 
a_{i  j} \hsy  x_i \hsy x_j.
\]

\qquad \textbullet \quad
\un{$A$} is \un{positive} if $\forall \ \un{x} \neq 0$, 
\[
\un{A} (\un{x})
\ > \ 
0.
\]
\\[-1.25cm]
\end{x}

\qquad 
{\small\bf FACT} \ 
$\un{A}$  is positive iff all successive principal minors of $\un{A}$ are positive, i.e., 
\[
a_{1 1}
\ > \ 0, 
\quad
\begin{vmatrix*}[l]
\ 
a_{1 1} 
&
a_{1 2}\ 
\\[4pt]
a_{2 1} 
&
a_{2 2}\ 
\end{vmatrix*}
\ > \ 0, 
\quad
\cdots, 
\quad
\begin{vmatrix*}[l]
\ 
a_{1 1} 
&
\cdots
&
a_{1 n} \hspace{0.25cm} { }
\\[4pt]
\
\cdots\cdot
&
\hspace{-.2cm}
\cdots\cdot
&
\hspace{-.2cm}
\cdots \cdot
\\[4pt]
a_{n 1} 
&
\cdots
&
a_{n n}\ 
\end{vmatrix*}
\ > \ 0.
\]
\\[-1.25cm]

\begin{x}{\small\bf SCHOLIUM} \ 
The zeros of $p$ are real and simple iff the quadratic form
\[
\sum\limits_{i, j = 0}^{n-1} \ 
s_{i + j} \hsy  x_i \hsy x_j
\]
is positive.
\\[-.5cm]
\end{x}

Put 
\[
s_k 
\ = \ 
\frac{1}{z_1^k} 
+ 
\frac{1}{z_2^k} 
+ \cdots + 
\frac{1}{z_n^k} 
\qquad (k = 1, 2, \ldots).
\]
\\[-1.cm]

\begin{x}{\small\bf LEMMA} \ 
There is an expansion
\[
-\frac{p^\prime (z)}{p (z)} 
\ = \ 
s_1 + s_2 \hsy z + s_3 \hsy z^2 + \cdots \hsy.
\]
\\[-1.5cm]
\end{x}

\qquad
{\small\bf \un{N.B.}} \ 
This is the point of departure for the ensuing extension of the theory.  
\\[-.5cm]

[Note: \ 
By way of reconciliation, observe that
\allowdisplaybreaks
\begin{align*}
\frac{p (z)}{a_0}\ 
&=\ 
\Big(1 - \frac{z}{z_1}\Big) \cdots \Big(1 - \frac{z}{z_n}\Big) 
\\[15pt]
&=\ 
\scalebox{1.11}{$e^{- s_{_1}  z}$} \ 
\prod\limits_{k = 1}^n \ 
\Big(1 - \frac{z}{z_k}\Big)
\hsx
e^{a / z_k}, 
\end{align*}
so the ''$b$'' below is, in fact, $-s_1$.]
\\[-.25cm]

Let 
$f \not\equiv 0$ be a transcendental real entire function with an infinity of zeros such that 
$f (0) \neq 0$: 
\allowdisplaybreaks
\begin{align*}
f (z) \ 
&=\ 
\sum\limits_{n = 0}^\infty \ 
c_n \hsy z^n
\\[15pt]
&=\ 
\sum\limits_{n = 0}^\infty \ 
\frac{\gamma_n}{n !} \hsy z^n
\qquad (\gamma_n = f^{(n)} (0)).
\end{align*}

Assume further that 
$f \in \sL - \sP$ 
$-$then in view of 10.19, $f$ has a representation of the form
\[
f (z) 
\ = \ 
C \hsy e^{a z^2 + b\hsy z } \ 
\prod\limits_{n = 1}^\infty \ 
\Big(1 - \frac{z}{\lambda_n}\Big)
\hsx
e^{z / \lambda_n}, 
\]
where $C \neq 0$ is real, $a$ is real and $\leq 0$, $b$ is real, 
the $\lambda_n$ are real with 
$
\
\ds
\sum\limits_{n = 1}^\infty \ 
\frac{1}{\lambda_n^2} < \infty
$. 
\\[-.25cm]

Consider now the expansion 
\allowdisplaybreaks
\begin{align*}
-\frac{f^\prime (z)}{f (z)} \ 
&=\ 
 - 2 \hsy a \hsy z \hsx - \hsx b 
\ +  \ 
\sum\limits_{n = 1}^\infty \ 
\Big(
\frac{1}{\lambda_n^{- z}}  - \frac{1}{\lambda_n} 
\Big)
\\[15pt]
&=\ 
- b - 2 \hsy a \hsy z
 \ +  \ 
\sum\limits_{n = 1}^\infty \ 
\Big(
\frac{z}{\lambda_n^2}   
\hsx + \hsx
\frac{z^2}{\lambda_n^3} 
\hsx + \hsx 
\cdots
\Big)
\\[15pt]
&=\ 
s_1 + s_2 \hsy z + s_3 \hsy z^2 + \cdots,
\end{align*}
thus
\[
s_1 \ = \ -b, \quad
s_2 \ = \ - 2 \hsy a \hsx + \hsx 
\sum\limits_{n = 1}^\infty \ 
\frac{1}{\lambda_n^2} 
\]
and 
\[
s_k 
\ = \ 
\sum\limits_{n = 1}^\infty \ 
\frac{1}{\lambda_n^k} 
\qquad (k \geq 3).
\]
\\[-1cm]


\begin{x}{\small\bf THEOREM} \ 
$\forall \ r \geq 0$, 
the quadratic form 
\[
\sum\limits_{i, j = 0}^r \ 
s_{2 + i + j} \hsy  x_i \hsy x_j
\] 
is positive.  
\\[-.5cm]

PROOF 
Inserting the data, consider
\[
- 2 a \hsy x_0^2 
\ + \ 
\sum\limits_{n = 1}^\infty \ 
\bigg(
\sum\limits_{i, j = 0}^r \ 
\frac{x_i \hsy x_j}{\lambda_n^{2 + i + j}}
\bigg)
\]
or still, 
\[
- 2 a \hsy x_0^2 
\ + \ 
\sum\limits_{n = 1}^\infty \ 
\frac{1}{\lambda_n^2} 
\hsx
\Big(
x_0 + \frac{x_1}{\lambda_n} + \cdots + \frac{x_r}{\lambda_n^r}
\Big)^2, 
\]
an expression in which each term is manifestly nonnegative.  
Suppose that $\exists$ 
$x_0^{(0)}$,  $x_1^{(0)}$, $\ldots, x_r^{(0)}$
such that 
\[
\sum\limits_{i, j = 0}^r \ 
s_{2 + i + j} \hsy  x_i^{(0)} \hsy x_j^{(0)} 
\ = \ 
0.
\]
Let
\[
P_r (x) 
\ = \ 
x_0^{(0)} + x_1^{(0)} x + \cdots + x_r^{(0)} x^r.
\]
Then
\[
P_r \Big(\frac{1}{\lambda_n}\Big) 
\ = \ 
0
\qquad (n = 1, 2, \ldots).
\]
But the number of distinct
$
\hsy 
\ds
\frac{1}{\lambda_n}
$
is infinite implying, therefore, that 
$P_r \equiv 0$, 
hence
$x_0^{(0)} = 0$, 
$x_1^{(0)} = 0$, 
\ldots, 
$x_r^{(0)} = 0$.
\\[-.25cm]
\end{x}


\begin{x}{\small\bf SCHOLIUM} \ 
If $f \not\equiv 0$ 
is a transcendental real entire function with an infinity of zeros such that 
$f (0) \neq 0$ 
and if 
$f \in \sL - \sP$, 
then the determinants
\[
D_r \ \equiv \ 
\begin{vmatrix*}[l]
\ 
s_2 
&
s_3
&\cdots
&s_{2 + r}\ 
\\[4pt]
\
s_3
&
s_4
&\cdots
&s_{2 + r  + 1}\ 
\\[4pt]
\
\vdots
&
\vdots
&
&
\vdots
\\[4pt]
\
s_{2 + r}
&
s_{2 + r  + 1}
&\cdots
&s_{2 + r  + r}\text{\hsy} 
\end{vmatrix*}
\qquad (r \geq 0)
\]
are positive.
\\[-.25cm]
\end{x}

\begin{x}{\small\bf EXAMPLE} \ 
Take $r = 0$ $-$then
\[
D_0 
\ = \ s_2 
\ = \ 
- 2 \hsy a 
\hsx + \hsx
\sum\limits_{n = 1}^\infty \ 
\frac{1}{\lambda_n^2} 
\ > \ 0.
\]

[Note: \ 
Assume that $c_0 = 1$, $-$then from the theory
\[
- 2 \hsy a 
\ = \ 
c_1^2 
\hsx - \hsx
2 \hsy c_2 
\hsx - \hsx
\sum\limits_{n = 1}^\infty \ 
\frac{1}{\lambda_n^2} 
\]
or still, 
\[
- 2 \hsy a 
\hsx + \hsx
\sum\limits_{n = 1}^\infty \ 
\frac{1}{\lambda_n^2} 
\ = \ 
c_1^2 
\hsx - \hsx
2 \hsy c_2
\]
or still, 
\[
s_2 
\ = \ 
c_1^2 
\hsx - \hsx
2 \hsy c_2
\qquad (\tcf. \ 43.13).]
\]
\\[-1.25cm]
\end{x}

\begin{x}{\small\bf EXAMPLE} \ 
Take $r = 1$  $-$then
\[
D_1 \ = \ 
\begin{vmatrix*}[l]
\ 
s_2 
&
s_3\ 
\\[4pt]
\
s_3 
&
s_4\text{\hsy}
\end{vmatrix*}
\ > \ 
0.
\]
\\[-1.cm]
\end{x}


\begin{x}{\small\bf LEMMA} \ 
We have
\allowdisplaybreaks
\begin{align*}
&
c_0 \hsy s_1 + c_1 \ = \ 0
\\[11pt]
&
c_0 \hsy s_2 + c_1 \hsy s_1  + 2 \hsy c_2\ = \ 0
\\[11pt]
&
c_0 \hsy s_3 + c_1 \hsy s_2 + c_2 \hsy s_1 + 3 \hsy c_3 \ = \ 0
\\[11pt]
&
c_0 \hsy s_4 + c_1 \hsy s_3 + c_2 \hsy s_2 + c_3 \hsy s_1 + 4 \hsy c_4\ = \ 0
\\[11pt]
&
\bcdot \hspace{-.04cm} 
\bcdot \hspace{-.04cm} 
\bcdot \hspace{-.04cm} 
\bcdot \hspace{-.04cm} 
\bcdot \hspace{-.04cm} 
\bcdot \hspace{-.04cm} 
\bcdot \hspace{-.04cm} 
\bcdot \hspace{-.04cm} 
\bcdot \hspace{-.04cm} 
\bcdot \hspace{-.04cm} 
\bcdot \hspace{-.04cm} 
\bcdot \hspace{-.04cm} 
\bcdot \hspace{-.04cm} 
\bcdot \hspace{-.04cm} 
\bcdot \hspace{-.04cm} 
\bcdot \hspace{-.04cm} 
\bcdot \hspace{-.04cm} 
\bcdot \hspace{-.04cm} 
\end{align*}
\\[-1.25cm]
\end{x}

\begin{x}{\small\bf APPLICATION} \ 
Suppose that $c_0$ is positive and $f$ is even $-$then 
$c_1 = 0$, 
$c_3 = 0$, 
\ldots \ 
and 
$s_1 = 0$, 
$s_3 = 0$, 
\ldots \hsy.  
Therefore
\[
s_2 
\ = \ 
- \frac{2 c_2}{c_0} 
\ > \ 
0
\qquad (\implies c_2 < 0)
\]
while 
\[
c_0 \hsy s_4 \hsx + \hsx c_2 \hsy \Big(-\frac{2 \hsy c_2}{c_0} \Big) \hsx + \hsx 4 \hsy c_4 
\ = \ 
0
\]
\qquad
\qquad 
$\implies$
\[
c_0 \hsy s_4
\ = \ 
\frac{2 \hsy c_2^2}{c_0} \hsx - \hsx 4 \hsy a_4
\quad \implies \quad
\frac{c_2^2}{c_0} \hsx - \hsx 2 \hsy c_4 
\ > \ 0.
\]
\\[-1.cm]
\end{x}

\begin{x}{\small\bf EXAMPLE} \ 
In the notation of \S41, take
\allowdisplaybreaks
\begin{align*}
f (z) \ 
&=\ 
\Sh (z) 
\\[15pt]
&=\ 
\frac{1}{8} \ \Xi \Big(\frac{z}{2}\Big)
\\[15pt]
&=\ 
\sum\limits_{k = 0}^\infty \ 
\frac{(-1)^k}{(2 k)!}
\hsx 
b_k \hsy z^{2 k}.
\end{align*}
Then $\Sh$ is even and under RH, 
$\Sh \in \sL - \sP$, 
thus the positivity of the $D_r$ $(r \geq 0)$ 
provides a countable set of necessary conditions for its validity.  
To illustrate, in the case at hand
\[
c_0 \ = \ b_0, \quad
c_1 \ = \ 0, \quad
c_2\ = \ \frac{1}{2 !} \hsy b_1, \quad
c_3\ = \ 0, \quad
c_4 \ = \ \frac{1}{4 !} \hsy b_2.
\]
Accordingly, 

\begin{align*}
\frac{c_2^2}{c_0} - 2 \hsy c_4 \ 
&=\ 
\frac{1}{b_0} \hsy \Big(- \frac{1}{2}\hsx b_1\Big)^2 - \frac{2}{24} \hsx b_2
\\[15pt]
&=\ 
\frac{1}{4} \hsx \frac{b_1^2}{b_0} - \frac{1}{12} \hsx b_2
\\[15pt]
&=\ 
\frac{1}{4 \hsy b_0} \hsx\Big(b_1^2 - \frac{1}{3} \hsx b_0 \hsx b_2 \Big).
\end{align*}
And
\[
b_1^2 \hsy - \hsy \frac{1}{3} \hsx b_0 \hsx b_2
\ = \ 
3. \ 588 \ 449 \ 148 \ldots
\ > \ 
0.
\]
\\[-1.25cm]
\end{x}

The central conclusion thus far is 43.9: \ 
If 
$f \in \sL - \sP$, 
then 
$\forall \ r \geq 0$, 
the quadratic form 
\[
\sum\limits_{i, j = 0}^r \ 
s_{2 + i + j} \hsy  x_i \hsy x_j
\] 
is positive.  
But this can be turned around.
\\[-.25cm]

\begin{x}{\small\bf THEOREM}\footnote[2]{\vspace{.11 cm}
J. Grommer,  \textit{J. Reine Angew. Math.} \textbf{144} (1914), pp. 114-166;
see also
N. Kritikos, \textit{Math. Annalen} \textbf{81} (1920), pp. 97-118.}
\ 
Suppose that 
\[
f (z) 
\ = \ 
C \hsy e^{a z^2 + b} \ 
\prod\limits_{n = 1}^\infty \ 
\Big(1 - \frac{z}{z_n}\Big)
\hsx
e^{z / z_n}
\]
is in 
$A - \sL - \sP$ 
(cf. 10.31).  
Assume: \ 
$\forall \ r \geq 0$, 
the quadratic form 
\[
\sum\limits_{i, j = 0}^r \ 
s_{2 + i + j} \hsy  x_i \hsy x_j
\] 
is positive $-$then 
$f \in \sL - \sP$.
\\[-.25cm]
\end{x}

Since 
\[
\Sh \in 1 - \sL - \sP,
\]
one approach to RH is potentially through 43.16.


\chapter{
$\boldsymbol{\S}$\textbf{44}.\quad  ONE EQUIVALENCE}
\setlength\parindent{2em}
\setcounter{theoremn}{0}
\renewcommand{\thepage}{\S44-\arabic{page}}


\qquad
There are a number of statements which are equivalent to the Riemann Hypothesis.  
What follows is one of them (of a semi-trivial nature \ldots).
\\[-.25cm]

Per \S41, 
\[
\Sh (z) 
\ = \ 
\sum\limits_{k = 0}^\infty \ 
\frac{(-1)^k}{(2 k)!} 
\ 
b_k \hsx z^k,
\]
where
\[
b_k 
\ = \ 
\int\limits_0^\infty \ 
t^{2 k}
\hsx 
\Phi (t) 
\ \td t
\qquad (k = 0, 1, \ldots).
\]
In particular: 
\[
b_0 
\ = \ 
\int\limits_0^\infty \ 
\Phi (t) 
\ \td t, 
\quad 
b_1
\ = \ 
\int\limits_0^\infty \ 
t^2
\hsx
\Phi (t) 
\ \td t.
\]
\\[-.5cm]

Let 
$0 < x_1 \leq x_2 \leq \cdots $ 
be the positive zeros of 
$\Sh$. 
\\[-.5cm]

Let $S = \{\rho\}$ be the set of nonreal zeros of $\Sh$ whose imaginary part is positive: \ 
\[
\rho 
\ = \ 
\alpha + \sqrt{-1} \hsx \beta
\qquad (0 < \beta < 1).
\]
\\[-1.5cm]

[Note: \ 
A sum over the empty set is 0 and a product over the empty set is 1.]
\\[-.25cm]

\begin{x}{\small\bf LEMMA} \ 
\[
\Sh (z) 
\ = \ 
\Sh (0) \ 
\prod\limits_{n = 1}^\infty \ 
\Big(
1 - \frac{z^2}{x_n}
\Big)
\ 
\prod\limits_{\rho \in S} \  
\Big(
1 - \frac{z^2}{\rho^2}
\Big).
\]
\\[-1.25cm]
\end{x}

\begin{x}{\small\bf LEMMA} \ 
\allowdisplaybreaks
\begin{align*}
\frac{\td}{\td z} \ 
\Big(\frac{\Sh^\prime (z) }{\Sh (z) }\Big) \ 
&=\ 
-
\sum\limits_{n = 1}^\infty \ 
\bigg(\frac{1}{(z - x_n)^2}
\hsx  + \hsx 
\frac{1}{(z + x_n)^2}\bigg)
\\[15pt]
&\hspace{2.5cm}
- 
\sum\limits_{\rho \in S} \ 
\bigg(\frac{1}{(z - \rho)^2}
\hsx  + \hsx 
\frac{1}{(z + \rho)^2}\bigg).
\end{align*}
Now evaluate the left hand side of 44.2 at $z = 0$: 
\allowdisplaybreaks
\begin{align*}
\frac{\td}{\td z} \ 
\Big(\frac{\Sh^\prime (z) }{\Sh (z)}\Big)\Big|_{z = 0} \ 
&=\ 
\Big(\frac{\Sh^\prime}{\Sh}\Big)^{\raisebox{.07cm}{\text{\scriptsize$\prime$}}} \hsy (0)
\\[15pt]
&=\ 
\frac{ \Sh (0) \hsy \Sh^{\prime\prime} (0) - \Sh^\prime (0)^2}{\Sh (0)^2}
\\[15pt]
&=\ 
\frac{\Sh^{\prime\prime} (0)}{\Sh (0)}\hsx .
\end{align*}
And
\[
\begin{cases}
\ 
b_0 \ = \ 
\Sh (0)
\\[4pt]
\ 
b_1 \ = \ 
- 
\Sh^{\prime\prime} (0)
\end{cases}
.
\]

[Note: \ 
$\Sh^\prime (0) = 0$ $(\Sh$ being even).]
\\[-.5cm]

On the other hand, the right hand side of 44.2 evaluated at $z = 0$ is

\[
- 2 \ 
\sum\limits_{n = 1}^\infty \ 
x_n^2
\ - \  2 \ 
\sum\limits_{\rho \in S}^\infty \ 
\frac{1}{\rho^2}\hsx .
\]
And
\allowdisplaybreaks
\begin{align*}
\frac{1}{\rho^2} \ 
&=\ 
\frac{1}{\alpha^2 - \beta^2 + 2 \hsy \sqrt{-1} \hsx \alpha \hsy \beta}
\\[15pt]
&=\ 
\frac
{\alpha^2 - \beta^2 - 2 \hsy \sqrt{-1} \hsx \alpha \hsy \beta}
{(\alpha^2 - \beta^2)^2 + 4 \hsy  \alpha^2 \hsy \beta^2}
\\[15pt]
&=\ 
\frac
{\alpha^2 - \beta^2 - 2 \hsy \sqrt{-1} \hsx \alpha \hsy \beta}
{\alpha^4 + 2 \hsy  \alpha^2 \hsy \beta^2 + \beta^4}\hsx .
\end{align*}

[Note: \ 
Working with $-\bar{\rho} = -\alpha + \sqrt{-1} \hsx \beta$ leads to 
\[
\frac
{\alpha^2 - \beta^2 + 2 \hsy \sqrt{-1} \hsx \alpha \hsy \beta}
{\alpha^4 + 2 \hsy  \alpha^2 \hsy \beta^2 + \beta^4},
\]
hence when summed the imaginary parts cancel out.]
\\[-.5cm]

Therefore
\[
\frac{b_1}{2 \hsy b_0}
\ = \ 
\sum\limits_{n = 1}^\infty \ 
\frac{1}{x_n^2} 
\ + \ 
\sum\limits_{\rho \in S} \ 
\frac
{\alpha^2 - \beta^2}
{\alpha^4 + 2 \hsy  \alpha^2 \hsy \beta^2 + \beta^4}
\hsx .
\]
\\[-.25cm]
\end{x}

\qquad
{\small\bf \un{N.B.}} \ 
$\forall \ \rho \in S$: 
\[
\begin{cases}
\ 
1 < \abs{\alpha}
\\[11pt]
\ 
0 < \beta < 1
\end{cases}
\implies 
\alpha^2 - \beta^2 > 0 \hsx .
\]
\\[-.25cm]

\begin{x}{\small\bf THEOREM} \ 
RH holds iff 
\[
\sum\limits_{n = 1}^\infty \ 
\frac{1}{x_n^2} 
\ = \ 
\frac{b_1}{2 \hsy b_0}.
\]

[The point is that if $S$ is not empty, then $\forall \ \rho \in S$, 
$\alpha^2 - \beta^2 > 0 \hsx$.]
\\[-.25cm]
\end{x}


\chapter{
$\boldsymbol{\S}$\textbf{45}.\quad  SUGGESTED READING}
\setlength\parindent{2em}
\setcounter{theoremn}{0}
\renewcommand{\thepage}{\S45-\arabic{page}}

\vspace{-.5cm}
\begin{rf}
Bhaskar Bagchi, On Nyman, Beurling and Baez-Duarte's Hilbert Space Reformulation of the Riemann Hypothesis, 
\textit{Proc. Indian Acad. Sci. $($Math. Sci.$)$}, 
\textbf{116} (2003), pp. 137-146.
\end{rf}

\begin{rf}
Michel Balazard, Completeness Problems and the Riemann Hypothesis: 
An Annotated Bibliography, In: 
\textit{Surveys in Number Theory}, 
A. K. Peters Ltd. (2003), pp. 1-28.
\end{rf}

\begin{rf}
Michel Balazard, 
Un Si\`ecle et Demi de Recherches sur L'Hypoth\`ese de Riemann, 
\textit{Gazette des Math\'ematiciens}  
\textbf{126} (2010), pp. 7-24.
\end{rf}

\begin{rf}
N. G. de Bruijn, The Roots of Trigonometric Integrals, 
\textit{Duke Math. J.} 
\textbf{17} (1950), pp. 197-226.
\end{rf}

\begin{rf}
Brian Conrey, The Riemann Hypothesis, 
\textit{Notices AMS}  
\textbf{50} (2003), pp. 341-353.
\end{rf}

\begin{rf}
George Csordas, Linear Operators, Fourier Transforms, and the Riemann $\xi$-Function, 
In: \ 
\textit{Some Topics on the Value Distribution and Differentiability in Complex and $p$-Adic Analysis}, 
Beijing Science Press (2008), pp. 188-218.
\end{rf}

\begin{rf}
Haseo Ki, The Zeros of Fourier Transforms, 
In: \ 
\textit{Fourier Series Methods in Complex Analysis}, 
Univ. Joensuu Dept. Math. Rep. Ser. 
\textbf{10}  (2006), pp. 113-127.
\end{rf}

\begin{rf}
M. G. Krein, Concerning a Special Class of Entire and Meromorphic Functions, 
In: \ 
\textit{Some Questions in the Theory of Moments}, 
AMS (1962), pp. 214-265.
\end{rf}


\begin{rf}
M. Marden, 
\textit{Geometry of Polynomials}, Math. Surveys no. 3, AMS, Providence, RI,
2nd ed., 1966.
\end{rf}

\begin{rf}
George P\'olya, 
\"Uber die Algebraisch-Funktionentheoretischen Untersuchungen von J. L. W. V. Jensen, 
\textit{Kgl. Danske Vid. Sel. Math.-Fys. Medd.} 
\textbf{7}  (1927), pp. 3-33.
\end{rf}

\begin{rf}
Richard S. Varga, Theoretical and Computational Aspects of the Riemann Hypothesis, 
In: \ 
\textit{Scientific Computation on Mathematical Problems and Conjectures}, 
SIAM (1990), pp. 39-63.
\end{rf}


\newpage
\setcounter{page}{1}
\renewcommand{\thepage}{Index-\arabic{page}}
\printindex
\end{document}